\documentclass{amsart}

\usepackage[all]{xy}
\usepackage[a4paper, top=2cm, bottom=2cm, left=3cm, right=3cm, heightrounded, marginparwidth=3.5cm, marginparsep=.2cm, centering]{geometry}

\usepackage{hyperref}
\usepackage[capitalise]{cleveref}
\usepackage{calligra,mathrsfs}

\usepackage{tikz-cd}
\usepackage{bbold}
\usepackage{amssymb}
\usepackage{amsaddr}
\usepackage{todonotes}
\usepackage[T1]{fontenc}

\usepackage{relsize}
\usepackage[bbgreekl]{mathbbol}
\usepackage{amsfonts,amsthm}
\usepackage{stmaryrd}

\newcommand{\n}[1]{\mathcal{#1}}
\newcommand{\bb}[1]{\mathbb{#1}}

\newcommand{\f}[1]{\mathfrak{#1}}
\newcommand{\s}[1]{\mathscr{#1}}
\newcommand{\cat}[1]{\mathsf{#1}}
\newcommand{\ob}[1]{\mathrm{#1}}

\newcommand{\proet}{\mathrm{pro\acute{e}t}}

\DeclareSymbolFontAlphabet{\mathbb}{AMSb}
\DeclareSymbolFontAlphabet{\mathbbl}{bbold}
\newcommand{\prism}{{\mathlarger{\mathbbl{\Delta}}}}

\numberwithin{equation}{subsection}

\DeclareRobustCommand{\SkipTocEntry}[5]{}

\theoremstyle{plain}
\newtheorem{theorem}{Theorem}[subsection]
\newtheorem{thm}{Theorem}[section]
\newtheorem{proposition}[theorem]{Proposition}
\newtheorem{prop}[thm]{Proposition}
\newtheorem{lemma}[theorem]{Lemma}
\newtheorem*{claim*}{Claim}
\newtheorem{corollary}[theorem]{Corollary}

\newtheorem{assumption}[theorem]{Assumption}

\theoremstyle{definition}
\newtheorem{construction}[theorem]{Construction}

\newtheorem{definition}[theorem]{Definition}
\newtheorem{defi}[thm]{Definition}
\newtheorem{example}[theorem]{Example}
\newtheorem{remark}[theorem]{Remark}
\newtheorem{warning}[theorem]{Warning}
\newtheorem{rem}[thm]{Remark}

\newtheorem{notation}[theorem]{Notation}

\newcommand{\A}{\mathbb{A}}
\newcommand{\B}{\mathbb{B}}
\newcommand{\C}{\mathbb{C}}

\newcommand{\F}{\mathbb{F}}
\newcommand{\N}{\mathbb{N}}
\newcommand{\Q}{\mathbb{Q}}
\newcommand{\R}{\mathbb{R}}
\newcommand{\Z}{\mathbb{Z}}

\newcommand{\HK}{\mathrm{HK}}

\newcommand{\Der}{\mathrm{Der}}

\newcommand{\Lie}{\mathrm{Lie}}

\newcommand{\lc}{\mathrm{lc}}

\newcommand{\an}{\mathrm{an}}

\newcommand{\AnSpec}{\mathrm{AnSpec}}
\newcommand{\cone}{\mathrm{cone}}
\newcommand{\GSpec}{\mathrm{GSpec}}

\newcommand{\ArcStk}{\cat{ArcStk}}

\newcommand{\arc}{\mathrm{arc}}
\newcommand{\BC}{\mathcal{BC}}

\newcommand{\cont}{\mathrm{cont}}
\newcommand{\Betti}{\mathrm{Betti}}
\newcommand{\BK}{\mathrm{BK}}
\newcommand{\Div}{\mathrm{Div}}
\newcommand{\cDiv}{\mathcal{D}\mathrm{iv}}
\newcommand{\real}{\mathrm{real}}

\newcommand{\GelfStk}{\cat{GelfStk}}

\newcommand{\cycl}{\mathrm{cycl}}
\newcommand{\cyc}{\mathrm{cyc}}

\newcommand{\dR}{\mathrm{dR}}

\newcommand{\et}{\mathrm{{\acute{e}t}}}

\newcommand{\VB}{\mathrm{VB}}
\newcommand{\fd}{\mathrm{qfd}}
\newcommand{\lfd}{\mathrm{qfd}}

\newcommand{\FF}{\mathrm{FF}}
\newcommand{\FFdisc}[1]{\mathbb{D}_{\FF,#1}}

\newcommand{\Tate}{\mathrm{Tate}}

\newcommand{\Fil}{\mathrm{Fil}}

\newcommand{\Gal}{\mathrm{Gal}}
\newcommand{\Gm}{\mathbb{G}_m}
\newcommand{\Ga}{\mathbb{G}_a}

\newcommand{\gr}{\mathrm{gr}}
\newcommand{\GL}{\mathrm{GL}}

\newcommand{\Hom}{\mathrm{Hom}}

\newcommand{\HT}{\mathrm{HT}}

\newcommand{\Id}{\mathrm{Id}}

\renewcommand{\inf}{\mathrm{inf}}
\newcommand{\id}{\operatorname{id}}

\newcommand{\Mod}{\mathrm{Mod}}
\newcommand{\Marc}{\mathcal{M}_{\ob{arc}}}
\newcommand{\la}{\mathrm{la}}

\newcommand{\op}{\mathrm{op}}
\newcommand{\Perf}{\mathcal{P}\mathrm{erf}}

\newcommand{\perf}{\mathrm{perf}}

\newcommand{\red}{\mathrm{red}}

\newcommand{\sm}{\mathrm{sm}}
\newcommand{\Spa}{\mathrm{Spa}}
\newcommand{\Spd}{\mathrm{Spd}}
\newcommand{\Spf}{\mathrm{Spf}}
\newcommand{\Spec}{\mathrm{Spec}}

\newcommand{\Sym}{\mathrm{Sym}}

\newcommand{\site}{\mathrm{site}}

\newcommand{\colim}{\mathrm{colim}~}
\newcommand{\Aff}{\mathrm{Aff}}

\newcommand{\qfd}{\mathrm{qfd}}
\newcommand{\Nil}{\mathrm{Nil}}

\newcommand{\sol}{{\scalebox{0.5}{$\square$}}}

\newcommand{\xto}[1]{\overset{#1}{\to}}

\newcommand*\isomarrow{%
	\xrightarrow{\raisebox{-0.35em}{\smash{\ensuremath{\sim}}}}
}

\begin{document}

	\title[Rational analytic prismatization]{Rational analytic prismatization}
	
	 \author[Ansch\"utz, Le Bras, Rodr\'iguez Camargo, Scholze]{Johannes Ansch\"utz, Arthur-C\'esar Le Bras, Juan Esteban Rodr\'iguez Camargo, Peter Scholze}

	\begin{abstract}
          We develop a general theory of analytic prismatization with $\Q_p$-coefficients in rigid analytic geometry over $\Q_p$.
          We prove an equivalence of perfect complexes on the perfect analytic prismatization modulo Frobenius with perfect complexes on the relative Fargues--Fontaine curve for any smooth partially proper rigid analytic variety.
          In particular, we explain how this equivalence  conceptualizes the  decompletion techniques by passing to locally analytic vectors which appear in (geometric) Sen theory. 
	\end{abstract}
	
	\maketitle
	\tableofcontents

\section{Introduction}\label{SectionIntroduction}

\textit{Prismatic cohomology} has been defined in \cite{bhatt_scholze_prisms_and_prismatic_cohomology} as a new cohomology theory for $p$-adic formal schemes, unifying the already known existing $p$-adic cohomology theories such as de Rham cohomology, crystalline cohomology (of the special fiber), $p$-adic \'etale cohomology (of the rigid generic fiber), etc.. Later, Bhatt--Lurie \cite{bhatt2022prismatization}, and independently Drinfeld \cite{drinfeld2020prismatization}, introduced the \textit{prismatization} $X^\prism$ of a $p$-adic formal scheme $X$ in such a way that the coherent cohomology of the stack $X^\prism$ recovers the prismatic cohomology of $X$. They also introduced various refinements of this construction, encoding extra structures on prismatic cohomology, such as the Nygaard filtration or the Frobenius.
A toy model and an inspiration for these developements is furnished by de Rham cohomology (in characteristic $0$): according to a celebrated result of Grothendieck \cite{grothendieck1968crystals}, it can be thought of as the cohomology of a site, the \textit{infinitesimal site}, and Simpson \cite{simpson1996hodge} has defined the \textit{de Rham space} (or \emph{stack}) of a scheme, with coherent cohomology given by de Rham cohomology. Moreover, he introduced a filtered variant of this construction, encoding the Hodge filtration on de Rham cohomology.

Our eventual goal is to provide an analogue of these works for (a suitable class of) analytic spaces over $\Z_p$, such as rigid spaces over a non-archimedean field in mixed characteristic $0$ and $p$,  or in equal characteristic $p$. In this paper, however, we only partially fulfill this goal: first, we restrict ourselves to the prismatization itself and not any of its refinements; second, we focus on the base change to $\Q_p$, \textit{both} in the geometry and in the coefficients, of what should ultimately be the analytic prismatization.
Before commenting on these limitations (see \Cref{rk-intro:limitations}), let us  describe more precisely what we achieve below.

\subsection{Rational $\delta$-spaces and rational analytic prismatization} \label{intro-sec-1}
The prismatization of \cite{bhatt2022prismatization} and \cite{drinfeld2020prismatization} is defined for $p$-adic formal schemes; it takes values in $p$-adic formal stacks, defined as functors on $p$-nilpotent rings. Since we aim at defining a version of the prismatization for a class of geometric objects of more analytic flavour, we need to use a different geometric setup. Such a geometric framework has been defined in \cite{drff_stacks}: Gelfand analytic geometry, i.e., the theory of Gelfand rings and Gelfand stacks. The category of \textit{Gelfand rings} is a subcategory of the category of solid $\Q_p$-algebras, which is stable under colimits. It contains Banach $\Q_p$-algebras and conversely any Gelfand ring $A$ admits a uniform completion $A^u$ which is a (uniform) Banach $\Q_p$-algebra. Let $\cat{GelfRing}$ be the category of Gelfand rings.  (For technical reasons, we always assume our Gelfand rings to be separable and quasi-finite dimensional, assumptions that will be implicit in the rest of this introduction.) When endowed with the descendable topology, a basis of $\cat{GelfRing}$ is formed by \textit{nilperfectoid rings}, which geometrically (on spectra) should be thought of as overconvergent thickenings of affinoid perfectoid spaces over $\Q_p$: a Gelfand ring $A$ is nilperfectoid if its quotient by its $\dagger$-nilradical is perfectoid. Finally, a \textit{Gelfand stack} is (in first approximation\footnote{The actual definition is more subtle, and involves a condition intermediate between the notion of sheaf and hypersheaf, best explained through the formalism of Gestalten.}) a sheaf of anima on $\cat{GelfRing}$, for the descendable topology.

The category $\cat{GelfStk}$ contains fully faithfully the category of $\dagger$-rigid spaces over $\Q_p$ and more generally the category of derived Berkovich spaces over $\Q_p$ from \cite{drff_stacks} (a derived generalization of Berkovich spaces\footnote{In the following, we will often drop the adjective ``derived''.}); as such it is a variant of the general notion of analytic stacks very well suited to our needs. 

\begin{rem}[Gelfand stacks and arc-stacks]
Any Gelfand stack gives rise to an arc-stack over $\Marc(\Q_p)$: the \textit{perfectoidization functor}
$$
(-)^{\diamond}: \cat{GelfStk} \to \cat{ArcStk}_{\Q_p}
$$
is a colimit-preserving functor sending $\GSpec(A)$ to $\Marc(A^u)$, for $A$ Gelfand. The reader can think of the perfectoidization of a Gelfand stack as analogous to taking the topological manifold underlying a complex manifold. The perfectoidization functor admits both left and right adjoints. Its right adjoint, the \textit{analytic de Rham stack functor} $(-)^\dR$, is especially interesting, and studied at length in \cite{RCAndeRham} (in a slightly different setup) and \cite{drff_stacks}. 

Although the main constructions below are modeled on the theory of prismatization of $p$-adic formal schemes and guided by the desire of ``analytifying'' them, the relation to arc-stacks provides a different perspective on Gelfand geometry, as an enhancement of perfectoid geometry allowing for non-trivial (overconvergent) thickenings. In fact, we will see that many of our definitions are extensions to nilperfectoid rings of already familiar constructions for perfectoid rings. 

As such, the theory of analytic prismatization reconciles two different and very successful approaches to $p$-adic Hodge theory: the theory of prismatic cohomology and the theory of diamonds. 
\end{rem}  

Like in the integral theory, we will both define the prismatic site of a Gelfand stack (e.g., a $\dagger$-rigid space over $\Q_p$) and its prismatization. This relies classically on the key definitions of \textit{$\delta$-rings} and \textit{prisms}. In our context, they take the following form. To make sense of the following definition, let $\cat{Perf}$ be the category of (qfd) perfectoid spaces in characteristic $p$ and recall that to any perfectoid space $S\in \cat{Perf}$, one can attach the Berkovich spaces\footnote{The notations $Y_S$, $Y_{S,[0,\infty)}, Y_{S,(0,\infty)}$ or $\mathcal{Y}_S$ are common, but due to their importance we use a different notation, and occasionally call them the (punctured) Fargues--Fontaine disc relative to $S$.}
$$
\FFdisc{S}, ~ \FFdisc{S}^{\times}= \FFdisc{S}\backslash V(p),
$$
where for $S=\Spa(A)$ affinoid perfectoid with pseudo-uniformizer $\varpi$, 
$$
\FFdisc{S}=\Spa(W(A^{\circ}))\backslash V([\varpi]).
$$

\begin{defi}[\Cref{Definition:Generalized-FF-Curvesv2}, \Cref{Definition-Analytic-prismatic-site}]
\begin{enumerate}
\item A \textit{rational  $\delta$-space} is a triple
$$
(S, \mathcal{X},\iota, \varphi_\mathcal{X}),
$$
 where $S\in \cat{Perf}$, $\mathcal{X}$ is a Berkovich space with an endomorphism $\varphi_\mathcal{X}$ (called ``the Frobenius of the $\delta$-space''), and $\iota\colon \FFdisc{S}^\times \to \mathcal{X}$ is a morphism of Berkovich spaces inducing an isomorphism of diamonds $\FFdisc{S}^{\times,\diamond}\xrightarrow{\sim} \mathcal{\mathcal{X}}^{\diamond}$ such that $\varphi_\mathcal{X}\circ \iota\cong \iota\circ \varphi$, where $\varphi\colon \FFdisc{S}^\times\to \FFdisc{S}^\times$ is the morphism induced by the Frobenius on $S$. We will often drop $\varphi_\mathcal{X}$ and $\iota$ from the notation. We let $\cat{Spc}^{\delta}$ denote the category of rational $\delta$-spaces. 
  
  \item  A morphism of rational  $\delta$-spaces $(S,\mathcal{X})\to (S',\mathcal{X}')$ is a morphism $\alpha\colon S \to S^\prime$ of arc-stacks, a $\varphi$-equivariant morphism $\mathcal{X}\to \mathcal{X}'$ of Gelfand stacks and a commutative diagram of Gelfand stacks
  \[
\begin{tikzcd}
\FFdisc{S}^\times \ar[r] \ar[d] & \mathcal{X} \ar[d] \\ 
\FFdisc{S'}^\times \ar[r] & \mathcal{X}^{\prime},
\end{tikzcd}
\]
where the left horizontal arrow $\FFdisc{S}^\times \to \FFdisc{S'}^\times$ is induced by $\alpha$.
\item A \textit{perfect rational $\delta$-space} is a rational  $\delta$-space $(S,\mathcal{X})$ such that $\varphi\colon \mathcal{X}\to \mathcal{X}$ is an automorphism.
  \item Let $(S,\mathcal{X})$ be a rational $\delta$-space. A (derived) Cartier divisor $(s\colon \mathcal{L}\to \s{O}_{\mathcal{X}})$ on $\mathcal{X}$ is of degree $1$ if its pullback to $\FFdisc{S}^\times$ along the natural map $\FFdisc{S}^\times\to \mathcal{X}$ is of degree $1$ (in the sense of \cite{fargues2021geometrization}: it represents the closed Cartier divisior of an untilt of $S$). A \textit{rational analytic prism} $(S,\mathcal{X},D)$ is a rational $\delta$-space $(S,\mathcal{X})$ together with a Cartier divisor $D\subseteq \mathcal{X}$ of degree $1$.\footnote{We write $D\subseteq \mathcal{X}$ to simplifiy the exposition, but we caution the reader that the morphism $D\to \mathcal{X}$ is not a monomorphism in general.}
\end{enumerate}
\end{defi}

\begin{rem}
The reader should think of rational $\delta$-spaces as analogous to (formal spectra of) $\delta$-rings and to rational analytic prisms as analogous to prisms. In the above, as well as in the very title of this paper, the adjective ``rational'' is here to recall that there exists an integral version of the above definitions, of which they are the base change to $\Q_p$. We will however from now on omit the adjectives ``rational'' and ``analytic'', to keep the terminology light. 
\end{rem}

The moduli problems we will consider will naturally be defined as sheaves on the category of $\delta$-spaces. We let 
$$\cat{AnStk}^\delta$$ 
be the localization of $\mathcal{P}(\cat{Spc}^\delta)$ at the $!$-equivalences (defined on underlying Berkovich spaces), and call objects therein \emph{rational analytic $\delta$-stacks} (or just \textit{$\delta$-stacks}). Any  $\delta$-stack can be realized as a Gelfand stack: the \emph{realization functor} on $\delta$-stacks is the colimit preserving functor
$$
(-)_\real\colon \cat{AnStk}^\delta\to \GelfStk,
$$
 sending $(S,\mathcal{X})$ to $\mathcal{X}$.

The next proposition makes the realization functor more explicit.

\begin{prop}[\Cref{Prop:Universal-FFNilPerfd2}, \Cref{PropAdjunctionFFStacks2}]
  Let $A$ be a nilperfectoid ring.
Let $T=\GSpec(A)$ and consider the category
$$
\cat{Spc}^\delta_{T/-}
$$
of  $\delta$-spaces $(S,\mathcal{X})$ with a morphism $T\to \mathcal{X}$ of Berkovich spaces.
Then $\cat{Spc}^\delta_{T/-}$ admits an initial object given by the pair $(T^\flat,\mathcal{X}_T)$, called the $\delta$-space associated to $T$, described as follows:
  \begin{enumerate}
 
  \item The perfectoid space $T^{\flat}$ over $\F_p$  is the tilt of $T^u:=\GSpec(A^u)$.
    In particular, we have a natural map of Berkovich spaces 
\begin{equation*}
f\colon T^u\to \FFdisc{T^\flat}^\times,
\end{equation*}
induced by the untilt $T^u$ of $T^\flat$.
We let $f_{\varphi}\colon \N\times T^u\to \FFdisc{T^\flat}^\times$ be the unique morphism whose restriction to $\{n\}\times T^u$ is $\varphi^n\circ f$.

\item The thickening  $\mathcal{\mathcal{X}}_T$ of $\FFdisc{T^\flat}^\times$ is the pushout in Berkovich spaces of the morphisms $f_\varphi\colon \N\times T^u\to \FFdisc{T^\flat}^\times$ and $\N\times T^u\to \N\times T$. 
\end{enumerate}
The same result holds for perfect $\delta$-spaces, replacing $\N$ by $\Z$ in the above construction. The quotient of the perfect $\delta$-space associated to $T$ by Frobenius will be denoted $\FF_T$ (this is the usual Fargues--Fontaine curve of the tilt when $T$ is perfectoid). 

Consequently, the realization functor $(-)_{\real}\colon \cat{AnStk}^\delta\to \GelfStk$ on  $\delta$-stacks sends an object $Y\in \cat{AnStk}^\delta$ to the sheafification in Gelfand stacks (i.e., with respect to $!$-equivalences) of the presheaf sending a $T=\GSpec(A)$ with $A$ nilperfectoid to $Y(T^\flat, \mathcal{X}_T)$. 
\end{prop}

  \begin{rem}[Witt vectors]
    \label{sec:thick-farg-font-relation-to-witt-vectors2-intro}
     Pursuing our analogy with the $p$-complete theory, the functor $T \mapsto (T^\flat, \mathcal{X}_T)$ should be thought of as analogous to the Witt vectors functor $R\mapsto W(R)$ on $p$-nilpotent rings. In fact, this can be made more precise. Let $A$ be a nilperfectoid $\Q_p$-algebra. We can describe $\mathcal{X}_T$, $T=\GSpec(A)$, as follows: let\footnote{For each light profinite set $S$, we have that $W(A^{\leq 1})(S)=W(A^{\leq 1}(S))$ are the Witt vectors of the animated ring $A^{\leq 1}(S)$, in the sense of \cite[Appendix A]{bhatt2022prismatization}.}
     $$
     W'(A^{\leq 1}):= W(A^{\leq 1}) \times_{W(A^{u,\circ})}  \varprojlim_{F} W(A^{u,\circ }) = W(A^{\leq 1}) \times_{W(A^{u,\circ})}  W(A^{u,\circ, \flat}).
     $$
    Then $\mathcal{X}_T$ is the Gelfand stack sending a Gelfand ring $B$ to the anima of maps $W'(A^{\leq 1})\to B^{\leq 1}$, such that its composition to any perfectoid $B$-algebra $C$ factors over the subspace $\FFdisc{A^u}^\times$.
      Indeed, this follows from the explicit description of $\mathcal{X}_T$ as a pushout and the fact that we have a distinguished triangle
    \[
W(\Nil^\dagger(A))\to W(A^{\leq 1})\to W(A^{u,\leq 1})= W(A^{u,\circ})
\]
(since $A$ is nilperfectoid) and that $W(\Nil^\dagger(A))\cong \prod\limits_{n\in \N} \Nil^\dagger(A)$ because the non-unital ring $\Nil^\dagger(A)$ is a $\Q_p$-algebra. This last statement critically uses that $A$ is a $\Q_p$-algebra and explains why with $\Q_p$-coefficients, the theory gets much simpler. For instance,  with $\Q_p$-coefficients the difference between ``analytic prismatization'' and ``perfect analytic prismatization'' is rather mild, and boils down to changing $\N$ to $\Z$ in the above description of Witt vectors of the $\dagger$-nilradical.

This observation also shows that in fact the theory of analytic prismatization is really analogous not to the theory of prismatization, but to the theory of \emph{sheared prismatization} currently developed by Bhatt--Kanaev--Mathew--Vologodsky--Zhang, cf., \cite{drinfeld2025ring}. The analogy is obscured by our restriction to $\Q_p$-coefficients, but would be cleaner if we had worked integrally. For the analytic prismatization with $\Z_p$-coefficients, the use of Witt vectors is critical, but we can and will avoid it in this paper.
  \end{rem}

 We can finally come to the main definition. 

\begin{defi}[\Cref{Definition-Analytic-prismatic-site}]
  Let $X$ be a Gelfand stack.
  \begin{enumerate}
  \item The \textit{rational analytic prismatic site of $X$}, denoted as $X_{\prism-\site}$, is the category of tuples $(S,\mathcal{X}, D,\alpha)$ where $(S,\mathcal{X},D)$ is a rational analytic prism and $\alpha\colon D\to X$ is a morphism of Gelfand stacks.
  Coverings are induced by $!$-covers in $(S,\mathcal{X})$.
   \item Let $\Q_{p,\prism}$ be the $\delta$-stack sending $(S,\mathcal{X})$ to the anima of degree $1$ Cartier divisors on $(S,\mathcal{X})$.
  \item Let $X_{\prism}$ be the $\delta$-stack over $\Q_{p,\prism}$ mapping $(S,\mathcal{X}, D)$ to the anima of maps $D\to X$.
 \item If $X$ is a Berkovich space, the \emph{rational analytic prismatization $X^\prism$ of $X$} is the realization of the  $\delta$-stack $X_{\prism}$.    Equivalently, it is the functor sending a nilperfectoid $\Q_p$-algebra $A$ over $\Q_p^\prism:=\GSpec(\Q_p)^\prism$ to the anima of maps $D\to X$, where $D\subseteq \mathcal{X}_{\GSpec(A)}$ is the degree $1$ divisor classified by the map $\GSpec(A)\to \Q_p^\prism$.

\item The \textit{rational perfect analytic prismatization $X^{\widetilde{\Div^1}}$ of $X$} is defined as $X^{\widetilde{\Div^1}}:=\varprojlim_{\varphi} X^{\prism}$.  We denote  $\widetilde{\Div^1}:= \varprojlim_{\varphi} \bb{Q}_p^{\prism}$.
 
  \item The \emph{rational perfect analytic prismatization of $X$ modulo Frobenius}, or alternatively the \textit{$p$-adic twistor stack of $X$}, is the Gelfand stack
    \[
      X^{\Div^1}:=X^{\widetilde{\Div^1}}/\varphi^\Z.
    \]
    We denote $\Div^1:=\widetilde{\Div^1}/\varphi$.
  \end{enumerate}  
\end{defi}

As before, we will drop the adjective ``rational'' in what follows. 

\begin{rem} The restriction to Berkovich spaces in the definition of the prismatization is here to ensure that the functor of points description on nilperfectoid rings holds as written, without any sheafification required (this is related to the considerations of \Cref{sec:thick-farg-font-relation-to-witt-vectors2-intro}). More generally, one could allow $X$ to be a $\dagger$-cohesive Gelfand stack  in the sense of \Cref{sec:push-deriv-berk-definition-cohesiveness}. In this introduction, we will restrict to Berkovich spaces, which is the case of interest for applications.
\end{rem}

For any Berkovich space $X$, there is a natural morphism 
$$
X^\prism \to \A^1/\Gm
$$
coming from pullback of Cartier divisors along $\{0\}\times T \to \N \times T \to \mathcal{X}_T$ for $T=\GSpec(A)$ nilperfectoid. This allows us to single out an interesting Cartier divisor of $X^\prism$.

\begin{defi}
Let $X$ be a Berkovich space.  The \emph{analytic Hodge-Tate stack of $X$} is the Gelfand stack
    \[
      X^\HT:=X^\prism\times_{\A^1/\Gm} B\Gm.
    \]
\end{defi}

\begin{rem}\label{rk-intro:explicit-description}
For any prism $(S,\mathcal{X},D)$, one has from the definition  natural maps of Gelfand stacks
$$
\mathcal{X} \to  D^{\prism} \to \Q_p^\prism
$$
mapping $D$ to $D^{\HT}\to \Q_p^\HT$, and where the first map corresponds to the identity of $D$, cf. \Cref{rk:explicit-description}.
For example, we have the morphism
$$
\mathcal{X}_{\rm cyc} \to \Q_p^\prism
$$
(factoring through $\widetilde{\Div^1}$) from the cyclotomic prism given by
$$(S,\mathcal{X}_{\rm cyc}):=(\Marc(\F_p((\pi^{1/p^{\infty}}))), \underset{q\mapsto q^p} \lim (1+ \mathbb{D}^{<1}) \backslash \{1\} ), $$
with the divisor corresponding to $[p]_q=(q^p-1)/q-1$, with $q$ the coordinate on the open unit disc $\mathbb{D}^{<1}$ over $\Q_p$ (and Frobenius given by $q\mapsto q^p$).
It is an analogue in the analytic theory of the $q$-de Rham prism.
\end{rem}

\begin{rem}[Relative versions]
One can define relative versions of the prismatization (relative to a fixed base prism), exactly like in the $p$-complete theory: let $(S,\mathcal{X},D)$ be a prism, $X$ a Berkovich space over $D$. Define
$$
X^{\prism/\mathcal{X}} = X^\prism \times_{D^\prism} \mathcal{X}
$$
(the left hand side really depends on $(S,\mathcal{X},D)$, not just on $\mathcal{X}$, but we follow an abuse of notation common in the $p$-complete theory to make the notation shorter), using the map $\mathcal{X}\to D^\prism$ from \Cref{rk-intro:explicit-description}. 
\end{rem}

\begin{rem}[Perfectoid case]
A difference with the $p$-complete theory is that the prismatization of a perfectoid space $X$ over $\Q_p$ is not simply given by $\FFdisc{X^\flat}^\times$. However, it is true that the Hodge--Tate stack $X^\HT$ of a perfectoid space $X$ is isomorphic to $X$: see \Cref{prop:ht-stack-of-perfectoids}. 
\end{rem}

In \Cref{Section:Analytic-Prism} and \Cref{sec:analyt-prism-ga}, we make explicit the (perfect) prismatizations of $\GSpec(\Q_p)$, the rigid analytic multiplicative group $\Gm$ and the analytic affine line. For $\GSpec(\Q_p)$, we show that the morphism $\mathcal{X}_{\rm cyc} \to \widetilde{\Div^1}$ from \Cref{rk-intro:explicit-description} induces isomorphisms\footnote{Given a $p$-adic Lie group $G$, we let $G^{\rm cont}$, $G^{\la}$ and $G^{\sm}$ denote the groups in Berkovich spaces given by the sheaf of continuous, locally analytic and locally constant functions on $G$ valued in $\Q_p$.}
$$
\widetilde{\Div^1} \cong \mathcal{X}_{\cyc}/\Z_p^{\times, \la}, ~~  \Div^1 \cong \mathcal{X}_{\cyc}/\Q_p^{\times, \la},
$$
with action of $a\in \Z_p^{\times, \la}$ given by $a.(q,q^{1/p},\dots):=(q^a, q^{a/p},\dots)$.  In particular, we have an isomorphism\footnote{Finding a theory of analytic prismatization outputing this explicit quotient stack as the Hodge-Tate stack of $\GSpec(\Q_p)$ has been one desideratum of the authors since the beginning of this project, inspired by \cite{colmez1994resultat}, \cite{ramsey2009bsen} and \cite{anschutz2022v}.
}
  \[
    \Q_p^\HT\cong \GSpec(\Q_p(\mu_{p^\infty}))/\Z_p^{\times,\la},
  \]
  where the action of $\Z_p^{\times,\la}$ on $\GSpec(\Q_p(\mu_{p^{\infty}}))$ factors through the natural Galois action of $\Z_p^{\times,\sm}\cong \Gal(\Q_p(\mu_{p^{\infty}})/\Q_p)^{\sm}$ provided by the cyclotomic character.  

  For the rigid analytic multiplicative group $\Gm$, and then also the analytic affine line, the explicit presentations involve a discussion of the prismatic logarithm (which takes a simpler form than the one of \cite{bhatt2022prismatization} and \cite{drinfeld2020prismatization}, since we only work rationally).
.

\subsection{The $p$-adic twistor stack and (geometric) Sen theory}
For a Berkovich space $X$, the analytic de Rham stack of $X^\prism$ is simply $(\FFdisc{X^\diamond}^\times)^\dR$. In particular, there is a radius map
$$
\kappa: X^\prism \to (0,\infty)
$$
intertwining the Frobenius with the multiplication by $p$ map
(here the interval $(0,\infty)$ is seen as a Gelfand stack over $\Q_p$ through the Betti stack functor). In most of the text and in the rest of this introduction, we will focus on the $p$-adic twistor stack $X^{\Div^1}$. This is justified by the following observation.
\begin{itemize}
\item Over $(1,\infty)$ (via $\kappa$), the morphism $X^\prism \to (\FFdisc{X^\diamond}^\times)^\dR$ is an isomorphism, thus the Frobenius induces an isomorphism between $X^\prism_{(1,\infty)}$ and $X_{(p,\infty)}^\prism$, and the quotient by Frobenius on this locus is the Hyodo--Kato stack of $X$.
\item Over $(0,p)$, the morphism $X^{\widetilde{\Div^1}} \to X^\prism$ is an isomorphism, thus the Frobenius also induces an isomorphism between $X^\prism_{(0,1)}$ and $X^\prism_{(0,p)}$ and the quotient by Frobenius on this locus is isomorphic on $X^{\Div^1}$.
\end{itemize}
Since the theory of Hyodo--Kato stacks has been already developed in \cite{drff_stacks}, this means the remaining work is in understanding the properties of $X^{\Div^1}$, the $p$-adic twistor stack of $X$.

For a Berkovich space $X$, there is a diagram of Gelfand stacks
\[
\FF_{X^{\diamond}}\xrightarrow{\alpha_X} X^{\Div^1}\xrightarrow{\beta_X} X^{\Div^1,\dR}=X^{\HK},
\]
where the composition is the natural map $\FF_{X^{\diamond}} \to X^{\HK}=(\FF_{X^{\diamond}})^\dR$. Hence, $X^{\Div^1}$ sits in between the relative Fargues-Fontaine curve over $X^{\diamond}$ (incarnated as a Gelfand stack) and its analytic de Rham stack. Note that in this diagram, the left and right terms only depend on the perfectoidization $X^{\diamond}$ of $X$ (its ``underlying topological space''), but $X^{\Div^1}$ itself does not. 

The next two results, \Cref{sec:defin-first-prop-2-first-properties-of-analytic-prismatizationintro} and \Cref{thm:equivalence-perfect-complexesintro}, make use of the morphisms $\beta_X$ and $\alpha_X$ respectively.

\begin{prop}[\Cref{sec:defin-first-prop-2-first-properties-of-analytic-prismatization}] \label{sec:defin-first-prop-2-first-properties-of-analytic-prismatizationintro}
Let $X$ be a Berkovich space. 
\begin{enumerate}
   \item The morphism $\beta_X: X^{\Div^1}\to \FF_{X^\diamond}^{\dR}=X^\HK$ is an isomorphism after base change to $\Div^1\setminus \Q_p^\HT\cong \Div^{1,\dR}\setminus \GSpec(\Q_p)$. 
   \item There exists a natural morphism of Gelfand stacks
  \[
    \pi_X\colon X^\HT\to X\times \Q_p^\HT.
  \]
  If $X$ is a smooth partially proper rigid space over $\Q_p$ (or more generally, a formally smooth Gelfand stack over $\Q_p$), $\pi_X$ realizes $X^\HT$ as a gerbe over $X\times \Q_p^\HT$, banded by $\mathbb{T}_{X/\Q_p}^{\dagger}\{1\}$, where $\{1\}$ refers to the Breuil--Kisin twist on $\Q_p^\HT$ (cf. \Cref{sec:breuil-kisin-twist-2-definition-bk-twist} and \Cref{ConstructionBKtwist}).
 \end{enumerate} 
\end{prop}

Hence, generically, $X^{\Div^1}$ coincides with the Hyodo--Kato stack of $X$, but interesting new phenomena happen on (the overconvergent neighborhood of) the Hodge-Tate locus. 

\begin{rem}[Relation to the Faltings' extension]
In fact the Hodge--Tate stack, and even the overconvergent neighborhood of $X^\HT$ inside $X^{\Div^1}$, is controlled by a familiar object in $p$-adic Hodge theory: the \emph{Faltings' extension} of $X$. This result admits a relative variant, over a perfectoid base, when $X$ lifts to Fontaine's ring $B_\dR^+$. See \Cref{sec:faltings-extension} for precise statements. 
\end{rem}

\begin{rem}[Archimedean/non-archimedean analogy] \label{rk-intro:analogy-real}
In \cite{real_local_langlands}, an archimedean version of this story has been developed. Let us briefly sketch the instructive analogy, since it motivates our choice of terminology. First, \cite{real_local_langlands} defines a notion of nilperfectoid $\C$-algebra. In this context, the perfectoidization functor is really taking the ``underlying topological space'' (or, rather, condensed anima). Moreover, for any nilperfectoid $\C$-algebra $A$, a twistor projective line relative to $A$ can be defined. If $X$ is a complex analytic space, \cite{real_local_langlands} then defines a stack on the category of nilperfectoid $\C$-algebras, denoted $X^{\rm tw}$ and called the \textit{twistor stack of $X$}, which lives over the twistor projective line $\widetilde{\mathbb{P}}_{\mathbb{R}}^1$. The  relative analytic de Rham stack  of $X^{\rm tw}$ over $\widetilde{\mathbb{P}}_{\mathbb{R}}^1$ is simply the base change of the Betti stack $X_{\rm Betti}$ of $X$ to $\widetilde{\mathbb{P}}_{\mathbb{R}}^1$. The morphism from $X^{\rm tw}$ to its relative de Rham stack is an isomorphism away from the point $\infty \in \widetilde{\mathbb{P}}_{\mathbb{R}}^1$; the fiber at $\infty$ of $X^{\rm tw}$ is (a version of) the Hodge stack of $X$. 

This is entirely parallel to the $p$-adic story, $X^{\rm tw}$, resp., $X_{\rm Betti}\times_{\AnSpec(\Z)} \widetilde{\mathbb{P}}_{\mathbb{R}}^1$, playing the role of the relative $p$-adic twistor stack $X^{\Div^1/\FF_C}$ of a partially proper rigid space $X$ over $\GSpec(C)$, $C$ complete algebraically closed non-archimedean extension of $\Q_p$, resp.\ of the relative Hyodo--Kato stack $X^{\HK/\FF_C}$. One important difference is that the Hodge stack is a split gerbe for the overconvergent tangent bundle of $X$, while the relative Hodge--Tate stack is a gerbe for the twisted overconvergent tangent bundle of $X$, as follows from a relative variant of \Cref{sec:defin-first-prop-2-first-properties-of-analytic-prismatizationintro}, but in general non-split. 
\end{rem}

The next result is the technical heart of this paper.

\begin{thm}[\Cref{sec:comp-proet-cohom}] \label{thm:equivalence-perfect-complexesintro}
Let $X$ be a partially proper smooth rigid space over $\GSpec(\Q_p)$. The morphism
    \[
      \alpha_X \colon \FF_{X^\diamond}\to X^{\Div^1}
    \]
    is prim with  trivial codualizing complex, the pullback $\alpha_X^*$ on quasi-coherent sheaves is fully faithful, and $\alpha^\ast_X$  induces an equivalence on perfect complexes 
    $$\alpha_X^\ast\colon \Perf(X^{\Div^1})\isomarrow{}  \Perf(\FF_{X^\diamond}).
    $$ 
    A similar result holds for the preimages in $X^{\Div^1}$ and $\FF_{X^\diamond}$ of any open arc-substack $U$ of $\FF_{X^\diamond}^\diamond$, as well as for the Hodge--Tate stack $X^\HT$ and $\widehat{X^\diamond}$. 
    
    Moreover, a variant of this result holds for the relative $p$-adic twistor stack of a smooth partially proper rigid space over a perfectoid base.
    \end{thm}

\begin{rem} \label{rk:significance-equivalence-perfect-complexes}
  We want to think of this result, at least morally, as a non-abelian $p$-adic Hodge theory statement.
  Over the Hodge-Tate locus, the equivalence restricts (for vector bundles) to an equivalence between \textit{twisted} Higgs bundles (aka vector bundles on the analytic Hodge--Tate stack, which is a \emph{non-split} gerbe) and v-vector bundles, thus is a weak form of the $p$-adic Simpson correspondence\footnote{Only a weak form, since arguably an important aspect of the $p$-adic Simpson correspondence is about providing conditions under which, after additional choices, one can undo the twist.}.
  Over the whole $p$-adic twistor stack, but still for vector bundles, it describes relative $(\varphi,\Gamma)$-modules on $X$ (aka vector bundles on $\FF_{X^\diamond}$) as vector bundles on $X^{\Div^1}$.
  Away from the Hodge-Tate locus, the latter are analytic $D$-modules on $\FF_{X^\diamond} \backslash X^\diamond$ and these are glued to twisted Higgs bundles which one sees on the Hodge-Tate locus.
  So at least intuitively speaking, the theorem says that relative $(\varphi,\Gamma)$-modules (e.g., pro-\'etale $\Q_p$-local systems on $X$) and their cohomology can be described in terms of some version of \textit{$p$-adic twistors} on $X$, cf., also \Cref{rk-intro:analogy-real}. 
\end{rem} 

\begin{rem}
In restriction to the Hodge--Tate locus and in the geometric case, \Cref{thm:equivalence-perfect-complexesintro} explains why v-vector bundles come equipped with a canonical Higgs field (\cite{camargo2022geometricsen}): indeed, they come from the (relative) Hodge--Tate stack, which is a gerbe for the overconvergent Breuil-Kisin twisted tangent bundle.
\end{rem}

\begin{rem}
Already when $X=\GSpec(K)$, for $K$ a finite extension of $\Q_p$, \Cref{thm:equivalence-perfect-complexesintro} is an interesting statement, encoding a decompletion result for $(\varphi,\Gamma)$-modules over the Robba ring. We refer to \Cref{sec:arithm-case:-q_ppr} for a discussion, including a reinterpretation of familiar functors of Fontaine for $p$-adic Galois representations in terms of the geometry of $\GSpec(K)^{\Div^1}$.
\end{rem}

\Cref{thm:equivalence-perfect-complexesintro} is a local statement on $X$. \'Etale locally, the partially proper smooth rigid space $X$ looks like $\Gm^d$ ($d$ the dimension of $X$), the $d$-dimensional rigid analytic multiplicative group over $\Q_p$, with coordinate $\underline{T}=(T_1,\dots,T_d)$. On the one hand, the explicit descriptions using the cyclotomic prism and prismatic logarithm mentioned at the end of \Cref{intro-sec-1} show that for $X=\Gm^d$, if 
$\widetilde{X}=\lim_{\underline{T} \mapsto \underline{T}^p} X$ denotes the perfected $d$-dimensional multiplicative group (in Gelfand stacks), then the perfect prismatization  can be described as 
$$
X^{\widetilde{\Div^1}}\cong (\widetilde{X} \times_{\GSpec(\Q_p)} \mathcal{X}_{\rm cyc})/\Gamma^\la_d,
$$
with $\Gamma_d:= \Z_p^{\times,\la}\ltimes \Z_p^{d, \la}$.
On the other hand, a rather standard geometric Sen theory computation shows that, after passing to the quotient by Frobenius, this description can equivalently be rewritten as
\begin{equation} \label{locan-description}
X^{\Div^1} = \FF_{\widetilde{X}_{\Q_p^{\rm cyc}}^{d,\diamond}}^\la /\Gamma_d^\la,
\end{equation}
where $\FF_{\widetilde{X}_{\Q_p^{\rm cyc}}^{d,\diamond}}^\la$ denotes the Berkovich space obtained by gluing the spectra of rings of $\Gamma_d$-locally analytic vectors in the rings of functions on $\Gamma_d$-stable affinoid pieces of $\FF_{\widetilde{X}_{\Q_p^{\rm cyc}}^{d,\diamond}}$. This reveals an interesting and perhaps unexpected from first principles relation between the theory of analytic prismatization and the theory of locally analytic vectors, that we will exploit to prove \Cref{thm:equivalence-perfect-complexesintro}.

Indeed, the above description has to be compared with the formula
$$
\FF_{X^\diamond} = \FF_{\widetilde{X}_{\Q_p^{\rm cyc}}^{d,\diamond}}/\Gamma_d^{\rm cont}.
$$
Taking this into account, our approach to \Cref{thm:equivalence-perfect-complexesintro} is by studying the more general problem of comparing quasi-coherent sheaves and perfect complexes on $Y/G^{\rm cont}$ and $Y^\la/G^\la$ (to be precisely defined), when $G$ is a $p$-adic Lie group acting on a Berkovich space $Y$.  This will be split into two closely intertwined questions: 
\begin{enumerate}
\item decompletions by (derived) locally analytic vectors, and
\item vanishing of higher locally analytic vectors.
\end{enumerate}
The usual axioms of Sen theory enter the discussion of (2); as we just mentioned, this vanishing will hold in the examples of interest. Regarding (1), we prove the following result, of independent interest.

\begin{thm}[\Cref{TheoDecompletion}] \label{TheoDecompletion-intro} 
  Let $Y$ be a Berkovich space endowed with a continuous action of a $p$-adic Lie group $G$, that is, endowed with an action of the Gelfand group stack $G^{\rm cont}$.
  Let $h\colon Y/G^{\rm cont}\to |Y|_{\Betti}/G^{\la}$ be the natural induced map of Gelfand stacks, where the action of $G^{\la}$ on $|Y|_{\Betti}$ arises from the morphism of groups $G^{\la}\to G^{\sm}$ and the natural continuous action of $G^{\sm}$ on $|Y|_{\Betti}$ induced by the action on topological spaces.
  The following hold:
\begin{enumerate}

\item The map $h$ is prim with trivial codualizing complex, and hence we have a fully faithful  pullback functor 
\[
\ob{Mod}_{h_* \s{O}_{Y/G^{\rm cont}}} \ob{D}(|Y|_{\Betti}/G^{\la}) \hookrightarrow \ob{D}(Y/G^{\rm cont}).
\] We denote $\s{O}_Y^{\la}:= h_* \s{O}_{Y/G^{\rm cont}}$  and call it the sheaf of locally analytic functions on $Y$ with respect to $G$.

\item Suppose that $Y=\GSpec(R)$ is affinoid. Assume moreover that $R$ is \textit{classically decompleted}, in the sense of \Cref{Definition-Classically-decompleted}. 
  Let $R^{\la}\in \ob{CAlg}(\ob{D}(BG^{\la}))$ denote the Gelfand algebra of locally analytic vectors of $R$.
  Then the natural fully faithful functor
\[
\ob{Mod}_{R^{\la}}(\ob{D}(BG^{\la}))\xrightarrow{\sim} \ob{Mod}_{\s{O}_Y^{\la}} \ob{D}(|Y|_{\Betti}/G^{\la})
\] 
is an isomorphism.

\item Assume that $Y$ is locally strongly decompleted, in the sense of \Cref{AssumptionKeyTheoremDecompletionPerfectComplexes}. Then the pullback functor 
\[
h^*\colon \ob{Mod}_{\s{O}_Y^{\la}}(\ob{D}(|Y|_{\Betti}/G^{\la}))\hookrightarrow \ob{D}(Y/G^{\rm cont})
\]
induces an equivalence of dualizable objects.
Hence, we have an equivalence of stable categories
\[
\ob{Mod}_{\s{O}_Y^{\la}}(\ob{D}(|Y|_{\Betti}/G^{\la}))^{\ob{dual}} \xrightarrow{\sim} \Perf(Y/G^{\rm cont}). 
\]
\end{enumerate}
\end{thm}

The notion of a \textit{classical decompletion} that appears in (2) is defined for an action of a $p$-adic Lie group on a Berkovich space and essentially asks for vanishing of the higher locally analytic vectors. The technical condition appearing in (3) is even stronger (it implies classical decompletion) and will be satisfied in the examples of interest. We refer to the main body for a precise discussion of this notion.

If $Y$ is (locally) classically decompleted, $\s{O}^\la_Y$ is locally on $|Y|_\Betti$ concentrated in degree $0$.  In this case, one can use the relative spectrum (cf. \Cref{sec:deriv-berk-spac-relative-gelfand-spectrum}) of $\s{O}^\la_Y$ to define the Berkovich space
  \[
    Y^\la:=\underline{\GSpec}_{|Y|_\Betti}(\s{O}^\la_Y).
  \]
Moreover, under the strong decompletion assumption, part (3) says that the natural pullback map of quasi-coherent sheaves along $f\colon Y/G^{\rm cont}\to Y^{\la}/G^{\la}$   gives rise to an equivalence of perfect complexes
 \[
 \Perf(Y^{\la}/G^{\la}) \xrightarrow{\sim} \Perf(Y/G^{\rm cont}).
 \]

\begin{rem}\label{CoroEquivalenceTorsor-intro}
\Cref{TheoDecompletion-intro} has another interesting consequence, \Cref{CoroEquivalenceTorsor}: informally stated, \Cref{CoroEquivalenceTorsor} says that if  $X$ is a partially proper smooth rigid space over $\Q_p$ and $G$ a  $p$-adic Lie group, pro-\'etale $\underline{G}$-torsors on $X^\diamond$ decomplete uniquely into $G^{\la}$-torsors on $X^\HT$, and furthermore any such torsor produces a natural $G^{\la}$-torsor over the  $p$-adic twistor stack $X^{\Div^1}$. 
\end{rem}

\begin{rem} \label{rk-intro:sen-op}
If $X$ is a  partially proper smooth rigid space over $\Q_p$, the Faltings' extension for $X$ can be described in terms of the cotangent complex of $X^\HT$\footnote{This is in line with the heuristic philosophy that the Hodge--Tate stack is a twisted version of the Hodge stack of $X$ over ``$\mathbb{F}_1$'', cf. \cite[Remark 2.1.3]{bhatt2023algebraic}, and the fact that the Faltings' extension is something like the cotangent complex of $X$ over ``$\mathbb{F}_1$''.} and the arithmetic and geometric Sen operators attached to a $\underline{G}$-torsor over $X^\diamond$, with $G$ a $p$-adic Lie group, find a natural interpretation in terms of the morphism on cotangent complexes induced by the morphism $X^\HT \to BG^\la$ from \Cref{CoroEquivalenceTorsor-intro} (1). See \Cref{sec:decompl-tors}. 
\end{rem}

\begin{rem}\label{cohomo-twistor-stack-perfectoid-cover}
Let $X$ be a partially proper smooth rigid space over $\Q_p$, let $G$ be a compact $p$-adic Lie group and let $\widetilde{X}^\diamond \to X^\diamond$ be a  perfectoid pro-\'etale $\underline{G}$-torsor. In the case where $\FF_{\widetilde{X}}$ is locally classically decompleted, it is natural to wonder whether there is an isomorphism of Gelfand stacks
 $$
 X^{\Div^1}=\FF_{\widetilde{X}}^{\la}/G^{\la}
 $$
 over $BG^{\la}$, that is, whether the natural $G^{\la}$-torsor over $X^{\Div^1}$ mentioned in \Cref{CoroEquivalenceTorsor-intro} is $\FF_{\widetilde{X}}^{\la}$, generalizing what was proved and used above for $X=\GSpec(\Q_p)$ or $X=\Gm$ and the cyclotomic or toric chart. We do not know if this holds. However, note that it at least follows from the above (cf. \Cref{TheoLocallyAnalyticComparisonGeneralbis}), even without assuming that $\FF_{\widetilde{X}}$ is locally classically decompleted, that\footnote{We implicitly work in the derived world throughout, so the global sections functor is always understood as being derived.}
 \[
 \Gamma(X^{\Div^1},\s{O}) = \Gamma(G^\la, \Gamma(\FF_{\widetilde{X}},\s{O})^{G-\la}).
 \]
 Moreover, for the Hodge--Tate stack, we can prove such a statement under a suitable assumption on the (arithmetic or geometric) Sen operators of the torsor (mentioned in \Cref{rk-intro:sen-op}, see \Cref{sec:decompl-tors}.
\end{rem}

\begin{rem}
  Besides the work of Bhatt--Lurie and Drinfeld, another major motivation for this paper has been the work of Pan \cite{MR4390302} \cite{pan2022locally}, which started geometric Sen theory, following earlier work (for $p$-adic fields) of Berger--Colmez \cite{BergerColmezII}.
  The results listed above offer a stacky interpretation of parts of his theory.
  The work of Pan goes, however, much beyond and his results are best interpreted and generalized in the framework of the geometrization of the locally analytic $p$-adic Langlands program, which will be the subject of future work of the authors.
\end{rem}

An immediate consequence of \Cref{thm:equivalence-perfect-complexesintro} is the fact that for smooth(oid) partially proper rigid spaces, the $p$-adic twistor stack geometrizes pro-\'etale $\Q_p$-cohomology, in the sense that
\begin{equation}
\label{proetale}
\Gamma(X^{\Div^1}, \s{O}_{X^{\Div^1}}) = \Gamma_{\mathrm{pro\acute{e}t}}(X,\Q_p).
\end{equation}
This is of course to be compared to the statement that the cohomology of the structure sheaf on the \textit{perfect} prismatic site of a $p$-adic formal scheme is related to the $p$-adic \'etale cohomology of its generic fiber. Note however that the two sides of \eqref{proetale} differ for an arbitrary, e.g., non-reduced, Berkovich space. 

\begin{rem}[Archimedean/non-archimedean analogy, continued]
We have just seen that if $X$ is a smooth and partially proper rigid space over a complete algebraically closed extension $C$ of $\Q_p$, the $p$-adic twistor cohomology of $X$, i.e., the coherent cohomology of $X^{\Div^1/\FF_C}$, is simply the pro-\'etale $\Q_p$-cohomology of $X$. Continuing the analogy of \Cref{rk-intro:analogy-real}, this makes it clear that the Betti--de Rham comparison over $\C$ is \textit{not} an analogue of the \'etale-de Rham comparison over $C$, but rather of the Hyodo--Kato-de Rham comparison (both statements being a direct consequence of the description of the fiber over $\infty$ of the relevant stacks as the relative de Rham stack  of $X$); and similarly the twistor-Betti and \'etale-Hyodo--Kato comparisons are parallel to each other, coming from the generic isomorphism between the twistor stack and its de Rham stack.\footnote{Another good reason for which the correct analogue of Betti cohomology is Hyodo--Kato cohomology, not $p$-adic \'etale cohomology, is the fact that the latter is not motivic, contrary to the first two.} 

This perspective also corrects some apparent asymetry between complex and $p$-adic geometry: complex manifolds are locally very simple from the point of view of Betti cohomology -- locally they are just discs, hence contractible -- while small affinoid subspaces of a smooth rigid space still have a large \'etale fundamental group and hence a non-trivial $p$-adic \'etale cohomology. Trying to kill this cohomology leads to the theory of perfectoid spaces: to get rid of these cohomology classes in positive degrees, one extracts a lot of $p$-th roots. But this is only a difference between the two worlds if one takes for granted the analogy between Betti and $p$-adic \'etale cohomology; if instead, as we suggest here, one considers that the correct complex analogue of $p$-adic \'etale cohomology is twistor cohomology, then indeed the twistor cohomology of a disc is non-trivial and one needs to pass to a perfectoid cover to kill it, cf. \cite[\S VIII.4]{real_local_langlands}.  
\end{rem}

As usual, transmuting a cohomology theory has the benefit of providing for free a category of coefficients for it, with six operations. The $p$-adic twistor stack is no exception. Even though we do not provide a thorough discussion of the $6$-functor formalism for the $p$-adic twistor stack, we at least prove the following results.\footnote{We only state them for the $p$-adic twistor stack in this introduction, but they admit a version for the whole prismatization.}

\begin{thm}[\Cref{sec:div1-analyt-prism-properness-smoothness-for-div-1}, \Cref{CorollaryPoincar\'e duality}]
  \label{CorollaryPoincar\'e duality-intro} 
  \begin{enumerate}
  \item  The morphism $\Div^1\to \GSpec(\Q_p)$ is prim with trivial codualizing sheaf and cohomologically smooth with dualizing complex given by $\s{O}\{1\}[2]$, the shift of the first Breuil--Kisin twist of $\Div^1$ (cf.   \Cref{sec:breuil-kisin-twist}), i.e., the $(\varphi,\Gamma)$-module associated with the cyclotomic character. 
  \item   Let $f\colon Y\to X$ be a rigid smooth morphism of Berkovich spaces of relative dimension $d$.
    Let $f^{\Div^1}\colon Y^{\Div^1}\to X^{\Div^1}$ be the induced morphism.
  Then $f^{\Div^1}$ is cohomologically smooth and (if in addition $X$ is a partially proper smooth rigid space over $\Q_p$) there is a natural isomorphism $f^{\Div^1,!}\s{O}= \s{O}\{d\}[2d]$, where $\{d\}$ refers to the $d$-th Breuil--Kisin twist on $\Div^1$. 
  \end{enumerate}
  \end{thm}
  
  Part (2) of \Cref{CorollaryPoincar\'e duality-intro} gives Poincar\'e duality for geometric pro-\'etale $\Q_p$-cohomology, and combined with part (1) gives Poincar\'e duality for arithmetic pro-\'etale $\Q_p$-cohomology.
 
 \begin{remark}
Moreover, regarding coefficients, we relate the $p$-adic Riemann--Hilbert functor \cite{scholze2013p} to the category of vector bundles on the overconvergent resp.\ formal neighborhood of $X^\HT$ in $X^{\Div^1}$, see \Cref{sec:relation-to-rh}.
\end{remark}

\begin{rem}[Relation to other work]
The idea of enlarging the category of perfectoid rings (or of Fargues--Fontaine curves over such) to allow for some nilpotent thickenings has also been considered by Howe in his theory of inscribed v-sheaves \cite{howe2025inscriptiontwistorspadicperiods}. The main difference is that he restricts the class of allowed thickenings by only allowing finite flat augmented nilpotent thickenings of Fargues--Fontaine curves. We do not attempt a precise comparison of the two set ups here.

Independently of our definition of the analytic Hodge--Tate stack, Bhatt and Zhang defined the \textit{Simpson gerbe} of a smooth rigid space over a perfectoid extension of $\Q_p$ as a $\Gm$-gerbe over the Tate-twisted cotangent bundle of $X$, cf. \cite{bhatt-aspects}. They prove an equivalent version of \Cref{thm:equivalence-perfect-complexesintro} in restriction to the Hodge--Tate locus, in the geometric case (i.e., over a perfectoid base). Their method of proof is however quite different from ours. They also investigate the relation to the work of Faltings and Heuer on the $p$-adic Simpson correspondence in the proper case \cite{heuer2025p}.
The relation between the Simpson gerbe and the analytic Hodge--Tate stack can be made precise via Cartier duality, cf. \cite{RCCartierGerbe}. See \Cref{sec:div1-analyt-prism-simpson-gerbe-for-q-p} for more details.  

Work in progress of Bhatt--Heuer--Petrov--Vologodsky investigates, in a different language, the geometry of the overconvergent Hodge--Tate stack. These authors have independently obtained results analogous to the ones of \Cref{sec:faltings-extension}; see \cite{bhatt_takagi} for an overview. 

Combined with the fact that the cohomology of the $p$-adic twistor stack of a smooth partially proper rigid space computes its pro-\'etale $\Q_p$-cohomology, the cohomological smoothness of $\Div^1$ and the description of its dualizing complex, \Cref{CorollaryPoincar\'e duality-intro} gives Poincar\'e duality for pro-\'etale $\Q_p$-cohomology of smooth partially proper rigid spaces. Such results have been previously obtained in \cite{colmez2023arithmetic}, \cite{colmez2024duality}, \cite{zhenghui_li_2024}, \cite{ALBM-FF}.

The recent paper \cite{RHRobbaLiuZhu} establishes a Riemann--Hilbert correspondence over the Robba ring. We (briefly) discuss the relation with our work in \Cref{relationrobaarh}.
\end{rem}

\begin{rem}[What we do not do in this paper] \label{rk-intro:limitations}
As mentioned at the beginning of this introduction, we only define the prismatization with $\Q_p$-coefficients, both in the geometry and in the coefficients. An integral version of these constructions will be studied in future work of the authors and Maximilian Hauck. The case where the coefficients are still rational but the geometry is in characteristic $p$ is related to the theory of rigid cohomology of schemes in characteristic $p$ and the subject of the forthcoming PhD thesis of Junhui Qin, see \cite{qin2026arithmetic}. We also postpone a more detailed investigation of the precise relation to the theory of (sheared) prismatization to the point where the integral theory will be available.

Finally, we do not consider refined variants of the prismatization here, such as the filtered prismatization or the syntomification. They have been investigated by Hauck in \cite{hauck}, at least again with $\Q_p$-coefficients. 
\end{rem}

\textbf{Plan of the paper.} \Cref{Section:Preliminaries} contains various preliminary results on the geometry of Gelfand stacks, as initiated in \cite{drff_stacks}. Worth pointing out here are a discussion of the analytification functor from algebraic stacks to Gelfand stacks (\Cref{sec:analyt-tann-stacks}) as well as a variant in the Gelfand context of the deformation to the normal bundle that will turn to be handy for analyzing the geometry of the prismatization ( \Cref{sec:cone-construction,sss:ConeFlatConnection}). \Cref{Section:Generalized-FF-curves,Section:rational_analytic_prismatization} contain our main definitions: (rational) $\delta$-spaces and analytic $\delta$-stacks, the (rational) analytic prismatization and its variants, such as the $p$-adic twistor stack. \Cref{sec:sen-theory} establishes our main locally analytic decompletion results for modules on quotients of Berkovich spaces by the action of a $p$-adic Lie group, \Cref{TheoDecompletion-intro}; the results of this chapter are independent of the theory of analytic prismatization and of independent interest. Combining the main results of \Cref{sec:sen-theory} with an explicit description of the $p$-adic twistor stack for $\GSpec(\Q_p)$ and the rigid multiplicative group, we establish \Cref{thm:equivalence-perfect-complexesintro} in \Cref{sec:comp-proet-cohom}. \Cref{sec:comp-proet-cohom} continues and ends with various applications of these results in (abelian and non-abelian) $p$-adic Hodge theory. Finally, \Cref{appendix-def-theory} discusses some basics on the cotangent complex and deformation theory for Gelfand stacks.
\\

\textbf{Notations and conventions.}
We denote by $\mathbb{A}^{1}$ the  Gelfand stack $A\in \cat{GelfRing} \mapsto A(\ast)$, and by $\mathbb{G}_m$ the  Gelfand stack $A\mapsto A(\ast)^\times$. We let
$$
\widehat{\mathbb{A}}^1 =\underset{n} \colim \GSpec(\Q_p[T]/T^n), ~~ \bb{A}^{1,\dagger}= \GSpec(\underset{n}\colim \Q_p\langle T/p^n \rangle)
$$
denote the formal and overconvergent completions of the affine line at the origin (in Gelfand stacks). We write $\Ga$, $\widehat{\mathbb{G}}_a$, $\Ga^\dagger$ when we see $\A^1, \widehat{\mathbb{A}}^1, \A^{1,\dagger}$ as group stacks. 

In this paper we will work only with induced solid ring structures over $\bb{Q}_p$.  Given a Gelfand ring $A$ and elements $f_1,\ldots, f_n,g\in A$ generating the unit ideal, the \textit{$\dagger$-rational localization of $A$} given by the locus $\{|f_i|\leq |g|\neq 0\}$ is the one of $\dagger$-spaces, namely,  
\[
A\langle \frac{f_1,\ldots, f_n}{g} \rangle_{\leq 1}:= A\langle T_1,\ldots, T_n \rangle_{\leq 1} /^L (gT_1-f_1,\ldots, gT_n-f_n),
\] 
where  $A\langle T_1,\ldots, T_n \rangle_{\leq 1}=A\otimes_{\Q_p} \Q_p\langle T_1,\ldots, T_n \rangle_{\leq 1}$. On the other hand, the   \textit{Tate-rational localization}  associated to the locus $\{|f_i|\leq |g|\neq 1\}$ is the one given by  
\begin{equation}\label{eqTateRatLocConventions}
A\langle \frac{f_1,\ldots, f_n}{g} \rangle:= A\langle T_1,\ldots, T_n \rangle /^L (gT_1-f_1,\ldots, gT_n-f_n),
\end{equation}
where $A\langle T_1,\ldots, T_n \rangle=A\otimes_{\Q_p} \Q_p\langle T_1,\ldots, T_n \rangle$.  

 Given a topos $\n{X}$ and a group $G$ on $\n{X}$ we let $BG$ denote the  classifying stack of $G$ in $\n{X}$. We think of $BG$ as a pointed object in $\n{X}$ via the presentation $BG=*/G$ as a right quotient of the point.
 
 We work derived throughout the paper: in particular, notations like $f_\ast, \Gamma(X,-), ...$ should implicitly be understood as derived. 
 
 We sometimes use (in proofs) the language of Gestalten (\cite{GestaltenScholze}). This is only for convenience, since this formalism allows to conveniently think geometrically about various categorical statements; all the arguments could be written  without appealing to the theory of Gestalten.
\\

\textbf{Acknowledgements.} The influence of the ideas of Bhatt--Lurie, Drinfeld and Pan on this work should be obvious to the reader. For useful discussions related to the content of this paper, we are grateful to Bhargav Bhatt, Guido Bosco, Maximilian Hauck, Ben Heuer, Sean Howe, Akhil Mathew and Lue Pan. The second author would like to thank Akhil Mathew and Matthew Emerton for giving him the opportunity to present the results of this paper at the University of Chicago, and all the participants there for their questions.

Some of the ideas in this paper were conceived during the Hausdorff trimester program ``The arithmetic of the Langlands program'' in 2023, and the authors would like to thank the HCM for making such an event possible. Work on this project was carried out at several institutions, including the Max Planck Institute for Mathematics in Bonn, the university of Bonn, the Institute for Advanced Study, Princeton University, l'IRMA (Universit\'e de Strasbourg), l'IH\'ES and LMO (Universit\'e d'Orsay); we thank these institutions for their hospitality and support. The third author would also like to thank Columbia University and the Simons Foundation for supporting him during the academic year 2023--2024 as a Junior Simons Fellow.
\\

\textbf{AI disclosure.} During the writing of this paper, we had several conversations with Claude, 55, owner of the neighborhood bar in Strasbourg. Unfortunately, none of them was of relevance to the topics of the paper.

\newpage
\section{Supplements on Analytic Geometry}\label{Section:Preliminaries}

The language used to develop the theory of analytic prismatization is based on the theory of analytic geometry of Clausen and Scholze \cite{clausenscholze_analytic_stacks}, more precisely, on the theory of Gelfand stacks introduced in \cite{drff_stacks}.
In this first section of the document we give a brief summary and introduction to the theory for the convenience of the reader, and add supplements to the theory that are necessary for later applications.

\emph{We note that starting in \Cref{sec:analyt-tann-stacks} we will always assume that Gelfand stacks/Gelfand rings/arc-stacks are quasi-finite dimensional (see \Cref{Subsection:Gelfand-Stacks} for the definition).}

\subsection{Gelfand stacks}\label{Subsection:Gelfand-Stacks}

Let $p$ be a prime number, and let $\ob{D}(\Q_p):=\ob{D}(\Q_{p,\sol})$ be the category of derived light solid $\Q_p$-vector spaces. 
 We fix the $p$-adic norm of $\Q_p$ such that $|p|=p^{-1}$.
 Let $\cat{GelfRing}$ be the category of animated separable Gelfand solid $\Q_{p}$-algebras in the sense of  \cite[Definition 3.3.1]{drff_stacks}.
 Instead of recalling the exact definition here, we list the main properties and structures of these.
   Given $r\in \mathbb{R}_{\geq 0}$ we let $A^{\leq r}\subset A$ be the solid $\Z_p$-submodule of elements of norm $\leq r$ as in \cite[Definition 2.2.8]{drff_stacks}.
    For $r=1$, $A^{\leq 1} $ has the structure of a solid animated $\Z_p$-algebra, and all the $A^{\leq r}$ are solid $A^{\leq 1}$-modules for $r\in \R_{\geq 0}$.
     As a consequence of $A$ being bounded (\cite[Definition 2.1.1]{drff_stacks}), we have that $A=\bigcup_{r\to \infty} A^{\leq r}$, and $A^{\leq 0}=\Nil^{\dagger}(A)$ is the \textit{spectral} or \textit{$\dagger$-nilradical} of $A$.
     Given $A$ a bounded $\Q_p$-algebra, we let
     $$
     A^{\dagger-\ob{red}}=A/A^{\leq 0}
     $$
     denote the \textit{$\dagger$-reduction} of $A$, and let
     $$
     A^{u}= \big(A^{\leq 1} \big)^{\wedge_p} [\frac{1}{p}]
     $$
      be its uniform completion; we recall that $A$  being a  separable Gelfand ring is equivalent to $A$ being a bounded $\Q_p$-algebra (with induced analytic ring structure from $\Q_p$) such that $A^u$ is a separable Banach algebra over $\Q_p$.
     The natural morphism of rings $A^{\dagger-\red}\to A^u$ is always an injection of static solid $\Q_p$-algebras (\cite[Proposition 2.3.6]{drff_stacks}).
     Each Gelfand ring $A$ has its \emph{Berkovich spectrum} $\mathcal{M}(A)$ (\cite[Definition 2.3.1]{drff_stacks}), which agrees with the classical Berkovich spectrum of the Banach algebra $A^u$.

 We let $\cat{GelfStk}^{\ob{aff}}=\cat{GelfRing}^{\op}$ be the category of affinoid Gelfand stacks endowed with the $!$-topology. 
 For $A$  a separable Gelfand ring we let $\GSpec(A)\in\cat{GelfStk}^{\ob{aff}} $ denote the affinoid Gelfand stack represented by $A$. 
 We let $\cat{GelfStk}$ be the category of Gelfand stacks, that is, the full  subcategory of $\cat{Ani}$-valued presheaves on $\cat{GelfStk}^{\ob{aff}}$ consisting of those presheaves that invert $!$-equivalences as in \cite[Definition 4.2.1]{drff_stacks}.   
 Let $\cat{GelfRing}^{\fd}\subset \cat{GelfRing}$ be the full subcategory of \textit{quasi-finite dimensional (qfd)} separable Gelfand rings $A$ as in \cite[Definition 4.2.14]{drff_stacks}.
 We let $\cat{GelfStk}^{\ob{aff},\fd}$ be the category of qfd affinoid Gelfand stacks and let $\cat{GelfStk}^{\fd}$ be the category of qfd Gelfand stacks (as in \cite[Definition 4.2.14]{drff_stacks}).
 We have geometric morphisms of $\infty$-topoi
\[
 \cat{AnStk} \xrightarrow{\beta} \cat{GelfStk}\xrightarrow{\alpha}\cat{GelfStk}^{\fd}, 
\]
where $\cat{AnStk}$ is the category of analytic stacks as in \cite{clausenscholze_analytic_stacks} (modified as in \cite{drff_stacks} using $!$-equivalences). 
These morphisms of topoi induce  left exact and colimit preserving functors  
\[
\cat{GelfStk}^{\fd}  \xrightarrow{\alpha^*}  \cat{GelfStk} \xrightarrow{\beta^*} \cat{AnStk}.
\]
We note though that $\alpha^\ast, \beta^\ast$ are probably not fully faithful.

Let $\cat{CondAni}$  be the category of condensed anima, i.e., the category of hypersheaves on the category of light profinite sets. 
We let 
\[
(-)_{\Betti}\colon \cat{CondAni}\to  \cat{GelfStk}^{\fd}
\]
 be the Betti realization of condensed anima  in $\fd$-Gelfand stacks \cite[Lemma 4.2.9, Remark 4.2.10]{drff_stacks}.

A major tool used in \cite{drff_stacks} and in the forthcoming sections is the theory of six-functor formalisms for Gelfand stacks. The six-functor formalism of  \textit{loc.\ cit.} was constructed from the one on analytic stacks as in \cite[Remark 4.2.8]{drff_stacks}, and is based on the general theory of abstract six-functor formalisms of  \cite{heyer20246functorformalismssmoothrepresentations} and \cite{scholze2022lectures}. However, it will be convenient for some of the arguments, in particular those appearing in \Cref{sec:p-adic-twistor}, to have an extended class of $!$-able maps. To construct this larger class of $!$-able maps,  we will take the (perhaps overkill) approach of \textit{Gestalten} of \cite{GestaltenScholze}.

\begin{remark}\label{Remark-Gestalten-Not-Necessary}
We must clarify to the reader that being familiar with the theory of Gestalten is not  necessary to follow any of the main constructions and arguments of the paper, and that working with the six-functor formalism from analytic stacks would be good enough for  all the important  results of this paper. We only use the extended six-functor formalism  on Gelfand stacks constructed  from Gestalten  in two points: 
\begin{enumerate}

\item  To trivialize the codualizing sheaf of $0$-proper maps (to be recalled in \Cref{Construction-Six-Functors-From-Gestalten}), this is then used in the identification of the dualizing sheaf of the prismatization in \Cref{sec:div1-analyt-prism-properness-smoothness-for-div-1,CorollaryPoincar\'e duality}; and 

\item To prove that enough morphisms of Berkovich spaces of finite presentation are $!$-able in \Cref{berkovich-smooth-qc-weakly-cohomologically-proper}. 

\end{enumerate}

The triviality of codualizing sheaves in (1) is used for convenience and to use the computations of \cite{ALBM-FF}. Nevertheless, there are more conceptual and direct arguments that are discussed and or mentioned in \Cref{alternative-computation-dualizing,sec:p-adic-twistor-1-remark-general-case}.  Concerning point (2), one could avoid using the theory of Gestalten and use the six-functor formalism from analytic stacks if one restricts to maps  of Berkovich spaces which are, locally in the $\dagger$-analytic topology, l.c.i. maps, that is, morphisms that factor as a regular Zariski closed immersion inside an overconvergent polydisc. In particular, this would suffice for smooth morphisms of Berkovich spaces and the discussion of Poincar\'e duality. 
\end{remark}

In the following, we recall some basic terminology of the theory of Gestalten, and construct the desired six-functor formalism on Gelfand stacks. 
  
 \begin{construction}\label{Construction-Six-Functors-From-Gestalten}
 
Let $\kappa$ be an uncountable regular cardinal that we fix once and for all and let $\cat{Gest}= (\varinjlim_{n} \ob{CAlg}(n\cat{Pr}_{\kappa}))^{\op}$  be the category of ($\kappa$-small) Gestalten. For example, the convention of \cite{GestaltenScholze}  is to take $\kappa=\aleph_1$, though all the theory developed in \textit{loc.\ cit.} can be adapted in an obvious way to any uncountable  cardinal $\kappa$.  We will need the following properties of morphisms of Gestalten in this paper:
\begin{enumerate}
\item $n$-affine morphisms as in \cite[Definition 6.5]{GestaltenScholze} for $n\geq 0$; and
\item $n$-proper morphisms as in \cite[Definition 6.12]{GestaltenScholze} for $n\geq 0$. 
\end{enumerate}
Other important properties of morphisms of Gestalten  introduced in \textit{loc.\ cit.} are \textit{$n$-\'etale, $n$-prim and $n$-suave} morphisms, but they will not be needed in this work. 

Let $\n{C}$ be a small category with finite  limits, for instance, the opposite category of $\kappa$-compact separable Gelfand rings.  Let $\ob{D}\colon \ob{Corr}(\n{C})\to \cat{Pr}_{\kappa}$ be a six-functor formalism on $\n{C}$ with all maps being $!$-able.  Let $*$ denote the final object of $\n{C}$. In \cite[Proposition 9.5 and Corollary 9.8]{GestaltenScholze} it is explained how to construct a ($\kappa$-small) colimit preserving functor commuting with pullbacks 
 \[
\ob{Gest}( (-)_{\ob{D}})\colon  \ob{Shv}_{\ob{D}}(\n{C})^{\kappa}\to \cat{Gest},
 \]
 where $\ob{Shv}_{\ob{D}}(\n{C})$ is the $\infty$-topos of presheaves obtained by inverting $!$-equivalences as in \cite[Definition 4.2.1]{drff_stacks}.  Set $S:=\ob{Gest}( (*)_{\ob{D}})$ be the Gestalt associated to the point, then the  natural restriction 
\[
\ob{Gest}( (-)_{\ob{D}})\colon  \ob{Shv}_{\ob{D}}(\n{C})^{\kappa}\to \cat{Gest}_{/S}
\]
is a morphism of ($\kappa$-small) $\infty$-topoi.  An important property of the functor  $\ob{Gest}( (-)_{\ob{D}})$ is that it is conservative by design, that is, the topos of presheaves on objects in $\n{C}$ obtained by inverting $!$-equivalences is precisely the one defined by the topology of  $\cat{Gest}$.  Given a morphism $f\colon Y\to X$ in $ \ob{Shv}_{\ob{D}}(\n{C})^{\kappa}$, we say that $f$ is $n$-affine, $n$-proper, $n$-\'etale, etc. if the morphism of Gestalten satisfies the same property.

In the rest of the discussion, we will omit the cardinal $\kappa$. Applying the previous construction to the category $\n{C}=\cat{GelfStk}^{\ob{aff}}$ of affinoid Gelfand stacks (where all arrows are $!$-able with $f_!=f_*$ by design), we get a conservative functor  of $\infty$-topoi
\[
\ob{Gest}( (-)_{\ob{D}})\colon \cat{GelfStk}\to \cat{Gest}_{/\ob{D}(\bb{Q}_{p})},
\]
where the target are Gestalten over $\ob{Gest}(\ob{D}(\bb{Q}_{p}))$.  Now, let $P=E$ denote the class of $0$-proper maps in $\cat{Gest}_{/\ob{D}(\bb{Q}_{p})}$ (for example, any map between affinoid Gelfand stacks is in $P$) and $I=\ob{iso}$ the class of isomorphisms. Thanks to \cite[Proposition 6.13]{GestaltenScholze} the class $E$ is stable under composition, base change, and has the left cancellation property, in particular, the pair $( \cat{Gest}_{/\ob{D}(\bb{Q}_{p})},0\textrm{-proper})$ is a geometric set up as in \cite[Definition 2.1.1 and Lemma 2.1.5]{heyer20246functorformalismssmoothrepresentations}.  Moreover, it is obvious that $(P,I)$ is a suitable decomposition of $E=P$.  Consider now the functor
\[
\s{O}(-)_1\colon \cat{Gest}_{/\ob{D}(\bb{Q}_{p})}\to \cat{Pr}
\]
 sending a Gestalt $X$ over $\ob{D}(\bb{Q}_{p})$ to its category of $1$-sheaves $\s{O}(X)_1$. By \cite[Definition 6.9]{GestaltenScholze} given $f\colon Y\to X$ a $0$-proper morphism of Gestalten, the functor $f^*\colon \s{O}(X)_1\to \s{O}(Y)_{1}$ admits a linear right adjoint which is compatible with  base change thanks to  \cite[Proposition 6.10 (iii)]{GestaltenScholze}. Then, thanks to \cite[Proposition 3.3.3]{heyer20246functorformalismssmoothrepresentations} the functor $X\mapsto \s{O}(X)_1$ extends to a  $3$-functor  formalism 
 \[
 \s{O}(-)_1\colon \ob{Corr}(\cat{Gest}_{/\ob{D}(\bb{Q}_{p})}, 0\textrm{-proper})\to \cat{Pr}
 \]
 which restricts to a $3$-functor formalism
\[
\ob{D}\colon \ob{Corr}(\cat{GelfStk},0\textrm{-proper})\to  \cat{Pr}. 
\]
Applying \cite{heyer20246functorformalismssmoothrepresentations}, there is a minimal class of arrows $E_{0\textrm{-proper}}$ in  $\cat{GelfStk}$ for which $\ob{D}$ extends uniquely to a six functor  formalism on $\ob{Corr}(\cat{GelfStk},E_{0\textrm{-proper}})$ and satisfying the properties of \cite[Theorem 3.4.11 and Lemma 3.4.13]{heyer20246functorformalismssmoothrepresentations}.  We will use this six-functor formalism on Gelfand stacks throughout the paper. 
 \end{construction}

\begin{remark}
  \label{sec:gelfand-stacks-cohomologically-smooth-proper-maps}
In order to apply the results obtained in \cite{drff_stacks} we must compare the six-functor formalism on Gelfand stacks of \Cref{Construction-Six-Functors-From-Gestalten} with the one used in \cite{drff_stacks}. Let $E_1$ denote the class of $!$-able maps on Gelfand stacks constructed in   \Cref{Construction-Six-Functors-From-Gestalten} and let $E_2$ denote the class of $!$-able maps for the six-functor formalism in \cite{drff_stacks}.

By construction, the class of arrows $E_2$ arises as pullback of the class of $!$-able arrows on analytic stacks via the map $\beta^*\colon \cat{GelfStk}\to \cat{AnStk}$. In particular, $E_2$ contains the class $E_0'$ of morphisms between affine Gelfand stacks. Let $E_0'\subset E_0$  denote the class of arrows in Gelfand stacks constructed via \cite[Theorem 3.4.11]{heyer20246functorformalismssmoothrepresentations}
 from the arrows $E_0'$. Then, by the stability properties established in the extension theorem, we have that  $E_0\subset E_1\cap E_2$. The previous implies the following:
 
 \begin{enumerate}
 
\item  The restrictions of $\ob{D}$ from $\ob{Corr}(\cat{GelfStk},E_1)$ or  $\ob{Corr}(\cat{GelfStk},E_2)$ to a six-functor formalism on  $\ob{Corr}(\cat{GelfStk},E_0)$ are naturally equivalent. 

\item  Any map that exists in $E_0$ and is suave, prim, cohomologically proper or cohomologically \'etale in $E_2$, has the same property when restricted to a six-functor formalism with $!$-able arrows in $E_0$, and also holds when considered as an arrow in the six-functor formalism with $!$-able arrows $E_1$. Indeed, by \cite[Lemmas 4.4.5 and 4.4.6]{heyer20246functorformalismssmoothrepresentations} for a map $f\colon Y\to X$ to satisfy these properties, it suffices to consider morphisms and functors constructed from the natural projections and (iterated) diagonal maps of finite products of $f$, whose $6$-functors are independent of $E_0,E_1$ or $E_2$. 
 \end{enumerate}

We finish this remark by pointing out some key examples of arrows of Gelfand stacks that are in $E_0$, and hence in the class $E_1$, that we will keep using from now on. 

\begin{itemize}

\item[(a)] Consider the pullback functor of topoi of condensed anima $(-)_{\Betti}\colon \cat{CondAni}\to \cat{GelfStk}$ from \cite[Lemma 4.2.9]{drff_stacks}, then all maps  of condensed anima which are $\bb{Q}$-fine in the sense of  \cite{heyer20246functorformalismssmoothrepresentations} are in $E_0$. Indeed, following \cite[Construction 3.5.16]{heyer20246functorformalismssmoothrepresentations} this class of arrows is constructed from the class $P$ of all maps between light profinite sets, which are mapped to morphisms between affinoid Gelfand stacks, and hence in $E_0$.

\item[(b)] Morphisms between derived Berkovich spaces $f\colon Y\to X$  (see \Cref{Subsection:Derived-Berkovich-Spaces}) with finite dimensional underlying topological space.

\item[(c)] De Rham stacks of morphisms of qfd arc-stacks over $\bb{Q}_p$ which are arc-locally on the target lqfd, see \cite[Proposition 5.1.4]{drff_stacks}. 

\end{itemize}
  \end{remark}
  
   \begin{remark}
   \label{sec:deriv-berk-spac-f-lower-shriek-via-0-affine-stuff} 
   By construction, all the $0$-proper maps $f\colon Y\to X$ in the six-functor formalism of Gelfand stacks with $!$-able arrows as in \Cref{Construction-Six-Functors-From-Gestalten} satisfy that $f_!=f_*$. More precisely, all these maps are prim by construction, and there is a natural trivialization of the codualizing sheaf $\delta_{f}\cong 1_{Y}$.
 \end{remark}

We introduce the following notions for Gelfand rings (cf., \cite[Definition 3.6.1]{drff_stacks}):

\begin{definition}
  \label{sec:gelfand-stacks-definition-totally-disc-...-gelfand-rings}
Let $A$ be a Gelfand ring. We say that $A$ is
\begin{enumerate}

\item  \textit{totally disconnected} if its Berkovich spectrum $\mathcal{M}(A)$ is profinite,

\item \textit{strictly totally disconnected} if $\mathcal{M}(A)$ is profinite and for each $x\in \mathcal{M}(A)$ its Berkovich residue field $k(x)_B$ (\cite[Definition 2.4.13]{drff_stacks}) is algebraically closed, and

\item  \textit{nilperfectoid} if $A^{\dagger-\red}\to A^u$ is an isomorphism and $A^u$ is a perfectoid ring.

\end{enumerate}
We let $\cat{GelfRing}^{\ob{td}}$, $\cat{GelfRing}^{\ob{std}}$ and $\cat{GelfRing}^{\ob{nilperfd}}$ be the full subcategories of $\cat{GelfRing}$ of totally disconnected, strictly totally disconnected and nilperfectoid rings, and $\cat{GelfStk}^{\ob{td}}$, $\cat{GelfStk}^{\ob{std}}$ and $\cat{GelfStk}^{\ob{nilperfd}}$ their opposite categories.
\end{definition}

\begin{remark}
Any strictly totally disconnected Gelfand ring $A$ with $A^{\dagger-\red}=A^u$ is nilperfectoid, cf.\ \cite[Lemma 3.6.2]{drff_stacks}.
\end{remark}

In \cite{drff_stacks} we proved the following statements.
The last two justify the restriction to qfd Gelfand rings.

\begin{theorem}
  \label{sec:gelfand-stacks-3-theorem-on-basis-for-gelfand-rings}
  \begin{enumerate}
  \item The category of separable Gelfand rings admits a basis for the $!$-topology given by nilperfectoid rings.
    More precisely, each separable Gelfand ring $A$ admits a descendable map $A\to B$ with $B$ nilperfectoid.
    
  \item The category of $\fd$ separable Gelfand rings admits a basis consisting of strictly totally disconnected nilperfectoid rings.
  \item Given an arc-hypercover $A\to A_\bullet$ of perfectoid rings with $A, A_\bullet$ qfd, then 
  \[
  \varinjlim_{\Delta}\GSpec(A_\bullet)\to \GSpec(A)
  \] is a $!$-equivalence.
  \end{enumerate}
\end{theorem}
\begin{proof}
Assertions (1) and (2) are proven in \cite[Proposition 4.6.4]{drff_stacks}. Point (3) is a consequence of \cite[Corollary 4.4.10]{drff_stacks}. 
\end{proof}

As recalled above, to any Gelfand stack $X$ is attached its category $\ob{D}(X)$ of quasi-coherent sheaves. If $X=\GSpec(A)$, $\ob{D}(X)$ contains as a full subcategory the category $\Perf(A)$ of perfect complexes of $A(\ast)$. It is shown in \cite[Lemma 4.2.11, Definition 4.2.12]{drff_stacks} that this association satisfies $!$-descent, defining for any Gelfand stack $X$  a full subcategory $\Perf(X)$ of $\ob{D}(X)$ which in fact is given exactly by dualizable objects in $\ob{D}(X)$.

Another subcategory of $\ob{D}(X)$ which will turn out to be useful when discussing cotangent complexes is the subcategory of \textit{locally eventually connective modules}. For a Gelfand stack $X$, we say that $M\in \ob{D}(X)$ is locally eventually connective if for each Gelfand ring $B$ with a morphism $\GSpec(B)\to X$ the pullback of $M$ is eventually connective, i.e., if it lies in $\ob{D}_{\geq -n}(B)$ for some $n\geq 0$ (depending on $B$). The following lemma shows that locally eventually connective modules satisfy descent for the $!$-topology.

\begin{lemma}\label{LemmaDescentlocallyAlmostConnective}
Let $\ob{D}_+(-)$ denote the right Kan extension from  $\cat{GelfStk}^{\ob{aff}}$  to presheaves on $\cat{GelfStk}^{\ob{aff}}$ of the functor 
\[
\ob{D}_{+}(-)\colon\cat{GelfStk}^{\ob{aff},\op}\to \cat{Cat}_{\infty}
\]
of eventually connective modules on Gelfand rings. Then $\ob{D}_+(-)$ inverts $!$-equivalences, and induces a sheaf of categories
\[
\ob{D}_+(-)\colon \cat{GelfStk}^{\op}\to \cat{Cat}_{\infty}. 
\]
\end{lemma}
\begin{proof}
We first show that $\ob{D}_+(-)$ is a sheaf for the $!$-topology. Let $A\to B$ be a descendable map of Gelfand rings and let $K=\ob{fib}(A\to B)$. Then by \cite[Proposition 2.6.5]{mann2022p} there is some $n\in \bb{N}$ such that the composite map 
\[
K^{\otimes_A n}\to K^{\otimes_A n-1}\to \cdots \to A
\]
vanishes.  This already holds in the stable category of eventually connective modules $\ob{D}_+(A)$. In particular $B$ is a descendable algebra in $\ob{D}_+(A)$, and by \cite[Proposition 3.22]{mathew2016galois} the natural map 
\[
\ob{D}_+(A)\xrightarrow{\sim} \ob{Tot}( \ob{D}_+(B^{\bullet}))
\]
is an equivalence, where $B^{\bullet}$ is the \v{C}ech co-nerve of $A\to B$, proving what we wanted.  It is left to see that $\ob{D}_+(-)$ inverts $!$-equivalences. By \Cref{LemmaCotopologicalLocalizations} down below, it suffices to show that if $Y\to X$ is a $!$-equivalence of $!$-sheaves, then the pullback functor $\ob{D}_+(X)\to \ob{D}_+(Y)$ is fully faithful, but this is clear since $\ob{D}_+(X)\subset \ob{D}(X)$ is a full subcategory and  $\ob{D}$ satisfies descent along $!$-equivalences by \cite[Remark 4.2.3]{drff_stacks}.
\end{proof}

We have used the following abstract lemma:

\begin{lemma}\label{LemmaCotopologicalLocalizations}
  Let $\n{C}$ be an $\infty$-topos.
  Let $\n{C}\to \n{C}^{\prime}$ be a cotopological localization and $F\colon \n{C}^{\op}\to \cat{Cat}$ a sheaf of categories. Then $F$ localizes through a sheaf on $\n{C}^{\prime}$ if and only if for each map $Y\to X$ in $\n{C}$ that is a $\n{C}^{\prime}$-isomorphism, the pullback functor $F(X)\to F(Y)$ is fully faithful.
\end{lemma}
\begin{proof}
  The condition is clearly necessary, suppose now that for a $\n{C}'$-isomorphism $Y\to X$ the map $F(X)\to F(Y)$ is fully faithful.
  By the definition of a cotopological localization, any $\n{C}'$-isomorphism is  $\infty$-connective which implies that all diagonals of $Y\to X$ are epimorphisms.
  Let $Y^{\bullet}\to X$ denote the \v{C}ech nerve.
  Since $F$ is a sheaf, the natural map $F(X)\xrightarrow{\sim} \ob{Tot}( F(Y^{\bullet}))$ is an equivalence.
  We claim that for all $[n]\to [m]$ in $\Delta$ the pullback map $F(Y^{n})\to F(Y^{m})$ is an isomorphism.
  If this holds we will have that $F(X)\to F(Y)$ is an isomorphism as wanted.

 To show the claim, by writing  $[n]\to [m]$  as a composite of boundary and degeneracy maps, one is reduced to show that the projection map $\pi_{n}\colon Y^{n+1}\to Y^{n}$ sending $(y_0,\ldots, y_{n+1})\mapsto (y_0,\ldots, y_n)$, and the degeneracy map $\sigma_{n}\colon Y^{n}\to Y^{n+1}$ given by $(y_0,\ldots, y_{n})\mapsto (y_0,\ldots, y_{n}, y_{n})$ induce isomorphisms when applying $F$. But both maps are $\n{C}'$-isomorphisms and the composite $\pi_n\sigma_{n}$ is the identity of $Y^n$, hence the composition $F(Y^{n})\xrightarrow{F(\pi_n)} F(Y^{n+1}) \xrightarrow{F(\sigma_n)} F(Y^n)$ consists of fully faithful maps whose composite is the identity. This implies that both $F(\sigma_n)$ and $F(\pi_n)$ must be isomorphisms, proving what we wanted. 
\end{proof}

As another application of \Cref{LemmaCotopologicalLocalizations} we can verify that nuclear modules form a sheaf of categories on Gelfand stacks.

\begin{lemma}
  \label{sec:gelfand-stacks-nuclear-modules-on-gelfand-stacks}
  The functor $R\mapsto \mathrm{Nuc}(\ob{D}(R))$ extends uniquely to a sheaf of categories on Gelfand stacks.
\end{lemma}
\begin{proof}
  By \Cref{LemmaCotopologicalLocalizations} and the fact that $R\mapsto \ob{D}(R)$ is a sheaf of categories on Gelfand stacks, it suffices to show descent of nuclear modules along a descendable map $R\to S$ of Gelfand rings.
  If $S_\bullet$ is the \v{C}ech conerve of $R\to S$ and $R_n:=\mathrm{Tot}_{\leq n}(S_\bullet)$, then the pro-object $\{R_n\}_{n}$ is pro-constant of value $R$ by descendability.  
  Let $N\in \ob{D}(R)$ and assume that $N\otimes_{R}S$ is nuclear.
  Let $P\in \ob{D}(R)$ be compact.
  We can deduce the following isomorphisms of pro-objects:  
\[
\begin{aligned}
 \{\Hom_{R}(P,R)\otimes_{R} N\}_n  \cong & \{  \Hom_{R}(P, R_n ) \otimes_R N\}_n  \\ 
 \cong  & \{ \ob{Tot}_{\leq n}( \Hom_{R}(P, S_{\bullet} )\otimes_{R} N )\}_n \\
 \cong & \{\ob{Tot}_{\leq n}( \Hom_{S_{\bullet}}(S_{\bullet}\otimes_{R} P , S_{\bullet}) \otimes_{S_{\bullet}} (S_{\bullet}\otimes_R N) ) \}_n \\ 
\cong & \{\ob{Tot}_{\leq n} \Hom_{S_{\bullet}}( S_{\bullet}\otimes_R P, S_{\bullet}\otimes_R N) \}_n \\ 
\cong & \{\Hom_R(P, \ob{Tot}_{\leq n} (S_{\bullet}\otimes_R N) )\}_n \\
\cong  &   \{\Hom_R(P, R_n\otimes_R N )\}_n \\ 
\cong & \{\Hom_R(P,N)\}_n.
\end{aligned}
\]
  Here, we used the pro-constancy of $\{R_n\}_n$ in the first and last equivalences. Nuclearity of $S_{\bullet}\otimes_R N$ was used in the fourth equivalence. All other equivalences are formal by $\otimes$-$\Hom$ adjunctions, the fact that all functors involved are exact functors of stable categories, and that $\ob{Tot}_{\leq n}$ is a finite limit. Taking limits of  pro-objects we get the desired nuclearity of $N$.  
\end{proof}

\subsection{Derived Berkovich spaces}\label{Subsection:Derived-Berkovich-Spaces}

Within the category of Gelfand stacks we find a generalization of Berkovich spaces that we call \textit{derived Berkovich spaces}, see \cite[Section 4.3]{drff_stacks}.
We recall their definition in this section.
Let $A$ be a Gelfand ring and let $X=\mathcal{M}(A)$ be its Berkovich spectrum \cite[Section 3.1]{drff_stacks}.
Given $f_1,\ldots, f_n,g$ elements in $A(\ast)$ generating the unit ideal, we let 
\[
A\langle \frac{f_1,\ldots, f_n}{g} \rangle_{\leq 1}:= A\langle T_1,\ldots, T_n\rangle_{\leq 1}/^{L}(gT_1-f_1,\ldots, gT_n-f_n)
\]
be the $\dagger$-rational localization of $A$ corresponding to the inequalities $\{|f_i|\leq |g|\neq 0\colon i=1,\ldots,n \}\subseteq \mathcal{M}(A)$ (here $A\langle T_1,\ldots, T_n\rangle_{\leq 1}$ is given by overconvergent functions on the $n$-dimensional unit disc over $A$).
Given a $\dagger$-rational localization $A\to  B=A\langle \frac{f_1,\ldots, f_n}{g} \rangle_{\leq 1}$ of Gelfand rings, the morphism of Berkovich spectra $\mathcal{M}(B)\to \mathcal{M}(A)$ is a closed immersion.
Its image is the classical rational subset $\{|f_i|\leq |g|\neq 0\colon i=1,\ldots,n \}$.
More precisely, the passage to the uniform completion induces an homeomorphism identifying $\dagger$-rational localizations $\mathcal{M}(A)=\mathcal{M}(A^u)$ (\cite[Proposition 2.4.9]{drff_stacks}).
Under this isomorphism $\mathcal{M}(B)$ is identified with the space of multiplicative seminorms $|-|_x\colon A^u\to \mathbb{R}_{\geq 0}$ such that $|f_i|_x\leq |g|_x\neq 0 $ for $ i=1,\ldots, n$ and $|p|_x=p^{-1}$.
 
The \textit{analytic topology} on Gelfand rings, equivalenly on affinoid Gelfand stacks, is the Grothendieck topology for which the covering sieves of $\GSpec(A)$ are those that contain a finite set of $\dagger$-rational localizations $\{\GSpec(A_i)\to \GSpec(A)\}_{i=1}^d$ such that $\bigsqcup_{i} \mathcal{M}(A_i)\to \mathcal{M}(A)$ is a \textit{strict} cover of topological spaces, that is, a closed cover that is refined by an open cover (\cite[Definition 4.3.1.(1)]{drff_stacks}). 

A \textit{derived Berkovich space} is a sheaf $X$ on $\cat{GelfStk}^{\ob{aff}}$ for the analytic topology such that there is an open cover $\bigsqcup_i U_i\to X$  for the analytic topology, with $U_i\subset \GSpec(A_i)$ an open subspace of an affinoid Gelfand stack (\cite[Definition 4.3.1]{drff_stacks}).
 We let $\cat{BerkSp}$ denote the category of derived Berkovich spaces.  
 By \cite[Proposition 4.3.3]{drff_stacks} the natural functor $\cat{BerkSp}\to \cat{GelfStk}$ obtained by left Kan extension from $\cat{GelfStk}^{\ob{aff}}\to \cat{GelfStk}$ is fully faithful.
 We let $\cat{BerkSp}^{\lfd}\subset \cat{BerkSp}$ be the full subcategory of \textit{locally quasi-finite dimensional} Berkovich spaces, that is, the full subcategory of derived Berkovich spaces admitting a $\dagger$-rational cover corepresented by $\fd$ Gelfand rings.
 We have a fully faithful embedding $\cat{BerkSp}^{\lfd}\subset \cat{GelfStk}^{\fd}$.
 From now on we will drop the ``derived'' adjective and simply call \textit{Berkovich spaces} the objects of $\cat{BerkSp}$.

 \begin{example}
   \label{sec:deriv-berk-spac-examples}
   $\dagger$-rigid spaces in the sense of \cite[Definition 4.3.6]{drff_stacks} are examples of Berkovich spaces.
   In particular, the structure sheaf of a Berkovich space can be ``overconvergent''.
 \end{example}

 \begin{example}
   \label{sec:deriv-berk-spac-rigid-spaces-as-berkovich-spaces}
   If $X$ is any rigid-analytic variety over a (qfd) Banach field over $\Q_p$, we can associate to it the Gelfand stack $X':=\varinjlim_{\Spa(A,A^+)\subseteq X \textrm{ open affinoid}} \GSpec(A)$, and the resulting functor $X\mapsto X'$ sends rational open covers (in the sense of adic spaces, also called Tate-rational localizations in this paper, cf. \eqref{eqTateRatLocConventions}) to epimorphisms of Gelfand stacks (using that Tate-rational localizations are descendable). In particular, the Gelfand stack $X'$ only depends on the compactification of $X$ as adic space over $\bb{Q}_p$.
   In this generality, it is not clear whether  $X'$ is a derived Berkovich space, as rigid-analytic varieties are glued along Tate-rational localizations while derived Berkovich are glued along $\dagger$-rational localizations.
   If $X$ is partially proper, this problem disappears and one obtains a fully faithful functor from partially proper rigid spaces to derived Berkovich spaces.
   If $X=\Spa(A,A^+)$ is affinoid, then its associated Gelfand stack is $\GSpec(A)$, which is a derived Berkovich space as well.
   
\emph{From now, we treat rigid spaces implicitly via their associated Gelfand stacks and add, whenever necessary, the assumption that they are derived Berkovich spaces, e.g., partially proper or affinoid}.
\end{example}

By a variant of \cite[Section 3.2]{drff_stacks}, given $A$ a $\fd$ Gelfand ring\footnote{The qfd assumption is used here to ensure the existence of a descendable cover $S_\Betti\to \mathcal{M}(A)_\Betti$ with $S$ light profinite.}  there is a natural morphism of Gelfand stacks (and hence also of associated analytic stacks) 
 \[
 \GSpec(A)\to \mathcal{M}(A)_{\Betti}.
 \]
 This functor extends to a proper morphism of Gelfand stacks
 \[
 X\to |X|_{\Betti}
 \]
 for a $\lfd$ Berkovich space $X$, where $|X|=\bigcup_{Z\subset X} \mathcal{M}(A)$ is the underlying Berkovich space of $X$.
 Here, $Z$ runs over affinoid closed subspaces of $X$. 
 (In general, given a Berkovich space $X$ we have a morphism of locales $\rho_X\colon \ob{Sm}(\ob{D}(X))\to |X|$ where $\ob{Sm}(\ob{D}(X))$ is the smashing spectrum of $\ob{D}(X)$, and $\rho_X$ is the map of locales whose pre-image of $|U|\subset |X|$ open is the open localization of symmetric monoidal stable categories $\ob{D}(X)\to \ob{D}(U)$.)

 We can use the Betti stacks to realize Berkovich spaces as ``topological (more precisely,  condensed) ringed  spaces''.
 This uses the relative Gelfand spectrum of a quasi-coherent algebra.
 
 \begin{definition}
   \label{sec:deriv-berk-spac-relative-gelfand-spectrum}
   Let $X$ be a qfd Gelfand stack and let $A\in \ob{CAlg}(\ob{D}(X))$ be a commutative algebra.
   We define $\underline{\GSpec}_{X}(A)$ as the functor on Gelfand rings $f\colon\GSpec(R)\to X$ over $X$ with value the anima $\mathrm{Hom}_{\ob{CAlg}(\ob{D}(R))}(f^\ast A,R)$.
 \end{definition}

 Here are some remarks about this construction.

 \begin{remark}
   \label{sec:deriv-berk-spac-remarks-relative-gelfand-spectrum}
   \begin{enumerate}
   \item Using the descent of $\ob{D}$, one can verify that $\underline{\GSpec}_{X}(A)$ is indeed a qfd Gelfand stack.
   \item If $g\colon Y\to X$ is a morphism of qfd Gelfand stacks, then $Y\times_X \underline{\GSpec}_{X}(A)\cong \underline{\GSpec}_{Y}(g^\ast A)$.
        \item If $f\colon Y=\GSpec(B)\to X=\GSpec(A)$ is a morphism of affinoid Gelfand stacks, then one verifies directly that $\underline{\GSpec}_X(f_\ast(\s{O}_Y))\cong Y$ as Gelfand stacks over $Y$.
   Indeed, given $h\colon \GSpec(R)\to \GSpec(A)$, we have\footnote{For the last equivalence to hold we recall that Gelfand rings have characteristic zero in our framework.}
   \[
     \underline{\GSpec}_X(f_\ast(\s{O}_Y))(R)=\Hom_{\ob{CAlg}(\ob{D}(R))}(R\otimes_AB,R)\cong \Hom_{\cat{GelfRings}}(B,R).
   \]
   \item Clearly, $A\mapsto \underline{\GSpec}_X(A)$ is a functor on $\ob{CAlg}(\ob{D}(X))$.
     In general, $\underline{\GSpec}_X(-)$ is not fully faithful, e.g., if $\pi_\ast A(\ast)$ contains an invertible homogeneous element of positive homological  degree, then $\underline{\GSpec}_X(A)=\emptyset$.
     Nonetheless, it is by descent fully faithful on those $A$, for which there exists a $!$-cover $\{g_i\colon \GSpec(R_i)\to X\}_{i\in I}$ such that $g_i^\ast A$ is a qfd Gelfand $R_i$-algebra.
   \item If $f\colon Y\to X$ is any morphism of Gelfand stacks, one can consider $A:=f_\ast(\s{O}_Y)\in \ob{D}(X)$, and set $Y^{\mathrm{aff}/X}:=\underline{\GSpec}_X(A)$ as the \emph{relative affinization} of $Y$ over $X$.
          If $f$ is prim, the relative affinization commutes with base change along morphisms $X'\to X$, i.e., $(Y\times_X X')^{\mathrm{aff}/{X'}}\cong X'\times_X(Y^{\mathrm{aff}/X})$.

          We note that there exists a canonical morphism $Y\to Y^{\mathrm{aff}/X}$ induced by the counit $f^\ast A=f^\ast f_\ast(\s{O}_Y)\to \s{O}_Y$.
          We say that $f$ is \textit{To\"en-affine} if $f$ is prim and the natural map $Y\to Y^{\ob{aff}/X}$ is an equivalence. 
   \item If $f\colon Y\to X$ is a morphism of Gelfand stacks and $X'\to X$ an epimorphism of Gelfand stacks with base change $f'\colon Y':=Y\times_X X'\to X'$, then $f$ is To\"en-affine if and only if $f'$ is To\"en-affine.
     Indeed, this follows from the fact that primness is local  on the target, that by primness of $f$ the formation of $h\colon Y\to Y^{\ob{aff}/X}$ commutes with base change, and that $h$ is an isomorphism if and only if $f$ is To\"en-affine.
   \end{enumerate}
 \end{remark}
 \begin{lemma}
   \label{sec:deriv-berk-spac-berkovich-spaces-as-relative-spectrum}
   Let $X$ be a qfd Berkovich space.
   Let $h\colon X\to |X|_\Betti$ be the natural morphism.
   Then the natural map $f\colon X\to Z:=X^{\mathrm{aff}/|X|_\Betti}=\underline{\GSpec}_{|X|_\Betti}(h_\ast(\s{O}_X))$ is an isomorphism, in fact, $h$ is a $0$-affine To\"en-affine morphism of Gelfand stacks (and in particular $0$-proper).
 \end{lemma}
 \begin{proof}
   We note that $h$ is $0$-affine, indeed, by \cite[Proposition 7.9]{GestaltenScholze} it suffices to prove this property locally on the target, and by covering $\bigsqcup_i S_{i,\Betti}\to |X|_{\Betti}$ with light profinite sets $S_i$, the base change $X\times_{|X|_{\Betti}} S_{i,\Betti}$ is represented by an affinoid Gelfand stack (thanks to the proof of \cite[Proposition 3.2.10]{drff_stacks}). Since $S_i=\GSpec(C^{\lc}(S_i, \bb{Q}_p))$ is also affinoid, the map $X\times_{|X|_{\Betti}} S_i\to S_i$ is a morphism of affinoid Gelfand stacks, and hence a $0$-affine map of Gestalten, this    yields the  $0$-afiness of $h$.   To\"en-affiness then follows from \Cref{sec:deriv-berk-spac-remarks-relative-gelfand-spectrum} (3) and (6).
   \end{proof}

 \begin{remark}
   \label{sec:deriv-berk-spac-different-notions-of-affineness}
   There are several options for the definition of an affine morphism $f\colon Y\to X$ of Gelfand stacks:
   \begin{enumerate}
   \item $g\colon Y\to Y^{\ob{aff}/X}$ is an isomorphism (this notion is not stable by base change, and so we avoid it),
   \item $f$ is To\"en-affine in the sense of \Cref{sec:deriv-berk-spac-remarks-relative-gelfand-spectrum} (this property is probably not stable under composition),
   \item $f$ is $0$-affine on the associated Gestalten, a notion which is stable under composition and base change, see \Cref{Construction-Six-Functors-From-Gestalten} and  \cite{GestaltenScholze},
     \item there exists an epimorphism $\bigsqcup_i \GSpec(A_i) \to X$ such that each base change $Y_{A_i}\to \GSpec(A_i)$ is represented by a Gelfand ring over $A_i$.
       If this condition is satisfied, $f$ is To\"en-affine (by \Cref{sec:deriv-berk-spac-remarks-relative-gelfand-spectrum} (6)),  $0$-affine by \cite[Proposition 7.9]{GestaltenScholze}.   We note that it is not clear whether this notion of ``affineness'' is  stable under composition.
           \end{enumerate}
 \end{remark}

\subsection{Analytic de Rham stacks}\label{Subsection:Analytic-de-Rham}

The main character of \cite{drff_stacks} is the analytic de Rham stack functor for qfd Gelfand stacks. 
Let $\cat{Perfd}^{\ob{aff}}$ be the category of separable affinoid perfectoid spaces, which are partially proper over $\Z_p$, seen as Berkovich spaces, and let $\cat{Perf}^{\ob{aff}}\subset \cat{Perfd}^{\ob{aff}}$ be the full subcategory of perfectoid spaces over $\F_p$. 
We let $\cat{Perfd}^{\ob{aff},\fd}$ and $\cat{Perf}^{\ob{aff},\fd}$ be the full subcategories of $\fd$ affinoid perfectoid spaces over $\Z_p$ and $\F_p$ respectively.
We endow $\cat{Perfd}^{\ob{aff}}$ with the arc-topology, that is the Grothendieck topology generated by surjections between affinoid perfectoid spaces.
The category of (light) arc-stacks is the category of arc-hypersheaves of anima 
\[
\cat{ArcStk}:=\cat{ArcStk}_{\F_p}:=\widehat{\ob{Shv}}(\cat{Perf}^{\ob{aff}}, \cat{Ani}),
\]
(\cite[Definition 4.1.5]{drff_stacks}).
Similarly, we let $\cat{ArcStk}^{\fd}=\widehat{\ob{Shv}}(\cat{Perf}^{\ob{aff},\ob{qfd}}, \cat{Ani})$ be the category of $\qfd$ arc-stacks.
Given a  separable Banach algebra $A$ we let $\Marc(A)$ denote the arc-sheaf associated to $A$, that is, the sheaf sending a separable perfectoid ring $R$ over $\F_p$ to the set of isomorphism classes of tuples $(R^{\sharp}, A\to R^{\sharp})$, where $R^{\sharp}$ is an untilt of $R$, and $A\to R^{\sharp}$ is a morphism of Banach algebras.
The tilting equivalence induces an equivalence of topoi
\[
\cat{ArcStk}_{\Q_p}:=\cat{ArcStk}_{/\Marc(\Q_p)}\cong \widehat{\ob{Shv}}(\cat{Perfd}^{\ob{aff}}_{\Q_p},\cat{Ani})
\]
and similarly for the $\qfd$ variant. 

Restriction along $\cat{Perfd}^{\ob{aff},\ob{td}}_{\bb{Q}_p}\subset \cat{GelfStk}^{\ob{aff}}$ to totally disconnected perfectoid rings induces the \textit{perfectoidization functor}
\[
(-)^{\diamond}\colon \cat{GelfStk}\to \cat{ArcStk}_{\Q_p},
\]
(\cite[Definition 4.5.1]{drff_stacks}, \cite[Proposition 4.5.4]{drff_stacks}).
Given $A$ a separable Gelfand ring, we let $\Marc(A)=\GSpec(A)^{\diamond}$ denote its perfectoidization.
Notice that the uniform completion gives rise to an equivalence of arc-stacks $\Marc(A^u)\xrightarrow{\sim} \Marc(A)$ (there is also no conflict of notation, besides that the perfectoidization is viewed as an arc-stack over $\Marc(\Q_p)$).

Since geometric realizations of arc-hypercovers of separable strictly totally disconnected perfectoid rings (or qfd separable perfectoid rings) are $!$-equivalences (cf.\ \cite[Lemma 4.4.3]{drff_stacks}), the functor  $(-)^{\diamond}$ is right adjoint to a functor $\widehat{(-)}\colon \cat{ArcStk}_{\Marc(\Q_p)} \to  \cat{GelfStk}$.
The functor $\widehat{(-)}$ preserves fiber products, but not the final object.
We have obvious variants for $\qfd$ arc and Gelfand stacks.

The functor $(-)^{\diamond}$ is itself a left adjoint in a morphism of topoi, cf.\ \cite[Proposition 3.3.3]{drff_stacks}.

\begin{definition}\label{def:de-Rham-stack}
The analytic de Rham stack is defined as the right adjoint
\[
(-)^{\dR}\colon \cat{ArcStk}^{\qfd}\to \cat{GelfStk}^{\qfd}
\]
of the perfectoidization functor.
More explicitly, given a qfd arc-stack $X$, its analytic de Rham stack $X^{\dR}$ is the functor whose value on a qfd nilperfectoid ring $A$ over $\Q_p$ is 
\[
X^{\dR}(A)= X(\Marc(A^u)).
\]
If $X$ is a qfd Gelfand stack, we set $X^\dR:=(X^\diamond)^\dR$.
\end{definition}

Thanks to \cite[Theorem 5.6.1]{drff_stacks} the functor $(-)^{\dR}$ satisfies arc-hyperdescent on qfd arc-stacks.

\begin{remark}
  \label{sec:prop-morph-gelf-2-description-of-de-rham-stack-in-qfd-setting}
  Let $X$ be a qfd Gelfand stack.
  As a feature of working in the qfd setting we note that for any qfd nilperfectoid $\Q_p$-algebra $A$ we have $\widehat{\mathcal{M}_\arc(A^u)}\cong \GSpec(A^u)$ and hence equivalences
  \[
    \Hom_{\cat{ArcStk}_{\Q_p}}(\Marc(A^u), X^\diamond)\cong \Hom_{\GelfStk}(\widehat{\Marc(A^u)},X)=\Hom_{\GelfStk}(\GSpec(A^u),X),
  \]
  i.e., $X^\dR(A)=X(A^u)$.
\end{remark}

\begin{definition}
  \label{sec:analytic-de-rham-locally-analytic-version}
  Let $X$ be a rigid-analytic variety over $\Q_p$ (seen as a Gelfand stack, cf.  \Cref{sec:deriv-berk-spac-rigid-spaces-as-berkovich-spaces}).
  We define
  \[
    X(\Q_p)^\la:=X\times_{X^\dR} (\underline{X(\Q_p)})^\dR,
  \]
  where $\underline{X(\Q_p)}$ is the locally profinite set of $\Q_p$-rational points of $X$.
\end{definition}

Of course, if $X$ has a group structure, then $X(\Q_p)^\la$ is a group object in Gelfand stacks.
This applies notably to $X$ the analytification of a reductive group over $\Q_p$. 
\\

\emph{From now on, and until the end of this paper, we will always assume that Gelfand stacks/Gelfand rings/arc-stacks are quasi-finite dimensional, unless explicitly mentioned. Consequently, we drop superscripts ``qfd'' from now on to lighten notation.}

\subsection{Analytifications and Tannakian stacks}
\label{sec:analyt-tann-stacks}

In this subsection, we develop a bit of material on the analytification of schemes. A recent reference on this matter (and generalizations) is \cite{wang2026relativegagatheoremapplication}.
In order to keep the paper self-contained, we will prove now the key analytification results that will be relevant for this paper.
We note that our results will be slightly different from those in \cite{wang2026relativegagatheoremapplication}.

Let $S=\GSpec(R)$ be an affinoid Gelfand stack over $\bb{Q}_p$. Given an affine (derived) scheme $\mathcal{Y}=\Spec(R')\to \mathcal{S}:=\Spec(R(\ast))$ over the underlying (animated) ring $R(\ast)$, we let
\[
  \mathcal{Y}^\an\in \GelfStk_{/S}
\]
be the Gelfand stack sending $\GSpec(A)\to S$ to
\[
\n{Y}^{\an}(A)=  \mathrm{Maps}_{R(\ast)}(R',A(\ast)).
\]
Using $\ast$-descent of the  formation of the  underlying discrete ring $A(*)$ along $!$-equivalences, it is easy to see that $\mathcal{Y}^\an$ is indeed a Gelfand stack. 

\begin{lemma}
  \label{sec:analyt-tann-stacks-1-factorization-over-etale-stacks}
  The functor $(-)^\an$ on the category $\Aff_{/\mathcal{S}}$ of affine (derived) schemes over $\mathcal{S}$ extends uniquely to a colimit preserving functor
  \[
    (-)^\an\colon \ob{Shv}_{\et}(\Aff_{/\mathcal{S}})\to \GelfStk_{/S}
  \]
  on \'etale sheaves, that is called the ``analytification''.
  Furthermore, $(-)^\an$ commutes with finite limits and hence is the left adjoint in a morphism of topoi. Concretely, if $\n{X}\in\ob{Shv}_{\et}(\Aff_{/\mathcal{S}}) $, then $\n{X}^{\an}$ is the $!$-sheafification of the presheaf on Gelfand rings  sending $\GSpec(A)\mapsto \n{X}(A(*))$.   We call $\mathrm{Shv}_{\et}(\Aff_{/\mathcal{S}})$ the category of \'etale algebraic stacks over $\mathcal{S}$.\footnote{The name ``algebraic stack'' is in the literature also used to mean Artin stack, but we warn the reader that here we only use it  to refer  that these are stacks are of algebro-geometric nature.} 
\end{lemma}
\begin{proof}
  We can formally extend $(-)^\an$ to a colimit preserving functor on presheaves of anima on $\Aff_{\mathcal{S}}$. By construction, this functor sends a presheaf $\n{X}$ on $\Aff_{\mathcal{X}}$ to the $!$-sheafification of $\GSpec(A)\mapsto \n{X}(A(*))$.
  Indeed, this follows from the fact that the analytification functor is constructed as a partial right adjoint on presheaves of the   functor 
  \[
  (-)^{\ob{alg}}\colon \cat{GelfStk}_{/S}^{\ob{aff}}\to \cat{Aff}_{/\n{S}}
  \]
mapping $\GSpec(A)\mapsto \Spec A(*)$.  This adjunction promotes to an adjunction on presheaves
\[
(-)^{\ob{alg}}\colon \n{P}(\cat{GelfStk}_{/S}^{\ob{aff}}) \leftrightarrows \colon \n{P}(\cat{Aff}_{/\n{S}}) \colon (-)^{\ob{an}},
\]
where both $(-)^{\ob{alg}}$ and $(-)^{\an}$ are left Kan extensions of their restrictions to affine objects. This makes the  left Kan extension 
\[
(-)^{\an}\colon \n{P}(\cat{Aff}_{/\n{S}})\to \cat{GelfStk}_{/S}
\]
the $!$-sheafification  of the functor sending $\n{X}$ to $\GSpec(A)\mapsto \n{X}(A(*))$, proving the second claim of the lemma.

The factorization over \'etale sheaves is implied if $\mathcal{Y}\mapsto \mathcal{Y}^\an$ on affine derived schemes over $\mathcal{S}$ sends \'etale covers to $!$-covers.
  If $f\colon \mathcal{Y}'\to \mathcal{Y}$ is an \'etale cover of schemes over $\mathcal{S}$, then we first reduce to the case that $\mathcal{Y}, \mathcal{Y}'$ are of finite presentation over $\mathcal{S}$ (as $!$-covers are stable under base change).
  But any \'etale morphism of schemes of finite presentation over $\mathcal{S}$ induces a rigid \'etale morphism of derived Berkovich spaces, see \Cref{sec:prop-morph-gelf-3-definition-of-rigid smooth}.
  Moreover, coverings are sent to coverings (as one only needs to check surjectivity on underlying topological spaces, cf. \Cref{sec:prop-morph-gelf-1-weakly-formally-smooth-is-epimorphism} below).
  This finishes the argument.
  To see that $(-)^\an\colon \mathrm{Sh}_{\et}(\Aff_{/\mathcal{S}})\to \GelfStk_{/S}$ preserves finite limits it is sufficient to note that the analytification of affine derived schemes over $\mathcal{S}$ preserves fiber products.
\end{proof}

For an animated ring $R$ we let $D(R)$ be the derived category of (discrete) $R$-modules.
Given an \'etale algebraic stack $\mathcal{X}$ we denote by $\ob{QCoh}(\mathcal{X})$ its category of quasi-coherent sheaves, i.e., the value of the right Kan extension of the functor $R\mapsto D(R)$ on (animated) affine schemes.
By Kan extending from the affine case, we can define the connective part $\ob{QCoh}_{\geq 0}(\mathcal{X})\subseteq \ob{QCoh}(\mathcal{X})$ and the full subcategories $\Perf(\mathcal{X}), \ob{PCoh}(\mathcal{X})\subseteq \ob{QCoh}(\mathcal{X})$ of perfect and pseudo-coherent (sometimes also called almost perfect or almost compact in the literature) complexes.
By flat descent (\cite[Proposition 2.8.4.2]{Luriea}), the property of being connective/perfect/pseudo-coherent can be tested on a flat cover $g\colon X\to \mathcal{X}$, i.e., flat morphism which is a surjection for the flat topology.

\begin{definition}\label{DefinitionTannakian}
Let $\n{S}=\Spec (R)$ be an (animated) affine scheme over $\bb{Q}$.

\begin{enumerate}

\item  An \'etale algebraic stack $\n{X}$ over $\n{S}$ is called \textit{Tannakian} if for all $R$-algebras $R'$, the natural map
\[
\ob{Map}_{/\Spec (R)}(\Spec (R'), \n{X})\to \ob{Map}^{\otimes}_{D(R)}(\ob{QCoh}(\n{X}),D(R')),\ f\mapsto f^\ast
\]
towards $D(R)$-linear symmetric monoidal functors $\ob{QCoh}(\n{X})\to D(R')$ is fully faithful, and its essential image consists of those symmetric monoidal functors preserving connective objects. 
\item An \'etale algebraic stack $\n{X}$ is called \textit{strongly Tannakian} if it is Tannakian and if $\ob{QCoh}_{\geq 0}(\mathcal{X})$ is generated under filtered colimits by connective pseudo-coherent objects.
\item An \'etale algebraic stack $\n{X}$ over $\n{S}$ is called \textit{perfect} if the natural map 
\[
\ob{Ind}(\cat{Perf}(\n{X}))\to \ob{QCoh}(\n{X})
\]
is an equivalence.
\end{enumerate}
\end{definition}

To construct examples, we will use the following simple observation.

\begin{lemma}
  \label{sec:analyt-tann-stacks-2-when-pseudo-coherence-is-not-necessary}
  Let $\mathcal{X}$ be a perfect \'etale algebraic stack. 
    If there exists a locally noetherian derived scheme $X$\footnote{This means that locally $X$ is $\Spec(R')$ with $R'$ a noetherian animated ring, i.e., $\pi_0(R')$ is noetherian and each $\pi_i(R'), i\geq 0,$ is a finitely generated $\pi_0(R')$-module.} and a flat morphism $X\to \mathcal{X}$, which is a surjection for the flat topology, then $\ob{QCoh}_{\geq 0}(\mathcal{X})$ is generated under filtered colimits by connective pseudo-coherent objects.
\end{lemma}
\begin{proof}
Let $M\in \ob{QCoh}_{\geq 0}(\mathcal{X})$.
  As $\mathcal{X}$ is perfect, we can write $M=\varinjlim_{i}P_i$ as a filtered colimit of perfect complexes $P_i$ on $\mathcal{X}$.
  Let $g\colon X\to \mathcal{X}$ be a flat morphism from a locally noetherian derived scheme, which is a surjection for the flat topology.
  We note that over the locally noetherian derived scheme $X$ an object $\ob{QCoh}(X)$ is pseudo-coherent if and only if each cohomology module is finitely generated.
  This implies that each $\tau_{\geq 0}g^\ast P_i$ is pseudo-coherent.
  By flatness of $X\to \mathcal{X}$ we can conclude that each $\tau_{\geq 0}P_i$ is pseudo-coherent.
  This in turn implies that
  \[
    M=\tau_{\geq 0} M=\tau_{\geq 0}\varinjlim\limits_{i}P_i=\varinjlim\limits_{i}\tau_{\geq 0}P_i
  \]
  is a filtered colimit of connective pseudo-coherent objects as desired (we use that filtered colimits commute with the truncation functor $\tau_{\geq 0}$).
\end{proof}

\begin{example}\label{ExamplePerfectlyTannakian}
Let $R$ be a ring over $\bb{Q}$ and $\n{S}=\Spec (R)$.  

\begin{enumerate}

\item  Let $\n{X}$ be an \'etale algebraic stack over $\n{S}$ admitting an fpqc surjection $\Spec (R')\to \n{X}$ from an affine scheme over $R$ (also called \textit{geometric  stack} in \cite[Section 1.1]{stefanich2023tannakaduality1affineness}). By Theorem 1.0.3 (1) of  \textit{loc.\ cit.}, if  $\n{X}$ is qcqs and has quasi--affine diagonal, it is Tannakian.

\item Following \cite{BenZviFrancisNadlerIntegralTransforms}, we have the following examples of perfect stacks. 

\begin{enumerate}

\item Quasi-compact derived schemes with affine diagonal, cf.  \cite[Proposition 3.19]{BenZviFrancisNadlerIntegralTransforms}.  

\item In characteristic zero, the quotient $\n{X}/G$ where $\n{X}$ is a quasi-projective derived scheme and $G$ an affine algebraic group acting on $\n{X}$, cf \cite[Corollary 3.22]{BenZviFrancisNadlerIntegralTransforms}. See also Proposition 3.21 of \textit{loc.\ cit.} for a more general class of examples of this kind.

\item Finite  products of perfect stacks are perfect, cf.\ \cite[Proposition 3.24]{BenZviFrancisNadlerIntegralTransforms}.
  More generally, if $\n{X}_1,\n{X}_2$ are perfect stacks over a stack $\n{Y}$ with affine diagonal, then $\n{X}_1\times_\n{Y} \n{X}_2$ is perfect.

\item Let $\n{X}$ be a perfect stack and $\Sigma$ a compact anima, then $\n{X}^{\Sigma}=\underline{\ob{Map}}(\Sigma,\n{X})$ is a perfect stack, cf.\ \cite[Corollary 3.25]{BenZviFrancisNadlerIntegralTransforms}. For example, $\n{X}^{S^1}=\n{X}\times_{\n{X}\times_{\n{S}} \n{X}} \n{X}$ where the fiber product is taken with respect to diagonal maps is perfect.

\end{enumerate}

\item Let $K$ be a field of characteristic zero.  Combining the examples of (1), (2) and \Cref{sec:analyt-tann-stacks-2-when-pseudo-coherence-is-not-necessary}, we deduce that the following stacks are strongly Tannakian:

\begin{enumerate}

\item  The quotient stack of schemes $\bb{A}^{1,\mathrm{sch}}_{K}/\bb{G}_{m,K}^{\mathrm{sch}}$ corepresenting virtual Cartier divisors, where $\bb{A}^{1,\mathrm{sch}}_{K}=\Spec(K[T])$ and $\bb{G}_{m,K}^{\mathrm{sch}}=\Spec(K[T^{\pm 1}])$ denote the schematic affine line and multiplicative group.\footnote{Recall that we reserve the notations $\A^1, \Gm$ for the analytic affine line and multiplicative group over $\Q_p$, seen as Gelfand stacks.}

\item The classifying stacks $BG$ where $G/\Spec K$ is an affine algebraic group.

\end{enumerate}

\end{enumerate}

\end{example}

What makes perfect strongly Tannakian stacks relevant to us is the fact that their analytification admits an explicit functor of points description. 

\begin{proposition}\label{PropositionFunctorOfPointsAnalytification} 
  Let $S=\GSpec (R)$ be an affinoid Gelfand stack with underlying affine scheme $\mathcal{S}=\Spec (R(*))$.
  Let $\n{X}\to \n{S}$ be a strongly Tannakian stack, which is perfect.
  Then the functor 
\[
\n{X}^{\mathrm{pre}-\an}\colon \cat{GelfRing}_{R}\to \cat{Ani}
\]
sending $R'$ to $\n{X}(\Spec (R'(*)))$ is a Gelfand stack.
In particular, $\n{X}^{\mathrm{pre}-\an}=\n{X}^{\an}$ is the analytification of \Cref{sec:analyt-tann-stacks-1-factorization-over-etale-stacks} and if $Y$ is a Gelfand stack, morphisms $Y\to \n{X}^\an$ are the same as $\Perf(R)$-linear symmetric monoidal functors $F\colon \Perf(\mathcal{X})\to \Perf(Y)$ whose extension to $\ob{QCoh}(\mathcal{X})=\mathrm{Ind}(\Perf(\mathcal{X}))\to \ob{D}(Y)$ sends connective pseudo-coherent objects to  locally connective objects (i.e. objects which are connective after base change to a cover given by  affinoid Gelfand stacks). 
\end{proposition}
\begin{proof}
  By the assumption that $\n{X}$ is Tannakian and perfect, morphisms from $\Spec(R'(\ast))$ to $\n{X}$ can be described as $R$-linear symmetric monoidal functors from perfect complexes on $\n{X}$ to perfect complexes on $R'(\ast)$, whose extension to $\ob{QCoh}(\mathcal{X})\to D(R'(\ast))$ sends connective objects on $\mathcal{X}$ to connective objects.
  As $\mathcal{X}$ is strongly Tannakian, it suffices to ask that connective pseudo-coherent objects map to connective objects in $D(R'(*))$ (or equivalently to connective objects in $\ob{D}(R')$).
  We note that each such functor automatically maps pseudo-coherent objects to pseudo-coherent modules over $R'(\ast)$ (as it is the pullback of a morphism $\Spec(R'(\ast))\to \mathcal{X}$).
  Since perfect complexes satisfy $!$-descent on Gelfand stacks (they even form a sheaf of categories on $\ob{GelfStk}$), cf.\ \Cref{Subsection:Gelfand-Stacks}, \cite[Lemma 4.2.11]{drff_stacks}, and connectiveness of discrete pseudo-coherent complexes can be checked $!$-locally (by applying \cite[Proposition VI.1.2]{real_local_langlands} to the top cohomology), we can deduce that $\mathcal{X}^{\mathrm{pre}-\an}=\mathcal{X}^{\an}$ and that, as claimed, morphisms $Y\to \mathcal{X}^\an$ from a Gelfand stack $Y$ over $\GSpec(R)$ identify with $\Perf(R)$-linear symmetric monoidal functors $\Perf(\mathcal{X})\to \Perf(Y)$ whose extension to $\ob{QCoh}(\mathcal{X})\to \ob{D}(Y)$ sends connective pseudo-coherent objects to locally connective objects. 
\end{proof}

In the above proof, we use that connectivity of \emph{discrete} pseudo-coherent objects can be checked $!$-locally.
For this reason, we use strongly Tannakian stacks, and not the perfectly Tannakian stacks from \cite[Definition 5.2.1]{wang2026relativegagatheoremapplication}.

\subsection{Generalized Cartier divisors on Gelfand stacks}
\label{sec:gener-cart divis}

In this subsection we will discuss the notion of (generalized) Cartier divisors on Berkovich spaces, that will enter prominently in the construction of the analytic prismatization.

\begin{notation}\label{NotationGmAffineLine}
Recall that we denote by  $\bb{G}_m$ and $\bb{A}^1$  the analytification as Gelfand stacks of the multiplicative group $\Gm^{\mathrm{sch}}$ and the affine line $\mathbb{A}^{1,\mathrm{sch}}$ over $\bb{Q}_p$. 
\end{notation}

We note that $\mathbb{A}^1/\Gm, B\Gm$ are the analytifications of the strongly Tannakian stacks $\mathbb{A}^{1,\mathrm{sch}}/\Gm^{\mathrm{sch}}, B\Gm^{\mathrm{sch}}$ over $\Q_p$.
Recall from \Cref{Subsection:Gelfand-Stacks} the sheaf of categories on Gelfand stacks $\Perf\colon \GelfStk\to \cat{Cat}_\infty$, whose value on $\GSpec(A)$ is the category of perfect complexes of $A(\ast)$.

We clarify now the notion of a line bundle, i.e., of vector bundles of (constant) rank $1$ (\cite[Definition 4.2.9]{drff_stacks}).

\begin{lemma}
  \label{sec:gener-cart divis-1-equivalent-description-of-line-bundles}
  Let $X$ be a Gelfand stack.
  The following anima are naturally equivalent:
  \begin{enumerate}
  \item The anima of maps $\ob{Map}_{\cat{GelfStk}}(X,B\Gm)$,
  \item The anima of line bundles in $\ob{D}(X)$, i.e., objects which $*$-locally on $X$ are represented by perfect complexes of perfect amplitude $0$ and rank $1$. 
  \end{enumerate}
\end{lemma}
\begin{proof}
This follows formally from the analogous description in algebraic geometry, \Cref{PropositionFunctorOfPointsAnalytification} and the examples of perfect Tannakian stacks in \Cref{ExamplePerfectlyTannakian}. We give a self contained argument for the convenience of the reader.  Let $\mathcal{Y}$ be the stack parametrizing objects as in (2). The trivial line bundle produces a map $f\colon *\to \mathcal{Y}$ which is surjective, namely, if $X=\GSpec(A)$ is a separable Gelfand ring and $M\in \ob{D}(X)$ is an object that $*$-locally is a perfect complex of perfect amplitude $0$ and rank $1$, then by \cite[Proposition 3.3.5]{drff_stacks}  and \cite[Corollary VI.1.4]{real_local_langlands}, $M$ is a projective $A$-module of rank $1$, i.e., a line bundle placed in degree $0$. After a further $\dagger$-rational localization one can then assume that $M$ is trivial proving the desired surjectivity.  Now, for the loops of $f$ to be $\mathbb{G}_m$ it suffices to see that for a Gelfand ring $A$, the automorphisms of the trivial line bundle are $A(*)^{\times}$. 
\end{proof}

Next, we turn to Cartier divisors.

\begin{proposition}\label{PropDivisorsBerkovichSpaces}
  Let $X$ be a  Gelfand stack.
  The following anima are naturally equivalent:
\begin{enumerate}
\item  The mapping anima $\ob{Map}_{\cat{GelfStk}}(X, \A^1/\mathbb{G}_m)$ where $\mathbb{G}_m$ acts by multiplication on $\A^1$. 

\item The anima of (generalized) Cartier divisors on $X$, that is, the anima of pairs $(\mathcal{L},s)$ given by a line bundle $\mathcal{L}$ and a morphism $s\colon\mathcal{L}\to \s{O}_X$ in $\ob{D}(X)$. 

\item (If $X$ is furthermore assumed to be a  Berkovich space.) The anima of morphisms of Berkovich spaces $f\colon Z\to X$ such that, locally in the analytic topology of $X$, $f$ is given by the (derived) zero locus of a section in $\s{O}_X$.
\end{enumerate}
The functor from (1) to (3) sends a map $g\colon X\to \A^1/\mathbb{G}_m$ to the pullback $f\colon Z\to X$ of $g$ along the zero map $\{0\}/\Gm\to \A^1/\mathbb{G}_m$.
The functor from (3) to (2) sends $f\colon Z\to X$ to the map $\mathcal{L}\to \s{O}_X$ where  $\mathcal{L}=\ob{fib}(\s{O}_X\to \s{O}_Z)$. The functor from (2) to (1) sends a morphism $\mathcal{L}\to \s{O}_X$ to the map $X\to \A^1/\mathbb{G}_m$ obtained by descent from the $\mathbb{G}_m$-equivariant map 
\[
\widetilde{f}\colon \mathbb{V}(\mathcal{L})^{\times}\to \A^1
\]
where $\mathbb{V}(\mathcal{L})^{\times}$ is the $\mathbb{G}_m$-torsor of trivializations of $\mathcal{L}$, and $\widetilde{f}$ is obtained by the function defined by the  map $\mathcal{L}\to \s{O}_X$ and the universal trivialization of $\mathcal{L}$ over $ \mathbb{V}(\mathcal{L})^{\times}$. 
\end{proposition}
\begin{proof}
This follows formally from the analogous description in algebraic geometry, \Cref{PropositionFunctorOfPointsAnalytification} and the examples of perfect Tannakian stacks in \Cref{ExamplePerfectlyTannakian}. We give a self-contained  argument for the convenience of the reader.

  We first show the equivalence between (1) and (2). Let $\mathcal{Y}$ be the Gelfand stack representing the functor in (2).
  Over $\A^1$ we have the natural section in $\s{O}_{\A^1}$ induced by the identity map, and taking $\mathcal{L}$ to be the trivial bundle this produces a map $f\colon \A^1\to \mathcal{Y}$.
  We claim that this map is surjective, namely, by \Cref{sec:gener-cart divis-1-equivalent-description-of-line-bundles} we know that any  Gelfand stack $X$ endowed with a line bundle $\mathcal{L}$ admits a $!$-cover $X'$ where $\mathcal{L}$ is isomorphic to the trivial line bundle $\s{O}_{X'}$, thus, any section $s\in \mathcal{L}$ arises $!$-locally from a global section of $\s{O}_X$.
  Next, we prove that the automorphisms of the map $f$ are given by the action of $\mathbb{G}_m$. For that, let $X=\GSpec(A)$ be an affinoid Gelfand stack   and consider a section $t\colon \s{O}_X\to \s{O}_X$. We are interested in the anima of commutative triangles 
\[
\begin{tikzcd}
\s{O}_X \ar[dr,"{t}"] \ar[dd,"g"] & \\ 
& \s{O}_X \\
\s{O}_X \ar[ru,"{t'}"] & 
\end{tikzcd}
\]
where $g$ is an isomorphism.  Since $g$ is an isomorphism, it is given by a unit $g\in \s{O}_X(*)^{\times}$. After fixing $g$, there is at most one of such triangles determined by $t'=g^{-1}t$. This implies that the anima of such triangles is given by $\s{O}_X(*)^{\times}$, and therefore that the natural map $\overline{f}\colon \A^1/\mathbb{G}_m \to \mathcal{Y}$ is an equivalence of Gelfand stacks.   It is clear from the construction that the isomorphism $\overline{f}$ is the inverse of the one constructed in the statement of the proposition. 

For the equivalence between (1),(2) with (3), it suffices to show that for $X$ a  Berkovich space and a map $h\colon X\to \A^1/\mathbb{G}_m$, the pullback at zero gives rise to a Berkovich space $Z\to X$ as described in (3). Indeed, the map $h$ corresponds to a map of line bundles $\mathcal{L}\to \s{O}_X$.
By the Fredholm property of Gelfand rings (cf.\ \cite[Proposition 3.3.5]{drff_stacks}  and \cite[Corollary VI.1.4]{real_local_langlands}) locally in the analytic topology of $X$ the sheaf $\mathcal{L}$ is given by an invertible $\s{O}_X$ sitting in degree zero and by further localizing it is given by a section $f\in \s{O}_X$. In that case the pullback of the zero section of $h$ is the (derived) zero locus of $f$ in $X$.  
\end{proof}

In this text, the words ``Cartier divisor'' will always be understood in the above sense, i.e., what would often be rather called a generalized or virtual Cartier divisor.

\subsection{Properties for morphisms of Gelfand stacks}
\label{sec:prop-morph-gelf}

In this subsection, we discuss several properties for morphisms of Gelfand stacks.

We note that the properties that we introduce here are stable under base change (via obvious arguments), unless we mention the contrary.

\begin{lemma}
  \label{sec:prop-morph-gelf-1-equivalent-characterizations-dagger-formally-etale-morphisms}
  Let $f\colon Y\to X$ be a morphism of Gelfand stacks.
  The following properties are equivalent:
  \begin{enumerate}
  \item The morphism $Y\to Y^{\dR/X}=X\times_{X^\dR} Y^\dR$ is an isomorphism.
  \item For any Gelfand ring $A$ the natural map $Y(A)\to X(A)\times_{X^\diamond(\Marc(A^{\dagger-\red}))}Y^\diamond(\Marc(A^{\dagger-\red}))$ is an isomorphism of anima.
  \item $f$ is $\dagger$-formally \'etale in the sense of \cite[Definition 4.7.2]{drff_stacks}, or equivalently by \Cref{sec:prop-morph-gelf-2-description-of-de-rham-stack-in-qfd-setting}, for any nilperfectoid $\Q_p$-algebra $A$ the natural map $Y(A)\to X(A)\times_{X(A^u)} Y(A^u)$
    is an isomorphism of anima.
  \end{enumerate}
  If these conditions are satisfied, we call $f$ \emph{$\dagger$-formally \'etale}.
\end{lemma}
\begin{proof}
  The equivalence of (1) and (3) follows by evaluating on nilperfectoids and using \Cref{sec:prop-morph-gelf-2-description-of-de-rham-stack-in-qfd-setting}.
  
  Clearly (2) implies (1).
  Thus, assume that (1) holds.
  Let $A$ be a Gelfand ring, and let $A\to B$ be a $!$-cover with $B$ nilperfectoid.
  We note that all terms $B_i:=B\otimes_A\ldots \otimes_A B$, $[i]\in \Delta$, in the \v{C}ech nerve of $A\to B$ are nilperfectoid.
  The assumption implies $Y(B_i)\cong X(B_i)\times_{X(B_i^u)}Y(B_i^u)$ for any $i$.
  Taking the limit over $\Delta$ yields the equivalence
  \[
    Y(A)\cong X(A)\times_{X^\diamond(\Marc(A^{\dagger-\red}))} Y^\diamond(\Marc(A^{\dagger-\red}))
  \]
  because $\Marc(B^u)\to \Marc(A^{\dagger-\red})$ is an arc-cover with \v{C}ech nerve $\Marc(B_i^u)$.
\end{proof}

For example, $\dagger$-rational localizations of Berkovich spaces are $\dagger$-formally \'etale, e.g., using \cite[Proposition 3.1.9]{drff_stacks}.

\begin{remark}
  \label{sec:prop-morph-gelf-1-dagger-formally-unramified}
  Replacing ``isomorphism'' by ``monomorphism'' in \Cref{sec:prop-morph-gelf-1-equivalent-characterizations-dagger-formally-etale-morphisms} one could also define $\dagger$-formally unramified morphisms of Gelfand stacks, but this is not necessary for this paper.

  We note that we will not consider the following stronger version of $\dagger$-formally \'etaleness of a morphism $f\colon Y\to X$: for \emph{any} Gelfand ring $A$ the natural map $Y(A)\to X(A)\times_{X(A^{\dagger-\red})}Y(A^{\dagger-\red})$ is an isomorphism.
\end{remark}

\begin{example}\label{ExampleFormallyEtale}
  Let $f\colon Y\to X$ be any map of Gelfand stacks. Then $Y^{\dR/X}\to X$ is $\dagger$-formally \'etale.
 \end{example}

Regarding notions of $\dagger$-formal smoothness the situation is more ambiguous.
We start by giving several possible definitions of $\dagger$-formal smoothness, in particular recalling the notion introduced in \cite[Definition 4.7.2]{drff_stacks}.

\begin{definition}
  \label{sec:prop-morph-gelf-1-definition-dagger-formally-etale}
  Let $f\colon Y\to X$ be a morphism of Gelfand stacks.
  We say that $f$ is
  \begin{enumerate}
  \item \textit{$\dagger$-formally} smooth if for any nilperfectoid $\Q_p$-algebra $A$ the natural map
    \[
      Y(A)\to X(A)\times_{X(A^u)} Y(A^u)
    \]
    is an epimorphism,
  \item \textit{weakly $\dagger$-formally smooth} if $Y\to Y^{\dR/X}=X\times_{X^\dR}Y^\dR$ is an epimorphism,
  \item \textit{strongly $\dagger$-formally smooth} if for any Gelfand ring $A$ the natural map $Y(A)\to X(A)\times_{X(A^{\dagger-\red})}Y(A^{\dagger-\red})$ is an epimorphism.
  \item \textit{$\dagger$-formally \'etale locally $\dagger$-formally smooth} if there exist $\dagger$-formally \'etale epimorphisms $Y'\to Y$ and $X'\to X$ and a commutative diagram
\[\begin{tikzcd}
	{Y'} & {X'} \\
	Y & X
	\arrow["{f'}", from=1-1, to=1-2]
	\arrow[from=1-1, to=2-1]
	\arrow[from=1-2, to=2-2]
	\arrow["f", from=2-1, to=2-2]
      \end{tikzcd}\]
    with $f'$ $\dagger$-formally smooth.
  \end{enumerate}
\end{definition}

Clearly we have the implications ``strongly $\dagger$-formally smooth $\Rightarrow$ $\dagger$-formally smooth $\Rightarrow$ $\dagger$-formally \'etale locally $\dagger$-formally smooth $\Rightarrow$ weakly $\dagger$-formally smooth''.

For a weakly $\dagger$-formally smooth morphism there is a convenient criterion to check that it is a cover.

\begin{lemma}
  \label{sec:prop-morph-gelf-1-weakly-formally-smooth-is-epimorphism}
  Let $f\colon Y\to X$ be a weakly $\dagger$-formally smooth morphism.
  Then $f$ is an epimorphism of Gelfand stacks if and only if $f^\diamond\colon Y^\diamond\to X^\diamond$ is an epimorphism of arc-stacks.
\end{lemma}
\begin{proof}
  As $(-)^\diamond$ is a left adjoint, the ``only if'' statement is clear.
  Conversely, note that both morphisms in the composition $Y\to Y^{\dR/X}\to X$ are epimorphisms: the first because $f$ is weakly $\dagger$-formally smooth, and the second because it is a base change of the epimorphism $(f^\diamond)^{\dR}$.
\end{proof}

More in line with classical notions, we discuss deformation theoretic properties of Gelfand stacks in \Cref{app:DefTheoryGelfand}. 

\begin{remark}
  \label{sec:prop-morph-gelf-1-formally-smooth-not-dagger-formally-smooth}
  An interesting example of a formally smooth (\Cref{sec:cotang-compl-morph-definition-having-a-cotangent-complex}) morphism, which is not $\dagger$-formally smooth is given by $\GSpec(\Q_p\langle T\rangle)\to \GSpec(\Q_p)$ (\cite[Example 4.7.7.]{drff_stacks}).
  The formal smoothness is implied by the calculation of the cotangent complex $\mathbb{L}_{\Q_p\langle T\rangle/\Q_p}=\Q_p\langle T\rangle dT[0]$, cf.\ \Cref{sec:cotang-compl-morph-cotangent-complex-via-p-completion}.
  
\end{remark}

Next we want to discuss notions of ``\'etale/smooth'' morphism on Berkovich spaces.
Again there are several options, coming from the fact that $\dagger$-rational localizations of Berkovich spaces are not open immersions.

\begin{remark}
  \label{sec:prop-morph-gelf-4-finite-etale-maps}
  There is no problem of defining \emph{finite} \'etale maps $f\colon Y\to X$ of Berkovich spaces.
  Namely, the following properties are equivalent (\cite[Remark 4.1.12]{drff_stacks}, \cite[Proposition 3.2.11]{drff_stacks}):
  \begin{enumerate}
  \item $f$ is $\dagger$-formally \'etale and the morphism $Y^\diamond\to X^\diamond$ is finite \'etale (in the sense of arc-stacks),
  \item for each morphism $\GSpec(A)\to X$ for a Gelfand ring $A$, one has $Y\times_X \GSpec(A)=\GSpec(B)$ with $B\cong A\otimes_{A(\ast)}B(\ast)$ and $B(\ast)$ is a finite \'etale $A(\ast)$-algebra (i.e., $A\to B$ is a finite \'etale morphism of Gelfand rings), and
    \item there exists a $!$-cover of $X$ by $\GSpec(A_i), i\in I$, such that for each $i\in I$, we have $Y\times_{X}\GSpec(A_i)\cong \GSpec(B_i)$ with $B_i$ finite \'etale over $A_i$.
    \end{enumerate}
    By descent, one can then also define finite \'etale morphisms of Gelfand stacks.
\end{remark}

In the following we mostly restrict our attention to Berkovich spaces.
We clarify the meaning of an open immersion in this context.

\begin{definition}
  \label{sec:prop-morph-gelf-1-definition-open-immersion}
  Let $f\colon Y\to X$ be a morphism of Berkovich spaces.
  Then $f$ is an \textit{open immersion} if $|Y|\to |X|$ is an open immersion of topological spaces and the natural map $Y\to |Y|_{\Betti}\times_{|X|_\Betti}X$ is an isomorphism. We define \textit{closed and locally closed immersions} in a similar way. 
\end{definition}

Equivalently, $f$ is an open immersion if $f$ is $\dagger$-formally \'etale and $f^{\diamond}\colon Y^{\diamond}\to X^{\diamond}$ is an open immersion in the sense of arc-stacks (and same for closed  and locally closed immersions).

\begin{remark}
  \label{sec:prop-morph-gelf-2-open-immersions-for-gelfand-stacks}
  Using that $!$-covers of Berkovich spaces induce quotient maps on topological spaces (even on arc-stacks), one checks that the property of being an open immersion is local on the target.
  In particular, we can define an open immersion of Gelfand stacks to be a morphism whose pullback to any Berkovich space is an open immersion. The same works for closed immersions. 
\end{remark}

We now define the notions of rigid \'etale and rigid smooth morphisms, cf.\ \cite[Definition 4.3.4]{drff_stacks}.

\begin{definition}
  \label{sec:prop-morph-gelf-3-definition-of-rigid smooth}
  Let $f\colon Y\to X$ be a morphism of Berkovich spaces.
  \begin{enumerate}
  \item $f$ is \emph{rigid \'etale} if it is locally in the analytic topology on $Y$ and $X$ the composition of an open immersion and a finite \'etale map.
  \item $f$ is \emph{rigid smooth} if it is locally in the analytic topology on $Y$ and $X$ the composition of a rigid \'etale morphism $Y\to \A^{n}_{X}:=\A^n_{\Q_p}\times_{\GSpec(\Q_p)}X$ with the projection $\A^n_X\to X$.
  \end{enumerate}
\end{definition}

For example, smooth morphisms of partially proper rigid-analytic varieties over a non-archimedean field $K/\bb{Q}_p$ yield rigid smooth morphisms on associated Berkovich spaces.

\begin{lemma}
  \label{sec:prop-morph-gelf-6-rigid smooth-implies-strongly-dagger-formally-smooth}
  Let $f\colon Y\to X$ be a rigid smooth morphism of Berkovich spaces.
  Then locally in the analytic topologies of $Y$ and $X$ the morphism $f$ is strongly $\dagger$-formally smooth.
\end{lemma}
\begin{proof}
  This follows from \cite[Proposition 3.7.5]{RCAndeRham}.
\end{proof}

We note that a $\dagger$-rational localization $\GSpec(A')\to \GSpec(A)$ is in general \emph{not} rigid smooth.
This motivates the following weakening.

\begin{definition}
  \label{sec:prop-morph-gelf-5-berkovich-etale-smooth} Let $f\colon Y\to X$ be a morphism of Berkovich spaces.
    \begin{enumerate}
    \item We say that $f$ is \emph{Berkovich-\'etale} if it is locally on a strict closed cover on $Y$ and $X$ the composition of a $\dagger$-rational localization and a finite \'etale map.
    \item We say that $f$ is \emph{Berkovich-smooth} if it is locally on a strict closed cover on $Y$ and $X$ the composition of a Berkovich-\'etale morphism $Y\to \A^{n}_X$ and the projection $\A^{n}_X\to X$.
  \end{enumerate}  
\end{definition}

Clearly, each rigid smooth morphism is also Berkovich-smooth, and each rigid \'etale/Berkovich-\'etale morphism is $\dagger$-formally \'etale.

\begin{warning}
   Continuing \Cref{sec:deriv-berk-spac-rigid-spaces-as-berkovich-spaces}, we note that our choice of terminology is that smoothness of rigid spaces still refers to the classical notion of smoothness, which does not map to Berkovich-smoothness in general. This problem disappears for smooth partially proper rigid spaces, which are rigid smooth, hence Berkovich-smooth, as Berkovich spaces. More generally, if a $\dagger$-rigid space is smooth in the sense of the theory of $\dagger$-rigid spaces, it is Berkovich-smooth when seen as a Berkovich space. 
   
   For example, the open unit disc is rigid smooth (hence also Berkovich-smooth).
   The overconvergent unit disc (a $\dagger$-rigid space) is Berkovich-smooth, but not rigid smooth.
   The closed unit disc, rather its canonical compactification seen as a Berkovich space, is not Berkovich-smooth nor rigid smooth.
    \end{warning}

\subsection{Pushout of Berkovich spaces and $\dagger$-cohesiveness}\label{ssec:pushoutsBerkovichSpaces}

In this subsection, we construct pushouts in Berkovich spaces, at least in a very limited case.

We recall the following definition.

\begin{definition}
  \label{DefinitionZariskiCLosed}
  A map $f\colon Z\to X$ of Berkovich spaces is said to be \textit{Zariski closed} (or a \textit{Zariski closed immersion}) if it is To\"en-affine in the sense of \Cref{sec:deriv-berk-spac-relative-gelfand-spectrum}, and locally on a strict closed cover $\{X_i=\GSpec(A_i)\to X\}_i$ of $X$ with base changes $f_i\colon Z_i\to X_i$ of $f$, the map $A_i\to f_{i,\ast}(\s{O}_{Z_i})$ in $\ob{D}(A_i)$ is surjective on $\pi_0$ and $f_{i,\ast}(\s{O}_{Z_i})$ is connective.
\end{definition}

\begin{example}
  \label{sec:push-berk-spac-non-shriek-descent-for-zariski-closed-immersion}
  Zariski closed immersions of Berkovich spaces do not satisfy $!$-descent: Let $Z\to X$ be the inclusion of a $\dagger$-rational subspace, say a disc of radius $1/2$ into one of radius $1$.
  Let $X'\to X$ be a descendable cover with $X'$ strictly totally disconnected.
  Now, $Z':=Z\times_X X'$ is Zariski closed in the sense of \Cref{DefinitionZariskiCLosed} because it is an inverse limit of inclusions of clopen subsets, but in general $Z\to X$ is not Zariski closed in the sense of \Cref{DefinitionZariskiCLosed}.
\end{example}

\begin{lemma}\label{LemmaPushoutNaiveBerkovich}
Let $f\colon Z\to Y$ and $g\colon Z\to X$ be Zariski closed morphism of Berkovich spaces. Then the pushout $X\bigsqcup_{Z} Y$ exists in the category $\cat{BerkSp}$ of Berkovich spaces, and satisfies the following properties: 
\begin{enumerate}

\item There is a natural equivalence of topological spaces
\[
|X|\bigsqcup_{|Z|} |Y|=|X\bigsqcup_{Z} Y|.
\]

\item Suppose that $X,Y,Z$ are affinoid and $\s{O}(X)\to \s{O}(Z), \s{O}(Y)\to \s{O}(Z) $ are surjective on $\pi_0$, then there is a natural isomorphism
\[
X\bigsqcup_{Z} Y= \GSpec(\s{O}(X)\times_{\s{O}(Z)} \s{O}(Y)). 
\]

\end{enumerate}
\end{lemma}
\begin{proof}
  We first prove the claim when all the spaces are affine with $\s{O}(X)\to \s{O}(Z), \s{O}(Y)\to \s{O}(Z)$ surjective on $\pi_0$.
  Write $Z=\GSpec(A)$, $X=\GSpec (B)$ and $Y=\GSpec (C)$ and define $D= B\times_A C$ and  $W=\GSpec (D)$.
  We first show that the natural map of compact Hausdorff spaces
\begin{equation}\label{eqo93jo3nmlawrqwr}
|X|\bigsqcup_{|Z|} |Y|\to|W|
\end{equation}
is an homeomorphism.  For that, let $x\in |W|$ be a point with associated Banach residue field $K$. This point arises from a map $f\colon D=B\times_A C\to K$. Since $B\to A$ is surjective, so is $D\to C$, and there is a natural equivalence 
\[
\ob{fib}(D\to C)\xrightarrow{\sim} \ob{fib}(B\to A)=I. 
\]
Suppose that $I$ does not vanish in $K$, then there is a $\dagger$-rational localization $D\to D_U$ such that $C\otimes_{D} D_U =0$ and $D\to K$ factors through $D_U$ (it suffices to localize such that an element in $I(\ast)$ that does not vanish in $K(\ast)$ gets invertible). Consequently, the map 
\[
D_U\to B\otimes_{D} D_U
\]
is an equivalence. We deduce that in the complement of the Zariski closed immersion $Y\to W$, the map \eqref{eqo93jo3nmlawrqwr} is a bijection.  Similarly, in the complement of the Zariski closed immersion $X\to W$ the map \eqref{eqo93jo3nmlawrqwr} is a bijection. It is left to show that the natural map $|Z|\to |X|\cap|Y|\subset |W|$ is a bijection, but this is the  vanishing locus of both surjective maps $D\to B$ and $D\to C$ which is precisely the vanishing  locus  of the map $D\to A$.
This proves the claim about topological spaces. 

It is clear that $W$ is a  pushout  when restricted to affinoid Berkovich spaces.  Let $S$ be a general Berkovich space with compatible maps $Y\to S$ and $X\to S$ along $Z$.  Let $W\backslash Z\subset W$ be the open complement of $W$, and let $Z\to V\subset W$ be an open neighbourhood. Then, we have a pullback square
\[
\begin{tikzcd}
\ob{Map}(W, S) \ar[r] \ar[d] & \ob{Map}(W\backslash Z, S ) \ar[d] \\ 
\ob{Map}(V,S) \ar[r] & \ob{Map}(V\cap(W\backslash Z), S) 
\end{tikzcd}
\]
We have an isomorphism $W\backslash Z=(Y\backslash Z) \bigsqcup (X\backslash Z)$. Hence, to show the pushout property of $W$, we are allowed to take strict closed neighbourhoods of $Z$ in $W$.

We can find $\dagger$-rational covers of $Y=\bigcup_i Y_i$ such that the maps $Y_i\to S$ factor through an affinoid subspace, and such that $Z_i=Z \times_Y Y_i\to Y_i$ arises from a surjective map of rings (this uses that morphisms of connective rings that are surjective on $\pi_0$ are stable under base change).
Similarly, we find a $\dagger$-rational cover $X=\bigcup_j X_j$ and $Z_j'=Z\times_{X}X_j\to X_j$ with the same properties.
We can find further $\dagger$-rational covers  by affinoids $Z_i\cap Z_j'=\bigcup Z_{ijk}'$. Since $\s{O}_{X_j}\to \s{O}_{Z_j'}$ is surjective such  $\dagger$-rational covers of $Z_j'$ can be lifted to $\dagger$-rational subspaces $X_{i,j,k}\subset X_{j}$ such that $Z_{i,j,k}=X_{i,j,k}\times_{X_{j}} Z_{j}'$, and similarly for $Y_{i}$.
Shrinking $Y$  and $X$ around $Z$ if necessary, we can  assume that we are given finite $\dagger$-rational covers $Y=\bigsqcup_s Y_s$ and $X=\bigsqcup_s X_s$ inducing the $\dagger$-rational cover of $Z$
\[
Z\times_Y Y_s =Z_s= Z\times_{X} X_s,
\]
and such that, for each $s$, the maps $X_s,Y_s\to S$ factor through the same closed affinoid subspace $S_s\subset S$.
By analytic descent, and the fact that $W_s=Y_s\bigsqcup_{Z_s} X_s$ is a pushout for affinoid Berkovich spaces, we deduce that $W$ is a pushout for all Berkovich spaces, proving what we wanted.

Finally, the construction and pushout property for general $X,Y,Z$ follow by analytic descent and the reasoning of  \cite[Tag 07RT]{StacksProjectAuthors2017}. 
\end{proof}

We have the following variant where we only ask  for an equivalence on uniform completions in one term. 

\begin{lemma}\label{LemmaPushoutsBerkovichSpaces}
  Consider a morphism $g\colon Z=\GSpec(A)\to Z'=\GSpec(A')$ of affinoid Gelfand stacks, such that $A'\to A$ induces an isomorphism on uniform completions.
  Consider a Zariski closed morphism $f\colon Z\to X$ of  Berkovich spaces.
  Then there exists the pushout $X\bigsqcup_{Z} Z'$ in $\cat{BerkSp}$ and it is given by the (qfd) Berkovich space with the same topological space $|X|$ and structural sheaf given by $\s{O}_X\times_{\s{O}_Z} \s{O}_{Z'}$.
  If $X=\GSpec(B)$ and $B\to A$ is surjective on $\pi_0$, then $X\bigsqcup_Z Z'=\GSpec(B\times_A A')$.
\end{lemma}
\begin{proof} 
  It suffices to construct the pushout in the case that $X=\GSpec(B)$ is affinoid with $B\to A$ surjective on $\pi_0$, the general construction follows from analytic descent e.g., by following the reasoning of \cite[Tag 07RT]{StacksProjectAuthors2017}.
  Then, we are in the situation where $X=\GSpec(B)$, $\pi_0(B)\to \pi_0(A)$ is surjective, and $A\to A'$ induces an isomorphism on uniform completions.
  We claim that the pushout is given by $\GSpec (B\times_A A')$.
  We first note that $(B\times_A A')^{\leq 1}=B^{\leq 1} \times_{A^{\leq 1}} A^{\prime \leq 1}$ (as follows from the definition of $(-)^{\leq 1}$ via morphisms from certain idempotent algebras), and consequently that the uniform completion of the bounded $\Q_p$-algebra $B\times_{A}A'$ is $B^u\times_{A^u} A^{\prime,u}\cong B^u$.
  In particular, $B\times_{A}A'$ is a qfd Gelfand ring.

  Now we check the desired universal property for the pushout. 
  Given compatible morphisms $\GSpec(B)\to Y,\ \GSpec(A')\to Y$ with $Y$ a Berkovich space, one finds a strict closed covering $\cup_{i\in I}\GSpec(A'_i)=\GSpec(A')$ by $\dagger$-rational subspaces, such that each $\GSpec(A'_i)$ maps into an affinoid subspace in $Y$.
  Using that $|\GSpec(A)|=|\GSpec(A')|$, one can reduce to the case that $Y$ is affinoid.
  In this case, the assertion is clear.
\end{proof}

We note the following lemma, see \cite[Theorem 7.2]{lurie_dag_IX} for an analogous statement.   
  
\begin{lemma}
  \label{sec:analyt-delta-spac-descent-for-fiber-product}
  Let $A\to B,\ C\to B$ be two morphisms of connective solid $\Q_p$-algebras with $A\to B$ surjective. Then $D:= A\times_{B} C$ is connective, and the following properties hold:
  \begin{enumerate}
  \item The natural functor $\ob{D}(D)\to \ob{D}(A)\times_{\ob{D}(B)} \ob{D}(C)$ obtained by base change is fully faithful.
  \item The natural functor $\ob{D}_+(D)\to \ob{D}_+(A)\times_{\ob{D}_+(B)}\ob{D}_+(C)$ is an equivalence.
    \item Suppose that $A,B$ and $C$ are Gelfand.  Then the natural functor $\Perf(D)\to \Perf(A)\times_{\Perf(B)}\Perf(C)$ is an equivalence.
  \end{enumerate}
 
\end{lemma}
\begin{proof}
The connectivity of $D$ is clear due to the surjectivity of $A\to B$.

\begin{enumerate}

\item    The functor $F\colon \ob{D}(D)\to \ob{D}(A)\times_{\ob{D}(B)} \ob{D}(C)$ has a right adjoint given by sending an element $(N,M)$ in the fiber product to the limit $N\times_{K} M$ where $K=N\otimes_A B= M\otimes_C B$. Note that all functors are linear over $\ob{D}(D)$ (as they consist of base changes  and a finite limit of  forgetful functors of $D$-algebras). Therefore, to prove that $F$ is fully faithful, it suffices to prove that the unit of the adjunction is an equivalence and  it suffices to prove it when evaluating at $D\in \ob{D}(D)$.
  This reduces to the fact that $D= A\times_{B} C$ which holds  by definition. This proves (1).

\item Next we prove (2). By (1) we know that the functor is fully faithful, so we only need to show essential surjectivity.  By shifting if necessary, it suffices to show that the fully faithful functor
\[
\Psi\colon \ob{D}_{\geq 0}(D)\to \ob{D}_{\geq 0}(A)\times_{\ob{D}_{\geq 0}(B)} \ob{D}_{\geq 0}(C)
\]
is an equivalence. Let $N\in \ob{D}_{\geq 0}(A)$, $M\in \ob{D}_{\geq 0}(C)$ together with an identification $\beta\colon B\otimes_{A} N \cong B\otimes_C M$. Let $\Psi^{R}$ denote the right adjoint of $\Psi$.

First, we show  that $K:=\Psi^R(N,M)$ is connective.
Indeed, since $A\to B$ is surjective the cofiber $\ob{cofib}(N\to B\otimes_A N)\cong \ob{cofib}(K\to M)$ is $1$-connective, and as $M$ is $0$-connective one deduces connectivity of $K$.
Next, we claim that the maps $A\otimes_D K\to N$ and $C\otimes_D K\to M$ are surjective on $\pi_0$.
Suppose that this holds and let $(N,M)$ be an object in the fiber product.
We want to see that $Q:=\ob{cofib}(\Psi(\Psi^{R}(N,M))\to (N,M))$ is zero.
Note that, by fully faithfulness of $\Psi$, we have that $\Psi^R(Q)=0$.
Thus, it suffices to show that if $Q\in \ob{D}_{\geq 0}(A)\times_{\ob{D}_{\geq 0}(B)} \ob{D}_{\geq 0}(C)$ is such that $\Psi^R(Q)=0$ then $Q=0$.
Since by our claim the map $0=\Psi(\Psi^{R}(Q))\to Q$ is surjective on $\pi_0$, we know that $Q$ is $1$-connective.
By an inductive argument one deduces that $Q$ is $n$-connective for all $n$ and therefore that $Q=0$ by left completeness of the fiber product.

We finish the proof of (2) by proving the claim, that is, for  $(N,M)\in \ob{D}_{\geq 0}(A)\times_{\ob{D}_{\geq 0}(B)} \ob{D}_{\geq 0}(C)$ with $K:=\Psi^{R}(N,M)$ the maps $A\otimes_{D} K\to N$ and $C\otimes_{D} K\to M$ are surjective on $\pi_0$. By construction, we have that $K= N\times_{B\otimes_{A} N} M$.
Let $S$ be a light profinite set and consider a map $f\colon S\to M$.
Since $A\to B$ is surjective, the map $N\to B\otimes_A N=B\otimes_C M$ is an epimorphism, and as $\Z_{\sol}[S]$ is  projective, there is a lift $S\to N$ of the composite $S\to M\to B\otimes_{C} M\cong B\otimes_A N$. In total, this produces a map $S\to K$ lifting the map $S\to M$. The previous implies that $K\to M$ is surjective, and in particular that $C\otimes_D K\to M$ is also surjective.  We note that, in particular, the map $D\to C$ is surjective. We claim now that the natural map $A\otimes_D K\to N$ is surjective, for that,  denote $Q=\ob{cofib}(K\to M)=\ob{cofib}(N\to B\otimes_{A} N)$.  The pullback square $K=N\times_{B\otimes_A N} M$ in $\ob{D}_+(\bb{Z}_{\sol})$ produces a natural morphism of fiber sequences 
\[
\begin{tikzcd}
K \ar[r] \ar[d]  & M \ar[d] \ar[r] & Q \ar[d]\\ 
N \ar[r] & B\otimes_A N  \ar[r]& Q
\end{tikzcd}
\]
that factors, by base change along $D\to A$, through a morphism of fiber sequences 
\[
\begin{tikzcd}
A\otimes_D K \ar[r] \ar[d]  & A\otimes_D M \ar[d] \ar[r] & A\otimes_D Q \ar[d]\\ 
N \ar[r] & B\otimes_A N  \ar[r]& Q.
\end{tikzcd}
\]
Since $B\otimes_A N = B\otimes_{C} M$ and the maps $A\to B$  and $D\to C$  are surjective,  the map $A\otimes_D M \to B\otimes_{C} M= B\otimes_A N $ is an isomorphism on $\pi_0$.
Indeed, we have
\[
  \pi_0(A\otimes_DM)\cong \pi_0((A\otimes_{D}C)\otimes_{C}M)\cong \pi_0(A\otimes_DC)\otimes_{\pi_0(C)}^{\mathrm{un}}\pi_0(M)\cong \pi_0(B)\otimes_{\pi_0(C)}^{\mathrm{un}} \pi_0(M)
\] with $\otimes^{\mathrm{un}}$ refering to the underived solid tensor product, as $\pi_0(A\otimes_{D}C)\cong \pi_0(A)/\mathrm{ker}(\pi_0(D)\to \pi_0(C))\cong \pi_0(B)$. 
As $K\to M$ is surjective on $\pi_0$, we have that $Q$ is $1$-connective and therefore that $\pi_0(Q)=\pi_0(A\otimes_D Q )=0$ and that the map   $\pi_1(A\otimes_D Q)\to \pi_1(Q)$ is surjective. Thus, by passing to long exact sequences, we get a morphism of exact sequences of solid abelian groups
\[
\begin{tikzcd}
\pi_1(A\otimes_D Q) \ar[r] \ar[d,twoheadrightarrow]& \pi_0(A\otimes_{D} K) \ar[r] \ar[d]& \pi_0(A\otimes_D M) \ar[r] \ar[d,"="]&  0 \ar[d] \\
\pi_1(Q)\ar[r]& \pi_0(N)\ar[r]& \pi_0(B\otimes_A N) \ar[r]& 0 
\end{tikzcd}
\]
The surjectivity of $\pi_0(A\otimes_{D} K) \to \pi_0(N)$ then follows from a diagram chase.

  \item The claim for perfect complexes follows from (2) by passage to dualizable objects, and the Fredholm property of Gelfand rings (\cite[Proposition 3.3.5]{drff_stacks}).\qedhere
  \end{enumerate} 
  \end{proof}

\begin{example} \label{ExampleCounterExample}
\Cref{sec:analyt-delta-spac-descent-for-fiber-product} does not hold for the full category of modules, cf.\ \cite[Warning 7.3]{lurie_dag_IX}. Let $A=C=\bb{Q}_p[T]$ be the polynomial algebra over $\bb{Q}_p$ and $B=\bb{Q}_p$. Consider the map $A\to B$ and $C\to B$ sending $T\mapsto 0$, and let $D=A\times_B C$. Then, let $M=\bigoplus_{n\in 2\bb{Z}} B[n]$ as an $A$-module and $N=\bigoplus_{n\in 2\bb{Z}+1} B[n]$ as a $C$-module. We have a natural equivalence as $B$-modules
\[
M\otimes_{A} B \cong  \bigoplus_{n\in \bb{Z}} B[n]  \cong N\otimes_{C} B. 
\]
Hence, $(M,N)\in \ob{D}(A)\otimes_{\ob{D}(B)} \ob{D}(C)$, but it does not belong to the essential image of the functor $\ob{D}(D)\to \ob{D}(A)\otimes_{\ob{D}(B)} \ob{D}(C)$. \end{example}

The key construction behind the analytic prismatization are $\dagger$-thickenings of Fargues--Fontaine curves. The key examples are obtained by pushouts of classical Fargues--Fontaine curves with thickenings of degree $1$ divisors. It is natural to ask for which stacks this construction still behaves as a pushout.
This motivates the following definition. 

  \begin{definition}
    \label{sec:push-deriv-berk-definition-cohesiveness}
\begin{enumerate}    
    \item A morphism $A\to B$ of Gelfand rings is said to be a \textit{$\dagger$-thickening} if it is surjective and induces an isomorphism on $\dagger$-reductions.   
        
\item    Let $f\colon Y\to X$ be a morphism of Gelfand stacks.  We say that $f$ is \textit{$\dagger$-cohesive} if  for a $\dagger$-thickening map of Gelfand rings $B\to A$ and a surjection $C\to A$ of Gelfand rings, the following diagram is cartesian

   \[
   \begin{tikzcd}
   Y(B\times_{A} C)  \ar[r] \ar[d] &  Y(B)\times_{Y(A)} Y(C) \ar[d] \\
  X(B\times_{A} C) \ar[r] & X(B)\times_{X(A)} X(C) 
   \end{tikzcd}
   \]
   \end{enumerate}
     \end{definition}

In the following we give a list of basic examples.

  \begin{example}
    \label{sec:push-deriv-berk-examples-for-dagger-cohesiveness}
    \begin{enumerate}
    \item Any morphism of Berkovich spaces is $\dagger$-cohesive by \Cref{LemmaPushoutsBerkovichSpaces}.
    
\item The functor $\Perf(-)\colon \cat{GelfRing}\to \cat{Cat}_{\infty}$ is $\dagger$-cohesive (in the obvious variant for sheaves valued in categories)    by \Cref{sec:analyt-delta-spac-descent-for-fiber-product}.

\item Let $R$ be a Gelfand ring and  $\n{X}$ be a perfect strongly Tannakian stack over $\Spec R(*)$.
  Then $\n{X}^{\ob{an}}\to \GSpec(R)$ is $\dagger$-cohesive.
  Indeed, this follows from (2) above and the Tannakian description of $\n{X}^{\ob{an}}$-points of \Cref{PropositionFunctorOfPointsAnalytification} (the connectivity assumption in \Cref{PropositionFunctorOfPointsAnalytification} can easily be verified, as it can be checked $!$-locally by \cite[Proposition VI.1.2]{real_local_langlands}).
  This applies in particular to classifying stacks of algebraic groups. 

    \end{enumerate}
  \end{example}
  
  The next lemma shows some basic stability properties of $\dagger$-cohesive maps. 
  
  \begin{lemma}
The following hold.  
  
\begin{enumerate}

\item Let $Z\xrightarrow{g} Y \xrightarrow{f} X$ be maps of Gelfand stacks with $f$ $\dagger$-cohesive. Then $g$ is $\dagger$-cohesive if and only if $g\circ f$ is $\dagger$-cohesive.

\item Let $f\colon Y\to X$ and $g\colon X'\to X$ be maps of Gelfand stacks with base change $f'\colon Y'\to X'$.  Suppose that $f$ is $\dagger$-cohesive, then $f'$ is $\dagger$-cohesive.

\item Let $f\colon Y\to X$ be a $\dagger$-formally \'etale map of Gelfand stacks, then $f$ is $\dagger$-cohesive. 

\end{enumerate}  
  \end{lemma}
  \begin{proof}
  \begin{enumerate}
  
	\item  Let $B\to A$ be a $\dagger$-thickening of Gelfand rings and $C\to A$ a surjection.   Denote $D=B\times_A C$. For a Gelfand stack $W$ let us define the anima  $W^{\dagger-\ob{coh}}(D):=W(B)\times_{W(A)} W(C)$. For a morphism of Gelfand stacks $W\to V$ write $W^{\dagger-\ob{coh}/V}(D):= W^{\dagger-\ob{coh}}(D)\times_{V^{\dagger-\ob{coh}}(D)} V(D)$.  By definition, $W\to V$ is $\dagger$-cohesive if for all rings $A,B,C$ as above the natural map $W(D)\to W^{\dagger-\ob{coh}/V}(D)$ is an equivalence of anima. 
	
	We have a diagram with pullback squares
	\[
	\begin{tikzcd}
	 Z(D) \ar[d] & &  \\ 
Z^{\dagger-\ob{coh}/Y}(D) \ar[r]  \ar[d] &  Y(D) \ar[d]&  \\
Z^{\dagger-\ob{coh}/X}(D)\ar[r] \ar[d] & Y^{\dagger-\ob{coh}/X}(D)\ar[r] \ar[d] &  X(D)\ar[d] \\
Z^{\dagger-\ob{coh}}(D) \ar[r]  & Y^{\dagger-\ob{coh}}(D) \ar[r] & X^{\dagger-\ob{coh}}(D)  
	\end{tikzcd}
	\]
	By assumption, the map $Y(D)\to Y^{\dagger-\ob{coh}/X}(D)$ is an equivalence, so that $Z^{\dagger-\ob{coh}/Y}(D)\to Z^{\dagger-\ob{coh}/X}(D)$ is also an equivalence. Hence, $Z(D)\to Z^{\dagger-\ob{coh}/Y}(D)$ is an equivalence if and only if $Z\to Z^{\dagger-\ob{coh}/X}(D)$ is so, which proves part (1).
	
	\item Keep the notation of part (1). The statement follows from the fact that
	\[
	Y^{\dagger-\ob{coh}}(D)\times_{X^{\dagger-\ob{coh}}(D)} X^{\prime,\dagger-\ob{coh}}(D) = Y^{\prime,\dagger-\ob{coh}}(D), 
	\]
	so that 
	\[
	\begin{aligned}
	 Y^{\prime,\dagger-\ob{coh}/X'}(D) & = Y^{\dagger-\ob{coh}}(D)\times_{X^{\dagger-\ob{coh}}(D)} X^{\prime,\dagger-\ob{coh}}(D) \times_{X^{\prime,\dagger-\ob{coh}}(D)} X'(D) \\ 
	 & = Y^{\dagger-\ob{coh}}(D)\times_{X^{\dagger-\ob{coh}}(D)} X'(D)  \\ 
	 & =  Y^{\dagger-\ob{coh}}(D)\times_{X^{\dagger-\ob{coh}}(D)} X(D)\times_{X(D)} X'(D) \\ 
                                           & = Y^{\dagger-\ob{coh}/X}(D)\times_{X(D)} X'(D) \\
                                           & = Y(D)\times_{X(D)}X'(D) \\
          & = Y'(D).
	\end{aligned}
	\]
	
	\item   Suppose that $f\colon Y\to X$ is $\dagger$-formally \'etale, that is, the map $Y\to Y^{\dR/X}$ is an equivalence.   By (2), it suffices to show that $Y^{\dR}\to X^{\dR}$ is $\dagger$-cohesive, and by (1) even just that $Y^{\dR}\to \GSpec(\Q_p)$ is $\dagger$-cohesive.  Let $A,B,C, D$ be as in (1), the statement then follows from the fact that  for $R'\to R$ a $\dagger$-thickening the natural map $Y^{\dR}(R')\to Y^{\dR}(R)$ is an equivalence, and the fact that $B\to  A$ and hence $D\to C$ are $\dagger$-thickenings. \qedhere
  \end{enumerate}
  \end{proof}

\subsection{The cone construction: geometry}
\label{sec:cone-construction}

In this section, we introduce the innocent looking \emph{cone construction} of a Gelfand stack.
It will play a prominent role in our discussion of the analytic prismatization. The main features of the cone construction in algebraic geometry can be found in \cite{hekking2025deformationnormalbundleblowups} (where in \textit{loc. cit.} this construction is called \textit{deformation to the normal bundle}).

\begin{definition}
  \label{sec:defin-first-prop-1-definition-cone-construction}
  Let $Y$ be a Gelfand stack.
  
  \begin{enumerate}
\item   We define \textit{the cone} $Y^\cone\to \A^1/\Gm$ of $Y$ to be the functor sending a Gelfand stack $Z$ with a Cartier divisor $D\subseteq Z$ to $Y(D)$.

\item We define the \textit{overconvergent cone} as $Y^{\dagger-\cone} = Y^\cone\times_{\A^1/\Gm} \A^{1,\dagger}/\bb{G}_m$. 

\item We define the \textit{formal cone} as $Y^{\widehat{\cone}} = Y^\cone\times_{\A^1/\Gm} \widehat{\mathbb{A}}^1/\bb{G}_m$.

\item   We define the \textit{special fiber of the cone} as $Y^{\cone}_{s}:=Y^{\cone}\times_{\A^1/\Gm} B\Gm$.   
\end{enumerate}

Given $f\colon Y\to X$ a morphism of Gelfand stacks, we define the \textit{relative cone of $f$} to be  $Y^{\cone/X}:=Y^{\cone}\times_{X^{\cone}} (X\times \bb{A}^1/\bb{G}_m)$.
Here, $X\times \bb{A}^1/\Gm\to X^\cone$ is the morphism sending a Gelfand stack $Z$ over $X$ with a Cartier divisor $D\to Z$, i.e., a Gelfand stack $Z$ over $X\times \bb{A}^1/\Gm$, to the composition $D\to Z\to X$.
\end{definition}

The cone construction is a well defined Gelfand stack. 

\begin{lemma}\label{LemmaConeConstructionWellDefined}
  Let $Y$ be a Gelfand stack, $\GSpec (R)\to \bb{A}^1/\bb{G}_m$ a Gelfand ring endowed with a Cartier divisor $I\to R$,  and let $S\to \GSpec(R)$ be a $!$-equivalence of presheaves on Gelfand rings over $\bb{A}^1/\bb{G}_m$.  Let $Z(S):= S\times_{\bb{A}^1/\bb{G}_m} B\bb{G}_m$ be the induced Cartier divisor of $S$.
  Then the natural map of anima 
\[
Y(R/I)\to Y(Z(S))
\]
is an equivalence.
In particular, $Y^{\cone}$ is a Gelfand stack.
\end{lemma}
\begin{proof}
This follows from the fact that $!$-equivalences are stable under base change, hence if $S\to \GSpec (R)$ is a $!$-equivalence of presheaves over $\bb{A}^1/\bb{G}_m$, then  $Z(S)\to \GSpec (R/I)$ is also a $!$-equivalence of presehaves over $B\bb{G}_m$ (and thus of presheaves over $\GSpec(\Q_p)$). 
\end{proof}

\begin{remark}\label{RemarkMapSpecialFiberConeToX}
  Let $X$ be a Gelfand stack, then there is a natural map $X^{\cone}_s\to X\times B\bb{G}_m$.
  Indeed, given $R$ a Gelfand ring and $I$ a line bundle on $\GSpec(R)$, the map above arises from  evaluating $X$ at the natural projection $R\oplus I[1]\to R$, where $R\oplus I[1]$ is the trivial square-zero extension.
\end{remark}

\begin{remark}\label{RemarkToposTeoreticDescriptionCone}
In the following we give a topos theoretic description of the cone construction.  Let $\s{X}=\cat{GelfStk}$ be the topos of Gelfand stacks. We have maps of Gelfand stacks
\[
\GSpec(\Q_p) \xleftarrow{f} B\bb{G}_m \xrightarrow{g} \bb{A}^1/\bb{G}_m.
\]
Then, by unraveling the definitions, the cone construction is nothing but the map
\[
g_*f^*\colon \s{X} \to \s{X}_{/(\bb{A}^1/\bb{G}_m)},
\]
where $f^*$ is the pullback along $B\bb{G}_m\to \GSpec(\Q_p)$, and $g_*$ is the right Kan extension along the map $B\bb{G}_m\to \bb{A}^1/\bb{G}_m$.  In particular, the formation of the cone commutes with limits of Gelfand stacks as $g_\ast$ and $f^\ast$ do (note that $f^\ast$ is just the restriction to a slice). 
\end{remark}

\begin{example}\label{ExamplesOfCones}
\begin{enumerate}

\item Consider the ring stack $\bb{G}_a$ over $\GSpec(\Q_p)$ and its base change $\bb{G}_{a,\bb{A}^1/\bb{G}_m}$ to $\A^1/\Gm$.
  Over  $\bb{A}^1/\bb{G}_m$, the universal Cartier divisor produces a map of $\bb{G}_{a,\bb{A}^1/\bb{G}_m}$-modules 
\[
\bb{G}_a(-1)\to \bb{G}_{a,\bb{A}^1/\bb{G}_m},
\]
where $\bb{G}_a(-1)$ is the vector bundle associated to the ideal $I\to R$ (for an $R$-point).
Then the cone of $\bb{G}_a$ is the quotient
\[
\bb{G}_a^{\cone}= \bb{G}_{a,\bb{A}^1/\bb{G}_m}/\bb{G}_a(-1). 
\]

\item Let $f\colon Z\to X$ be a Zariski closed immersion of Berkovich spaces. Then \cite[Remark 7.7.17]{solidnotes} implies that $Z^{\cone/X}$ is the analytification of the \textit{deformation to the normal bundle} of $f$ given by blowing up (in an appropriate sense for solid rings) the Zariski closed map $Z\times B\bb{G}_m\to X\times \bb{A}^1/\bb{G}_m$.  Indeed, suppose for simplicity that $Z\to X$ is locally given by an lci Zariski closed immersion.  As we will see in \Cref{sec:cone-construction-1-geometry-of-cone-construction}, over the locus $\bb{G}_m/\bb{G}_m$ the stack $Z^{\cone/X}$ is isomorphic to $X$, and the pullback along the zero map $B\bb{G}_m\to \bb{A}^1/\bb{G}_m$ is the shifted tangent bundle $\bb{T}_{Z/X}[1](-1)$ (with a twist by $-1$ for the $\bb{G}_m$-action), equivalently after forgetting the twist, it is the \textit{normal cone} $N_{Z/X}$ given by the  (analytification of the) relative spectrum of the symmetric algebra $\Sym_{\s{O}_Z} (N_{Z/X}^*)$ with $N_{Z/X}^*=\bb{L}_{Z/X}[-1]$ being the conormal cone, or shifted cotangent complex, of $Z\to X$ (for cotangent complexes for morphisms of Gelfand stacks, see \Cref{app:DefTheoryGelfand}). 
\end{enumerate}
\end{example}

The slightly strange geometry of the cone construction is described in \Cref{sec:cone-construction-1-geometry-of-cone-construction}.

\begin{lemma}
  \label{sec:cone-construction-1-geometry-of-cone-construction}
  Let $f\colon Y\to X$ be a  morphism of Gelfand stacks. The following hold:
  \begin{enumerate}
  \item The natural map  $Y^{\cone/X}\times_{\A^1/\Gm}\Gm/\Gm\to X$ is an isomorphism.
  \item Assume that $f$ admits a cotangent complex, see \Cref{sec:cotang-compl-morph-definition-having-a-cotangent-complex}.   Then there is a natural equivalence  of stacks over $Y\times B\bb{G}_m$
  \[
  Y^{\cone/X}_s= \bb{T}_{Y/X}[1](-1)=\bb{T}_{Y/X}[1]\otimes_{\bb{G}_a} \bb{G}_a(-1),
  \]
  where $ \bb{T}_{Y/X}[1]$ is the $1$-shifted tangent bundle of $f$ (cf.\ \Cref{ConstructionTangentSpaces}), and the $(-1)$ refers to the action of $\bb{G}_m$ on the shifted tangent space via multiplication by the character $t\mapsto t^{-1}$. 
  Equivalently, $\bb{T}_{Y/X}[1](-1)$ is the twist of the tangent space by the universal ideal defining the Cartier divisor over $\bb{A}^1/\bb{G}_m$.
  \end{enumerate}
\end{lemma}
\begin{proof}
\begin{enumerate}
\item  It suffices to prove the statement when $X=\GSpec(\Q_p)$, where $Y^{\cone/X}=Y^\cone$.
  If $\GSpec(R)\to \A^1/\Gm$ factors over $\Gm/\Gm$, then the classified Cartier divisor $s\colon I \to R$ satisfies that $s$ is an isomorphism.
  Thus, $R/I=0$ which implies $Y(R/I)=\ast$ as desired.

\item  If $\GSpec(R)\to \A^1/\Gm$ is classifying $s\colon I\to R$, then a factorization over $B\Gm$ yields naturally a zero homotopy of $s$, making the quotient $R/I= R\oplus I[1]$ a trivial square-zero extension of $R$ by $I[1]$.
  The morphism of rings $R\oplus I[1]\to R$ induced by the first projection yields the natural map $Y_s^{\cone/X}\to Y\times B\bb{G}_m$.

Furthermore, we have the following functor of points description: given $\eta\colon \GSpec(R) \to X\times B\bb{G}_m$ one has 
\[
\begin{aligned}
Y^{\cone/X}_s(\GSpec(R)\to X\times B\bb{G}_m)  & = Y(R\oplus I[1])_{/\eta} \\
&  = \tau_{\geq 0} \Hom_{R}(\eta^* \bb{L}_{Y/X}, I[1]) \\ 
&  = \tau_{\geq -1} \Hom_{R}(\eta^* \bb{L}_{Y/X}, R) \otimes_R I  \\ 
& = \bb{T}_{Y/X}[1](-1)(\GSpec(R)\to X\times B\bb{G}_m) 
\end{aligned} 
\]
  proving the desired equivalence.  \qedhere
  \end{enumerate}
  \end{proof}

Overconvergent cones have the following behaviour with respect to $\dagger$-formally \'etale maps. 

\begin{lemma}\label{LemmaFormalConePullback}
Let $Y\to X$ be a  $\dagger$-formally \'etale map of Gelfand stacks. Then the natural commutative diagram
\[
\begin{tikzcd}
Y^{\dagger-\ob{cone}}  \ar[r] \ar[d]& X^{\dagger-\ob{cone}} \ar[d] \\ 
(Y^{\dagger-\ob{cone}} )^{\dR}\ar[r] \ar[d] &   (X^{\dagger-\ob{cone}} )^{\dR} \ar[d]\\ 
Y^{\dR} \ar[r] &  X^{\dR}
\end{tikzcd}
\]
has cartesian squares. 
\end{lemma}
The same assertion is true for the formal cones, by the same argument (where $I\to R$ is then actually a nilpotent ideal).
\begin{proof}
Let $R$ be a nilperfectoid ring endowed with a Cartier divisor $I\to R$ landing in $\ob{Nil}^{\dagger}(R)$. We have to show that the following diagram has cartesian squares
\[
\begin{tikzcd}
Y(R/I) \ar[r] \ar[d] & X(R/I) \ar[d] \\ 
Y(R^u/(I\otimes_{R} R^u))  \ar[d] \ar[r] & X(R^u/(I\otimes_{R} R^u))  \ar[d] \\ 
Y(R^u) \ar[r] & X(R^u). 
\end{tikzcd}
\]
But this follows from $\dagger$-formally \'etaleness of $Y\to X$, and the fact that $R^u=(R/I)^u=(R^u/ I\otimes_{R} R^u)^u$ as $I\to R$ is $\dagger$-nilpotent. 
\end{proof}

We can define cones relative to Cartier divisors. First, we make the following observation.

\begin{remark} \label{RemarkConeAdjunction}
Let $W\to \bb{A}^1/\bb{G}_m$ be a Gelfand stack with associated Cartier divisor $D\to W$.  Then there is a natural map 
\[
W\to D^{\cone}
\]
 over $\bb{A}^1/\bb{G}_m$.
Concretely, the map $W\to D^{\cone}$ over $\bb{A}^1/\bb{G}_m$ sends a morphism $Z\to W$ to the map $Z\times_{\bb{A}^1/\bb{G}_m} B\bb{G}_m = Z\times_{W} D \to D$ given by the projection to $D$. 
\end{remark}

\begin{definition}\label{DefinitionRelativeConesCartierDivisors}
Let $W\to \bb{A}^1/\bb{G}_m$ be a Gelfand stack endowed with a Cartier divisor $D\to W$, and  let $Y\to D$ be a morphism of Gelfand stacks. The \textit{relative cone of $Y$ over $(D,W)$} is defined as 
\[
Y^{\cone/(D,W)}:= Y^{\cone}\times_{D^{\cone}} W,
\] 
where the map $W\to D^{\cone}$ is as in \Cref{RemarkConeAdjunction}.  We define the overconvergent, formal and special fiber cone in a similar fashion, by base changing the respective locus in $\A^1/\Gm$. 
\end{definition}

\begin{remark}
  \label{sec:cone-constr-geom-relation-of-relative-cone-constructions}
  We note that the relative cone $Y^{\cone/X}$ for a morphism $Y\to X$ in \Cref{sec:defin-first-prop-1-definition-cone-construction} agrees with the relative cone $(Y\times B\bb{G}_m)^{\cone/(X\times B\bb{G}_m, X\times \A^1/\Gm)}$ with respect to the divisor $X\times B\bb{G}_m$.
\end{remark}

The geometry of the relative cone construction with respect to a Cartier divisor is similar to \Cref{sec:cone-construction-1-geometry-of-cone-construction}. Moreover, the special fiber of the  cone relative to a divisor  is intimately related with the problem of lifting a Gelfand stack to the natural square-zero thickening of the divisor. In order to make the previous ideas precise, we need the following definition.

\begin{definition}\label{DefinitionsSquareZeroExtension}
\begin{enumerate}
\item Let $Y\to X$ be a morphism of Gelfand stacks admitting a cotangent complex  (\Cref{sec:cotang-compl-morph-definition-having-a-cotangent-complex}). 
  Let $M\in \ob{D}_+(Y)$ be a locally eventually connective sheaf.
  The anima of \textit{square-zero extensions of $Y$ over $X$ with ideal $M$} is defined as $\ob{Sqz}(Y/X;M):=\tau_{\geq 0} \Hom_{Y}(\bb{L}_{Y/X}, M[1])$. 

\item  Let $Y\xrightarrow{f} X \to S$ be morphisms of Gelfand stacks admitting a cotangent complex, $M\in \ob{D}_+(X)$, and let $\eta\in \pi_0(\Hom_{X}(\bb{L}_{X/S}, M[1]))$ be a square-zero extension of $X$ over $S$. The anima of \textit{square-zero lifts of $Y$ over $\eta$} is the fiber  
\[
\ob{Sqz(Y/S;\eta)}=\ob{fib}(\tau_{\geq 0}\Hom_{Y}(\bb{L}_{Y/S}, f^*M[1])\to \tau_{\geq 0}\Hom_{Y}(f^*\bb{L}_{X/S}, f^*M[1]))
\]
at the induced point by $f^\ast\eta$ on the right term (we note that despite the notation, this anima does not depend on $S$).

\end{enumerate}

\end{definition}

\begin{example}\label{ExampleDivisorCase}
\begin{enumerate}

\item  Suppose that $Y\to X$ is a map of affinoid Gelfand stacks with associated morphism of rings $A\to B$, and let $N\in \ob{D}_{\geq 0}(B)$ be a connective $B$-module.
  Then $\ob{Sqz}(Y/X;M)=\tau_{\geq 0}(\Hom_B(\bb{L}_{B/A}, M[1]))$ is precisely the anima of square-zero extensions $\widetilde{B}\to B$ as $A$-algebras of $B$ with ideal $M$ (see \cite[Section 5.1]{cesnavicius2019purity}).
  Here, $\widetilde{B}=B\times_{B\oplus M[1]}B$ is the fiber product of the zero morphism $(\mathrm{Id}_B,0)\colon B\to B\oplus M[1]$ with the morphism classified by the given morphism $\bb{L}_{B/A}\to M[1]$.

\item Let $Y\to X \to S$ be maps of Gelfand stacks admitting a cotangent complex, let $M\in \ob{D}_+(X)$, and let $\eta\in \ob{Sqz}(X/S;M)$ be a square-zero extension of $X$ with ideal $M$ over $S$.  We have a fiber sequence of spectra
\[
\Hom_{Y}(\bb{L}_{Y/X}, f^\ast M[1])\to \Hom_{Y}(\bb{L}_{Y/S}, f^\ast M[1])\to \Hom_Y(f^*\bb{L}_{X/S}, f^\ast M[1]).
\]
Thus, the anima $\ob{Sqz}(Y/S;\eta)$ of square-zero lifts of $Y$ over $\eta$  is non-empty if and only if the image of $f^\ast\eta$ along 
\[
\pi_{-1}(\Hom_{Y}(f^*\bb{L}_{X/S}, f^*M))\to \pi_{-2}(\Hom_{Y}(\bb{L}_{Y/X}, f^*M))
\]
is zero.
If this is the case, then the anima of lifts is  naturally a $\tau_{\geq 0} \Hom_{Y}(\bb{L}_{Y/X}, f^\ast M[1])$-torsor.
In general, we say that  $\ob{Sqz}(Y/S;\eta)$ is a \textit{principally homogeneous space} of the group  $\tau_{\geq 0} \Hom_{Y}(\bb{L}_{Y/X}, f^\ast M[1])$. 

If $M$ is connective, $Y=\GSpec(B)\to X=\GSpec(A)\to S=\Spec(R)$ are affinoid and $\eta$ is given by the square-zero extension $\widetilde{A}\to A$ (as an $R$-algebra) with ideal $M$, then the anima $\ob{Sqz}(Y/S;\eta)$ is given by square-zero lifts $\widetilde{B}\to B$ (as an $R$-algebra) with ideal $B\otimes_A M$ with a morphism $(\widetilde{A}\to A)\to (\widetilde{B}\to B)$ of square-zero extensions of $R$-algebras, which induces the natural map $M\to B\otimes_AM$ on ideals.
Here, a morphism of square-zero extensions is the datum of a commutative square of Gelfand rings 
\[\begin{tikzcd}
	A & B \\
	{A\oplus M[1]} & {B\oplus (B\otimes_A M)[1],}
	\arrow[from=1-1, to=1-2]
	\arrow["{\delta_{\widetilde{A}}}"', from=1-1, to=2-1]
	\arrow["{\delta_{\widetilde{B}}}", from=1-2, to=2-2]
	\arrow[from=2-1, to=2-2]
      \end{tikzcd}\]
    where $\delta_{\widetilde{A}}, \delta_{\widetilde{B}}$ are the morphisms classifying the square-zero extensions $\widetilde{A}$ and $\widetilde{B}$.
 \item The map $s\colon  B\bb{G}_m\to \bb{A}^1/\bb{G}_m$ has a natural square-zero extension given by the relative spectrum of $\bb{Q}_p[t]/t^2$, where $t$ is the coordinate ot $\bb{A}^1/\bb{G}_m$ (which has weight $-1$ for the $\bb{G}_m$ action). Concretely, this square-zero extension is associated to the identity map $\eta\colon N_{s}^*\to N_s^*$ of the conormal cone $N_s^*=\bb{L}_{B\bb{G}_m/(\bb{A}^1/\bb{G}_m)}[-1]$, which is isomorphic to $\s{O}(-1)$. In particular, if $g\colon W\to \bb{A}^1/\bb{G}_m$ is a Cartier divisor on a Gelfand stack $W$, the divisor $D\to W$ has a natural square-zero extension $g^*\eta$ that is associated to the first nilpotent thickening of $D\to W$ given by $D^{(1)}=\underline{\GSpec}_{\bb{A}^1/\bb{G}_m}(\bb{Q}_p[t]/t^2)\otimes_{\bb{A}^1/\bb{G}_m} W$. 

\end{enumerate}
\end{example}

Under some reasonable geometric conditions, we can produce an \textit{honest} geometric square-zero thickening from a square-zero extension as in \Cref{ExampleDivisorCase}.  We first deal with the case of Berkovich spaces.

\begin{lemma}\label{LemmaBerkovichSpacesSquareZero}
 Let  $f\colon Y\to X$ be a morphism of Berkovich spaces, and let $M\in \ob{D}_+(Y)$ be a module that, locally on a strict closed cover of $Y$ by affinoids, becomes connective. Let $\eta\in \ob{Sqz}(Y/X;M)$ be a square-zero extension of $Y$ with ideal $M$ over $X$. Let  $(\id,d)\colon \s{O}_{Y}\to \s{O}_{Y}\oplus M[1]$ be the associated morphism of rings on $|Y|_{\Betti}$ associated to the map $\eta$, and let $\s{O}_{\widetilde{Y}_{\eta}}:= \s{O}_{Y}\times_{\s{O}_{Y}\oplus M[1]}\s{O}_Y$ as algebras over $|Y|_{\Betti}$, where the left map  in the fiber product corresponds to the trivial derivation and the right map to $d$. Then, locally in the analytic topology of $Y$, the ring $\s{O}_{\widetilde{Y}_{\eta}}$ is connective and has a natural structure of animated $\s{O}_X$-algebra. We call  $\widetilde{Y}_{\eta}:=\underline{\GSpec}_{|Y|_{\Betti}}(\s{O}_{\widetilde{Y}_{\eta}})$ the \textit{underlying square-zero thickening of $\eta$}.
\end{lemma}
\begin{proof}
In case $Y\to X$ is a morphism of affinoid Gelfand stacks for which $M$ is connective, then the statement about $\s{O}_{\widetilde{Y}_{\eta}}$ is clear, that is, it is an $\s{O}_X$-algebra given by the underlying algebra of a square-zero extension of $\s{O}_Y$ by $M$. To deduce the construction for general Berkovich spaces, it suffices to show that the formation of $\widetilde{Y}_{\eta}$ is compatible with analytic descent in both $Y$ and $X$. This is clear since $\dagger$-rational  localizations are formally \'etale. 
\end{proof}

The following lemma shows some stability properties of the formation of the underlying square-zero thickening for Berkovich spaces. 

\begin{lemma}\label{LemmaStabilitySquareZeroGeometrically}
  Let $f\colon Y\to X$ be a map of Berkovich spaces, $M\in \ob{D}_+(Y)$ a module that locally on a strict cover of $Y$ is connective, and $\eta\colon \s{O}_Y\to M[1]$ a square-zero extension of $Y$ over $X$ with ideal $M$.
  Let $Y\to \widetilde{Y}_{\eta} \to X$ be the underlying square-zero thickening of $Y$ associated to $\eta$. The following hold: 
  \begin{enumerate}

\item  Let $g\colon X'\to X$ be a morphism with base change $f'\colon Y'\to X'$ and $g'\colon Y'\to Y$. Let $g^{\prime,*}M$ and $g^{\prime,*}\eta$ be the pullback to $Y'$ of the module and the square-zero extension. Then the following natural square is cartesian
\[
\begin{tikzcd}
\widetilde{Y'}_{g^{\prime,*}\eta} \ar[r] \ar[d] & X' \ar[d] \\ 
\widetilde{Y}_\eta \ar[r] & X.
\end{tikzcd}
\]

\item  Let $g\colon Y'\to Y$ be a $\dagger$-formally \'etale morphism, and consider the square-zero extension $g^*\eta$ of $Y'$ by $g^*M$ over $X$. Then the following natural diagram has cartesian squares
\[
\begin{tikzcd}
Y' \ar[r]\ar[d] & Y\ar[d]  \\
\widetilde{Y'}_{g^*\eta} \ar[r] \ar[d]& \widetilde{Y}_{\eta} \ar[d] \\  
Y^{\prime,\dR} \ar[r] & Y^{\dR} 
\end{tikzcd}
\]

\end{enumerate}
\end{lemma}
\begin{proof}
\begin{enumerate}

\item Arguing as in \Cref{LemmaStabilitySquareZeroGeometrically} one reduces to the fact that square-zero extensions $\widetilde{B}\to B$ over $A$ of  maps of Gelfand rings  $A\to B$ are compatible with base change along $A\to A'$, which is obvious from construction.

\item Similarly as before, one can reduce to the affine case.
  Recall that  $\dagger$-formally \'etale maps are in particular formally \'etale.  Thus, one reduces the statement to the fact that the formation of a square-zero extension $\widetilde{B}\to B$ over $A\to B$  is compatible with base change along formally \'etale maps maps $B\to B'$. This proves the fact that the first square is cartesian. Now, since $Y'\to  Y$ is  $\dagger$-formally \'etale the outer square is cartesian. Furthermore, the map $Y\to \widetilde{Y}_{\eta}$ is a $!$-cover being locally a nilpotent thickening (and hence descendable). The pasting lemma then implies that the bottom square is cartesian. \qedhere
\end{enumerate}
\end{proof}

Going back to the relation between the cone construction and the problem of lifting a morphism, we have the following proposition.

\begin{proposition}
  \label{sec:cone-construction-1-geometry-of-relative-cone-construction}
  Let $S$ be a Gelfand stack, let  $D\to W$ be a Cartier divisor on a Gelfand stack $W$  over $S$,  and $f\colon Y\to D$ a morphism of Gelfand stacks over $S$. Consider the relative cone $Y^{\cone/(D,W)}\to W$. The following hold:
  \begin{enumerate}
  \item The map  $Y^{\cone/(D,W)}\to W$ is an isomorphism on the complement $W\backslash D$ of the Cartier divisor.

  \item  Assume that the maps $Y\to D\to W\to S$ are infinitesimally cohesive (\Cref{sec:cotang-compl-morph-definition-having-a-cotangent-complex}) and admit a cotangent complex. Then there is a natural pullback square 
  \[
  \begin{tikzcd}
  Y^{\cone/(D,W)}_{s} \ar[r] \ar[d] &  D \ar[d] \\ 
  		\bb{T}_{Y/S}[1](-1) \ar[r]			& \bb{T}_{D/S}[1](-1),
  \end{tikzcd}
  \]
where the map  $\eta\colon D\to \bb{T}_{D/S}[1](-1)$ corresponds to the  natural map $\bb{L}_{D/S}\to \bb{L}_{D/W}$ under the equivalence   
\[
\Hom_{D}(\bb{L}_{D/S}, \bb{L}_{D/W})= \Hom_{D}(\bb{L}_{D/S},g^* \n{I}[1])
\] with $\n{I}\to \s{O}_W$ the underlying ideal of $D$ in  $W$ and $g\colon D\to W$.

\item     The natural map $Y_{s}^{\cone/(D,W)}\to Y$ is a principally homogeneous stack for $\bb{T}_{Y/D}[1]\otimes_{\bb{G}_a} N_{D/W}^*$. Moreover, for $j\colon U\to Y$ a $\dagger$-formally \'etale map   the anima of sections  $U\dashrightarrow Y_{s}^{\cone/(D,W)}$ is precisely the anima $\ob{Sqz}(U/S, j^*\eta)$ of square-zero lifts of $U$ over $S$ from $D$ to the first nilpotent thickening $D^{(1)}$ (cf. \Cref{ExampleDivisorCase}).
 
  \end{enumerate}
\end{proposition}
The morphism $Y^{\cone/(D,W)}_s\to Y$ is the obstruction class for lifting $Y$ to the first infinitesimal neighborhood of $D$ in $W$ in the sense that sections $Y\to Y^{\cone/(D,W)}_s$ correspond to the liftings.
\begin{proof}
Part (1) is a direct consequence of \Cref{sec:cone-construction-1-geometry-of-cone-construction} as the locus $W\backslash D$ is by definition the pullback $W\times_{\bb{A}^1/\bb{G}_m} \bb{G}_m/\bb{G}_m$.  

For part (2), the pullback description follows directly from the definitions,  the fact that the divisor of $W$ is given by a map $W\to S\times \bb{A}^1/\bb{G}_m$, and from \Cref{sec:cone-construction-1-geometry-of-cone-construction}.  The fact that $D\to D^{\cone/S}_{s}= \bb{T}_{D/S}[1](-1)$ is given by the natural map   $\bb{L}_{D/S}\to g^* \n{I}[1]=\bb{L}_{D/W}$  of cotangent complexes on $D$  follows from the description of the natural square-zero thickening $D\to D^{(1)}\to W$ as in \Cref{ExampleDivisorCase} (3). For the claim about being a principally homogeneous stack over $Y$, we note that there is a fiber sequence of shifted tangent spaces over $Y\times B\bb{G}_m$
\[
 \bb{T}_{Y/D}[1](-1)\to  \bb{T}_{Y/S}[1](-1)\to f^*\bb{T}_{D/S}[1](-1)
\]
produced from the fiber sequence of cotangent complexes  $f^*\bb{L}_{D/S}\to \bb{L}_{Y/S}\to \bb{L}_{Y/D}$.  This yields the principally homogeneous property of $Y_{s}^{\cone/(D,W)}$ after base change along $\eta$.

Finally, for part (3), the above discussion shows that  the $Y$-points of $Y_{s}^{\ob{cone}/(D,W)}$ are precisely the fiber 
\[
\ob{fib}(\tau_{\geq 0}\Hom_{Y}(\bb{L}_{Y/S}, f^*N_{D/W}^*[1])\to \Hom_{Y}(f^*\bb{L}_{D/S}, f^*N_{D/W}^*[1])) 
\]
at the point $\eta$ corresponding to the pullback of the natural map of cotangent complexes $\bb{L}_{D/S}\to \bb{L}_{D/W}=N_{D/W}^*[1]$. This is, by definition,  the anima of square-zero extensions of $Y$ over $\eta$ to $D^{(1)}$. If $U\to Y$ is $\dagger$-\'etale, the cartesian diagram of \Cref{LemmaFormalConePullback} implies that the anima of  $U$-points of  $Y_{s}^{\ob{cone}/(D,W)}$ is the same as the anima of $U$-points  of $U_{s}^{\ob{cone}/(D,W)}$, and so it is equal to the anima of square-zero extensions of $U$ over $\eta$ to $D^{(1)}$. 
\end{proof}

\begin{remark}\label{RemarkConeTwoExamples}
Let $D\to W$ be a Cartier divisor and $X\to D$ a map admitting a cotangent complex. One can take two possible cone constructions, namely, $X^{\cone/D}$ and $X^{\cone/(D,W)}$. Generically, the former is $D$ while the last is $W$. The special fiber $X^{\cone/D}_s= \bb{T}_{X/D}[1](-1)$ is a shifted and twisted tangent space, but the special fiber $X^{\cone/(D,W)}_s$ is only a principally homogeneous $\bb{T}_X[1](-1)$-space over $X\times B\bb{G}_m$. The key difference between the two constructions are the difference between the maps $W\to D^{\cone}$ and $D\times \bb{A}^1/\bb{G}_m \to D^\cone$; the former is associated to the divisor $D\to W$, while the latter (after pullying back along $*\to B\bb{G}_m$) corresponds to the ``trivial divisor'' $D \xrightarrow{0} D\times \bb{A}^1$. In particular, $X^{\cone/D}_s$ corepresents deformations of $X$ to $D\times \GSpec (\bb{Q}_p[t]/t^2)$ (where one always has the trivial deformation), while $X^{\cone/(D,W)}_s$ corepresents deformations of $X$ to the square-zero extension $D^{(1)}$  associated to $D\to W$ (which might be empty). 
\end{remark}

\subsection{The cone construction: sheaf theory}\label{sss:ConeFlatConnection}
The goal of this subsection is to establish that if $X$ is a smooth  rigid variety, sheaves on $X^\cone$ are naturally related with flat connections, as well as a variant relative to a Cartier divisor.
For that, we shall follow the discussion of \cite[Section 5.2.5]{RCBBCourse2026}. Let $\bb{A}^{1,\ob{alg}}$ and $\bb{G}_m^{\ob{alg}}$ denote the algebraic affine line and multiplicative group seen as schemes over $\bb{Q}_p$.    Following \cite[Section 2.2]{bhatt2022prismatic}, the classical derived category $\ob{QCoh}(\bb{A}^{1,\ob{alg}}/\bb{G}_m^{\ob{alg}})$ is equivalent to the category of filtered $\bb{Q}_p$-modules with decreasing filtration.
Taking formal completions at $0$, quasi-coherent sheaves on $\ob{QCoh}(\widehat{\bb{A}}^1/\bb{G}_m^{\ob{alg}})$ are equivalent to complete filtered $\bb{Q}_p$-modules.   We want to have analogous descriptions for the base change to general analytic rings (or even analytic stacks!). This is justified in the next remark.

\begin{remark}\label{RemarkKunneth}
We want to see algebraic stacks as analytic stacks, for that we follow \cite[Section 6.4]{solidnotes}. Let $\cat{AlgStk}$ denote the category of functors on affine schemes satisfying descent along $!$-equivalences for the $\ob{QCoh}$-six-functor formalism.  There is a pullback map of topoi 
\[
(-)^{\ob{cond}}\colon \cat{AlgStk}\to \cat{AnStk}
\]
sending $\Spec A\mapsto \AnSpec( A^{\ob{cond}})$ where $A^{\ob{cond}}$ is the trivial analytic ring structure on $A$. By \cite[Remark 6.4.1]{solidnotes}, for $X\in \cat{AlgStk}$,  we have an equivalence of sheaves
\[
\ob{D}(X^{\ob{cond}}) = \ob{Fun}_{\ob{QCoh}(\bb{Z})}(\ob{QCoh}^!(X), \ob{D}(\bb{Z}^{\ob{cond}})),
\]
where $\ob{QCoh}^!(X)$ is the category of quasi-coherent  $!$-sheaves.  By design,  $\ob{QCoh}^!(X)$ is a pre-dual of $\ob{QCoh}(X)$, that is, we have 
\[
\ob{QCoh}(X)=\ob{Fun}_{\ob{QCoh}(\bb{Z})}(\ob{QCoh}^!(X),\ob{QCoh}(\bb{Z})).
\]

Now, if $X$ is $!$-able over a base ring $A$, $\ob{QCoh}(X)$ is selfdual and we have that  
\[
\begin{aligned}
\ob{D}(X^{\ob{cond}}) & = \ob{Fun}_{\ob{QCoh}(\bb{Z})}(\ob{QCoh}^!(X), \ob{D}(\bb{Z}^{\ob{cond}})) \\ 
& = \ob{QCoh}(X)\otimes_{\ob{QCoh}(\bb{Z})} \ob{D}(\bb{Z}^{\ob{cond}}). 
\end{aligned}
\]
The previous produces a natural morphism of $6$-functor formalisms 
\[
\ob{QCoh}(-)\to \ob{D}((-)^{\ob{cond}})
\]
on the geometric set up $(\cat{AlgStk},E)$ where  $E$ is the class of $!$-able maps obtained by applying \cite[Theorem 3.4.11]{heyer20246functorformalismssmoothrepresentations}
 to the class of morphisms  of affine schemes.  Thus, if $X$ is a $!$-able algebraic stack, the map $X^{\ob{cond}}\to \AnSpec (\bb{Z}^{\ob{cond}})$ is  $!$-able as analytic stack and therefore $\ob{D}(X^{\ob{cond}})$ is selfdual over $\ob{D}(\bb{Z}^{\ob{cond}})$. We apply the previous observation to the following cases:

\begin{itemize}

\item  $X=B\bb{G}_m^{\ob{alg}}$, where for $Y$ an analytic stack $\ob{D}(B\bb{G}_m^{\ob{alg}}\times Y )$  is the category of graded modules on $Y$. 

\item  $X=\bb{A}^{1,\ob{alg}}/\bb{G}_m^{\ob{alg}}$, where for $Y$ an analytic stack $\ob{D}(\bb{A}^{1,\ob{alg}}/{\bb{G}_m^{\ob{alg}}}\times Y )$  is the category of filtered modules on $Y$. 

\item  $X=\widehat{\bb{A}}^{1}/\bb{G}_m^{\ob{alg}}$, where for $Y$ an analytic stack $\ob{D}(\widehat{\bb{A}}^{1}/\bb{G}_m^{\ob{alg}}\times Y )$  is the category of complete filtered modules on $Y$.

\end{itemize}
\end{remark}

From now on we work with analytic rings over $\bb{Q}_{p,\sol}$.
The map  $\AnSpec( \bb{Q}_{p,\sol})\to \widehat{\bb{A}}^1$ is an epimorphism for the $!$-topology, thus we have an equivalence 
\[
\widehat{\bb{A}}^1= B \bb{G}_a[-1],
\]
where $\bb{G}_a[-1]= \AnSpec( \bb{Q}_{p,\sol})\times_{\bb{G}_a} \AnSpec (\bb{Q}_{p,\sol}) = \AnSpec ( \bb{Q}_{p,\sol}[\epsilon_1])$ for $\epsilon_1$ a free generator in degree $1$ (so $\epsilon_1^2=0$).
Therefore, we have an isomorphism of stacks
\[
\widehat{\bb{A}}^1/\bb{G}^{\ob{alg}}_m= B \bb{G}_a[-1]/  \bb{G}^{\ob{alg}}_m,
\]
where $\epsilon_1$ has degree $1$ and weight $-1$ with respect to the $\bb{G}_m$-action.

\begin{remark}\label{RemarkDiscussionJustification}
 Sheaves on $B \bb{G}_a[-1]/  \bb{G}^{\ob{alg}}_m$ are precisely chain complexes in $\ob{D}(\bb{Q}_p)$ with a differential shifted by $1$, that is, a sheaf on  $B \bb{G}_a[-1]/  \bb{G}^{\ob{alg}}_m$ is a graded  solid $\bb{Q}_p$-module $M=\bigoplus_n M_n$ together with maps $d_n\colon M_n\to M_{n-1}[1]$ satisfying $d_{n+1}\circ d_n=0$ together with higher coherences. Following \cite[Example 2.2.4 (5)]{bhatt2022prismatic}, the category $\ob{D}(B \bb{G}_a[-1]/  \bb{G}^{\ob{alg}}_m)$ has a natural $t$-structure (the Beilinson $t$-structure) making an object $M=\bigoplus_n M_n$ connective (resp.\ coconnective) if $\bigoplus_{n} M_n[n]$ is connective (resp.\ coconnective). The heart of the Beilinson $t$-structure is equivalent to the symmetric monoidal category of chain complexes $\ob{Ch}(\bb{Q}_{p,\sol})$ of solid $\bb{Q}_p$-vector spaces.  

 All the previous facts about sheaves on $B \bb{G}_a[-1]/  \bb{G}^{\ob{alg}}_m$ follow from a  shearing  as in \cite[Section 6]{benzvi2024relativelanglandsduality}, after taking a base change along $\ob{QCoh}(\bb{Q})\to \ob{D}(\bb{Q}_{p})$ in $\cat{Pr}_{\ob{QCoh}(\bb{Q})}$. One has obvious linear variants of the previous discussion over any symmetric presentable $\ob{QCoh}(\bb{Q})$-linear category by base change. 
 \end{remark}

The next proposition describes sheaves on the formal cone  for formally smooth maps of Berkovich spaces. For the formal algebraic filtered de Rham stack, a similar result has been independently observed by Jiang, \cite[Lemma 2.7.18]{jiang2026differential}.

\begin{proposition}\label{PropositionConeSheaves}
  Let $f\colon Y\to X$ be a morphism of Berkovich spaces such that $\Omega^1_{Y/X}=\bb{L}_{Y/X}$ is a vector bundle, i.e., $f$ is formally smooth in the sense of \Cref{sec:cotang-compl-morph-definition-having-a-cotangent-complex}.
  Consider the natural map of analytic stacks over $\bb{Q}_{p,\sol}$ 
\[
h\colon Y^{\widehat{\ob{cone}}/X} \to |Y|_{\Betti}\times \widehat{\bb{A}}^1/\bb{G}_m^{\ob{alg}}.
\]
 The following hold:
\begin{enumerate}

\item  The map $h$ is prim. In particular, the natural morphism 
\begin{equation}\label{eqFullyFaithfulMapConeFormal}
\ob{Mod}_{h_*1}(\ob{D}(|Y|_{\Betti}\times \widehat{\bb{A}}^1/\bb{G}_m^{\ob{alg}}))\hookrightarrow \ob{D}(Y^{\widehat{\ob{cone}}/X}) 
\end{equation}
is fully faithful. 

\item  The underlying graded solid sheaf   $h_*1$ on $|Y|_{\Betti}$ is the exterior algebra $\bigoplus_i \bigwedge^{i} \Omega_{Y/X}^i[-i](i)$, where $\Omega_{Y/X}^i[-i](i)$ is shifted by $-i$ and has weight $i$ for the $\bb{G}_m$-action.  Hence, following  \Cref{RemarkDiscussionJustification}, $h_*1$ is a promotion of $\bigoplus \bigwedge^{i} \Omega_{Y/X}^i[-i](i)$ to a commutative differential graded  algebra in $\ob{D}(|Y|_{\Betti})$. 

\item Suppose that the map $f$ factors through a formally \'etale map $Y\to \bb{A}^n_{X}$.
  Then, after a shearing as in \Cref{RemarkDiscussionJustification},  there is a natural equivalence  of commutative differential graded  algebras 
\[
h_*1 \cong (\Omega^{\bullet}_{\bb{A}^n/\bb{Q}_p},d)\otimes_{\s{O}_{\bb{A}^n_{\bb{Q}_p}}} \s{O}_{Y},
\]
where $(\Omega^{\bullet}_{\bb{A}^n/\bb{Q}_p},d)$ is the de Rham complex of $\bb{A}^n_{\bb{Q}_p}$. 

\end{enumerate}
\end{proposition}
\begin{proof}
\begin{enumerate}

\item Since $\AnSpec \bb{Q}_{p,\sol}\to \widehat{\bb{A}}^1$ is an epimorphism of analytic stacks, it suffices to show that $h$ is prim on special fibers. Then, it suffices to show that 
\[
Y/ \bb{T}_{Y/X}\rtimes \bb{G}_m\to |Y|_{\Betti}/ \bb{G}_m^{\ob{alg}}
\] 
is prim. This map factors through 
\[
Y/ (\bb{T}_{Y/X}\rtimes \bb{G}_m)\to Y/ \bb{G}_m\to  |Y|_{\Betti}/\bb{G}_m\to |Y|_{\Betti}/ \bb{G}_m^{\ob{alg}}. 
\]
Each of the following maps is prim by the following reasonings:

\begin{enumerate}

\item For an analytic vector bundle $V$ on a Gelfand stack $X$, the map $h\colon X/V\to X$ is prim.
  That is, by Cartier duality \cite[Theorem 3.2.27 (3)]{RCCartierGerbe} it suffices to show that $h^{\vee}\colon X\to V^{\vee,\dagger}$ is suave, where $V^{\vee,\dagger}\subset V^{\vee}$ is the overconvergent neighbourhood at $1$.
  Indeed, being prim means that the pullback  map $h^*\colon [X]\to [X/V]$ in the kernel category over $X$ admits a linear right adjoint, and Cartier duality identifies the previous map with $h^{\vee}_!\colon [X]\to [V^{\vee,\dagger}]$.
  This map admits a right adjoint in the kernel category of $X$ provided that $h^{\vee}$ is suave.
  Now, to prove that $h^{\vee}\colon X\to V^{\vee,\dagger}$ is suave, it suffices to consider the universal case over $X=B\GL_n$ and by descent for $X=\GSpec(\Q_p)$. In that case, $V\cong \bb{G}_a^n$ and suaveness follows form the fact that the zero section $0\colon \GSpec(\Q_p)\to \bb{G}_a^{\dagger,n}$ is a regular Zariski closed map and \cite[Lemma 6.3.5]{solidnotes}.

\item The map $Y/\bb{G}_m\to  |Y|_{\Betti}/\bb{G}_m$ is prim by \Cref{sec:deriv-berk-spac-berkovich-spaces-as-relative-spectrum}.

\item For the map of analytic stacks  $|Y|_{\Betti}/\bb{G}_m\to |Y|_{\Betti}/\bb{G}_m^{\ob{alg}}$, it suffices to show that $B\bb{G}_m\to B\bb{G}_m^{\ob{alg}}$ is prim, and by descent that $\bb{G}_m^{\ob{alg}}/\bb{G}_m\to \AnSpec \bb{Q}_{p,\sol}$ is prim. Since  $\bb{G}_m^{\ob{alg}}$ is represented by the ring $\bb{Q}_{p,\sol}[T^{\pm 1}]$, the map  $\bb{G}_m^{\ob{alg}}/\bb{G}_m\to B\bb{G}_m$ is prim, and it suffices to show that $B\bb{G}_m\to \AnSpec (\bb{Q}_p)$ is prim (we recall that $\bb{G}_m$ is the rigid analytic multiplicative group over $\bb{Q}_p$).
  Let $g\colon \GSpec(\Q_p)\to B\bb{G}_m$.
  The map $g$ is suave being represented in smooth partially proper rigid spaces, and therefore $g_{\sharp}1\in \ob{D}(B\bb{G}_m)$ is a prim object over $\GSpec(\Q_p)$. It suffices to show that $1$ is a direct summand of     $g_{\sharp}1$. By \cite[Proposition 4.4.8]{jacinto2023solid} the $*$-pushforward along $\GSpec(\Q_p)\to B\bb{G}_m$ produces an equivalence of categories
\[
\ob{D}(B\bb{G}_m)\xrightarrow{\sim} \ob{LMod}_A(\ob{D}(\bb{Q}_{p})),
\]
where $A$ is the dual nuclear Fr\'echet algebra of sequences $(a_n)_{n\in \bb{Z}}$ in $\bb{Q}_p$ such that there is some $r>0$ with $|a_{n}|r^{|n|}\to 0$ as $|n|\to \infty$, and with component-wise  multiplication. Under this equivalence, $g_{\sharp} 1 $ is precisely $A$ seen as a left $A$-module, and the unit of $B\bb{G}_m$ is $\Q_p$ endowed with the multiplication by $A$ given by the projection to the $0$-th component. It is clear that $\Q_p$ is a direct summand of $A$ as an $A$-module, this yields primness of $B\bb{G}_m\to B\bb{G}_m^{\ob{alg}}$ as wanted. 
\end{enumerate}

The last claim about the fully faithfulness of the map  \eqref{eqFullyFaithfulMapConeFormal} follows formally from the projection formula for $h_*$ due to the primness of $h$.

\item The underlying graded module is the $*$-pushforward of the unit along the map 
\[
Y/(\bb{T}_{Y/X}\rtimes \bb{G}_m)\to Y/\bb{G}_m\to |Y|_{\Betti}/\bb{G}_m. 
\]
It suffices to compute the $*$-pushforward of $Y/(\bb{T}_{Y/X}\rtimes \bb{G}_m)\to Y/\bb{G}_m$.
This is the classifying stack of a relative vector bundle over $Y/\bb{G}_m$ and can be computed via Cartier as in \cite[Theorem 3.2.37 (3)]{RCCartierGerbe}.  That is, the Cartier dual of $Y/(\bb{T}_{Y/X}\rtimes \bb{G}_m)$ over $Y/\bb{G}_m$ is $\bb{T}_{Y/X}^{*,\dagger}(-1)$, the overconvergent neighbourhood of the zero section of the cotangent bundle of $Y$ over $X$  with weight $-1$ for the $\bb{G}_m$-action. Cartier duality switches $*$ and $!$-functors, thus, if $g\colon Y/(\bb{T}_{Y/X}\rtimes \bb{G}_m)\to Y/\bb{G}_m$ is the structural map and $\iota\colon Y/\bb{G}_m\to \bb{T}_{Y/X}^{*,\dagger}(-1)$ the zero section, we have that 
\[
g_*1= g_*g^*1 = \iota^!\iota_! 1.
\]
This algebra agrees with the endomorphism algebra of $\s{O}_{Y/\bb{G}_m}$ as a $\ob{Sym}_{\s{O}_{Y/\bb{G}_m}} (\Omega^{1,\vee}_{Y/X}(1))$-module. This follows from the fact that the analytification of a vector bundle is an open immersion of the associated algebraic vector bundle. One can reduce this problem to the universal case over the stack of vector bundles which then follows by an explicit computation, see \cite[Proposition 5.2.40 (2)]{RCBBCourse2026}.

\item  By localizing in the analytic topology we can assume without loss of generality that $X=\GSpec (A)$ and  $Y=\GSpec (B)$  are affinoid  and that the $\dagger$-formally \'etale map $Y\to \bb{A}^n_{X}$ factors through $\bb{D}^{\leq 1,n}_X$.   By \Cref{LemmaFormalConePullback} (rather its variant for the formal cone) we have a Cartesian diagram 
\[
\begin{tikzcd}
 Y^{\widehat{\ob{cone}}/X} \ar[r] \ar[d]  & \bb{D}^{\leq 1,n,\widehat{\ob{cone}}/X}_{X}  \ar[d] \\ 
Y^{\dR} \ar[r]  & \bb{D}^{\leq 1,n,\dR}_{X}  
\end{tikzcd}
\]
Taking global sections and using (2), we deduce a natural map of commutative differential graded  algebras over $A$
\[
(\Omega_{\bb{D}^{\leq 1,n}_A}^{\bullet}, \partial)\to (\Omega^{\bullet}_{Y/X}, d)
\]
such that  $\Omega_{\bb{D}^{\leq 1,n}_A}^{\bullet}\otimes_{\s{O}_{\bb{D}^{\leq 1,n}_A}} B\to \Omega^{\bullet}_{Y/X}$ is an equivalence. Therefore, it suffices to show that the commutative differential graded  algebra defined by  $\bb{D}^{\leq 1,n}_{\bb{Q}_p}$ is the de Rham complex of the closed polydisc, i.e., that $\partial$ is the de Rham differental.
By K\"unneth formula, it suffices to deal with a single copy $\bb{D}^{\leq 1}_{\bb{Q}_p}$.
This special case is done by an explicit computation as that of \cite[Theorem 5.2.41 (2)]{RCBBCourse2026}. \qedhere

\end{enumerate}
\end{proof}

  We note that we do not define a de Rham complex for a general formally smooth morphism $f\colon Y\to X$ of Berkovich spaces, that could be equivalent to $h_\ast(1)$ in \Cref{PropositionConeSheaves} (3).
  In fact, any reasonable definition should agree with $h_\ast(1)$, and so one could take $h_\ast(1)$ as a definition of the de Rham complex. Moreover, using the topos theoretic interpretation of the cone \Cref{RemarkToposTeoreticDescriptionCone}, one can show by abstract nonsense that this construction of the de Rham complex produces the initial cdga over $|Y|_{\Betti}$, linear over  $\s{O}(X)$, and  admitting a map of graded algebras from $\s{O}_Y$ in degree $0$. 
 We will not show this fact in this paper, instead, we will content ourselves by  proving that  $h_\ast(1)$ recovers the classical de Rham complex for smooth rigid spaces.

\begin{corollary}\label{CorollaryConeSmoothRigidSpace}
  Let  $K$ be a Banach field extension of $\bb{Q}_p$ and $X$ a smooth rigid space over $K$, which is a Berkovich space (e.g., partially proper or affinoid, cf.\ \Cref{sec:deriv-berk-spac-rigid-spaces-as-berkovich-spaces}). 
  Let $h\colon X^{\widehat{\ob{cone}}/K}\to |X|_{\Betti}\times \widehat{\bb{A}}^1/\bb{G}_m^{\ob{alg}}$.
  Then $h_*1$ sits in the heart of Beilinson $t$-structure and, after a shearing as in \Cref{RemarkDiscussionJustification},  it is naturally isomorphic  to the de Rham complex $(\Omega^{\bullet}_{X/K},d)$ as a  commutative differential graded  algebra. 
\end{corollary}
\begin{proof}
  By \Cref{PropositionConeSheaves} (2), and after applying a shearing as in \Cref{RemarkDiscussionJustification},  the algebra $h_*1$ is a suitable commutative differential graded  algebra structure $(\Omega^{\bullet}_{X/K},\partial)$ on the Hodge complex  of $X$.
  By the universal property of the de Rham complex, there is a unique morphism of  commutative differential graded  algebras $(\Omega^{\bullet}_{X/K},d)\to (\Omega^{\bullet}_{X/K},\partial)$ that induces the identity on $\s{O}_X$. We claim that this map is an isomorphism.
  This can be proven locally on \'etale charts to affine spaces, in which case it follows from \Cref{PropositionConeSheaves} (3).
\end{proof}

\begin{example}\label{ExampleFilteredVectorBundles}
  Let $K$ be a non-archimedean field over $\bb{Q}_p$ and $X$ a smooth rigid space over $K$, which is a Berkovich space.
  Under the equivalence between differential graded modules and complete filtered modules of \Cref{RemarkKunneth}, the de Rham complex $(\Omega^{\bullet}_X,d)$ obtained as the relative cohomology of the formal cone of \Cref{CorollaryConeSmoothRigidSpace} gives rise to the filtered  commutative algebra $(\Fil^\bullet (\ob{dR}(X/K)))$ given by the de Rham complex
\[
\ob{dR}(X/K)=   \s{O}_X\xrightarrow{d} \Omega^1_X\xrightarrow{d} \Omega^2_X\xrightarrow{d} \cdots \xrightarrow{d} \Omega^{\dim X}_X
\]
endowed with the Hodge filtration 
\[
\Fil^i(\ob{dR}(X/K))=    \cdots  0 \to  \Omega^i_X\xrightarrow{d} \Omega^{i+1}_X\xrightarrow{d} \cdots \xrightarrow{d} \Omega^{\dim X}_X. 
\]

Let  $(\Fil^{\bullet} V, V,\nabla)$ be a filtered flat connection over $X$ satisfying Griffiths transversality, i.e.,  $\nabla(\Fil^i(V))\subseteq \Fil^{i-1}(V) \otimes_{\s{O}_X} \Omega_X^1$ for $i\in \Z$.
We can construct a complete filtered module over $\Fil^\bullet(\ob{dR}(X/K))$ as
 \[
 \ob{dR}(V,\nabla)=( V\xrightarrow{\nabla} V\otimes_{\s{O}_X} \Omega^1_X\xrightarrow{\nabla} V\otimes_{\s{O}_X} \Omega^2_X \xrightarrow{\nabla} \cdots  \xrightarrow{\nabla} V\otimes_{\s{O_X}} \Omega^{\dim X}_X)
 \]
 with filtration given by 
 \[
 \ob{Fil}^i(\ob{dR}(V,\nabla)) = \Fil^i V \xrightarrow{\nabla} \Fil^{i-1 }V \otimes_{\s{O}_X} \Omega^1_X \xrightarrow{\nabla} \Fil^{i-2}V\otimes_{\s{O}_X} \Omega^2_X \xrightarrow{\nabla} \cdots  \xrightarrow{\nabla} \Fil^{i-\dim_X} V \otimes_{\s{O}_X} \Omega^{\dim X}_X
 \]
 (cf.\ \cite[Theorem 3.3]{faltings_chai_degenerations_of_abelian_varieties}). 
 
Applying the global sections functor $\ob{Mod}_{\ob{dR}(X/K)}(\ob{D}(|X|_{\Betti}\times \widehat{\bb{A}}^1/\bb{G}_m^{\ob{alg}}))\to \ob{D}(|X|_{\Betti}/ \bb{G}_m^{\ob{alg}})$, the filtered $\ob{dR}(X/K)$-module $ \ob{dR}(V,\nabla)$ is sent to the graded sheaf over $|X|_{\Betti}$ given by the  filtered de Rham cohomology 
 \[
 \bigoplus_{i}  \ob{Fil}^i(\ob{dR}(V,\nabla)).
 \]
 
 On the other hand, using the fully faithful map $\ob{Mod}_{{\ob{dR}(X/K)}}(\ob{D}(|X|_\Betti\times \widehat{\bb{A}}^1/\bb{G}_m^{\ob{alg}}))\hookrightarrow \ob{D}(X^{\widehat{\cone}/K})$ of \Cref{PropositionConeSheaves}, we obtain a fully faithful embedding of the category of filtered flat connections on $X$ into modules over the formal cone $X^{\widehat{\ob{cone}}/K}$.
\end{example}

\begin{remark}
  \label{sec:cone-constr-sheaf-argument-via-filtered-de-rham-stack}
  An alternative way to describe vector bundles on the formal cone is via the filtered de Rham stack, cf. \cite[Lemma 2.7.13]{jiang2026differential} in the algebraic and logarithmic setting or the upcoming work of the third author with Joaqu\'in Rodrigues Jacinto on quasi-coherent sheaves on de Rham stacks.
  From the cartesian square \Cref{sec:cone-construction-1-proposition-filtered-de-rham-stack}, base changed to $\widehat{\mathbb{A}}^1/\mathbb{G}_m$,
  one deduces that the pullback along
  \[
    \alpha \colon Y^{\dR,+/X}\to Y^{\widehat{\cone}/X}
  \]
  is fully faithful (it is a torsor under the de Rham stack of a geometric vector bundle, see the proof of \Cref{fg-suave} for a similar argument).
  Now, quasi-coherent sheaves on $Y^{\dR,+/X}$ identify with (an analytic version of) filtered $\mathcal{D}$-modules, and on perfect complexes $\alpha^\ast$ is an equivalence.
  Namely, for vector bundles the existence of a filtration by subbundles on which the connection satisfies Griffiths transversality automatically implies the nilpotence of the Higgs field on the associated graded.
\end{remark}

\begin{remark}\label{RemarkAlgebraicConeGeneralAnalyticStacks}
\begin{enumerate}

\item  The cone construction can be naturally extended to all analytic stacks as a functor sending $X$ to $X^{\ob{cone-alg}}\to \bb{A}^{1,\ob{alg}}/\bb{G}_m^{\ob{alg}}$ by taking the sheafification of the functor sending an affinoid analytic stack $\AnSpec (A)\to \bb{A}^{1,\ob{alg}}/\bb{G}_m^{\ob{alg}}$ endowed with a Cartier divisor  $I\to A$ to $X(A/I)$.  One has the obvious relative variants as in \Cref{DefinitionRelativeConesCartierDivisors}, as well as the formal cone $X^{\widehat{\ob{cone-alg}}}$, and the special fiber of the cone $X^{\ob{cone-alg}}_s$.

If $f\colon Y\to X$ is a morphism of analytic stacks over $\bb{Q}$ that admits a cotangent complex and is infinitesimally cohesive (for a suitable basis of analytic rings), one has that the special cone
\[
Y^{\ob{cone-alg}/X}_s = \bb{T}_{Y/X}(-1)[1]
\]
is the shifted and twisted tangent bundle of $Y\to X$. If $f$ arises from a morphism of Gelfand stacks over $\bb{Q}_p$, then we have a natural commutative square  of analytic stacks
\[
\begin{tikzcd}
Y^{\cone/X} \ar[r] \ar[d] & Y^{\ob{cone-alg}/X} \ar[d]  \\ 
\bb{A}^1/\bb{G}_m \ar[r] & \bb{A}^{1,\ob{alg}}/\bb{G}_m^{\ob{alg}}
\end{tikzcd}
\]
making $Y^{\cone/X}$ the \textit{Gelfandification of $Y^{\ob{cone-alg}/X} $} as explained in \cite[Section 4]{drff_stacks}.

Assume now that $f\colon Y\to X$ is a morphism of Berkovich spaces which is formally smooth and $\bb{T}_{Y/X}$ is a vector bundle. Then the exact same argument as in \cite[Theorem 5.2.41]{RCBBCourse2026} shows that the natural map $f\colon Y^{\widehat{\ob{cone-alg}}/X}\to |Y|_{\Betti}\times \widehat{\bb{A}}^1/\bb{G}_m^{\ob{alg}}$ is $0$-affine, given by the relative  spectrum of the de Rham  complex $(\Omega^{\bullet}_{Y/X}, d)$ considered as a commutative differential graded  algebra over $|Y|_{\Betti}$.
Hence, the proof of \Cref{PropositionConeSheaves} shows that the natural map $g\colon Y^{\widehat{\cone}/X}\to Y^{\widehat{\ob{cone-alg}}/X}$ is prim and $1\xrightarrow{\sim} g_*1$ is an equivalence.
As a further consequence of $0$-affineness, we see that
\[
  \ob{D}(Y^{\widehat{\ob{cone-alg}}/X})\cong \mathrm{Mod}_{(\Omega^\bullet_{Y/X},d)}(\ob{D}(|Y|_\Betti\times \widehat{\bb{A}}^1/\bb{G}_m^{\ob{alg}})).
  \]

\item Suppose that  $f\colon Z\to X$ is a $!$-Zariski closed map of Gelfand stacks, that is $f$ is To\"en-affine and $\s{O}_X\to \s{O}_Z$ is a $!$-locally a surjection of connective rings (this is different from \Cref{DefinitionZariskiCLosed}, cf.\ \Cref{sec:push-berk-spac-non-shriek-descent-for-zariski-closed-immersion}).
  Then the formal cone $Z^{\widehat{\ob{cone-alg}}/X}\to X\times \bb{A}^{1,\ob{alg}}/\bb{G}_m^{\ob{alg}}$ is given by the relative spectrum of the filtered ring $\ob{Fil}^{\bullet} \s{O}_X$ given by the adic filtration of $\s{O}_X\to \s{O}_{Z}$ (which is defined by descent from the case of affinoid stacks), see \cite[Remark 7.7.17]{solidnotes}.
\end{enumerate}
\end{remark}

\begin{remark}
  \label{sec:cone-constr-sheaf-natural-surjection-to-cone}
  In the notations of \Cref{PropositionConeSheaves} and \Cref{RemarkAlgebraicConeGeneralAnalyticStacks}, there is a natural surjection
  \[
    g\colon Y\times \widehat{\bb{A}^1}/\bb{G}_m^{\ob{alg}}\to Y^{\widehat{\ob{cone-alg}}/X}
  \]
  between $0$-affine stacks over $|Y|_\Betti\times \widehat{\bb{A}^1/\bb{G}_m^{\ob{alg}}}$.
  Under the equivalence,
  \[
    \ob{D}(Y^{\widehat{\ob{cone-alg}}/X})\cong \mathrm{Mod}_{(\Omega^\bullet_{Y/X},d)}(\ob{D}(|Y|_\Betti\times \widehat{\bb{A}}^1/\bb{G}_m^{\ob{alg}})),\]
  the morphism $g$ is induced by the morphism
  \[
    (\Omega^\bullet_{Y/X},d)\to \mathscr{O}_Y=\mathrm{gr}^0(\Omega^\bullet_{Y/X},d)
  \]
  of filtered algebras (whenever the de Rham complex is defined).
\end{remark}

In \Cref{PropositionConeSheaves}, we have discussed the relative formal cone for a morphism $f\colon Y\to X$ of Berkovich spaces. We will now discuss the formal cone relative to a Cartier divisor, cf.\ \Cref{DefinitionRelativeConesCartierDivisors}.
Fix a Gelfand stack $W$ with a Cartier divisor $D\subseteq W$, and assume that $W$ is complete along $D$, i.e., that $D\to W$ is an epimorphism.
Let $f\colon X\to D$ be a morphism of Gelfand stacks, and fix a lift $\widetilde{f}\colon \widetilde{X}\to W$ of $f$.
Furthermore, we have a natural equivalence 
\[
  X^{{\cone}/(D,W)}\cong \widetilde{X}^{{\cone}/W}\times_{W\times \bb{A}^1/\bb{G}_m} W
\]
where the morphism $W\to W\times \mathbb{A}^1/\Gm$ is the identity on the first factor and the map classifying $D$ on the second.
Indeed, we have a commutative diagram  
\[ \begin{tikzcd}
	& {X^{\cone/(D,W)}} & W & \\
	{\widetilde{X}^{\cone/W}} & {X^\cone} & {D^\cone} & {W\times \mathbb{A}^1/\mathbb{G}_m} \\
	& {\widetilde{X}^\cone} & {W^\cone}
	\arrow[from=1-2, to=1-3]
	\arrow[from=1-2, to=2-1]
	\arrow[from=1-2, to=2-2]
	\arrow[from=1-3, to=2-3]
	\arrow[from=1-3, to=2-4]
	\arrow[bend left=15, from=2-1, to=2-4]
	\arrow[from=2-1, to=3-2]
	\arrow[from=2-2, to=2-3]
	\arrow[from=2-2, to=3-2]
	\arrow[from=2-3, to=3-3]
	\arrow[from=2-4, to=3-3]
	\arrow[from=3-2, to=3-3]
\end{tikzcd}\]
with all squares cartesian.

As we assumed that $W$ is complete along $D$, the morphism $W\to \mathbb{A}^1/\mathbb{G}_m$ factors over $\widehat{\mathbb{A}}^1/\mathbb{G}_m$, and we deduce
\[
    X^{\cone/(D,W)}= X^{\widehat{\cone}/(D,W)}\cong \widetilde{X}^{\widehat{\cone}/W}\times_{W\times \widehat{\mathbb{A}}^1/\Gm}W.
\]

From \Cref{PropositionConeSheaves} we can derive the following consequence.

\begin{proposition}
  \label{sec:cone-constr-sheaf-sheaves-relative-cone-divisor-case}
  Assume furthermore that $f\colon X\to D$ is a formally smooth morphism of Berkovich spaces.
  Consider the natural map of Gelfand stacks over $\bb{Q}_{p}$ 
\[
h\colon X^{\widehat{\ob{cone}}/(D,W)} \to |X|_{\Betti}\times W.
\]
 The following hold:
\begin{enumerate}

\item  The map $h$ is prim.
  In particular, the natural morphism 
\[
\ob{Mod}_{h_*1}(\ob{D}(|X|_{\Betti}\times W))\hookrightarrow \ob{D}(X^{\widehat{\ob{cone}}/(D,W)}) 
\]
is fully faithful. 
\item Suppose that $\widetilde{f}$ factors through a formally \'etale map $\widetilde{X}\to \bb{A}^n_{W}$ (and that $D$ is static).
  Then, after a shearing as in \Cref{RemarkDiscussionJustification},  there is a natural equivalence  of commutative differential graded  algebras 
\[
h_*1 \cong (\Omega^{\bullet}(\bullet D)_{\bb{A}^n_W/W},\widetilde{d})\otimes_{\s{O}_{\bb{A}^n_{W}}} \s{O}_{X},
\]
where $(\Omega^{\bullet}(\bullet D)_{\bb{A}^n_W/W},\widetilde{d})$ is the $D$-de Rham complex of $\bb{A}^n_{\bb{Q}_p}$, i.e.,
\[
    \mathscr{O}_{\bb{A}^n_W}\overset{\widetilde{d}}{\to} \Omega^1_{\bb{A}^n_W/W}(D)\overset{\widetilde{d}}{\to} \Omega^2_{\bb{A}^n_W/W}(2\cdot D)\overset{\widetilde{d}}{\to} \ldots
\]
with $\widetilde{d}$ induced by the composition $\mathscr{O}_{\bb{A}^n_W}\overset{d}{\to} \Omega^1_{\bb{A}^n_W/W}\to \Omega^1_{\bb{A}^n_W/W}(D)$.
\end{enumerate}  
\end{proposition}
In (2) the assumption that $D$ is static is only added to be able to speak about the de Rham complex in a classical way.
\begin{proof}
  We note that $\widetilde{X}\to W$ is formally smooth, as this holds after base change along the cover $D\to W$.
  Using the isomorphism
  \[
      X^{\widehat{\cone}/(D,W)}\cong \widetilde{X}^{\widehat{\cone}/W}\times_{W\times \widehat{\mathbb{A}}^1/\Gm}W,
    \]
    that was explained before, the statements follow by base change from \Cref{PropositionConeSheaves} (applied to the formally smooth morphism $\widetilde{X}\to W$) along the map $W\to W\times \widehat{\mathbb{A}}^1/\mathbb{G}_m^{\ob{alg}}$
    (with this argument we can actually deduce that $\ob{Mod}_{h_*1}(\ob{D}(|X|_{\Betti}\times W))$ is equivalent to quasi-coherent sheaves on the formal algebraic relative cone of $X$ over $(D,W)$).
    The appearance of the $D$-de Rham complex comes as follows:
    applying the Rees construction to the algebra $(\Omega^{\bullet}_{\bb{A}^n/\bb{Q}_p},d)$ with its Hodge filtration yields the $\mathbb{G}_m^{\ob{alg}}$-equivariant de Rham complex
    \[
      \mathscr{O}_{\mathbb{A}^1_{W[t]}}\overset{d}{\to}t^{-1}\Omega^1_{\mathbb{A}^1_{W[t]/}/W[t]}\overset{d}{\to} t^{-2}\Omega^2_{\mathbb{A}^1_{W[t]/}/W[t]}\overset{d}{\to}\ldots
    \]
    with $t$ of weight $-1$ for the $\mathbb{G}_m^{\ob{alg}}$-action.
    Pulling back along the morphism $W\to W\times \widehat{\mathbb{A}}^1/\mathbb{G}_m$ yields the $D$-de Rham complex.
\end{proof}

\begin{remark}
  \label{sec:cone-constr-sheaf-relation-to-t-connections}
  There exists a natural morphism
  \[
    \gamma\colon \widetilde{X}\to X^{\widehat{\cone-\ob{alg}}/(D,W)}
  \]
  of $0$-affine stacks over $|X|_\Betti\times W$, and in fact $\gamma$ is induced by the natural morphism $\mathcal{A}^\bullet:=(\Omega^\bullet_{\widetilde{X}/W}(\bullet D),\widetilde{d})\to \mathscr{O}_{\widetilde{X}}$ (whenever the $D$-de Rham complex is defined).
  Using $\gamma^\ast$, one checks that vector bundles on $X^{\widehat{\cone}/(D,W)}$ identify with $D$-connections whose reduction modulo $D$ is nilpotent.
  More precisely, assume that $X,D$ are static in the following and that the (classical) $D$-de Rham complex of $\widetilde{X}$ relative to $W$ is defined, e.g., $W=\mathrm{Spf}(B^+_\dR)$ and $\widetilde{X}$ is an inverse limit of smooth analytic adic spaces over $B^+_\dR/t^n$.
  Set
  \[
    \mathcal{D}_{\widetilde{X},D}:=\Hom_{\mathcal{A}^\bullet}(\mathscr{O}_{\widetilde{X}},\mathscr{O}_{\widetilde{X}}),
  \]
  the algebra of $D$-differential operators.
  In coordinates $x_1,\ldots, x_n$ on $\widetilde{X}$, and after fixing a generator $t$ of $\mathcal{O}_W(-D)$, it is isomorphic to the power series algebra $\mathscr{O}_{\widetilde{X}}[[t\frac{\partial}{\partial x_1},\ldots,t\frac{\partial}{\partial x_n}]]$.
  In particular, given $M\in \Mod_{h_\ast(1)}(\ob{D}(|X|_\Betti\times W))$ then $\Phi(M):=\mathscr{O}_{\widetilde{X}}\otimes_{\mathcal{A}^\bullet} M$ acquires a natural $\mathcal{D}_{\widetilde{X},D}$-module structure.
  By base change to $D$, one verifies that if $\Phi(M)$ is a vector bundle on $\widetilde{X}$ the $\mathcal{D}_{\widetilde{X},D}$-module structure defines a $D$-connection on $\Phi(M)$ whose pullback to $D$ is nilpotent.
  One verifies (again by reduction modulo $D$) that this yields an equivalence between vector bundles on $X^{\widehat{\cone}/(D,W)}$ with $D$-connections on $\widetilde{X}$ whose reduction is nilpotent.
\end{remark}

\begin{remark}
  \label{sec:cone-constr-sheaf-argument-via-a-relative-de-rham-stack}
  Similar to \Cref{sec:cone-constr-sheaf-argument-via-filtered-de-rham-stack}, one can argue via filtered de Rham stacks (again relying on upcoming work of the third author with Joaqu\'in Rodrigues Jacinto on quasi-coherent sheaves on de Rham stacks).
  Namely, vector bundles on $\widetilde{X}^{\dR,+}\times_{W\times \A^1/\Gm} W$ identify with $D$-connections on $\widetilde{X}$ and (if $W$ is complete along $D$) the pullback along
  \[
    \widetilde{X}^{\dR,+}\times_{W\times \A^1/\Gm} W\to X^{\widehat{\cone}/(D,W)}\cong \widetilde{X}^{\widehat{\cone}/W}\times_{W\times \A^1/\Gm}W
  \]
  is fully faithful with essential image exactly given by those $D$-connections whose reduction modulo $D$ is nilpotent.
\end{remark}

\newpage
\section{Analytic $\delta$-stacks}\label{Section:Generalized-FF-curves}

In the theory of prismatization for $p$-adic formal schemes the key algebraic objects are $\delta$-rings, that is, rings $A$ together with a Frobenius lift $\varphi$ (in the animated sense).
In the theory of analytic prismatization over $\Q_p$, the analogues of these objects will have a more geometric  flavour and will typically be given by $\dagger$-nilthickenings of Fargues-Fontaine discs.
In this section we define precisely these objects and a category of stacks on them, and prove some basic properties.

\subsection{A toy model: the filtered de Rham stack}\label{sec:filtered-de-rham-stack}
Our definition of the analytic version of $\delta$-rings can be motivated by Witt vectors considerations parallel to the theory of $p$-adic formal schemes, cf. \Cref{sec:thick-farg-font-relation-to-witt-vectors2-intro}. To provide another motivation for looking at thickenings of Fargues--Fontaine discs, we start by briefly revisiting the theory of filtered de Rham stacks. As will be made precise later in this text, the filtered de Rham stack shares similarities with the (to be defined) $p$-adic twistor stack (i.e., the perfect analytic prismatization modulo Frobenius). 

 The filtered de Rham stack geometrizes de Rham cohomology endowed with the Hodge filtration; it interpolates between the de Rham stack and the Hodge stack. In the algebraic context, its definition goes back to Simpson, see \cite[\S 2.3]{bhatt2022prismatic} for a discussion. In the $p$-adic analytic context relevant to us, cf.\ \cite[Definition 5.1.1.(1)]{RCAndeRham}. We recall the definition of the filtered de Rham stack.

\begin{definition}
Let $X$ be a Gelfand stack, its \textit{filtered de Rham stack} is the (sheafification of the) functor
\[
X^{\dR,+}\to \bb{A}^1/\bb{G}_m
\]
sending a Gelfand ring $R$ endowed with a Cartier divisor $I\to R$ to the anima 
\[
X^{\dR,+}(I\to R)= X(R/(I\otimes_R\ob{Nil}^{\dagger}(R))).
\]
Note that the projection map $R\to R/(I\otimes_R\ob{Nil}^{\dagger}(R))$ produces a natural morphism of stacks $X\times \bb{A}^1/\bb{G}_m\to X^{\dR,+}$. 

For a morphism $Y\to X$ we define the \textit{relative filtered de Rham stack} to be 
\[
Y^{\dR,+/X}:= Y^{\dR,+}\times_{X^{\dR,+}} (X\times \bb{A}^1/\bb{G}_m). 
\]

More generally, for $W\to \bb{A}^1/\bb{G}_m$ a Gelfand stack endowed with a Cartier divisor $D\to W$, and $X\to W$ a Gelfand stack, we define the \textit{filtered de Rham stack of $X$ over $(D,W)$} to be 
\[
X^{\dR,+/(D,W)}:=X^{\dR,+}\times_{W^{\dR,+}} W,
\]
where $W\to W^{\dR,+}$ arises from the diagonal map $W\to W\times \bb{A}^1/\bb{G}_m$ (with morphism to $\bb{A}^1/\mathbb{G}_m$ classified by $D$) and the natural projection to $W^{\dR,+}$. 
\end{definition}

 \begin{remark}\label{RemarkFilteredDeRhamRingStack}
Recall that $R/(I\otimes_R\ob{Nil}^{\dagger}(R))$ has a natural structure as an animated ring since it fits in the pullback diagram
\[
\begin{tikzcd}
R/I\otimes_R\ob{Nil}^{\dagger}(R) \ar[r] \ar[d] & R/I \ar[d] \\ 
R^{\dagger-\ob{red}} \ar[r]& R^{\dagger-\ob{red}} /  (R^{\dagger-\ob{red}}\otimes_{R} I).
\end{tikzcd}
\]
In particular, by taking the underlying pullback diagram of discrete animated rings, we obtain a pullback square of ring stacks
\[
\begin{tikzcd}
\bb{G}_a^{\dR,+}  \ar[r] \ar[d]& \bb{G}_a^{\cone} \ar[d] \\ 
\bb{G}_a^{\dR}\ar[r] &  (\bb{G}_a^{\cone})^{\dR}
\end{tikzcd}
\]
proving in particular that the presheaf defining $\bb{G}_a^{\dR,+}$  is already a sheaf for the $!$-topology when restricted to nilperfectoids. 
 \end{remark}

Thanks to \Cref{RemarkFilteredDeRhamRingStack} we can think of the filtered de Rham stack as a \textit{relative de Rham stack over the cone}, more precisely, we have the following proposition:

\begin{proposition}
  \label{sec:cone-construction-1-proposition-filtered-de-rham-stack}
  Let $Y\to X$ be a $\dagger$-cohesive morphism of Gelfand stacks. Then the natural commutative diagram
  \[
  \begin{tikzcd}
  Y^{\dR,+/X} \ar[r] \ar[d] & Y^{\cone/X} \ar[d] \\ 
   (Y\times \bb{A}^1/\bb{G}_m)^{\dR} \ar[r] &   (Y^{\cone/X})^{\dR}
  \end{tikzcd}
  \] 
  produced by the pullback diagram of Gelfand rings of \Cref{RemarkFilteredDeRhamRingStack} is cartesian. 
\end{proposition}
\begin{proof}
This follows directly from  the Cartesian diagram of rings stacks and $\dagger$-cohesivity as $R/I\to  R^{\dagger-\red}/ (R^{\dagger-\red}\otimes_{R} I)$ is a $\dagger$-thickening and $R^{\dagger-\red}\to R^{\dagger-\red}/ (R^{\dagger-\red}\otimes_{R} I)$ is a surjection. 
\end{proof}

\begin{remark}
  \label{sec:toy-model:-filtered-geometry-of-filtered-de-rham-stack}
  Let $X$ be a Gelfand stack, which admits a cotangent complex.
  Over $\mathbb{G}_m/\mathbb{G}_m$ the filtered de Rham stack $X^{\dR,+}$ is given by $X^\dR$, but over $B \mathbb{G}_m$ it is the Hodge stack, i.e., the classifying stack of the overconvergent shifted tangent bundle $\mathbb{T}_{X/\Q_p}^\dagger[1](-1):=\mathrm{fib}(\mathbb{T}_{X/\Q_p}[1](-1)\to \mathbb{T}_{X/\Q_p}[1](-1)^\dR\times_{X^\dR}X)$.
\end{remark}

We now explain a different description of the filtered de Rham stack, using thickened affine lines. If $A$ is a nilperfectoid $\Q_p$-algebra, let $\tilde{\A}_A^1$ be the pushout in Berkovich spaces (cf.\ \Cref{LemmaPushoutsBerkovichSpaces})

\[
\tilde{\A}_A^1 = \A^1_{A^u} \sqcup_{\GSpec(A^u)} \GSpec(A),
\]
where on the right  the morphism $\GSpec(A^u)\xto{0} \A_{A^u}^1$ is the zero-section. Informally, $\tilde{\A}^1_A$ sits in between the constant families of affine lines $\A_{A^u}^1$ and $\A_A^1$; contrary to the latter, which is a $\dagger$-thickening of $\A_{A^u}^1$ at every point, $\tilde{\A}_A^1$ is only thickened at the origin. 

\begin{definition}
\begin{enumerate}
\item Let $A$ be a nilperfectoid $\Q_p$-algebra. A \textit{degree $1$ divisor} on $\tilde{\A}^1_A$ is a generalized Cartier divisor on $\tilde{\A}_A^1$ whose pullback to $\A_{A^u}^1$ comes from a section $\GSpec(A^u) \to \A_{A^u}^1$ of the projection. 
\item Let $\Div_{\A^1}^1$ be the Gelfand stack sending a nilperfectoid ring $A$ to the anima of degree $1$ divisors on $\tilde{\A}_A^1$. 
\item Let $X$ be a Gelfand stack. Define $X^{\Div_{\A^1}^1}$ as the  Gelfand stack  over $\Div_{\A^1}^1$  (after sheafification)  sending a nilperfectoid ring $A$ with a degree $1$ divisor $D$ on $\tilde{\A}_A^1$ to the anima of maps $D\to X$.  Note that there is a natural map $\bb{A}^1_X:=\mathbb{A}^1\times_{\GSpec(\Q_p)}X\to X^{\Div_{\A^1}^1}$ sending a nilperfectoid ring $\GSpec(A)\to X$ together with a section $a\in A$ to the degree $1$ divisor $D=V(T-a)$ obtained by pullback along $\widetilde{\bb{A}}^1_A\to \bb{A}^1_A$, and to the morphism   $D=\GSpec(A)\to X$ defined by $A$.  

\item  More generally, for $Y\to X$ a morphism of Gelfand stacks define $Y^{\Div_{\A^1}^1,/X}:= Y^{\Div_{\A^1}^1}\times_{X^{\Div_{\A^1}^1}} \bb{A}^1_X$.  
\end{enumerate}
\end{definition}

 The pushout definition of $\tilde{\A}_A^1$ gives rise to a pullback diagram
  \begin{equation}\label{eqPullbackDiv1FilteredDeRham}
\begin{tikzcd}
\Div_{\A^1}^1 \ar[r] \ar[d] & \A^{1,\dR} \ar[d] \\ 
\A^1/\Gm \ar[r] & (\A^1/\Gm)^\dR,
\end{tikzcd}
\end{equation}
where the maps correspond to pulling back divisors along the natural maps $\GSpec(A) \to \tilde{\A}_A^1$, $\A_{A^u}^1 \to \tilde{\A}_A^1$, $\GSpec(A^u) \to \GSpec(A)$, $\GSpec(A^u) \xto{0} \A_{A^u}^1$. As a consequence, we deduce an isomorphism
\[
\Div_{\A^1}^1 \cong \A^1/\Gm^\dagger.
\]

\begin{remark}\label{RemarkVectorBundlesStacks}
Vector bundles on $ \A^1/\Gm^\dagger$  are the same as vector bundles on $\A^1$ equipped with a $T$-connection, i.e.\ (by dividing by $T$) connections on $\A^1$ with a logarithmic singularity at $0$.  Indeed, sheaves on $\A^1/\widehat{\bb{G}}_m$ are precisely all $T$-connections on $\A^1$.  The $*$-pullback along $f\colon \A^1/\widehat{\bb{G}}_m \to \A^1/\Gm^\dagger$  is fully faithful as the map is prim and $f_*1=1$ (this is a consequence of the contractibility of $\bb{G}_m^{\dagger}$ for de Rham cohomology).  The essential image of $f^*$ are those $T$-connections on $\A^1$ for which the action of $T\partial_T$ converges as an holomorphic function on $\bb{A}^1$ (this is a formal consequence of Cartier duality for analytic vector bundles, see \cite{RCCartierGerbe}). For vector bundles with $T$-connection, this convergence condition is immediate as,  locally on rigid affinoid spaces, they are Banach spaces and $T\partial_T$ is a continuous operator. 
\end{remark}

Let $A$ be a nilperfectoid ring and $D\to \widetilde{\bb{A}}^1_{A}$ a degree $1$ divisor. The pushout  description of $\widetilde{\bb{A}}^1_{A}=\bb{A}^1_{A^u}\bigsqcup_{\GSpec(A^u)} \GSpec(A)$ produces a pushout description of $D$
\[
D= \overline{D}\bigsqcup_{D_{A^u}} D_A.
\]
The description of $\Div^1_{\bb{A}^1}$ in \eqref{eqPullbackDiv1FilteredDeRham} shows that $D$ is  the spectrum of the pullback  $A/I \times_{A^u/(A^u\otimes_A I)} A^u$, where $I\to A$ is the pullback of $D$ to $A$.
Thus, $D$ is simply $\GSpec (A/ (I\otimes_A \ob{Nil}^{\dagger}(A)))$. From this, we deduce that for a $\dagger$-cohesive Gelfand stack $X$ there is a natural equivalence
\[
X^{\Div^1_{\A^1}}= X^{\dR,+}\times_{\bb{A}^1/\bb{G}_m} \A^1/\bb{G}_m^{\dagger}. 
\]
In other words, after base changing along $\A^1/\Gm^\dagger\to \A^1/\Gm$, the filtered de Rham stack of $X$ can be interpreted as the moduli of degree $1$ divisors on thickened affine lines with a map to $X$.
Furthermore, vector bundles on $X^{\Div_{\A^1}^1}$ are connections on $X\times \A^1$ with logarithmic singularities at $T=0$.

The construction of the $p$-adic twistor stack explained later in this text will be along similar lines. In fact, continuing the analogy mentioned in \Cref{rk-intro:analogy-real}, we note that over the complex numbers, (the complex version of) the stack $X^{\Div_{\A^1}^1}$ naturally appears when describing a base change of a variant of the complex twistor stack, cf.\ \cite[Chapter VIII]{real_local_langlands}. To pass from the filtered de Rham stack to the $p$-adic twistor stack, one roughly has to replace in the above the constant families of affine lines over perfectoid rings by perfectoid families of Fargues--Fontaine curves. We therefore start by recalling in the next subsection some basic properties of such families.

\subsection{Classical Fargues-Fontaine curves}
\label{sec:farg-font-curv}
Let $\cat{Perfd}^{\ob{aff}}$ and $\cat{Perf}^{\ob{aff}}$  be  the categories of (light, qfd, partially proper) affinoid perfectoid spaces over $\Z_p$ and $\F_p$ respectively. Given $X=\Marc(A)$ in $\cat{Perfd}^{\ob{aff}}$ we let $X^{\flat}$ denote its tilt over $\F_p$ and we let $\A_{\inf}(X)=\A_{\inf}(A)= W(A^{\flat,\circ})$ be the ring of Witt vectors of $A^{\flat,\circ}$.\footnote{We note that $A^\circ=A^{\leq 1}$ for each Banach $\Q_p$-algebra $A$, seen as a Gelfand ring.} 
We denote by  
\[
\FFdisc{A}=\FFdisc{X}=\Spa(\A_{\inf}(X))\backslash V([\varpi])
\]
the relative \textit{Fargues-Fontaine  disc  over $X$} seen as a Berkovich space, where $\varpi$ is a pseudo-uniformizer of $X$.\footnote{We note that $\FFdisc{A}$ is a qfd Berkovich space as $A$ is qfd.}
Concretely, if $\varpi$ is a fixed pseudo-uniformizer of $A$, we can write 
\[
\FFdisc{X}=\bigcup_{s} \GSpec (\mathbb{B}_{[0,s]}^{\dagger}(X))
\]
where $[0,s]\subset [0,\infty)$ runs over closed intervals containing $0$, and  $\mathbb{B}_{[0,s]}^{\dagger}(X)=\varinjlim_{s'>s}\bb{B}_{[0,s']}(X)$  is the rational overconvergent localization of the locus $\{|p|^s\leq |[\varpi]|\neq 0\} $ in $\Spa( \A_{\inf}(X))$.   
Similarly, we let 
\[
\FFdisc{X}^{\times}= \FFdisc{X}\backslash V(p)
\]
be the \textit{punctured Fargues-Fontaine disc over $X$} obtained by removing the locus $p=0$ of $\FFdisc{X}$.
On associated arc-stacks, i.e., after perfectoidization, we have that 
\[
\FFdisc{X}^{\diamond}=X\times_{\Marc(\F_p)} \Marc(\Z_p) \mbox{ and }  
\FFdisc{X}^{\times, \diamond}=X\times_{\Marc(\F_p)} \Marc(\Q_p).  
\]
The space $\FFdisc{X}$ admits a Frobenius action $\varphi_X$ leaving the locus where $p=0$ and the open subspace $\FFdisc{X}^{\times}$ stable.
Given a pseudo-uniformizer $\varpi$, the space $\FFdisc{X}$ admits a radius function, i.e., a morphism of Gelfand stacks
\[
\kappa_{\varpi}\colon \FFdisc{X} \to [0,\infty)_{\Betti}
\]  
such that for $x\in |\FFdisc{X}|$ a point in the Berkovich spectrum we have that $\kappa_{\varpi}(x)= -\frac{\log{|[\varpi^{\flat}]|_x}}{\log |p|_{x}}$ (where $\varpi^{\flat}$ is any pseudo-uniformizer in $A^{\flat}$ determined up to a norm-$1$-unit by the fact that $|\varpi^{\flat}(x)|=|\varpi(x)|$ for $x\in \mathcal{M}(X^{\flat})\cong \mathcal{M}(X)$).
In particular, $\varphi_X$  commutes with the multiplication by $p$ in $[0,\infty)_{\Betti}$,  $\kappa^{-1}_{\varpi}(0)= \{|p|=0\}$ is the locus where $|p|=0$ (i.e., the overconvergent neighborhood of $p=0$),  and $\kappa^{-1}_{\varpi}((0,\infty))=\FFdisc{X}^{\times}$.

	\begin{remark}
	If $T$ is an affinoid perfectoid space over $\Q_p$ and $X=T^\flat$, we can make a canonical choice of radius map on $\FFdisc{X}$: if $T$ is a $\Q_p(p^{1/p^\infty})$-algebra, take $\varpi$ to be a choice of $p^\flat$. Since any two choices differ by a compatible system of $p$-th roots of unity, which has norm $1$, the associated radius map does not depend on this choice and descends to the case where $T$ is any affinoid perfectoid space over $\Q_p$. In what follows, we will always be given a perfectoid space over $\Q_p$, so we fix this choice of radius map for the associated Fargues--Fontaine disc of its tilt and we will denote it  simply by $\kappa$.
	\end{remark}

We note that the functors $X\mapsto \FFdisc{X},\ X\mapsto \FFdisc{X}^\times$ globalize to all (light) perfectoid spaces over $\F_p$ in such a way that $X\mapsto \FFdisc{X}$ and $X\mapsto \FFdisc{X}^\times$ preserve Tate-analytic covers.

The following property is important.

\begin{lemma}\label{LemmaCategoryFFCurves}
Let $X$ be a perfectoid space over $\Q_p$ with tilt $X^{\flat}$, and let $Y$ be a perfectoid space over $\F_p$.
 Then any morphism of Gelfand stacks $X\to \FFdisc{Y}^{\times}$ extends uniquely to a $\varphi$-equivariant morphism of Gelfand stacks $\FFdisc{X^{\flat}}^{\times}\to \FFdisc{Y}^{\times}$ making the following diagram commutative 
\[
\begin{tikzcd}
X \ar[r] \ar[rd] & \FFdisc{X^{\flat}}^\times \ar[d] \\
			& \FFdisc{Y}^\times,
\end{tikzcd}
\]
where $X\to  \FFdisc{X^{\flat}}^{\times} $ corresponds to the untilt $X$ of $X^{\flat}$, and the morphism  $\FFdisc{X^{\flat}}^{\times}\to \FFdisc{Y}^{\times}$ arises from the unique morphism of perfectoid spaces $X^{\flat}\to Y$ in characteristic $p$ corresponding to the projection 
\[
X^{\diamond}=X^{\flat,\diamond}\to  \FFdisc{Y}^{\times\diamond} =Y^\diamond\times \Marc(\Q_p)\to Y^{\diamond}.
\]
\end{lemma}
\begin{proof}
This follows easily from the isomorphism $\FFdisc{Y}^{\times,\diamond}\cong Y\times_{\Marc(\F_p)}\Marc(\Q_p)$ and the fact that, because $X$ is perfectoid, morphisms $X\to \FFdisc{Y}^{\times}$ correspond to morphisms $X^\diamond\to \FFdisc{Y}^{\times,\diamond}$ over $\Marc(\Q_p)$ (this uses $\widehat{X^\diamond}=X$). 
\end{proof}

\begin{remark}
  \label{sec:class-farg-font-analogy-to-universal-property-of-witt-vectors}
  \Cref{LemmaCategoryFFCurves} should be seen as a (perfect, analytic) analogue of the following well-known property of Witt vectors:
  for any rings $A,R$ any morphism $W(R)\to A$ of rings factors uniquely as $W(R)\to W(A)\to A$, where $W(R)\to W(A)$ is a morphism of $\delta$-rings and $W(A)\to A$ the canonical projection.
\end{remark}

By \cite[Corollary 3.4.10]{drff_stacks}, the functor $\Marc(A)\mapsto \FFdisc{A}^{\times}$  sends geometric realizations of arc-hypercovers of separable totally disconnected perfectoid rings over $\F_p$ to $!$-equivalences of Gelfand stacks.
By descent, this produces a left exact and colimit preserving functor
\begin{equation}\label{eqo23kfo23d}
\FFdisc{-}^{\times}\colon \cat{ArcStk}\to \cat{GelfStk}
\end{equation}
called the \textit{relative punctured Fargues-Fontaine disc}.

\begin{remark}
Beware that, if $X\in \cat{ArcStk}$ is a non-empty arc-stack, the natural map $\widehat{\FF_{X}^\diamond} \to \FF_X$ is not an isomorphism. Indeed, for $X$ a perfectoid space, let $\FFdisc{X}^{\times}$ denote the punctured Fargues--Fontaine disc of $X$ and $\widetilde{\FFdisc{X}}^{\times}$ the  base change along the map $\bb{Q}_p\to \bb{Q}_p(\zeta_{p^{\infty}})$ and $\widetilde{\FFdisc{X}}^{\times,u}$ its uniform completion. The object $\widetilde{\FFdisc{X}}^{\times,u}$ is a perfectoid space, and if we denote $\widetilde{\FF}_X:= \widetilde{\FFdisc{X}}^{\times}/\varphi$ and $\widetilde{\FF}_X^u$ its uniform completion, we have that 
\[
\widehat{\FF_{X}^\diamond}= \widetilde{\FFdisc{X}}^{\times,u}/\bb{Z}_p^{\times,\ob{cont}} 
\]
while 
\[
\FF_{X}= \widetilde{\FFdisc{X}}^{\times}/  \bb{Z}_p^{\times,\sm},
\]
where we recall that $\bb{Z}_p^{\times,\ob{cont}} =\GSpec (C^{\rm cont}(\bb{Z}_p^{\times} , \bb{Q}_p))$  (resp. $\bb{Z}_p^{\times,\sm}= \GSpec( C^{\lc}(\bb{Z}_p^{\times} , \bb{Q}_p))$ ) is the analytic spectrum of the $\bb{Q}_p$-valued continuous (resp. locally constant) functions of $\bb{Z}_p^{\times}$. In particular,   the map 
\begin{equation}\label{eqMapsComparisonFFStacksHat}
\widehat{\FF_{X}^\diamond}\to \FF_{X}
\end{equation} is not an equivalence. Passing to arc-stacks via colimits, one has a similar description for $\widehat{\FF_X^{\diamond}}$ and $\FF_{X}$ compatible with base change with respect to maps of arc-stacks. This implies that \eqref{eqMapsComparisonFFStacksHat} is not an isomorphism unless $X=\emptyset$.

 The reader is invited to pay attention to the difference between $\FF_{Y^\diamond}$ and $\FF_{Y^\diamond}^\diamond$, if $Y$ is a Gelfand stack (the former is a Gelfand stack, the latter is an arc-stack), despite the notational proximity!
\end{remark}

\begin{remark}
  \label{sec:class-farg-font-comparison-with-albm}
  The theory of quasi-coherent sheaves of \eqref{eqo23kfo23d} encodes a variant of the theory of solid quasi-coherent sheaves on the Fargues-Fontaine curve of \cite{ALBM-FF} (restricted to partially proper spaces), more precisely of the formalism $\mathcal{D}_{(0,\infty)}(-)$ of \textit{loc.\ cit.}.
\end{remark}

Finally, if $X$ is an arc-stack over $\F_p$, we let $\FF_X=\FFdisc{X}^{\times}/\varphi^{\Z}$ be the quotient by the Frobenius of $X$.

\begin{definition}
  \label{sec:class-farg-font-FF-curve}
  We call $\FF_X$ the \textit{relative Fargues-Fontaine curve of $X$} (note that $\FF_X$ is a Gelfand stack).
\end{definition}

\begin{example}
  \label{sec:class-farg-font-ff-curve-for-q-p}
  If $X=\mathcal{M}_\arc(\Q_p)$ with perfectoid cover $\mathcal{M}_\arc(\Q_p^\cycl)$, then $\FF_X\cong \FF_{\Q_p^\cycl}/\Gamma^{\rm cont}$, where $\Gamma:=\Gal(\Q_p^\cycl/\Q_p)\cong \Z_p^\times$ and $\Gamma^{\rm cont}:=\GSpec(C(\Gamma,\Q_p))$ acting naturally on $\FF_{\Q_p^\cycl}$. 
  \end{example}

\subsection{Rational  $\delta$-spaces}\label{section-Thickened-FF}
The following definition of a ``rational  $\delta$-space''  plays an important role in this paper.

\begin{definition}\label{Definition:Generalized-FF-Curvesv2} 
\begin{enumerate}
\item A \textit{rational  $\delta$-space} is the datum $(S, \mathcal{X},\iota, \varphi_\mathcal{X})$ where $S\in \cat{Perf}$, $\mathcal{X}$ is a Berkovich space with an endomorphism $\varphi_\mathcal{X}$ (called ``the Frobenius of the $\delta$-space''), and $\iota\colon \FFdisc{S}^\times \to \mathcal{X}$ is a morphism of Berkovich spaces inducing an isomorphism of diamonds $\FFdisc{S}^{\times,\diamond}\xrightarrow{\sim} \mathcal{X}^{\diamond}$ such that $\varphi_\mathcal{X}\circ \iota\cong \iota\circ \varphi$, where $\varphi\colon \FFdisc{S}^\times\to \FFdisc{S}^\times$ is the morphism induced by the Frobenius on $S$.\footnote{We note that due to the fact that $\FFdisc{S}^\times$ is locally the spectrum of a \emph{static} Banach $\Q_p$-algebra, this is a property and not a datum.}
  
  We often omit $\iota, \varphi_\mathcal{X}$ from the notation and write $\varphi=\varphi_\mathcal{X}$.
  We sometimes informally say that $\mathcal{X}$ is a \textit{thickening} of $\FFdisc{S}^\times$. 
\item  A morphism of rational  $\delta$-spaces $(S,\mathcal{X})\to (S',\mathcal{X}')$ is a morphism $\alpha\colon S \to S^\prime$ of arc-stacks, a $\varphi$-equivariant morphism $\mathcal{X}\to \mathcal{X}'$ of Gelfand stacks\footnote{The $\varphi$-equivariance is an additional datum.} and a commutative\footnote{The commutativity is a property, not a datum as the Berkovich space $\FFdisc{S}^\times$ is locally given a static Banach $\Q_p$-algebra.} diagram of Gelfand stacks
  \[
\begin{tikzcd}
\FFdisc{S}^\times \ar[r] \ar[d] & \mathcal{X} \ar[d] \\ 
\FFdisc{S'}^\times \ar[r] & \mathcal{X}^{\prime},
\end{tikzcd}
\]
where the left horizontal arrow $\FFdisc{S}^\times \to \FFdisc{S'}^\times$ is induced by $\alpha$.
\item A \textit{perfect rational $\delta$-space} is a rational  $\delta$-space $(S,\mathcal{X})$ such that $\varphi\colon \mathcal{X}\to \mathcal{X}$ is an automorphism.
\item We let $\cat{Spc}^\delta$ be the category of rational  $\delta$-spaces, and we let $\cat{PerfSpc}^{\delta}$ be the full subcategory of  perfect rational $\delta$-spaces.  
\end{enumerate} 
\end{definition}

\emph{From now on, to lighten notation, we will always omit the adjective ``rational'' and only say ``(perfect) $\delta$-space''.}

\begin{remark}
  \label{sec:rati-analyt-prisms-ff-curves}
  If $(S,\mathcal{X})$ is a perfect $\delta$-space, then the $\varphi$-action on $\mathcal{X}$ is totally discontinuous, i.e., each point in the underlying topological space $|\mathcal{X}|$ has a neighborhood disjoint from its non-trivial Frobenius translates. In particular, we can without any harm consider the quotient $\mathcal{X}/\varphi^\Z$, that we call a \emph{thickened Fargues--Fontaine curve}.
  
  If $S\in \cat{Perf}$ and $\mathcal{Y}$ is a Berkovich space with a map $\iota: \FF_S \to \mathcal{Y}$ such that $\iota^\diamond$ is an isomorphism, then $(S,\mathcal{X})$ is a perfect $\delta$-space, for $\mathcal{X} = \mathcal{Y} \times_{|\mathcal{Y}|_{\rm Betti}} |\FFdisc{S}^\times|_{\rm Betti}$ (where the map $|\FFdisc{S}^\times|_{\rm Betti} \to |\mathcal{Y}|_{\rm Betti}$ comes from the identification $|\mathcal{Y}|_{\rm Betti}=|\FF_S|_{\rm Betti}$ induced by $\iota^\diamond$), with its obvious map from $\FFdisc{S}^\times$. For example, $\mathcal{Y}$ could be (locally) a $\dagger$-thickening of $\FF_S$. But not all perfect $\delta$-spaces are $\dagger$-thickenings: see for instance \Cref{ExamplePerfectoidDisc} below.
\end{remark}

\begin{remark}\label{RemarkInscriptionsHowe}
  In \cite{howe2025inscriptiontwistorspadicperiods} Howe introduces the theory of \textit{inscribed $v$-sheaves}.
  This is ressembling a theory of sheaves on thickened Fargues--Fontaine curves $(S,\mathcal{X})$ where one restricts the class of allowed thickenings by imposing $\mathcal{X}$ to be a finite flat augmented nilpotent thickening over $\FF_S$.
\end{remark}

\begin{example}[Classical $\delta$-spaces]\label{ExampleClassicalFFCurves}
  If $S$ is a perfectoid space, the pair $(S,\FFdisc{S}^\times)$ is a perfect $\delta$-space.
  Objects in $\cat{PerfSpc}^{\delta}$ obtained in this way are called \textit{classical} and we denote the category of such by $\cat{Spc}^{\delta, \mathrm{cl}}$. 
\end{example}

\begin{example}\label{ExamplePerfectoidDisc}
  Let $S=\Marc(\F_p((\pi^{1/p^{\infty}})))$. Then $\mathbb{D}^{\times}_{\FF,S}$ is given by the open punctured (pre)perfectoid disc over $\Q_p$.
  Consider the ``decompleted'' perfected punctured open unit disc $\mathbb{D}^{\times, <1,\perf}_{\Q_p}$ with Frobenius given by $\varphi(x)=x^p$.
  More precisely, $\mathbb{D}^{\times,<1, \perf}_{\Q_p}$ is the inverse limit over the Frobenius $x\mapsto x^p$ of the punctured open unit disc $\mathbb{D}^{\times,<1}_{\Q_p}$ over $\Q_p$, and we note that due to the definition of Gelfand rings this inverse limits is calculated locally on affinoids as an uncompleted colimit on rings.
  
 We have a natural map $\FFdisc{S}^\times\to \mathbb{D}^{\times,<1, \perf}_{\Q_p}$ inducing an isomorphism on diamonds.
  Hence, $(S, \mathcal{X}:=\mathbb{D}^{\times,<1,\perf}_{\Q_p})$ is a perfect $\delta$-space, without $\mathcal{X}$ being a $\dagger$-thickening of $\FFdisc{S}^\times$. 
\end{example}

\begin{remark}
  \label{sec:analyt-delta-spac-notation-thick-ff-curves}
  In order to lighten our terminology, we will often drop the ``thickened'' when dealing with thickened Fargues--Fontaine curves.
  If we want to highlight that a Fargues--Fontaine curve is associated with a perfectoid space (and thus not thickened), then we say that it is a classical Fargues--Fontaine curve following \Cref{ExampleClassicalFFCurves}.
\end{remark}

\subsection{Realizations of (perfect) analytic $\delta$-stacks}
\label{sec:real-farg-font}

In this subsection we will consider stacks on $\delta$-spaces resp.\ perfect $\delta$-spaces, and ``realize'' them to Gelfand stacks.

This process is not strictly necessary for the aims of this paper, but we develop this viewpoint here for future applications, and for a slightly deeper structure on the analytic prismatization.
Strictly speaking, the results of this paper only need the construction of $\delta$-spaces/perfect $\delta$-spaces (or thickened Fargues--Fontaine curves) for nilperfectoids, cf.\ \Cref{Prop:Universal-FFNilPerfd2}.\footnote{In the analogy explained in \Cref{sec:real-farg-font-1-functors-on-delta-rings}, these correspond to $\delta$-rings of the form $W(R)$ for an animated $p$-nilpotent ring $R$, which is semiperfectoid.}

\begin{remark}
  \label{sec:real-farg-font-1-functors-on-delta-rings}
  A  $\delta$-space can be thought of as an analytic analogue of the spectrum of a ($p$-complete, animated) $\delta$-ring.
  We note that several important constructions in the theory of prismatic cohomology or prismatization are naturally functors on
  such $\delta$-rings $A$:
  \begin{itemize}
  \item we can send $A$ to invertible ideals $I\subseteq A$ such that $(A,I)$ is a prism,
  \item given a $p$-adic formal scheme $X$, we get the functor sending a $\delta$-ring $A$ to prism structures $(A,I)$ on $A$ with a map $\Spf(A/I)\to X$, and
  \item we can send (perfect) $\delta$-rings $A$ to vector bundles on $A[1/p]$.
  \end{itemize}

  Evaluating at Witt vectors $W(R)$ for (animated $p$-nilpotent) rings $R$ ``realizes'' these functors to:
  \begin{itemize}
  \item the prismatization $\Z_p^\prism$ of $\Z_p$,
  \item the prismatization $X^\prism$ of $X$, and
  \item the stack of isocrystals on perfect rings.
  \end{itemize}
\end{remark}

We construct now the ``realization'' of stacks on  $\delta$-spaces resp.\ Fargues--Fontaine curves in detail.

\begin{definition}
  \label{sec:thick-farg-font-definition-shriek-equivalence-for-ff-stacks2}
  \begin{enumerate}
  \item We let
  \[
(-)_\real^p\colon \mathcal{P}(\cat{Spc}^\delta)\to \GelfStk
\]
be the unique  colimit preserving functor sending a $\delta$-space $(S,\mathcal{X})$ to $\mathcal{X}$. 
\item A morphism $f\colon Y\to X$ in $\mathcal{P}(\cat{Spc}^\delta)$ is a $!$-equivalence if its image in $\GelfStk$ under $(-)_{\real}^p$ is an isomorphism (note that this topology is subcanonical).
\item We let $\cat{AnStk}^\delta$ be the localization of $\mathcal{P}(\cat{Spc}^\delta)$ at the $!$-equivalences, and call objects therein \emph{rational analytic $\delta$-stacks}.
\item The \emph{realization functor} on rational $\delta$-stacks is the induced  colimit preserving functor $(-)_\real\colon \cat{AnStk}^\delta\to \GelfStk$.
  \end{enumerate}
\end{definition}

\begin{remark}
  \label{sec:real-farg-font-1-frobenius-on-realization}
  We note that the Frobenius on any  $\delta$-space is naturally an endomorphism of  $\delta$-spaces.
  In particular, any rational $\delta$-stack acquires naturally a Frobenius endomorphism, which also induces a Frobenius on its realization. 
\end{remark}

\begin{remark}
  \label{sec:real-farg-font-2-ff-version-of-realization}
  Replacing the category $\cat{Spc}^\delta$ by its full subcategory $\cat{PerfSpc}^{\delta}$ we can make the same definitions as in \Cref{sec:thick-farg-font-definition-shriek-equivalence-for-ff-stacks2}, and get the category $\cat{PerfAnStk}^\delta$ of \emph{perfect rational analytic $\delta$-stacks} as the localization of $\mathcal{P}(\cat{PerfSpc}^{\delta})$ at the $!$-equivalences, and their realizations.
  Of course, in this case the natural Frobenius endomorphism of a perfect analytic $\delta$-stack is an automorphism.
\end{remark}

\emph{As above, we will from now on always drop the adjective ``rational'' when speaking about (perfect) analytic $\delta$-stacks.}

\begin{remark}
  \label{sec:real-farg-font-1-realization-is-conservative}
  By construction, we see that the realization functor on analytic $\delta$-stacks (resp.\ perfect analytic $\delta$-stacks) is conservative and commutes with   colimits.  
\end{remark}

\begin{remark}
  \label{sec:real-farg-font-1-criticism-for-endomorphism}
  One might object that we do not ask any condition on the endomorphism $\varphi_\mathcal{X}\colon \mathcal{X}\to \mathcal{X}$ of a  $\delta$-space $(S,\mathcal{X})$, and we don't have a good answer for this matter.
  We remark however that the final class of  $\delta$-spaces associated with nilperfectoids, cf.\ \Cref{sec:real-farg-font-2-definition-rational-analytic-prism-for-uniformly-perfectoid}, carries a rather reasonable endomorphism, which makes the theory well-behaved.
\end{remark}

Our next aim will be to construct a (partial) left adjoint to $(-)_\real$, cf.\ \Cref{Prop:Universal-FFNilPerfd2}, that can be thought of as an analytic variant of the functor of Witt vectors.

The following proposition constructs the necessary  $\delta$-spaces for a large class of Gelfand rings.

\begin{proposition}\label{Prop:Universal-FFNilPerfd2}
  Let $A$ be a Gelfand ring whose uniform completion $A^u$ is perfectoid, i.e., $A$ is uniformly perfectoid.
Let $T=\GSpec(A)$ and consider the category
$$
\cat{Spc}^\delta_{T/-}
$$
of  $\delta$-spaces $(S,\mathcal{X})$ with a morphism $T\to \mathcal{X}$ of Berkovich spaces.
Then $\cat{Spc}^\delta_{T/-}$ admits an initial object given by the pair $(T^\flat,\mathcal{X}_T)$, described as follows:
  \begin{enumerate}
 
  \item The perfectoid space $T^{\flat}$ over $\F_p$  is the tilt of $T^u:=\GSpec(A^u)$.
    In particular, we have a natural map of Berkovich spaces 
\begin{equation*}
f\colon T^u\to \FFdisc{T^\flat}^\times
\end{equation*}
induced by the untilt $T^u$ of $T^\flat$.
We let $f_{\varphi}\colon \N\times T^u\to \FFdisc{T^\flat}^\times$ be the unique morphism whose restriction to $\{n\}\times T^u$ is $\varphi^n\circ f$.

\item The thickening  $\mathcal{X}_T$ of $\FFdisc{T^\flat}^\times$ is the pushout in Berkovich spaces of the morphisms $f_\varphi\colon \N\times T^u\to \FFdisc{T^\flat}^\times$ and $\N\times T^u\to \N\times T$.
\end{enumerate}
\end{proposition}
\begin{proof}
  Let $(S,\mathcal{X})$ be a  $\delta$-space and let $T\to \mathcal{X}$ be a morphism.
  We note that $T^u=\widehat{T^{u,\diamond}}$, which implies that the composition
  \[
    T^u\to T\to \mathcal{X}
  \]
  is uniquely determined by the morphism $T^{u,\diamond}\to \mathcal{X}^\diamond=\FFdisc{S}^{\times,\diamond}$, and thus factors as $T^u\to \FFdisc{S}^{\times}\to \mathcal{X}$.
  By \Cref{LemmaCategoryFFCurves} the morphism $T^u\to \FFdisc{S}^{\times}$ factors canonically as $T^u\xto{f} \FFdisc{T^\flat}^\times\to \FFdisc{S}^\times$, with the latter morphism induced by $T^\flat\to S$ (and hence it is $\varphi$-equivariant).
  Clearly, we get a unique extension of $T^u\to \FFdisc{T^\flat}\to \FFdisc{S}^\times\to \mathcal{X}$ to a $\varphi$-equivariant morphism $\N\times T^u\to \mathcal{X}$ factoring over $f_\varphi$.
  Similarly, we get a unique $\varphi$-equivariant morphism $\N\times T\to \mathcal{X}$.
  Combing these with \Cref{LemmaPushoutsBerkovichSpaces}, we see that we get a unique morphism of  $\delta$-spaces
  \[
    (T^\flat,\mathcal{X}_T=(\N\times T)\sqcup_{(\N\times T^u)} \FFdisc{T^\flat}^\times)\to (S,\mathcal{X})
  \]
  as desired.
  \end{proof}

  \begin{definition}
    \label{sec:real-farg-font-2-definition-rational-analytic-prism-for-uniformly-perfectoid}
    For $T=\GSpec(A)$ uniformly perfectoid as in \Cref{Prop:Universal-FFNilPerfd2} we call $(T^\flat, \mathcal{X}_T)$ the  $\delta$-space associated with $T$. 
      \end{definition}
  
  \begin{remark}\label{rk:radius-discs}
In the notations of \Cref{Prop:Universal-FFNilPerfd2}, the de Rham stack of $\mathcal{X}_T$ is the same as the de Rham stack of $\FFdisc{S}^\times$, so we get a radius map
\[
\kappa: \mathcal{X}_T \to (0,\infty).  \qedhere
\]
\end{remark}

Using the radius map, one can nicely vizualise the $\delta$-space $\mathcal{X}_T$ from \Cref{Prop:Universal-FFNilPerfd2}, see Figure 1. In particular,  the preimage of $(0,\infty)\backslash \{1,p,p^2,\ldots\}$ agrees with the same region on the  classical Fargues--Fontaine curve of $\GSpec(A^u)$. The $n$-th copy of $\GSpec(A)$ in the thickened curve maps to $p^n$ via $\kappa$. However, $\kappa^{-1}(p^n)$ is much larger than $\GSpec(A)$, and is equal to the spectrum of the ring $A\times_{A^u} \bb{B}_{[p^n,p^n],A^u}^{\dagger}$, where  $\bb{B}_{[p^n,p^n],A^u}^{\dagger}$ is the colimit of rings $\bb{B}_{[s,r],A^u}$  defined by intervals with ends  $s<p^n<r$ converging to $p^n$. 

      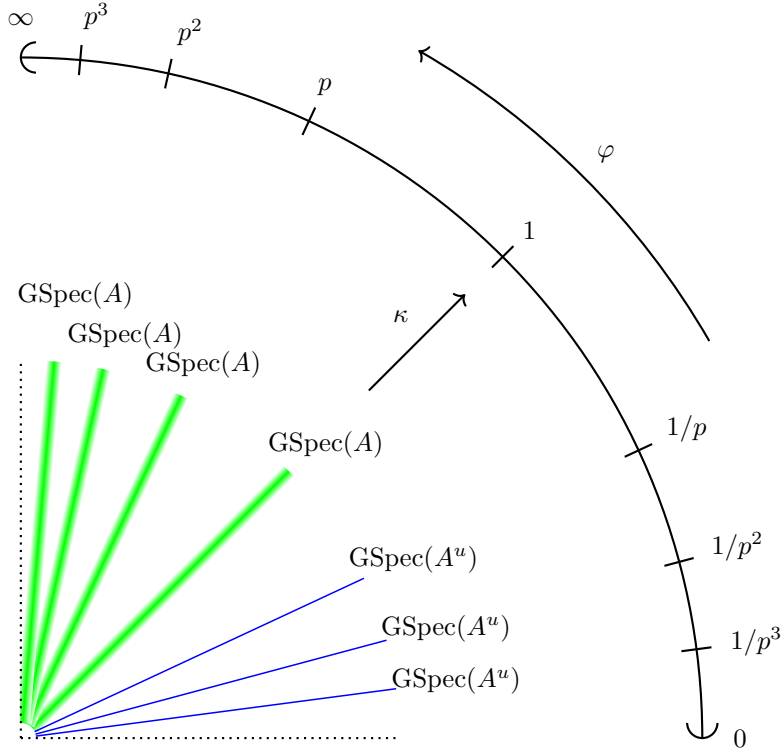
\begin{figure}[h!t]  \label{Figure1FFDiscs}
\begin{center}
\begin{tikzpicture}
  \def\r{9}
  \def\shortr{5} % Slightly shorter radius
  \def\longr{9.5}
  \def\vlongr{10.5}

  % Draw the figure
  \draw[thick] (\r,0) arc[start angle=0, end angle=90, radius=\r];
  
  % Draw the arrow at the center of the arc, slightly outwards
  \draw[->, thick] (0.866*\vlongr, 1/2*\vlongr) arc[start angle=30, end angle=60, radius=\vlongr];
  
  % Label for the angle
  \node at (1.4142/2*\vlongr + 0.3, 1.4142/2*\vlongr + 0.3) {$\varphi$};

  \def\shorttr{\r-0.2}
  \def\longgr{\r+0.2}

  \draw[thick] (\shorttr,0.2) arc (-180:0:0.2);
  \node at (\longr, 0) {$0$};
  
  \draw[thick] (7.5:\shorttr) -- (7.5:\longgr);
  \node at (7.5:\longr+0.3) {$1/p^3$};

  \draw[thick] (15:\shorttr) -- (15:\longgr);
  \node at (15:\longr+0.3) {$1/p^2$};

  \draw[thick] (25:\shorttr) -- (25:\longgr);
  \node at (25:\longr+0.2) {$1/p$};

  \draw[thick] (45:\shorttr) -- (45:\longgr);
  \node at (45:\longr) {$1$};

  \draw[thick] (65:\shorttr) -- (65:\longgr);
  \node at (65:\longr) {$p$};

  \draw[thick] (77.5:\shorttr) -- (77.5:\longgr);
  \node at (76.5:\longr+0.1) {$p^2$};
  
  \draw[thick] (85:\shorttr) -- (85:\longgr);
  \node at (84:\longr+0.1) {$p^3$};

  \draw[thick] (0.2, \longgr) arc (90:270:0.2);
  \node at (0, \longr) {$\infty$};

  \begin{scope}

  \clip (0,0) rectangle (\longr,\longr);
  
  % Draw rotated lines with the middle part darker (last one drawn should be on top)
  \begin{scope}[transform canvas={rotate=45}]
        \shade[shading=axis, bottom color=white, top color=white, middle color=green, shading angle=0] (0.2, -0.1) rectangle (\shortr, 0.1);
  \end{scope} 

  % Second layer, drawn on top of the previous one
  \begin{scope}[transform canvas={rotate=65}]
        \shade[shading=axis, bottom color=white, top color=white, middle color=green, shading angle=0] (0.2, -0.1) rectangle (\shortr, 0.1);
  \end{scope} 

  % Third layer, drawn on top of the previous one
  \begin{scope}[transform canvas={rotate=77.5}]
        \shade[shading=axis, bottom color=white, top color=white, middle color=green, shading angle=0] (0.2, -0.1) rectangle (\shortr, 0.1);
  \end{scope}

  % Fourth layer, drawn on top of the previous one
  \begin{scope}[transform canvas={rotate=85}]
        \shade[shading=axis, bottom color=white, top color=white, middle color=green, shading angle=0] (0.2, -0.1) rectangle (\shortr, 0.1);
  \end{scope}

  \draw[line width=0.2mm, color=blue] (7.5:0.2) -- (7.5:\shortr);

  \draw[line width=0.2mm, color=blue] (15:0.2) -- (15:\shortr);

  \draw[line width=0.2mm, color=blue] (25:0.2) -- (25:\shortr);

  \end{scope}
  
  \def\shortishr{\shortr+0.5}
  
  \node at (7.5:\shortishr+0.3) {$\GSpec(A^u)$};
  
  \node at (14.5:\shortishr+0.3) {$\GSpec(A^u)$};
  
  \node at (24.5:\shortishr+0.2) {$\GSpec(A^u)$};
  
  \node at (44:\shortishr+0.1) {$\GSpec(A)$};
  
  \node at (64:\shortishr) {$\GSpec(A)$};
  
  \node at (75.5:\shortishr) {$\GSpec(A)$};
  
  \node at (83:\shortishr+0.4) {$\GSpec(A)$};
  
  \draw[thick, dotted] (0,0) -- (\shortr,0);
  \draw[dotted, thick] (0,0) -- (0,\shortr); % Dotted vertical radius
  
  \draw[->, thick] (45:\shortishr+1) -- (45:\shorttr-0.5);
  \node at (48:7.5) {$\kappa$};

\end{tikzpicture}
\end{center}
  \caption{\footnotesize{A picture of $\mathcal{X}_T$, stolen from \cite{hauck}.}}  
\end{figure}

  \begin{remark}
    \label{sec:real-farg-font-1-perfected-version}
    \Cref{Prop:Universal-FFNilPerfd2} holds similarly for perfect $\delta$-spaces, i.e., by restricting to perfect  $\delta$-spaces.
    In this case, the initial perfect $\delta$-space associated with a uniformly perfectoid $T=\GSpec(A)$ is given by the \emph{relative perfect $\delta$-space $\FFdisc{A}^{\times}$ for $A$} (we also write $\FFdisc{T}^\times$), and it has the description
    \[
     (T^\flat, \FFdisc{A}^\times)=(T^\flat,(\Z\times T)\sqcup_{(\Z\times T^u)} \FFdisc{T^\flat}^\times).
   \]
   We will also use the \emph{Fargues--Fontaine curve $\FF_A$} associated with $A$, i.e., the quotient
   \[
     \FF_A:=\FFdisc{A}^\times/\varphi^\Z,
   \]
   which is a Berkovich space.
  \end{remark}

  \begin{remark}
    \label{sec:thick-farg-font-relation-to-witt-vectors2}
    The functor $T \mapsto (T^\flat, \mathcal{X}_T)$ should be thought of as analogous to the Witt vectors functor $R\mapsto W(R)$ on $p$-nilpotent rings. We refer to \Cref{sec:thick-farg-font-relation-to-witt-vectors2-intro} for a description of $\mathcal{X}_T$ via Witt vectors. We simply remark here that to make this description true one needs to take the pushout in Berkovich spaces rather than in Gelfand stacks; this explains our restriction to Berkovich spaces in the definition of $\delta$-spaces.
    Witt vectors do not play an important role in this paper, but this changes (for a more complicated version of them) for the analytic prismatization with $\Z_p$-coefficients.
      \end{remark}

  We can deduce now a useful description of the realization functor for analytic $\delta$-stacks (and similar for perfect analytic $\delta$-stacks).

  \begin{proposition}\label{PropAdjunctionFFStacks2}
The realization functor $(-)_{\real}\colon \cat{AnStk}^\delta\to \GelfStk$ on analytic $\delta$-stacks commutes with finite limits and it is thus a left adjoint in a morphism of topoi.
Explicitly, it sends an object $Y\in \cat{AnStk}^\delta$ to the sheafification in Gelfand stacks (i.e.\ with respect to $!$-equivalences) of the presheaf sending a uniformly perfectoid space $T=\GSpec(A)$ to $Y(T^\flat, \mathcal{X}_T)$. 
\end{proposition}

\begin{proof}
  By \cref{Prop:Universal-FFNilPerfd2}, there is a partial left adjoint of $(-)_\real$ when restricted to  $\delta$-spaces and uniformly perfectoid affinoid spaces.
  At the level of presheaves this induces an adjunction
  \[
    \Phi'\colon\mathcal{P}(\cat{BerkSp}^{\ob{up}}) \rightleftharpoons \mathcal{P}(\cat{Spc}^\delta)\colon \Psi',
  \]
  where $\cat{BerkSp}^{\ob{up}}$ is the category of affinoid uniformly perfectoid spaces.
  Here $\Psi'$ sends some $(S,\mathcal{X})\in \cat{Spc}^\delta$ to the functor represented by $\mathcal{X}$ on affinoid uniformly perfectoids.
  By definition of the category of analytic $\delta$-stacks, the functor $\Psi'$ sends $!$-equivalences to $!$-equivalences,  and the functor $\Psi^\prime$ commutes with limits being a right adjoint. Therefore, the realization functor $(-)_{\real}$ is the localization along $\n{P}(\cat{Spc}^{\delta})\to \cat{AnStk}^{\delta}$ of the composite functor
  \[
  \n{P}(\cat{Spc}^{\delta})\xrightarrow{\Psi^{\prime}} \n{P}(\cat{BerkSp}^{\ob{up}})\to \cat{GelfStk}
  \]
  where the second arrow is the sheafification functor. Since sheafification commutes with colimits and finite limits, we deduce that $(-)_{\real}\colon \cat{AnStk}^{\delta}\to \cat{GelfStk}$ is a colimit and finite limit preserving functor. In addition, the above factorization produces the explicit description of $(-)_{\real}$ as stated in the proposition, finishing the proof. 
  \end{proof}

\begin{remark}
  \label{sec:real-farg-font-1-analog-of-realization-for-ff-stacks}
  Of course, \Cref{PropAdjunctionFFStacks2} holds similarly for the perfect version with  $\delta$-spaces replaced by perfect $\delta$-spaces.
  Due to its importance, we state this analogue explicity:
  If $Y\in \cat{PerfAnStk}^\delta$ is a perfect analytic $\delta$-stack, then its realization is the sheafification for $!$-equivalences of the functor
  \[
    A\mapsto Y(\Marc(A)^\flat, \FFdisc{A}^\times)
  \]
  on uniformly perfectoid Gelfand rings $A$ (or even nilperfectoids).
\end{remark}

Note that it is not even clear that for $Y\in \cat{AnStk}^\delta$ the functor
\[
  T=\GSpec(A)\mapsto Y(T^\flat, \mathcal{X}_T)
\]
on affinoid uniformly perfectoid Gelfand stacks is a sheaf for the $!$-topology.\footnote{The issue is that $T\mapsto \mathcal{X}_T$ does not preserve fiber products, even over a perfectoid base. One can still show that if $T'\to T$ is descendable, then $\mathcal{X}_{T'}\to \mathcal{X}_T$ is descendable.}

\begin{remark}\label{RemarkRelationPerfectNonPerfectVariants}
 The  inclusion $\cat{PerfSpc}^{\delta}\to \cat{Spc}^{\delta}$ has a right adjoint sending a   $\delta$-space $(S,\n{X})$ to $(S,\varprojlim_{\varphi} \n{X})$. We get an adjunction 
 \[
 \Phi\colon \cat{PerfStk}^{\delta} \rightleftharpoons \cat{AnStk}^{\delta} \colon (-)^{\perf},
 \]
where $\Phi$ is the left Kan  extension from the inclusion of $\delta$-spaces, and  $(-)^{\perf}$ is the  \textit{perfection} functor sending an analytic $\delta$-stack $Y$ to its restriction to perfect $\delta$-spaces  (which is already a sheaf for the $!$-topology).  For  an analytic $\delta$-stack, $X$ the previous produces a natural functor $\Phi(X)^{\perf}\to X$ which factors through $\Phi (X)^{\perf}\to \varprojlim_{\varphi} X$.   

It is clear  from the construction that we have a commutative diagram of realization functors
\[
\begin{tikzcd}
\cat{PerfStk}^{\delta}  \ar[rd,"(-)_{\real}"'] \ar[r,"\Phi"]& \cat{AnStk}^{\delta} \ar[d,"(-)_{\real}"] \\ 
& \cat{GelfStk}
\end{tikzcd}
\]
\end{remark}

In order to adress the (potentially annoying) sheafification at $!$-equivalences in \Cref{PropAdjunctionFFStacks2} we introduce the following technical terminology. 

\begin{definition}
  \label{sec:real-farg-font-1-admissible-delta-stacks-etc}
  We call $Y\in \cat{AnStk}^\delta$ (resp.\ $Y\in \cat{PerfAnStk}^\delta$) admissible if $Y_\real$ is the functor $T=\GSpec(A)\mapsto Y(T^\flat, \mathcal{X}_T)$ (resp.\ $T=\GSpec(A)\mapsto Y(T^\flat, \FFdisc{A}^\times)$) on nilperfectoid $\Q_p$-algebras. 
  \end{definition}

We can record the following properties.

\begin{lemma}
  \label{sec:real-farg-font-2-properties-of-admissible-delta-stacks}
  \begin{enumerate}
  \item Admissible analytic $\delta$-stacks resp.\ perfect analytic $\delta$-stacks are stable under limits.
  \item Assume that $Y\in \cat{AnStk}^\delta$ is admissible, and let $Y^{\perf}\in \cat{PerfAnStk}^\delta$ be its perfection.
      Then the canonical morphism $(Y^{\perf})_\real\to \varprojlim_{\varphi}Y_\real$ is an equivalence. 
  \end{enumerate}
\end{lemma}
\begin{proof}
  The first assertion is clear because limits in analytic $\delta$-stacks and Gelfand stacks are computed pointwise.
  For the second assertion, let $A$ be a nilperfectoid $\Q_p$-algebra, and $T=\GSpec(A)$.
  Then we calculate
  \begin{align*}
    \varprojlim\limits_{\varphi}Y_\real(T) = & \varprojlim\limits_{\varphi}Y(T^\flat, \mathcal{X}_T)\\
    = & Y(\varinjlim\limits_{\varphi}(T^\flat, \mathcal{X}_T)) \\
    = & Y(T^\flat, \FFdisc{T}^\times)\\
    = & (Y^{\perf})_\real(T)
  \end{align*}
  Indeed, the first equality comes from the fact that limits of Gelfand stacks are computed pointwise, the second comes from the fact that the Frobnius of $Y$ induces on $Y(T^\flat, \mathcal{X}_T)$ canonically the same map as the Frobenius of $(T^\flat, \mathcal{X}_T)$, and the third equality follows from the explicit descriptions in \Cref{Prop:Universal-FFNilPerfd2} resp.\ \Cref{sec:real-farg-font-1-perfected-version}. 
\end{proof}

In other words, the second part of \Cref{sec:real-farg-font-2-properties-of-admissible-delta-stacks} implies that on realizations the perfection of an admissible analytic $\delta$-stack can be calculated conveniently by restrictions to perfect $\delta$-spaces.

\begin{remark}
  \label{sec:real-farg-font-1-different-incarnations-of-ff-curves}
  Different incarnations of perfect $\delta$-spaces and Fargues--Fontaine curves play a role in this paper, and we apologize to the reader if this creates confusion.
  For example, if $A$ is a uniformly perfectoid ring, one has the following versions:
  \begin{enumerate}
  \item the arc-stack $\FFdisc{\Marc(A^u)}^{\times, \diamond}=\Marc(A^u)\times \Marc(\Q_p)$, this makes sense for $\Marc(A^u)$ replaced by any arc-stack,
  \item the Gelfand stack $\FFdisc{\Marc(A)}^\times=\FFdisc{\Marc(A^u)}^\times $, which is the classical punctured Fargues--Fontaine disc, i.e., the classical adic space realized as a Gelfand stack, and this construction extends to any arc-stack, 
  \item the Berkovich space $\FFdisc{A}^\times$ from \Cref{sec:real-farg-font-1-perfected-version}, this is only defined for $A$ uniformly perfectoid (e.g., nilperfectoid).
  \end{enumerate}
  These objects are linked as follows: the perfectoidizations of $\FFdisc{\Marc(A)}^\times$ and $\FFdisc{A}^\times$ are isomorphic to $\FFdisc{\Marc(A^u)}^{\times,\diamond}$, and by adjunction we get from here a natural morphism
  \[
    \FFdisc{\Marc(A^u)}^\times\to \FFdisc{A}^\times,
  \]
  which is an isomorphism if and only if $A$ is perfectoid.
\end{remark}

We note two general constructions for (perfect) analytic $\delta$-stacks.

\begin{definition}
  \label{RemarkNotationDescent}
  Let $T$ be a Gelfand stack (or any functor $\GelfStk\to \cat{Cat}_\infty$). One defines an analytic $\delta$-stack $T_{\cat{Spc}^\delta}$ (or a functor $\cat{Spc}^\delta\to \cat{Cat}_\infty$) as the functor
  \[
    (S,\mathcal{X})\mapsto T(\mathcal{X}).
  \]
  Similarly, one defines the perfect analytic $\delta$-stacks $T_{\cat{PerfSpc}^{\delta}}$  as before by evaluating only on perfect  $\delta$-spaces. One also defines  $T_{\mathcal{FF}}$ as the  the functor sending a perfect  $\delta$-space 
  \[
    (S,\mathcal{X})\mapsto T(\mathcal{X}/\varphi^\Z).
  \]
  \end{definition}
  We note that in general we do not expect to have a direct description of the realizations of $T_{\cat{Spc}^{\delta}}$, $T_{\cat{PerfSpc}^\delta}$ or  $T_{\mathcal{FF}}$.
  We finally remark that the Frobenius automorphism of $T_{\mathcal{FF}}$ is canonically trivialized.

\subsection{Relation to Hyodo--Kato stacks}
\label{sec:hyodo-kato-stacks-1}

We saw that the theory of perfect analytic $\delta$-stacks mimics the theory of sheaves on perfect $\delta$-rings.
But it also resembles the theory of arc-stacks, since nilperfectoid rings are (obviously) closely related to perfectoid rings.
In this subsection we make this idea more precise, by relating our theory of functors on Fargues--Fontaine curves to the Hyodo-Kato stacks introduced in \cite{drff_stacks}.

We recall from \Cref{ExampleClassicalFFCurves} the terminology of classical Fargues--Fontaine curves.

\begin{definition}
  \label{sec:relation-hyodo-kato-reduced-ff-curves}
  We let $\cat{Spc}^{\delta, \mathrm{cl}}\subseteq \cat{Spc}^{\delta}$ be the full subcategory of classical $\delta$-spaces $(S,\mathcal{X})$, namely, those $\delta$-spaces for which the given map $\FFdisc{S}^\times\to \mathcal{X}$ is an isomorphism, i.e., $\mathcal{X}$ is a classical Fargues--Fontaine disc.
\end{definition}

In fact, this category is equivalent to the category of (qfd) perfectoid spaces in characteristic $p$.

\begin{lemma}\label{LemmaEquivalencePerfectoidFFstk}
  The functor  $\FFdisc{(-)}^\times\colon \cat{Perf}\to \cat{Spc}^{\delta, \mathrm{cl}}$  from the category of perfectoid spaces over $\F_p$ to $\cat{Spc}^{\delta, \mathrm{cl}}$ is an equivalence of categories.
  Furthermore, the inclusion $\cat{Perf}\cong \cat{Spc}^{\delta, \mathrm{cl}} \hookrightarrow \cat{Spc}^{\delta}$ has a right adjoint sending $(S,\mathcal{X})$ to $S$. 
\end{lemma}
\begin{proof}
  By definition objects in $\cat{Spc}^{\delta}$ are pairs $(S,\mathcal{X})$ where $S$ is perfectoid over $\F_p$ and $\mathcal{X}$ is a Berkovich space with endomorphism $\varphi$ together with a $\varphi$-equivariant map $\FFdisc{S}^\times\to \mathcal{X}$ inducing an isomorphism at the level of diamonds.
  Morphisms from $(S,\mathcal{X})$ to $(S',\mathcal{X}')$ consist of maps $S\to S'$ of perfectoid spaces and $\varphi$-equivariant maps $\mathcal{X}\to \mathcal{X}'$ compatible with the morphisms of Fargues--Fontaine discs.
  From this it is clear that the functor $\cat{Perf}\to \cat{Spc}^{\delta, \mathrm{cl}}$ sending $S\mapsto (S,\FFdisc{S}^\times)$ is an equivalence of categories.
  On the other hand, let $S\in \cat{Perf}$ and let $(S',\mathcal{X}')\in \cat{Spc}^{\delta}$.
  Then any map $(S,\FFdisc{S}^\times)\to (S',\mathcal{X}')$ will factor uniquely as a map $(S,\FFdisc{S}^\times)\to (S',\FFdisc{S'}^\times)\to (S',\mathcal{X}')$ by definition of $\cat{Spc}^{\delta}$.
  This proves the desired adjunction.  
\end{proof}

Replacing perfectoid spaces over $\Q_p$ by classical perfect $\delta$-spaces, we arrive at a formalism very similar to the adjunctions
\[
  \widehat{(-)}\vdash (-)^\diamond\vdash (-)^\dR
\]
from \Cref{Subsection:Analytic-de-Rham}.

\begin{proposition}
  \label{sec:relation-hyodo-kato-adjunction-for-ff-stacks2}
  The functor
  \[
    (-)^{\diamond_\FF}\colon \cat{PerfAnStk}^\delta\to \ArcStk,\ Y\mapsto (S\mapsto Y(S,\FFdisc{S}^\times))
  \]
  is well-defined, and admits both a left adjoint $\FFdisc{(-)}^\times$ and a right adjoint $\FFdisc{(-)}^{\times,\dR}$.
  More precisely, we have the following:
  \begin{enumerate}
  \item $\FFdisc{(-)}^\times$ is the unique colimit preserving functor sending $S$ to $\FFdisc{S}^\times$.
    Moreover, the functor $\FFdisc{(-)}^\times\colon \ArcStk\to \cat{PerfAnStk}^\delta$ is fully faithful and preserves finite limits.

  \item $\FFdisc{(-)}^{\times,\dR}$ sends some $Y\in \ArcStk$ to the perfect analytic $\delta$-stack $(S,\mathcal{X})\mapsto Y(S)$.
    Moreover, $\FFdisc{(-)}^{\times,\dR}$ is fully faithful and takes images in admissible perfect analytic $\delta$-stacks.
  \end{enumerate}
 \end{proposition}
 \begin{proof}
   We need to see that $(-)^{\diamond_\FF}$ is well-defined.
   From \Cref{LemmaEquivalencePerfectoidFFstk} we formally get functors   
     \begin{align*}
       j^\ast\colon \mathcal{P}(\cat{Spc}^{\delta})\to \mathcal{P}(\cat{Perf}),\ & \mathcal{G}\mapsto (S\mapsto \mathcal{G}(S,\FFdisc{S}^\times)), \\
       j_!\colon \mathcal{P}(\cat{Perf})\to \mathcal{P}(\cat{Spc}^{\delta}),\ & S\in \cat{Perf}\mapsto (S,\FFdisc{S}^\times),  \\
       j_\ast\colon \mathcal{P}(\cat{Perf})\to \mathcal{P}(\cat{Spc}^{\delta}),\ & Z\mapsto ((S,\mathcal{X}),\mapsto Z(S,\FFdisc{S}^\times)),       
     \end{align*}
   such that $j_!$ is left adjoint to $j^\ast$ and $j^\ast$ is left adjoint to $j_\ast$.
   We need to see that $j^\ast$ and $j_\ast$ restrict to functors between arc-stacks and perfect analytic $\delta$-stacks.
   We first claim that sheaves for the arc-topology resp.\ $!$-topology are preserved.
   For $j_\ast$ this is clear because any $!$-cover on perfect $\delta$-spaces induces an arc-cover on the underlying classical perfect $\delta$-spaces.
   For $j^\ast$ this follows from the fact that arc-covers on perfectoid spaces in characteristc $p$ induce $!$-covers on perfect $\delta$-spaces (\Cref{Subsection:Analytic-de-Rham}).

   The passage to hypersheaves (resp.\ localization at $!$-equivalences) follows as in \cite[Lemma 4.5.1]{drff_stacks}.
   
   We can conclude that $j^\ast$ restricts to $(-)^{\diamond_\FF}$, and deduce the desired descriptions of its left and right adjoint (and the fact that $\FFdisc{(-)}^{\times,\dR}$ takes images in admissible perfect analytic $\delta$-stacks).
   The fully faithfulness of (1) and (2) follows as in \cite[Proposition 4.5.7]{drff_stacks} (replacing perfectoid $\Q_p$-algebras by classical perfect $\delta$-spaces).
   
   Finally, it is left to show that $\bb{D}^{\times}_{\FF,(-)}$ commutes with finite limits. Let $\cat{Perf}$ denote the category of perfectoid spaces in characteristic $p$, by construction the functor $\bb{D}^{\times}_{\FF,(-)}$ is the left Kan extension to its restriction to  $\cat{Perf}$, hence it suffices to show that 
   \[
\bb{D}^{\times}_{\FF,(-)} \colon    \cat{Perf} \to \cat{PerfAnStk}^{\delta}
   \]
   commutes with finite limits.
   It suffices to show that the passage to the Fargues--Fontaine disc commutes with fiber products and finite products of perfectoid spaces.
   The case of fiber products is clear since this is already true for the $\A_{\inf}$-construction of affinoid rings.
   For the case of absolute products, by descent, it suffices to prove it for a product of perfectoid spaces $S_1$ and $S_2$ which are affinoid and admit a pseudo-uniformizer $\pi$, that is, both $S_i$ live over $K=\bb{F}_p((\pi^{1/p^{\infty}}))$. Since we already know the compatibility with fiber products, it  suffices to consider the case $S_1=S_2=\Marc(K)$.
   In this case, we have that 
   \[
   \Marc(K)\times \Marc(K) = \bb{D}^{\times,<1}_{K,\perf}
   \]
   is isomorphic to the open punctured  perfectoid unit disc over $K$, with coordinate $X^{\flat}$. By taking the coordinate $[X^{\flat}]$ on Witt vectors, we have an isomorphism
   \[
   \FFdisc{\Marc(K)\times \Marc(K)}^{\times} = \bb{D}_{\perf}^{u,\times,<1}\times_{\GSpec(\Q_p)} \FFdisc{K}^{\times},
   \]
where $\bb{D}_{\perf}^{u,\times,<1}$ is the uniformly completed pre-perfectoid open punctured unit disc over $\bb{Q}_p$, which is precisely isomorphic to $\FFdisc{K}$. One deduces that 
\[
  \FFdisc{\Marc(K)\times \Marc(K)}^{\times}= \FFdisc{K}\times_{\GSpec(\Q_p)} \FFdisc{K}
\]
   proving what we wanted. 
  \end{proof}

 After realization, the left adjoint becomes the functor $\FFdisc{-}^\times$ from \Cref{sec:farg-font-curv} and the right adjoint becomes $(\FFdisc{-}^{\times,\diamond})^\dR = (- \times_{\Marc(\mathbb{F}_p)} \Marc(\Q_p))^\dR$. In particular, if $Z$ is any object of $\cat{AnStk}^\delta$, we have a composition of morphisms of Gelfand stacks
 \begin{equation} \label{composition}
\FFdisc{Z^{\diamond_\FF} }^\times \to Z_{\rm real} \to (\FFdisc{Z^{\diamond_\FF}}^\times)^\dR
\end{equation}
factoring the natural map from $\FFdisc{Z^{\diamond_\FF} }^\times$ to its de Rham stack.

\begin{definition}
  \label{sec:hyodo-kato-stacks-2-de-rham-stack-for-ff-stack}
  Let $Y\in \ArcStk$.
  We define its \textit{de Rham Fargues--Fontaine stack} $\FFdisc{Y}^{\times,\dR}\in \cat{PerfAnStk}^\delta$ as the functor
  \[
(S,\mathcal{X})\mapsto Y(S,\FFdisc{S}^\times),
\]
i.e., as the right adjoint to $(-)^{\diamond_\FF}$ from \Cref{sec:relation-hyodo-kato-adjunction-for-ff-stacks2}.
Given any $Y\in \cat{PerfAnStk}^\delta$, we define its \textit{de Rham stack} as $Y^\dR:=\FFdisc{Y^{\diamond_\FF}}^{\times,\dR}$.
\end{definition}

\begin{remark}
  \label{RemarkCompatibilityNotation}
  The notation of \Cref{sec:relation-hyodo-kato-adjunction-for-ff-stacks2} and  \Cref{sec:hyodo-kato-stacks-2-de-rham-stack-for-ff-stack} for de Rham stacks are compatible.
  That is, for $Y$ an arc-stack the de Rham stack of its perfect $\delta$-space $\FFdisc{Y}^\times$ is precisely $(\FFdisc{Y}^{\times})^{\dR}=\FFdisc{Y}^{\times,\dR}$. 
\end{remark}

\begin{remark}
  \label{sec:relation-hyodo-kato-relation-to-hk-stacks}
  If $Y\in \cat{ArcStk}$, then we can see that
  \[
    (\FFdisc{Y}^{\times,\dR})_\real=(Y\times \Marc(\Q_p))^\dR.
  \]
  Indeed, if $A$ is nilperfectoid, $T=\GSpec(A)$, then we get
  \[
  (\FFdisc{Y}^{\times,\dR})_\real(A)=  (\FFdisc{Y}^{\times,\dR})(T^\flat, \mathcal{X}_T)=Y(T^\flat),
    \]
    which agrees with $\mathrm{Maps}_{\Marc(\Q_p)}(T^u,Y\times \Marc(\Q_p))=(Y\times \Marc(\Q_p))^\dR(A)$.

    We note that this implies that the Hyodo-Kato stack $Y^\HK$ of $Y$ (as defined in \cite{drff_stacks}) is given by $Y^\HK=\FFdisc{Y/\varphi^\Z}^{\times,\dR}$.
\end{remark}

Similar as for de Rham stacks (see \cite[Theorem 5.6.6]{drff_stacks}), the formation of de Rham Fargues-Fontaine stacks commutes with colimits. 

\begin{lemma}\label{Lemma-Colimit-de-Rham-FF}
The functor $\FFdisc{(-)}^{\times,\dR}\colon \ArcStk\to \cat{PerfAnStk}^\delta$ commutes with colimits. 
\end{lemma}
\begin{proof}
  By \Cref{sec:real-farg-font-1-realization-is-conservative} we can check this on realizations.
  By \Cref{sec:relation-hyodo-kato-relation-to-hk-stacks} we need to see that the composition 
  \[
    \ArcStk\xto{(-)\times\Marc(\Q_p)} \ArcStk_{\Q_p}\xto{(-)^\dR} \GelfStk
  \]
  commutes with colimits.
  For the first functor this is clear, and for the second functors this is \cite[Theorem 5.6.6]{drff_stacks}.
 \end{proof}

\begin{remark}
  \label{sec:relation-hyodo-kato-realization-commutes-with-de-rham-stacks}
  \begin{enumerate}
\item   We also note that the realization functor functor $(-)_\real\colon \cat{PerfAnStk}^\delta\to \GelfStk$ commutes with the passage to de Rham stacks.  Indeed, by definition $Y_{\real}$ is the sheafification of the functor $Y_{\real}^{\ob{pre}}\colon A\mapsto Y(\Marc(A), \FFdisc{A}^{\times})$, this functor is already a sheaf when restricting to perfectoid rings as the formation of $A^u\mapsto \FFdisc{A^u}^{\times}$ satisfies arc hyperdescent. Thus, $(Y_{\real}^{\ob{pre}})^{\dR}=(Y_{\real})^{\dR}$. It is clear from this functor of points description that $(Y^{\dR})_{\real}=(Y_{\real})^{\dR}$.

\item We can define de Rham stacks of analytic $\delta$-stacks in a totally similar way. That is, we have a Fargues--Fontaine disc functor
\[
\bb{D}^{\times}_{\FF,(-)}\colon \cat{ArcStk}\to \cat{AnStk}^{\delta}
\]
given by the composite of the Fargues--Fontaine disc functor to $\cat{PerfAnStk}^{\delta}$ and the left Kan extension from perfect analytic $\delta$-stacks to analytic $\delta$-stacks. This functor  commutes with colimits and finite limits, and has as a right adjoint the perfectoidization functor
\[
(-)^{\diamond}\colon \cat{AnStk}^{\delta}\to \cat{ArcStk},
\]
which is the restriction of an analytic $\delta$-stack to classical Fargues--Fontaine curves. The exact same argument as in \Cref{sec:relation-hyodo-kato-adjunction-for-ff-stacks2} implies that the perfectoidization functor admits an additional right adjoint given by the de Rham  Fargues--Fontaine stack functor. This functor commutes with colimits  by the same argument as that of \Cref{Lemma-Colimit-de-Rham-FF}, and even factors as the composite 
\[
\cat{ArcStk} \xrightarrow{\bb{D}^{\times,\dR}_{\FF,(-)}} \cat{PerfAnStk}^{\delta}\to \cat{AnStk}^{\delta},
\]
where the first map is that of \Cref{sec:relation-hyodo-kato-adjunction-for-ff-stacks2}, and the second map a left Kan extension. Both the Fargues--Fontaine disc functor and the de Rham Fargues--Fontaine stack functor with  values in $\cat{AnStk}^{\delta}$ are fully faithful

\item  Let $Y$ be an analytic $\delta$-stack.
  Then the passage to the de Rham stack and the realization functor agree by the exact same argument as in (1). Furthermore, the natural maps of Gelfand stacks 
\[
\FFdisc{Y^{\diamond}}^{\times}\to Y_{\real}\to (Y^{\dR})_{\real}
\]
induce equivalences after passing to analytic de Rham stacks. In particular,
\[
(Y_{\real})^{\dR}=\FFdisc{Y^{\diamond}}^{\times,\dR}. 
\]

\end{enumerate} 
\end{remark}

\subsection{Perfect complexes on Fargues--Fontaine curves}
\label{sec:perf-compl-farg}

Many important constructions in prismatization or geometrization deal with perfect complexes on Fargues--Fontaine curves, e.g., Cartier divisors.

\begin{definition}
  \label{DefinitionFFStk}
  Let $[n,m]\subseteq [-\infty, \infty]$ be any interval with integral boundary (possibly including $\pm \infty$ as the boundary)
  \begin{enumerate}
  \item We let
    \[
      \Perf_{[n,m]}\colon \GelfStk\to \cat{Cat}_\infty
    \]
    be the sheaf of categories on Gelfand stacks whose value on an affinoid Gelfand stack $\GSpec(A)$ is $\Perf_{[n,m]}(A)$.
  \item We let
    \[
      \Perf_{[n,m],\mathcal{FF}}\colon \cat{PerfSpc}^\delta\to \cat{Cat}_{\infty}
    \]
    be the sheaf of categories on perfect analytic $\delta$-stacks whose value on a perfect $\delta$-space $(S,\mathcal{X})$ is $\Perf_{[n,m]}(\mathcal{X}/\varphi^\Z)$.
  \end{enumerate}
  When $[n,m]=[-\infty, \infty]$ we simply write $\Perf(-)$ resp.\ $\Perf_{\mathcal{FF}}$, when $n=m=0$ we simply write $\ob{VB}(-)$ resp.\ $\ob{VB}_{\mathcal{FF}}$ (in which we speak of the \textit{stack of vector bundles on Fargues--Fontaine curves}).
\end{definition}

Of course, one can also consider perfect complexes on perfect $\delta$-spaces (or even  $\delta$-spaces), i.e., the functor
    \[
      \Perf_{\cat{PerfSpc}^\delta} \colon \cat{PerfSpc}^\delta\to \cat{Cat}_{\infty},\ (S,\mathcal{X})\mapsto \Perf(\mathcal{X}),
    \]
    but for simplicity of exposition we concentrate on perfect complexes on Fargues--Fontaine curves in this section.

    \begin{remark}
      \label{sec:perf-compl-farg-2-descent-for-perfect-complexes}
      Perfect complexes (of fixed amplitude) form a sheaf of categories on Gelfand stacks, cf.\ \Cref{sec:gener-cart divis} and \ \cite[Lemma 4.2.11]{drff_stacks}, and hence the same is true for $\Perf_{[n,m],\mathcal{FF}}(-)$.
      Restricting to the core categories, i.e., the maximal subanima, yields the perfect analytic $\delta$-stack $\Perf_{[n,m],\mathcal{FF}}(-)^{\simeq}$.
    \end{remark}

The realization of $\Perf_{\mathcal{FF}}$ has a very convenient description (in particular, $\Perf_{\mathcal{FF}}^\simeq(-)$ is admissible in the terminology of~\Cref{sec:real-farg-font-2-properties-of-admissible-delta-stacks}).

\begin{proposition}
  \label{sec:perf-compl-farg-1-perfect-complexes-on-ff-curves}
  Let $A$ be a nilperfectoid $\Q_p$-algebra (or just $A$ uniformly perfectoid).
  The diagram of categories 
\[\begin{tikzcd}
	{\Perf(\FF_A)} & {\Perf(\FF_{A^u})}  \\
	{\Perf(\GSpec(A))}  & {\Perf(\GSpec(A^u))}
	\arrow[from=1-1, to=1-2]
	\arrow[from=1-1, to=2-1]
	\arrow[from=1-2, to=2-2]
	\arrow[from=2-1, to=2-2]
      \end{tikzcd}\]
    is cartesian.
    In particular, the functor $A\mapsto \Perf(\FF_A)$ is a sheaf of categories on the $\infty$-topos of Gelfand stacks and we have a cartesian diagram of Gelfand stacks 
    \[
    \begin{tikzcd}
(\Perf_{\mathcal{FF}}^{\simeq})_{\real}  \ar[r] \ar[d] & (\Perf_{\mathcal{FF}}^{\simeq})^{\dR}_{\real}   \ar[d] \\
   \Perf^{\simeq}    \ar[r]   & \Perf^{\simeq, \dR}.
    \end{tikzcd}
    \]
  \end{proposition}
  \begin{proof}
    By covering $\FF_A$ by $\dagger$-rational localizations and using descent on perfect complexes, the pullback diagram of categories follows from  \Cref{sec:analyt-delta-spac-descent-for-fiber-product}. 

    To check that $A\mapsto \Perf(\FF_A)$ is a sheaf of categories on nilperfectoid $\Q_p$-algebras, we can use \cite[Proposition 4.4.7.]{drff_stacks} to conclude that $A\mapsto \Perf(\FF_{A^u})$ and $A\mapsto \Perf(A^u)$ are sheaves of categories on nilperfectoid $\Q_p$-algebras.
    Similarly, $A\mapsto \Perf(A)$ is a sheaf of categories, so that we can conclude.
    For the last assertion we note that by virtue of $A\mapsto \Perf(\FF_A)$ being a sheaf of categories on $\GelfStk$ the sheafification for $!$-equivalences in \Cref{PropAdjunctionFFStacks2} is not necessary, and thus the realization of the perfect analytic $\delta$-stack $\Perf_{\mathcal{FF}}^{\simeq}$ is the functor $A\mapsto \Perf(\FF_A)^\simeq$.
  \end{proof}
  
  \begin{remark}
    \label{sec:perf-compl-farg-1-complexes-in-d-+}
  In fact, the conclusion of \Cref{sec:perf-compl-farg-1-perfect-complexes-on-ff-curves} holds for $\ob{D}_+$ by  \Cref{sec:analyt-delta-spac-descent-for-fiber-product}.  In particular, since the perfect amplitude of a perfect complex is preserved under base change, \Cref{sec:perf-compl-farg-1-perfect-complexes-on-ff-curves} also holds for the stacks $\Perf_{[n,m]}$. 
  \end{remark}

  \begin{corollary}
    \label{sec:perf-compl-farg-1-calculation-for-hom}
    Let $A$ be a nilperfectoid $\Q_p$-algebra, and $\mathcal{E}\in \Perf(\FF_A)$ with base changes $\mathcal{E}_A\in \Perf(A),\ \mathcal{E}_{A^u}\in \Perf(A^u),\ \mathcal{E}_{\FF_{A^u}}\in \Perf(\FF_{A^u})$.
    Then we have a pullback diagram
 \[\begin{tikzcd}
	\Gamma(\FF_A, \mathcal{E}) & {\Gamma(\FF_{A^u},\mathcal{E}_{\FF_{A^u}})} \\
	{\Gamma(\GSpec(A),\mathcal{E}_A)} & {\Gamma(\GSpec(A^u),\mathcal{E}_{A^u})}.
	\arrow[from=1-1, to=1-2]
	\arrow[from=1-1, to=2-1]
	\arrow[from=1-2, to=2-2]
	\arrow[from=2-1, to=2-2]
      \end{tikzcd}\]
   \end{corollary}
\begin{proof}
  This follows from the fully faithfulness of $\Perf(\FF_A)\to \Perf(\FF_{A^u})\times_{\Perf(A^u)}\Perf(A)$  in \Cref{sec:perf-compl-farg-1-perfect-complexes-on-ff-curves}.
\end{proof}

\begin{example}
  \label{sec:perf-compl-farg-1-global-section-on-ff-curves}
  \Cref{sec:perf-compl-farg-1-calculation-for-hom} is already interesting for $\mathcal{E}=\s{O}_{\FF_A}$: if $A$ is strictly totally disconnected nilperfectoid,  then $\Gamma(\FF_{A^u},\s{O}_{\FF_{A^u}})\cong \underline{\Q_p}^\dR(A)$ (using \cite[Chapter II]{fargues2021geometrization} and the fact that $A$ is strictly totally disconnected to ensure that $\Gamma(\FF_{A^u},\s{O}_{\FF_{A^u}})$ is concentrated in degree $0$), which implies that
  \[
\Gamma(\GSpec(A),\s{O}_{\FF_A})\cong A\times_{A^u}\underline{\Q_p}^\dR(A)=\Q_p^\la(A),
\]
cf.\ \Cref{sec:analytic-de-rham-locally-analytic-version}.
We note that this example convinced us that our approach of using Fargues--Fontaine curves for nilperfectoids, or functors on Fargues--Fontaine curves, might be relevant for a geometrization of the locally analytic local Langlands correspondence over $\Q_p$. 
\end{example}

\begin{remark}
  \label{sec:perf-compl-farg-1-general-case}
  As remarked in the beginnning of this section, the results of this section hold for the variants
  \[
    \Perf_{\cat{PerfSpc}^\delta}\colon \cat{PerfAnStk}^\delta\to \cat{Cat}_\infty,\ (S,\n{X})\mapsto \Perf(\n{X}),
  \]
  and
  \[
    \Perf_{\cat{Spc}^\delta}\colon \cat{AnStk}^\delta\to \cat{Cat}_\infty,\ (S,\n{X})\mapsto \Perf(\n{X}).
  \]
  In particular, these are admissible and for each nilperfectoid ring $A$ over $\Q_p$ with $T:=\GSpec(A)$  we have equivalences
  \[
    \Perf_{\cat{PerfSpc}^\delta,\real}(T)=\Perf(\FFdisc{T}^\times)\cong \Perf(\FFdisc{\Marc(T)}^\times)\times_{\Perf(\Z\times T^u)} \Perf(\Z\times T)
  \]
  and
  \[
    \Perf_{\cat{Spc}^\delta,\real}(T)=\Perf(\n{X}_T)\cong \Perf(\FFdisc{\Marc(T)}^\times)\times_{\Perf(\N\times T^u)} \Perf(\N\times T).
  \]
\end{remark}

\subsection{Locally analytic Banach--Colmez spaces}\label{subsec-loc-an-BC-spaces}

Let $A$ be a nilperfectoid $\Q_p$-algebra.
Using the uniformization $\FFdisc{A}^\times\to \FF_A$ one can define plenty of vector bundles on $\FF_A$.

\begin{definition}
  \label{sec:perf-compl-farg-1-definition-vector-bundles-on-ff-a} 
  \begin{enumerate}
\item   For $\lambda\in \Q$, let $(D_\lambda,\varphi_{D_\lambda})$ be the standard isocrystal over $\Q_p$ of slope $\lambda$.
  Let $A$ be a nilperfectoid $\Q_p$-algebra.
  We define
  \[
    \s{O}_{\FF_A}(\lambda)\in \VB(\FF_A)
  \]
  as the vector bundle obtained by descent along $\FFdisc{A}^\times\to \FF_A$ of the $\varphi$-equivariant vector bundle $(\varphi_{\FFdisc{A}^\times}\otimes_{\Q_p} D_{-\lambda},\varphi \otimes \varphi_{D_{-\lambda}})$.

 \item Let $\lambda\in \Q\cap [0,\infty)$, we define the \textit{absolute Banach--Colmez space} $\BC(\s{O}(\lambda))$ to be the Gelfand stack sending a nilperfectoid ring $A$ to $\Gamma(\FF_A, \s{O}_{\FF_A}(\lambda))$. 
 
 \item Let $\lambda\in \Q\cap (-\infty,0)$, we define the \textit{absolute Banach--Colmez space} $\BC(\s{O}(\lambda)[1])$ to be the Gelfand stack sending a nilperfectoid ring $A$ to $\Gamma(\FF_A, \s{O}_{\FF_A}(\lambda)[1])$. 
 
 \end{enumerate}
\end{definition}

\begin{remark}\label{Remark-BC-Spaces-FFStk}

  The constructions of \Cref{sec:perf-compl-farg-1-definition-vector-bundles-on-ff-a} have an incarnation in the category $\cat{PerfAnStk}^\delta$.
  More precisely, let $\s{O}(\lambda)\in \VB(B\varphi^{\Z})$ be the vector bundle defined by the $\varphi$-equivariant vector bundle $D_{-\lambda}$.
  For $(S,X)$ a (thickened) perfect $\delta$-space with structural map $f \colon X/\varphi^\Z\to B\varphi^{\Z}$ we define $\s{O}_{X/\varphi^\Z}(\lambda):=f^* \s{O}(\lambda)$. 

  Let $\lambda\geq 0$ be a rational number (the discussion for $\lambda<0$ is analogous).
  Let $\mathcal{BC}_{\mathcal{FF}}(\lambda)\in  \cat{PerfAnStk}^\delta$ be the  sheaf  on  perfect $\delta$-spaces  sending $(S,X)$ to $\Gamma(X/\varphi^\Z, \s{O}_{X/\varphi^\Z}(\lambda))$.
  Since Fargues--Fontaine curves of the form $(\Marc(A), \FF_A)$ for a nilperfectoid $\Q_p$-algebra $A$ form a basis for the $!$-topology on Fargues--Fontaine curves, locally in the $!$-topology the spectral valued sheaf $\mathcal{BC}_{\mathcal{FF}}(\lambda)$ is actually connective.

  We can deduce (using \Cref{sec:perf-compl-farg-1-perfect-complexes-on-ff-curves}) that one has a natural isomorphism of Gelfand stacks over $\Q_p$
  \[
  \BC_{\mathcal{FF}}(\s{O}(\lambda))_{\real}= \BC(\s{O}(\lambda)). 
  \]
\end{remark}

\begin{example}
  \label{sec:perf-compl-farg-2-example-bc-o-1}
  Let $\mathbb{D}_{\Q_p}$ be the open unit disc over $\Q_p$.
  We write
  \[
    1+\widetilde{\mathbb{D}}^{<1}_{\Q_p}:=\varprojlim_{x\mapsto x^p}(1+\mathbb{D}^{<1}_{\Q_p}),
  \]
  where the limit is computed in Gelfand stacks, i.e., at the level of functions on affinoids is an \textit{uncompleted} filtered colimit.
  
  Using \cite[Proposition II.2.2]{fargues2021geometrization} we see that there is a natural isomorphism 
  \[
   1+{\widetilde{\mathbb{D}}_{\Q_p}}^{\dR}\xrightarrow{\sim} \mathcal{BC}(\s{O}(1))^\dR
  \]
  such that for a $A$ nilperfectoid $\Q_p$-algebra, the composition with the natural map $\mathcal{BC}(\s{O}(1))^\dR(A)\to \GSpec(A^u)$ is given by the composition of the logarithm $1+A^{u,<1}\to A^u$ and the projection $\varprojlim_{x \mapsto x^p} (1+A^{u,<1})\to  1+A^{u,<1}$ to the first component.
  We note that the logarithm $(1+A^{u,<1})\to A^u$ lifts naturally to a map\footnote{E.g., the natural transformation $1+\mathbb{D}_{\Q_p}^{<1}\to \A^1_{\Q_p}$ of Gelfand stacks can be defined via the usual series of the logarithm evaluated on \emph{static} Gelfand rings, and then a posteriori evaluated on \emph{animated} Gelfand rings.}
  \[
    \log\colon 1+A^{<1}\to A.
  \]
  From  the pullback diagram of \Cref{sec:perf-compl-farg-1-calculation-for-hom} we get a natural map
  \[
   1+\widetilde{\mathbb{D}}_{\Q_p}^{<1}\to \mathcal{BC}(\s{O}(1)),
  \]
  which has to be an isomorphism (even of abelian group objects in Gelfand stacks). Indeed, this follows from the fact that the logarithm map  $\log\colon 1+\bb{D}_{\bb{Q}_p}^{<1} \to \bb{A}^1_{\bb{Q}_p}$ is  an \'etale map of partially proper rigid spaces, and that the map $ 1+\widetilde{\mathbb{D}}_{\Q_p}^{<1}\to  1+\mathbb{D}_{\Q_p}^{<1}$ is $\dagger$-formally \'etale being a limit of finite \'etale maps. Thus, we have a cartesian diagram
  \[
  \begin{tikzcd}
   1+\widetilde{\mathbb{D}}_{\Q_p}^{<1}  \ar[r] \ar[d,"\log"]& ( 1+\widetilde{\mathbb{D}}_{\Q_p}^{<1})^{\dR} \ar[d,"\log"] \\ 
   \bb{A}^1_{\bb{Q}_p} \ar[r] &  \bb{A}^{1,\dR}_{\bb{Q}_p}. 
  \end{tikzcd}
  \]
\end{example}

In general, absolute Banach--Colmez spaces have the following form.

\begin{proposition}\label{Prop-Description-BC-Spaces}
The following hold:
\begin{enumerate}

\item Let $\lambda\geq 0$ be a rational number. Then there is a pullback square of Gelfand stacks over $\Q_p$
\[
\begin{tikzcd}
\BC(\s{O}(\lambda)) \ar[r]\ar[d] & \BC(\s{O}(\lambda))^{\dR} \ar[d] \\ 
\mathbb{G}_a^{\an}\otimes_{\Q_p} D_{-\lambda} \ar[r] & \mathbb{G}_a^{\an,\dR}\otimes_{\Q_p} D_{-\lambda}
\end{tikzcd}
\]
In particular, we have $\BC(\s{O})=\Q_p^{\la}$, and for any $\lambda\geq 0$ we have an extension of $\Q_p^{\la}$-modules
\[
\mathbb{G}_{a}^{\dagger}\otimes_{\Q_p} D_{-\lambda}\to \BC(\s{O}(\lambda))\to \BC(\s{O}(\lambda))^{\dR}. 
\]

\item  Let $\lambda<0$ be a rational number.
  Then there is a pullback square of Gelfand stacks over $\Q_p$ 
\[
\begin{tikzcd}
\BC(\s{O}(\lambda)[1]) \ar[r] \ar[d] & \BC(\s{O}(\lambda)[1])^{\dR} \ar[d] \\ 
B(\mathbb{G}_a^{\dR}\otimes_{\Q_p} D_{-\lambda}) \ar[r]  & B(\mathbb{G}_a^{\an,\dR}\otimes_{\Q_p} D_{-\lambda}) .
\end{tikzcd}
\]
In particular, we have an extension of $\Q_p^{\la}$-modules 
\[
B(\mathbb{G}_a^{\dagger}\otimes_{\Q_p} D_{-\lambda})  \to \BC(\s{O}(\lambda)[1])\to \BC(\s{O}(\lambda)[1])^{\dR}. 
\]
\end{enumerate}
\end{proposition}
\begin{proof}
Let $A$ be a nilperfectoid ring and let $\mathcal{E}$ denote either $\s{O}(\lambda)$ for $\lambda\geq 0$ or $\s{O}(\lambda)[1]$ for $\lambda<0$. By \Cref{sec:perf-compl-farg-1-calculation-for-hom} we  have a pullback diagram of solid $\Q_p$-modules
\[
\begin{tikzcd}
\BC(\mathcal{E}) \ar[r] \ar[d] & \Gamma(\FF_{A^u}, \mathcal{E}_{\FF_{A^u}}) \ar[d] \\ 
\Gamma(\GSpec(A), \mathcal{E}_A) \ar[r] & \Gamma(\GSpec(A^u), \mathcal{E}_{A^u}). 
\end{tikzcd}
\]
The identification of the lower left term with $\mathbb{G}_a^{\an}\otimes D_{-\lambda}$ for $\lambda\geq 0$ (resp.\ for  $B(\mathbb{G}_a^{\an}\otimes D_{-\lambda})$ for $\lambda<0$) is clear from the definition of $\s{O}(\lambda)$.
It is also clear that the lower right term is the de Rham stack of the lower left term.
Finally, it is also clear from the definition that the upper right term is the de Rham stack of the upper left term, and that it agrees with the de Rham stack of the classical absolute  Banach--Colmez spaces of \cite[Chapter II]{fargues2021geometrization} (after base change to $\Marc(\Q_p)$). 
\end{proof}

We can note the following consequence.

\begin{corollary}\label{Corollary-Epi-BC-Spaces}
  Let $\lambda\geq 0$ be a rational number.
  Then the natural map
  $$
  \BC_{\mathcal{FF}}(\s{O}(\lambda))\to \BC_{\mathcal{FF}}(\s{O}(\lambda))^{\dR}
  $$
   in $\cat{PerfAnStk}^\delta$ is an epimorphism (where de Rham stacks are taken as in \Cref{sec:hyodo-kato-stacks-2-de-rham-stack-for-ff-stack}).
  The analogous  statement holds for $\BC_{\mathcal{FF}}(\s{O}(\lambda)[1])$ when $\lambda<0$.
\end{corollary}
\begin{proof}
By \Cref{sec:real-farg-font-1-realization-is-conservative}, it suffices to prove that this map is an epimorphism after applying the realization functor, in this case the claim follows from \Cref{Prop-Description-BC-Spaces} as the map $\mathbb{G}^{\an}_{a}\to \mathbb{G}_{a}^{\an,\dR}$ is an epimorphism of Gelfand stacks. 
\end{proof}

\newpage
\section{The rational analytic prismatization}
\label{Section:rational_analytic_prismatization}

\subsection{Main definitions}\label{main def}
We recall from \Cref{sec:gener-cart divis} that the Gelfand stack $\A^1/\Gm$ parametrizes Cartier divisors on Gelfand rings, i.e., one has the natural equivalence
\[
  \mathrm{Maps}(X,\A^1/\Gm)\cong \{ \textrm{ Cartier divisors } (s\colon \mathcal{L}\to \s{O}_X) \textrm{ on $X$ }\}.
\]
Similarly, we see that $\A^{1,\dR}/\Gm^\dR\cong (\A^{1}/\Gm)^\dR$ is equivalent to the functor
\[
  A\mapsto \{ \textrm{ Cartier divisors } (s\colon L\to A^u) \textrm{ on $\GSpec(A^u)$ }\}
\]
on nilperfectoid $\Q_p$-algebras $A$.

\begin{definition}\label{Definition-Analytic-prismatic-site}
  Let $X$ be a Gelfand stack.
  \begin{enumerate}
  \item  If $(S,\mathcal{X})$ is a $\delta$-space, a Cartier divisor $(s\colon \mathcal{L}\to \s{O}_{\mathcal{X}})$ on $\mathcal{X}$ is of degree $1$ if its pullback to $\FFdisc{S}^\times$ along the natural map $\FFdisc{S}^\times\to \mathcal{X}$ is of degree $1$ in the sense of Fargues--Scholze, i.e. if it represents the closed Cartier divisor $S^\sharp\to \FFdisc{S}^\times$ of an untilt $S^\sharp$ of $S$ (cf.\ \cite[Proposition II.1.18]{fargues2021geometrization}). A \textit{rational analytic prism} $(S,\mathcal{X},D)$ is a  $\delta$-space $(S,\mathcal{X})$ together with a Cartier divisor $D\subseteq \mathcal{X}$ of degree $1$.

  \item The \textit{rational analytic prismatic site of $X$}, denoted as $X_{\prism-\site}$, is the category of tuples $(S,\mathcal{X}, D,\alpha)$ where $(S,\mathcal{X},D)$ is a rational analytic prism and $\alpha\colon D\to X$ is a morphism of Gelfand stacks.
  Coverings are induced by $!$-covers in $(S,\mathcal{\mathcal{X}})$.
   \item Let $\Q_{p,\prism}$ be the analytic  $\delta$-stack sending $(S,\mathcal{X})$ to the anima of degree $1$ Cartier divisors on $(S,\mathcal{X})$.
  \item Let $X_{\prism}$ be the analytic $\delta$-stack over $\Q_{p,\prism}$ mapping $(S,\mathcal{X}, D)$ to the anima of maps $D\to X$.
 \item The \emph{rational analytic prismatization $X^\prism$ of $X$} is the realization of the analytic $\delta$-stack $X_{\prism}$.    Equivalently, it is the sheafification of the functor sending a nilperfectoid $\Q_p$-algebra $A$ over $\Q_p^\prism$ to the anima of maps $D\to X$, where $D\subseteq \mathcal{X}_{A}$ is the degree $1$ divisor classified by the map $\GSpec(A)\to \Q_p^\prism$.
    \item The \emph{rational perfect analytic prismatization of $X$} is the Gelfand stack
    \[
      X^{\widetilde{\Div^1}}:=\varprojlim\limits_{\varphi}X^\prism.
    \]
  \item The \emph{rational perfect analytic prismatization of $X$ modulo Frobenius}, or alternatively the \textit{$p$-adic twistor stack of $X$}, is the Gelfand stack
    \[
      X^{\Div^1}:=X^{\widetilde{\Div^1}}/\varphi^\Z.
    \]
  \end{enumerate}  
\end{definition}

As before, we will drop the adjective ``rational'' in what follows.
Some remarks are in order.

\begin{remark}
  \label{sec:rati-analyt-prisms-2-remarks-definition-prism}
  \begin{enumerate}
    \item Of course, \Cref{Definition-Analytic-prismatic-site} resembles the definition of a prism in \cite{bhatt_scholze_prisms_and_prismatic_cohomology} because locally an untilt of $S$ is defined by an element of degree $1$ in the respective Witt vectors.

    \item    We have a natural colimit preserving functor
\[
G\colon \cat{AnStk}^\delta_{/ \Q_{p,\prism}}\to \GelfStk,\ (S,\mathcal{X},D)\mapsto D
\]
given by the left Kan extension of the functor sending an analytic prism $(S,\mathcal{X},D)$ with underlying degree $1$ divisor $D\to \mathcal{X}$ to the Berkovich space $D$ seen as a Gelfand stack.
The functor $G$ has a right adjoint sending a Gelfand stack $X$ to the analytic $\delta$-stack $X_{\prism}$ over $\Q_{p,\prism}$ mapping $(S,\mathcal{X}, D)$ to the anima of maps $D\to X$ (up to inverting $!$-equivalences).
We note that the analytic prismatic topos of $X$ is $\cat{AnStk}^\delta_{/X_{\prism}}$.

It is clear from this construction that the functor $X\mapsto X_{\prism}$ commutes with fiber products. 
    \item Prism structures satisfy rigidity: if $(S',\mathcal{X}',D')\to (S,\mathcal{X},D)$ is a morphism of rational analytic prisms (i.e., one is given a factorization of $D'\to \mathcal{X}'\to \mathcal{X}$ over $D$), then the map $D'\to D\times_{\mathcal{X}} \mathcal{X}'$ is an isomorphism.
      Indeed, one needs to see that one has an isomorphism on invertible ideal sheaves, but this can be done on the uniform completions, where it reduces to the classical rigidity for degree $1$ divisors on $\FFdisc{S'}^\times$.
      
      \item Let $X$ be a Gelfand stack, one can  define  the \textit{perfect  analytic prismatic site} $X_{\widetilde{\cDiv^1}-\ob{site}}$ by restricting the prismation site to perfect analytic prisms (that is, analytic prisms with underlying perfect $\delta$-space). Equivalently, one can define the perfect analytic $\delta$-stack $X_{\widetilde{\cDiv^1}}\in \cat{PerfStk}^{\delta}$ as in \Cref{Definition-Analytic-prismatic-site} (4), we set $X_{\cDiv^1}=X_{\widetilde{\cDiv^1}}/\varphi$.  By construction, one has that $X_{\widetilde{\cDiv^1}}=X_{\prism}^{\perf}$ and therefore a  natural functor of Gelfand stacks 
      \begin{equation}\label{eq-MapRealizationsPerfStacksDiv1}
      (X_{\widetilde{\cDiv^1}})\to X^{\Div^1}.
      \end{equation}
      By \Cref{sec:real-farg-font-2-properties-of-admissible-delta-stacks} and \Cref{LemmaDaggerCohesivePrismatization} down below, the map \eqref{eq-MapRealizationsPerfStacksDiv1} is an equivalence if $X$ is $\dagger$-cohesive. 
  \end{enumerate}
\end{remark}

\begin{remark}\label{RemarkInscriptionstwo}
We continue the discussion of \Cref{RemarkInscriptionsHowe}. Compared to the theory of inscriptions of Howe, the functor sending $X\mapsto X_{\cDiv^1}$ is the analogue in our setting of the functor $X\mapsto X^{\diamond, \ob{lf}}$ of \cite[\S 4.4]{howe2025inscriptiontwistorspadicperiods} for the pair appearing in Proposition 4.3.1 (1) with $i=1$ in \textit{loc.\ cit.}.\footnote{Indeed, the other variants of the functor $X\mapsto X^{\diamond, \ob{lf}}$  of Howe is the analogue of the $X\mapsto X_{\cDiv^1}$ where instead of degree $1$ divisors one considers degree $d$ divisors, quotients by powers of the divisor, the completion at the divisor, the complement of the divisor, etc.}
\end{remark}

For any  Gelfand stack $X$, there is a natural morphism 
$$
X^\prism \to \A^1/\Gm
$$
coming from pullback of Cartier divisors along $\{0\}\times T \to \N \times T \to \mathcal{X}_T$ for $T=\GSpec(A)$ nilperfectoid. This allows us to single out an interesting substack of $X^\prism$.

\begin{definition}
Let $X$ be a  Gelfand stack.  The \emph{analytic Hodge-Tate of $X$} is the Gelfand stack
    \[
      X^\HT:=X^\prism\times_{\A^1/\Gm} B\Gm.
    \]
\end{definition}

\begin{remark}
  \label{sec:rati-analyt-prisms-1-warning-about-distinguished-locus}
  We warn the reader that the natural morphism $\{0\}\times T\to \mathcal{X}_T$ is \emph{not} a Cartier divisor.
  If $D\in \Q_p^{\HT}(A)$, then by definition one is given a morphism $T\to D$, which is a derived thickening (indeed, one has fiber sequence $\Nil^\dagger(A)[1]\to \s{O}_D\to A$).
  For a general $D\in \Q_p^\prism(A)$ one can consider the Cartier divisor $T\cap D$ on $T$.
  We have that $D$ factors over the open complement of $\Q_p^\prism\setminus \Q_p^\HT$ of $\Q_p^\HT$ if and only if $T\cap D=\emptyset$.
\end{remark}

\begin{remark}\label{rk:explicit-description}
Let $(S,\mathcal{X},D)$ be an analytic prism. By definition of $\Q_{p,\prism}$, the prismatic structure defined by $D$ on the  $\delta$-space $(S,\mathcal{X})$ corresponds to a morphism of analytic $\delta$-stacks
$$
(S,\mathcal{X}) \to \Q_{p,\prism}
$$
which factors naturally through $(S,\mathcal{X}) \to D_{\prism}\to \Q_{p,\prism}$, where the left arrow corresponds to the identity of $D$.  
Passing to realizations, we get  morphisms of Gelfand stacks 
$$
\mathcal{X} \to D^{\prism}\to  \Q_p^{\prism}. 
$$
This is similar to the morphism $\Spf(A) \to \Z_p^\prism$ in the $p$-complete theory, attached to a prism structure on the $\delta$-ring $A$. If $\mathcal{X}$ is perfect, the morphism factors through $D^{\widetilde{\Div^1}}\to \bb{Q}_p^{\widetilde{\Div^1}}$.
\end{remark}

\begin{example}\label{ExamplesRelativePrismatizations}
\begin{enumerate}
\item Let $\mathcal{X}_{\rm cyc}=\big((\varprojlim_{x\mapsto x^p} 1+{\bb{D}}_{\Q_p}^{<1})\backslash \{1\} \big)$. Then $$(S=\mathcal{M}(\F_p((t^{1/p^\infty}))), \mathcal{X}_{\rm cyc})$$ defines a (perfect) $\delta$-space, called the \textit{cyclotomic $\delta$-space}.  Then, the divisor $(\frac{q^p-1}{q-1})$ gives rise to a degree $1$ divisor on $\mathcal{X}_{\rm cyc}$ which is the natural map $\GSpec (\Q_p(\zeta_{p^{\infty}}))\to \mathcal{X}$. The resulting prism, called the \textit{cyclotomic prism}, is an analogue of the $q$-de Rham prism. This produces a cartesian diagram of Gelfand stacks 
\[
\begin{tikzcd}
\GSpec(\Q_p(\zeta_{p^{\infty}}))  \ar[r] \ar[d] & \mathcal{X}_{\rm cyc} \ar[d] \\ 
\Q_p^{\HT}\ar[r] & \Q_p^\prism.
\end{tikzcd}
\]
\item Let $\bb{Q}_p^{\ob{Kum}}$ be the completion of $\bb{Q}_p(p^{1/p^{\infty}})$, and let  $\mathcal{X}=\big( \varprojlim_{x\mapsto x^p} {\bb{D}}_{\Q_p}^{<1} \backslash \{0\} \big)$  be the decompletion of $\FF_{\Q_p^{\ob{Kum}}}$  obtained by the Breuil--Kisin $\delta$-ring $\Z_p[[X]]$ with $\varphi(X)=X^p$. The divisor $(X-p)$ on $\Z_p[[X]]$ gives rise to a degree $1$ divisor on $\mathcal{X}$ which is the natural map $\GSpec (\Q_p(p^{1/p^{\infty}}))\to \mathcal{X}$. This produces a cartesian diagram of Gelfand stacks 
\[
\begin{tikzcd}
\GSpec (\Q_p(p^{1/p^{\infty}}))  \ar[r] \ar[d] &  \varprojlim_{x\mapsto x^p} {\bb{D}}_{\Q_p}^{<1}\backslash \{0\}  \ar[d] \\ 
\Q_p^{\HT}\ar[r] & \Q_p^\prism.
\end{tikzcd}
\]

\end{enumerate}
\end{example}

As in the $p$-complete theory, one can define relative versions of the above definitions. 

\begin{definition}[Relative prismatization]\label{RelativePrismatization}
  Let $(S,\mathcal{X},D)$ be an analytic prism and consider the natural map $\mathcal{X}\to \Q_p^\prism$ arising from the divisor $D$.
  Let $Y$ be a  Gelfand stack over $D$. 
\begin{enumerate}

\item  The \textit{relative analytic  prismatization of $Y$ over $(S,\mathcal{X},D)$} is the pullback map 
\[
Y^{\prism/\mathcal{X}}:= Y^{\prism}\times_{D^{\prism}} \mathcal{X},
\] 
where the map $\mathcal{X}\to D^{\prism}$ arises from the identity on $D$. Similarly, if $(S,\mathcal{X},D)$ is a perfect analytic prism, namely a perfect $\delta$-space endowed with a degree $1$ divisor, we define the \textit{relative perfect analytic  prismatization $Y^{\widetilde{\Div^1}/\mathcal{X}}$ of $Y$ over $(S,\mathcal{X},D)$}, resp.\ the \textit{relative perfect analytic prismatization $Y^{\Div^1/\mathcal{X}}$ of $Y$ over $(S,\mathcal{X},D)$ modulo Frobenius} (or the \textit{relative $p$-adic twistor stack over $(S,\n{X},D)$}).

\item The \textit{relative Hodge--Tate stack of $Y$ over $(S,\mathcal{X},D)$} is the pullback
\[
Y^{\HT/D}:= Y^{\prism/\mathcal{X}}\times_{\mathcal{X}} D. 
\]

\end{enumerate}
\end{definition}

\subsection{Geometry of the prismatization}
We start with the following easy statement describing, for a Gelfand stack $X$, the underlying arc-stack of $X_\prism$.

\begin{lemma}
  \label{sec:defin-first-prop-3-de-rham-stack-of-y-subscript-prism}
  Let $X$ be a Gelfand stack, and $X_\prism$ as in \Cref{Definition-Analytic-prismatic-site}.
  Then one has an isomorphism of arc-stacks
  $$
(X_\prism)^{\diamond_{\FF}} = X^\diamond
$$
and thus an isomorphism
  \[
    X_{\prism}^{\dR}=\FFdisc{X^\diamond}^{\times,\dR}
  \]
  of analytic $\delta$-stacks over $\Q_{p,\prism}^\dR$.
\end{lemma}
\begin{proof}
  This immediately follows by unravelling the definitions and \Cref{sec:relation-hyodo-kato-realization-commutes-with-de-rham-stacks} (3). 
\end{proof}

\begin{remark}\label{RemarkNaturalMapsXprism}
From \Cref{sec:defin-first-prop-3-de-rham-stack-of-y-subscript-prism} and the adjunctions of \Cref{sec:relation-hyodo-kato-realization-commutes-with-de-rham-stacks} (2) we have natural maps of analytic $\delta$-stacks
\[
  \FFdisc{X^{\diamond}}^{\times} \to  X_{\prism} \to    \FFdisc{X^{\diamond}}^{\times, \dR}
\]
inducing equivalences on de Rham stacks.  Passing to realizations, we have morphisms of Gelfand stacks
\[
\FFdisc{X^{\diamond}}^{\times} \to  X^{\prism} \to \FFdisc{X^{\diamond}}^{\times, \dR}
\]
also inducing equivalences on de Rham stacks. Hence, $X^{\prism}$ sits in between the relative Fargues-Fontaine disc over $X^{\diamond}$ (incarnated as a Gelfand stack) and its analytic de Rham stack. Note that in this diagram, the left and right terms only depend on the perfectoidization $X^{\diamond}$ of $X$ (its ``underlying topological space''), but $X^{\prism}$ itself does not.

Since the left term is perfect, we can in fact replace $X^\prism$ by $X^{\widetilde{\Div^1}}$ in this diagram.
\end{remark}

\begin{proposition}
  \label{sec:defin-first-prop-2-perfection-of-analytic-prismatization}
  Let $X$ be a  Gelfand stack.
  The natural maps $X^{\widetilde{\Div^1},\dR}\to X^{\prism,\dR}$ and $\varinjlim\limits_{\varphi} X^\prism\to X^{\prism,\dR}$ are isomorphisms, and $X^{\widetilde{\Div^1},\dR}\cong \FFdisc{X^{\diamond}}^{\times,\dR}$. 
  \end{proposition}
\begin{proof}
  The isomorphisms $X^{\widetilde{\Div^1},\dR} \cong X^{\prism,\dR}\cong \FFdisc{X}^{\times,\dR}$ are immediate by evaluating on nilperfectoids.
  For the rest of the first statement, we notice that $X^{\prism}$ is the sheafification for the $!$-topology of the functor $X^{\prism-\ob{pre}}$ sending a nilperfectoid $A$ with $T:=\GSpec(A)$ to the anima of degree $1$ divisors $D\to \n{X}_Y$ and maps $D\to X$. Since sheafification commutes with colimits, it suffices to show that the natural map of presheaves
  \[
  \varinjlim_{\varphi }X^{\prism-\ob{pre}}\to \FFdisc{X^{\diamond}}^{\dR}
  \]
  is an equivalence. By definition, the image of a point $x\in  X^{\prism,\ob{pre}}(A)$ along Frobenius is the datum of the degree $1$ divisor $\varphi^*D\to \n{X}_T$ and the map $\varphi^* D\to D\to X$. The claim  then follows from the fact that for each degree $1$ divisor $D\subseteq \mathcal{X}_T$ its pullback along $\varphi^n\colon \mathcal{X}_T\to \mathcal{X}_T$ for some $n\gg 0$ is contained in the locus $\mathcal{X}_T\setminus  (\N\times T)\cong \FFdisc{T^\flat}^\times\setminus (\N\times T^u)$. 
\end{proof}

\begin{proposition}\label{prop:ht-stack-of-perfectoids}
Let $Y$ be a perfectoid space over $\Q_p$. The natural morphism $\pi_Y: Y^\HT \to Y$ from \Cref{sec:defin-first-prop-2-first-properties-of-analytic-prismatization} (2) is an isomorphism.
\end{proposition}
\begin{proof}
  Let $A$ be a nilperfectoid ring over $\Q_p$.
  Let $\GSpec(A) \to Y^\HT$ be an $A$-point of $Y^\HT$, i.e., a degree $1$ divisor $D$ on $\FF_A$ with a factorization $\GSpec(A) \to D \to \FF_A$, and a map $D\to Y$.
  The morphism $\pi_Y$ sends this datum to the composition $f_D: \GSpec(A) \to D \to Y$.

  Conversely, let $f: \GSpec(A) \to Y$. We get an induced map on Fargues-Fontaine curves $\FF_A \to \FF_Y$ and thus a degree $1$ divisor $D_f = Y \times_{\FF_Y} \FF_A$ on $\FF_A$ with a morphism $\GSpec(A)\to D$ (note that $\FF_Y$ makes sense for $Y$ nilperfectoid, but the fact that $Y \to \FF_Y$ is a degree $1$ divisor makes use of the perfectoidness assumption on $Y$) with a morphism $D_f\to Y$, defining an $A$-point of $Y^\HT$.
  This yields a section $Y\to Y^\HT$ of $\pi_Y$.

For $D$ as above, we claim that $D$ is naturally isomorphic to $D_{f_D}$. By rigidity (\Cref{sec:rati-analyt-prisms-2-remarks-definition-prism}
(4)), it is enough to produce a morphism of prism structures $D\to D_{f_D}$. By assumption we have morphisms $D \to Y$ and $D\to \FF_A$, thus a morphism $D \to Y \times_{\GSpec(\Q_p)} \FF_A$. The composition of $D\to \FF_A$ with the morphism $\FF_A \to \FF_Y$ induced by $f_D$ on Fargues-Fontaine curves factors by construction through $\FF_D$. But the composition $D \to \FF_D \to \FF_Y$ identifies with $D\to Y \to \FF_D$, by the construction of nilperfectoid Fargues-Fontaine curves. Hence we get the desired factorization $D \to D_{f_D}$.
This finishes the proof.
\end{proof}

In particular, for a  Gelfand stack over a perfectoid base, its relative and absolute Hodge--Tate stacks are the same. We warn the reader that however, and in contrast with the $p$-complete prismatization, the stacks $Y^{\Div^1}$ and $\FF_Y$ are not isomorphic for $Y$ perfectoid.

Analytic prismatization behaves well with respect to \'etale, even $\dagger$-formally \'etale localization.
We note that in \Cref{subsec:desc-perf-analyt} we will prove stronger descent statements for the analytic prismatization.

  \begin{lemma}
    \label{sec:defin-first-prop-5-prismatization-for-dagger-formally-etale-maps}
   Let $X'\to X$ be a $\dagger$-formally \'etale map.
      Then the natural map
      \[
        (X')^\prism\to X^\prism\times_{\FFdisc{X^\diamond}^{\times, \dR}}\FFdisc{X^{\prime,\diamond}}^{\times,\dR}
      \]
      is an isomorphism.
      In particular, the functor  $X\mapsto X^\prism$ satisfies descent for $\dagger$-formally \'etale hypercoverings on  Gelfand stacks.
  \end{lemma}
  \begin{proof}
    Let $A$ be nilperfectoid, with $T:=\GSpec(A)$.
    Let $D\subseteq \mathcal{X}_T$ be a degree $1$ Cartier divisor.
Let $\overline{D}:= D_{|\FFdisc{T^\flat}^\times}$.
      We note that $\GSpec(\overline{D})\to \GSpec(D)$ is the $\dagger$-reduction.
    This implies by $\dagger$-formally \'etaleness of $X^\prime\to X$ that
    \[
      X^\prime(D)\cong X(D)\times_{X(\overline{D})} X'(\overline{D}),
    \]
 and then, after passing to sheafifications for the $!$-topology,  that $X^{\prime,\prism}\cong X^\prism\times_{\FFdisc{X^\diamond}^{\times,\dR}}\FFdisc{X^{\prime,\diamond}}^{\times,\dR}$
    as desired.

    The descent statement follows by base change along $X^\prism\to X^{\prism,\dR}\cong \FFdisc{X^\diamond}^{\times,\dR}$ and the arc-hyperdescent of the functor $\FFdisc{(-)}^{\times,\dR}$.
  \end{proof}

Let $A$ be a nilperfectoid ring with $T=\GSpec(A)$.
By construction, the $\delta$-space associated to $A$ is the pushout of Berkovich spaces
 \[
 \n{X}_T= \FFdisc{T^{\diamond}}\bigsqcup_{\bb{N}\times  T^{u}} (\bb{N} \times T).
 \]
Let $D\to \n{X}_A$ be a degree $1$ divisor, we have a pushout description  of divisors
\begin{equation}\label{eqPushoutDegree1divisor}
D=\overline{D}\bigsqcup_{D|_{(\bb{N}\times T^u)}} (D|_{\bb{N}\times T}),
\end{equation}
where the intersection refers to a fiber product over $\n{X}$. Thus, for a Gelfand stack $X$, we have a natural commutative diagram 
\begin{equation}\label{eqDiagramPrismatizationCones}
\begin{tikzcd}
X^{\prism} \ar[r] \ar[d] & \FFdisc{X^{\diamond}}^{\times,\dR} \ar[d] \\ 
\prod_{n\in \bb{N}} X^{\cone} \ar[r] &  \prod_{n\in \bb{N}}( X^{\cone})^{\dR},
\end{tikzcd}
\end{equation}
where the map $X^{\prism}\to \prod_{n\in \bb{N}} X^{\cone}$ arises by the restriction of $D\to X$ to each intersection with $\bb{N}\times T$. It is natural to ask under which conditions the diagram \eqref{eqDiagramPrismatizationCones} becomes cartesian; this is precisely what the $\dagger$-cohesive assumption guarantees:

\begin{lemma}\label{LemmaDaggerCohesivePrismatization}
Let $X$ be a $\dagger$-cohesive Gelfand stack. Then for $A$ a a nilperfectoid $\Q_p$-algebra with $Y=\GSpec (A)$ and $D\to \n{X}_T$ a degree $1$ divisor, the natural commutative diagram  of anima
\[
\begin{tikzcd}
X(D) \ar[r] \ar[d]& X(\overline{D})  \ar[d] \\ 
X(D|_{(\bb{N}\times T)})  \ar[r] & X(D|_{(\bb{N}\times T^u)}) 
\end{tikzcd}
\]
induced by \eqref{eqPushoutDegree1divisor} is cartesian.  In particular, $X_{\prism}$ is an admissible $\delta$-stack, and the diagram of Gelfand stacks \eqref{eqDiagramPrismatizationCones} is cartesian. 
\end{lemma}
\begin{proof}
 First, note that by quasi-compactness of $D$, only finitely many intersections with $\bb{N}\times T$ are non-empty.  The lemma  follows directly from the definition of $\dagger$-cohesivity and the fact that the pushout decomposition of \eqref{eqDiagramPrismatizationCones} is precisely as in  \Cref{LemmaPushoutsBerkovichSpaces} with surjective morphisms of rings. 
\end{proof}

Let us rewrite the description of \Cref{LemmaDaggerCohesivePrismatization} in a slightly nicer way. For this, we need one more definition. 

\begin{definition}[Breuil--Kisin ideal]
  \label{sec:breuil-kisin-twist-2-definition-bk-twist}
  Let $\mathcal{I}_\prism$ be the ideal sheaf of the Hodge--Tate locus. The \textit{Breuil--Kisin ideal} $\n{I}_{\BK}$ is the line bundle on $\Q_p^\prism$ defined by
  \[
    \bigotimes\limits_{n\in \N} (\varphi^n)^\ast(\mathcal{I}_\prism). 
  \]
We let 
\[
\mu_{\BK}\colon \bb{Q}^{\prism}\to \bb{A}^1/\bb{G}_m
\]
be the resulting map giving rise to the divisor $\mathrm{can}\colon \mathrm\n{I}_{\BK}\to \s{O}_{\bb{Q}^{\prism}}$ that we call the \textit{Breuil--Kisin divisor}. 
\end{definition}

\begin{remark}\label{RemarkProductMakesSense}
 \Cref{sec:breuil-kisin-twist-2-definition-bk-twist} has to be understood in the following way: let  $A$ be a nilperfectoid $\bb{Q}_p$-algebra with a morphism $f\colon T:=\GSpec(A)\to \Q_p^\prism$ given by a degree $1$ divisor $D$.  Let  $U\subseteq \mathcal{X}_T$ be a rational affinoid subspace in the $\delta$-space of $T$,  and consider the composite 
 \begin{equation}\label{eqProductDivisors}
 U\to  \mathcal{X}_T \to \prod_{n\in \bb{N}} \bb{A}^1/\bb{G}_m
 \end{equation}
 obtained by the pullback along $U\to \mathcal{X}_T $ of the Frobenius translates of the degree $1$ divisor $D$.
 The map \eqref{eqProductDivisors} has all but finitely many projections that factor through the trivial divisor $\GSpec (\bb{Q}_p) \xrightarrow{1} \bb{A}^1/\bb{G}_m$. Therefore, if  $(\prod_{n\in \bb{N}} \bb{A}^1/\bb{G}_m)^{\ob{res}}\subset \prod_{n\in \bb{N}} \bb{A}^1/\bb{G}_m$  denotes the locus where all but finitely many divisors are trivial, we have a factorization  $U\to (\prod_{n\in \bb{N}} \bb{A}^1/\bb{G}_m)^{\ob{res}}$ and hence of  $\mathcal{X}_T \to (\prod_{n\in \bb{N}} \bb{A}^1/\bb{G}_m)^{\ob{res}}$. 
 
Composing with $T\to \n{X}_T$, the previous implies that the map $\bb{Q}_p^{\prism}\to \prod_{n\in \bb{N}} \bb{A}^1/\bb{G}_m$ also factors through  $(\prod_{n\in \bb{N}} \bb{A}^1/\bb{G}_m)^{\ob{res}}$. Now, there is a well defined product map 
\[
(\prod_{n\in \bb{N}} \bb{A}^1/\bb{G}_m)^{\ob{res}} \xrightarrow{\otimes} \bb{A}^1/\bb{G}_m
\]
as only finitely many divisors are non  trivial. We deduce a product map
\[
\mu_{\BK}\colon \bb{Q}_p^{\prism}\to  \bb{A}^1/\bb{G}_m
\]
which is the Breuil--Kisin divisor 
\[
\n{I}_{\BK}= \bigotimes_{n\in \bb{N}} \varphi^{n,*} \n{I}_{\prism}. 
\]
 \end{remark}

\begin{remark}\label{RemarkBreuilKisinTwist}
The Breuil--Kisin ideal of  \Cref{RemarkProductMakesSense} will be isomorphic \textit{as a plain line bundle} to the analogue of the Breuil--Kisin twist $\s{O}\{1\}$ in our theory of the analytic prismatization.  There is a reason why we have not called the object $\n{I}_{\BK}$ the Breuil--Kisin twist, and this is due to the correct definition of the Frobenius. We will elaborate on this point later in \Cref{sec:analyt-prism-ga}. 

However, we  highlight that in the analytic prismatization (with $\Q_p$-coefficients!) the intricate arguments in \cite{bhatt2022absolute} or \cite{drinfeld20211}  for the existence of the Breuil-Kisin twist are not necessary. The reason for this simplification in the rational theory is the simple observation   that the infinite tensor product of the ideals  $\varphi^{n,*}\n{I}_{\prism}$ makes perfect sense on rational Fargues--Fontaine curves, but that its extension to $p=0$ requires a delicate control of the convergence of this tensor product. In particular, constructing the Breuil--Kisin twist in the theory of the \textit{integral analytic prismatization} will be equally difficult as it is in the prismatization of formal schemes. 
\end{remark}

  Using the Breuil--Kisin divisor, one can rewrite the  cartesian diagram \eqref{eqDiagramPrismatizationCones} in a more pleasant way.   We need the following easy lemma.

  \begin{lemma}\label{LemmaObservationBKIdeal}
  Let $A$ be a nilperfectoid ring with $T=\GSpec(A)$. Let $D\to \n{X}_T$ be a degree $1$ divisor and $D_{\BK}\to \n{X}_T$ the Breuil--Kisin divisor defined by the map $\mu_{\BK}$, as follows: $D$ defines a map $T\to \Q_p^\prism$ which lifts to a map $\mathcal{X}_T \to \Q_p^\prism$; its composition with $\mu_\BK$ gives the Cartier divisor $D_\BK$ on $\mathcal{X}_T$. Then there is a natural isomorphism of stacks over  $\n{X}_T$
\[
D_{\BK} = \bigsqcup_{n\in \bb{N}} \varphi^{n,*}D.
\]
  \end{lemma}

  \begin{proof}
  We first observe that given  $n\neq m$ the divisors $\varphi^{*,n}D$ and $\varphi^{*,m}D$ are disjoint in $\n{X}_{T}$, that is, $\varphi^{*,n}D\times_{\n{X}_T} \varphi^{*,m}D=\emptyset$. Indeed, this follows from the simple observation that their support in the associated Berkovich space $|\n{X}_T|= |\FFdisc{B^u}^{\times}|$ are  disjoint due to the compatibility of the radius map from \Cref{sec:farg-font-curv} with the Frobenius. 
  
  Then the statement of the lemma follows from the Frobenius-equivariance of the morphism $\mathcal{X}_T\to \Q_p^\prism$ (induced by $D$) and the fact that under this morphism, the divisor $\Q_p^\HT$ pulls back to $D$.
    \end{proof}

  Note first that, by construction of the Breuil--Kisin divisor, there exists for each $n$ a morphism of Cartier divisors $( \mathrm\n{I}_{\BK} \xto{\mathrm{can}} \s{O}_{\Q_p^\prism}) \to ((\varphi^n)^\ast(\mathcal{I}_\prism) \to \s{O}_{\Q_p^\prism})$, leading for any Gelfand stack $X$ to a morphism 
  $$
  X^{\rm cone} \times_{\A^1/\Gm,\mu_\BK} \Q_p^\prism \to \prod_{n\in \N} X^{\cone}.
  $$
Moreover, given a Gelfand stack $X$, the divisor $\mu_{\BK}$ of \Cref{sec:breuil-kisin-twist-2-definition-bk-twist} extends to a natural morphism 
\[
\mu_{\BK}\colon X^{\prism}\to X^{\cone}\times_{\A^1/\Gm,\mu_{\BK}} \Q_p^{\prism}
\]
constructed as follows: let $T=\GSpec(A)$ be a nilperfectoid ring with Fargues--Fontaine curve $\n{X}_A$. Let $D\to \n{X}_A$ be a degree $1$-divisor, then we have a natural maps of Gelfand stacks 
\[
D_{\BK}\to D
\]
whose restriction to $\varphi^{n,*}D$ is the natural Frobenius map $\varphi^{n}\colon \varphi^{n,*}D\to D$. Now, taking a map $D\to X$ corresponding to a point $\GSpec(A)\to X^{\prism}$, and restricting to the map
\[
D_{\BK}\times_{\n{X}_A} \GSpec(A)\to D_{\BK}\to D\to X , 
\]
we obtain the desired point $\GSpec(A)\to X^{\cone}\times_{\A^1/\Gm,\mu_{\BK}} \Q_p^{\prism}$. 

\begin{proposition}\label{prop:pullback-description-twistor-stack1}
Let $X$ be a $\dagger$-cohesive Gelfand stack. Then the natural commutative diagram 
\begin{equation}\label{eqBKIdealPullback}
\begin{tikzcd}
X^{\prism} \ar[d,"{\mu_{\BK}}"] \ar[r] & {\FFdisc{X^{\diamond}}^{\times,\dR}}  \ar[d, "{\mu_{\BK}^{\dR}}"] \\ 
X^{\cone} \ar[r] &  (X^{\cone})^{\dR}
\end{tikzcd}
\end{equation}
is cartesian.  In particular, $\mu_{\BK}\colon X^{\prism}\to  X^{\cone}$ is $\dagger$-formally \'etale. 
\end{proposition}

\begin{proof} 
Let $A$ be a nilperfectoid ring with $T=\GSpec(A)$, $T^u=\GSpec (A^u)$. Let $D\to \n{X}_T$ be a degree $1$ divisor.  The pushout diagram in Berkovich spaces 
\[
\n{X}_{T}= \FFdisc{A^u}^{\times}\bigsqcup_{\bb{N}\times T^u}  \bb{N}\times T
\]
induces a pushout diagram 
\[
D = \overline{D}\bigsqcup_{D\times_{\n{X}_T}   (\bb{N}\times T^u)} D\times_{\n{X}_T}   (\bb{N}\times T)
\]
with  $\overline{D}$ the reduction of the divisor $D$.

On the other hand,  the pullback of $D$ along the $n$-th component of  $\bb{N}\times T\to \n{X}_{T}$ is precisely  the divisor $\varphi^{n,*}D\times_{\n{X}_T} T$, namely, the Frobenius morphism induces an isomorphism $\varphi^{n}\colon 0\times T \xrightarrow{\sim} n\times T$. This produces a natural identification 
\[
D\times_{\n{X}_T}   (\bb{N}\times T) \cong \bigsqcup_{n\in \bb{N}} \varphi^{n,*}D\times_{\n{X}} T
\]
and similarly for $T^u$. By \Cref{LemmaObservationBKIdeal}  we  know that  $\bigsqcup_{n\in \bb{N}} \varphi^{n,*}D\times_{\n{X}} T = D_{\BK}\times_{\n{X}_T} T$.

 In total, we get a pushout diagram of Berkovich spaces
\[
D= \overline{D} \bigsqcup_{D_{\BK}\times_{\n{X}_T} T^u} D_{\BK}\times_{\n{X}_T} T.
\]
as $X$ is $\dagger$-cohesive, we get a cartesian diagram 
\[
\begin{tikzcd}
X(D) \ar[r] \ar[d] & X(\overline{D}) \ar[d] \\ 
X(D_{\BK}\times_{\n{X}_T} T) \ar[r] & X(D_{\BK}\times_{\n{X}_T} T^u).
\end{tikzcd}
\]
This is precisely the value at $A$ of the commutative diagram \eqref{eqBKIdealPullback}, finishing the proof of the proposition. 
\end{proof}

\begin{remark} \label{rk:radius-map}
As another consequence of the morphism $X^{\prism}\to (\FFdisc{X^{\diamond}}^\times)^\dR$ and the radius map from \Cref{sec:farg-font-curv}, we obtain a radius map
$$
\kappa: X^\prism \to (0,\infty).
$$
\end{remark}

Using the radius map, we can analyze the geometry of $X^\prism$ (cf. Figure 2).  In particular, the preimage of $(0,\infty)\backslash \{\cdots, p^{-2},p^{-1},1\}$ is precisely the same region on the Fargues--Fontaine de Rham stack of $X$, and the $n$-th Frobenius inverse of the Hodge--Tate stack maps to $p^{-n}$. Note however, that $\kappa^{-1}(p^{-n})$ is much larger than $X^{\HT}$.

  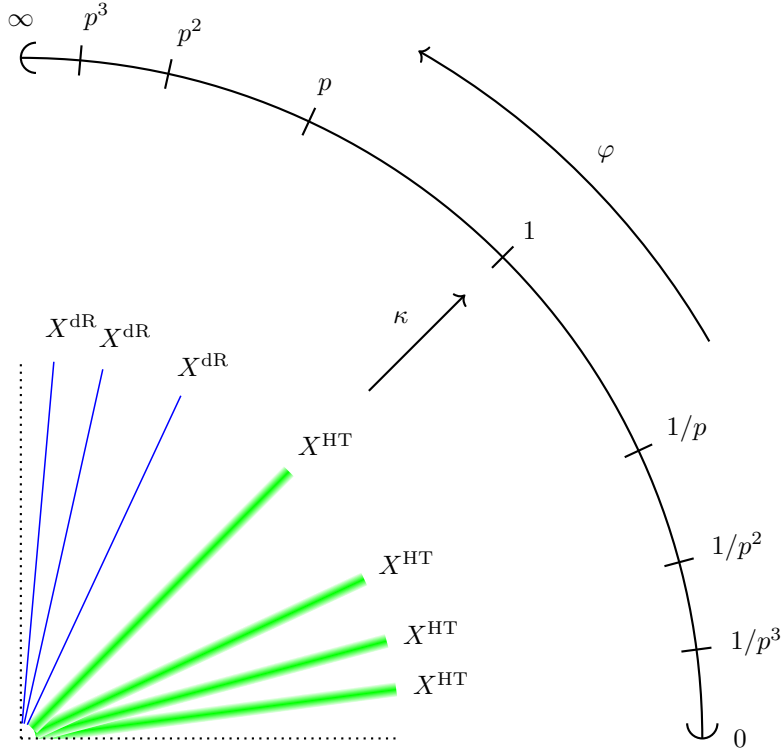
\begin{figure}[h!t]  
\begin{center}
\begin{tikzpicture}
  \def\r{9}
  \def\shortr{5} % Slightly shorter radius
  \def\longr{9.5}
  \def\vlongr{10.5}

  % Draw the figure
  \draw[dotted, thick] (0,0) -- (0,\shortr); % Dotted vertical radius
  \draw[thick] (\r,0) arc[start angle=0, end angle=90, radius=\r];
  
  % Draw the arrow at the center of the arc, slightly outwards
  \draw[->, thick] (0.866*\vlongr, 1/2*\vlongr) arc[start angle=30, end angle=60, radius=\vlongr];
  
  % Label for the angle
  \node at (1.4142/2*\vlongr + 0.3, 1.4142/2*\vlongr + 0.3) {$\varphi$};

  \def\shorttr{\r-0.2}
  \def\longgr{\r+0.2}

  \draw[thick] (\shorttr,0.2) arc (-180:0:0.2);
  \node at (\longr, 0) {$0$};
  
  \draw[thick] (7.5:\shorttr) -- (7.5:\longgr);
  \node at (7.5:\longr+0.3) {$1/p^3$};

  \draw[thick] (15:\shorttr) -- (15:\longgr);
  \node at (15:\longr+0.3) {$1/p^2$};

  \draw[thick] (25:\shorttr) -- (25:\longgr);
  \node at (25:\longr+0.2) {$1/p$};

  \draw[thick] (45:\shorttr) -- (45:\longgr);
  \node at (45:\longr) {$1$};

  \draw[thick] (65:\shorttr) -- (65:\longgr);
  \node at (65:\longr) {$p$};

  \draw[thick] (77.5:\shorttr) -- (77.5:\longgr);
  \node at (76.5:\longr+0.1) {$p^2$};
  
  \draw[thick] (85:\shorttr) -- (85:\longgr);
  \node at (84:\longr+0.1) {$p^3$};

  \draw[thick] (0.2, \longgr) arc (90:270:0.2);
  \node at (0, \longr) {$\infty$};

  \begin{scope}

  \clip (0,0) rectangle (\longr,\longr);
  
  % Draw rotated lines with the middle part darker (last one drawn should be on top)
  \begin{scope}[transform canvas={rotate=7.5}]
        \shade[shading=axis, bottom color=white, top color=white, middle color=green, shading angle=0] (0.2, -0.1) rectangle (\shortr, 0.1);
  \end{scope} 

  % Second layer, drawn on top of the previous one
  \begin{scope}[transform canvas={rotate=15}]
        \shade[shading=axis, bottom color=white, top color=white, middle color=green, shading angle=0] (0.2, -0.1) rectangle (\shortr, 0.1);
  \end{scope} 

  % Third layer, drawn on top of the previous one
  \begin{scope}[transform canvas={rotate=25}]
        \shade[shading=axis, bottom color=white, top color=white, middle color=green, shading angle=0] (0.2, -0.1) rectangle (\shortr, 0.1);
  \end{scope}

  % Fourth layer, drawn on top of the previous one
  \begin{scope}[transform canvas={rotate=45}]
        \shade[shading=axis, bottom color=white, top color=white, middle color=green, shading angle=0] (0.2, -0.1) rectangle (\shortr, 0.1);
  \end{scope}

  \draw[line width=0.2mm, color=blue] (0.4266*0.2,0.9063*0.2) -- (0.4226*\shortr, 0.9063*\shortr);

  \draw[line width=0.2mm, color=blue] (0.21644*0.2,0.9763*0.2) -- (0.21644*\shortr, 0.9763*\shortr);

  \draw[line width=0.2mm, color=blue] (0.087156*0.2,0.9962*0.2) -- (0.087156*\shortr, 0.9962*\shortr);

  \end{scope}
  
  \def\shortishr{\shortr+0.5}
  
  \node at (7.5:\shortishr+0.1) {$X^\HT$};
  
  \node at (14.5:\shortishr+0.1) {$X^\HT$};
  
  \node at (24.5:\shortishr+0.1) {$X^\HT$};
  
  \node at (44:\shortishr+0.1) {$X^\HT$};
  
  \node at (64:\shortishr) {$X^\dR$};
  
  \node at (75.5:\shortishr) {$X^\dR$};
  
  \node at (83:\shortishr) {$X^\dR$};
  
  \draw[thick, dotted] (0,0) -- (\shortr,0);
  
  \draw[->, thick] (45:\shortishr+1) -- (45:\shorttr-0.5);
  \node at (48:7.5) {$\kappa$};

\end{tikzpicture}
\end{center}
  \caption{\footnotesize{A picture of $X^\prism$, stolen from \cite{hauck}.}}  
\end{figure}

\begin{proposition}\label{PropositionDescriptionPrismatizationHKDivtilde}
Let $X$ be a $\dagger$-cohesive Gelfand stack.
\begin{enumerate}
\item The map $X^{\prism} \to (\FFdisc{X^{\diamond}}^\times)^\dR$ constructed in \Cref{sec:defin-first-prop-2-perfection-of-analytic-prismatization} is an isomorphism over $(1,\infty)$ under $\kappa$. Thus, the Frobenius induces an isomorphism $\varphi: X_{(1,\infty)}^\prism \cong X_{(p,\infty)}^\prism$ and 
$$
\mathrm{coeq}(\varphi, \mathrm{id}: X_{(1,\infty)}^\prism \to X_{(1,\infty)}^\prism)\cong X^\HK.
$$
\item The natural map $X^{\widetilde{\Div^1}} \to X^\prism$ is an isomorphism over $(0,p)$ under $\kappa$. Thus, the Frobenius induces an isomorphism $\varphi: X_{(0,1)}^\prism \cong X_{(0,p)}^\prism$ and 
$$
\mathrm{coeq}(\varphi, \mathrm{id}: X_{(0,1)}^\prism \to X_{(0,p)}^\prism)\cong X^{\Div^1}.
$$
\end{enumerate}
\end{proposition}
\begin{proof}
The key (and somewhat confusing) point of the proof is the following.  Let $A$ be a nilperfectoid ring, $T=\GSpec(A)$. Let $D\to \n{X}_T$ be a degree $1$ divisor. 

There is the radius map of $\ob{rad}_T\colon \n{X}_T\to (0,\infty)$ as explained in \Cref{rk:radius-discs}. Under the pushout  presentation $\n{X}_T=\FFdisc{T}^{\times}\bigsqcup_{(\bb{N}\times T^u)} (\bb{N} \times T)$, the radius map sends $n\times T$ to $p^n$.  On the other hand, the map $X^{\prism}\to \FFdisc{X^{\diamond}}^{\times,\dR}$ produces another radius map $\ob{rad}_X\colon X^{\prism}\to (0,\infty)$ that arises from the map of diamonds 
\[
\FFdisc{X^{\diamond}}^{\times,\diamond} = X^{\diamond}\times_{\Marc(\bb{F}_p)} \Marc(\bb{Q}_p) \to  \Marc(\bb{Q}_p)\times_{\Marc(\bb{F}_p)} \Marc(\bb{Q}_p)
\]
and the natural radius map of $\Marc(\bb{Q}_p)\times \Marc(\bb{Q}_p)$. 

The two radius maps are related as follows: given a degree $1$ divisor $D\to \n{X}_T$ and a map $D\to X$ giving rise to a point $\GSpec(A)\to X^{\prism}$, the divisor $D\to \n{X}_T$ lands in the locus $\n{X}_{T, (r,s)}$ if and only if the map $\GSpec(A)\to X^{\prism}$ lands in $X^{\prism}_{(1/s,1/r)}$, see \cite[Proposition 3.23]{hauck} for a discussion. Indeed, we can assume without loss of generality that $A=K$ is an algebraically closed  perfectoid field over $\bb{Q}_p$, for convenience in the notation, let us write $\pi\in K$ a pseudo-uniformizer with the same norm of $p$.  let $\pi^{\flat}=(\pi^{1/p^n})_{n\geq 0}$ denote a compatible sequence of $p$-power roots of $\pi$ and $[\pi^{\flat}]$ its Teichm\"uller lift. Then, the radius map 
\[
\ob{rad}\colon \FFdisc{K^{\flat}}^{\times}\to (0,\infty)
\]
is precisely $-\log_{|\pi^{\flat}|}(|p|)$ where for a Gelfand ring $B$ and $b\in K$, the map $|b|\colon \GSpec (B) \to [0,\infty)$ is the composite of $\GSpec(B)\to \n{M}(B)_{\Betti}$ and the norm map $|b|\colon \n{M}(B)_{\Betti}\to [0,\infty)$.  Let $f\colon \GSpec (K^{\sharp}) \to  \FFdisc{K^{\flat}}^{\times}$ be a degree $1$ divisor corresponding to an untilt $K^{\sharp}$ of $K^{\flat}$. We have the diagram
\[
\GSpec (K^{\sharp}) \to  \FFdisc{K^{\flat}}^{\times}  \cong   \FFdisc{K^{\sharp, \flat}}^{\times}  \leftarrow \GSpec (K).
\]
Then the behaviour of the radius map for $\n{X}_T$ and $X^{\prism}$  is just the fact that for  $\n{X}_T$ we fix the pseudo-uniformizer of $K^{\flat}$ to be $\pi^{\flat}$, while for $X^{\prism}$ we use a pseudo-uniformizer $\pi'\in K^{\sharp}$ (which is $p$ up to a unit of norm $1$) and a tilt $\pi^{\prime,\flat}\in K^{\flat}$ instead. Then, the relation follows from the fact that  in the disc $\FFdisc{K^{\flat}}^{\times}$ one has
\[
|\pi^{\prime,\flat}|\leq |\pi^{\flat}|^{r} \mbox{ if and only if } |\pi^{\flat}|^{1/r} \leq |\pi^{\prime,\flat}|
\]
with respect to the $p$-adic norm, and that at $x=\GSpec (K^{\sharp})$ (resp.\ $y=\GSpec (K)$) one has $|p|_x=|[\pi^{\prime,\flat}]|_x$ (resp.\ $|p|_y=|[\pi]|_y$) by design.

Thus,  the description of $\delta$-spaces of \Cref{Figure1FFDiscs} produces natural isomorphisms $\varphi\colon \n{X}_{T,(1/p,\infty)}\xrightarrow{\sim} \n{X}_{T,(1,\infty)}$ and $\varphi\colon \n{X}_{T,(0,1/p)}\xrightarrow{\sim} \n{X}_{T,(0,1)}$, furthermore, in the locus $(0,1)$ the disc is reduced. By the previous discussion on radius  maps, we deduce isomorphisms $\varphi\colon X^{\prism}_{(0,1)}\xrightarrow{\sim} X_{T,(0,p)}$ and $\varphi\colon X^{\prism}_{(1,\infty)}\xrightarrow{\sim} X^{\prism}_{ (p,\infty)}$ respectively.

Finally, the identification of the coequalizers of $\ob{id},\varphi$ with $X^{\Div^1}$ and $X^{\HK}$ follow from the previous discusion on radius maps and the  fact that the natural maps
\[
\n{X}_{T,(1/p,\infty)}\xrightarrow{\sim} \FFdisc{A,(1/p,\infty)}^{\times} \mbox{ and }  \FFdisc{A^{u,\flat},(0,1)}^{\times} \xrightarrow{\sim} \n{X}_{T, (0,1)}
\]
are isomorphisms (recall $\FFdisc{A}^{\times} $ is the perfect thickened Fargues--Fontaine disc, see \Cref{sec:real-farg-font-1-perfected-version}). 
\end{proof}

\begin{corollary}\label{cor:prismatization-push-out}
Let $X$ be a $\dagger$-cohesive Gelfand stack. Via the isomorphism $X^{\widetilde{\Div^1}}_{(1,p)} \cong X_{(1,p)}^\prism \cong (\FFdisc{X^{\diamond}}^\times)_{(1,p)}^\dR$, the prismatization $X^\prism$ is a pushout
$$
X^\prism \cong X^{\widetilde{\Div^1}}_{(0,p)} \bigsqcup_{X^{\widetilde{\Div^1}}_{(1,p)}} (\FFdisc{X^{\diamond}}^{\times,\dR})_{(1,\infty)}.
$$
\end{corollary}
\begin{proof}
We have on Betti stacks
\[
(0,\infty) \cong (0,p) \bigsqcup_{(1,p)} (1,\infty).
\]
The statement follows from this,  \Cref{PropositionDescriptionPrismatizationHKDivtilde}, and the fact that in Gelfand stacks (like in any topos), colimits commute with pullbacks.
\end{proof}

Consequently, understanding $X^\prism$ boils down to understanding $(\FFdisc{X^{\diamond}}^\times)^\dR$ and $X^{\widetilde{\Div^1}}$, or what is the same $X^\HK$ and $X^{\Div^1}$.  We will focus on the latter, since the former has already been analyzed in \cite{drff_stacks} and only depends on the perfectoidization of $X$.

\begin{remark}\label{sec:breuil-kisin-twist-1-modification-associated-to-bk-twist}
  As a consequence of \Cref{cor:prismatization-push-out}, for $X$ a $\dagger$-cohesive Gelfand stack, the category formed by vector bundles $\mathcal{E}$ on $X^\prism$ together with an isomorphism $\varphi_E\colon \varphi^\ast \mathcal{E}[1/\mathcal{I}_\prism] \cong \mathcal{E}[1/\mathcal{I}_\prism]$ (where we recall that $\mathcal{I}_\prism$ denotes the ideal sheaf of $X^\HT$ inside $X^\prism$) could be understood as \textit{rational shtukas},  i.e.\ triples formed by  a vector bundle $\mathcal{F}_1$ on $X^{\Div^1}$, a vector bundle $\mathcal{F}_2$ on $X^\HK$ and a meromorphic isomorphism $\mathcal{F}_{1}\dasharrow \mathcal{F}_2$ outside the image of $X^\HT$.
  We don't develop this theory in detail as a powerful theory of coefficients for rational analytic syntomic cohomology has been defined in \cite{hauck}. 
\end{remark}

\subsection{The $p$-adic twistor stack} \label{sec:p-adic-twistor-stack}
 Let $X$ be a  Gelfand stack. Passing to the quotient by Frobenius in \Cref{RemarkNaturalMapsXprism}, we get for the $p$-adic twistor stack of $X$ a diagram
\begin{equation} \label{composition-for-twistor}
\FF_{X^\diamond} \to X^{\Div^1}\to \FF_{X^\diamond}^\dR= X^\HK.
\end{equation}

From \Cref{prop:pullback-description-twistor-stack1}, and its obvious variant for $X^{\widetilde{\Div^1}}$, we deduce the following very useful statement:

\begin{proposition}\label{prop:pullback-description-twistor-stack}
Let $X$ be a Gelfand stack. Restriction of degree $1$ divisors on $\FF_A$ to $\GSpec(A)$ and $\FF_{A^u}$ gives rise to a commutative square
\[\begin{tikzcd}
	{X^{\Div^1}} \ar[r] \ar[d] & {X^\HK}  \ar[d] \\
	{X^{\rm cone}} \ar[r] & {(X^{\rm cone})^\dR.}
      \end{tikzcd}\]
      If $X$ is $\dagger$-cohesive, this square is cartesian.
\end{proposition}
\begin{proof}
This follows formally from the pushout diagram of thickened Fargues--Fontaine curves $\FF_T= \FF_{T^u}\bigsqcup_{T^u} T$ where $T=\GSpec(A)$ is the spectrum of a nilperfectoid ring $A$. 
\end{proof}

Outside the Hodge--Tate stack the $p$-adic twistor stack looks like a Hyodo--Kato stack.

\begin{proposition}\label{PropTwistorOutsideHTisHK}
Let $X$ be a Gelfand stack. The morphism $ X^{\Div^1}\to X^{\HK}$ is an isomorphism after base change to $\Div^1\setminus \Q_p^\HT\cong \Div^{1,\dR}\setminus \GSpec(\Q_p)$. 
\end{proposition}
\begin{proof}
  This follows directly from the definitions and the observation that for $A$ nilperfectoid, the map of  Fargues--Fontaine curves $\FF_{A^{u,\flat}}\to \FF_A$ is an isomorphism outside the locus $\GSpec(A)$.
\end{proof}

\subsection{The Hodge--Tate stack}\label{sec:hodge-tate-stack}
We have seen that outside the Hodge--Tate locus, the $p$-adic twistor stack is just a Hyodo--Kato stack. Next, we describe the geometry of the Hodge--Tate stack.
For that, we make the following definition. 

\begin{definition}\label{DefinitionIdealDiv1HTDivisor}
We let $\n{I}_{\Div^1}\to \s{O}_{\Div^1}$ denote the ideal defined by the Hodge--Tate divisor $\bb{Q}_p^{\HT}\to \Div^1$. Its base change to $\bb{Q}_p^{\HT}$ will be denoted by $\overline{\n{I}}_{\Div^1}$. 
\end{definition}

The following lemma produces the expected description of the Hodge--Tate stack for $\bb{G}_a$ as a square-zero extension. 

\begin{lemma}\label{LemmaGerbeDescriptionGa}
  Let $A$ be a nilperfectoid ring with $T=\GSpec(A)$ and let $D\to \FF_T$ be a degree $1$ divisor.
  Suppose we are given a factorization $T\to D\to \FF_T$ that corresponds to a factorization $T\to \bb{Q}_p^{\HT}\to \Div^1$.   The following hold:

\begin{enumerate}

\item  The natural square-zero thickening $D^{(1)}\to \FF_{T}$ arising from the divisor $D\to \FF_T$ gives rise to natural derivations 
\[
d\colon A\to  A\oplus  \n{I}_{A}[1]
\]
and 
\[
\overline{d}\colon A^u\to A^u \oplus \n{I}_{A^u}[1].
\]
Furthermore,  the derivation $\overline{d}$ corresponds to the square-zero extension $\bb{B}^+_{\dR}(A^u)/\xi^2\to A^u$, where $\xi$ is any generator of the kernel of Fontaine's map $\theta\colon \bb{B}^+_{\dR}(A^u)\to A^u$.

\item Let $\overline{d}_A\colon A\to A\oplus \n{I}_{A^u}[1]$ be the derivation on $A$ induced by $\overline{d}$ by composing with the projection $A\to A^u$. Then we have a natural derivation 
\[
\widetilde{d}\colon A\to A\oplus \ob{Nil}^{\dagger}(A)\otimes_A \n{I}_A[2]
\]
obtained by composing $\overline{d}_A$ with the connecting map of the fiber sequence 
\begin{equation}\label{eqFiberSequenceDerivationNilDagger}
\n{I}_{A}[1]\to \n{I}_{A^u}[1] \to \ob{Nil}^{\dagger}(A)\otimes_A \n{I}_A[2].
\end{equation}
In particular, the derivation $d\colon A\to A\oplus \n{I}_A[1]$ arises from $\widetilde{d}$ and the pullback square
 \begin{equation}\label{eqPullbackSquaresqextHTGa}
 \begin{tikzcd}
A\oplus \n{I}_A[1] \ar[d] \ar[r] & A \ar[d,"0"] \\
A\oplus \n{I}_{A^u}[1]   \ar[r]  &  A\oplus  \ob{Nil}^{\dagger}(A)\otimes_A \n{I}_A[2]
\end{tikzcd}
 \end{equation}

\item The ring  $\s{O}(D)$ naturally fits in a a commutative diagram with pullback squares
\[
\begin{tikzcd}
\s{O}(D) \ar[r] \ar[d] &  A\oplus \n{I}_A[1] \ar[r] \ar[d] & A \ar[d,"0"] \\
A \ar[r,"\overline{d}_A"] \ar[d] & A\oplus \n{I}_{A^u}[1] \ar[r] \ar[d]  & A\oplus \ob{Nil}^{\dagger}(A)\otimes_A \n{I}_A [2] \\ 
A^u \ar[r,"\overline{d}"] & A^u \oplus \n{I}_{A^u}[1] & 
\end{tikzcd}
\]
In particular, $\s{O}(D)\to A$ is a square-zero extension with ideal $\ob{Nil}^{\dagger}(A)\otimes_{A} \n{I}_A[1]$. 

\end{enumerate}

\end{lemma}
\begin{proof}

In the following we denote $T^u:=\GSpec(A^u)$,   $\overline{D}:=D\times_{\FF_T} \FF_{T^u}$ and for $\GSpec(B)\to \FF_T$ we set $D_B:=D\times_{\FF_T} \GSpec(B)$.

\begin{enumerate}

\item  By construction, see \Cref{ExampleDivisorCase} (3), the square-zero extension $D^{(1)}$ defined by the divisor $D\to \FF_T$ corresponds to the derivation of $D$ associated to the second projection $\pi_2\colon D\times_{\FF_T} D\to D$ where $D\times_{\FF_T} D$ is considered as a pointed stack over $D$ via the first projection $\pi_1\colon D\times_{\FF_T} D\to D$ and the  diagonal map $\Delta\colon D\to D\times_{\FF_T} D$.  The associated derivations of $A$ and $A^u$ are obtained by base change of $\pi_2$ along the maps $T^u\to T\to D$. In particular, for $T^u$, the square-zero thickening is, by \Cref{ExampleDivisorCase} (3), the square-zero extension given by $\bb{B}^+_{\dR}(A^u)/(\xi^2)$ as $\bb{B}^{+}_{\dR}(A^u)$ is precisely the formal completion of $\FF_{T^u}$ at the divisor  $\overline{D}$. 

\item The construction of $\widetilde{d}$ is  formal from the fiber sequence \eqref{eqFiberSequenceDerivationNilDagger}, as well as the  the pullback square \eqref{eqPullbackSquaresqextHTGa}. The fact that $d\colon A\to A\oplus \n{I}_A(A)[1]$ is recovered from the derivations $\overline{d}_A\colon A\oplus \n{I}_{A^u}(A)[1]$  and $\widetilde{d}$ via the pullback diagram \eqref{eqPullbackSquaresqextHTGa} follows from the simple observation that \eqref{eqFiberSequenceDerivationNilDagger} is a fiber sequence and that the composite $A\oplus \n{I}_A[1]\to A\oplus \n{I}_{A^u}[1]$ sends $d$ to $\overline{d}_A$.

\item We first construct the left vertical squares and prove that they are cartesian. Indeed, we have a pushout diagram of Berkovich spaces $\FF_T = \FF_{T^u}\bigsqcup_{T^u} T$ which produces a pushout diagram $D= \overline{D}\bigsqcup_{D_{T^u}} D_T$. Taking the underlying rings, and noticing the $\overline{D}=\GSpec(A^u)$ as we are over the Hodge--Tate point, this produces the cartesian diagram 
\begin{equation}\label{eqCartesianBasicPresentationOD}
\begin{tikzcd}
\s{O}(D) \ar[r] \ar[d] & A\oplus \n{I}_A[1] \ar[d] \\ 
A^u  \ar[r]& A^u \oplus \n{I}_{A^u} [1]
\end{tikzcd}
\end{equation}
Now, since $D$ sits in the Hodge--Tate locus, we have a morphism $T^u\to T\to D$, so that the map $\s{O}(D)\to A^u$ will factor through $\s{O}(D)\to A$. This produces the commutative diagram 
\begin{equation}\label{eqCartesianODvsAAI1}
\begin{tikzcd}
\s{O}(D) \ar[r] \ar[d] & A\oplus \n{I}_A[1] \ar[d] \\ 
A  \ar[r]& A \oplus \n{I}_{A^u} [1]
\end{tikzcd}
\end{equation}
which must be cartesian as \eqref{eqCartesianBasicPresentationOD} is cartesian and $A= A^u\times_{A^u\oplus \n{I}_{A^u}[1]} (A\oplus \n{I}_A[1])$.

On the other hand, part (2) produces for us the cartesian diagram  \eqref{eqPullbackSquaresqextHTGa}. Combining with \eqref{eqCartesianODvsAAI1} we get the full commutative diagram of part (3), and by using the pasting lemma we deduce that all squares are cartesian, as wanted. \qedhere
\end{enumerate}

\end{proof}

\Cref{LemmaGerbeDescriptionGa} gives a precise way of viewing the ring stack $\mathbb{G}^{\HT}_a$ over $\Q_p^\HT$ as a square-zero extension of $\Ga\times_{\GSpec(\Q_p)}\Q_p^\HT$ by $\Ga^{\dagger}\{1\}[1]$.
We note that this is very much in line with the prismatization of $p$-adic formal schemes, cf.\ \cite[Construction 5.10]{bhatt2022absolute}, if one replaces $\Ga^{\sharp}$ by $\Ga^\dagger$.

\begin{proposition}\label{sec:defin-first-prop-2-first-properties-of-analytic-prismatization}
Let $X$ be a  Gelfand stack. 
\begin{enumerate}

   \item There exists a natural commutative diagram of Gelfand stacks
\begin{equation}\label{EquationConeConstructionvsHodgeTate}
\begin{tikzcd}
X^{\HT}  \ar[r] \ar[d] & X^{\dR} \ar[d] \\ 
X^{\ob{cone}}_s  \ar[r] & (X^{\ob{cone}}_s)^{\dR},
\end{tikzcd}
\end{equation}
   where   $X^{\ob{cone}}_s \to X\times B\bb{G}_m$ is the special fiber of the cone of X, see \Cref{sec:defin-first-prop-1-definition-cone-construction}. In particular, we have a natural map     
  \[
    \pi_X\colon X^\HT\to X\times \Q_p^\HT,
  \]
  where the morphism $X^\HT\to X$ is induced by $X^\cone_s\to X$, cf.\ \Cref{RemarkMapSpecialFiberConeToX}.

  \item   Assume that $X$ is $\dagger$-cohesive, then  \eqref{EquationConeConstructionvsHodgeTate} is cartesian. If in addition $X$ admits a cotangent complex (cf.\ \Cref{sec:cotang-compl-morph-definition-having-a-cotangent-complex}), then the map $X^{\HT}\to X\times \bb{Q}_p^{\HT}$ is a principal homogeneous space for the overconvergent shifted and twisted tangent complex  $(\bb{T}_X[1]\{1\})^{\dagger}$  (see \Cref{RemarkdeRhamStackTangentComplex}), where $\{1\}$ refers to a twist by the ideal $\overline{\n{I}}_{\Div^1}$ (see \Cref{DefinitionIdealDiv1HTDivisor}).

 \item Suppose that $X$ is $\dagger$-cohesive and  admits a cotangent complex.  We have a pullback square
  \[
 \begin{tikzcd}
 X^{\HT} \ar[r,"\pi_{X}"]  \ar[d] & X\times \bb{Q}_p^{\HT} \ar[d,"0"] \\ 
 X\times \bb{Q}_p^{\HT} \ar[r,"\widetilde{d}_X"] & (\bb{T}_{X}\{1\}[2])^{\dagger}
 \end{tikzcd}
 \]
 
 Thus, if $X$ is formally smooth with $\bb{T}_X$ a vector bundle with overconvergent neighbourhood at zero $\bb{T}_X^{\dagger}\subset \bb{T}_{X}$, the map $\pi_X$ is an epimorphism and has the structure of a  $B\bb{T}_X^{\dagger}\{1\}$-torsor, i.e., a gerbe banded by $\bb{T}_X^{\dagger}\{1\}$. 
\end{enumerate} 
\end{proposition}
\begin{proof}

\begin{enumerate}

\item The existence of the commutative diagram follows from \Cref{prop:pullback-description-twistor-stack} by taking pullbacks along $B\bb{G}_m\to \bb{A}^1/\bb{G}_m$ and $\Div^1\to \bb{A}^1/\bb{G}_m$. This produces a natural map 
\[
X^{\HT}\to X^{\cone}_s\times_{B\bb{G}_m} \bb{Q}_p^{\HT} \to X\times \bb{Q}_p^{\HT},
\]
where the second map arises form the natural morphism $X^{\cone}_s\to X\times B\bb{G}_m$ of \Cref{RemarkMapSpecialFiberConeToX}.

\item Assume that $X$ is  $\dagger$-cohesive. Let $A$ be a nilperfectoid ring  with $T=\GSpec(A)$ and $D\to \FF_T$ a degree $1$ divisor. Suppose that there is a factorization $\GSpec(A)\to D \to \FF_T$ corresponding to a factorization $\GSpec(A)\to \bb{Q}_p^{\HT}\to \Div^1$.  The commutative diagram \eqref{EquationConeConstructionvsHodgeTate} appears by evaluating $X$ at  \eqref{eqCartesianBasicPresentationOD} which remains cartesian as $X$ is $\dagger$-cohesive.  Taking base change to $X\times \bb{Q}_p^{\HT}$, we have a cartesian diagram
\begin{equation}\label{eqHTCartesianDiagram1}
\begin{tikzcd}
 X^{\HT} \ar[r] \ar[d] & X \times \bb{Q}_{p}^{\HT}  \ar[d] \\ 
X_s^{\cone}\times_{X\times B\bb{G}_m} X\times \bb{Q}_p^{\HT} \ar[r]  &  (X_s^{\cone})^{\dR}\times_{X^{\dR}\times B\bb{G}_m} X\times \bb{Q}_p^{\HT}
\end{tikzcd}
\end{equation}
But now, as $X$ admits a cotangent complex, \Cref{sec:cone-construction-1-geometry-of-cone-construction} yields the equivalence $X^{\cone}_s = \bb{T}_{X}[1]\{1\}$ (where in \textit{loc.\ cit.} the $(-1)$ twist corresponds to the $\{1\}$-twist here). Therefore, the   bottom map of \eqref{eqHTCartesianDiagram1} is the map of abelian groups over $X\times \bb{Q}_p^{\HT}$
\[
\bb{T}_{X}[1]\{1\}\times_{X\times B\bb{G}_m} X \times \bb{Q}_p^{\HT}\to  (\bb{T}_{X}[1]\{1\})^{\dR}\times_{X\times B\bb{G}_m} X \times \bb{Q}_p^{\HT},
\]
in particular it has fiber $(\bb{T}_{X}[1]\{1\})^{\dagger}$. This implies that $X^{\HT}\to X\times \bb{Q}_p^{\HT}$ is a principal homogeneous space for $(\bb{T}_{X}[1]\{1\})^{\dagger}$ as desired. 

\item  The pullback diagram follows from $\dagger$-cohesivity and the pullback square  
\[
\s{O}(D)=A\times_{A\oplus \ob{Nil}^{\dagger}(A)\otimes_{A} \overline{\n{I}}_{\Div^1}[2]} A
\]
of  \Cref{LemmaGerbeDescriptionGa} (3).
The last assertion for formally smooth stacks is then clear as $(\bb{T}_{X}\{1\}[2] )^{\dagger}$ is nothing but $B^2 \bb{T}_X^{\dagger} \{1\}$ as a $\bb{G}_a$-module over $X\times \bb{Q}_p^{\HT}$ (without the formal smoothness, the object $(\bb{T}_{X}\{1\}[2] )^{\dagger}$ remains more inexplicit).
Alternatively, by formal smoothness $\mathbb{L}_{X/\Q_p}$ is a vector bundle, which implies that the bottom line in \eqref{eqHTCartesianDiagram1} is an epimorphism of Gelfand stacks.\qedhere
 \end{enumerate}

\end{proof}

\Cref{prop:pullback-description-twistor-stack,sec:defin-first-prop-2-first-properties-of-analytic-prismatization} admit relative variants that are important for studying geometric aspects of $p$-adic twistor stacks.  

\begin{proposition}\label{PropRelativeDescriptionPrismatization}
Let $(S, \n{X},D)$ be a perfect analytic prism and let $X\to D$ be a  $\dagger$-cohesive morphism of Gelfand stacks. The following hold:

\begin{enumerate}

\item We have  a  commutative diagram with cartesian squares
\[
\begin{tikzcd}
X^{\Div^1/\n{X}} \ar[r] \ar[d] & X^{\HK/\n{X}}  \ar[r] \ar[d] & X^{\HK} \ar[d] \\ 
X^{\ob{cone}/(D,\n{X}/\varphi^\Z)}  \ar[r] & (X^{\ob{cone}/(D,\n{X}/\varphi^\Z)})^{\dR/\n{X}/\varphi^\Z} \ar[r]  &  (X^{\ob{cone}/(D,\n{X}/\varphi^\Z)})^{\dR},
\end{tikzcd}
\]
where $X^{\ob{cone}/(D,\n{X}/\varphi)}$ is as in \Cref{DefinitionRelativeConesCartierDivisors}.

\item Suppose that $X$ admits a cotangent complex, then there is a cartesian square
\[
\begin{tikzcd}
X^{\HT/D} \ar[r] \ar[d] & X \ar[d,"0"] \\  
 X \ar[r,"\widetilde{d}_X"] &  (\bb{T}_{X/D}\{1\}[2])^{\dagger},
 \end{tikzcd}
\]
where $\widetilde{d}_X$ is the derivation that arises from the map $\widetilde{d}$ of \Cref{LemmaGerbeDescriptionGa}.  Here $\{1\}$ refers to a twist by the pullback of the ideal $\n{I}_D$ that defines $D\to \n{X}$.

In particular, if $X\to D$ is formally smooth with $\bb{T}_{X/D}$ a vector bundle, and $\bb{T}^{\dagger}_{X/D}$ is the overconvergent neighbourhood at $0$, the map $X^{\HT/D}\to X$ is a $B\bb{T}^{\dagger}_{X/D}\{1\}$-torsor, aka, a $\bb{T}^{\dagger}_{X/D}\{1\}$-gerbe.

\end{enumerate}
\end{proposition}
\begin{proof}
This follows from the exact same arguments as those of \Cref{prop:pullback-description-twistor-stack,sec:defin-first-prop-2-first-properties-of-analytic-prismatization}. 
\end{proof}

\begin{remark}
  \label{sec:stacks-ydiv1-yht-relative-case}
  Keep the notation of \Cref{PropRelativeDescriptionPrismatization} and suppose that $X\to D$  is $\dagger$-cohesive and admits a cotangent complex. By  \Cref{sec:cone-construction-1-geometry-of-relative-cone-construction}  the  special fiber of the relative cone construction  $X^{\ob{cone}/(D,\mathcal{X}/\varphi^\Z)}_s\to X$ has sections  given by square-zero lifts of $X$ to the first nilpotent thickening $D^{(1)}\to \n{X}$ of $D$.
  In the case $X\to D$ is formally smooth (thus with $\bb{T}_{X/D}$ a vector bundle), the map 
  \[
    X^{\HT/D}\to X^{\ob{cone}/(D,\n{X}/\varphi^\Z)}
 \]
 is then a $\bb{T}_{X/D}^{\dR}\{1\}$-torsor.
 Moreover the map is equivariant for the map of groups $B\bb{T}_{X/D}^{\dagger}\{1\}\to B\bb{T}_{X/D}\{1\}$. In other words, the pushout of the Hodge--Tate stack $X^{\HT/D}\to X$ along the map of groups  $B\bb{T}_{X/D}^{\dagger}\{1\}\to B\bb{T}_{X/D}\{1\}$ is precisely the stack of square-zero lifts of $X$ to $D^{(1)}$. 
\end{remark}

\begin{remark}[Relation to the work of Bhatt--Zhang]
  \label{sec:div1-analyt-prism-simpson-gerbe-for-q-p}
  Let $Y\in \GelfStk$, and let $V\in \VB(Y)$ be a (geometric) vector bundle on $Y$.
  Let $V^\dagger$ be the overconvergent neighborhood of the zero section of $V\to Y$.
  Assume that $Z\to Y$ is a gerbe banded by $V^\dagger$, which is classified by a morphism $Y\to B^2V^\dagger$.
  We can make the following general construction in order to produce a $\Gm^\dagger$-gerbe over the dual $V^\vee$ of $V$. In the following all products and classifying stacks are taken relative to $Y$.

  We have the natural evaluation map $\ob{ev}\colon V^\vee\times B^2V^\dagger\to B^2\Gm^\dagger\times B^2V^{\dagger}$ (given by the natural pairing to $B^2\bb{G}_a^{\dagger}$ as $\bb{G}_a$-modules, composed with the second delooping of the exponential $B^2\exp\colon B^2\bb{G}_a^{\dagger}\to B^2\Gm^{\dagger}$), which defines a morphism of animated group stacks over $B^2V^\dagger$.
  Then we can define the \emph{abstract $\dagger$-Simpson gerbe associated $V$} as the pullback 
  \[
  \begin{tikzcd}
\mathcal{S}_V^\dagger \ar[r] \ar[d] & B^2V^{\dagger}  \ar[d,"0"] \\ 
V^\vee\times B^2V^\dagger \ar[r,"\ob{ev}"] & B^2\Gm^\dagger\times B^2V^{\dagger},
\end{tikzcd}
\]
where the right vertical map is the zero section of the group $B^2\Gm^{\dagger}$ relative to $B^2V^{\dagger}$.

The group $\mathcal{S}_V^\dagger$ sits in a fiber/cofiber sequence  of groups over $B^{2}V^{\dagger}$
\[
B\Gm^{\dagger}\times B^{2}V^{\dagger}\to \mathcal{S}^{\dagger}_V\to V^{\vee}\times B^2V^{\dagger},
\]
which in particular says that $S^{\dagger}_V\to V^{\vee}\times B^2V^{\dagger}$ is a $\Gm^{\dagger}$-gerbe.
Note that we can pushforward $\mathcal{S}_V^\dagger$ along $B\Gm^\dagger\to B\Gm$ to obtain a $\Gm$-gerbe $\mathcal{S}_V$ over $V^\vee\times B^2V^\dagger$ (the \emph{abstract Simpson gerbe associated with $V$}).
Pulling $\mathcal{S}_V$ back along $Y\to B^2V^\dagger$ defines the \emph{abstract Simpson gerbe associated with the $V^\dagger$-gerbe $Z$.}

Let $X$ be a smooth  rigid space over a perfectoid space $D$.
Applying this construction to the $\bb{T}_{X/D}^{\dagger}\{1\}$-gerbe $X^{\HT/D}=X^\HT$ over $X$ (\Cref{PropRelativeDescriptionPrismatization}\footnote{More precisely, we note that by \Cref{sec:desc-analyt-prism-descent-for-hodge-tate-stack} the formation of the Hodge-Tate stack satisfies descent along formally \'etale $!$-covers, e.g., Tate-rational localizations. This allows to propagate the geometry of the Hodge-Tate stack from the smooth affinoid case to the general smooth rigid case.}), one obtains a $\Gm$-gerbe over the twisted cotangent bundle $\bb{T}_{X/D}^*\{-1\}$, which is the Simpson gerbe of Bhatt--Zhang \cite{bhatt-aspects}, as proved in \cite{RCCartierGerbe}. The paper \cite{RCCartierGerbe} also shows that the Simpson gerbe is Cartier dual to (an extended version of) the Hodge--Tate stack in this situation, so that weight $1$ sheaves on the Simpson gerbe are equivalent to sheaves on the Hodge--Tate stack.

 This construction can also be applied to the $\Gm^\dagger$-gerbe $\pi_{\Q_p}\colon \Q_p^\HT\to \GSpec(\Q_p)$ and we obtain a $\Gm$-gerbe $\mathcal{S}_{\Q_p}\to \A^1$ (admitting a reduction to a $\Gm^\dagger$-gerbe).
Following \cite{bhatt-aspects,RCCartierGerbe} it has the following features:
\begin{enumerate}
\item $\mathcal{S}_{\Q_p}$ classifies line bundles on $\Q_p^\HT$,
\item $\mathcal{S}_{\Q_p}$ is Cartier dual to $\Q_p^\HT$, and weight $1$ sheaves on $\mathcal{S}_{\Q_p}$ are equivalent to quasi-coherent sheaves on $\Q_p^\HT$, and 
  \item $\mathcal{S}_{\Q_p}$ is \'etale locally split over $\A^1$.
  \end{enumerate}
  The gerbe $\mathcal{S}_{\Q_p}$ (and generalizations for tori over $\Q_p$) will be studied in more depth in upcoming work of Hellmann and Heuer.
  \end{remark}

  \subsection{Formal and overconvergent Hodge--Tate stacks}
  \label{sec:form-overc-hodge}

 Given \Cref{sec:defin-first-prop-2-first-properties-of-analytic-prismatization} and the results of \Cref{sec:hodge-tate-stack}, the natural next step is to analyze the geometry of the overconvergent or formal neighbourhood of $X^{\HT}$ in $X^{\Div^1}$. These stacks will play a role later in relation to the Riemann--Hilbert correspodence (\Cref{sec:relation-to-rh}), so they deserve their own name.

  \begin{definition}\label{OverconvergentHTLocus}
\begin{enumerate}

\item Let $X$ be a Gelfand stack. Its \textit{overconvergent  Hodge--Tate stack} is the pullback 
\[
X^{\dagger-\HT}:= X^{\Div^1}\times_{\bb{A}^1/\bb{G}_m} \bb{A}^{1,\dagger}/\bb{G}_m.
\]
Similarly, we define the \textit{formal Hodge--Tate stack} to be 
\[
X^{\widehat{\HT}}:=X^{\Div^1}\times_{\bb{A}^1/\bb{G}_m} \widehat{\mathbb{A}}^1/\bb{G}_m.
\]

\item Let $(S,\n{X},D)$ be a perfect analytic prism and let $X\to D$ be a morphism of Gelfand stacks. Let $D^\wedge\subset D^{\dagger}\subset \n{X}$ be the formal completion and overconvergent neighbourhood of $D$ in $\n{X}$ respectively. We define the \textit{relative overconvergent (resp.\ formal) Hodge--Tate stack} of $X$ over $(S,\n{X},D)$ to be  $X^{\dagger-\HT/D^{\dagger}}:=X^{\Div^1/\n{X}}\times_{\n{X}/\varphi^\Z} D^{\dagger}$ (resp.\ $X^{\widehat{\HT}/D^{\dagger}}:=X^{\Div^1/\n{X}}\times_{\n{X}/\varphi^\Z} D^\wedge$).

\end{enumerate}
  \end{definition}

 The overconvergent and formal Hodge--Tate stack satisfy the following formal properties.
   
\begin{proposition}\label{PropOverconvergentHT}
  \begin{enumerate}
 Let $X$ be a $\dagger$-cohesive Gelfand stack. 
\item     The natural commutative diagram
\[
  \begin{tikzcd}
X^{\widehat{\HT}} \ar[r] \ar[d]  &  X^{\dagger-\HT} \ar[r] \ar[d] &  X^{\dR}  \ar[d]    \\ 
  X^{\widehat{\cone}} \ar[r]  &		X^{\dagger-\cone}	 \ar[r]	& (X^{\cone}_s)^{\dR} 
  \end{tikzcd}
  \]
has cartesian squares, where $	X^{\dagger-\cone}$ and $X^{\widehat{\cone}}$  are as in \Cref{sec:defin-first-prop-1-definition-cone-construction}.

\item  Suppose that $X$ admits a cotangent complex and is formally smooth. Then $X^{\dR}\to (X_s^{\cone})^{\dR}=(X/\bb{T}_{X}(-1)\rtimes \bb{G}_m)^{\dR}$ is a $(\bb{T}_{X}(-1)\rtimes \bb{G}_m)^{\dR}$-torsor where $\bb{G}_m$ acts by weight $-1$ on the tangent space. In particular, the same holds for the maps $X^{\widehat{\HT}}\to X^{\widehat{\cone}}$ and $X^{\dagger-\HT}\to X^{\dagger-\cone}$ respectively. 
  
  \end{enumerate}

\end{proposition}
\begin{proof}
  Part (1) follows from \Cref{sec:defin-first-prop-2-first-properties-of-analytic-prismatization} (2) after taking pullbacks along the maps $\widehat{\mathbb{A}}^1/\bb{G}_m\to \A^{1,\dagger}/\bb{G}_m\to \bb{A}^1/\bb{G}_m$, using that $(\widehat{\mathbb{A}}^1/\bb{G}_m)^{\dR}=(\A^{1,\dagger}/\bb{G}_m)^{\dR}= B\bb{G}_m^{\dR}$ to identify the pullback of the Hyodo--Kato stack with $X^{\dR}$, and therefore that 
  \[
  (X^{\widehat{\cone}})^{\dR}= (X^{\dagger-\cone})^{\dR}= (X^{\cone}_{s})^{\dR}. 
  \]
    Part (2)  follows  formally from the fact that if $X$ is formally smooth then $X^{\cone}_s= X/\bb{T}_X(-1)\rtimes \bb{G}_m$. 
\end{proof}

  \Cref{PropOverconvergentHT} admits a relative version as follows.

  \begin{proposition}\label{PropRelativeVersionOverconvergentHTStacks}
Let $(S,\n{X}, D)$ be a perfect analytic prism and $X\to D$ a Gelfand stack.

\begin{enumerate}

\item Suppose that  $X\to D$ is $\dagger$-cohesive, then there are natural cartesian diagrams
\[
  \begin{tikzcd}
X^{\widehat{\HT}/D^\wedge} \ar[r] \ar[d]  &  X^{\dagger-\HT/D^{\dagger}} \ar[r] \ar[d] &  X^{\dR/D^{\dagger}}  \ar[d] \ar[r]  & X^{\dR}   \ar[d] \\ 
  X^{\widehat{\cone}/(D,\n{X}/\varphi^\Z)} \ar[r]  &		X^{\dagger-\cone/(D,\n{X}/\varphi^\Z)}	 \ar[r]	& (X^{\cone/(D,\n{X}/\varphi^\Z)}_s)^{\dR/D^{\dagger}} \ar[r] & (X_s^{\cone/(D,\n{X}/\varphi^\Z)})^{\dR}.
  \end{tikzcd}
  \]

\item Suppose that $X\to D$ is $\dagger$-cohesive and  that $X\to D$ is formally smooth.  Then the map $X^{\dR}\to (X_{s}^{\cone/(D,\n{X}/\varphi^\Z)})^{\dR}$ is a $\bb{T}_{X/D}\{1\}^{\dR}$-torsor, where $\{1\}$ refers to a twist by the invertible ideal defining the divisor  $D$. In particular the same holds for the other vertical maps in the diagram of part (1).
\end{enumerate}
  \end{proposition}
  \begin{proof}
  Part (1) follows from \Cref{PropRelativeDescriptionPrismatization} after taking pullbacks along $D^\wedge\to D^{\dagger}\to \n{X}/\varphi^\Z$.  For part (2), it suffices to show that the map  of  arc-stacks $\pi_X\colon X^{\diamond}\to (X_s^{\cone/(D,\n{X}/\varphi^\Z)})^{\diamond}$ is a $\bb{T}_{X/D}\{1\}^{\diamond}$-torsor.   We note that, since $\n{X}$ is a Berkovich space, both $D$ and $\n{X}$ admit a cotangent complex over $\bb{Q}_p$, and since by assumption the map $X\to D$ admits a cotangent complex, then $X$ also admits a cotangent complex over $\bb{Q}_p$. 
  
  Thus, as $X$ is formally smooth over $D$,   \Cref{sec:cone-construction-1-geometry-of-relative-cone-construction} shows that  the natural map
  $$
  X^{\cone/(D,\n{X}/\varphi^\Z)}_s\to X
  $$
  is the  $B\bb{T}_{X/D}\{1\}$-torsor  of square-zero lifts of $X$ to $D^{(1)}$. Passing to arc-stacks we see that  
  \[
  (X^{\cone/(D,\n{X}/\varphi^\Z)}_s)^{\diamond} \to X^{\diamond}
  \]
  is then a $B\bb{T}_{X/D}\{1\}^{\diamond}$-torsor, and the natural section $X^{\diamond}\to (X^{\cone/(D,\n{X}/\varphi^\Z)}_s)^{\diamond}$ produced by the Hyodo--Kato stack  trivializes this torsor. In particular, this makes the map $X^{\diamond}\to  (X^{\cone/(D,\n{X}/\varphi^\Z)}_s)^{\diamond}$ a $\bb{T}_{X/D}\{1\}^{\diamond}$-torsor as wanted. 
  \end{proof}

  \begin{remark}\label{TorsorExistsOnlyDiamond}
    The proof of \Cref{PropRelativeVersionOverconvergentHTStacks} shows that if $X\to D$ is formally smooth, then $(X_{s}^{\cone/(D,\n{X}/\varphi^\Z)})^{\diamond}\to X^{\diamond}$ admits a section making the perfectoidization of the cone a trivial $B\bb{T}_{X/D}\{1\}^{\diamond}$-torsor.
    Thanks to \Cref{sec:cone-construction-1-geometry-of-relative-cone-construction},  the map $X_{s}^{\cone/(D,\n{X}/\varphi^\Z)}\to X$ is the stack of square-zero lifts of $X$ to $D^{(1)}$. In particular, if $X$ is a Berkovich space,  a section $X\to X_{s}^{\cone/(D,\n{X}/\varphi^\Z)}$ is equivalent to a lift $X^{(1)}\to D^{(1)}$ (see \Cref{LemmaStabilitySquareZeroGeometrically}).
    In other words, even though the map  $X_{s}^{\cone/(D,\n{X}/\varphi^\Z)}\to X$ might not admit a section, it does when passing to arc-stacks. Furthermore, by construction, this section corresponds to the natural square-zero extension $X^{\diamond}\to X^{\diamond,(1)}$ given by the first nilpotent thickening of the divisor $X^{\diamond}\to \FF_{X^{\diamond}}$. 
  \end{remark}

  \subsection{Faltings' extension}
  \label{sec:faltings-extension}
  
  \Cref{sec:defin-first-prop-2-first-properties-of-analytic-prismatization} (1) produces a map $X^{\dR}\to (X_s^{\cone})^{\dR}$ which comes from the natural square-zero extension produced by the $\bb{B}_{\dR}^+$-construction.  For $X$ formally smooth, this map is equivalent to a $(\bb{T}_X\{-1\}\rtimes \bb{G}_m)^{\diamond}$-torsor on $X^{\diamond}$.
  In the case of smooth rigid spaces over $\Q_p$  we can relate this torsor to the so-called \emph{Faltings' extension}.
  Furthermore, if $X$ is a smooth rigid space over an affinoid   perfectoid space $D$ (a smoothoid space in the sense of \cite{MR5022120}), we have a similar interpretation for the  $\bb{T}_{X/D}\{1\}^{\diamond}$-torsor over $X^{\diamond}$   produced by a lift $\widetilde{X}$ of $X$ from $\s{O}(D)$ to $\bb{B}_{\dR}^+(D)/\xi^2$, and the associated map $X^{\diamond}\to (X^{\cone/(D,\FF_D)}_s)^{\diamond}$ of \Cref{PropRelativeVersionOverconvergentHTStacks}. Bhargav Bhatt informed us that he has obtained analogous results (in a different language) in his joint work in progress with Heuer, Petrov and Vologodsky (cf. \cite{bhatt_takagi}). 

In order to make the previous discussion precise, we will \textit{define} the Faltings' extension via the aforementioned torsors,  and then compare this construction with other definitions in the literature.

  \begin{definition}
    \label{sec:p-adic-twistor-1-definition-faltings-extension}
    \begin{enumerate}
    \item Let $X$ be a smooth rigid space over $\Q_p$.
      Its Faltings' extension $\mathrm{FE}_X$ is the extension
      \[
        0\to \s{O}_{\widehat{X^\diamond}}\{1\}\to \mathrm{FE}_X\to \Omega_X^1\otimes_{\s{O}_X} \s{O}_{\widehat{X^\diamond}}\to 0
      \]
      on $\widehat{X^\diamond}$ whose associated $\bb{T}_X\{1\}$-torsor $\mathcal{L}_{\mathrm{FE}_X}$ over $\widehat{X^\diamond}$ is given by $\widehat{X^\diamond}\times_{X^{\cone}_s}(X\times B\Gm)$, where the morphism $X\times B\Gm\to X^\cone_s$ is the one from \Cref{sec:cone-construction-1-geometry-of-cone-construction} and the other one is the canonical one $\widehat{X^\diamond}\to X^\HT\to X^\cone_s$ .
    \item Let $D$ be an affinoid\footnote{We only added this assumption for simplicity of exposition, it can be removed easily.} perfectoid space  over $\bb{Q}_p$ with degree $1$ divisor $D\to \FF_{D}$. Let $D^\wedge$ denote the formal completion of $\FF_D$ along $D$: it is the Gelfand stack given by 
    \[
    D^\wedge= \varinjlim_n \GSpec \bb{B}_{\dR}^+(D)/\xi^{n},
    \]
where $\xi$ is a generator of $\bb{B}_{\dR}^+(D)\to \s{O}(D)$. 
    
Let $X\to D$ be a smoothoid space over $D$ and   let  $\widetilde{X}$ be a lift of $X$ to $D^{(1)}=B^+_\dR(D)/t^2$. Let $a\colon X\to X_s^{\cone/(D,\FF_D)}$ be the map $a$ that corresponds to the square-zero extension $\widetilde{X}$ (cf.\ \Cref{sec:cone-construction-1-geometry-of-relative-cone-construction}).

The Faltings' extension $\mathrm{FE}_X^{\widetilde{X}}$ associated with $\widetilde{X}$ is the extension
      \[
        0\to \s{O}_{\widehat{X^\diamond}}\{1\}\to \mathrm{FE}_X^{\widetilde{X}}\to \Omega_X^1\otimes_{\s{O}_X} \s{O}_{\widehat{X^\diamond}}\to 0
      \]
      on $\widehat{X^\diamond}$ whose associated $\bb{T}_{X/D}\{1\}$-torsor $\mathcal{L}_{\mathrm{FE}^{\widetilde{X}}_X}$ over $\widehat{X^\diamond}$ is given by
      \[
        \widehat{X^\diamond}\times_{X^{\cone/(D,\FF_D)}_s}X. 
      \]

    \end{enumerate}
    
In the previous definition  $\{1\}$ refers to a Tate twist, or equivalently, a twist by the ideal defining the square-zero extension $D\to \FF_D$.   
  \end{definition}

\begin{remark}[Faltings' extension as a difference of sections: geometric case] \label{RemarkFaltingsObstructionGeometric}
  Let $X$ be a smoothoid space over an affinoid perfectoid space $D$ over $\bb{Q}_p$.
  In the following we identify the geometric Faltings' extension with a difference of sections of $\widehat{X^\diamond} \to X_s^{\cone/(D,\FF_D)}$.
  This description of the Faltings' extension is analogue to that of Wang in \cite[Section 2.2]{wang2023p}.
  
Let $Y\to X\to D\to S$ be  morphisms  of Gelfand stacks, suppose they are infinitesimally cohesive and admit a cotangent complex. We have a commutative diagram of shifted tangent spaces over $Y\times \bb{G}_m$ with cartesian squares
\[
\begin{tikzcd}
\bb{T}_{Y/S}[1](-1) \ar[r]  \ar[d] & \bb{T}_{X/S}[1](-1)\times_X Y  \ar[r] \ar[d]& \bb{T}_{D/S}[1](-1)\times_{D} Y  \ar[d] \\ 
Y\times B\bb{G}_m \ar[r,"0"] & \bb{T}_{Y/X}[2](-1) \ar[d] \ar[r] & \bb{T}_{Y/D}[2]  (-1) \ar[d] \\ 
& Y\times B\bb{G}_m \ar[r,"0"] & \bb{T}_{X/D}[2](-1)\times_X Y 
\end{tikzcd}
\]
 Recall from \Cref{sec:cone-construction-1-geometry-of-cone-construction} that the upper diagram is the natural map of special fibers of relative cones
 \[
 Y^{\cone/S}_s\to X^{\cone/S}_s\times_{X} Y \to D^{\cone/S}_s\times_{D} Y. 
 \]

Let us now take $S=\GSpec (\bb{Q}_p)$, $D$ a perfectoid space, $X\to D$ a smoothoid  space over $D$, and $Y=\widehat{X^{\diamond}}$.  In this case, we have that $\bb{T}_{\widehat{X}^{\diamond}/D}=0$ (cf. \Cref{sec:tang-cotang-spac-cotangent-complexes-of-diamonds}) and obtain a commutative diagram with pullback squares
\begin{equation}\label{eqDiagramPullbackTangentSpacesComputation}
\begin{tikzcd}
\bb{T}_{\widehat{X^{\diamond}}}[1](-1) \ar[r]  \ar[d] & \bb{T}_{X}[1](-1)\times_X \widehat{X^{\diamond}}  \ar[r] \ar[d]& \bb{T}_{D}[1](-1)\times_{D} \widehat{X^{\diamond}}  \ar[d] \\ 
\widehat{X^{\diamond}}\times B\bb{G}_m \ar[r,"0"] & \bb{T}_{\widehat{X^{\diamond}}/X}[2](-1) \ar[d] \ar[r,"0"] & \widehat{X^{\diamond}}\times B\bb{G}_m  \ar[d,"0"]  \\ 
& \widehat{X^{\diamond}}\times B\bb{G}_m \ar[r,"0"] & \bb{T}_{X/D}[2](-1)\times_X \widehat{X^{\diamond}},
\end{tikzcd}
\end{equation}
which identifies $\bb{T}_{\widehat{X}^{\diamond}/X}[2](-1)\cong \bb{T}_{X/D}[1](-1)\times_X \widehat{X^{\diamond}}$ over $\widehat{X^{\diamond}}\times B\bb{G}_m$.  Now, the square-zero-extension $D^{(1)}=\bb{B}_{\dR}^+(D)/\xi^2$ of $D$ produces a natural section of the cone 
\[
D\xrightarrow{[D^{(1)}]} D^{\cone}_s,
\]
which  then produces a commutative diagram with pullback squares
\[
\begin{tikzcd}
\widehat{X^{\diamond}}^{\cone/(D,\FF_D)}_s \ar[r] \ar[d] & X^{\cone/(D,\FF_D)}_s  \ar[r] \ar[d] & D \ar[d]\\ 
\widehat{X^{\diamond}}^{\cone}_s \ar[r]  & X^{\cone}_s  \ar[r] & D^{\cone}_s,
\end{tikzcd}
\] 
which in combination with  \eqref{eqDiagramPullbackTangentSpacesComputation} gives rise to a commutative diagram with pullback squares
\[
\begin{tikzcd}
\widehat{X^{\diamond}}^{\cone/(D,\FF_D)}_s \ar[r] \ar[d] & X^{\cone/(D,\FF_D)}_s\times_{X}\widehat{X^{\diamond}}  \ar[r] \ar[d] & \widehat{X^{\diamond}} \ar[d,"{[D^{(1)}]}"]\\  
\bb{T}_{\widehat{X^{\diamond}}}[1](-1) \ar[r]  \ar[d] & \bb{T}_{X}[1](-1)\times_X \widehat{X^{\diamond}}  \ar[r] \ar[d]& \bb{T}_{D}[1](-1)\times_{D} \widehat{X^{\diamond}}  \ar[d] \\ 
\widehat{X^{\diamond}}\times B\bb{G}_m \ar[r,"0"] & \bb{T}_{\widehat{X^{\diamond}}/X}[2](-1) \ar[d] \ar[r,"0"] & \widehat{X^{\diamond}}\times B\bb{G}_m  \ar[d,"0"]  \\ 
& \widehat{X^{\diamond}}\times B\bb{G}_m \ar[r,"0"] & \bb{T}_{X/D}[2](-1)\times_X \widehat{X^{\diamond}}.
\end{tikzcd}
\]
From this big diagram we can deduce two important facts: 
\begin{enumerate}

\item  We have a cartesian diagram 
\[
\begin{tikzcd}
\widehat{X^{\diamond}}^{\cone/(D,\FF_D)}_s \ar[r] \ar[d] & \widehat{X^{\diamond}}  \ar[d,"\s{O}_{X^{\diamond}}(1)"] \\ 
\widehat{X^{\diamond}} \times B\bb{G}_m \ar[r,"\ob{id}"] & \widehat{X}^{\diamond}\times B\bb{G}_m,
\end{tikzcd}
\]
where the right map is the identity on $\widehat{X^{\diamond}}$ and the Tate twist on the projection to $B\bb{G}_m$.
 This implies that the natural map  $\widehat{X^{\diamond}}^{\cone/(D,\FF_D)}_s\to \widehat{X^{\diamond}}$ is an equivalence. Interpreting the cone as the natural obstruction class to lift $\widehat{X^{\diamond}}\to D$ to a morphism over $D^{(1)}$ via \Cref{sec:cone-construction-1-geometry-of-relative-cone-construction}, this says that there is a unique such lift given by $\bb{B}_{\dR}^+(\widehat{X^{\diamond}})/\xi^{2}$.

\item The previous identification produces a natural section  of the map 
\[
X_s^{\cone/(D,\FF_D)}\times_X \widehat{X^{\diamond}}\to \widehat{X^{\diamond}}.
\]
This section arises from the class 
\[
\ob{cl}\colon \widehat{X^{\diamond}}\to  \bb{T}_{\widehat{X}^{\diamond}}[1](-1) \to \bb{T}_{X}[1](-1)\times_X \widehat{X^{\diamond}},
\]
which is the image of the square-zero lift $[\bb{B}_{\dR}^+(\widehat{X^{\diamond}})]\colon \widehat{X^{\diamond}} \to \bb{T}_{\widehat{X}^{\diamond}}[1](-1)$. Since $X_s^{\cone/(D,\FF_D)}\to X$ is a $\bb{T}_{X/D}[1]\{-1\}$-torsor (with $\bb{T}_{X/D}[1]\{-1\}\to X$ being the shifted tangent space with a twist by the conormal cone of $D\to \FF_D$), the class $\ob{cl}$ gives rise to an identification 
\[
\alpha_{\ob{cl}}\colon  \bb{T}_{X/D}[1]\{-1\}\times_X \widehat{X^{\diamond}}   \xrightarrow{\sim} X_s^{\cone/(D,\FF_D)}\times_X \widehat{X^{\diamond}}.
\]
\end{enumerate}
Now,  by \Cref{sec:cone-construction-1-geometry-of-relative-cone-construction} the datum of a square-zero lift $\widetilde{X}\to D^{(1)}$ of $X\to D$ is equivalent to a section 
\[
[\widetilde{X}]\colon X\to X^{\cone/(D,\FF_D)}_s\to X,
\]
or equivalently, to a class $\widetilde{\ob{cl}}\colon X\to \bb{T}_{X}[1](-1)$ that lifts the class $X\xrightarrow{[D^{(1)}]} \bb{T}_{D}[1](-1)\times_D X$. Now, taking base change to $\widehat{X^{\diamond}}$, we have two maps 
\[
\ob{cl},\widetilde{\ob{cl}}\colon  \widehat{X^{\diamond}}\to  X^{\cone/(D,\FF_D)}_s\times_{X}\widehat{X}^{\diamond}.
\]
By trivializing the cone via the class $\ob{cl}$, we can think of the base change 
\[
\widehat{X^{\diamond}} \times_{ [D^{(1)}],X^{\cone/(D,\FF_D)}_s, [\widetilde{X}] } X     \xleftarrow{\alpha_{\ob{cl}}}  \widehat{X^{\diamond}}    \times_{ 0,  \bb{T}_{X/D}[1]\{-1\}\times_X, \widetilde{\ob{cl}}-\ob{cl} } \widehat{X^{\diamond}}    
\]
as the $\bb{T}_{\widehat{X^{\diamond}}/D}[1]\{1\}\cong \bb{T}_{X/D}\{1\}\times_{X}\widehat{X^{\diamond}}$-torsor   over  $\widehat{X^{\diamond}} $ associated to the class  $\widetilde{\ob{cl}}-\ob{cl}$.  By \Cref{sec:p-adic-twistor-1-definition-faltings-extension} (2), the class $\widetilde{\ob{cl}}-\ob{cl}$ is precisely the Faltings' extension $\n{L}_{\mathrm{FE}_{X}^{\widetilde{X}}}$.

\end{remark}

    \begin{remark}[Faltings' extension as an obstruction class]
    \label{sec:p-adic-twistor-3-functor-of-points-description}
    Let $f\colon Y\to X$ be a morphism of Gelfand stacks over a perfectoid space $D$.
    Assume given a section of the natural map $X_s^{\cone/(D,\FF_D)} \to X$, i.e., a square-zero lift $\widetilde{X}$ of $X$ to $D^{(1)}=\GSpec(B_\dR^+(D)/t^2)$, the first-order thickening of $D$ in $\FF_D$. Then the sections of $Y_s^{\cone/(D,\FF_D)} \times_{X_s^{\cone/(D,\FF_D)}} X$ (where the second map in the fiber product is defined by $\widetilde{X}$) are the same thing as square-zero lifts $\widetilde{Y}$ of $Y$ to $D^{(1)}$ together with a lift $\widetilde{f}\colon \widetilde{Y}\to \widetilde{X}$. 
    
    If $X \to D$ is a smoothoid space and $Y=\widehat{X^\diamond}$ with the natural map $f\colon \widehat{X^\diamond}\to X$, we have seen in \Cref{RemarkFaltingsObstructionGeometric} that $Y_s^{\cone/(D,\FF_D)}= (\widehat{X^\diamond})_s^{\cone/(D,\FF_D)}\cong \widehat{X^\diamond}$, so in this case the fiber product $Y_s^{\cone/(D,\FF_D)} \times_{X_s^{\cone/(D,\FF_D)}} X$ is our definition of the Faltings' extension. In particular, we deduce the following functor of points descriptions for the geometric Faltings' extension: if $A$ is a perfectoid ring over $\widehat{X^\diamond}$, $A$-points of $\mathcal{L}_{\mathrm{FE}_X^{\widetilde{X}}}$ correspond to morphisms $\GSpec(B_\dR^+(A)/t^2) \to \widetilde{X}$ over $D^{(1)}$, lifting the morphism $\GSpec(A)\to \widehat{X^\diamond}\to X$.\footnote{For a discussion in the same spirit (for different, integral, statements), see \cite[\S II.10]{abbes2016p}.} 
    
    In the arithmetic case, the same argument for viewing the Faltings' extension as an obstruction class goes through if one works relative to $\widehat{\Q_p^\diamond}$, i.e., we consider the cone construction relative to $(D,\FF_D):=(\widehat{\Q_p^\diamond},\FF_{\widehat{\bb{Q}_p^{\diamond}}})$.
    Let $X$ be a smooth rigid space over $\Q_p$.
    The canonical section $X\times B\Gm\to X^{\cone}_s$ can be pulled back to the special fiber of the cone construction relative to $(D,\FF_D)$.
    This yields a section $X\times D\to (X\times  D)^{\cone/(D,\FF_D)}$ which is classifying the canonical lift $X\times \widehat{B^+_\dR(\Q_p)}/t^2$ of $X\times \widehat{\Q_p^\diamond}$.
    The Faltings' extension is now the obstruction to lifting $\widehat{X^\diamond}\to X\times \widehat{\Q_p^\diamond}$ to $X\times \widehat{B^+_\dR/t^2}$.
    More precisely, the same argument as before applies  and we have
        \[
      (\widehat{X^\diamond})_s^{\cone/(D,\FF_D)} \times_{(X\times D)^{\cone/(D,\FF_D)}_s} (X\times D) \cong \widehat{X^\diamond}\times_{(X\times D)^{\cone/(D,\FF_D)}_s} (X\times D).
    \]
    The right hand side can be rewritten as the fiber product $\widehat{X^\diamond}\times_{X^{\cone}_s} (X\times B\Gm)$ as both fiber products are given by a $\bb{T}_{X\times D/D}\{1\}\cong (\bb{T}_{X/\Q_p}\times_{X}(X\times D))\{1\}$-torsor over $\widehat{X^\diamond}$.
    But this is exactly the Faltings' extension as defined in \Cref{sec:p-adic-twistor-1-definition-faltings-extension}.

    We note that the arithmetic case is slightly confusing because the Faltings' extension is not the obstruction class of deforming $\widehat{X^\diamond}\to X\times B\Gm$.
    More precisely, the argument breaks down because $\mathbb{L}_{\widehat{X^\diamond}/B\Gm}$ is not zero and the obstruction class, i.e.,
    \[
(\widehat{X^\diamond})^\cone_s\times_{(X\times B\Gm)^\cone_s}(X\times B\Gm),
\]
differs from the Faltings' extension, i.e., from $\widehat{X^\diamond}\times_{X^\cone_s}(X\times B\Gm)$.
\end{remark}

  In the arithmetic case, for $X$ a smooth rigid space over $\Q_p$, the Faltings' extension has been previously defined and we now show that our definition agrees with the already existing one. The standard definition of the Faltings' extension, e.g.\ \cite{scholze2013p}, is via period sheaves: one defines the Faltings' extension of $X$ as $\mathrm{gr}^1 \s{O}\mathbb{B}_{\dR,X}^+$, so that the associated torsor is the analytification of 
  \[
  \AnSpec(\s{O}\mathbb{C}_X) = \AnSpec(\mathrm{gr}^0 \s{O}\mathbb{B}_{\dR,X}).
  \]
  Our next goal, obtained in \Cref{cor:identification-fe} below, is to see that this coincides with $\widehat{X^\diamond}\times_{X^{\cone}_s}(X\times B\Gm)$. For that, we will use the algebraic cone construction as discussed in \Cref{RemarkAlgebraicConeGeneralAnalyticStacks}.
  
Let $\Spf(\bb{B}_{\dR,X}^+)$ denote the formal completion of the divisor $\widehat{X}^{\diamond}\to \FFdisc{X^{\diamond}}^{\times}$, equivalently, we write 
\[
\Spf(\bb{B}_{\dR,X}^+)=\varinjlim_n \underline{\GSpec}_{\FFdisc{X^{\diamond}}^{\times}} (\bb{B}^+_{\dR,X}/(\ker \theta)^n),
\]
where $\theta\colon \bb{B}_{\dR,X}^+\to \widehat{\s{O}}_X$ is Fontaine's map.   The natural Cartier divisor  $\Spf(\bb{B}_{\dR,X}^+)\to \widehat{\bb{A}}^1/\bb{G}_m$ corresponds to the  decreasing $\bb{Z}$-filtered ring given by $\Fil^i \bb{B}_{\dR,X}^+= \xi^i \bb{B}_{\dR,X}^+$ where $\xi$ is a local generator of $\ker \theta$. 

On the other hand, let $ \s{O}\mathbb{B}_{\dR,X}^+$ be the big de Rham period sheaf endowed with its natural decreasing $\bb{N}$-filtration given by the adic filtration along the kernel of $\theta\colon \s{O}\mathbb{B}_{\dR,X}^+\to \widehat{\s{O}}_{X}$, and let $\Spf ( \s{O}\mathbb{B}_{\dR,X}^+)$ denote its formal spectrum. The following statement provides a geometric construction of this Gelfand stack.

\begin{lemma}\label{LemmaReconstructionOBdR} 
  Let $X$ be a smooth rigid space over $\bb{Q}_p$.
  Let $X_{\proet}^{\omega_1}$ denote the pro-\'etale site of $X$ (in the sense of \cite{scholze2013p}) restricted to countable projective limits, and given an object $Y=\varprojlim_i Y_i$ in $X_{\proet}^{\omega_1}$ we let $Y$ denote its limit in Gelfand stacks (with each $Y_i$ viewed as a Gelfand stack via \Cref{sec:deriv-berk-spac-rigid-spaces-as-berkovich-spaces}), and $Y^{\diamond}$ its associated diamond.   For $Z$ an arc-stack over $\bb{Q}_p$ we let $\FF_Z^\wedge$ denote the formal completion of $\FF_Z$ at $\widehat{Z^{\diamond}}\to \FF_Z$.  The following hold:

\begin{enumerate}

\item  For $Y\in X_{\proet}^{\omega_1}$ the  natural diagram
\[
\begin{tikzcd}
\widehat{Y^{\diamond}}  \ar[r] \ar[d] & \FF_{Y^{\diamond}}^\wedge \times_{Y^{\dR}} Y  \ar[d] \\ 
\widehat{X^{\diamond}} \ar[r] & \FF_{X^{\diamond}}^\wedge \times_{X^{\dR}} X
\end{tikzcd}
\]
is cartesian. 

\item  Suppose that $X$ is affinoid and that  there is an \'etale map $f\colon X\to \bb{T}^d_{\bb{Q}_p}$.
  Let $\bb{T}^{d}_{\bb{Q}_p(\zeta_{p^{\infty}}),\perf}\to \bb{T}^d_{\bb{Q}_p}$ be the uncompleted perfectoid torus and $X_{\infty}\to X$ its pullback along $f$.
  Let $T_1,\ldots, T_n$ denote the coordinates induced by $f$ and set $T_i^{\flat}=(T_i^{1/p^n})_n$.
  Let $\epsilon= (\zeta_{p^n})_n$ be a compatible sequence of  $p$-power roots of unity.
  Then the map 
\begin{equation}\label{eqFormalCompletionOBdRplus}
\widehat{X}_{\infty} \to \FF_{X^{\diamond}_{\infty}}^\wedge \times_{X^{\dR}} X_{\infty}
\end{equation}
is a regular immersion with kernel generated by the ideal $(\log[\epsilon], [T_1^{\flat}]\otimes 1-1\otimes T_1,\ldots, [T_d^{\flat}]\otimes 1-1\otimes T_d )$.  In particular, the formal completion of \eqref{eqFormalCompletionOBdRplus} is precisely the formal spectrum of the ring $\s{O}\mathbb{B}_{\dR,X}^+(X_{\infty})$. 

\item The map $\widehat{X}^{\diamond}\to \FF_{X^{\diamond}}^\wedge \times_{X^{\dR}} X$ is, locally on $X_{\proet}^{\omega_1}$, a regular closed immersion, with formal completion $\Spf  (\s{O}\mathbb{B}_{\dR,X}^+)$.

\end{enumerate}
\end{lemma}
\begin{proof}

\begin{enumerate}

\item  This follows by taking the limit of the following diagram of cartesian  squares
\[
(\widehat{Y^{\diamond}} = \widehat{X^{\diamond}} \times_{\FF_{X^{\diamond}}^\wedge} \FF_{Y^{\diamond}}^\wedge)\to (Y=X\times_{X^{\dR}} Y^{\dR})  \leftarrow (Y=X\times_{X} Y)
\]
where the left cartesian square follows from the fact that $(\widehat{Y}^{\diamond}\to \FF_{Y^{\diamond}})\to (\widehat{X}^{\diamond}\to \FF_{X^{\diamond}})$ is a morphism of degree $1$ divisors in Fargues--Fontaine curves, and the middle square is cartesian as $Y\to X$ is $\dagger$-formally \'etale (as an inverse limit of Berkovich-\'etale maps). 

\item Let $X=\GSpec(A)$ and $X_{\infty}=\GSpec(A_{\infty})$ with completion the perfectoid ring $A_{\infty}^u$.  We let $t=\log [\epsilon]$.

 By definition, the product $\FF_{X_{\infty}^{\diamond}}^\wedge \times X_{\infty}$ is the formal analytic stack given by the colimit 
\[
\FF_{X_{\infty}^{\diamond}}^\wedge \times X_{\infty}=\varinjlim_n \GSpec (\bb{B}^{+}_{\dR}(A^u_{\infty})/t^n\otimes_{\bb{Q}_p} X_{\infty}),
\]
that is, the formal analytic stack provided by the $t$-adic completion  of  $\bb{B}^{+}_{\dR}(A^u_{\infty})\otimes_{\bb{Q}_{p}} A_{\infty}$.   By assumption we have \'etale maps 
\[
\bb{Q}_p\langle T_1,\ldots, T_d\rangle \to A_0\to A_1\to A_2\to \cdots,
\]
and we have a formal lift $\bb{Q}_p\langle T_1,\ldots, T_d\rangle\to \bb{B}_{\dR}^+(A_{\infty})$ given by $T_i\mapsto [T_i^{\flat}]$, and by formal \'etaleness  we can find compatible lifts $A_n\to  \bb{B}_{\dR}^+(A_{\infty})$, this produces a factorization 
\[
 \widehat{X^{\diamond}_{\infty}}\to \FF_{X^{\diamond}_{\infty}}^\wedge \xrightarrow{a} X_{\infty}
\]
whose maps induce equivalences on de Rham stacks.
We have a pullback diagram
\begin{equation}\label{eqCartesianDeRhamFFOB}
  \begin{tikzcd}
{  \FF_{X^{\diamond}_{\infty}}^\wedge  }\ar[r,"{(\ob{id},a)}"] \ar[d,"a"] & {\FF_{X^{\diamond}_{\infty}}^\wedge \times_{X^{\dR}_{\infty}} X_{\infty}  }\ar[d,"{a\times \ob{id}} "] \\ 
  X_{\infty} \ar[r,"\Delta"] &  X_{\infty} \times_{X^{\dR}_{\infty}} X_{\infty} 
  \end{tikzcd}
\end{equation}
Since the map $\widehat{X_{\infty}^{\diamond}}\to \FF_{X^{\diamond}_{\infty}}^\wedge$ is already the regular Zariski closed immersion determined by $t=0$, it suffices to prove that the upper horizontal map of \eqref{eqCartesianDeRhamFFOB} is a regular immersion.
For that, it suffices to show that the lower horizontal map is a regular immersion,  but this follows from the fact that it is the pullback of the regular Zariski closed immersion $X\to X\times_{X^{\dR}} X=(X\subset X\times X)^{\dagger}$ (as $X$ is a smooth rigid space), along the morphism of arrows
\[
(X\to X\times_{X^{\dR}} X)\to (X^{\dR}\to  X^{\dR}) \leftarrow (X_{\infty}^{\dR} \to X_{\infty}^{\dR}),
\]
namely, this is a formal consequence of the fact that $X_{\infty}\to X$ is $\dagger$-formally \'etale, i.e., $X_\infty\cong X_\infty^\dR\times_{X^\dR}X$. 

A bookkeeping of the construction shows  that $(t, [T_1^{\flat}]\otimes 1-1\otimes T_1,\ldots, [T_d^{\flat}]\otimes 1-1\otimes T_d )$ is the desired regular sequence.

\item This follows directly from (2) and (1) as there is a basis of $X^{\omega}_{\proet}$ consisting of perfectoids which are pro-\'etale over perfectoid tori. 
\qedhere

\end{enumerate}

\end{proof}

  \begin{remark}\label{rk:faltings-extension2}
  Let
    $$
    i: \widehat{X^\diamond} \to \Spf(\s{O}\mathbb{B}_{\dR,X}^+),
    $$
      Since $\mathbb{L}_{\FF_{X^\diamond}/\GSpec(\Q_p)}=0$ (by the usual argument that the $p$-completed cotangent complex of $A_{\mathrm{inf}}$ over $\Z_p$ vanishes by perfectness mod $p$, \Cref{sec:cotang-compl-morph-cotangent-complex-via-p-completion} and \cite[Lemma 4.4.9]{drff_stacks} to pass to affinoids on the curve), and $\bb{L}_{X^{\dR}/\GSpec(\Q_p)}=0$, we see from \Cref{LemmaReconstructionOBdR} that
    $$
    \mathbb{L}_i = \mathbb{L}_{\widehat{X^\diamond}/X}.
    $$
    On the other hand, $i$ is a regular  Zariski closed immersion by \Cref{LemmaReconstructionOBdR}, and thus its cotangent complex is the suspension of $\ker(\theta)/\ker(\theta)^2$, where $\ker(\theta)=\ker(\s{O}\mathbb{B}_{\dR,X}^+ \to \widehat{\s{O}}_X)$ is the ideal sheaf defining the closed immersion. Hence,
        $$
    \mathbb{L}_i = \mathrm{gr}^1(\s{O}\mathbb{B}_{\dR,X}^+)[1].
    $$
    Putting things together, we get 
     $$
    \bb{L}_{\widehat{X^\diamond}/X} \cong \mathrm{gr}^1(\s{O}\mathbb{B}_{\dR,X}^+)[1].
    $$
   recovering a result of Guo--Li (\cite[Remark 4.22]{guo2021period})). 
  \end{remark}

\begin{construction}\label{ConstructionNfilteredInsideZfiltered}
  We can think of decreasing $\bb{N}$-filtered modules as the full subcategory of decreasing $\bb{Z}$-filtered modules of filtered objects $\Fil^i$ (equivalently sheaves on $\bb{A}^{1,\ob{alg}}/\bb{G}_m^{\ob{alg}}$) such that $\Fil^{0}\to \Fil^{-i}$ is an equivalence for all $i\geq 0$; moreover, this identification is also symmetric monoidal.
  Therefore, we can consider $\bb{B}_{\dR,X}^+$ both as an $\bb{N}$-filtered and  $\bb{Z}$-filtered ring.
  The description as an $\bb{N}$-filtered ring is equivalent to the completion at $\widehat{\bb{A}}^1/\bb{G}_m^{\ob{alg}}$ of the deformation to the normal bundle of $\theta\colon \bb{B}_{\dR,X}^+\to \widehat{\s{O}}_X$; this follows for example from the discussion of \cite[Remark 7.7.17]{solidnotes}. We let $\Spf_{\geq 0} \bb{B}_{\dR,X}^+$ (resp. $\Spf_{\geq 0}(\s{O}\mathbb{B}_{\dR,X}^+)$) denote the spectrum of the $\bb{N}$-filtration of $\bb{B}_{\dR,X}^+$  (resp. $\s{O}\mathbb{B}_{\dR,X}^+$) considered as a sheaf over $\FFdisc{X^{\diamond}}^{\times}\times \widehat{\bb{A}}^1/\bb{G}_m^{\ob{alg}}$ (resp. over $\FF_{X^{\diamond}}^{\wedge}\times_{X^{\dR}}X \times \widehat{\A}^1/\Gm^{\ob{alg}}$).
  The description as a $\bb{Z}$-filtered ring is equivalent to the datum of the Cartier divisor $\ker \theta\to \bb{B}_{\dR,X}^+$, this corresponds to the map $\Spf( \bb{B}_{\dR,X}^+)\to \widehat{\bb{A}}^1/\bb{G}_m^{\ob{alg}}$.  In particular, we have  a natural map  over $\FFdisc{X^{\diamond}}^{\times}\times \widehat{\bb{A}}^1/\bb{G}_m^{\ob{alg}}$
\[
\Spf (\bb{B}_{\dR,X}^+)\to  \Spf_{\geq 0} (\bb{B}_{\dR,X}^+). 
\]

Let $\s{O}\mathbb{B}_{\dR,X}$ be the  completion of the convolution filtration  of $\s{O}\mathbb{B}_{\dR,X}^+[\frac{1}{\xi}]$,  we define the formal spectrum $\Spf  (\Fil^{\bullet} \s{O}\mathbb{B}_{\dR,X})$   seen as a stack over $\widehat{\bb{A}}^1/\bb{G}_m^{\ob{alg}}$  to be 
\[
\Spf (\Fil^{\bullet}\s{O}\mathbb{B}_{\dR,X}):=  \Spf (\bb{B}_{\dR,X}^+)\times_{\Spf_{\geq 0} (\bb{B}_{\dR,X}^+)}  \Spf_{\geq 0} (\s{O}\mathbb{B}_{\dR,X}^+).
\]
\end{construction}

\begin{remark}\label{RemarkRefinementConeConstruction}
Let $X/\bb{Q}_p$ be a smooth rigid space as above, and let $\widehat{X^{\diamond}}\to X$ be the natural map.  Consider the morphism of normal cones
\[
\Spf (\bb{B}_{\dR,X}^+) \to  \widehat{X^{\diamond}}^{\widehat{\cone}} \to X^{\widehat{\cone}}\to X^{\widehat{\ob{cone-alg}}},
\]
where the first map arises from the Cartier divisor $\widehat{X^{\diamond}}\to \Spf (\bb{B}_{\dR,X}^+)$. 

By \Cref{RemarkAlgebraicConeGeneralAnalyticStacks} the map $X^{\widehat{\ob{cone-alg}}}\to |X|_{\Betti}\times \widehat{\bb{A}}^1/\bb{G}_m^{\ob{alg}}$  is the relative spectrum of the de Rham complex $(\Omega^{\bullet}_{X},d)$ seen as a cdga. In particular it is $\bb{N}$-filtered, which yields a factorization 
\[
\Spf( \bb{B}_{\dR,X}^+) \to \Spf_{\geq 0} (\bb{B}_{\dR,X}^+ ) \to  X^{\widehat{\ob{cone-alg}}}. 
\]
Geometrically, the previous diagram arises from the natural morphisms
\[
\Spf( \bb{B}_{\dR,X}^+)\to   (\widehat{X^{\diamond}})^{\widehat{\ob{cone-alg}}/\Spf (\bb{B}^{+}_{\dR,X})}  \to  X^{\widehat{\ob{cone-alg}}}
\]
(using  \cite[Remark 7.7.17]{solidnotes} for the middle term, see \Cref{RemarkAlgebraicConeGeneralAnalyticStacks} (2)).
\end{remark}

  The following result is stronger than what we need for the comparison to the Faltings' extension. It can be interpreted as a geometric incarnation of \cite[Theorem 1.1, Remark 4.14]{guo2021period}, where one should replace relative Hodge--completed de Rham complexes by the formal cone construction (see \Cref{RemarkAlgebraicConeGeneralAnalyticStacks}).  A related result has been proved independently by Jiang (see \cite[\S 2.10]{jiang2026differential}).

 \begin{proposition} \label{prop:fiber-product-in-terms-of-faltings-extension}
 Let $X$ be a smooth rigid space over $\Q_p$ and $X^{\widehat{\ob{cone-alg}}}\to \widehat{\bb{A}}^1/\bb{G}_m^{\ob{alg}}$ its algebraic cone over $\bb{Q}_p$. Then the natural map  
 \[
 \Spf_{\geq 0} (\bb{B}_{\dR,X}^+ ) \to  \widehat{X^{\diamond}}^{\widehat{\ob{cone-alg}}}
 \]
 is an equivalence.

  Furthermore,  there is a natural  equivalence of analytic stacks over $\widehat{\bb{A}}^1/\bb{G}_m^{\ob{alg}}$
 \begin{equation}\label{BaseChangeOBdRplus}
 \Spf_{\geq 0}( \bb{B}_{\dR,X}^+) \times_{X^{\widehat{\ob{cone-alg}}}} (X\times \widehat{\bb{A}}^1/\bb{G}_m^{\ob{alg}}) \cong  \Spf_{\geq 0} (\s{O}\mathbb{B}_{\dR,X}^+).
 \end{equation}
 Similarly, we have an equivalence 
 \begin{equation} \label{BaseChangeOBdR}
 \Spf(\bb{B}_{\dR,X}^+) \times_{X^{\widehat{\ob{cone-alg}}}} (X\times \widehat{\bb{A}}^1/\bb{G}_m^{\ob{alg}}) \cong  \Spf 
 (\Fil^{\bullet} \s{O}\mathbb{B}_{\dR,X}).
 \end{equation}
\end{proposition}

 \begin{proof}  It is clear that \eqref{BaseChangeOBdRplus} implies \eqref{BaseChangeOBdR}.  We recall from \Cref{RemarkAlgebraicConeGeneralAnalyticStacks} (1) that if $W\to V$ is a formally \'etale map of Gelfand stacks, then the natural map of Gelfand stacks 
\[
W\times \widehat{\bb{A}}^1/\bb{G}_m^{\ob{alg}}\to W^{\widehat{\ob{cone-alg}}/ V} 
\]
is an equivalence.
In particular, this applies to the map $\FF_{X^{\diamond}}^{\wedge} \times_{X^{\dR}} X \to X$ (note that $X^{\dR}\to \GSpec (\bb{Q}_p)$ and $\FF_{X^{\diamond}}^{\wedge}\to \GSpec (\bb{Q}_p)$ are formally \'etale).
Thus, we deduce that the following commutative diagram 
\begin{equation}\label{eqCommutativeDiagramConesObdR1}
\begin{tikzcd}
(\FF_{X^{\diamond}}^{\wedge} \times_{X^{\dR}} X )\times  \widehat{\bb{A}}^1/\bb{G}_m^{\ob{alg}} \ar[r] \ar[d]&  X\times  \widehat{\bb{A}}^1/\bb{G}_m^{\ob{alg}}     \ar[d] \\ 
(\FF_{X^{\diamond}}^{\wedge} \times_{X^{\dR}} X )^{\widehat{\ob{cone-alg}}}  \ar[r]  &  X^{\widehat{\ob{cone-alg}}}
\end{tikzcd}
\end{equation}
is cartesian.
On the other hand, by \Cref{LemmaReconstructionOBdR} the stack $\Spf (\s{O}\bb{B}_{\dR,X}^+)$ is obtained as the formal neighborhood of the Zariski closed immersion $g\colon \widehat{X}^{\diamond}\to \FF_{X^{\diamond}}^{\wedge} \times_{X^{\dR}} X$. Moreover, by \Cref{RemarkAlgebraicConeGeneralAnalyticStacks} (2), the relative cone of  $g$ is the spectrum of the adic filtration of $ \s{O}\bb{B}_{\dR,X}^+$, that is $\Spf_{\geq 0}(\s{O}\bb{B}_{\dR,X}^+)$. We deduce a pullback square 
\begin{equation}\label{eqCommutativeDiagramConesObdR2}
\begin{tikzcd}
\Spf_{\geq 0}(\s{O}\bb{B}_{\dR,X}^+) \ar[r] \ar[d] & (\FF_{X^{\diamond}}^{\wedge} \times_{X^{\dR}} X )\times  \widehat{\bb{A}}^1/\bb{G}_m^{\ob{alg}}   \ar[d] \\ 
\widehat{X^{\diamond}}^{\widehat{\ob{cone-alg}}}  \ar[r]&  (\FF_{X^{\diamond}}^{\wedge} \times_{X^{\dR}} X )^{\widehat{\ob{cone-alg}}}
\end{tikzcd}
\end{equation}
Combining \eqref{eqCommutativeDiagramConesObdR1} and \eqref{eqCommutativeDiagramConesObdR2}, we deduce a cartesian square
\[
\begin{tikzcd}
\Spf_{\geq 0}(\s{O}\bb{B}_{\dR,X}^+)\ar[r]  \ar[d] &   X\times  \widehat{\bb{A}}^1/\bb{G}_m^{\ob{alg}} \ar[d] \\
\widehat{X^{\diamond}}^{\widehat{\ob{cone-alg}} } \ar[r]& X^{\widehat{\ob{cone-alg}}} .
\end{tikzcd}
\]

Therefore, to conclude the proof of the proposition it suffices to show that the natural map
   \[
    \Spf_{\geq 0} (\bb{B}_{\dR,X}^+ ) \to  \widehat{X^\diamond}^{\widehat{\ob{cone-alg}}}
   \]
is an equivalence. For that, we notice that since $\FF^{\wedge}_{X^{\diamond}}\to \GSpec (\bb{Q}_p)$ is formally \'etale, the natural map of formal cones 
\[
\widehat{X^{\diamond}}^{\widehat{\ob{cone-alg}}/\FF^{\wedge}_{X^{\diamond}}}\to \widehat{X^{\diamond}}^{\widehat{\ob{cone-alg}}}
\]
is an equivalence, but by \Cref{RemarkAlgebraicConeGeneralAnalyticStacks} (2) the former is precisely the spectrum of the filtered ring given by the adic filtration of $\bb{B}_{\dR,X}^+$, namely, the spectrum $\Spf_{\geq 0}(\bb{B}_{\dR,X}^+)$ as desired. 
\end{proof}  

  \begin{remark}
  As a side remark, note that the argument proving \Cref{prop:fiber-product-in-terms-of-faltings-extension} gives the following description of Fontaine's period ring $B_\dR^+$, seen as an $\N$-filtered ring with filtration given by powers of $\ker(\theta)$: it is the pushforward of $\s{O}$ along $\GSpec(\C_p)^{\widehat{\cone-\mathrm{alg}}}\to \widehat{\A^1}/\Gm^{\mathrm{alg}}$.
  \end{remark}

\begin{corollary} \label{cor:identification-fe}
Let $X$ be a smooth rigid space over $\Q_p$. There is a natural equivalence of Gelfand stacks 
\[ \widehat{X^\diamond}\times_{X^{\cone}_s}(X\times B\Gm) \cong (\AnSpec(\s{O}\mathbb{C}_X))_{\mathrm{Gelf}},
\]
where $(-)_{\mathrm{Gelf}}$ is the ``Gelfandification'' functor of \cite[Section 4]{drff_stacks}. In other words, the Faltings' extension of $X$ as defined in \Cref{sec:p-adic-twistor-1-definition-faltings-extension} coincides with the usual Faltings' extension.
\end{corollary}
\begin{proof}
Base changing along $B\bb{G}_m^{\ob{alg}} \to \widehat{\bb{A}}^1/\bb{G}_m^{\ob{alg}}$ the equivalence \eqref{BaseChangeOBdR} of analytic stacks over $\widehat{\bb{A}}^1/\bb{G}_m^{\ob{alg}}$ 
 \[ \Spf(\bb{B}_{\dR,X}^+) \times_{X^{\widehat{\ob{cone-alg}}}} X\times \widehat{\bb{A}}^1/\bb{G}_m^{\ob{alg}} \cong  \Spf 
 (\Fil^{\bullet} \s{O}\mathbb{B}_{\dR,X})
\]
one gets an equivalence of stacks over $B\Gm^{\mathrm{alg}}$
\[ 
\widehat{X^\diamond}\times_{X^{\cone}_s}(X\times B\Gm^{\mathrm{alg}}) \cong \Spf 
 (\Fil^{\bullet} \s{O}\mathbb{B}_{\dR,X}) \times_{\widehat{\bb{A}}^1/\bb{G}_m^{\ob{alg}}} B\bb{G}_m^{\ob{alg}}.
\]
The stack $\Spf (\Fil^{\bullet} \s{O}\mathbb{B}_{\dR,X}) \times_{\widehat{\bb{A}}^1/\bb{G}_m^{\ob{alg}}} B\bb{G}_m^{\ob{alg}}$ as a stack over $B\Gm^{\mathrm{alg}}$ corresponds to the graded ring $\mathrm{gr}^{\bullet} \s{O}\mathbb{B}_{\dR,X}$. Since $\mathrm{gr}^{\bullet} \s{O}\mathbb{B}_{\dR,X}= \oplus_i \mathrm{gr}^{0} \s{O}\mathbb{B}_{\dR,X}(i)$, we have
\[
\Spf (\Fil^{\bullet} \s{O}\mathbb{B}_{\dR,X}) \times_{\widehat{\bb{A}}^1/\bb{G}_m^{\ob{alg}}} B\bb{G}_m^{\ob{alg}} \cong \AnSpec(\mathrm{gr}^{0} \s{O}\mathbb{B}_{\dR,X})= \AnSpec(\s{O}\mathbb{C}_X).
\]
One concludes the statement of the corollary by applying the functor $(-)_{\mathrm{Gelf}}$.
 \end{proof}

  With these preliminaries out of the way, we specialize \Cref{PropOverconvergentHT,PropRelativeVersionOverconvergentHTStacks} to  the promised (almost tautological) description of the formal (and overconvergent) the Hodge--Tate stack for smooth rigid spaces, first in the absolute case.

  \begin{corollary}\label{restriction-ht-overconvergent}
    Let $X$ be a smooth rigid space over $\bb{Q}_p$.
    There is a cartesian square
\[\begin{tikzcd}
	{X^{\widehat{\HT}}} & {X^\dR} \\
	{X^{\widehat{\cone}}} & {(X^{\cone}_s)^\dR = (X^{\cone})^\dR \times_{(\A^1/\Gm)^\dR} B\Gm^\dR,}
	\arrow[from=1-1, to=1-2]
	\arrow[from=1-1, to=2-1]
	\arrow[from=1-2, to=2-2]
	\arrow[from=2-1, to=2-2]
      \end{tikzcd}\]
    where  the bottom right term naturally identifies with the classifying stack $B_{X^\dR \times B\Gm^\dR} \bb{T}_X\{1\}^\dR$, the map $X^{\dR}\to B\bb{G}_m^{\dR}$ arising from the second projection of the right vertical map corresponds to the Tate twist,  and the $\bb{T}_X\{1\}^\diamond$-torsor on $X^\diamond$ corresponding to the right vertical map\footnote{Recall that if $Y, Z$ are Gelfand stacks, a morphism $Y^\dR \to Z^\dR= (Z^\diamond)^\dR$ is the same thing by adjunction as a morphism  $(Y^\dR)^\diamond= Y^\diamond \to Z^\diamond$ of arc-stacks over $\Q_p$.} is the Faltings' extension $\mathcal{L}_{\mathrm{FE}_X}$.
  \end{corollary}
  If in addition $X$ is a Berkovich space, e.g., $X$ is partially proper or affinoid, then the formal completion along the Hodge-Tate stack can be replaced by the overconvergent neighborhood of it.
\begin{proof}
  As $X_s^{\cone}\to X^{\widehat{\cone}}$ is an epimorphism, it suffices to show the same for the Hodge-Tate stack.
  By descent along formally \'etale $!$-covers, e.g., Tate-rational localizations, as will be proven in \Cref{sec:desc-analyt-prism-descent-for-hodge-tate-stack}, we can then reduce to the case that $X$ is affinoid, and thus $\dagger$-cohesive.
  Then the existence of the cartesian diagram is clear.
\end{proof}

Next, we state the relative variant.

 \begin{corollary}\label{restriction-ht-overconvergent-relative}
Let $D$ be a perfectoid space over $\Q_p$ and $X$ be a smoothoid space over $D$. Let $D^{\wedge}\subset \FF_D$ denote the formal neighbourhood of $D$ in $\FF_D$.  There is a cartesian square
\[\begin{tikzcd}
	{X^{\widehat{\HT} /D^{\wedge}}} & {X^{\dR/D^{\wedge}}} \\
	{X^{\widehat{\cone}/(D,\FF_D)}} & {(X^{\cone/(D,\FF_D)}_s)^{\dR/D^{\wedge}}.}
	\arrow[from=1-1, to=1-2]
	\arrow[from=1-1, to=2-1]
	\arrow[from=1-2, to=2-2]
	\arrow[from=2-1, to=2-2]
      \end{tikzcd}\]
      Assume that $X$ has a smooth lift $\widetilde{X}$ to $B_\dR^+(D)/t^2$, depicted by a map $X\to X_{s}^{\cone/(D,\FF_D)}$. Then $\widetilde{X}$ determines an isomorphism of the bottom right term with the classifying stack $B_{X^{\dR/D^{\wedge}}} \bb{T}_{X/D}\{1\}^{\dR/D^{\wedge}}$ and the $\bb{T}_{X/D}\{1\}^\diamond$-torsor on $X^\diamond$ corresponding to the right vertical map through this isomorphism is the Faltings' extension $\mathcal{L}_{\mathrm{FE}_X^{\widetilde{X}}}$.
    \end{corollary}
    Again, if $X$ is in addition a Berkovich space, the analog holds for the overconvergent neighborhood (replacing relative de Rham stacks over $D^{\wedge}$ by relative de Rham stacks over the overconvergent neighborhood $D^\dagger$ of $D$ in $\FF_D$).
    \begin{proof}
      The proof is similar to the one of \Cref{restriction-ht-overconvergent}.
    \end{proof}

Both corollaries follow immediately from the definition of the Faltings' extension in \Cref{sec:p-adic-twistor-1-definition-faltings-extension}.

   \begin{remark}[Overconvergent Hodge--Tate stack and filtered de Rham stack] \label{rk:overconvergent-ht-stack-vs-filtered-dr-stack}
     Going back to the analogy with the filtered de Rham stack mentioned in \Cref{sec:filtered-de-rham-stack}, \Cref{restriction-ht-overconvergent} can be compared to the analogous description for the pullback of the filtered de Rham stack of $X$ to $\A^{1,\dagger}/\Gm$.
     Namely, we have a cartesian square
      \[\begin{tikzcd}
	{X^{\dagger-\dR, +}} & {X^\dR} \\
	{X^{\dagger-\cone}} & {(X^{\cone})^\dR \times_{(\A^1/\Gm)^\dR} B\Gm^\dR,}
	\arrow[from=1-1, to=1-2]
	\arrow[from=1-1, to=2-1]
	\arrow[from=1-2, to=2-2]
	\arrow[from=2-1, to=2-2]
      \end{tikzcd}\]     
      where now the right vertical map is the unit section of the classifying stack, corresponding to the trivial $\bb{T}_X(-1)$-torsor (where the $(-1)$ refers to a twist by the universal ideal of $\bb{A}^1/\bb{G}_m$, which agrees with the $\{1\}$-twist on the Hodge--Tate stack). Hence as a consequence of \Cref{restriction-ht-overconvergent} we see that if $\widetilde{X}\to X$ is a formally \'etale map such that the Faltings' extension trivializes after pullback to $\widetilde{X}^\diamond$, then the choice of a trivialization of the Faltings'  extension  produces an isomorphism  $\widetilde{X}^{\dagger-\HT}\cong \widetilde{X}^{\dagger-\dR,+}$. This applies for example when $\widetilde{X}\to X$ is $\GSpec(\Q_p(\mu_{p^\infty}))\to \GSpec(\Q_p)$ or $\varprojlim_{x\mapsto x^p} \bb{G}_{m,D}\to \bb{G}_{m,D}$ with $D$ perfectoid. 
\end{remark}

\begin{remark}
The left vertical map in \Cref{restriction-ht-overconvergent}, from $X^{\dagger-\HT}$ to $X^{\dagger-\cone}$, or rather a formal variant of it, will be used again in \Cref{sec:relation-to-rh}: we will see, based on the results established in this subsection, that pullback along this map recovers on vector bundles the $p$-adic Riemann--Hilbert functor of \cite{scholze2013p}.

We will also see the Faltings' extension reappear in \Cref{sec:mixed-case}, when discussing the cotangent complex of the Hodge--Tate stack. 
\end{remark}

\subsection{Descent for the analytic prismatization}
\label{subsec:desc-perf-analyt}

In this section we study more seriously the descent properties of the analytic prismatization and the  Hodge--Tate stack.
The first kind of descent we are able to prove involves covers which are $\dagger$-formally smooth, cf.\ \Cref{sec:prop-morph-gelf-1-definition-dagger-formally-etale}.  

More precisely, we have the following theorem.

\begin{theorem}\label{TheoremDescentPrismatization}
  Let $f\colon Y\to X$ be a map of Gelfand stacks.
  \begin{enumerate}
  \item If $f$ is $\dagger$-formally \'etale locally $\dagger$-formally smooth, then the maps $f^{\prism}\colon Y^{\prism}\to X^{\prism}$, $f^{\widetilde{\Div^1}}\colon Y^{\widetilde{\Div^1}}\to X^{\widetilde{\Div^1}}$ and $f^{\Div^1}\colon Y^{\Div^1}\to X^{\Div^1}$ are weakly $\dagger$-formally smooth.
  \item If in addition $f$ is an arc-cover, then $f^\prism, f^{\widetilde{\Div^1}}, f^{\Div^1}$ are epimorphisms.
  \item If $f$ is $\dagger$-formally smooth and $Y,X$ are $\dagger$-cohesive, then $f^\prism, f^{\widetilde{\Div^1}}, f^{\Div^1}$ are $\dagger$-formally smooth.
  \end{enumerate}
  \end{theorem}
\begin{proof}
  We give the argument for the map $Y^{\prism}\to X^{\prism}$, the other cases being totally analogous.
  Let $A$ be a nilperfectoid ring with $T=\GSpec(A)$, and let $D\to \n{X}_T$ a  degree $1$ divisor on the thickened Fargues--Fontaine disc of $T$.
  Then the space $D$ is an affinoid nilperfectoid  space, with reduction $\overline{D}$, and the map 
\[
Y(D)\to Y(\overline{D})\times_{X(\overline{D})} X(D)
\]
is an epimorphism if $Y\to X$ is $\dagger$-formally smooth.
Passing to sheafifications for the $!$-topology this implies that $Y^{\prism}\to Y^{\prism,\dR}\times_{X^{\prism,\dR}} X^{\prism}$ is an epimorphism in this case, i.e., weakly $\dagger$-formally smooth.
In case both $X$ and $Y$ are in addition $\dagger$-cohesive, the sheafification is not necessary and $f^{\prism}$ is already $\dagger$-formally smooth. 
 By \Cref{sec:defin-first-prop-5-prismatization-for-dagger-formally-etale-maps},  we can use $\dagger$-formally \'etale localizations to reduce to the case that $f$ is $\dagger$-formally smooth.

To prove the second claim about descent, it suffices to note that since $Y\to X$ is an arc-cover the map $Y^{\prism,\dR}\to X^{\prism,\dR}$ is an epimorphism, and hence it suffices to prove that the natural map
\[
Y^{\prism}\to Y^{\prism,\dR}\times_{ X^{\prism,\dR}} X^{\prism}
\]
is an epimorphism, which is exactly implied by the first point.
\end{proof}

For the Hodge-Tate stack, we get a stronger descent statement.

\begin{theorem}\label{sec:desc-analyt-prism-descent-for-hodge-tate-stack}
  Let $f\colon Y\to X$ be an epimorphism of Gelfand stacks which is  infinitesimally cohesive  and admits a cotangent complex with $\bb{L}_{Y/X}$ a perfect module with amplitude $[-1,0]$ (i.e.\ $f$ is quasi-syntomic with perfect cotangent complex). Then $f^{\HT}\colon Y^{\HT}\to X^{\HT}$ is an epimorphism. 
\end{theorem}
\begin{proof}
Let $A$ be a nilperfectoid $\Q_p$-algebra. Let $g: \GSpec(\mathbb{G}_a^\HT(A))\to X$ be given. We want to lift $g$ along $f$, up to possibly replacing $A$ by a descendable cover.  The ring $\mathbb{G}_a^\HT(A)$ is a square-zero extension of $A$ by $\mathrm{Nil}^\dagger(A)\{1\}[1]$,  thus it fits in a pullback square
\[
\begin{tikzcd}
\bb{G}_a^{\HT}(A)  \ar[r] \ar[d] & A \ar[d,"d"] \\
A \ar[r,"0"] & A\oplus \mathrm{Nil}^\dagger(A)\{1\}[2], 
\end{tikzcd}
\]
where $0$ is the trivial derivation and $d$ the derivation that realizes $\bb{G}_{a}^{\HT}(A)$. Since the map $f\colon Y\to X$ is infinitesimally cohesive, the natural commutative square
\[
\begin{tikzcd}
Y(\bb{G}_a^{\HT}(A)) \ar[r] \ar[d]& X(\bb{G}_a^{\HT}(A))  \ar[d]\\ 
 Y(A)\times_{Y(A\oplus \mathrm{Nil}^\dagger(A)\{1\}[2] )} Y(A)  \ar[r] &  X(A)\times_{X(A\oplus \mathrm{Nil}^\dagger(A)\{1\}[2] )} X(A)
\end{tikzcd}
\]
is cartesian.  Therefore, to produce a lift $\GSpec(A) \dashrightarrow Y^{\HT}$ it suffices to find a lift to $Y(A)\times_{Y(A\oplus \mathrm{Nil}^\dagger(A)\{1\}[2] )} Y(A) $ compatible with the map $g$. 

Since $f$ is assumed to be a $!$-cover, we can already lift the composition $\overline{g}\colon \GSpec(A) \to \mathbb{G}_a^\HT(A) \xto{g} X$ to $h: \GSpec(A)\to Y$, up to replacing $A$ by a descendable cover. Consider the commutative diagram 
\[
\begin{tikzcd} 
\ob{fib}(r)  \ar[r]    & \ob{fib}(b) \ar[d] \ar[r,"r"]  & \ob{fib}(a)  \ar[d]\\
  & Y(A) \times_{Y(A\oplus \mathrm{Nil}^\dagger(A)\{1\}[2] ) } Y(A)  \ar[d,"b"]  \ar[r] &   X(A) \times_{X(A\oplus \mathrm{Nil}^\dagger(A)\{1\}[2] )} X(A) \ar[d,"a"] \\ 
& Y(A) \ar[r] & X(A), 
\end{tikzcd}
\]
where the fiber of $a$ is at  the point $\overline{g}$, the fiber at $b$ is  at the point $h$, the fiber at $r$  is at the point $x$ induced by $g$. We want to see that the map $r$ is surjective, for that it suffices to see that $\ob{fib}(r)$ is non-empty.  But then, in the terminology of \cite[Definition 7.5]{MR4176782} and the proof of \Cref{PropFinalDescentDeformationTheory}, the map $r$ is the morphism of anima of paths from $0\rightsquigarrow d$
\[
\Omega_{0,d}\ob{Der}_{Y}(A,\mathrm{Nil}^\dagger(A)\{1\}[2])\to \Omega_{0,d} \ob{Der}_X(A,\mathrm{Nil}^\dagger(A)\{1\}[2] ).
\]
Thus,  we have that 
\[
\ob{fib}(r)= \Omega_{0,d} \ob{Der}_{Y/X}(A,\mathrm{Nil}^\dagger(A)\{1\}[2] ) = \Omega_{0,d}  \tau_{\geq 0} \Hom_{A}(\overline{g}^*\bb{L}_{Y/X}, \mathrm{Nil}^\dagger(A)\{1\}[2]). 
\]
But now, $\overline{g}^*\bb{L}_{Y/X}$ is perfect of homological perfect amplitude $[-1,0]$, and this implies that the anima $\tau_{\geq 0} \Hom_{A}(\overline{g}^*\bb{L}_{Y/X}, \mathrm{Nil}^\dagger(A)\{1\}[2])$ is connective and that there is therefore a path from $0$ to $d$, proving that $\ob{fib}(r)\neq \emptyset$ as wanted.  
\end{proof}

\begin{remark}
We note that the above proof shows quasi-syntomic descent more generally for the transmutation of any ring stack which is a square-zero extension of $\mathbb{G}_a$ by a $1$-connective module. 
\end{remark}

\begin{remark}\label{RemarkDescentFiltereddeRhamStacks}
  Essentially the same arguments in this subsection can also be applied to the filtered de Rham stack obtaining similar descent results.
  Indeed, following \Cref{sec:filtered-de-rham-stack} the definitions of the filtered de Rham stack for a $\dagger$-cohesive Gelfand stack $Y^{\dR,+}\cong Y^{\cone}\times_{Y^{\cone,\dR}} (Y^\dR\times \A^{1,\dR}/\Gm^\dR)$ and of the perfect analytic prismatization modulo Frobenius $Y^{\Div^1}\cong Y^\cone\times_{Y^{\cone,\dR}} Y^{\HK}$ are very similar - for the analytic prismatization one only replaces $Y^{\dR}\times \A^{1,\dR}/\Gm^\dR$ by $Y^{\HK}$.
\end{remark}

\subsection{The analytic prismatization of $\Q_p$}
\label{Section:Analytic-Prism}

The first concrete example for the analytic prismatization to consider is of course the example $X=\GSpec(\Q_p)$. We (obviously) have $\Q_p^\prism=\GSpec(\Q_p)^\prism$ and we write $\widetilde{\Div^1}=\GSpec(\Q_p)^{\widetilde{\Div^1}}$, $\Div^1=\GSpec(\Q_p)^{\Div^1}$.

Let $1+\widetilde{\bb{D}}_{\bb{Q}_p}^{<1}:=\varprojlim_{x\mapsto x^p} (1+\bb{D}^{<1}_{\bb{Q}_p})$ denote the perfect analytic $\delta$-stack given by the uncompleted perfectoid disc  with coordinate $q$ and Frobenius $\varphi(q)=q^p$. We consider the cyclotomic prism $(S_{\cyc},\n{X}_{\cyc},I):=( \Marc(\bb{F}_p((q^{1/p^{\infty}-1}))),(1+\widetilde{\bb{D}}_{\bb{Q}_p}^{<1}) \backslash \{1\},I)$ of  \Cref{ExamplesRelativePrismatizations} (1) with divisor generated by $[p]_q=\frac{q^p-1}{q-1}$.

\begin{proposition}\label{prop:explicit-div-1}
\begin{enumerate}
 \item Let $\mu_{\BK}\colon \bb{Q}_p^{\prism}\to \bb{A}^1/\bb{G}_m$ be the Breuil--Kisin ideal of \Cref{sec:breuil-kisin-twist-2-definition-bk-twist}.    The natural diagram
\[\begin{tikzcd}
	{\Q_p^\prism} & {\Q_p^{\prism,\dR}} \\
	{\A^1/\Gm} & {\A^{1,\dR}/\Gm^\dR}
	\arrow[ from=1-1, to=1-2]
	\arrow["\mu_{\BK}"', from=1-1, to=2-1]
	\arrow["{\mu_{\BK}^\dR}", from=1-2, to=2-2]
	\arrow[from=2-1, to=2-2]
      \end{tikzcd}\]
    is a fiber square.
  \item
    The diagram
\[\begin{tikzcd}
	{\Q_p^\HT} & {\GSpec(\Q_p)} \\
	{B\Gm} & {B\Gm^\dR}
	\arrow[from=1-1, to=1-2]
	\arrow["\iota", from=1-1, to=2-1]
	\arrow["\iota^\dR", from=1-2, to=2-2]
	\arrow[from=2-1, to=2-2]
      \end{tikzcd}\]
    is a fiber square. 
    In particular, the natural map
    \[
\pi_{\Q_p}\colon \Q_p^\HT\to \GSpec(\Q_p),
\]
is a $B\Gm^{\dagger}$-torsor and therefore a gerbe banded by $\Gm^\dagger$.

\item Consider the morphism $\n{X}_{\cyc}\to \bb{Q}_p^{\prism}$ produced from the cyclotomic prism (cf. \Cref{rk:explicit-description}). Then the pullback of the Breuil--Kisin divisor $I_{\BK}\to \s{O}_{\bb{Q}_p^{\prism}}$ is the ideal in $\n{X}_{\cyc}$ generated by the function  $\log(q)/(q-1)$.  Similarly, the pullback of $\n{I}_{\prism}$ along the composite $\n{X}_{\cyc}\to \widetilde{\Div^1}\to \Div^1$ is generated by the function $\log (q)$.

\item The morphism 
$$
 \mathcal{X}_{\rm cyc} \to \widetilde{\Div^1}
$$
from the cyclotomic prism (cf. \Cref{ExamplesRelativePrismatizations}  \Cref{rk:explicit-description})  induces an isomorphism
$$
\widetilde{\Div^1} \cong \mathcal{X}_{\rm cyc}/\Z_p^{\times, \la},
$$
with $\Z_p^{\times, \la}$ acting by $a.(q,q^{1/p},\dots) \mapsto (q^a, q^{a/p},\dots)$.  In particular, we have isomorphisms
$$
\Div^1 \cong \mathcal{X}_{\rm cyc} /\Q_p^{\times, \la}
$$
and
  \[
    \Q_p^\HT\cong \GSpec(\Q_p(\mu_{p^\infty}))/\Z_p^{\times,\la},
  \]
  where the action of $\Z_p^{\times,\la}$ on $\GSpec(\Q_p(\mu_{p^{\infty}}))$ factors through the natural Galois action of $\Z_p^{\times,\sm}\cong \Gal(\Q_p(\mu_{p^{\infty}})/\Q_p)^{\sm}$ provided by the cyclotomic character.  

The cartesian square of \Cref{prop:pullback-description-twistor-stack} becomes in this case:
\[\begin{tikzcd}
	{((1+\widetilde{\bb{D}}_{\bb{Q}_p}^{<1})\backslash \{1\} )/\Q_p^{\times, \la}} & {((1+\widetilde{\bb{D}}_{\bb{Q}_p}^{<1})^{\dR}\backslash \{1\} )/\Q_p^{\times, \sm}} \\
	{\A^1/\Gm} & {\A^{1,\dR}/\Gm^\dR}
	\arrow[ from=1-1, to=1-2]
	\arrow["\log"', from=1-1, to=2-1]
	\arrow["{\log^\dR}", from=1-2, to=2-2]
	\arrow[from=2-1, to=2-2]
      \end{tikzcd}\]
      \end{enumerate}
\end{proposition}
\begin{proof}
Points (1) and (2) are simply an explicit rewriting of \Cref{prop:pullback-description-twistor-stack}.

For (3), we note that by definition the pullback of the ideal $\n{I}_{\prism}$ to $\n{X}_{\cyc}$ is  generated by $[p]_q$. Hence, we only need to see that the function $\frac{\log(q)}{q-1}$ has simple zeros at $[p]_{q^{p^n}}=0$ for $n\geq 0$. It suffices to show this for the open unit disc $1+ \bb{D}^{<1}_{\bb{Q}_p}$ with coordinate $q$.
In that case $\log(q)$ has simple roots at $p$-power roots of unity $\mu_{p^{\infty}}$, while the function $[p]_{q^n}$ has simple roots at primitive $p^{n+1}$-th roots of unity.
This implies that $\frac{\log(q)}{q-1}$ has the same zeros at $\{[p]_{q^n}: n\geq 0\}$ as desired.   For the claim about the pullback of $\n{I}_{\prism}$ along the composite $\n{X}_{\cyc}\to \widetilde{\Div^1}\to \Div^1$, we note that this ideal coincides with the ideal obtained as the perfection of $\log(q)/(q-1)$. Using that Frobenius acts totally discontinuously on $\n{X}_{\cyc}$, this ideal is the limit of the ideal generated by the function $\lim_n \log(q)/(q^{p^{1/n}}-1) $ as $n\to\infty$.
The zeros of the denominator approach more and more to the boundary of the disc, which makes the zeros of the limit equal to the zeros of $\log(q)$, proving what we wanted. 

For (4), one way to argue is to use \Cref{prop:pullback-description-twistor-stack} and the fact that $\Q_p^\HK= (1+\widetilde{\bb{D}}_{\bb{Q}_p}^{<1})^\dR$ (which follows from the same statement on arc-stacks). This gives rise to a map $((1+\widetilde{\bb{D}}_{\bb{Q}_p}^{<1})\backslash \{1\} )/\Q_p^{\times, \la} \to \Div^1$, which is an isomorphism (since the logarithm map is $\dagger$-formally \'etale).

Another way to argue is to observe, using the pushout description of Fargues--Fontaine curves, that 
$$
\BC(\s{O})=\Q_p^\la, ~~ \BC(\s{O}(1)) = 1+\widetilde{\bb{D}}_{\bb{Q}_p}^{<1}.
$$
Moreover, a degree $1$ divisor on Fargues--Fontaine curves is the same as a degree $1$ line bundle with a fiberwise non-zero section. We claim that any such a global section is $!$-locally isomorphic to $\log(q) \colon \s{O}\to \s{O}(1)$.
Indeed, let $A$ nilperfectoid strictly totally disconnected and $\s{O}_{\FF_A}\to \mathcal{L}$ be a degree $1$ divisor on $\FF_A$. One can pick an isomorphism of the pullback $\bar{\mathcal{L}}$ of $\mathcal{L}$ to $\FF_{A^u}$  with $\s{O}(1)$,  after a Berkovich-pro-\'etale cover of $A^u$, which lifts uniquely to a Berkovich-pro-\'etale cover $A'$ of $A$ by \cite[Proposition 2.5.9]{drff_stacks}. Since $A$ is strictly totally disconnected nilperfectoid, we can by \cite[Proposition 4.8.5]{drff_stacks} assume that $A\to A'$ is descendable, up to replacing the original $A'$ by another ring with the same uniform completion. We can also assume that $A'$ is nilperfectoid, up to taking a further descendable cover.  Thus, we have proven that any degree $1$-divisor on a thickened Fargues--Fontaine curve arises locally in the $!$-topology from a global section of $\s{O}(1)$.
This proves that the natural map 
\[
\n{BC}(\s{O}(1))\backslash \{0\}\to \Div^1
\]
is surjective. Since the automorphism group of $\s{O}(1)$ is $\n{BC}(\s{O})^{\times}$, we deduce that the natural map 
\[
(\n{BC}(\s{O}(1))\backslash \{0\})/\n{BC}(\s{O})^{\times} \xrightarrow{\sim} \Div^1
\]
is an equivalence.
Finally, the desired description  of the global section then follows from the computations in \Cref{sec:perf-compl-farg-1-global-section-on-ff-curves}, \Cref{sec:perf-compl-farg-2-example-bc-o-1}, where the isomorphism $\n{BC}(\s{O}(1))\cong (1+\widetilde{\bb{D}}_{\bb{Q}_p}^{<1})$   is obtained precisely by the logarithm map, see \cite[Proposition II.2.2]{fargues2021geometrization}.
\end{proof}

\begin{remark}
  \label{sec:div1-analyt-prism-remarks-analytic-ht-stack}
  The construction of the analytic Hodge--Tate stack of $\Q_p$ as part of a general theory of analytic prismatization, which in the end takes the simple form
  \[
    \Q_p^\HT\cong \GSpec(\Q_p(\mu_{p^\infty}))/\Z_p^{\times,\la},
  \]
  has been a major desideratum in our project.
  We note that this formula should imply that\footnote{To be precise, one should take the limit in the category of Gestalten.}
  \[
    \Q_p^\HT\cong \varprojlim\limits_{\Q_p\subseteq K\subseteq \Q_p(\mu_{p^\infty}) \textrm{ Galois}} \Spf(\s{O}_K)^\HT_{\eta}/\Gal(K/\Q_p)
  \]
  with $(-)_\eta$ the generic fiber of the Hodge--Tate stack for the formal scheme $\Spf(\s{O}_K)$,  as anticipated in \cite{anschutz2022v}. Since formalizing this idea requires technical preliminaries such as  realizing $p$-adic formal stacks into Gelfand stacks via a generic fiber functor, we do not discuss it in further detail here.
\end{remark}

\begin{remark}
  \label{sec:div1-analyt-prism-remark-fiber-square-for-analytic-ht-stack}
  The description
  \[
    \Q_p^\HT\cong B\Gm\times_{B\Gm^\dR} \GSpec(\Q_p)
  \]
  from \Cref{prop:explicit-div-1} (2)
  might look surprisingly simple.
  However, the morphism $\iota^\dR\colon \GSpec(\Q_p)\to B\Gm^\dR$ is very interesting.
  Namely, for each Tate-Huber ring $B$ over $\Q_p$ the groupoid of morphisms
  \[
    \GSpec(B)\to B\Gm^\dR
  \]
  is equivalent to $\widehat{\s{O}}$-line bundles on the pro-\'etale site of $\Spa(B)$.
  Indeed, this follows by descent from the perfectoid case, where the vanishing of $\Gm^\dagger$-cohomology (\cite[Proposition 4.6.7]{drff_stacks}) implies that on perfectoids $\Gm^\dR$-torsors are equivalent to $\Gm$-torsors.
  Alternatively, one can use adjunction to see that maps $\GSpec(B)\to B\Gm^\dR$ are equivalent to maps $\Marc(B)\to B\Gm^\diamond$.
  
  Applied to $B=\Q_p$, we can see that the morphism $\iota^\dR\colon \GSpec(\Q_p)\to B\Gm^\dR$ is classified by the Tate twist $\widehat{\s{O}}\otimes_{\Q_p}\Q_p(1)\cong \mathrm{gr}^1 \mathbb{B}^+_\dR$ on $\Spa(\Q_p)_{\proet}$.\footnote{The identification with the Tate twisted $\widehat{\s{O}}$ follows from the following fact: if $(s\colon \mathcal{L}\to \s{O})$ is a Cartier divisor on some Gelfand stack $X$, classified by a morphism $X\to \A^1/\Gm$, then the morphism $Z:=X\times_{\A^1/\Gm} B\Gm \to B\Gm$ classifies the line bundle $\mathcal{L}_{|Z}$ on the Gelfand stack $Z$.}

  In more casual terms one could say that $\Q_p^\HT$ is the terminal Gelfand stack, such that the Tate twist (a pro-\'etale vector bundle on $\Spa(\Q_p)$) becomes an honest line bundle.
\end{remark}

\begin{example}\label{ExampleFrobeniusTwistCyc}
Let $\bb{Q}_p^{\cyc}$ denote the cyclotomic extension of $\bb{Q}_p$ and let $\FFdisc{\bb{Q}_p^{\cyc}}^{\times}$ denote its Fargues--Fontaine disc. Let $\epsilon=(\zeta_{p^n})_{n\geq 0}\in \bb{Q}_p^{\cyc,\flat}$ be a compatible sequence of $p$-power roots of unity, and $[\epsilon]$ its Teichm\"uller lift seen as a global section of $\FFdisc{\bb{Q}_p^{\cyc}}^{\times}$. There is a natural morphism of rational $\delta$-spaces
\begin{equation}\label{ComparisonCyclotomicPrismwithFFCurve}
\FFdisc{\bb{Q}_p^{\cyc}}^{\times}\to (1+\widetilde{\bb{D}}_{\bb{Q}_p}^{<1}) \backslash \{1\}
\end{equation}
sending $q\mapsto [\epsilon^{1/p}]$ (note the Frobenius inverse!) making the source the uniform completion of the target (locally on affinoids). This map is compatible with the natural action of $\bb{Z}_p^{\times}$ on the right hand side and the Galois action of $\Gamma=\ob{Gal}(\bb{Q}_p^{\cyc}/\bb{Q}_p)$  on the left hand side via the cyclotomic character $\chi\colon \Gamma\xrightarrow{\sim} \bb{Z}_p^{\times}$. Moreover, this map is compatible with Frobenius and sends the prism structure $([p]_q)$ to the ideal generated by $\frac{[\epsilon]-1}{[\epsilon^{1/p}]-1}$, that is, to the ideal associated to the untilt  $\bb{Q}_p^{\cyc}$ of $\bb{Q}_p^{\cyc,\flat}$.  In \Cref{subsec-Sen-theory-revisited} we will see that  \eqref{ComparisonCyclotomicPrismwithFFCurve} is the \textit{locally analytic decompletion} of $\FFdisc{\bb{Q}_p^{\cyc}}^{\times}$ for the action of $\Gamma$.  
\end{example}

\subsection{The Breuil--Kisin twist}
\label{sec:breuil-kisin-twist}
Let us define an important line bundle with Frobenius structure on $\Q_p^\prism$: the Breuil--Kisin twist, continuing \Cref{sec:breuil-kisin-twist-2-definition-bk-twist}. For that, let us first introduce a piece of notation: consider the natural $\varphi$-equivariant map $\bb{Q}_{p}^{\prism}\to \FFdisc{\bb{Q}_p^{\diamond}}^{\times,\dR}$. For $\lambda\in \bb{Q}$ we have the vector bundle $\s{O}\langle \lambda\rangle$ over $\bb{Q}_p^{\HK}$ associated to the isocrystal  of slope $-\lambda$, see \cite[Example 7.1.3]{drff_stacks}.  We will denote the pullback to $\bb{Q}_p^{\prism}$ of the vector bundle associated to  the isocrystal $(\bb{Q}_p,p^{-1}\varphi)$ by  $\s{O}_{\Q_p^\prism}\langle 1 \rangle$, seen as a $\varphi$-equivariant module. We shall use a similar notation $\langle - \rangle$ for variants of $\Q_p^\prism$ such as $\widetilde{\Div^1}$ or $\Div^1$, to indicate a twist of the Frobenius slope.

\begin{construction}\label{ConstructionBKtwist}
Consider the ideal $\n{I}_{\prism}$ of $\bb{Q}_p^{\prism}$. Recall from \Cref{sec:breuil-kisin-twist-2-definition-bk-twist} that the Breuil--Kisin ideal is given by the tensor product $\n{I}_{\BK}:=\bigotimes_{n\in \bb{N}} \varphi^{n,*} \n{I}_{\prism}$.  We have a natural equivalence $\varphi^* \n{I}_{\BK}\cong \n{I}_{\prism}^{-1}\n{I}_{\BK}$ which can be considered as a $\varphi$-equivariant map $\varphi_{\n{I}_{\BK}}\colon \n{I}_{\BK}\to \n{I}_{\prism}^{-1}\n{I}_{\BK}$.  

On the other hand, consider the natural $\varphi$-equivariant map $\bb{Q}_{p}^{\prism}\to \FFdisc{\bb{Q}_p^{\diamond}}^{\times,\dR}$. We define the \textit{Breuil--Kisin twist} $\s{O}_{\bb{Q}_p^{\prism}}\{1\}$ to be the $\varphi$-equivariant module
\[
\s{O}_{\bb{Q}_p^{\prism}}\{1\}:=\n{I}_{\BK}\otimes_{\s{O}} \s{O}_{\Q_p^\prism}\langle 1 \rangle. 
\]
\end{construction}

\begin{remark}[Comparison with the integral definition]\label{RemarkFrobeniusTwistBKT}
The reader might ask about the non-standard definition of the Breuil--Kisin twist of \Cref{ConstructionBKtwist} using the combination of the Breuil--Kisin ideal and the Frobenius twist. In the integral theory, given a prism $(A,I)$ the Breuil--Kisin twist is (heuristically) the infinite tensor product of Frobenius twists $\bigotimes_n \varphi^{n,*}I$; however, one has to give a precise meaning of this infinite tensor as it is not obvious that it converges in any reasonable sense. To make this precise, for example following \cite[Section 2.2]{bhatt2022absolute}, one forms the finite products $I_k=\bigotimes_{i=0}^{k-1} \varphi^{i,*} I_{\prism}$ and studies the partial quotients $I_k/I_k^{2}$ and the natural maps among them. We note that, in the \textit{rational theory}, the natural map $r_{k+1}\colon I_{k+1}/I_{k+1}^2\to I_k/I_k^2$ already induces a natural equivalence $I_{k+1}/I_k I_{k+1}\xrightarrow{\sim} I_k/I_k^2$. However, in the \textit{integral theory} the natural map   $r_{k+1}$ is divisible by $p$ and $p^{-1}r_{k+1}$ induces an isomorphism $p^{-1}r_{k+1}\colon I_{k+1}/I_k I_{k+1}\xrightarrow{\sim}  I_k/I_k^2$, see \cite[Corollary 	2.2.10]{bhatt2022absolute}. One then defines the Breuil--Kisin twist to be the limit 
\[
A\{1\}= \varprojlim_{k}  I_k/I_k^2
\]
 along the maps $p^{-1}r_{k}$, cf.\ \cite[Construction 2.2.11]{bhatt2022absolute}. 
 
 This modification in the transition maps  in the definition of $A\{1\}$ seems inoffensive at first sight, however it does an important modification to the natural Frobenius structure of  $A\{1\}$, analogous to the twist by $\s{O}_{\Q_p^\prism}\langle 1\rangle$.
 Indeed, for each $k$ we have an isomorphism  $\varphi^{*}  I_{k}\cong  I^{-1} I_{k+1}$.
 This produces a natural isomorphism $\varphi^{*}( I_{k}/ I_{k}^2)\cong  I^{-1} I_{k+1}/ ( I^{-1} I_{k+1})^2$, composing with the map $ p^{-1}r_{k+1}$, we see that the Frobenius  $\varphi_{A\{1\}}$ is the limit of the maps 
\[
I_{k}/I_{k}^2 \xrightarrow{\varphi} ( I^{-1}I_{k+1})/( I^{-1} I_{k+1})^{2} \xrightarrow{p^{-1}r_{k+1}} ( I^{-1} I_{k})/( I^{-1}I_{k})^2 .
\]
In terms of the rational theory, the first map is the obvious Frobenius of $\n{I}_{\BK}$, while the second map corresponds to the twist by $\s{O}_{\Q_p^\prism}\langle 1\rangle $, that is, multiplication by $p^{-1}$ in the Frobenius.

The most important justification for us to add the Frobenius twist in \Cref{ConstructionBKtwist} is the construction of the prismatic logarithm map in \Cref{sec:analyt-prism-ga} below. This is analogous to  \cite[Proposition 2.5.18]{bhatt2022absolute}, where the prismatic logarithm $\log_{\prism}\colon (1+I)_{\ob{rk}=1}\to A\{1\}$ satisfies $\varphi_{A\{1\}}(\log_{\prism}(u))=\log_{\prism}(u)$. Rationally, using cyclotomic coordinates,  the prismatic logarithm is equivalent to the usual logarithm up to multiplying by an element of  $\bb{Q}_p[[q-1]]$, see \cite[Proposition 2.6.10]{bhatt2022absolute} and \cite[Lemma 4.6]{anschutz2020p}, and we have that $p^{-1}\varphi(\log(q)) =\log(q)$.

In more explicit terms, let $(A,I):=(A_\inf(\Q_p^\cycl),(\widetilde{\xi}))$ with $\widetilde{\xi}=[p]_q$ be the completed  cyclotomic prism.
We have $A\{1\}=\frac{1}{\mu}A_\inf\otimes_{\Z_p}\Z_p(1)$, cf.\ \cite[Example 4.24]{Bhatt2018}, where $\mu=q-1$ and $q=[\varepsilon]$ for a generator $\varepsilon\in \Z_p(1)(\Q_p^\cycl)$.
After multiplying with $\log(q)$, we get a map
\[
  \log(q)A\{1\}\otimes_{A}\s{O}_{\FFdisc{\Q_p^{\cycl,\flat}}^\times}\to \s{O}_{\FFdisc{\Q_p^{\cycl,\flat}}^\times},\ a\otimes \varepsilon\mapsto a 
\]
with image the Breuil-Kisin ideal $I_{\BK}=(\frac{\log(q)}{\mu})$.
However, the multiplication by $\log(q)$ changes the Frobenius structure exactly by multiplication by $p$.
\end{remark}

\begin{remark}\label{rk:line-bundles-on-qp-prism-adn-div1}
To avoid confusion, let us recall several line bundles on $\Q_p^\prism, \Div^{1}$ or $\Div^{1,\dR}$ that appeared before and their relations.
\begin{enumerate}
    \item The isocrystal $(\Q_p,p^{-1}\varphi)$ pulls back to the line bundle $\s{O}_{\HK}\langle 1 \rangle$ on $\Div^{1,\dR}$. Note that (by definition) $\s{O}_{\Div^1}\langle 1 \rangle= f^\ast \s{O}_{\HK} \langle 1 \rangle$, if $f: \Div^1\to \Div^{1,\dR}$.

\item The Breuil-Kisin twist $\s{O}_{\Q_p^\prism}\{1\}$ on $\Q_p^\prism$ just defined.
\item On $\widetilde{\Div^1}$ one can consider a perfected version of the Breuil-Kisin twist defined as
  \[
    \s{O}_{\widetilde{\Div^1}}\{1\}=\prod\limits_{n\in \Z} \varphi^n(\mathcal{I}_\prism),
  \]
 cf.\ \Cref{quotient-perfect}. It is equipped with a Frobenius structure $\varphi_{\s{O}_{\widetilde{\Div^1}}\{1\}}: \varphi^\ast \s{O}_{\widetilde{\Div^1}}\{1\} \cong  \s{O}_{\widetilde{\Div^1}}\{1\}$, and thus descends to the stack $\Div^1$. 
  \item The line bundle $\s{O}_{\widetilde{\Div^1}}\{1\}$ descends to a line bundle $\s{O}_{\Div^1}\{1\}$ on $\Div^1$, which is $\mathcal{I}_{\Div^1}\langle 1\rangle$, the tensor product of $\s{O}_{\Div^1}\langle 1 \rangle$ by  the ideal sheaf $\mathcal{I}_{\Div^1}$ of the Hodge-Tate point on $\Div^1$. 
  \item By the fundamental exact sequence, the pullback of $\s{O}_{\Div^1}\{1\}$ along $\FF_{\Q_p^\diamond} \to \Div^1$ is the first Tate twist, seen as a pro-\'etale $\Q_p$-local system on $\Q_p$. The pullback of $\s{O}_{\Div^1}\langle 1 \rangle$ along $\FF_{\Q_p^\diamond} \to \Div^1$ is instead the line bundle with Frobenius slope $1$ (which would usually be denoted $\s{O}_{\FF_{\Q_p^\diamond}}(1)$, but which in this context should perhaps rather be denoted $\s{O}_{\FF_{\Q_p^\diamond}}\langle 1 \rangle$, for consistency).  

     \end{enumerate}
\end{remark}

\subsection{The  prismatic logarithm and the analytic prismatizaton of $\bb{G}_m$}
\label{sec:analyt-prism-ga}
The goal of this section is to give an  explicit description of the prismatization of $\bb{G}_m$ as a commutative group (along the way, we will also get a description of the prismatization of the affine line). For that, we follow the strategy of \cite{drinfeld20211} recalled below and cover $\bb{G}_m^{\prism}$ by a suitable perfected variant of $\bb{G}_m$, in this case the base change to $\bb{Q}_p^{\prism}$ of the  \textit{perfection} (in Gelfand stacks) $\widetilde{\bb{G}}_m:=\varprojlim_{x\mapsto x^p} \bb{G}_m$.   The map  $\widetilde{\bb{G}}_{m,\bb{Q}_p^{\prism}}\to \bb{G}_m^{\prism}$  will be constructed using a suitable Teichm\"uller lift for thickened Fargues--Fontaine discs. We want to give a precise description of the group $\n{G}:=\ob{ker}(\widetilde{\bb{G}}_{m,\bb{Q}_p^{\prism}}\to \bb{G}_m^{\prism})$, analogue to the Cartier dual of the $1$-dimensional formal group of \cite{drinfeld20211}. With this in mind, we construct a \textit{prismatic logarithm}  $\log_{\prism}\colon \n{G}\to \s{O}\{1\}$ with values in the  the Breuil--Kisin twist, and show that the map $\log_{\prism}$ is $\dagger$-formally \'etale. In order to have similar descriptions of $\bb{G}_m^{\widetilde{\Div^1}}$ and $\bb{G}_m^{\Div^1}$, we need to pay special attention to the Frobenius structures of all the objects discussed above.

Let us first recall Drinfeld's description of the prismatization of the $p$-adic formal multiplicative group, \cite{drinfeld20211}.

\begin{remark}
  \label{sec:analyt-prism-ga-1-canonical-morphism-for-delta-ring-in-classical-prismatization}
  If $A$ is a $p$-complete $\delta$-ring, then there exists a canonical morphism
  \[
    \Spf(A)\times \Z_p^\prism\to A^\prism,
  \]
  where $(-)^\prism$ denotes in this remark the (classical) prismatization of $p$-adic formal schemes.
  Indeed, given a $p$-nilpotent ring $R$ with a map $A\to R$ and a Cartier-Witt divisor $(s\colon I\to W(R))$, one gets a canonical morphism $A\to W(R)$ of $\delta$-rings, and thus a morphism $A\to W(R)\to \cone(s)$, i.e., a point in $A^\prism(R)$.
  The same construction works for $\delta$-stacks like $\A^1/\Gm$ (with Frobenius induced by the $p$-power map).

  We view now the $p$-adic formal scheme $\Gm=\Spf(\Z_p\langle T^{\pm 1}\rangle)$ as a $\delta$-stack by setting $\delta(T)=0$, i.e., $\varphi(T)=T^p$.
  The resulting morphism
  \[
    \rho\colon \Gm\times \Z_p^\prism\to \Gm^\prism
  \]
  sends a pair $(x\in R^\times, (\alpha\colon I\to W(R)))$ to the image of the Teichm\"uller lift $[x]\in W(R)^\times$ under the canonical morphism $W(R)\to \cone(I\to W(R))$.
  The morphism $\rho$ is a morphism of group stacks over $\Z_p^\prism$ and a faithfully flat affine surjection.
  Set
  \[
    \mathcal{G}_{\mathrm{Dr}}:=\mathrm{fib}(\Gm\times \Z_p^\prism\xto{\rho} \Gm^\prism),
  \]
  so that $\Gm^\prism\cong \cone(\mathcal{G}_{\mathrm{Dr}}\to \Gm\times \Z_p^\prism)$.
  Drinfeld describes the structure of $\mathcal{G}_{\mathrm{Dr}}$ explicitly: it is the Cartier dual of a $1$-dimensional commutative formal group $H_{\mathrm{Dr}}$ over $\Z_p^\prism$.
 After pullback along the surjection $\Spf(\Z_p[[q-1]])\to \Z_p^\prism$ the formal group $H_{\mathrm{Dr}}$ is isomorphic to the formal group with formal group law $x+y+(q-1)xy$.

  Moreover, Drinfeld interprets the prismatic logarithm as a morphism of group schemes $\log_\prism\colon \mathcal{G}_{\mathrm{Dr}}\to \s{O}_{\BK}\{1\}$, where $\s{O}_{\BK}\{1\}$ is the integral Breuil-Kisin twist on $\Z_p^\prism$ (seen here as a geometric line bundle).
\end{remark}

We switch back to the analytic setup and use the following notations for the rest of this subsection:
\begin{itemize}
\item $\widetilde{\mathbb{G}}_m=\varprojlim_{x\mapsto x^p} \Gm$, the (decompleted) perfected multiplicative group over $\Q_p$, we endow $\widetilde{\mathbb{G}}_m$ with the automorphism $\varphi(T)=T^p$, 
\item $\widetilde{\mathbb{G}}_m^u$ the uniformly completed version of $\widetilde{\mathbb{G}}_m$, or more precisely, the preperfectoid multiplicative group obtained as the uniform completion of the limit of $\bb{G}_m$ along raising to the $p$-th power maps,
\item $\Z_p(1)^\sm:=\mathrm{fib}(\widetilde{\mathbb{G}}_m\to \Gm)\cong \GSpec(\Q_p[T^{1/p^\infty}]/(T-1))$,
\item for a Gelfand stack $Y$ we denote by $\mathbb{G}_{m,Y},\widetilde{\mathbb{G}}_{m,Y},\ldots$ the base change of $\Gm,\widetilde{\mathbb{G}}_m,\ldots$ to $Y$.
\end{itemize}

We want to describe the analytic prismatization $X^\prism=\Gm^\prism$ in a similar way as in \Cref{sec:analyt-prism-ga-1-canonical-morphism-for-delta-ring-in-classical-prismatization}.
We start by constructing an analogue of the morphism $\rho$, which amounts to the construction of a Teichm\"uller lift.\footnote{We content ourselves with the explicit construction for the case of $\widetilde{\mathbb{G}}_m$, but we note the following general construction: if $\mathcal{X}$ is an analytic  $\delta$-stack, then one gets a natural morphism $\mathcal{X}\times \Q_{p,\prism}\to (\mathcal{X}_{\real})_\prism$ by sending a  $\delta$-space $(S,Y)$ with a degree $1$ Cartier divisor $D\subseteq Y$ and a morphism $(S,Y)\to \mathcal{X}$ to the composition $D\to Y=(S,Y)_{\real}\to \mathcal{X}_{\real}$.}

\begin{lemma}
  \label{sec:analyt-prism-gm-description-of-gm-perf}
  The natural map
  \[
    \widetilde{\mathbb{G}}_m \to \widetilde{\mathbb{G}}_m^{\dR}\times_{\Gm^\dR}\Gm
  \]
  is an isomorphism.
  In particular, for any nilperfectoid $\Q_p$-algebra $A$ the natural map
  \[
    \varprojlim\limits_{x\mapsto x^p}A^\times\to   \varprojlim\limits_{x\mapsto x^p} A^{u,\times}   \times_{A^{u,\times}} A^\times
  \]
  is an isomorphism.
\end{lemma}
\begin{proof}
  This follows from the fact that (over $\Q_p$) the morphism $\Gm\to \Gm,\ x\mapsto x^p$ is \'etale, which implies that $\widetilde{\mathbb{G}}_m\to \Gm$ is $\dagger$-formally \'etale.
\end{proof}

\begin{construction}[Teichm\"uller morphism]\label{ConstructionPrismaticLog}
 We construct now a $\varphi$-equivariant  morphism of commutative group stacks $\rho\colon \widetilde{\mathbb{G}}_{m,\Q_p^\prism} \to \Gm^{\prism}$ over $\Q_p^\prism$ as follows: 

  Let $A$ be a nilperfectoid ring over $\Q_p$, with $T:=\GSpec(A)$ and associated  $\delta$-space $\n{X}_T$.
  Let $D\to \n{X}_T$ be a degree $1$ divisor and let $\overline{D}\to \bb{D}^{\times}_{\FF,A^u}$ be the pullback of $D$.
  We have to produce a natural map  of commutative groups
\begin{equation}\label{eq0j0o3mllafdaf}
\rho_A\colon\varprojlim_{x\mapsto x^p} A^{\times} \to \s{O}(D)^{\times}.
\end{equation}
For that, we will construct a Teichm\"uller morphism of commutative monoids
\begin{equation}\label{eqTeichmullerMapPrisms}
[-]\colon \varprojlim_{x\mapsto x^p} A \to \s{O}(\n{X}_T)^{\varphi(-)=(-)^p}
\end{equation}
and define the map $\rho_A$ \eqref{eq0j0o3mllafdaf} as the composite of $[-]$ with the natural map $ \s{O}(\n{X}_T)\to  \s{O}(D)$.
For constructing the Teichm\"uller map, we observe the $\varphi$-equivariant  pullback isomorphism of solid $\Q_p$-algebras 
\[
\s{O}(\n{X}_T)=\s{O}(\bb{D}_{\FF,A^u}^{\times})\times_{\prod_{n\in \N} A^{u}} \prod_{n\in \N} A 
\]
induced by the description in \Cref{Prop:Universal-FFNilPerfd2} (here the Frobenius on $\prod_{n\in \N} A$ resp.\ $\prod_{n\in \N}A^u$ is the map $(a_n)_{n\in \N}\mapsto (a_{n+1})_{n\in \N}$).
Hence, it suffices to produce a map from $\varprojlim_{x\mapsto x^p} A$ to $\s{O}(\bb{D}_{\FF,A^u}^{\times})$ and $\prod_{n\in \N} A$ compatible with the natural map to $\prod_{n\in \N} A^{u}$ and Frobenius.
We do this as follows: 
\begin{itemize}
\item[(a)] The map to  $\s{O}(\bb{D}_{\FF,A^u}^{\times})$ factors through the classical Teichm\"uller map 
\[
\varprojlim_{x\mapsto x^p} A\to \varprojlim_{x\mapsto x^p} A^{u} = A^{u,\flat}\xrightarrow{[-]}  \s{O}(\bb{D}_{\FF,A^{u}}^{\times})^{\varphi(-)=(-)^p}.
\]

\item[(b)] Let $\theta\colon \s{O}(\bb{D}^{\times}_{\FF,A^u})\to A^u$ be Fontaine's map, then the map to $\prod_{n\in \N} A^u$  is the composition of the Teichm\"uller map of (a) with the map 
\[
 \s{O}(\bb{D}_{\FF,A^{u}}^{\times})\xrightarrow{x\mapsto \prod_{n\in \N} \theta(\varphi^{n}(x))} \prod_{n\in \N} A^u. 
\]

\item[(c)] Finally, the map $\varprojlim_{x\mapsto x^p} A \to \prod_{n\in \N} A$ is the composition of the first projection $\varprojlim_{x\mapsto x^p}A\to A$ with the map
  \[
    A\to \prod_{n\in \N} A,\ a\mapsto (a,a^p,a^{p^2},\ldots).
  \] 
\end{itemize}

We note that by construction $[-]$ takes its image naturally in the multiplicative monoid $\s{O}(\n{X}_T)^{\varphi(-)=(-)^p}$ of elements $x$ such that $\varphi(x)=x^p$.

Furthermore, we can factor $[-]$ over the natural map $\s{O}(\FFdisc{A}^\times)\to \s{O}(\n{X}_T)$ by similarly constructing the map
\[
 [-]^\perf\colon \varprojlim_{x\mapsto x^p}A\to \s{O}(\FFdisc{A}^\times)^{\varphi(-)=(-)^p}, 
\]
by using in (c) instead the map
\[
  \varprojlim\limits_{x\mapsto x^p}A\to \prod\limits_{n\in \Z} A,\ (a_0,a_1,\ldots)\mapsto (\ldots, a_1,a_0,a_0^p,\ldots)
\]
(we note that $\s{O}(\FFdisc{A}^\times)=\s{O}(\bb{D}_{\FF,A^u}^{\times})\times_{\prod_{n\in \Z} A^{u}} \prod_{n\in \Z} A $).
\end{construction}

The following lemma gives a first description of the analytic prismatization of $\Gm$ (and of $1+\mathbb{D}^{<1}_{\bb{Q}_p}$).

\begin{lemma}
  \label{sec:analyt-prism-gm-analytic-prismatization-of-gm}
  \begin{enumerate}
  \item The morphism $\rho\colon \widetilde{\mathbb{G}}_{m,\Q_p^\prism}\to \Gm^{\prism}$ is a $\varphi$-equivariant  epimorphism of commutative group Gelfand stacks over $\Q_p^\prism$.\footnote{With commutative group Gelfand stacks meaning functors on Gelfand rings factoring over the category of animated abelian groups whose underlying anima-valued functor is a Gelfand stack.}
    In particular, $\Gm^{\prism}\cong \widetilde{\mathbb{G}}_{m,\Q_p^\prism}/\mathcal{G}_{\Q_p^\prism}$, where $\mathcal{G}_{\Q_p^\prism}:=\mathrm{fib}(\rho)$ is the fiber at $1$.
  \item We have a $\varphi$-equivariant commutative diagram with cartesian squares
\[\begin{tikzcd}
	{\mathcal{G}_{\Q_p^\prism}} & {(1+\widetilde{\mathbb{D}}^{<1})\times_{\GSpec(\Q_p)} \Q_p^\prism} & {(1+\mathbb{D}^{<1})^\prism} \\
	{\mathcal{G}_{\Q_p^\prism}} & {\widetilde{\mathbb{G}}_{m,\Q_p^\prism}} & {\Gm^\prism}
	\arrow[from=1-1, to=1-2]
	\arrow["=", from=1-1, to=2-1]
	\arrow["\rho_{\mathbb{D}}" , from=1-2, to=1-3]
	\arrow[from=1-2, to=2-2]
	\arrow[from=1-3, to=2-3]
	\arrow[from=2-1, to=2-2]
	\arrow["\rho", from=2-2, to=2-3]
\end{tikzcd}\]
      \end{enumerate}
\end{lemma}

\begin{proof}
  Let $A$ be a nilperfectoid $\Q_p$-algebra, $T:=\GSpec(A)$ and $D\to \n{X}_T$ a degree $1$ divisor on the associated  $\delta$-space.
  For (1), we want to see that, locally in the $!$-topology on $T$, the map 
\begin{equation}\label{eqoj93oamasbduq3}
\varprojlim_{x\mapsto x^p} A^{\times}\to \s{O}(D)^{\times} 
\end{equation}
is an epimorphism.
Let $A^{\flat,\sharp}$ be the untilt of $A^{\flat}$ attached to the  classical degree $1$ divisor $\overline{D}\to \bb{D}^{\times}_{\FF,A^u}$ associated with $D$.
The map \eqref{eqoj93oamasbduq3} produces a natural map
\begin{equation}\label{eq0jo3awlmeawsd}
\varprojlim_{x\mapsto x^p} A^{u,\times} \to  \s{O}(\overline{D})^{\times}= A^{\flat, \sharp,\times}. 
\end{equation}
This map is surjective locally in the arc-topology of $A$ (even in the pro-finite \'etale topology).
Hence, we can assume without loss of generality that \eqref{eq0jo3awlmeawsd} is an epimorphism, e.g., by assuming that $A^u$ is strictly totally disconnected.
Now, by taking a further analytic cover of $A$ and a Frobenius twist, we can assume without loss of generality that the intersection of $\Marc(A^{\flat,\sharp})$ and $\bigsqcup_{n\in \N} \varphi^{-n}(\Marc(A^u))$ in $\bb{D}^{\times,\diamond}_{\FF,A^u}$ is contained in $\Marc(A)$.
We can also assume that the ideal of $D$ is principal generated by some element $\xi$.
We have a fiber sequence 
\[
\ob{Nil}^{\dagger}(A)/\xi\to \s{O}(D)\to A^{\flat,\sharp}
\]
so that we have a fiber sequence of commutative groups 
\[
1+(\ob{Nil}^{\dagger}(A)/\xi)\to \s{O}(D)^\times \to A^{\flat, \sharp,\times}
\]
(here we use that $\ob{Nil}^\dagger(A)/\xi$ is contained in the $\dagger$-nilradical of $\s{O}(D)$, so that $1+\ob{Nil}^\dagger(A)/\xi$ lies in $\s{O}(D)^\times$).
Thus, the map \eqref{eqoj93oamasbduq3} sends $1+\Nil^{\dagger}(A)=\varprojlim_{x\mapsto x^p} (1+\Nil^{\dagger}(A))\subset \varprojlim_{x\mapsto x^p} A$ to $(1+\Nil^{\dagger}(A)/\xi)\subset \s{O}(D)$, given by the obvious quotient map
\[
1+\Nil^{\dagger}(A)\to 1+(\Nil^{\dagger}(A)/\xi)
\]
thanks to \Cref{ConstructionPrismaticLog}.  Summarizing, we have a morphism of fiber sequences of commutative groups 
\[
\begin{tikzcd}
(1+\Nil^{\dagger}(A)) \ar[d] \ar[r] & \varprojlim_{x\mapsto x^p} A^{\times} \ar[d] \ar[r] &   \varprojlim_{x\mapsto x^p} A^{u,\times} \ar[d] \\
 (1+\Nil^{\dagger}(A)/\xi) \ar[r] & \s{O}(D)^{\times} \ar[r] & A^{\flat,\sharp,\times}
\end{tikzcd}
\]
where the right horizontal arrows and the left and right vertical arrows are epimorphisms. This implies that the middle vertical  arrow  is an epimorphism, proving what we wanted.

Assertion (2) can be deduced from \Cref{sec:defin-first-prop-5-prismatization-for-dagger-formally-etale-maps} and the same assertion on de Rham stacks.
\end{proof}

Our next task will be to describe the group stack $\mathcal{G}_{\Q_p^\prism}$ explicitly.
We do this first using a version of the prismatic logarithm (which in our case is easy to construct via the usual logarithm).

In the following we will denote by $\bb{G}_{a,\Q_p^\prism}\{1\}$ the geometric line bundle associated with the Breuil-Kisin twist, i.e.\ the Gelfand stack over $\Q_p^\prism$ whose $A$-valued sections (for $A$ a Gelfand ring over $\Q_p^\prism$) are given by $\Gamma(\GSpec(A),\s{O}_{\Q_p^\prism}\{1\}_A)$, where $\s{O}_{\Q_p^\prism}\{1\}_A$ is the pullback of $\s{O}_{\Q_p^\prism}\{1\}$ to $\GSpec(A)$.

Let
\[
  \mu_{\BK}\colon \Q_p^\prism\to \A^1/\Gm
\]
be the morphism classifying the Breuil--Kisin divisor $\n{I}_{\BK}\to \s{O}$ on $\Q_p^\prism$.
Twisting by the Frobenius twist $\s{O}_{\Q_p^\prism} \langle 1 \rangle $ we get a $\varphi$-equivariant section $\s{O}\{1\}\to \s{O}_{\Q_p^\prism} \langle 1 \rangle $ which geometrically we write as $\bb{G}_{a,\bb{Q}_p^{\prism}}\{1\}\to \bb{G}_{a,\bb{Q}_p^{\prism}}\langle 1 \rangle$. We note that we still have a natural identification 
\begin{equation}\label{eqConeBKidealHKTwisted}
\bb{G}_a^{\cone}\times_{\bb{A}^1/\bb{G}_m, \mu_{\BK}} \bb{Q}_p^{\prism} \cong  \ob{cofib}(\bb{G}_{a,\bb{Q}_p^{\prism}}\{1\} \to \bb{G}_{a,\bb{Q}_p^{\prism}}\langle 1 \rangle)
\end{equation}
as the Frobenius twist can be naturally trivialized after pullback to degree $1$ divisors. 
\begin{construction}[Prismatic logarithm]\label{ConstructionTwists}
  The \textit{prismatic logarithm} is the unique $\varphi$-equivariant morphism of commutative group Gelfand stacks 
\begin{equation}\label{eqomapoewnpoawfraew}
\log_{\prism}\colon \mathcal{G}_{\Q_p^\prism}\to \bb{G}_{a,\bb{Q}_p^{\prism}}\{1\}
\end{equation}
fitting into the following  commutative diagram of fiber sequences of commutative group Gelfand stacks over $\Q_p^\prism$:

\[\begin{tikzcd}
	{\mathcal{G}_{\Q_p^\prism}} & {(1+\widetilde{\mathbb{D}}^{<1})_{\Q_p^\prism}} & {(1+\mathbb{D}^{<1})^\prism} \\
	\bb{G}_{a,\bb{Q}_p^{\prism}}\{1\} & \bb{G}_{a,\bb{Q}_p^{\prism}}\langle 1 \rangle & {\Ga^\cone\times_{\A^1/\Gm,\mu_{\BK}} \Q_p^\prism}
	\arrow[from=1-1, to=1-2]
	\arrow["{\log_\prism}"', from=1-1, to=2-1]
	\arrow["{\rho_{\mathbb{D}}}", from=1-2, to=1-3]
	\arrow["\beta", from=1-2, to=2-2]
	\arrow["\gamma", from=1-3, to=2-3]
	\arrow[from=2-1, to=2-2]
	\arrow["\alpha", from=2-2, to=2-3]
      \end{tikzcd}\]
    where: 
\begin{itemize} 

\item The map $\alpha$ comes from  \eqref{eqConeBKidealHKTwisted}.

\item The map $\beta$ is the base change to $\Q_p^\prism$ of the composition
  \[
    1+\widetilde{\mathbb{D}}^{<1}\to 1+\mathbb{D}^{<1}\xto{\log} \Ga,
  \]
  with the first morphism being the projection to the first component; we note that to make this map $\varphi$-equivariant one must add the Hyodo--Kato twist.\footnote{Indeed, let $D_{1}$ be the isocrystal $(\Q_p,p\varphi)$ with dual $D_{-1}$. We have \[\underline{\GSpec}_{\GSpec(\Q_p)}(\mathrm{Sym}^\bullet_{\Q_p} D_1)\times_{\GSpec(\Q_p)} \Q_p^\prism\cong \mathbb{G}_{a,\Q_p^\prism}\langle 1\rangle\] with diagonal Frobenius action, according to the construction of the geometric vector bundle for $\s{O}_{\Q_p^\prism}\langle 1 \rangle$ (i.e., the geometric vector bundle is the relative spectrum of the symmetric algebra on the dual, and we note that $D_{-1}$ pulls back to $\s{O}_{\Q_p^\prism}\langle 1\rangle$ by definition). Now, the map $\Q_p[T]\cong \mathrm{Sym}^\bullet_{\Q_p} D_1\to \s{O}(1+\widetilde{\mathbb{D}}),\ T\mapsto \mathrm{log}([-])$ is $\varphi$-equivariant because $\varphi(T)=pT$ and $\varphi(\log([-]))=\log([-]^p)=p\log([-])$.} 
\item The map $\gamma$ is the composition
  \[
    (1+\mathbb{D}^{<1})^\prism\to (1+\mathbb{D}^{<1})^\cone\times_{\A^1/\Gm,\mu_{\BK}}\Q_p^\prism\xto{(\log)^\cone} \Ga^\cone\times_{\A^1/\Gm,\mu_{\BK}} \Q_p^\prism,
  \]
  where the base changes to $\Q_p^\prism$ are along the morphism $\mu_{\BK}\colon \Q_p^\prism\to \A^1/\Gm$.
\end{itemize}

The (natural) commutativity of the right square is easy to check. More precisely, let $A$ be a nilperfectoid ring with thickened Fargues--Fontaine disc $\n{X}_A$, and let $D\to \n{X}_A$ be a degree $1$-divisor. Let $\n{I}_A$ denote the pullback of $\n{I}_{\prism}$ to $A$.   The commutativity of the right square of the diagram then boils down to the commutativity of 
\[
\begin{tikzcd}
\varprojlim_{x\mapsto x^p} (1+ A^{<1} ) \ar[d,"{\log([-])}"] \ar[r]  & 1+\s{O}(D)^{<1}  \ar[d,"\log(-)"] \\ 
 \Gamma(\n{X}_A,\s{O}(1)) \ar[r] &  \prod_{n\in \N}  \varphi^{n,*}(A/\n{I}_A) = A/\n{I}_{\BK,A}
\end{tikzcd}
\]
which follows from the naturality of the logarithm map along morphisms of Gelfand rings.  This shows  that $\log_\prism$ is well-defined.

\end{construction}

We can now describe $\mathcal{G}_{\Q_p^\prism}$ explicitly as a twisted form of $\Z_p^\la$.\footnote{Strictly speaking, it is a twisted form only on the region $(0,p)$ corresponding to the perfect prismatization, see \Cref{RemarkTwistedFormZpla}.}

\begin{theorem}
  \label{sec:analyt-prism-gm-description-of-g-qp-prism}
  The morphism $\log_\prism\colon \mathcal{G}_{\Q_p^\prism}\to \bb{G}_{a,\Q_p^\prism}\{1\}$ is $\dagger$-formally \'etale.
  In particular, there is a fiber square
\[\begin{tikzcd}
	{\mathcal{G}_{\Q_p^\prism}} & { \bb{G}_{a,\bb{Q}_p^{\prism}}\{1\} } \\
	{ \FFdisc{\Marc(\bb{Q}_p)}^{\times,\dR} \times_{\GSpec(\Q_p)} \Z_p(1)^{\sm} } &{ (\bb{G}_{a,\bb{Q}_p^{\prism}}\{1\})^{\dR} }
	\arrow["{\log_\prism}", from=1-1, to=1-2]
	\arrow[from=1-1, to=2-1]
	\arrow[from=1-2, to=2-2]
	\arrow["{\log_\prism^\dR}", from=2-1, to=2-2]
\end{tikzcd}\]
  and an exact sequence
  \[
    0\to \Ga^\dagger\otimes_{\Ga} \bb{G}_{a,\bb{Q}_p^{\prism}}\{1\}  \to \mathcal{G}_{\Q_p^\prism}\to \Z_p(1)^\sm \to 0.
  \]
\end{theorem}
\begin{proof}
  To show that $\log_\prism$ is $\dagger$-formally \'etale, it suffices to show that the morphisms $\gamma$ and $\beta$ in \Cref{ConstructionTwists} are $\dagger$-formally \'etale.
  For $\beta$ this is clear as $1+\widetilde{\mathbb{D}}^{<1}\to 1+\mathbb{D}^{<1}$ and $\log\colon 1+\mathbb{D}^{<1}\to \Ga$ are $\dagger$-formally \'etale.
  For $\gamma$, we can use \Cref{prop:pullback-description-twistor-stack1} to see that $(1+\mathbb{D}^{<1})^\prism\to (1+\mathbb{D}^{<1})^\cone\times_{\A^1/\Gm,\mu_{\BK}} \Q_p^\prism$ is $\dagger$-formally \'etale, and then again that $\log\colon 1+\mathbb{D}^{<1}\to \Ga$ is \'etale.
  
  Finally, we note that $\mathcal{G}_{\Q_p^\prism}^\dR\cong \Z_p(1)^\sm\times \FFdisc{\Marc(\bb{Q}_p)}^{\times}  $, namely,  it follows by passing to de Rham stacks in \Cref{sec:analyt-prism-gm-analytic-prismatization-of-gm} and using that $(1+\mathbb{D}^{<1})^{\prism,\dR}\cong (1+\widetilde{\mathbb{D}}^{<1})^{\prism,\dR}/\Z_p(1)^\sm$.
\end{proof}

\begin{remark}\label{quotient-perfect}
The results above admit an analogue for the perfect prismatization: $\Gm^{\widetilde{\Div^1}}$ is the quotient of $\widetilde{\mathbb{G}}_{m,\widetilde{\Div^1}}$ by a commutative group Gelfand stack\footnote{Note that $\mathcal{G}_{\widetilde{\Div^1}}$ is \textit{not} the base change of $\mathcal{G}_{\widetilde{\Div^1}}$ along $\widetilde{\Div^1}\to \Div^1$. We hope that our choice of notation does not create too much confusion.} $\mathcal{G}_{\widetilde{\Div^1}}$, which admits a $\dagger$-formally \'etale map
$$
\log_{\prism,\mathrm{perf}}: \mathcal{G}_{\widetilde{\Div^1}} \to \bb{G}_{a,\widetilde{\Div^1}}\{1\},
$$
where $\bb{G}_{a,\widetilde{\Div^1}}\{1\} $ denotes the (underlying line bundle of the) perfect version of the first Breuil--Kisin twist  $\s{O}_{\widetilde{\Div^1}}\{1\}$ (which after forgetting Frobenius is the ideal sheaf of $\sqcup_{n \in \Z} \varphi^n(\Q_p^\HT)$), and similarly for $1+\mathbb{D}^{<1}$.

Moreover, the constructions are Frobenius-equivariant, so one also deduces analogous results for the $p$-adic twistor stack: $\Gm^{\Div^1}$ is the quotient of $\widetilde{\mathbb{G}}_{m,\Div^1}/\varphi^\Z$ by a commutative group Gelfand stack $\mathcal{G}_{\Div^1}$, which admits a $\dagger$-formally \'etale map
$$
\log_{\Div^1}: \mathcal{G}_{\Div^1} \to \bb{G}_{a,\Div^1}\{1\},
$$
where $\s{O}_{\Div^1}\{1\}$ denotes the descent by Frobenius of $\s{O}_{\widetilde{\Div^1}}\{1\}$, which is given by the tensor product $\n{I}_{\Div^1}\otimes \s{O}_{\Div^1}\langle 1 \rangle$ of  the ideal sheaf and  the Frobenius twist, so a line bundle of slope zero. Later in \Cref{sec:arithmetic-case:-p-adic-fields} we will see that $\s{O}_{\Div^1}\{1\}$ is  the \textit{decompletion} of the Tate twist from $\FF_{\Marc(\bb{Q}_p)}$ to $\Div^1$. 
\end{remark}

\begin{remark}\label{RemarkTwistedFormZpla}
We comment on the structure of  $\n{G}_{\bb{Q}_p^{\prism}}$ as a twisted form of $\bb{Z}_p^{\la}$.  At the level of arc-stacks,  we see that $\n{G}_{\bb{Q}_p^{\prism}}^{\diamond}$ is the kernel of the map 
\[
\widetilde{\mathbb{G}}_m^\diamond\times_{\bb{F}_p^{\diamond}} \bb{Q}_p^{\diamond}\to \bb{Q}_p^{\diamond} \times_{\bb{F}_p^{\diamond}} \Gm^\diamond
\]
that sends a triple $(R,R^{\sharp}, a^{\flat})$   consisting of a perfectoid  ring $R/\bb{Q}_p$, an untilt $R^{\sharp}$ of $R^\flat$, and a sequence $a^{\flat}\in \varprojlim_{x\mapsto x^p} R^{\sharp}$, to the triple $(R,R^{\sharp}, \theta([a^{\flat}]))$ where $(R,R^{\sharp})$ are as above, and $\theta\colon \s{O}(\FFdisc{R^{\flat}}^{\times})\to R$ is Fontaine's $\theta$-map with values in $R$.\footnote{In the previous, the ring $R/\bb{Q}_p$ corresponds to the untilt given by the second projection to $\bb{Q}_p^{\diamond}$ while $R^{\sharp}/\bb{Q}_p$ corresponds to the first projection which factors through $\widetilde{\bb{G}}_m^{\diamond}\to \bb{Q}_p^{\diamond}$.} 
 The tilting equivalence produces a natural equivalence of arc-stacks
\[
\widetilde{\mathbb{G}}_m^\diamond\times_{\bb{F}_p^{\diamond}} \bb{Q}_p^{\diamond}\cong  \bb{Q}_p^{\diamond} \times_{\bb{F}_p^{\diamond}} \widetilde{\mathbb{G}}_m^\diamond
\]
as $\varprojlim_{x\mapsto x^p} R^{\sharp}\cong \varprojlim_{x\mapsto x^p} R$.
This makes $\n{G}^{\diamond}_{\bb{Q}_p^{\prism}}$ isomorphic  to the base change to $\bb{Q}_p^{\prism}$ of the inclusion of the Tate module 
\[
\bb{Z}_p(1)=\ob{ker}( \widetilde{\mathbb{G}}_m^\diamond\to \Gm^\diamond). 
\]
We note that $\bb{Z}_p(1)\subset (1+\widetilde{\bb{D}}^{<1}_{\bb{Q}_p})^{\diamond}$.

Then, the logarithm of the Teichm\"uller lift produces a factorization of commutative group  diamonds over $\FFdisc{\bb{Q}_p^{\diamond}}^{\times,\diamond}$
\begin{equation*}\label{eqLogarithmFactorizationTateTwistAndHKTwist}
\log([-])\colon  \bb{Q}_p^{\diamond} \times_{\bb{F}_p^{\diamond}} \bb{Z}_p(1) \hookrightarrow  \bb{G}_{a,\FFdisc{\bb{Q}_p^{\diamond}}^{\times,\diamond}}(1)_{\ob{Tate}}^{\diamond}\to  \bb{G}_{a,\FFdisc{\bb{Q}_p^{\diamond}}^{\times,\diamond}}\langle 1 \rangle^{\diamond},
\end{equation*}
where  $ \bb{G}_{a,\FFdisc{\bb{Q}_p^{\diamond}}^{\times,\diamond}}(1)_{\ob{Tate}}^{\diamond}$ is the Tate twist obtained as the  base change of $\bb{Q}_p^{\diamond} \times_{\bb{F}_p^{\diamond}} \bb{Z}_p(1) $ along $\bb{Z}_p\to \bb{G}_a^{\diamond}$. Using \Cref{sec:analyt-prism-gm-description-of-g-qp-prism} and \Cref{quotient-perfect}, we can  identify the Tate twist with the perfectoidization of the perfect Breuil--Kisin twist
\[
 \bb{G}_{a,\FFdisc{\bb{Q}_p^{\diamond}}^{\times,\diamond}}(1)_{\ob{Tate}}^\diamond  \cong \bb{G}_{a,\widetilde{\Div^1}}\{1\}_{\perf}^{\diamond}. 
\]
This produces a factorization  of $\bb{G}_a^{\diamond}$-line bundles on $\FFdisc{\bb{Q}_p}^{\times,\diamond}$
\begin{equation}\label{eqMapLineBundlesDifferentTwists1}
\bb{G}_{a,\FFdisc{\bb{Q}_p^{\diamond}}^{\times,\diamond}}(1)_{\ob{Tate}}^{\diamond}\to \bb{G}_{a,\bb{Q}_p^{\prism}}\{1\}^{\diamond} \to  \bb{G}_{a,\FFdisc{\bb{Q}_p^{\diamond}}^{\times,\diamond}}\langle 1 \rangle^{\diamond}. 
\end{equation}

Let $\kappa\colon \bb{Q}_p^{\prism}\to (0,\infty)$ be the radius map of the prismatization.  By \Cref{ConstructionBKtwist}, left map of \eqref{eqMapLineBundlesDifferentTwists1} is an equivalence in the region $(0,p)$, resp.\ the right map is an equivalence in the region $(1,\infty)$.

We have the following descriptions of the group $\n{G}_{\bb{Q}_p^{\prism}}$. 
 
\begin{itemize}

\item  In the  region $(0,p)$, we have a pullback square 
\[
\begin{tikzcd}
\n{G}_{\bb{Q}_p^{\prism}}|_{(0,p)}  \ar[r] \ar[d] & \bb{G}_{a,\bb{Q}_p^{\prism}}\{1\}|_{(0,p)}  \ar[d] \\ 
{ \FFdisc{\Q_p,(0,p)}^{\times,\dR} \times_{\GSpec(\Q_p)} \Z_p(1)^{\sm} }  \ar[r] & \bb{G}_{a,\FFdisc{\bb{Q}_p,(0,p)}^{\times}}(1)_{\Tate}^{\dR}, 
\end{tikzcd}
\]
where the vertical arrows are the natural maps towards the de  Rham stacks.
Since the Tate twist is the base change along $\bb{Z}_p\to \bb{G}_a^{\diamond}$ of the Tate module  $\bb{Z}_p(1)$ via the logarithm map of Teichm\"uller lift, we see that after fixing a trivialization $\bb{Z}_p(1)\cong \bb{Z}_p$,  $\n{G}_{\bb{Q}_p^{\prism}}|_{(0,p)}$ is isomorphic to a base change of $\bb{Z}_p^{\la}= \bb{Z}_p^{\sm}\times_{\bb{G}_a^{\dR}} \bb{G}_a$.

\item  In the region $(1,\infty)$, it is not true that $\n{G}_{\bb{Q}_p^{\prism}}|_{(1,\infty)}$ is a twisted form of $\bb{Z}_p^{\la}$, instead the group law is additionally twisted by the inclusion of the Tate module into the Frobenius twist.  Indeed, $\n{G}_{\bb{Q}_p^{\prism}}|_{(1,\infty)}$   is the extension 
\[
0\to \bb{G}_{a,\FFdisc{\bb{Q}_p,(1,\infty)}^{\times,\dR}}^{\dagger} \to  \n{G}_{\bb{Q}_p^{\prism}}|_{(1,\infty)} \to \FFdisc{\bb{Q}_p^{\diamond}, (1,\infty)}^{\times,\dR} \times_{\GSpec(\Q_p)} \Z_p(1)^{\sm}  \to 0
\]
produced by the composite 
\[
 \FFdisc{\bb{Q}_p^{\diamond}, (1,\infty)}^{\times,\dR} \times_{\GSpec(\Q_p)} \Z_p(1)^{\sm}  \hookrightarrow \bb{G}_{a, \FFdisc{\bb{Q}_p^{\diamond}}^\times, (1,\infty)}(1)_{\ob{Tate}}^{\dR} \to \bb{G}_{a, \FFdisc{\bb{Q}_p^{\diamond}}^\times, (1,\infty)}\langle 1 \rangle^{\dR} \to B \bb{G}_{a,\FFdisc{\bb{Q}_p,(1,\infty)}^{\times,\dR}}^{\dagger}. 
\]

\end{itemize}

We will make the previous description of $\n{G}_{\bb{Q}_p^{\prism}}$ explicit in \Cref{ExampleCyclotomicPrism} for the case of the cyclotomic prism. 
\end{remark}

\begin{example}[Breuil--Kisin twist in the cyclotomic case]\label{ExampleCyclotomicPrism}
Consider the cyclotomic prism $(S,\n{X}_{\cyc},I)$ of \Cref{ExamplesRelativePrismatizations} (1). By \Cref{prop:explicit-div-1} (3) the Breuil--Kisin ideal $\n{I}_{\BK}$ is generated by the function $\frac{\log(q)}{q-1}$, similarly, the perfect Breuil--Kisin ideal is generated by the function $\log(q)$. Therefore, the sequence \eqref{eqMapLineBundlesDifferentTwists1} is the inclusion of line bundles
\[
\bb{G}_{a,\n{X}_{\cyc}}\cdot \log(q) \to  \bb{G}_{a,\n{X}_{\cyc}}\cdot (\frac{\log q}{q-1}) \to \bb{G}_{a,\n{X}_{\cyc}} 
\]
after taking a $p^{-1}$-Frobenius twist. The left map is an equivalence in the region $(0,p)$, as then $q-1$ becomes invertible, and the right map is an equivalence in the region $(1,\infty)$, as then $\log q$ and $(q-1)$ have the same zeroes with multiplicity in this region.

Over the cyclotomic prism $\n{X}_{\cyc}$ we have a trivialization of the Tate twist $q^{\bb{Z}_p^{\sm}}= \n{X}_{\cyc}\times \bb{Z}_p(1)^{\sm}$.   The perfect group $\n{G}_{\widetilde{\Div^1}}$ is then the pullback
\begin{equation}\label{PullbackSquarePrecise}
\begin{tikzcd}
\n{G}_{\widetilde{\Div^1}} |_{\n{X}_{\cyc}} \ar[r] \ar[d] &  \bb{G}_{a,\n{X}_{\cyc}}\cdot \log(q) \ar[d] \\ 
q^{\bb{Z}_p^{\sm}} \ar[r,"\log(-)"] &  \bb{G}_{a,\n{X}_{\cyc}}^{\dR/\n{X}_{\cyc}}\cdot \log(q),
\end{tikzcd}
\end{equation}
or equivalently, the group $\n{G}_{\widetilde{\Div^1}} |_{\n{X}_{\cyc}}=(\bb{Z}_p^{\la})_{\n{X}_{\cyc}}\cdot \log(q)$.  This identifies the group $\n{G}_{\bb{Q}_p^{\prism}}|_{\n{X}_{\cyc}, (0,p)}$ with a twisted form of $\bb{Z}_p^{\la}$ with basis $\log q$. 

On the other hand, in the region $(1,\infty)$ we have a natural map of groups
\[
(\bb{Z}_p^{\la})_{\n{X}_{\cyc}, (1,\infty)}\cdot  \log q \to \n{G}_{\bb{Q}_p^{\prism}}|_{\n{X}_{\cyc}, (1,\infty)} 
\]
that induces an isomorphism on de Rham stacks, and whose overconvergent neighbourhoods at $1$ is the morphism 
\[
\bb{G}_{a,\n{X}_{\cyc}}^{\dagger} \cdot \log q\to \bb{G}_{a,\n{X}_{\cyc}}^{\dagger}. 
\]
Thus, we have a pushout description 
\[
\n{G}_{\bb{Q}_p^{\prism}}|_{\n{X}_{\cyc}, (1,\infty)} = (\bb{Z}_p^{\la})_{\n{X}_{\cyc}, (1,\infty)} \cdot \log q \bigoplus_{ \bb{G}_{a,\n{X}_{\cyc}}^{\dagger} \cdot \log q} \bb{G}_{a,\n{X}_{\cyc}}^{\dagger}. 
\]

\end{example}

  With the analytic prismatization of $\Gm$ or $1+\mathbb{D}^{<1}$ at hand, it is also easy to describe $\Ga^\prism$ as a commutative group stack (but not as a \emph{ring} stack).

  \begin{lemma}
    \label{sec:analyt-prism-gm-description-for-g-a}
    The logarithm $1+\mathbb{D}^{<1}\to \Ga$ induces a commutative diagram
\[\begin{tikzcd}
	{\mathcal{G}_{\Q_p^\prism}} & {(1+\widetilde{\mathbb{D}}^{<1})_{\Q_p^\prism}} & {(1+\mathbb{D}^{<1})^\prism} \\
	{\mathcal{G}_{\Q_p^\prism}[1/p]} & {(1+\widetilde{\mathbb{D}}^{<1})_{\Q_p^\prism}} & {\Ga^{\prism}}
	\arrow[from=1-1, to=1-2]
	\arrow[from=1-1, to=2-1]
	\arrow["\rho", from=1-2, to=1-3]
	\arrow["{=}", from=1-2, to=2-2]
	\arrow["{\log^\prism}", from=1-3, to=2-3]
	\arrow[from=2-1, to=2-2]
	\arrow[from=2-2, to=2-3]
      \end{tikzcd}\]
    of commutative group stacks over $\Q_p^\prism$ with rows given by fiber sequences of spectra. 
\end{lemma}
\begin{proof}
By construction the upper horizontal diagram is a fiber sequence of spectral objects in analytic stacks over $\bb{Q}_{p}^{\prism}$. The map $\log \colon 1+\bb{D}^{<1}\to \bb{G}_a$ has $p$-power torsion kernel, and induces an equivalence after inverting $p$. On the other hand,   $1+\widetilde{\bb{D}}^{<1}$ is a $\bb{Q}_p^{\la}$-vector space so multiplication by $p$ is an isomorphism. The lemma follows. 
\end{proof}

\begin{remark}
  \label{sec:analyt-prism-gm-fundamental-exact-sequence}
  \Cref{sec:analyt-prism-gm-description-for-g-a} can be seen as a rational prismatic variant of the fundamental exact sequence in $p$-adic Hodge theory. Indeed, after  passing to diamonds, \textit{loc.\ cit.} becomes the short exact sequence over $\bb{Q}_p^{\diamond}$
  \[
  0\to \bb{Q}_p(1)^{\diamond}_{\Tate} \to (1+\widetilde{\bb{D}}_{\bb{Q}_p}^{<1})^{\diamond}=\n{BC}(\s{O}(1))^{\diamond}\to \bb{G}_{a,\bb{Q}_p}^{\diamond}\to 0. 
  \]
\end{remark}

\begin{remark}
With some additional effort, one can deduce Cartier self duality of $\mathbb{G}_a^\prism$  (up to a twist and shift)  from \Cref{sec:analyt-prism-gm-description-for-g-a}. More precisely, using the general formalism of Cartier duality in Gestalten \cite{GestaltenScholze},  the basic Cartier duality between discrete and profinite groups, and the Cartier duality between $\bb{Z}_p^{\la}$ and $1+\bb{D}^{<1}_{\bb{Q}_p}$ of \cite{RCCartierGerbe}, one can extend Cartier duality to general Banach--Colmez spaces. We do not try to formalize this idea in this paper. 
\end{remark}

Using \Cref{sec:analyt-prism-gm-description-of-g-qp-prism} and \Cref{ExampleCyclotomicPrism}  we get an explicit description of the base change of $\Gm^{\widetilde{\Div^1}}$ in terms of the cyclotomic  prism.

\begin{theorem}
  \label{sec:analyt-prism-gm-description-over-cyclotomic-prism}
Consider  $\Gm^{\widetilde{\Div^1}/\n{X}_{\cyc}}$,  the perfect prismatization of $\bb{G}_m$  over the cyclotomic prism.  Let 
\[
0\to \n{G}_{\widetilde{\Div^1}} |_{\n{X}_{\cyc}}  \to \widetilde{\mathbb{G}}_{m,\n{X}_{\cyc}} \to \Gm^{\widetilde{\Div^1}/\n{X}_{\cyc}} \to 0
\]
be the (perfection of the) fiber  sequence of \Cref{sec:analyt-prism-gm-analytic-prismatization-of-gm} (2) where  $\widetilde{\mathbb{G}}_{m,\n{X}_{\cyc}}$ has coordinate $T$.  

The following statements hold.
  \begin{enumerate}
 \item  The section $(q,q^{1/p}, q^{1/p^2},\ldots)$ of $\widetilde{\mathbb{G}}_{m,\n{X}_{\cyc}}\to \n{X}_{\cyc}$ produces a natural equivalence 
 \[
 \bb{Z}_{p,\n{X}_{\cyc}}^{\la}\xrightarrow {\sim} \n{G}_{\widetilde{\Div^1}} |_{\n{X}_{\cyc}}  \to  \widetilde{\mathbb{G}}_{m,\n{X}_{\cyc}}. 
 \]
By an abuse of notation, we write  $\n{G}_{\widetilde{\Div^1}} |_{\n{X}_{\cyc}} = q^{\bb{Z}_p^{\la}}$.

\item The natural action of $\bb{Z}_{p,\n{X}_{\cyc}}^{\la}\cong q^{\bb{Z}_p^{\la}}$ on $\widetilde{\mathbb{G}}_{m,\n{X}_{\cyc}}$ is given by
\[
a\cdot T^{1/p^n}= q^{a/p^n}T^{1/p^n}
\]
for $a\in \bb{Z}_p^{\la}$. In particular, we have an equivalence 
\[
\Gm^{\widetilde{\Div^1}/\n{X}_{\cyc}} \cong \widetilde{\mathbb{G}}_{m,\n{X}_{\cyc}}/\bb{Z}_p^{\la}. 
\]

\item We have  an equivalence of Gelfand stacks
\[
\Gm^{\widetilde{\Div^1}} \cong \widetilde{\mathbb{G}}_{m,\n{X}_{\cyc}}/(\bb{Z}_p^{\times, \la}\ltimes \bb{Z}_p^{\la} )
\] 
where $\bb{Z}_p^{\times,\la}$ acts on $\bb{Z}_p^{\la}$ by multiplication, and under the fixed coordinates $(q^{\flat}=(q^{1/p^n}), T^{\flat}=(T^{1/p^n}))_{n\geq 0}$ the action of $(a,b)\in \bb{Z}_p^{\times, \la}\ltimes \bb{Z}_p^{\la} $ is via 
\begin{equation}\label{eqMapCocycleGmCoordinates}
(a,b)\cdot (q^{\flat}, T^{\flat})= (q^{\flat, a}, q^{\flat,b}T^{\flat}). 
\end{equation}
       \end{enumerate}
  \end{theorem}

  \begin{proof}
    For part (1), it suffices to prove the claim for $1+\widetilde{\bb{D}}^{<1}$ instead of $\widetilde{\mathbb{G}}_m$. In that case the statement follows by taking prismatic logarithms  and from  the pullback square \eqref{PullbackSquarePrecise}.  The explicit description of the action in part (2) follows directly from (1) since the action of $\n{G}_{\widetilde{\Div^1}} |_{\n{X}_{\cyc}} $ on $\widetilde{\mathbb{G}}_{m,\n{X}_{\cyc}}$ is by multiplication.
    The explicit presentation of the perfect  prismatization  of $\Gm$ over $\n{X}_{\cyc}$ follows from \Cref{sec:analyt-prism-gm-analytic-prismatization-of-gm} (1).  

Finally, we prove the presentation of (3).  We use the notation of  \Cref{ExampleFrobeniusTwistCyc}. We have a natural morphism of Fargues--Fontaine discs
\[
\FFdisc{\bb{Q}_p^{\cyc}}^{\times}\to \n{X}_{\cyc}
\]
sending $q\mapsto \epsilon^{1/p}$ and the ideal $(\frac{q^p-1}{q-1})$ to the ideal defining the untilt  $\bb{Q}_p^{\cyc}$ of $\bb{Q}_p^{\cyc,\flat}$.  This morphism of $\delta$-spaces induces an equivalence on uniform completions when restricted to affinoids, and therefore is an equivalence on diamonds. Furthermore, the natural action of $\bb{Z}_p^{\times,\la}$ on the $\n{X}^{\cyc}$ is compatible with the natural  $\Gamma$-action on $\FFdisc{\bb{Q}_p^{\cyc}}^{\times}$ via the cyclotomic character $\chi_{\cyc}\colon \Gamma\xrightarrow{\sim} \bb{Z}_p^{\times}$. 

Similarly, the obvious map of $\delta$-spaces  
\[
\widetilde{\mathbb{G}}_{m,\FFdisc{\bb{Q}_p^{\cyc}}^{\times}}^u\to \widetilde{\mathbb{G}}_{m,\n{X}_{\cyc}}
\]
is an equivalence after passing to uniform completions on affinoids, and equivariant for the action of $\bb{Z}_p$ and the Galois action of the map  $\widetilde{\mathbb{G}}_{m,\bb{Q}_p^{\cyc}}\to \mathbb{G}_{m,\bb{Q}_p^{\cyc}}$.
The previous implies that if the action of $\bb{Z}_p^{\times,\la}$ on $\n{X}_{\cyc}$ and the action of $\bb{Z}_p^{\la}$ on  $\widetilde{\mathbb{G}}_{m,\n{X}_{\cyc}}$ lift to an action of $\bb{Z}_p^{\times,\la}\ltimes \bb{Z}_p^{\la}$  compatible with that on $\widetilde{\mathbb{G}}_{m,\FFdisc{\bb{Q}_p^{\cyc}}^{\times}}^u$,  it must be given by the formula \eqref{eqMapCocycleGmCoordinates}, namely, this follows from the fact that  the map of (simplicial) Berkovich spaces
\[
 \widetilde{\mathbb{G}}_{m,\FFdisc{\bb{Q}_p^{\cyc}}^{\times}}^u \rtimes (\bb{Z}_p^{\times,\la}\ltimes \bb{Z}_p^{\la})^{u,\bullet}\to \widetilde{\mathbb{G}}_{m,\n{X}_{\cyc}} \times (\bb{Z}_p^{\times,\la}\ltimes \bb{Z}_p^{\la})^{\bullet} 
\]
  induces an injection of structural  sheaves on the underlying topological space.

It is left to show that $\bb{Z}_p^{\times,\la}\ltimes \bb{Z}_p^{\la}$ acts on $\widetilde{\mathbb{G}}_{m,\n{X}_{\cyc}}$ compatible with the natural actions on classical Fargues--Fontaine curves, but this follows from the fact that   $\n{X}_{\cyc}\to \widetilde{\Div^1}$ is a $\bb{Z}_p^{\times,\la}$-torsor by \Cref{prop:explicit-div-1} (4), and that the isomorphism of  $\bb{Z}_p^{\la}\times \n{X}_{\cyc}\cong q^{\bb{Z}_p^{\la}}= \n{G}_{\widetilde{\Div^1}} |_{\n{X}_{\cyc}}$  endows $\bb{Z}_p^{\la}$ with the $\bb{Z}_p^{\times,\la}$-action given by multiplication.  Thus, the equivalence 
\[
\widetilde{\mathbb{G}}_{m,\widetilde{\Div^1}}\times_{\widetilde{
\Div}^1}\n{X}_{\cyc}= \widetilde{\mathbb{G}}_{m,\n{X}_{\cyc}}
\]
endows $\widetilde{\mathbb{G}}_{m,\n{X}_{\cyc}}$ with an action of $\bb{Z}_p^{\times}$ trivial on the copy $\widetilde{\mathbb{G}}_m$,  equivariant for the projection $\widetilde{\mathbb{G}}_{m,\n{X}_{\cyc}}\to \n{X}_{\cyc}$, and equivariant for the action of $\bb{Z}_p^{\la}\times \n{X}_{\cyc}$ with respect to the $\bb{Z}_p^{\times,\la}$-multiplication on  $\bb{Z}_p^{\la}$. The theorem follows.
    \end{proof}
   
 \begin{remark}
    \label{sec:analyt-prism-gm-equivariance-for-frobenius}
    By the construction of the Teichm\"uller map, the morphism
    \[
      \widetilde{\mathbb{G}}_{m,\mathcal{X}_{\cyc}}\to \Gm^{\widetilde{\Div^1}}\times_{\widetilde{\Div^1}} \mathcal{X}_{\cyc}
    \]
    is Frobenius equivariant, if we set $\varphi_{\widetilde{\mathbb{G}}_{m,\mathcal{X}_{\cyc}}}\colon \widetilde{\mathbb{G}}_{m,\mathcal{X}_{\cyc}}\to \widetilde{\mathbb{G}}_{m,\mathcal{X}_{\cyc}}$ as $([p],\varphi_{\mathcal{X}_{\cyc}})$, where $[p]\colon \widetilde{\mathbb{G}}_m\to \widetilde{\mathbb{G}}_m$ is the multiplication by $p$.
    We see furthermore that the map 
    \[
    \Z_p^\la\times_{\GSpec(\Q_p)}\mathcal{X}_{\cyc} \to \widetilde{\mathbb{G}}_{m,\mathcal{X}_{\cyc}}
    \] is $\varphi$-equivariant if $\varphi$ acts trivially on $\Z_p^\la$.
    We can deduce that
    \[
      \Gm^{\Div^1}\cong (\widetilde{\mathbb{G}}_{m,\mathcal{X}_{\cyc}})/(\varphi^\Z_{\widetilde{\mathbb{G}}_{m,\mathcal{X}_{\cyc}}}\times \Z_p^{\times,\la}\ltimes \Z_p^\la).
    \]
  \end{remark}

\begin{remark}[Towards (geometric) Sen theory]
The definition of the analytic prismatization above has been heavily inspired by the constructions of Bhatt--Lurie and Drinfeld for $p$-adic formal schemes. But the explicit description of \Cref{sec:analyt-prism-gm-description-over-cyclotomic-prism} lets appear a relation with a different chapter of $p$-adic Hodge theory: (geometric) Sen theory. Indeed, let $X=\Gm^d$, $d\geq 0$, be the $d$-dimensional rigid analytic multiplicative group and $\widetilde{X}$ the perfected (in Gelfand stacks) $d$-dimensional rigid analytic multiplicative group. The $p$-adic twistor stack of $X$ can be described as 
$$
X^{\Div^1}\cong (\widetilde{X}_{\mathcal{X}_{\cyc}})/(\varphi^\Z_{\widetilde{X}_{\mathcal{X}_{\cyc}}}\times \Gamma_d^\la),
$$
with $\Gamma_d:= \Z_p^{\times,\la}\ltimes \Z_p^{d, \la}$ (for the action of $\Z_p^{\times,\la}$ on $\Z_p^{d,\la}$ by multiplication). Using computations in geometric Sen theory (e.g.\ \cite[Theorem C]{PoratLocAnVetorsFF} when $d=0$), the right hand side can also be rewritten
$$
(\widetilde{X}_{\mathcal{X}_{\cyc}})/(\varphi^\Z_{\widetilde{X}_{\mathcal{X}_{\cyc}}}\times \Gamma_d^\la) = \FF_{\widetilde{X}_{\Q_p^{\rm cyc}}^{d,\diamond}}^\la /\Gamma_d^\la,
$$
even though we postpone the definition of the right hand side to the next section (basically, $\FF_{\widetilde{X}_{\Q_p^{\rm cyc}}^{d,\diamond}}^\la$ is obtained by gluing spectra of rings of locally analytic vectors on affinoid pieces of $\FF_{\widetilde{X}_{\Q_p^{\rm cyc}}^{d,\diamond}}$). This is to be compared to the description 
$$
\FF_{X^\diamond} = \FF_{\widetilde{X}_{\Q_p^{\rm cyc}}^{d,\diamond}}/\Gamma_d^{\rm cont}
$$
and (passing to de Rham stacks)
$$
\FF_{X^\diamond}^\dR = \FF_{\widetilde{X}_{\Q_p^{\rm cyc}}^{d,\diamond}}^\dR/\Gamma_d^{\rm sm}.
$$
The diagram \eqref{composition-for-twistor} becomes the diagram
$$
\FF_{\widetilde{X}_{\Q_p^{\rm cyc}}^{d,\diamond}}/\Gamma_d^{\rm cont} \xto{\alpha_X} \FF_{\widetilde{X}_{\Q_p^{\rm cyc}}^{d,\diamond}}^\la /\Gamma_d^\la \xto{\beta_X}  \FF_{\widetilde{X}_{\Q_p^{\rm cyc}}^{d,\diamond}}^\dR/\Gamma_d^{\rm sm}.
$$
We will make heavy use of this explicit description of $\alpha_X$ in the next section. This description also allows to analyze the difference between the source and target of $\beta_X$ in terms of the actions of Lie algebroids as follows. 
  
  The Lie algebra attached to the quotient $\FF_{\widetilde{X}_{\Q_p^{\rm cyc}}^{d,\diamond}}^\la \to \FF_{\widetilde{X}_{\Q_p^{\rm cyc}}^{d,\diamond}}^\dR$  is its tangent space $\mathcal{T}_{\FF_{\widetilde{X}_{\Q_p^{\rm cyc}}^{d,\diamond}}^\la}$, which has a basis given by the derivations $q\partial_q, T_1\partial_{T_1},\ldots, T_d\partial_{T_d}$.
On the other hand, the action of the Lie algebra of $\Gamma_d^\la$ has a basis given by the tangent vectors $(\log q) q\partial_q$ , $ (\log q)T_1\partial_{T_1},\ldots, (\log q)T_d\partial_{T_d}$ (using the explicit description in \Cref{sec:analyt-prism-gm-description-over-cyclotomic-prism}, the factor $\log q$ appears as the derivative of $a\mapsto q^a$). 
As we see, there is a short exact sequence  of sheaves  on $\FF_{\widetilde{X}_{\Q_p^{\rm cyc}}^{d,\diamond}}^\la$
  \[
  0\to \Lie \Gamma_d \otimes_{\Q_p} \s{O}_{\FF_{\widetilde{X}_{\Q_p^{\rm cyc}}^{d,\diamond}}^\la}\to \mathcal{T}_{\FF_{\widetilde{X}_{\Q_p^{\rm cyc}}^{d,\diamond}}^\la}\to \Lie\Gamma_d \otimes \s{O}_{\FF_{\widetilde{X}_{\Q_p^{\rm cyc}}^{d,\diamond}}^\la}/ (\log q)  \to 0. 
  \]
  In particular, the action of $\Lie \Gamma_d$ has logarithmic poles at the Hodge--Tate point, and its reduction modulo $\log q$ is trivial.
  This residual linear action on the Hodge--Tate stack is the local  arithmetic  and geometric Sen operators, as will be discussed in more detail in \Cref{sec:comp-proet-cohom}.
\end{remark}

\newpage
\section{Sen theory and decompletions by locally analytic vectors}
\label{sec:sen-theory}
In the next section, \Cref{sec:comp-proet-cohom}, we will prove \Cref{thm:equivalence-perfect-complexesintro}. As already sketched in the introduction, the proof will be by reduction, using the local explicit descriptions obtained at the end of the previous section, to the two central questions of (geometric) Sen theory:
\begin{enumerate}
\item[(a)] decompletions by (derived) locally analytic vectors, and
\item[(b)] vanishing of higher locally analytic vectors.
\end{enumerate}

This section investigates these questions and establishes results of interest independently of the theory of analytic prismatization. Given its technical nature, the reader interested only in the theory of prismatization is invited to skip this section and take its main results as a blackbox.

\subsection{Main result}

Given a $p$-adic Lie group $G$ we let $G^{\rm cont}$, $G^{\la}$ and $G^{\sm}$ denote the groups in Berkovich spaces given by the sheaf of continuous, locally analytic and locally constant functions on $G$ valued in $\Q_p$. We have the following decompletion theorem for Berkovich spaces with an action of $G^{\rm cont}$.

\begin{theorem}\label{TheoDecompletion}
  Let $X$ be a Berkovich space endowed with a continuous action of a $p$-adic Lie group $G$, that is, endowed with an action of the Gelfand group stack $G^{\rm cont}$.
  Let $h=h_X\colon X/G^{\rm cont}\to |X|_{\Betti}/G^{\la}$ be the natural induced map of Gelfand stacks, where the action of $G^{\la}$ on $|X|_{\Betti}$ arises from the morphism of groups $G^{\la}\to G^{\sm}$ and the natural action of $G^{\sm}$ on $|X|_{\Betti}$ induced by the action on topological spaces.\footnote{That this action has to be smooth is implied by the fact that $|X|_\Betti=|X^\dR|_\Betti$ and $G^\sm=\underline{G}^\dR$.}
  The following hold:
  
\begin{enumerate}

\item The map $h$ is $0$-proper,  and hence we have a fully faithful  pullback functor 
\[
\ob{Mod}_{h_* \s{O}_{X/G^{\rm cont}}} \ob{D}(|X|_{\Betti}/G^{\la}) \hookrightarrow \ob{D}(X/G^{\rm cont}).
\] We denote $\s{O}_X^{\la}:= h_* \s{O}_{X/G^{\rm cont}}$  and call it the sheaf of locally analytic functions on $X$ with respect to $G$.
\item Suppose that $X=\GSpec(R)$ is affinoid. Assume moreover that $R$ is \textit{classically decompleted}, in the sense of \Cref{Definition-Classically-decompleted}.
  Let $R^{\la}\in \ob{CAlg}(\ob{D}(BG^{\la}))$ denote the algebra of locally analytic vectors of $R$ (in particular, $R^\la$ is a Gelfand algebra with $R^{\la,u}=R^u$).
  Then the natural map $R^{\la}\to \s{O}_{X}^{\la}$ of commutative algebras in $|X|_{\Betti}$ is an equivalence. In particular, the natural pullback functor
\[
\ob{Mod}_{R^{\la}}(\ob{D}(BG^{\la}))\xrightarrow{\sim} \ob{Mod}_{\s{O}_X^{\la}} \ob{D}(|X|_{\Betti}/G^{\la})
\] 
is an isomorphism.  Furthermore, any $G^{\cont}$-equivariant affinoid Berkovich-\'etale map $\GSpec (R')\to \GSpec (R)$ is classically decompleted. 

\item Assume that $X$ is  locally strongly decompleted, in the sense of \Cref{AssumptionKeyTheoremDecompletionPerfectComplexes}. Then, for any $G^{\ob{cont}}$-equivariant Berkovich-\'etale morphism $X'\to X$, $X^{\prime}$ admits a  strict closed  affinoid  cover  $\{Z_i\}_{i\in I}$ which is classically decompleted with respect to the action of the (necessarily open) stabilizers $G_{Z_i}\subset G$.  In particular, $\s{O}^{\la}_{X^{\prime}}\in \ob{CAlg}(\ob{D}(|X^{\prime}|)$ is the structural sheaf of a Berkovich space $X^{\prime,\la}$ with underlying topological space $|X^{\prime,\la}|=|X^{\prime}|$. 
Furthermore, the pullback functor 
\[
h^*_{X'}\colon \ob{D}(X^{\prime,\la}/G^{\la})\hookrightarrow \ob{D}(X'/G^{\rm cont})
\]
induces an equivalence of perfect complexes 
\[
\Perf(X^{\prime,\la}/G^{\la})\xrightarrow{\sim} \Perf(X^{\prime}/G^{\rm cont}). 
\]
\end{enumerate}
\end{theorem}

 The notion of a \textit{classical decompletion} that appears in (2) is defined for an action of a $p$-adic Lie group on a Berkovich space and essentially asks for vanishing of the higher locally analytic vectors, together with density of locally analytic vectors of  the underlying classical ring. 
 The assumption in (3) is stronger than admitting a classical decompletion. Hence, we see that our discussion of decompletions and vanishing of higher locally analytic vectors are closely intertwined.

In  \Cref{subsec-Sen-theory-revisited}, we verify that for the key examples appearing in $p$-adic Hodge theory, the classical Sen method using traces shows that the assumption of \Cref{TheoDecompletion} (3) is satisfied. \Cref{thm:summary-main-results-sen} summarizes our conclusions and will be the relevant statement for the applications given in \Cref{sec:comp-proet-cohom}.

\begin{remark} \label{rk:reformulation-theo-decompletion}
We explain in more informal terms the content of \Cref{TheoDecompletion}.  If $X$ is (locally) classically decompleted, $\s{O}^\la_X$ is locally on $|X|_\Betti$ represented by a Gelfand ring.  In this case, one can use the relative spectrum (cf. \Cref{sec:deriv-berk-spac-relative-gelfand-spectrum}) of $\s{O}^\la_X$ to define the Berkovich space
  \[
    X^\la:=\underline{\GSpec}_{|X|_\Betti}(\s{O}^\la_X).
  \]
 Then part (1) of \Cref{TheoDecompletion} can be rephrased as saying that $f\colon X/G^\cont \to X^\la/G^\la$ is $0$-proper and that the pullback of quasi-coherent sheaves along $f$ is fully faithful.  Part (2) says that $X^{\la}$ and $X$ have the same underlying topological space, and that  part (1) holds for  Berkovich-\'etale localizations of $X$. Finally,   if $X$ is locally strongly decompleted, part (3) says that $X$ is locally classically decompleted and that  pullback of quasi-coherent sheaves along $f$   gives rise to an equivalence of perfect objects
 \[
 \Perf(X^{\la}/G^{\la})  \xrightarrow{\sim} \Perf(X/G^{\rm cont}).
 \]
\end{remark}

\begin{remark}\label{Remark-IntegralLocallyAnalytic}
  The whole discussion in this section can be adapted to mixed characteristic locally analytic representations, giving an analogue of  \Cref{TheoDecompletion} for  the integral theory.
  We decided to restrict to the characteristic zero situation in order to keep working with the theory of Gelfand stacks over $\Q_p$ following \cite{drff_stacks}. 
\end{remark}

\begin{remark}
We could only prove the equivalence of dualizable objects in \Cref{TheoDecompletion} (3) under the assumption that $X$ is locally strongly decompleted. We do not know if the result could be true for any continuous action of a $p$-adic Lie group on a Berkovich space, without any assumption. 
\end{remark}

The proof of \Cref{TheoDecompletion} spreads out over various subsections below: see \Cref{PropositionPart1MainTheorem}, \Cref{CorollaryPart2MainDecompletionTheorem}, \Cref{TheoremDecompletionPerfectComplexesII}. In \Cref{subsub-semilinear-solid-rep} we clarify notions surrounding locally analytic representations. \Cref{subsec-Sen-theory-revisited} revisits the examples coming from (geometric) Sen theory. In \Cref{sec:decompl-tors-1}, we discuss more generally the decompletion of torsors over $X/G^\cont$ under $p$-adic Lie groups.

\subsection{Semilinear solid representations and locally analytic representations}\label{subsub-semilinear-solid-rep}
Let $G$ be a  $p$-adic Lie group.
In this subsection, we clarify different notions regarding solid and locally analytic representations of $G$.
We let $\cat{Rep}_{\Z}^\sol(G)=\ob{D}(\Z_{\sol}[G])=\mathrm{LMod}_{\Z_{\sol}[G]}(\ob{D}(\Z))$ be the category of $\Z_{\sol}[G]$-modules on solid $\Z$-modules for the solid tensor product $-\otimes_{\Z_{\sol}}-$.
We call objects in $\ob{D}(\Z_{\sol}[G])$ \emph{$G$-equivariant solid $\Z$-modules}, or just \emph{solid $G$-representations}, and equip $\ob{D}(\Z_{\sol}[G])$ with the symmetric monoidal structure $-\otimes_{\Z_{\sol}}-$ (using that $\Z_{\sol}[G]$ is a Hopf algebra), so that the forgetful functor $\ob{D}(\Z_\sol[G])\to \ob{D}(\Z_\sol)$ is symmetric monoidal.

Let $R$ be a (not necessarily connective) commutative $G$-equivariant solid $\Z$-algebra, i.e., $R\in \ob{CAlg}(\ob{D}(\Z_{\sol}[G]))$.
We let $\cat{Rep}_R^{\sol}(G)=\mathrm{Mod}_R(\ob{D}(\Z_{\sol}[G]))$ denote the category of $R$-semilinear solid $G$-representations on solid $\Z$-modules and $\Lambda_R(G)=R\otimes_{\Z} \Z_{\sol}[G]  $.
We equip $\cat{Rep}_R^\sol(G)$ with the symmetric monoidal structure $-\otimes_{R}-$. The following construction endows $\Lambda_{R}(G)$ with a natural associative algebra structure which generalizes the classical semilinear group-algebra such that $(a[g])(b[h])= ag(b) [gh]$ for $a,b\in R$ and $g,h\in G$.

\begin{construction}\label{LemmaGroupAlgebra}
By definition, the $G$-equivariant ring  $R$ gives rise to a commutative algebra object in $\cat{Rep}^{\sol}_{\Z}(G)$, and one has $\cat{Rep}^{\sol}_{R}(G)=\mathrm{Mod}_R(\cat{Rep}^{\sol}_{\Z}(G))$.
  Note also that $\cat{Rep}^{\sol}_{\Z}(G)=\mathrm{LMod}_{\Z_{\sol}[G]}(\ob{D}(\Z_{\sol}))$.
Consider the composite of forgetful functors 
\begin{equation}\label{eqComposite1LambdaGForgetful}
F\colon \cat{Rep}^{\sol}_{R}(G)\to \cat{Rep}^{\sol}_{\Z}(G)\to \ob{D}(\Z_{\sol}).
\end{equation}
This composite preserves limits and colimits,  is $\Z_{\sol}$-linear and conservative, and by Lurie's monadicity theorem we have a natural equivalence 
\begin{equation*}\label{eqForgetfulEquivariantGtoQ_p}
 \mathrm{Rep}^{\sol}_{R}(G)\xrightarrow{\sim} \ob{LMod}_A(\ob{D}(\Z))
\end{equation*}
where $A\in \mathrm{Alg}(\ob{D}(\Z))$ is an associative algebra in solid abelian groups. The underlying solid abelian group of $A$ is computed as the monad of $F$ evaluated at $\Z$, and is given by $R\otimes_{\Z}\Z_{\sol}[G]=\Lambda_R(G)$. The previous produces a natural associative algebra structure  on $A=\Lambda_R(G)$.  The expression of $F$ as composite of forgetful functors \eqref{eqComposite1LambdaGForgetful} produces a natural morphism of algebras  $\Z_{\sol}[G]\to \Lambda_{R}(G)$.  Similarly, the presentation of $F$ as composite of forgetful functors
\[
F\colon \cat{Rep}^{\sol}_{R}(G)\to \ob{Mod}_R(\ob{D}(\Z))\to \ob{D}(\Z)
\]
produces the morphism of algebras $R\to \Lambda_{R}(G)$.  Finally, the classical compatibility of the $R$-linear structure and the action of $G$ depicted by the formula $(a[g])(b[h])= ag(b)[gh]$ is implied by the fact that, by construction,  we have a $G$-equivariant isomorphism $\Lambda_{R}(G)\cong R\otimes_{\Z} \Z_{\sol}[G]$ where  $\Z_{\sol}[G]$ is endowed with the left $G$-action, and that the left $R$-linear structure on  $\Lambda_{R}(G)$ arises precisely from the morphism of algebras $R\to \Lambda_{R}(G)$ constructed above.  We observe that the previous discussion can be applied to any (light) locally profinite group $\Pi$. 
\end{construction}

\subsubsection{$\Pi^{\mathrm{cont}}$-representations}

Let $\Pi$  be a (light) locally profinite group and  let $\Pi^{\ob{cont}}$ be the Berkovich space defined by the sheaf of $\bb{Q}_p$-valued continuous functions of $\Pi$.  A natural category of $\Pi$-representations that appears in \cite{ALBM-FF} is the category $\ob{D}(B\Pi^{\ob{cont}})$ of solid quasi-coherent sheaves on the classifying stack $B\Pi^\cont=\GSpec(\Q_p)/\Pi^\cont$.

\begin{lemma}\label{LemmaFirstPropertiesGcont}
Let $\Pi$ be a (light) locally profinite group as before and denote $f\colon  \GSpec(\bb{Q}_p) \to B\Pi^\cont$. The following holds:
\begin{enumerate}

\item  Suppose that $\Pi$ is profinite, the  $*$-pullback along $f$ induces an equivalence of categories 
\[
\ob{D}(\GSpec (B\Pi^\cont) \xrightarrow{\sim} \ob{cRMod}_{C(\Pi,\bb{Q}_p)}(\ob{D}(\bb{Q}_p)),
\]
where the right hand side is the category of right $C(\Pi,\bb{Q}_p)$-comodules in solid $\bb{Q}_p$-vector spaces. 
 
\item The category $\ob{D}(B\Pi^\cont)$ has a natural $t$-structure making $f^*$ a  $t$-exact functor. Furthermore,  $\ob{D}(B\Pi^\cont)$ is the left completed derived category of its heart. 

\item Suppose that $\Pi$ is profinite and has finite $p$-cohomological dimension. Then the map $f\colon \GSpec (\bb{Q}_p)\to B\Pi^\cont$ is proper and  descendable. In particular, $B\Pi^{\ob{cont}}$ is proper over $\GSpec (\bb{Q}_p)$. 

\item Suppose that $\Pi$ has virtually  finite $p$-cohomological dimension,\footnote{That is, there is an open subgroup $\Pi_0\subset \Pi$ with finite $p$-cohomological dimension.} then $B\Pi^{\ob{cont}}\to \GSpec(\bb{Q}_p)$ is $!$-able. 

\end{enumerate}
\end{lemma}
\begin{proof}
\begin{enumerate}

\item This is a special case of  \cite[Proposition 3.3.13 (3)]{RCBBCourse2026}.

\item  Let $\Pi_0\subset \Pi$ be a compact open subgroup.  The space $C(\Pi_0,\bb{Q}_p)$ is a flat solid $\bb{Q}_p$-module which makes the comonad   of $f^*$ flat.  Then, part (2) for $\Pi_0$ follows from  \cite[Lemma 4.3.5 (3)]{jacinto2023solid}. To prove the statement for $\Pi$, consider the map of classifying stacks $g\colon B\Pi_0\to  B\Pi$, the map $g$ is \'etale which makes $g^*$ monadic with monad given by $T=g^*g_{\sharp} $.  Suppose that the monad is $t$-exact, then part (2) will follow from \cite[Lemma 4.3.5 (2)]{jacinto2023solid} and the case of $\Pi_0$.  To show $t$-exactness of $T$, we use proper base change along the pullback diagram 
\[
\begin{tikzcd}
\Pi_0\backslash \Pi/\Pi_0 \ar[r,"{\ob{pr}_2}"] \ar[d,"{\ob{pr}_1}"]  & B\Pi_0 \ar[d,"g"]  \\
B\Pi_0 \ar[r,"g"] & B\Pi.
\end{tikzcd}
\]
Indeed, we have  $g^*g_{\sharp}\cong \ob{pr}_{2,\sharp}\ob{pr}_{1}^{*}$, and to prove that $T$ is $t$-exact  it suffices to show that   $ \ob{pr}_{2,\sharp}$ and  $\ob{pr}_{1}^{*}$ are $t$-exact. We prove that $\ob{pr}_{2}^*$ and $\ob{pr}_{2,\sharp}$ are $t$-exact, the claim for $\ob{pr}_1$ being  symmetric. The map $\ob{pr}_{2}$ can be written as a disjoint union of double cosets
\[
\bigsqcup_{g\in  [\Pi_0 g\Pi_0] } \GSpec(\bb{Q}_p)/ \big((g^{-1}\Pi_0g)\cap \Pi_0 \big)^{\ob{cont}} \to B\Pi_0^{\ob{cont}}.
\] 
Thus, flatness of $\ob{pr}_{2}^*$ and $\ob{pr}_{2,\sharp}$ are reduced to flatness of the analogue functors along the map $h\colon B\Pi^{\ob{cont}}_1\to B\Pi^{\ob{cont}}_0$ where $\Pi_1\subset \Pi_0$ is an open compact subgroup. In this case $h^*$ is a restriction functor which is clearly $t$-exact, and $h_{\sharp}=h_!$ is compact induction which is also clearly $t$-exact. This finishes the proof of part (2).

\item  If $f$ is proper and descendable, it will follow from \cite[Lemma 4.6.3 (ii')]{heyer20246functorformalismssmoothrepresentations} that $f\colon B\Pi^{\ob{cont}}\to \GSpec(\Q_p)$ is proper.  Properness of  $f$ is proper follows from the fact that its fiber is the map $\Pi^{\ob{cont }}\to \GSpec(\bb{Q}_p)$ which is $0$-proper and $0$-truncated.  It is left to prove that $f_* \bb{Q}_p$ is descendable. This is shown in the proof of \cite[Proposition 4.3.9]{ALBM-FF}.

\item Let $\Pi_0\subset \Pi$ be an open subgroup with finite cohomological dimension. The map $B\Pi_0^{\ob{cont}}\to B\Pi^{\ob{cont}}$ is an \'etale cover. By (3) the stack $B\Pi_0^{\ob{cont}}$ is $!$-able; $!$-ability of $B\Pi$ follows from \cite[Theorem 3.4.11 and Lemma 4.7.1]{heyer20246functorformalismssmoothrepresentations}. \qedhere
\end{enumerate}
\end{proof}

Next, we construct a symmetric monoidal functor from $\Pi^{\ob{cont}}$-representations to solid $\Pi$-representations which is compatible with forgetful functors to solid $\bb{Q}_p$-modules.

\begin{construction}\label{ConstructionSymmetricMonoidalFunctorContSolidReps}
Let $\Pi$ be a (light) locally profinite group. 
By \cite[Lemma 10.3]{mann2022zl} we can identify solid $\bb{Q}_p$-sheaves on $(B\Pi)_{\proet}$ with $\cat{Rep}^{\sol}_{\bb{Q}_p}(\Pi)$; this presentation also produces the natural symmetric monoidal structure of $\cat{Rep}^{\sol}_{\bb{Q}_p}(\Pi)$ as sheaves on the classifying stack $(B\Pi)_{\proet}$.    Consider the \v{C}ech nerve $\Pi^{\bullet}$ of $*\to */\Pi$. By \cite[Lemma 3.29]{anschutz2024descent}, for a light profinite set $S$ there is an equivalence of categories 
\[
\ob{D}_{\ob{nuc}}(S,\bb{Q}_p)\cong \ob{D}_{\ob{nuc}}( C(S,\bb{Q}_p)) =: \ob{Nuc}(S^{\ob{cont}})
\]
between nuclear solid $\bb{Q}_p$-valued sheaves on $S_{\proet}$ and nuclear solid $ C(S,\bb{Q}_p)$-modules. Taking  left Kan extensions along open immersions, for $S$ a locally profinite set we have a natural equivalence of categories 
\[
\ob{D}_{\ob{nuc}}(S,\bb{Q}_p)\cong \ob{Nuc}(S^{\ob{cont}})
\]
where $\ob{Nuc}(S^{\ob{cont}})\subset \ob{D}(S^{\ob{cont}})$ is the full subcategory of sheaves whose pullback to compact open subspaces are nuclear.   Applying this to $\Pi^{\bullet}$,   we have a   fully faithful symmetric monoidal morphism   of cosimplicial categories
\begin{equation}\label{eqInclusionNuclearComods1}
\ob{D}_{\ob{nuc}}(\Pi^{\bullet},\bb{Q}_p) \subset \ob{D}(\Pi^{\bullet,\ob{cont}})
\end{equation}
inducing a natural equivalence of cosimplicial symmetric monoidal categories 
\[
\ob{D}_{\ob{nuc}}(\Pi^{\bullet},\bb{Q}_p) \otimes_{\ob{Nuc}(\bb{Q}_p)} \ob{D}(\bb{Q}_p)\xrightarrow{\sim} \ob{D}(\Pi^{\bullet,\ob{cont}}).
\] 
 Since $\ob{D}(\Pi^{\bullet}, \bb{Q}_p)$ is a cosimplicial $\ob{D}(\bb{Q}_p)$-algebra in $\cat{Pr}$, the inclusion $\ob{D}_{\ob{nuc}}(\Pi^{\bullet},\bb{Q}_p)\subset \ob{D}(\Pi^{\bullet},\bb{Q}_p)$ produces a morphism of cosimplicial symmetric monoidal categories
\begin{equation}\label{eqSymmetricMonoidalDifferentGroups}
\ob{D}(\Pi^{\bullet,\ob{cont}})\to  \ob{D}(\Pi^{\bullet}, \bb{Q}_p)
\end{equation}
whose totalization is a symmetric monoidal functor 
\[
\Phi\colon \ob{D}(B\Pi^{\ob{cont}})\to \cat{Rep}^{\sol}_{\bb{Q}_p}(\Pi).
\] 
By construction, the functor $\Phi$ is compatible with the symmetric monoidal  forgetful functors to solid $\bb{Q}_p$-vector spaces, and its formation is natural on group homomorphisms.

We observe that \eqref{eqInclusionNuclearComods1} also produces a symmetric monoidal fully faithful functor 
\[
\ob{D}_{\ob{nuc}}(*/\Pi, \bb{Q}_p)\subset \ob{D}(B \Pi^{\rm cont})).
\]
Finally, by  \cite[Lemma 10.5]{mann2022zl}, the natural map  $\Phi\colon \ob{D}(\GSpec (\Q_p) / \Pi^{\rm cont})\to \cat{Rep}^{\sol}_{\bb{Q}_p}(\Pi)$ is fully faithful when restricted to $\ob{D}_{\ob{nuc}}(*/\Pi, \bb{Q}_p)$ and its essential image consists of those solid $G$-representations over $\bb{Q}_p$ whose underlying solid vector space is nuclear. 
\end{construction}

\begin{remark}\label{Remark-Difference-Actions}
\begin{enumerate}
\item The functor 
\[
\Phi\colon \ob{D}(\GSpec(\bb{Q}_p)/\Pi^{\ob{cont}})\to \cat{Rep}^{\sol}_{\bb{Q}_p}(\Pi)
\]
is not an equivalence of categories, for example,  the regular representation $\Q_{p,\sol}[\Pi]$ is not in the essential image as it is not a $C(\Pi_0,\bb{Q}_p)$-comodule for any open compact subgroup $\Pi_0\subset \Pi$.

\item For a general group $\Pi$  the functor $\Phi$ is not fully faithful.  For example, if $\Pi=\prod_{\bb{N}} \bb{Z}_p$ and $N=  (\prod_{\N} \bb{Z}_p)[\frac{1}{p}]$ we have that $H^1(\Pi, N)=\Hom(\Pi, N)$ is the space of condensed group homomorphisms $\Pi\to N$, which agrees with the space of solid $\bb{Z}_p$-linear maps. However, the image of  the map  $H^1(\Pi^{\ob{cont}}, N)\to H^1(\Pi, N)$  consits of those $\bb{Z}_p$-linear maps $\Pi\to N$ which arise from the image of 
\[
\Hom_{\bb{Z}_p}(\Pi, \bb{Z}_p)\otimes_{\bb{Z}_p} N \to \Hom_{\bb{Z}_p}(\Pi, N).
\]
In particular, the natural inclusion $\Pi\to N$ induced by the identity on  $\prod_{\bb{N}} \bb{Z}_p$ is not in the essential image.

\item We will see in \Cref{LemmaContinuousRepsvsSolid} that under a technical condition that holds in many cases in practice, the functor $\Phi$ is fully faithful.  In particular, this will be  the case if $\Pi=G$ is a $p$-adic Lie group.  
\end{enumerate}
\end{remark}

\begin{lemma}\label{LemmaContinuousRepsvsSolid}
Let $\Pi$ be a (light) locally profinite group and let 
\[
\Phi\colon \ob{D}(B\Pi^{\ob{cont}})\to \cat{Rep}^{\sol}_{\bb{Q}_p}(\Pi)
\] 
be the symmetric monoidal functor of \Cref{ConstructionSymmetricMonoidalFunctorContSolidReps}. 
Suppose in addition that there is a compact open subgroup $\Pi_0\subset \Pi$ such that $\bb{Q}_p$ is a perfect $\bb{Q}_{p,\sol}[\Pi_0]$-module.  Then the following holds: 
\begin{enumerate}

\item  The functor $\Phi$ admits a  $\ob{D}(B\Pi^{\ob{cont}})$-linear right adjoint $\Phi^{R}$. 

\item The composite of $\Phi^{R}$ with the forgetful functor to solid $\bb{Q}_p$-vector spaces is the functor  sending $M\in \cat{Rep}^{\sol}_{\bb{Q}_p}(\Pi)$ to  the solid group cohomology
\[
\Phi^{R}(M)= \Gamma(\Pi_0, M\otimes_{\bb{Q}_p} C(\Pi_0,\bb{Q}_p))
\]
for $\Pi_0\subset \Pi$ a compact open subgroup. 

\item The natural map $\bb{Q}_p\to \Phi^{R}(\bb{Q}_p)$ is an equivalence. In particular, the functor $\Phi$ is fully faithful.

\end{enumerate}
Furthermore, the hypothesis of the lemma holds for $\Pi=G$ a $p$-adic Lie group. 
\end{lemma}
\begin{proof}
The last claim about $p$-adic Lie groups follows from Lazard's resolution, see \cite[Theorem 5.7]{jacinto2021solid}.

In order to give a more conceptual proof  we will use the language of Gestalten of \cite{GestaltenScholze}.  Given a $6$-functor formalism $\ob{D}\colon \ob{Corr}(\n{C})\to \cat{Pr}$ and $X\in \n{C}$ we let $[X]_D\in \cat{Gest}$ denote the Gestalt associated to $X$ by $D$, see \cite[Proposition 9.5]{GestaltenScholze}. We recall that by design, the Gestalt $[X]_{D}$ correspond to the sequence of Stefanich rings 
\begin{equation}\label{eqStefanichRingSixFunctorF}
(\ob{D}(X), \cat{Pr}_{\ob{D},X}, 2\cat{Pr}_{\ob{D},X},\cdots)
\end{equation}
of higher categories of kernels constructed from $D$. 

The functor $\Phi$ involves two different $6$-functor formalisms, namely, that of Gelfand stacks (for $B\Pi^{\ob{cont}}$) which we write $\ob{D}$,  and that of $\bb{Q}_p$-solid pro-\'etale sheaves on condensed anima (for $\cat{Rep}^{\sol}_{\bb{Q}_p}(\Pi)$)  which we write as $\ob{D}_{\ob{proet}}$.

If $X=\GSpec (A)$ is an affinoid Gelfand stack, then $[X]_{\ob{D}}=\ob{Gest}(\ob{D}(A))$ is $1$-affine, namely, this follows form the  construction of the Stefanich ring of the  Gestalt $[X]_{\ob{D}}$ as in  \eqref{eqStefanichRingSixFunctorF}, from  the  categorical K\"unneth formula \cite[Proposition 4.1.14]{solidnotes} and from the fact that $\cat{Pr}_{\ob{D},X}=\ob{Mod}_{\ob{D}(A)}(\cat{Pr})$ thanks to \cite[Lemma 2.2.4]{RCBBCourse2026}.  We deduce that $[\GSpec(\bb{Q}_p)]_{D}=\ob{Gest}(\ob{D}(\bb{Q}_p))$ is a $1$-affine Gestalt, that $[\Pi]_{\ob{D}}=\bigsqcup_{g\in \Pi/\Pi_0} [g\Pi_0]_{\ob{D}}$ is a disjoint union of $1$-affine Gestalten, and that $[B\Pi]_{\ob{D}} = [\GSpec(\bb{Q}_p)]_{\ob{D}}/[\Pi]_{\ob{D}}$.  

We now explain how the $6$-functor formalism of pro-\'etale sheaves on condensed anima is constructed. Let $\n{X}:=\cat{CondAni}$ denote the category of condensed anima, one has a natural $0$-\'etale    $6$-functor formalism on the   topos $\n{X}$ given by sending $S\mapsto \n{X}_{/S}$, see \cite[Construction  3.2.1]{RCBBCourse2026}. Taking base change along the natural symmetric monoidal functor 
\[
\bb{Q}_p[-]\colon \cat{CondAni}\to \ob{D}(\bb{Q}_p)
\] 
given by free solid $\bb{Q}_p$-vector spaces, one obtains the $0$-\'etale $6$-functor formalism of $\bb{Q}_p$-solid pro-\'etale sheaves $\ob{D}_{\proet}(-,\bb{Q}_p)$ (which we  simply denote as $\ob{D}_{\proet}$).

From the previous constructions, we see that  the $1$-affinization of $[*]_{\ob{D}_{\proet}}$ is $\ob{Gest}(\ob{D}(\bb{Q}_p)) = [\GSpec(\bb{Q}_p)]_{\ob{D}}$. Thus, we have a natural morphism of Gestalten
\[
f\colon [*]_{\ob{D}_{\proet}}\to [\GSpec(\bb{Q}_p)]_{\ob{D}}. 
\]
To simplify notation we shall write $S:=[*]_{\ob{D}_{\proet}}$ and   let $\Pi^{\ob{cont}}_S$ denote the pullback of $[\Pi^{\ob{cont}}]_{\ob{D}}\to [\GSpec(\bb{Q}_p)]_{\ob{D}}$ along $f$. Denote $\Pi^{\sol}:=[\Pi]_{\ob{D}_{\proet}}$, the symmetric monoidal morphism of cosimplicial categories \eqref{eqSymmetricMonoidalDifferentGroups} gives rise to a morphism of group-Gestalten over $S$
\[
\Pi^{\sol}\to \Pi^{\ob{cont}}_{S}
\]
and hence a morphism of classifying stacks over  $S$
\[
g\colon B\Pi^{\sol}\to B(\Pi^{\ob{cont}}_S).
\]

By construction, the functor $g^*_1\colon \s{O}(B(\Pi^{\ob{cont}}_S))_1\to \s{O}(B\Pi^{\sol})_1$ on $1$-sheaves is precisely the functor $\Phi$ of \Cref{ConstructionSymmetricMonoidalFunctorContSolidReps}. 
Note that for $\Pi_0\subset \Pi$ an open compact subgroup one has a pullback square
\begin{equation}\label{eqPullbackSquareSolidVsContinuousStacks}
\begin{tikzcd}
B\Pi_0^{\sol} \ar[r,"g_0"] \ar[d] &  B(\Pi^{\ob{cont}}_S) \ar[d]\\ 
B\Pi^{\sol}  \ar[r,"g"] &  B(\Pi^{\ob{cont}}_S).
\end{tikzcd}
\end{equation}
Indeed, this follows from the fact that 
\[
\Pi^{\sol}/\Pi_{0}^{\sol}=  \Pi^{\ob{cont}}_S/\Pi^{\ob{cont}}_{0,S} = \Pi/\Pi_0\times S,
\]
where $\Pi/\Pi_0$ is a discrete set.

\begin{claim*}
The map $g$ is $0$-prim. 
\end{claim*}

\begin{proof}[Proof of the Claim] By \cite[Proposition 7.9]{GestaltenScholze} we can show the claim locally on the target, and by \eqref{eqPullbackSquareSolidVsContinuousStacks} we can then assume without loss of generality that $\Pi=\Pi_0$ is compact and that  $\bb{Q}_p$ is a perfect $\bb{Q}_{p,\sol}[\Pi_0]$-module.  In such case, the diagonal map of the Gelfand stack $B\Pi$ is $0$-proper (even $0$-affine) as it is given by 
 \[
 \Delta\colon B\Pi\cong \Pi\backslash \Pi/\Pi\to \Pi\backslash \GSpec(\bb{Q}_p) /\Pi
 \]
 and $\Pi$ is profinite. This implies that the target of $g$ has $0$-proper (even $0$-affine) diagonal.  Thus,  the $0$-primness of $g$ will follow if  we show that $B\Pi^{\sol}$ is $0$-prim by \Cref{LemmaDiagonalPrimGestalten} down below.  For this, we note that, by construction, any map of condensed anima $Y\to X$ is $0$-\'etale for the pro-\'etale $6$ functor-formalism. In particular, the map $h\colon *\to B\Pi$ is suave and $h_{\sharp} 1 = \bb{Q}_{p,\sol}[\Pi]\in \ob{D}_{\proet}(B\Pi,\bb{Q}_p)$ is a prim object  by \cite[Lemma 4.4.9 (ii)]{heyer20246functorformalismssmoothrepresentations}.
Since prim objects of $B\Pi$ form a stable category stable under retracts, cf. \cite[Corollary 4.4.13]{heyer20246functorformalismssmoothrepresentations}, the fact that $\bb{Q}_p$ is a perfect $\bb{Q}_{p,\sol}[\Pi]$-module implies that it is a prim object in $\ob{D}_{\proet}(B\Pi,\bb{Q}_p)$. Going back to Gestalten, we deduce that $B\Pi^{\sol}\to S$ is $0$-prim, which finishes the proof of the claim. 
\end{proof}

We finish the proof of the lemma using the claim. Since $g$ is $0$-prim, $g_{1}^*=\Phi$ will admit a linear right adjoint $g_{1,*}$ proving  part (1). This right adjoint is compatible with proper base change, so by considering the pullback square 
\[
\begin{tikzcd}
\Pi_{0,S}^{\ob{cont}}/\Pi_0^{\sol} \ar[r,"g'"]  \ar[d,"q'"]& S \ar[d,"q"] \\
 B\Pi^{\sol} \ar[r,"g"] & B\Pi_S^{\ob{cont}}, 
\end{tikzcd}
\]
we have an equivalence of functors $q^*_1g_{1,*}\cong g_{1,*}^{\prime}q^{\prime,*}_1$ from $\cat{Rep}^{\sol}_{\bb{Q}_p}(\Pi)=\s{O}(B\Pi^{\sol})_1\to \s{O}(S)_1=\ob{D}(\bb{Q}_p)$. The functor $q^{\prime,*}_1$ is the functor sending $M\mapsto M\otimes_{\bb{Q}_p} C(\Pi_0,\bb{Q}_p)$ endowed with the diagonal action as a solid $\Pi_0$-representation, and the functor $g^{\prime}_{1,*}$ is just solid $\Pi_0$-cohomology. This proves part (2) of the lemma. Finally,  for part (3) we use (2), namely that 
\[
\Phi^{R}(\bb{Q}_p)= \Gamma(\Pi_0, C(\Pi_0, \bb{Q}_p))= \bb{Q}_p
\]
where the second equivalence follows from the fact that  $C(\Pi_0, \bb{Q}_p)=\Hom_{\bb{Q}_p}(\bb{Q}_{p,\sol}[\Pi_0], \bb{Q}_p)$ is an induced representation. The fully faithfulness of $\Phi$ then follows formally from the projection formula. 
\end{proof}

We used the following lemma:

\begin{lemma}\label{LemmaDiagonalPrimGestalten}
Let $Y\to X$ be a morphism of Gestalten over $S$. Let $n\geq 0$ be an integer and suppose that the diagonal of $X\to S$  and that $Y\to S$ are $n$-prim. Then $Y\to X$ is $n$-prim. The same holds for $n$-properness, $n$-\'etalenss and $n$-suaveness.  
\end{lemma}
\begin{proof}
We can write $f\colon Y\to X$ as the composite of 
\[
Y\xrightarrow{\Gamma_f} Y\times_S X \to X,
\]
where the first map is the graph, and the second the projection. The first map is the base change of the diagonal $X\to X\times_S X$ along $Y\times_S X\to X\times_S X$, and the second a base change of the map $Y\to S$ along $X\to S$. The lemma follows since $n$-prim maps are stable under base change and composition, see \cite[Proposition 6.10 (iii) and (iv)]{GestaltenScholze}.
\end{proof}

Finally, we record the following projection formula for solid cohomology.

\begin{lemma}
  \label{sec:local-locally-analyt-nuclear-projection-formula-for-uniform-groups}
Supose that $\Pi$ is a profinite group and suppose that $\bb{Z}_p$ is a perfect $\bb{Z}_{p,\sol}[\Pi]$-module.  Then the cohomology functor
  \[
    \mathrm{Rep}^\sol_{\Z_p}(\Pi)\to \ob{D}(\Z_p)
  \]
  is $\ob{D}(\Z_p)$-linear.  In particular, this holds if    $\Pi=G$ is  a uniform pro-$p$ group.
\end{lemma}
\begin{proof}
 Since $\bb{Z}_p$ is a perfect $\bb{Z}_{p,\sol}[\Pi]$-module, the natural morphism of functors  $\mathrm{Rep}^\sol_{\Z_p}(\Pi)\to \ob{D}(\Z_p)$
  \[
  \Hom_{\bb{Z}_{p,\sol}[\Pi]}(
  \bb{Z}_p,\bb{Z}_{p,\sol}[\Pi])\otimes_{\bb{Z}_{p,\sol}[\Pi]} (-)\to   \Hom_{\Z_{p,\sol}[G]}(\Z_p,-)= \Gamma(\Pi,-)
    \]
    is an equivalence.
    The left functor is clearly $\ob{D}(\bb{Z}_{p})$-linear, proving the first claim of the lemma.
    
    Finally,  if $G$ is uniform-pro-$p$, the trivial representation $\Z_p$ has a finite resolution by finite free $\Z_{p,\sol}[G]$-modules, cf.\ \cite[Theorem 5.7]{jacinto2021solid}, proving the second assertion of the lemma. 
    \end{proof}

\subsubsection{Locally analytic representations} \label{sss-loc-an-Reps}
Next, we introduce the relevant notation for the theory of locally analytic representations.
We follow the setup of \cite{jacinto2021solid,jacinto2023solid}. Let $G$ be a compact $p$-adic Lie group and $G_0\subset G$ an open normal subgroup given by a uniform pro-$p$-group.
The group $G_0$ admits a basis $g_1,\ldots, g_d$ which yields an isomorphism of $p$-adic Lie manifolds $G_0\cong \bb{Z}_p^d$, and a $p$-valuation $w\colon G_0\backslash\{1\}\to (1/(p-1),\infty)$ (cf.\ \cite[Section 4]{zbMATH01965478}).
After shrinking $G_0$ if necessary, we can assume that the group law of $G_0$  extends to a group law on the affinoid unit polydisc $\bb{D}^{d}_{\bb{Q}_p}=\GSpec(\bb{Q}_p\langle T_1,\ldots, T_d\rangle)$ with coefficients in $\Z_p$, and that the action of $G$ on $G_0$ by conjugation is also given, termwise on $G$, by power series  with coefficients in $\Z_p$.
We let $\bb{G}_{0}(1)$ denote the rigid analytic group structure on  $\bb{D}^{d}_{\bb{Q}_p}$ defined by $G_0$, and for $r\in (0,1]\cap \bb{Q}$ a rational number we let $\bb{G}_{0}(r)\subset \bb{G}_{0}(1)$ denote the affinoid polydisc of radius $r$.
We let $G_r\subset G_0$ denote the compact open subgroup of $\bb{Q}_p$-points of $\bb{G}_{0}(r)$.
Since the group law of $\bb{G}_{0}(1)$ has integral coefficients it restricts to a group law on $\bb{G}_{0}(r)$ and, similarly, since the action of $G$ by conjugation on $G_{0}$ is given by power series with integral coefficients, one can construct a unique rigid analytic group structure on $\bb{G}_{r}:=\bigsqcup_{g\in G/ G_r} g\bb{G}_{0}(r)$ compatible with the group structure of $G$.
In the notation of  \cite[Definition 2.1.2 and 2.1.4]{jacinto2023solid}, if  $r=|p^h|$ for $h\geq 0$, one has  $\bb{G}^h= \bb{G}_{r}$.
The group $G^{\la}$ is defined as the intersection $G^{\la}=\bigcap_{r\to 0} \bb{G}_{r}$, or equivalently, it is the affinoid group object in Gelfand stacks whose algebra of functions is given by 
\[
\s{O}(G^{\la})= C^{\la}_G:=\varinjlim_{r\to 0} \s{O}(\bb{G}_{r}). 
\]
We note that  $C_{G}^{r}:=\s{O}(\bb{G}_{r})$ is denoted as $C^{h}(G,\bb{Q}_p)$ in  \cite[Definition 2.1.2]{jacinto2021solid} for $r=|p^h|$.
Finally, we observe that by construction one has
\[
\s{O}(\bb{G}_r)= \ob{Ind}_{G_r}^G(\s{O}(\bb{G}_0(r))). 
\] 

To study locally analytic representations one must consider the following groups built from the $\bb{G}_{r}$'s. 

\begin{definition}\label{DefinitionOpenOverconvergentGroups}
  For $r\in (0,1]$ we let $\bb{G}_{<r}\subset \bb{G}_{r}$ denote the open subgroup given by the union $\bb{G}_{<r}=\bigcup_{r'<r} \bb{G}_{r'}$.
  Similarly, for $r\in (0,1)$, we let $\bb{G}_{\leq r}\subset \bb{G}_1$ denote the subgroup given by the intersection $\bb{G}_{\leq r}:= \bigcap_{r'>r} \bb{G}_{r'}$.
\end{definition}

The groups $\bb{G}_{<r}$ consist of  finite disjoint unions of open polydiscs of radius $r$, while the groups $\bb{G}_{\leq r}$ consist of finite disjoint unions of overconvergent closed discs of radius $r$. The group $\bb{G}_{<r}$ is denoted by $\bb{G}^{(h^+)}$ in \cite{jacinto2021solid} for $r=|p^h|$. The groups $\bb{G}_{\leq r}$ are not considered neither in \cite{jacinto2021solid} nor in \cite{jacinto2023solid}, but they play an important role in \cite{mikami2025finiteness}.  We shall write
$$C^{<r}_G:= \s{O}(\bb{G}_{<r}), ~ C^{\leq r}_G:=\s{O}(\bb{G}_{\leq r}).$$
 By \cite[Remark 2.1.5]{jacinto2023solid} the ind-systems $\{C^{r}_{G}\}_{r\to 0}$, $\{C^{<r}_{G}\}_{r\to 0}$ and $\{C^{\leq r}_{G}\}_{r\to 0}$ are ind-equivalent.
In particular, we have that
\[
C^{\la}_G=\varinjlim_{r\to 0} C^{r}_G =\varinjlim_{r\to 0} C^{<r}_G =\varinjlim_{r\to 0} C^{\leq r}_G.  
\]
Following \cite{jacinto2023solid}, we let $\n{D}^{<r}_G :=\s{O}(\bb{G}_{<r})^{\vee}$ denote the $(<r)$-distribution algebra (what is denoted $\n{D}^h(G, \bb{Q}_p)$ in \cite[Definition 2.1.2]{jacinto2023solid}), obtained as the underived dual of the global functions of $\bb{G}_{<r}$.
We let $\n{D}^{\la}_G:=\varprojlim_{r\to 0} \n{D}^{<r}_G$ be the locally analytic distribution algebra of $G$.

The following notion is studied in depth in \cite{jacinto2021solid,jacinto2023solid}.

\begin{definition}
  \label{sec:semil-solid-repr-definition-locally-analytic-representation}
  An object $M\in \cat{Rep}_{\Q_p}^\sol(G)$ is called a \emph{locally analytic representation} if the natural map in $\cat{Rep}^\sol_{\Q_p}(G)$
  \[
     M^\la:=\Gamma(G,C^\la_G\otimes_{\Q_p}M)\to M
  \]
  is an isomorphism, where $\Gamma(G,-)\colon \cat{Rep}^\sol_{\Q_p}(G)\to D(\Q_p)$ is the functor of (derived, solid) invariants.\footnote{The $G$-action on $C^\la_G\otimes_{\Q_p}M$ is diagonal with the  left regular action on $C^\la_G$. The right regular action on $C^\la_G$ equips $M^\la$ with a $G$-action.}

  We let $\cat{Rep}^\la_{\Q_p}(G)\subseteq \cat{Rep}^\sol_{\Q_p}(G)$ be the full subcategory of locally analytic representations.
\end{definition}

We record some properties of the category of locally analytic representations.

\begin{lemma} \label{Lemma-Prim-padic}
\begin{enumerate}
\item Let $G$ be a (non necessarily compact)  $p$-adic Lie group.  Then the morphism $s\colon BG^{\rm cont}\to BG^\la$ is $0$-proper, and $s_\ast(\s{O})=\s{O}$. In particular,   $*$-pullback along  $s$ is fully faithful.

\item Consider the composite of fully faithful symmetric monoidal functors 
\[
\Psi\colon \ob{D}(BG^{\la})\xrightarrow{s^*} \ob{D}(BG^{\ob{cont}})\xrightarrow{\Phi} \cat{Rep}^{\sol}_{\bb{Q}_p}(G)
\]
where $\Phi$ is as in \Cref{ConstructionSymmetricMonoidalFunctorContSolidReps}. Then $\Psi$ admits a $\ob{D}(BG^{\la})$-linear right adjoint.  
\item Suppose that $G$ is compact.  Then the  composite of $\Psi^{R}$ with the $*$-pullback functor along $\GSpec(\bb{Q}_p)\to BG^{\la}$ is the functor $(-)^{\la}$ of \Cref{sec:semil-solid-repr-definition-locally-analytic-representation}. In particular, the essential image of $\Psi$ is the full subcategory $\cat{Rep}^{\la}_{\bb{Q}_p}(G)\subset \cat{Rep}^{\sol}_{\bb{Q}_p}(G)$ of solid locally analytic representations. 
\end{enumerate}
\end{lemma}

\begin{proof}
\begin{enumerate}
\item By descent on the target we can assume without loss of generality that $G$ is compact.  To see that the map  $s\colon BG^{\rm cont}\to BG^{\la}$ is $0$-proper, it suffices to prove it after taking the pullback along $\GSpec(\Q_p)\to BG^{\la}$ in which case one has $G^{\la}/G^{\rm cont}\to \GSpec(\Q_p)$ which is $0$-proper as $G^{\la}$ is affine (hence $0$-proper) and $BG^{\rm cont}$ is proper by \Cref{LemmaFirstPropertiesGcont} (3).    The computation $s_*\Q_p=\Q_p$ follows form the fact that $s_*$ computes locally analytic vectors, see \cite[Proposition 6.1.5]{jacinto2023solid}.

\item By \Cref{LemmaContinuousRepsvsSolid} (1) the functor $\Phi^{R}$ is $\ob{D}(BG^{\ob{cont}})$-linear, while by part (1) of this lemma the functor $s^*$ is $\ob{D}(BG^{\la})$-linear. One deduces that $\Psi^{R}$ is $\ob{D}(BG^{\la})$-linear.   

\item To see that the functor $\Psi^{R}$ is $(-)^{\la}$, one can argue exactly as in the proof of \Cref{LemmaContinuousRepsvsSolid} (2). \qedhere

\end{enumerate}

\end{proof}

We can define categories of   $(<r)$ and $(\leq r)$-analytic representations of $G$ using the groups $\bb{G}_{<r}$ and $\bb{G}_{\leq r}$.

 \begin{definition}\label{DefinitionhAnRep}
Let $r\in (0,1)$, we define the category of $(< r)$-analytic representations of $G$ (resp. the category of $(\leq r)$-analytic representations) to be $\cat{Rep}^{<r}_{\bb{Q}_p}(G):= \ob{D}(\GSpec( \bb{Q}_p)/ \bb{G}_{<r})$ (resp. $\cat{Rep}^{\leq r}_{\bb{Q}_p}(G):= \ob{D}(\GSpec (\bb{Q}_p)/ \bb{G}_{\leq r})$). If $R$ is an algebra in  $\cat{Rep}^{<r}_{\bb{Q}_p}(G)$, we let $\cat{Rep}^{<r}_{R}(G)$ denote the category of left $R$-modules on $\cat{Rep}^{<r}_{\bb{Q}_p}(G)$ (and similarly for $\leq r$). 
\end{definition}

The following lemma establishes some basic facts  of the categories of locally analytic, $(<r)$-analytic and $(\leq r)$-analytic representations. 

\begin{lemma}\label{LemmaAnalyticRepresentationProperty}
Let $r\in (0,1)$. The following hold.

\begin{enumerate}

\item The map $\GSpec(\Q_p)\to BG^{\la}$ is $0$-proper and descendable. In particular,  $BG^{\la}$ is $0$-proper over $\GSpec(\Q_p)$.  Moreover, $BG^{\la}$ is suave.

\item The map $\GSpec(\Q_p)\to B\bb{G}_{\leq r}$ is $0$-proper and descendable. In particular,  $B\bb{G}_{\leq r}$ is $0$-proper over $\GSpec(\Q_p)$. Moreover, $B\bb{G}_{\leq r}$ is suave.

\item  The map $\GSpec(\Q_p)\to B\bb{G}_{<r}$ is cohomologically smooth, in particular  $B\bb{G}_{<r}$ is $!$-able over $\GSpec(\Q_p)$ and cohomologically smooth. Moreover, $B\bb{G}_{<r}$ is prim. 

\item The map $f\colon BG^{\la}\to B\bb{G}_{\leq r}$ is $0$-proper  and the natural map $1\to f_*1$ is an equivalence. In particular the pullback functor $f^*\colon \cat{Rep}^{\leq r}_{\sol}(G)\to  \cat{Rep}^{\la}_{\sol}(G)$ is fully faithful.

\item The map $f\colon BG^{\la}\to B\bb{G}_{<r}$ is suave and the natural morphism $f_{\sharp} 1\to 1$ is an equivalence. In particular the pullback functor $f^*\colon \cat{Rep}^{< r}_{\sol}(G)\to  \cat{Rep}^{\la}_{\sol}(G)$ is fully faithful.

\item The $*$-pullback along $\AnSpec (\bb{Q}_p)\to B\bb{G}_{<r}$ produces a natural equivalence of categories 
\[
\ob{D}(B\bb{G}_{<r})\xrightarrow{\sim} \ob{LMod}_{\n{D}^{<r}_G}(\ob{D}(\bb{Q}_p)).
\]

\end{enumerate}
\end{lemma}
\begin{proof}
\begin{enumerate}

\item Let $G^{\dagger}\subset G^{\la}$ denote the overconvergent neighbourhood at $1$. We have a cartesian square of Gelfand stacks
\[
\begin{tikzcd}
BG^{\dagger} \ar[r] \ar[d] & BG^{\la} \ar[d] \\ 
\GSpec(\Q_p)  \ar[r] & BG^{\sm}.
\end{tikzcd}
\]
Since $G$ is a compact $p$-adic Lie group, the map $h\colon \GSpec(\Q_p)\to BG^{\sm}$ is proper  and $h_* \bb{Q}_p$ is the regular $\bb{Q}_p$-linear smooth representation of $G$.
A Haar measure (whose existence is guaranteed from the fact that we are working in characteristic zero) yields a retract of $\bb{Q}_p\to h_*\bb{Q}_p$ as a smooth $G$-representation, proving that $h$ is descendable and therefore that $BG^{\sm}$ is proper by \cite[Corollary 4.7.5 (iii)]{heyer20246functorformalismssmoothrepresentations}.
Thus, it suffices to show that $g\colon \GSpec(\Q_p)\to BG^{\dagger}$ is $0$-proper and descendable; $0$-properness is clear as $G^{\dagger}$ is an affinoid Gelfand stack.
For descendability, we can  use the relative de Rham complex of $g$ endowed with its Hodge filtration.
Namely, as $G^{\dagger}$ has trivial de Rham cohomology, the relative de Rham complex of $g$ produces an increasing $\bb{N}$-filtration of the trivial representation $\bb{Q}_p$ with $i$-th graded piece given by $g_* \bigwedge^i \bb{L}_{g}[-i]= g_* \bigwedge^i \f{g}^{\vee}[-i]$ with $\f{g}=\ob{Lie} G$.
This filtration is supported in degrees $[0,d]$ with $d=\ob{dim} G$, which makes the map $g$ descendable with index of descendability bounded by $\leq d+2$.
 Finally, for suaveness of $BG^{\la}$, we note that $BG^{\sm}\to \GSpec(\Q_p)$ is suave as:
\begin{itemize}

\item  The map $h\colon \GSpec(\Q_p)\to BG^{\sm}$ is prim so that $h_*\bb{Q}_p\in \ob{D}(BG^{\sm})$ is suave. 
  
\item  The trivial representation $\bb{Q}_p$ is a  direct summand of $h_*\bb{Q}_p$ and therefore it is suave. 

\end{itemize}

Hence, it suffices to see that $BG^{\dagger}\to\GSpec(\Q_p)$ is suave, but the unit of $BG^{\dagger}$ admits a finite filtration whose graded pieces are $g_* \bigwedge^i \f{g}^{\vee}[-i]$ with $g\colon \GSpec(\Q_p)\to BG^{\dagger}$.
As $g$ is prim,  $g_* \bigwedge^i \f{g}^{\vee}\in \ob{D}(BG^{\dagger})$ is suave for all $i$ and then  so is $1\in \ob{D}(BG^{\dagger})$.
This finishes the proof of (1).

\item   This follows from a similar argument as for (1).
  It is clear that the map $f\colon \GSpec(\Q_p)\to B\bb{G}_{\leq r}$ is $0$-proper as $\bb{G}_{\leq r}$ is an affinoid Gelfand stack. To prove that  $f$ is descendable, by passing to the connected component of $1$, and since all maps $\GSpec(\Q_p)\to BH$ are proper and descendable for $H$ a finite group, we can assume without loss of generality that $\bb{G}_{\leq r}$ is an overconvergent polydisc of radius $r$. Now, since $\bb{G}_{\leq r}\to \bb{G}_{\leq r,\dR}$ is surjective, we have a pullback square
\[
\begin{tikzcd}
{BG^{\dagger} }\ar[r] \ar[d]& {B\bb{G}_{\leq r}} \ar[d] \\ 
{\GSpec(\Q_p) }\ar[r] & {B\bb{G}_{\leq r}^{\dR}}
\end{tikzcd}
\] 
As the map $\GSpec(\Q_p)\to BG^{\dagger}$ is $0$-proper and descendable by the proof of (1), it suffices to show that $g\colon \GSpec(\Q_p) \to B\bb{G}_{\leq r}^{\dR}$ is proper and descendable; properness follows form \cite[Theorem 1.3.2 (1)]{drff_stacks}, descendability follows from the fact that $g_*1 =1$ as  $\bb{G}_{\leq r}^{\dR}$ has contractible cohomology.  Finally, for suaveness of $B\bb{G}_{\leq r}$, as $BG^{\dagger}\to \AnSpec (\bb{Q}_p)$ is suave it suffices to show that $B\bb{G}^{\dR}_{\leq r}\to \AnSpec (\bb{Q}_p)$ is suave, but the map $g\colon \GSpec(\Q_p) \to B\bb{G}_{\leq r}^{\dR}$ is proper and $g_*1=1$, this makes the unit of $B\bb{G}_{\leq r}^{\dR}$ suave, finishing the proof of (2).

\item This is proven in a similar way as (2), we will just mention the main changes.
  Suaveness of $\GSpec(\Q_p)\to B\bb{G}_{<r}$ follows from suaveness of the smooth partially proper rigid space $\bb{G}_{<r}$.
  By \cite[Lemma 4.4.10 (i)]{heyer20246functorformalismssmoothrepresentations} we deduce that $B\bb{G}_{<r}$ is suave over the point.
  Next we show primness.
  By passing to the connected component of the unit, we can assume that $\bb{G}_{<r}$ is an open polydisc of radius $r$. We have a cartesian diagram
\[
\begin{tikzcd}
{BG^{\dagger} }\ar[r] \ar[d]& {B\bb{G}_{< r}} \ar[d] \\ 
{\GSpec(\Q_p) }\ar[r] & {B\bb{G}_{< r}^{\dR}}
\end{tikzcd}
\] 
The map  $h\colon \GSpec(\Q_p) \to B\bb{G}_{< r}^{\dR} $ is suave thanks to \cite[Theorem 1.3.2 (2)]{drff_stacks} as it is represented in the de Rham stack of a  partially proper smooth rigid space. Furthermore, as  the de Rham homology of the open disc is contractible, we have that $h_{\sharp} 1=1$.   We deduce that $h'\colon BG^{\dagger}\to B\bb{G}_{<r}$ is suave and that $h'_{\sharp} 1=1$.
By \cite[Lemma 4.5.16]{heyer20246functorformalismssmoothrepresentations} we deduce that the unit of $B\bb{G}_{<r}$ is prim, finishing the proof of (3).

\item  The fact that $f\colon BG^{\la}\to B\bb{G}_{\leq r}$ is $0$-proper follows from the fact that  both stacks are $0$-proper. To show that $f_* 1=1$, it suffices to prove that the global sections of $\bb{G}_{\leq r}/G^{\la}$ are given by $\bb{Q}_p$. By taking the connected component of $1$ in $\bb{G}_{\leq r}$, and passing to an open compact subgroup of $G$, we can assume that $\bb{G}_{\leq r}$ is a closed overconvergent polydisc. In that case, we have that 
\[
\bb{G}_{\leq r}/G^{\la} = (\bb{G}_{\leq r}/G^{\dagger})/G^{\sm} = \bb{G}_{\leq r}^{\dR}/G^{\sm}.
\]
The acyclicity then follows from the fact that  $\bb{G}_{\leq r}^{\dR}$ has trivial de Rham cohomology, and that taking smooth $G$-invariants is a $t$-exact functor. 

\item Suaveness of $f$  follows formally from suaveness of $BG^{\la}$ and the fact that the diagonal of $B\bb{G}_{<r}$ is suave as  $\bb{G}_{<r}$ is a smooth partially proper rigid space.
  The rest of the proof is  similar to that  of (4).
  By passing to a compact open subgroup we can assume that $\bb{G}_{<r}$ is connected and that it is an open polydisc of radius $r$.
  The statement then follows from the fact that the de Rham homology of  $\bb{G}_{<r}$ is contractible, and again the $t$-exactness of invariants for smooth representations of $G$.

\item  Consider the map $f\colon \GSpec(\Q_p)\to B\bb{G}_{<r}$.
  It is suave and the $*$-pullback is conservative and preserves both limits and colimits.
  By the monadicity theorem, we have that 
\[
\ob{D}(B\bb{G}_{<r})=\ob{LMod}_T(\ob{D}(\bb{Q}_p))
\]
where $T$ is the monad $T=f^*f_{\sharp}$. Now, both functors $f^*$ and $f_{\sharp}$ are $\ob{D}(\bb{Q}_p)$-linear, which implies that $T$ arises from the monad of an associative algebra $\n{D}\in \ob{Alg}(\ob{D}(\bb{Q}_p))$. This shows that 
\[
\ob{D}(B\bb{G}_{<r})=\ob{LMod}_{\n{D}}(\ob{D}(\bb{Q}_p)).
\]
Now, one computes the  algebra $\n{D}$. By proper base change applied to the square 
\[
\begin{tikzcd}
\bb{G}_{<r}  \ar[r] \ar[d] & \GSpec(\Q_p) \ar[d] \\ 
\GSpec(\Q_p)  \ar[r]& B\bb{G}_{<r}
\end{tikzcd}
\]
we see that the underlying $\bb{Q}_p$-vector space of $\mathcal{D}$ is the homology of  $\bb{G}_{<r}$ which is precisely the underived dual of $\s{O}(\bb{G}_{<r})$, that is, $\n{D}^{<r}_G$.
It is left to identify the algebra structure.
From the description of the category $\ob{D}(B\bb{G}_{<r})$  as the totalization of the \v{C}ech nerve of $\GSpec(\Q_p)\to B\bb{G}_{<r}$, the multiplication map of $\n{D}$  arises from the homology of the multiplication map 
\[
m\colon \bb{G}_{<r}\times \bb{G}_{<r} \to \bb{G}_{<r}
\]
after applying the K\"unneth formula, namely, from the map 
\[
m_{\sharp}\colon \Gamma_{\sharp}(\bb{G}_{<r}\times \bb{G}_{<r})=\Gamma_{\sharp}(\bb{G}_{<r})\otimes \Gamma_{\sharp}(\bb{G}_{<r}) \to  \Gamma_{\sharp}(\bb{G}_{<r})
\] 
(where for a suave map $X\to \GSpec(\Q_p)$  we let $\Gamma_{\sharp} (X)$ denote its homology). This morphism is the underived dual of the morphism 
\[
m^*\colon \s{O}(\bb{G}_{<r})\to \s{O}(\bb{G}_{<r}\times \bb{G}_{<r}).
\]
Thus, $m_\sharp$ is the multiplication map of $\n{D}^{<r}_G$. \qedhere
\end{enumerate}
\end{proof}

\begin{remark}\label{RemarkRAnalyticFailure}
  It is natural to ask whether the affinoid group $\bb{G}_r$ also gives rise to a reasonable category of $r$-analytic representations.
  Naturally, we can define it as  $\ob{D}(B\bb{G}_r)$, and one might  ask whether  the $*$-pullback map along $g\colon BG^{\la}\to B\bb{G}_r$ is fully faithful.
  We claim that this is never the case when $G$ is infinite.
  First, note that $g$ is $0$-proper. Indeed, the quotient $\bb{G}_r/G^{\la}$ is $0$-proper over $\AnSpec (\bb{Q}_p)$ as  $\bb{G}_r$ is affinoid with trivial analytic ring structure, and $BG^{\la}$ is $0$-proper by \Cref{LemmaAnalyticRepresentationProperty} (1).
  Therefore, the pullback of  $g$ is fully faithful if and only if $1\to g_*1$ is an equivalence.
  Let us now compute $g_*$. Let $\bb{G}_r^0\subset \bb{G}_{r}$ denote the connected component of $1$, by construction we have that $\bb{G}_r/\bb{G}_r^0=G/G_r$.
  Thus, we have an equivalence of stacks
\[
BG^{\la}_r= BG^{\la}\times_{B\bb{G}_r} B\bb{G}_r^{0},
\]
and by shrinking $G$, we can assume without loss of generality that $\bb{G}_r$ is an affinoid polydisc. 

Now, the underlying $\bb{Q}_p$-vector space of $g_*1$ is precisely the cohomology of $\bb{G}_r/G^{\la}$.
Let $G^{\dagger}\subset G^{\la}$ denote the overconvergent neighbourhood at $1$.
We have a short exact sequence of groups
\[
1\to G^{\dagger}\to G^{\la}\to G^{\sm}\to 1,
\]
so that 
\[
\bb{G}_r/G^{\la} = (\bb{G}_r/G^{\dagger})/G^{\sm}
\]
and $g_*1$ computes the $G^{\sm}$-invariants of the cohomology of $\bb{G}_r/\bb{G}^{\dagger}$. But $\Gamma(\bb{G}_r/\bb{G}^{\dagger} ,\s{O})$ is the de Rham cohomology of the affinoid disc $\bb{G}_r$, on which the action of $G^{\sm}$ is trivial.
Therefore $g_*1=\Gamma_{\dR}(\bb{G}_r)$.
Now, as $\bb{G}_r$ is an affinoid polydisc, its de Rham cohomology is not concentrated in degree $0$, and the higher cohomologies consist of non quasi-separated solid $\bb{Q}_p$-vector spaces. 
\end{remark} 

\subsubsection{$\ob{Bin}_{r,\varpi}$-analytic representations}\label{sss:BinomReps} For proving \Cref{TheoDecompletion} we will need to analyse some integral variant of the theory of analytic representations.
This will require working with some analytic groups defined by $\bb{Q}_p$-Banach Hopf-algebras admitting natural $\bb{Z}_p$-formal models.
In view of  \Cref{RemarkRAnalyticFailure}, we must modify  the groups $\bb{G}_r$ in order to have a reasonable rational theory of $r$-analytic representations. To define the modified groups  one can find inspiration from the description of locally  analytic functions by Amice. In the following, given a tuple $\alpha\in \bb{N}$ with $\alpha=(\alpha_1,\ldots, \alpha_d)$, we denote $\binom{\mathbf{x}}{\alpha}=\binom{x_1}{\alpha_1}\cdots \binom{x_d}{\alpha_d}$ and $|\alpha|=\sum_i \alpha_i$. 

\begin{theorem}[{\cite[Corollary I.3.9]{zbMATH03225421}}]
  \label{sec:obbin_r-varpi-analyt-amice-theorem}
Let $G$ be a uniform pro-$p$-group with coordinates $\psi\colon G\xrightarrow{\sim} \bb{Z}_p^d$ and $M$ a $\bb{Q}_p$-Banach space with norm $\|-\|$. A function $f\colon G\to M$ is locally analytic if and only if its binomial expansion 
\[
f(\underline{x})=\sum_{\alpha\in \bb{N}^d} c_{\alpha}  \binom{\mathbf{x}}{\alpha} 
\]
satisfies the following condition: there is $r>1$ such that $\|c_{\alpha}\| r^{|\alpha|}\to 0$ as $|\alpha|\to \infty$, where $|\alpha|=\sum_i \alpha_i$. 
\end{theorem}

The following discussion is taken from  \cite[Section 4]{zbMATH01965478}, see also \cite[Section 3.4]{jacinto2023solid}.  Let $G$ be a uniform pro-$p$-group with fixed coordinates $\psi\colon G\xrightarrow{\sim} \bb{Z}_p^d$ and $p$-valuation $w\colon G\backslash \{1\}\to (1/(p-1),\infty)$.

 Let $g_1,\ldots, g_d\in G$ be the basis defined by $\psi$ and set $w_i:= w(g_i)$. For a tuple $\alpha\in \bb{N}^d$ we define $\tau \alpha= \sum_i w_i \alpha_i$.  By modifying the $p$-valuation, we can assume that the $w_i$ are rational numbers (cf. \cite[III.3.1.11]{zbMATH03225421}). Consider the Iwasawa algebra $\bb{Z}_{p,\sol}[G]$, and let $b_i:=g_i-1\in \bb{Z}_{p,\sol}[G]$. The Mahler basis  of $C(G,\bb{Z}_p)$ produces a power series description
\[
\bb{Z}_{p,\sol}[G]=\prod_{\alpha\in \bb{N}} \bb{Z}_p \mathbf{b}^{\alpha}
\]
where $\mathbf{b}^{\alpha}=b_1^{\alpha_1}\cdots b_d^{\alpha_d}$. The algebra $\bb{Z}_{p,\sol}[G]$ is endowed with a filtration induced by $w$, and independent of the coordinates $\psi$, such that 
\[
\mathrm{Fil}^s(\bb{Z}_{p,\sol}[G]) =\{ \sum_{\alpha\in \bb{N}^d} a_{\alpha} \mathbf{b}^{\alpha} \; : \; \sup_{\alpha}( |a_{\alpha}|p^{-\tau\alpha}) \leq p^{-s}\}. 
\]
We endow $\bb{Z}_p$ with the $p$-adic filtration, and let $\bb{F}_p[\epsilon]=\gr^{\bullet} \bb{Z}_p$ its graded algebra, where $\epsilon=\overline{p}$ is the residue class of $p$ with weight $1$. By  \cite[Theorem III.2.3.3]{zbMATH03225421}, the inclusion $G\to \bb{Z}_{p,\sol}[G]$ produces a natural isomorphism of  graded rings over  $\bb{F}_p[\epsilon]$
\begin{equation}\label{eqGradedRingsIsoLieALggLazard}
U(\mathrm{gr}^{\bullet}(G))\xrightarrow{\sim} \mathrm{gr}^{\bullet}(\bb{Z}_{p,\sol}[G])
\end{equation}
where:

\begin{itemize}

\item[i.] $\mathrm{gr}^{\bullet}(G)$ is the Lie algebra of the $p$-valued group $G$ as in \cite[II.1.1.7]{zbMATH03225421} and $U(\mathrm{gr}^{\bullet}(G))$ its enveloping algebra over $\bb{F}_p[\epsilon]$.  By modifying the valuation,   we can and will assume that $\mathrm{gr}^{\bullet}(G)$ is abelian, cf. \cite[III.3.1.12]{zbMATH03225421}.

\item[ii.]  $\mathrm{gr}^{\bullet}(G)$ is a finite free Lie algebra over $\bb{F}_p[\epsilon]$ by \cite[III.2.1.3]{zbMATH03225421} with basis given by the residue classes of the elements $g_1,\ldots, g_d$ with weights $w_1,\ldots, w_d$ respectively (cf.  \cite[Proposition III.2.2.5]{zbMATH03225421}).

\item[iii.]  The isomorphism \eqref{eqGradedRingsIsoLieALggLazard} sends $\overline{g}_i$ to $x_i:=\overline{b}_i$. Thus, we have an isomorphism of graded $\bb{F}_p[\epsilon]$-algebras $\gr^{\bullet}(\bb{Z}_{p,\sol}[G])\cong \bb{F}_p[\epsilon][x_1,\ldots, x_d]$ with $\epsilon$ of weight $1$, and $x_i$ of weight $w_i$.  

\end{itemize}

Following \cite[Section 4]{zbMATH01965478}, we introduce the following  Banach algebras:

\begin{definition}\label{DefinitionAnalyticAlgebrasBinomial}
  Let $\Lambda$ be a uniform Banach $\bb{Q}_p$-algebra with pseudo-uniformizer $\varpi$. Let $|-|_{\varpi}$ be the spectral norm of $\Lambda$ normalized as follows: for $\lambda\in \Lambda$, $|\lambda|_{\varpi} =\sup_{x\in \n{M}(\Lambda)} |\lambda|_{\varpi,x}$, where $|-|_{\varpi,x}$ is the unique multiplicative norm on the residue field $K(x)$ such that $|\varpi|_{\varpi,x}=1/2$. 
  For $r>1$ we define the $\varpi$-complete module 
\begin{equation}\label{eqBasisBinomFunctionsCont}
C_{r,\varpi}(G, \Lambda^{\circ}):= \{\sum_{\alpha\in \bb{N}^d} c_{\alpha} \binom{\mathbf{x}}{\alpha} \; : \; c_{\alpha}\in \Lambda,\; |c_{\alpha}|_{\varpi}r^{\tau \alpha}\leq 1\  \forall  \alpha \mbox{ and } |c_{\alpha}|_{\varpi}r^{\tau \alpha}\to 0 \mbox{ as }  \tau \alpha\to \infty\}. 
\end{equation}
We set $C_{r,\varpi}(G, \Lambda):=C_{r,\varpi}(G, \Lambda^{\circ})[\frac{1}{\varpi}]$. 
\end{definition}

By \Cref{sec:obbin_r-varpi-analyt-amice-theorem},  $C_{r,\varpi}(G,\Lambda)$ embeds into $C^\la(G,\Q_p)\otimes_{\Q_p}\Lambda$ and the solid $G$-representation $C_{r,\varpi}(G,\Lambda)$ is locally analytic.

\begin{lemma}\label{LemmaConvergengeBinomials}
Keep the notation of \Cref{DefinitionAnalyticAlgebrasBinomial}. Then, for $r>1$ the $\varpi$-complete submodule $C_{r,\varpi}(G, \Lambda^{\circ})\subset C(G, \Lambda^{\circ})$ is an algebra.  
\end{lemma}
\begin{proof}
We endow $\bb{N}^{d}$ with the partial order making $\alpha\leq \beta$ if and only if $\alpha_{i}\leq \beta_i$ for all $i=1,\ldots, d$.  In particular, if $\alpha\leq \beta$ then $\tau\alpha\leq \tau \beta$.  For $\alpha,  \beta\in \bb{N}^d$ we have 
\[
\binom{\mathbf{x}}{\alpha}\binom{\mathbf{x}}{\beta} = \sum_{\gamma\leq   \alpha+ \beta} n_{\alpha \beta, \gamma} \binom{\mathbf{x}}{\gamma}
\]
with $n_{\alpha \beta,  \gamma}\in \bb{Z}$.  Thus, if $f(x)= \sum_{\alpha} c_{\alpha} \binom{\mathbf{x}}{\alpha}$ and $g(x)= \sum_{\beta} d_{\beta} \binom{\mathbf{x}}{\beta}$ are elements in $C_{r,\varpi}(G, \Lambda^{\circ})$, we deduce
\[
f(x)g(x)= \sum_{\gamma}\bigg( \sum_{\gamma\leq \alpha+\beta} c_{\alpha} d_{\beta} n_{\alpha \beta, \gamma}  \bigg) \binom{\mathbf{x}}{\gamma}.
\]
Then, for $\gamma\leq \alpha+\beta$ we have that  $|c_{\alpha} d_{\beta} n_{\alpha \beta, \gamma}|_{\varpi}\leq r^{-(\tau \alpha+\tau \beta)}\leq r^{-\tau  \gamma}$.
Similarly, we have that $|c_{\alpha} d_{\beta} n_{\alpha\beta,\gamma}|_{\varpi} r^{\tau \gamma} \leq |c_{\alpha} d_{\beta}|_{\varpi} r^{\tau\alpha+\tau\beta} \to 0$ as $\tau \gamma\to \infty$.
\end{proof}

The following lemma shows that, for suitable choices of $\varpi$ and $r$,  the Hopf algebra structure of $C(G,\Lambda^{\circ})$ restricts to a Hopf algebra structure on $C_{r,\varpi}(G,\Lambda^{\circ})$.

\begin{lemma}\label{LemmaHopfALgebraStructureCBinom}
Suppose that $ |p|_{\varpi}\leq 1/2$ and let $2>r>1$.  Then the Hopf $\Lambda^{\circ}$-algebra structure on $C(G, \Lambda^{\circ})$ restricts to a Hopf algebra structure on $C_{r,\varpi}(G, \Lambda^{\circ})$.
\end{lemma}
\begin{proof}
Consider the $\Lambda^{\circ}$-dual $\n{D}_{r,\varpi}(G,\Lambda^{\circ})$ of $C_{r,\varpi}(G, \Lambda^{\circ})$.  In terms of the Mahler basis, $\n{D}_{r,\varpi}(G,\Lambda^{\circ})$  is the solid $\Lambda^{\circ}$-module given by the power series 
\begin{equation}\label{eqBinomDistAlgebrasFormulas33}
\n{D}_{r,\varpi}(G,\Lambda^{\circ})=\{\sum_{\alpha} a_{\alpha} \mathbf{b}^{\alpha} \; : \; |a_{\alpha}|_{\varpi}r^{-\tau \alpha}\leq 1 \},
\end{equation}
 where $\mathbf{b}^{\underline{\alpha}}= b_1^{\alpha_1}\cdots b_d^{\alpha_d}$, and endowed with the product topology with respect to the projections to the coefficients to each monomial. Writing $r^{-1}=|\varpi|^\delta_{\varpi}$ for some $\delta>0$, if  the root $\varpi^{\delta}$ exists in $\Lambda$ we have the explicit description 
 \[
 \n{D}_{r,\varpi}(G,\Lambda^{\circ}) = \prod_{\alpha} \Lambda^{\circ}  \frac{\mathbf{b}^{\alpha}}{\varpi^{ \delta \tau \alpha} }. 
 \]

  The space $C_{r,\varpi}(G,\Lambda^{\circ})$ can be recovered as the topological $\Lambda^{\circ}$-dual of $\n{D}_{r,\varpi}(G,\Lambda^{\circ})$, hence, to show that $C_{r,\varpi}(G,\Lambda^{\circ})$ is a Hopf algebra it suffices to show that  $\n{D}_{r,\varpi}(G,\Lambda^{\circ})$ has a (necessarily unique)  $\Lambda^{\circ}$-algebra structure compatible with that of $\bb{Z}_{p,\sol}[G]$, and that the antipode of $\bb{Z}_{p,\sol}[G]$ extends to an antipode on  $\n{D}_{r,\varpi}(G,\Lambda^{\circ})$.

We begin with the multiplication. We have to show that the multiplication map $m\colon \bb{Z}_{p,\sol}[G]\times \bb{Z}_{p,\sol}[G]$ extends to a    multiplication map $m\colon \n{D}_{r,\varpi}(G,\Lambda^{\circ})\times \n{D}_{r,\varpi}(G,\Lambda^{\circ})\to \n{D}_{r,\varpi}(G,\Lambda^{\circ})$  of solid $\Lambda^{\circ}$-modules. For this, it suffices to see that  $m$ extends to a morphism 
\begin{equation}\label{eqMultipplicationMapExtensionBinomDist}
m\colon \n{D}_{r,\varpi}(G\times G,\Lambda^{\circ})\to  \n{D}_{r,\varpi}(G,\Lambda^{\circ})
\end{equation}
where the left term consists of power series $\sum_{\alpha,\beta} a_{\alpha,\beta} \mathbf{b}^{\alpha} \widetilde{\mathbf{b}}^{\beta}$ such that $|a_{\alpha,\beta}|r^{-\tau\alpha-\tau \beta}\leq 1$ for all $\alpha,\beta\in \bb{N}^d$,  with respect to the product topology on the coefficients (where $\mathbf{b}$ are the coordinates of the Iwasawa algebra of the first term copy of $G$ while $\widetilde{\mathbf{b}}$ refers to the coordinates of the second copy).

For $\alpha, \beta$ multi-indices, write $\mathbf{b}^{\alpha} \widetilde{\mathbf{b}}^{\beta}=\sum_{\gamma} h_{\alpha  \beta,\gamma} \mathbf{b}^{\gamma}$ with $h_{\alpha \beta, \gamma} \in \bb{Z}_p$.   By \cite[Lemma 4.1]{zbMATH01965478} we have that 
 \begin{equation}\label{eqBoundCoeffMultiplicationIwasawaAlg}
 v_{p}(h_{\alpha\beta,\gamma})\geq \max(0,\tau\alpha+\tau \beta-\tau\gamma)
 \end{equation}
 where $v_p$ refers to the $p$-adic valuation. Hence, for  a power series $f(\mathbf{b},\widetilde{\mathbf{b}})= \sum_{\alpha,\beta} a_{\alpha,\beta} \widetilde{\mathbf{b}}^{\alpha} \mathbf{b}^{\beta}$, we get that 
 \[
 m(f)=\sum_{\gamma}\bigg( \sum_{\alpha,\beta} a_{\alpha,\beta} h_{\alpha\beta,\gamma} \bigg) \mathbf{b}^{\gamma}.  
 \]

 Therefore, \eqref{eqBoundCoeffMultiplicationIwasawaAlg}  implies that 
 \[
 |h_{\alpha \beta,\gamma} |_{\varpi}\leq \min(1,|p|_{\varpi}^{\tau \alpha +\tau \beta - \tau \gamma}) \leq  \min(1,2^{-(\tau \alpha+\tau \beta-\tau \gamma}))\leq \min(1, r^{-(\tau\alpha +\tau\beta-\tau \gamma}))
 \] 
 where the minimum is $1$ if $\tau\alpha +\tau \beta < \tau \gamma$, and it is equal to  the right terms otherwise.  The previous yields the following:
 
 \begin{itemize} 
 
 \item[i.] for $\gamma$ fixed,   $|a_{\alpha,\beta} h_{\alpha\beta, \gamma}|_{\varpi}\to 0$  as $\max(\tau \alpha,\tau \beta)\to \infty$,  so that the coefficients of the power series expansion of $m(f)$ converge. Indeed, we have that  
 \[
 |a_{\alpha,\beta}  h_{\alpha \beta, \gamma}|_{\varpi} \leq (2^{\tau  \gamma})( \frac{r}{2})^{\tau \alpha +\tau\beta} 
 \]  
 which converges to $0$ when $\max(\tau  \alpha, \tau  \beta )\to \infty$ as $r<2$. 
 \item[ii.] 
 For all $ \alpha, \beta$ and $ \gamma$ one has 
 \[
 |a_{\alpha,\beta}  h_{\alpha \beta, \gamma} |_{\varpi}\leq \min(r^{\tau \alpha+\tau \beta} , r^{\tau \gamma})\leq r^{\tau \gamma},
 \]
 proving that  $m(f)\in \n{D}_{r,\varpi}(G,\Lambda^{\circ})$ as desired.

\item[iii.] We have to show that the map \eqref{eqMultipplicationMapExtensionBinomDist} is continuous for the product topology on the distribution algebra.  For $f\in \n{D}_{r,\varpi}(G\times G,\Lambda^{\circ})$ write $m(f)=\sum_{\gamma} c_{\gamma}(f) \mathbf{b}^{\gamma}$. It suffices to prove that, for $\gamma$ fixed and $\epsilon>0$, there is some $N_0\gg 0$ such that $|c_{\gamma}(\mathbf{b}^{\alpha}\mathbf{b}^{\beta})|r^{\tau\alpha+\tau\beta-\tau\gamma}<\epsilon$ if $\max(\tau\alpha,\tau\beta)\geq N_0$. This last claim  is a direct consequence of \eqref{eqBoundCoeffMultiplicationIwasawaAlg}, and the inequalities $|p|_{\varpi}\leq 1/2$ and $1/2< r^{-1}<1$. Indeed, one has that 
\[
|c_{\gamma}(\mathbf{b}^{\alpha} \mathbf{b}^{\beta})|r^{\tau \alpha+ \tau \beta -\tau \gamma} \leq  \big(\frac{r}{2})^{-\tau \gamma} \big(\frac{r}{2}\big)^{\tau \alpha+\tau \beta}
\]  
and $\big(\frac{r}{2}\big)^{\tau \alpha+\tau \beta}$ decreases exponentially to $0$ as $\max(\tau\alpha,\tau\beta)\to \infty$.

 \end{itemize}
 
 Finally, we prove the claim for the antipode. We note that 
 \[
 g^{-1}_i= 1/(1+ (g_i-1)) = \sum_{n\in \bb{N}} (-1)^n b_i^n
 \]
 and that 
 \[
 \iota(b_i)= g_i^{-1}-1 =\sum_{n\geq 1} (-1)^n b_i^n= -b_i +\mathcal{O}(2)
 \]
 where for $k\in \bb{N}$ the notation  $\mathcal{O}(k)$ refers to a power series on the $b_i's$ with monomials of degree $\geq 2$ and integer coefficients. 
 
In particular, if $\iota$ denotes the antipode of $\bb{Z}_{p,\sol}[G]$, we see that 
\[
\iota(\mathbf{b}^{\alpha})=(-1)^{|\alpha|} \mathbf{b}^{\alpha} + \mathcal{O}(|\alpha|+1). 
\]
From the previous bound, an explicit estimate shows that if $f(\mathbf{b})\in \n{D}_{r,\varpi}(G,\Lambda^{\circ})$ then  $f(\iota(\mathbf{b}))$ also converges to an element in $\n{D}_{r,\varpi}(G,\Lambda^{\circ})$. It is immediate to see that the antipode is continuous with respect to the product topology on the coefficients.  This finishes the proof of the lemma. 
\end{proof}

We can now define the following groups in Gelfand stacks.

\begin{definition}\label{DefinitionGroupsBinoms}
Keep the notation of \Cref{DefinitionAnalyticAlgebrasBinomial} and suppose that $|p|_{\varpi}\leq 1/2$. For  $2>r>1$ we define the analytic group in Gelfand stacks $G_{\ob{Bin}_r,\varpi}:= \GSpec (C_{r,\varpi}(G,\Lambda))$. For a commutative algebra $R\in \ob{CAlg}(\ob{D}(B G_{\ob{Bin}_r,\varpi}))$ we define the category of $R$-semilinear $\ob{Bin}_{r,\varpi}$-analytic representations $\cat{Rep}_{R}^{\ob{Bin}_{r,\varpi}}(G):= \ob{Mod}_{R}(\ob{D}(B G_{\ob{Bin}_r,\varpi}))$.
\end{definition}

  The trivial representation $\bb{Z}_p$ admits a resolution by finite free $\bb{Z}_{p,\sol}[G]$  called the \textit{Lazard resolution}, see \cite[Definition V.2.2.2.1]{zbMATH03225421} of the form 
\begin{equation}\label{eqLazardResolutionZp22}
0\to X_{d}\to X_{d-1}\to \cdots \to X_1\to X_0\to \bb{Z}_p\to 0
\end{equation}
where $X_i\cong \bb{Z}_{p,\sol}[G]^{\binom{d}{i}}$.  We let  $\ob{LZ}^{\bullet}(\bb{Z}_{p,\sol}[G])$  denote the acyclic complex \eqref{eqLazardResolutionZp22}. The base change along $\bb{Z}_p\to \Lambda^{\circ}$ of the $\bb{Z}_p$-dual of \eqref{eqLazardResolutionZp22} produces an acyclic complex  of $C(G,\Lambda^{\circ})$-comodules
\begin{equation}\label{eqLazardDualContFunctions}
0\to \Lambda^{\circ}\to Y_{0}\to Y_{1}\to \cdots \to Y_{d}\to 0
\end{equation}
where $Y_{i}\cong C(G,\Lambda^{\circ})^{\binom{d}{i}}$. We denote this acyclic complex by $\ob{LZ}_{\bullet}(C(G,\Lambda^{\circ}))$. 

For the decompletion statements by locally analytic vectors in \Cref{sec:gener-lemma-decompl}, we will need to consider variants of the Lazard resolution \eqref{eqLazardDualContFunctions} for the binomial functions $C_{r,\varpi}(G,\Lambda^{\circ})$. In this case, the associated complex will not be exact, but we can manage to  control the failure  of this exactness. The ideas presented next are essentially those encoded in the proof of \cite[Theorem 4.4]{MR2774315}, and are  basically taken line to line from the proof of \cite[Proposition 6.6]{PoratSolidLocAn} with some minor modifications.

By definition, the Lazard resolution is a complex of filtered $\bb{Z}_{p,\sol}[G]$-modules. Furthermore, by the lemma of Serre used in the construction \cite[Lemme V.2.1.1]{zbMATH03225421}, the complex admits a $\bb{Z}_p$-linear homotopy $s^{\bullet}\colon \ob{LZ}^{\bullet}(\bb{Z}_{p,\sol}[G]) \to \ob{LZ}^{\bullet+1}(\bb{Z}_{p,\sol}[G])$ compatible with the filtrations.   Moreover, the graded complex \eqref{eqLazardResolutionZp22} is precisely the Koszul complex $\mathrm{Kos}^{\bullet}( \bb{F}_p[\epsilon][x_1,\ldots, x_n])$ of the ring $\bb{F}_p[\epsilon]$ seen as $\gr^{\bullet}(\bb{Z}_{p,\sol}[G])\cong \bb{F}_p[\epsilon][x_1,\ldots, x_n]$-algebra, and $\mathrm{gr}(s^{\bullet})$ is a homotopy between $1$ and $0$ preserving the graded pieces. Thus, if we write $X_i =\bigoplus_{|J|=i}  \bb{Z}_{p,\sol}[G]e_J$ where $J\subset \{1,\ldots, d\}$, and $e_J$ is a lift of weight $w_J:=\sum_{j\in J} w_j$ of the basis element $x_J:= \bigwedge_{j\in J} x_{j}$ (up to a sign) in the Koszul complex $\mathrm{Kos}^{\bullet}( \bb{F}_p[\epsilon][x_1,\ldots, x_n])$, the map $s^{i}$ must restrict to 
\begin{equation}\label{eqMapsFilteredModulesTechnical233ojaosd}
s^i\colon \bigoplus_{|I|=i} \ob{Fil}^s(\bb{Z}_{p,\sol}[G]) e_{I} \to   \bigoplus_{|I|=i+1} \ob{Fil}^{s-\omega}(\bb{Z}_{p,\sol}[G]) e_{I}
\end{equation}
with $\omega:= \sum_{I,J\subset \{1,\ldots, d\}} \{|w_J-w_I|\}$.

The following lemma gives the precise torsion control on the Lazard resolution for the spaces of binomial functions.

\begin{lemma}\label{LemmaExtensionLazardResolutionDist}
Suppose that $|p|_{\varpi}<\frac{1}{2}$ and let $2>r>1$ be a real number. Suppose that there is a topologically nilpotent element $\pi\in \Lambda^{\circ}$ such that for all $x\in \n{M}(\Lambda)$ one has $|\pi(x)|_{\varpi,x}=r^{-\omega}$.\footnote{A pseudo-uniformizer $\pi$ such that $|\pi(x)|_{\varpi,x}\leq r^{-\omega}$ for all $x\in \n{M}(\Lambda)$ suffices.}
The following hold:

\begin{enumerate}

\item The Lazard complex $\ob{LZ}^{\bullet}(\bb{Z}_{p,\sol}[G])$ extends uniquely to a $\Lambda^{\circ}$-linear complex  of $\n{D}_{r,\varpi}(G,\Lambda^{\circ})$-modules
\[
0\to X_{d}^{r}\to X_{d-1}^{r}\to \cdots \to X_1^{r}\to X_0^r\to \Lambda^{\circ}\to 0,
\]
where $X_{i}^r\cong \n{D}_{r,\varpi}(G,\Lambda^{\circ})$ are the distribution algebras \eqref{eqBinomDistAlgebrasFormulas33}.  We let $\ob{LZ}^{\bullet}(\n{D}_{r,\varpi}(G,\Lambda^{\circ}))$ denote this  complex.  

\item The Lazard complex $\ob{LZ}_{\bullet}(C(G,\Lambda^{\circ}))$ restricts uniquely to a complex of $C_{r,\varpi}(G,\Lambda^{\circ})$-comodules
\[
0\to \Lambda^{\circ}\to Y_{0}^r\to Y_1^r\to \cdots\to Y_{d}^r\to 0
\] 
with $Y_i^r\cong C_{r,\varpi}(G,\Lambda^{\circ})^{\binom{d}{i}}$. We let $\ob{LZ}_{\bullet}(C_{r,\varpi}(G,\Lambda^{\circ}))$ denote this  complex.

\item The homotopy $s^{\bullet}$ of $\ob{LZ}^{\bullet}(\bb{Z}_{p,\sol}[G])$ extends uniquely to a homotopy
\[
s^{\bullet}\colon \ob{LZ}^{\bullet}(\n{D}_{r,\varpi}(G,\Lambda^{\circ}))\to \pi^{-1} \ob{LZ}^{\bullet+1}(\n{D}_{r,\varpi}(G,\Lambda^{\circ})).
\]
Taking continuous $\Lambda^{\circ}$-duals, the homotopy $s^{\bullet}$ then induces a homotopy of $\Lambda^{\circ}$-modules
\[
s^{\bullet,\vee}\colon \ob{LZ}_{\bullet+1}(C_{r,\varpi}(G,\Lambda^{\circ}))\to \pi^{-1} \ob{LZ}_{\bullet}(C_{r,\varpi}(G,\Lambda^{\circ})).
\]

\end{enumerate}

In particular, $\Gamma(G, C_{r,\varpi}(G, \Lambda^{\circ}))$ has $\pi$-torsion higher cohomology groups and  $\Gamma(G, C_{r,\varpi}(G, \Lambda))=\Lambda$. 
\end{lemma}
\begin{proof}
By taking topological $\Lambda^{\circ}$-duals (where $C_{r,\varpi}(G,\Lambda^{\circ})$ is a Banach $\Lambda^{\circ}$-module with basis $\binom{\mathbf{x}}{\alpha}$ as in \eqref{eqBasisBinomFunctionsCont}, and $\n{D}_{r,\varpi}(G,\Lambda^{\circ})$ is its $\Lambda^{\circ}$-dual endowed with the product topology for the basis $\mathbf{b}^{\alpha}$ as in \eqref{eqBinomDistAlgebrasFormulas33}), it suffices to show part (1), and part (3) for $\n{D}_{r,\varpi}(G,\Lambda^{\circ})$.  The extension of $\ob{LZ}^{\bullet}(\bb{Z}_{p,\sol}[G])$ to $\ob{LZ}^{\bullet}( \n{D}_{r,\varpi}(G,\Lambda^{\circ}))$  in part (1) is clear by taking a base change along the morphism of rings $\bb{Z}_{p,\sol}[G]\to \n{D}_{r,\varpi}(G,\Lambda^{\circ})$, it is only left to prove the extension of the homotopy $s^{\bullet}$ in (3). For that, consider the homotopy 
\[
s^i\colon X_i\cong \bigoplus_{|I|=i} \bb{Z}_{p,\sol}[G] e_{I}\to  \bigoplus_{|J|=i+1}  \bb{Z}_{p,\sol}[G] e_{J}.
\]
Wee need to show that the restriction $s^{i}_{I,J}\colon  \bb{Z}_{p,\sol}[G] e_{I}\to \bb{Z}_{p,\sol}[G] e_{J}$ extends to a (necessarily unique) continuous map of $\Lambda^{\circ}$-modules $s^i_{I,J}\colon \n{D}_{r,\varpi}(G,\Lambda^{\circ}) e_{I}\to \pi^{-1} \n{D}_{r,\varpi}(G,\Lambda^{\circ}) e_{I}$. To simplify the notation, we will omit the basis elements $e_I,e_J$ and write $s:= s^i_{I,J}$.  Given $\alpha\in \bb{N}$ we write $s(\mathbf{b}^{\alpha})=\sum_{\beta} d_{\alpha,\beta} \mathbf{b}^{\beta}$. Thanks to \eqref{eqMapsFilteredModulesTechnical233ojaosd}, for  $\tau \beta< \tau\alpha-\omega$ we have that $v_p(d_{\alpha,\beta})\geq \tau \alpha-\tau\beta-\omega$, where $v_p$ is the $p$-adic valuation.  Thus, for all $\alpha,\beta$ we have 
\begin{equation}\label{eqBoundsDsLazardResolution}
|d_{\alpha,\beta}|\leq \min(1,|p|_{\omega}^{\tau \alpha-\tau\beta -\omega}) \leq \min(1,2^{-(\tau \alpha-\tau\beta -\omega)}) \leq \min(1, r^{-(\tau \alpha-\tau\beta -\omega)}),
\end{equation}
where if $\tau \beta\geq \tau \alpha -\omega$ the minimum above is $1$, and if $\tau\beta<\tau\alpha-\omega$ the minimum are the right terms.

 Let $f(\mathbf{b})=\sum_{\alpha\in \bb{N}^d} a_{\alpha} \mathbf{b}^{\alpha}\in \n{D}_{r,\varpi}(G,\Lambda^{\circ})$, we have that 
\[
s(f(\mathbf{b}))=\sum_{\beta} \bigg( \sum_{\alpha} a_{\alpha} d_{\alpha,\beta}  \bigg) \mathbf{b}^{\beta}=\sum_{\alpha} e_{\beta} \mathbf{b}^{\beta}. 
\] 
To prove the lemma, we need to prove the following facts:

\begin{enumerate}

\item For fixed $\beta$, the norms $|a_{\alpha}d_{\alpha,\beta}|_{\varpi}\to 0$ as $\tau\alpha\to \infty$. In particular, the series defining $e_{\beta}$ converge in $\Lambda$. 

\item  For all $\beta$, we have $|a_{\alpha}d_{\alpha,\beta}|_{\varpi}\leq r^{\tau \beta+\omega}$. This implies that $|e_{\beta}|\leq r^{\tau\beta+\omega}$ and that $s(f)$ converges to an element  in $\pi^{-1} \n{D}_{r,\varpi}(G,\Lambda^{\circ})$.

\item For all $\beta$, the $\Lambda^{\circ}$-linear maps $\n{D}_{r,\varpi}(G,\Lambda^{\circ})\to \pi^{-1}\Lambda^{\circ}$ defined by $f\mapsto e_{\beta}$ are continuous. This will imply that $s\colon \n{D}_{r,\varpi}(G,\Lambda^{\circ})\to \pi^{-1} \n{D}_{r,\varpi}(G,\Lambda^{\circ})$ is continuous.

\end{enumerate}

Recall that by \eqref{eqBinomDistAlgebrasFormulas33}  $|a_{\alpha}|_{\varpi}\leq r^{\tau \alpha}$.  This implies that, fixing $\beta$, we have  
\[
|a_{\alpha} d_{\alpha,\beta}|_{\varpi}\leq  r^{\tau \alpha} 2^{-\tau \alpha} 2^{\tau\beta +\omega} =  (\frac{r}{2})^{\tau \alpha} 2^{\tau\beta +\omega} 
\]
which converges to $0$ as $\tau\alpha\to \infty$ since $r<2$. This proves (1).  For (2), we also deduce from \eqref{eqBoundsDsLazardResolution} that 
\[
|a_{\alpha} d_{\alpha,\beta}|_{\varpi} \leq   |a_{\alpha}| r^{-\tau \alpha}  r^{\tau\beta +\omega} \leq   r^{\tau\beta +\omega} 
\]
which yields $|e_{\beta}|\leq r^{\tau\beta+\omega}$.  Finally, for (3), it suffices to show that the map $f\mapsto \pi e_{\beta}$ from $\n{D}_{r,\varpi}(G,\Lambda^{\circ})\to \Lambda^{\circ}$ is continuous, that is, for all $k\geq 1$ there is some $N\gg 0$ such that the composite 
\[
\Psi\colon \n{D}_{r,\varpi}(G,\Lambda^{\circ})\to \Lambda^{\circ} \to \Lambda^{\circ}/\varpi^k
\]
factors through the projection  $\n{D}_{r,\varpi}(G,\Lambda^{\circ})\to M_{N}$, with $M_{N}$ the $\varpi$-Banach $\Lambda^{\circ}$-module given by  
\[
M_{N}=\{\sum_{\tau \alpha\leq N} a_{\alpha} \mathbf{b}^{\alpha } \; : \; a_{\alpha} \in \Lambda, \;\; |a_{\alpha}|r^{-\tau \alpha}\leq 1 \}. 
\]
It suffices to show that $|\Psi(\mathbf{b}^{\alpha})|_{\varpi} \leq 2^{-k}$ for $\tau \alpha \geq N$. Take $N$ such that  $\big(\frac{r}{2}\big)^{N}  < 2^{-k -(\tau \beta + \omega)}$ and  $N-\tau \beta -\omega >0$. Then, for $\tau\alpha\geq N$ we have that  $\tau \alpha -\tau \beta -\omega  \geq  N-\tau \beta -\omega \geq 0$ an by \eqref{eqBoundsDsLazardResolution} that 
\[
|a_{\alpha} d_{\alpha,\beta}|_{\varpi} \leq |a_{\alpha}|_{\varpi} 2^{-(\tau \alpha -\tau \beta -\omega)}\leq  \bigg(\frac{r}{2}\bigg)^{\tau \alpha}   2^{\tau \beta + \omega} \leq \bigg(\frac{r}{2}\bigg)^{N} 2^{\tau \beta + \omega}  < 2^{-k}. 
\]
This finishes the proof of the lemma.
\end{proof}

\begin{remark}\label{RemarkBetterBetterBounds}
  Note that the $\varpi$-norm of $ \pi$ in \Cref{LemmaExtensionLazardResolutionDist} converges to $1$ as $r\to 1$. This implies that the higher cohomology groups of $C_{r,\varpi}(G,\Lambda^{\circ})$ are of smaller torsion as $r$ approaches $1$.
  In particular, if $\Lambda$ contains enough roots of $\varpi$, taking colimits the higher $G$-cohomology groups of $\varinjlim_{r\to 1} C_{r,\varpi}(G,\Lambda^{\circ})$ are almost zero, with respect the the almost setup $(\Lambda^{\circ}, \Lambda^{\circ\circ})$.
\end{remark}

 The following lemma implies good properties of the groups $G_{\ob{Bin}_r,\varpi}$. 
 
\begin{lemma}\label{LemmaGoodShierkPropertiesGBinom}
Let $2>r>1$. The natural map $f\colon BG^{\la}\times_{\GSpec(\Q_p)} \GSpec (\Lambda) \to BG_{\ob{Bin}_r,\varpi}$ is $0$-proper and satisfies $f_*\Lambda=\Lambda$. In particular $*$-pullback along $f$ is fully faithful.
\end{lemma}
\begin{proof}
  The natural morphism of groups $G^{\la}\times_{\GSpec(\Q_p)} \GSpec (\Lambda) \to G_{\ob{Bin}_r,\varpi}$ arises from the natural inclusion of Hopf algebras $C_{r,\varpi}(G, \Lambda)\subset C^{\la}(G, \Lambda)$.
  The fact that $f$ is $0$-proper follows formally from the fact that $G_{\ob{Bin}_r,\varpi}$ is $0$-proper over $\GSpec (\Lambda)$ (being affinoid), and that $BG^{\la}$ is $0$-proper by \Cref{LemmaAnalyticRepresentationProperty} (1). Finally,  the claim   $f_*\Lambda=\Lambda$ is equivalent to the fact that the locally analytic cohomology, i.e., the cohomology of $BG^\la$, of the $G^\la$-representation $C_{r,\varpi}(G,\Lambda)$ is just $\Lambda$.
  Since this cohomology can be computed simply as solid group cohomology, that is via taking Hom out of the Lazard resolution, this follows from \Cref{LemmaExtensionLazardResolutionDist}.
\end{proof}

\begin{remark}
  \label{sec:obbin_r-varpi-analyt-integral-torsion}
  We note that the Hopf algebra $C_{r,\varpi}(G,\Lambda)$ defining $G_{\ob{Bin}_r,\varpi}$ admits the $\varpi$-complete model $C_{r,\varpi}(G,\Lambda^\circ)$ and by \Cref{LemmaExtensionLazardResolutionDist}, the higher solid cohomology (as a $\Z_{p,\sol}[G]$-representation) of $C_{r,\varpi}(G,\Lambda^\circ)$ is torsion by \Cref{RemarkBetterBetterBounds}.
  This should be compared to \Cref{RemarkRAnalyticFailure}, where the natural integral models for $\mathbb{G}_r$ have higher non-torsion cohomology.
\end{remark}

\subsubsection{Relative decompletions} Let $G$ be a compact $p$-adic Lie group as in \Cref{sss-loc-an-Reps}.  If $r<1$, the natural morphisms of groups on Gelfand stacks $G^{\la}\to \bb{G}_{<r}\to \bb{G}_{\leq r}$ produce fully faithful pullback maps at the level of their representation categories.  Moreover, $\cat{Rep}^{\la}_{\bb{Q}_p}(G)$ embeds fully faithfully in $\cat{Rep}^{\sol}_{\bb{Q}_p}(G)$. Hence,  by an abuse of notation,  for $M$  a $\bb{Q}_p$-solid  $G$-representation    we identify the $<r$ and $\leq r$-analytic vectors of $M$, denoted as  $M^{< r-\an}$ and $M^{\leq r-\an}$ respectively, as objects in $\cat{Rep}^{\sol}_{\bb{Q}_p}(G)$, $\cat{Rep}^{\la}_{\bb{Q}_p}(G)$, $\cat{Rep}^{<r}_{\bb{Q}_p}(G)$ or $\cat{Rep}^{\leq r}_{\bb{Q}_p}(G)$ accordingly.

  If $R\in \ob{CAlg}(\cat{Rep}^{\sol}_{\bb{Q}_p}(G))$,  we have fully  faithful embeddings
\begin{equation}\label{eqChainFullyFaithfulEmbeddingsReps}
\cat{Rep}^{\leq r}_{R^{\leq r-\an}}(G) \to \cat{Rep}^{\la}_{R^{\la}}(G) \to \cat{Rep}_{R}^{\sol}(G). 
\end{equation}
Indeed, all inclusion are formal consequences of the projection formulas satisfied by the functors of  locally analytic and $\leq r$-analytic vectors.

\begin{definition}\label{DefDecompletionLocAn}
Consider the chain of fully faithful functors \eqref{eqChainFullyFaithfulEmbeddingsReps}. 
 We say that an object $M\in \cat{Rep}^{\sol}_{R}(G)$ admits a \textit{locally analytic decompletion} (resp.\ an \textit{$\leq r$-analytic decompletion})  as an $R$-module if it belongs to the essential image of $\cat{Rep}^{\la}_{R^{\la}}(G)$ (resp.\ of $\cat{Rep}^{\leq r}_{R^{\leq r-an}}(G)$).
Equivalently, $M$ admits a locally analytic  decompletion as an $R$-module  if and only if the natural map $M^{\la}\otimes_{R^{\la}} R\to M$  is an equivalence (and similarly for $\leq r$). 
\end{definition}

\begin{remark}\label{RemarkGContVersionDecompletion}
Given a commutative algebra $R\in  \ob{CAlg}(\ob{D}(BG^{\ob{cont}}))$ one can define a notion of relative decompletion using the $*$-push along $BG^{\ob{cont}}\to BG^{\la}$. Thanks to \Cref{Lemma-Prim-padic} this notion is equivalent to that of \Cref{DefDecompletionLocAn} after taking composition with  the fully faithful functor $\Phi\colon \ob{D}(BG^{\ob{cont}})\to \cat{Rep}^{\sol}_{\bb{Q}_p}(G)$.  
\end{remark}

\begin{remark}
  \label{sec:morphism-bgrm-contto-locally-analytic-vectors-for-products}
  Let $G$ and $H$ be $p$-adic Lie groups and $M\in \ob{D}(B(G^\cont\times H^\cont))$, then  $M^{(G\times H)-\la}\cong (M^{G-\la})^{H-\la}$. Indeed,  this follows by \Cref{Lemma-Prim-padic} applied to the composite  of morphisms  Gelfand stacks $B(G^\cont\times H^\cont)\to B(G^\la\times H^\cont)\to B(G^\la\times H^\la)$,  and by  proper base change to describe pushforward along the first morphism. 
\end{remark}

The following lemma shows that descendability is independent of the action of the group. 

\begin{lemma}\label{LemmaDescendable-Independent}
Let $f\colon A\to B$ be a morphism of  $G$-equivariant solid $\Q_p$-algebras with $G$ a compact $p$-adic Lie group.    Then  $f$ is descendable in the symmetric monoidal category $\cat{Rep}^{\sol}_{\Q_p}(G)$ if and only if   the underlying map of $\Q_p$-algebras is descendable.   Furthermore, if $R$ is a commutative algebra on $\cat{Rep}^{\sol}_{\Q_p}(G)$, then the morphism $R^{\la}\to R$ is descendable.   The same holds for a morphism of algebras in $\ob{D}(BG^{\rm cont})$ with $BG^{\rm cont}=\GSpec(\Q_p)/G^{\rm cont}$. 
\end{lemma}
\begin{proof}
  Since the forgetful functor $F\colon \cat{Rep}^{\sol}_{\Q_p}(G)\to \ob{D}(\Q_p)$ is symmetric monoidal, if $f$ is descendable as $G$-equivariant algebras then so is as just $\Q_p$-algebras. Conversely, the right adjoint  $F_*$ sending a $\Q_{p,\sol}$-vector space $M$ to $ C(G,M):=\underline{\Hom}_{\Q_p}(\Q_{p,\sol}[G], M)$ is lax symmetric monoidal, and therefore it sends descendable maps to descendable maps.
  Thus, $C(G,A)\to C(G,B)$ is a descendable map of $G$-equivariant algebras.
  Hence, it suffices to show that the map $A\to C(G,A)$ is descendable as $G$-equivariant algebras, but this follows from the Lazard resolution of $G$ as then $A$ is a finite limit of finite copies of $C(G,A)$ seen as a $G$-equivariant $A$-module.
  The same argument holds for algebras in $\ob{D}(BG^{\rm cont})$ where the right adjoint of the pullback functor to $\GSpec(\Q_p)$ sends $M$ to the induced representation $C(G,\Q_p)\otimes_{\Q_p} M$.  

Finally, the last claim about descendability of $R^{\la}\to R$ follows from the fact that $R^{\la}= \Gamma(G, C^{\la}(G,\Q_p)\otimes R)$ and that by Lazard's resolution the group cohomology is a finite limit of finite copies of $C^{\la}(G,\Q_p)\otimes R$ as $R^{\la}$-module, so that $R^{\la}$ belongs to the thick tensor ideal generated by $R$ as $R^{\la}$-module. 
\end{proof}

We note the following invariance properties of being locally analytic.

\begin{lemma}\label{Lemma-Decompletions}
Let $A\to B$ be a map of $G$-equivariant solid $\Q_p$-algebras. 

\begin{enumerate}
\item Suppose that $B$ admits a decompletion by locally analytic vectors as a $G$-equivariant $A$-module and let $M$ be a $G$-equivariant $B$-module.  Then there is a natural isomorphism 
\[
M^{\la}\otimes_{A^{\la}} A = M^{\la}\otimes_{B^{\la}} B.
\]
In particular, $M$ admits a decompletion as an $A$-module if and only if it does as a $B$-module.

\item Let $N$ be a $G$-equivariant $A$-module admitting a decompletion by locally analytic vectors.
  Then $N\otimes_A B$ admits a decompletion by locally analytic vectors as a $B$-module and $(N\otimes_AB)^\la\cong N^\la\otimes_{A^\la}B^\la$. 
  \item Conversely, assume that $A\to B$ is descendable and that $B$ admits a decompletion by locally analytic vectors as a $G$-equivariant $A$-module. Let $N$ be a $G$-equivariant $A$-module, such that $N\otimes_A B$ admits a locally analytic decompletion as a $B$-module. Then $N$ admits a locally analytic decompletion as an $A$-module.
\end{enumerate}
\end{lemma}
\begin{proof}
\begin{enumerate}
\item By assumption the natural map $B^{\la}\otimes_{A^{\la}} A\to B$  is an isomorphism. Let $M$ be a $G$-equivariant  $B$-module, we have that 
\[
M^{\la}\otimes_{A^{\la}} A= M^{\la}\otimes_{B^{\la}} (B^{\la}\otimes_{A^{\la}} A) = M^{\la}\otimes_{B^{\la}} B
\]

\item Let $N$ be an $A$-module admitting a decompletion by locally analytic vectors as $A$-module, then the natural map $N^{\la}\otimes_{A^\la} A\to N$ is an isomorphism and we have 
\[
(N\otimes_A B)^{\la}\otimes_{B^{\la}} B= (N^{\la}\otimes_{A^{\la}} B)^{\la}\otimes_{B^{\la}} B = N^{\la}\otimes_{A^{\la}}  B^{\la} \otimes_{B^{\la} } B= (N^{\la}\otimes_{A^{\la}} A)\otimes_A B = N\otimes_A B 
\]
using the projection formula for $(-)^\la$ in the second equality.

\item By assumption, we have
\[ 
(N\otimes_A B)^\la \otimes_{B^\la} B = N\otimes_A B.
\]
We can rewrite the left-hand side as
\[ 
(N\otimes_A B)^\la \otimes_{B^\la} B = (N\otimes_{A^\la} B^\la)^\la \otimes_{B^\la} B = N^\la \otimes_{A^\la} B^\la \otimes_{B^\la} B = N^\la \otimes_{A^\la}  B =(N^\la \otimes_{A^\la} A) \otimes_{A} B,
\]
where in the first equality we used the assumption that $B$ admits a locally analytic decompletion as an $A$-module and the projection formula in the second equality.
Comparing the two formulas, we get
\[
(N^\la \otimes_{A^\la} A) \otimes_{A} B = N\otimes_A B.
\]
Since $A\to B$ is descendable, the natural map $N^\la \otimes_{A^\la} A\to N$ must be an isomorphism. 
\qedhere
\end{enumerate}
\end{proof}

Locally analytic decompletion on semidirect products can be checked individually on each factor in an appropriate sense. 

\begin{lemma}\label{LemmaDecompletionsSemidirectProducrs}
  Let $G=H\ltimes N$ be a compact $p$-adic Lie group written as a semidirect product. Let $R$ be a  solid $\bb{Q}_p$-algebra  where $X=\GSpec (R)$ is endowed with an action of $G^{\ob{cont}}$, and let $M\in \ob{D}(X/G^{\ob{cont}})$.
  Suppose that $M$ admits a locally analytic decompletion over $R$ as an $N$-representation and that $M^{N-\la}$ admits a decompletion  over $R^{N-\la}$ as an $H$-representation.
  Then $M$ admits a decompletion over $R$ as a $G$-representation and it is given by 
\[
M^{G-\la}= (M^{N-\la})^{H-\la}.
\] 
\end{lemma}
\begin{proof}
Recall that for a $p$-adic Lie group $K$ one has a natural morphism of groups $K^{\ob{cont}}\to K^{\la}$.   Consider the following composition of Gelfand stacks 
\[
B(H^{\ob{cont}}\ltimes N^{\ob{cont}})\xrightarrow{f} B(H^{\ob{cont}}\ltimes N^{\la}) \xrightarrow{h} B(H^{\la}\ltimes N^{\la}). 
\]
Both the maps $f$ and $h$ are $0$-proper, thus, by proper  base change, $*$-pushforward along $f$ is the functor $(-)^{N-\la}$ while $*$-pushforward along $h$ is the functor $(-)^{H-\la}$.  In particular, we have a natural equivalence  $(-)^{G-\la}=((-)^{N-\la})^{H-\la}$ of functors $\ob{D}(BG^{\ob{cont}})\to \ob{D}(\bb{Q}_p)$.  Applying this to $R$, we get fully faithful embeddings 
\[
\ob{Mod}_{R^{G-\la}}\ob{D}(B(H^{\la}\ltimes N^{\la}))\hookrightarrow  \ob{Mod}_{R^{N-\la}}\ob{D}(B(H^{\ob{cont}}\ltimes N^{\la})) \hookrightarrow  \cat{Rep}^{\sol}_{R}(G).
\]
Then, by proper base change applied to $BH^{\ob{cont}}= BH^{\la}\times_{B(H^{\la}\ltimes N^{\la})} B(H^{\ob{cont}}\ltimes N^{\la})$, the essential image of the left inclusion is given precisely by  those $R^{N-\la}$-modules with an action of $H^{\ob{cont}}\ltimes N^{\la}$ which admit a decompletion as $H$-representations. The lemma follows.  
\end{proof}

\begin{remark}\label{lemmaSolidEquivariantVariant}
  Using the languange of Gestalten in the proof of \Cref{LemmaContinuousRepsvsSolid} one can prove a variant of \Cref{LemmaDecompletionsSemidirectProducrs} for solid $G$-representations with exactly the same argument.
  This will  not be necessary in the future sections, so we leave the details to the reader. 
\end{remark}

We will in particular be interested in $G^{\rm cont}$-actions on Gelfand stacks.

\begin{lemma}\label{LemmaDecompletionsIndependent}
  Let $X=\GSpec (R)$ be an affinoid Gelfand stack endowed with an action of $G^{\rm cont}$ with $G$  a (non necessarily compact) $p$-adic Lie group.
  Let $f\colon X/G^{\rm cont}\to BG^{\la}$ be the natural projection map, and let $\Phi\colon \ob{D}(X/G^{\rm cont})\to \cat{Rep}^\sol_R(G)$ be the natural functor. 
  For $M\in \ob{D}(X/G^\cont)$, there is a natural isomorphism $f_* M =\Phi(M)^{\la}$ where $\Phi(M)^{\la}\in \cat{Rep}^\la_{\Q_p}(G)$ is the locally analytic decompletion of $\Phi(M)$.
  In particular, $f_\ast(\s{O}_X)=R^\la$ as commutative algebras in $\cat{Rep}^\la_{\Q_p}(G)$.
\end{lemma}
To ease notation, we will abbreviate the notation and write $(-)^\la$ instead of $\Phi(-)^\la$.
\begin{proof}
This follows formally form the fully faithfulness of $\Phi$ (cf. \Cref{LemmaContinuousRepsvsSolid}) and  from \Cref{Lemma-Prim-padic}. 
\end{proof}

In fact, we see that in \Cref{LemmaDecompletionsIndependent} the representation  $M^\la$ is naturally an object in $\cat{Rep}^\la_{R^\la}(G)$.
If $R^\la$ is a Gelfand ring, in particular connective, then we see that the functor $f_\ast$ lifts to a natural functor
\[
  \ob{D}(X/G^{\rm cont})\to \ob{D}(\GSpec(R^\la)/G^\la),
\]
which is given by $(-)^\la$.

\subsection{Part (1) of \Cref{TheoDecompletion}}
\label{sec:part-1-crefth-part1}

We now prove part (1) of \Cref{TheoDecompletion}. First,   we prove in \Cref{Lemma-Etale-localization} that actions of $p$-adic Lie groups lift along Berkovich-\'etale maps. 

We also note that if $H\subseteq G$ is an open subgroup, then the natural maps $H^{\ob{cont}}\to G^{\rm cont}\times_{G^\la}H^\la\to G^{\rm cont}\times_{G^\sm}H^\sm$ are isomorphisms.
Using the cover $|X|_\Betti/H^\la\to |X|_\Betti/G^\la$ reduces part (1) of the theorem   to the case that $G$ is compact.

We now check in greater generality that we can localize on $X$.

\begin{lemma}\label{Lemma-Etale-localization}
  Let $R$ be a Gelfand $\Q_p$-algebra and $\Pi$ a (light) profinite group.
  Consider an action of $\Pi^{\rm cont}$ on $\GSpec(R)$.
  Let $R\to R'$ be a Berkovich-\'etale morphism (\Cref{sec:prop-morph-gelf-5-berkovich-etale-smooth}).
\begin{enumerate}

\item The action of $\Pi$  induces a continuous action on the Banach algebra $R^u$.

\item Consider the action of $\Pi$ on $R^u$, then there is an open subgroup $\Pi_0\subset \Pi$ and an extension of the action of $\Pi_0$ to $R^{\prime,u}$.
  Furthermore, if $\Pi_0,\Pi_1\subset \Pi$ are open subgroups whose actions extend to $R^{\prime,u}$, there is an open subgroup $\Pi_2\subset \Pi_0\cap \Pi_1$ such that both actions restricted to $\Pi_2$ coincide.

\item Let $\Pi_0\subset \Pi$ be an open subgroup and consider an extension to an action $\rho$ on $R^{\prime,u}$.
  Then there is a natural extension of the action of  $\Pi_0$  from $R$ to $R^{\prime}$ completely determined by $\rho$.
  More precisely, this action is depicted by a pullback diagram 
\begin{equation}\label{eq0o2awjioqwe}
\begin{tikzcd}
\GSpec (R')/\Pi^{\rm cont}_0 \ar[r] \ar[d] & \Marc(R')^{\dR}/ \Pi_0^{\sm} \ar[d] \\
\GSpec (R)/\Pi^{\rm cont} \ar[r] & \Marc(R)^{\dR}/\Pi^{\sm} . 
\end{tikzcd}
\end{equation}
\end{enumerate}
\end{lemma}
\begin{proof}

  For part (1), let $R$ be a Gelfand ring.
  We claim that for $S$ a profinite set, the uniform completion of $C(S,\Q_p)\otimes_{\Q_p} R$ is given by $C(S,R^u)$.
  Indeed, by  \cite[Corollary 2.3.7]{drff_stacks} the functor $A\mapsto A^u$ is the left adjoint of the fully faithful embedding of uniformly completed bounded algebras into all bounded $\Q_p$-algebras. Thus, $(C(S,\Q_p)\otimes_{\Q_p} R)^u$  is also given by the uniform completion of $C(S,\Q_p)\otimes_{\Q_p} R^u=C(S,R^u)$ which is itself.
  Hence, if $(C(\Pi^{\bullet},\Q_p)\otimes R)_{[\bullet]\in \Delta}$ is the cosimplicial diagram that depicts the $C(\Pi,\Q_p)$-comodule structure on $R$, passing to uniform completions we obtain a cosimplicial diagram of Banach algebras $(C(\Pi^{\bullet},R^u))_{\Delta}$ endowing $R^u$ with a natural $C(\Pi,\Q_p)$-comodule structure, i.e., we can see $R^u$ as an object in $\ob{D}(\GSpec (\Q_p)/\Pi^{\rm cont})$.
  It is clear by construction that this $\Pi^{\rm cont}$-action is compatible with the algebra structure of $R^u$, so that it gives rise to an action of $\Pi^{\rm cont}$ on $\GSpec( R^u)$ as wanted.

  Part (2) follows from the same argument as in \cite[Lemma 2.3]{ScholzeLubinTate}.
  Indeed, suppose first that $R\to R'$ is a $\dagger$-rational localization.
  Then there is an open subgroup $\Pi_0\subset \Pi$ that leaves $\GSpec (R')\subset \GSpec (R)$ stable, namely, using adic spaces write the associated  Tate-rational localization and its adic compactification $U:=\Spa (R^{\prime,u},R^{\prime,u,+}) \to Y:=\Spa(R^u,R_{\min}^{u,+})$  where $R_{\min}^u\subset R^{u,\leq 1}$ is the smallest integrally closed subring containing $R^{u,<1}$.
  Then $\Pi$ gives rise to an action in the adic space $Y=\Spa(R^u,R_{\min}^{u,+})$ and one can take $\Pi_0$ as the stabilizer of the open subspace $U$.
  Note that in this case, the action of $\Pi_0$ on $R^{\prime,u}$ us uniquely determined by that on $R^u$.
  Thus, for $R\to R'$ general, it suffices to prove the statement locally in the analytic topology of $R'$, and so we can assume that $R'$ is a composite of finite \'etale maps  and $\dagger$-rational localizations of $R$.
  We already dealt with $\dagger$-rational localizations, and for finite \'etale maps we use the argument of \cite[Lemma 2.3]{ScholzeLubinTate} (for which the finite type assumption is not relevant thanks to \cite[Lemma 15.6]{scholze_etale_cohomology_of_diamonds}).

Finally for part (3), it suffices to note  that we have a pullback square 
\[
\begin{tikzcd}
\GSpec (\Q_p)/\Pi^{\rm cont}_0  \ar[r] \ar[d] &  \GSpec (\Q_p)/\Pi^{\sm}_0 \ar[d] \\
\GSpec (\Q_p)/\Pi^{\rm cont} \ar[r] & \GSpec (\Q_p)/\Pi^{\sm}
\end{tikzcd}
\]
which follows from the fact that $\Pi^{\rm cont}\to \Pi^{\sm}$ is an epimorphism of groups.
Indeed, if $Z$ is the pullback of \eqref{eq0o2awjioqwe}, we have a natural map $Z\to (B\Pi^{\rm cont})\times_{B\Pi^{\sm}} ( B\Pi_0^{\sm})= B\Pi^{\rm cont}_0$, that produces a $\Pi^{\rm cont}_0$-torsor $\widetilde{Z}$ over $Z$.
The previous map factors through 
\[
Z\to \GSpec (R)/\Pi^{\rm cont}_0\to B\Pi^{\rm cont}_0
\]
and therefore 
\[
\widetilde{Z}= Z\times_{\GSpec (R)/\Pi^{\rm cont}_0} \GSpec (R) =\GSpec (R')
\]
where the last equivalence follows from \eqref{eq0o2awjioqwe} and the fact that $R'$ is $\dagger$-formally \'etale over $ R$ so that we have a fiber product 
\[
\GSpec(R')=\GSpec (R) \times_{\Marc(R)^{\dR}/\Pi_0^{sm}} \Marc(R')^{\dR}/\Pi_{0}^{\sm}. 
\]
This endows $\GSpec(R')$ with a $\Pi^{\rm cont}_0$-action such that $\GSpec(R')/\Pi^{\rm cont}_0=Z$ as wanted. 
\end{proof}

We can now give the the proof of the first assertion in \Cref{TheoDecompletion}.

\begin{proposition}\label{PropositionPart1MainTheorem}
Let $X$ be a Berkovich space endowed with an action of $G^{\cont}$. Write $h\colon X/G^{\ob{cont}}\to |X|_{\Betti}/G^{\la}$. Then the map $h$ is $0$-proper. Denote $\s{O}^{\la}_{X}:= h_{*} \s{O}_{X/G^{\ob{cont}}}$, then the $*$-pullback induces a fully faithful embedding
\[
\ob{Mod}_{\s{O}^{\la}_{X}}(\ob{D}(|X|_{\Betti}/G^{\la})) \hookrightarrow \ob{D}(X/G^{\cont}). 
\] 
\end{proposition}
\begin{proof}
  Keep the notation as in the theorem.
  We have a map $h\colon X/G^{\rm cont}\to |X|_{\Betti}/G^{\la}$, we want to see that $h$ is $0$-proper. For that, by \Cref{Lemma-Etale-localization} and \Cref{Lemma-Prim-padic}, we can restrict the action from $G$ to a smaller compact open subgroup, and then localize to closed $\dagger$-rational subspaces of $|X|_{\Betti}$.
  Hence, we can assume without loss of generality that $X=\GSpec(R)$ is affine.
  In that case, $0$-properness of the map $h$ follows from the $0$-properness of the stacks $X/G^{\rm cont}$ and $|X|_{\Betti}/G^{\la}$, which follows from the fact that $X$ is affine, $|X|_{\Betti}$ is a finite dimensional compact Hausdorff space, and that the stacks $BG^{\rm cont}$ and $BG^{\la}$ are $0$-proper by \Cref{LemmaFirstPropertiesGcont} (3) and  \Cref{Lemma-Prim-padic} (1).
  It is then a formal consequence of the projection formula that if $h\colon Y\to X$ is a prim map between analytic stacks, then the pullback map gives rise to a fully faithful functor $\ob{Mod}_{h_*\s{O}_Y}(\ob{D}(X))\hookrightarrow \ob{D}(Y)$.
  \end{proof}

\subsection{Localizations of locally analytic vectors, and part (2) of \Cref{TheoDecompletion}}
\label{sec:part-2-crefth-part-2}

We continue with the proof of part (2) of the theorem (cf.\ \Cref{CorollaryPart2MainDecompletionTheorem}). For that, we need a preliminary discussion on decompletions of Gelfand rings by locally analytic vectors.

\begin{lemma}\label{LemmaGelfandGroupCohomology}
  Let $R$ be a Gelfand $\Q_p$-algebra endowed with an action of a compact $p$-adic Lie group $G$.
  Then the connective cover $R^{\mathrm{cla}}:=\tau_{\geq 0} R^{\la}$ of the locally analytic vectors of $R$ is a Gelfand $\Q_p$-algebra.
\end{lemma}
\begin{proof}
  We begin by proving that $R^{\mathrm{cla}}$ is bounded. For that, let $S$  be a light profinite set and consider a map $f\colon S\to R^{\mathrm{cla}}$.
  By adjunction this gives rise to a $G$-equivariant map $ \Q_{p,\sol}[ S]\to C^{\la}_G\otimes R$ where the source has the trivial $G$-action. Since $R$ and $C^{\la}_{G}$ are  bounded, after normalizing by a power of $p$ this map naturally extends to a $G$-equivariant map $\Q_{p,\sol}\langle \N[S] \rangle_{\leq 1}\to  C^{\la}_G\otimes R$ and by taking invariants and the connective cover we get a natural extension of $f$ to a map $\Q_{p,\sol}\langle \N[S]\rangle_{\leq 1}\to R^{\mathrm{cla}}$, proving the bounded condition. 

  Next, we prove that $R^{\mathrm{cla}}$ is Gelfand. By  \cite[Definition 3.1.1]{drff_stacks}, we want to show that $R^{\mathrm{cla},\leq 1}/ pR^{\mathrm{cla},\leq 1}$ is a discrete condensed set.
  For that, recall that $R^{\mathrm{cla},\leq 1}$ is the subspace of $R^{\mathrm{cla}}$ for which maps from light profinite sets $S\to R^{\mathrm{cla}}$ extend to the algebra $\Q_{p,\sol}\langle \N[S]\rangle_{\leq 1}$. One deduces that 
\[
R^{\mathrm{cla},\leq 1}= \tau_{\geq 0} \Gamma(G, (C^{\la}_{G}\otimes R )^{\leq 1}).
\]
Hence, we have a fiber sequence 
\[
\Gamma(G, (C^{\la}_{G}\otimes R )^{\leq 1})\xrightarrow{p} \Gamma(G,(C^{\la}_{G}\otimes R )^{\leq 1}) \to  \Gamma(G, (C^{\la}_{G}\otimes R )^{\leq 1}/p).
\]
Since $(C^{\la}_{G}\otimes R )^{\leq 1}/p$ is static and discrete, we deduce that $R^{\mathrm{cla},\leq 1}/p$ is a subspace of $H^0(G, (C^{\la}_{G}\otimes R )^{\leq 1}/p)$ and hence that it is discrete as desired.  
\end{proof}

\begin{definition}\label{Definition-Classically-decompleted}
\begin{enumerate}
\item Let $G$ be a compact $p$-adic Lie group and $R$  a Gelfand $\bb{Q}_p$-algebra endowed with an action of $G$. 
  We say that $R$ is \textit{classically decompleted} if the following properties hold: 
\begin{enumerate}

\item $R^{\la}$ is a connective $\Q_p$-algebra, and hence a Gelfand ring by \Cref{LemmaGelfandGroupCohomology}.

\item The map $R^{\la}\to R$ induces an isomorphism on uniform completions. 

\end{enumerate}

\item Let $X$ be a Berkovich space endowed with an action of $G^{\rm cont}$ for a $p$-adic Lie group $G$.
  We say that $X$ is \textit{locally classically decompleted} if there is a strict closed cover $X=\bigcup_i \GSpec( R_i)$ by  affinoid subspaces, with $\GSpec (R_i)$ stabilized by an open compact subgroup $G_i\subset G$, such that $R_i$ is classically decompleted as in (1).
  We let $X^{\la}$ denote the associated Gelfand stack obtained as the union $X^{\la}= \bigcup_i \GSpec (R_i^{\la})$, and say that $X$ \textit{admits locally a classical decompletion $X^{\la}$ with respect to $G$} . 
\end{enumerate}
\end{definition}

\begin{remark}
  \label{sec:local-locally-analyt-x-la-derived-berkovich-space}
  If as in \Cref{Definition-Classically-decompleted}, $X$ is locally classically decompleted, then $X^\la$ is in fact a Berkovich space.
  Indeed, as $X$ is glued along $\dagger$-rational localizations of the $\GSpec(R_i)$ and by assumption for each $i$, $\dagger$-rational localizations of $\GSpec(R_i)$ are pulled back from $\GSpec(R_i^\la)$ (a consequence of $R^u_i\cong R^{\la,u}_i$, cf.\ \cite[Lemma 3.1.8.]{drff_stacks}), the Gelfand stack $X^\la$ is glued from the $\GSpec(R_i^\la)$ along $\dagger$-rational localizations.
\end{remark}

\begin{remark}\label{RemarkClassicalDecompletionGelfandRingEtaleSites}
Let $G$ be a compact $p$-adic Lie group acting on a Gelfand $\bb{Q}_p$-algebra $R$. Suppose that $R$ is classically decompleted. Then, \cite[Lemma 3.2.8 and Proposition 3.2.11]{drff_stacks} yields natural equivalences between the analytic, finite \'etale, and Berkovich-\'etale sites of $R^{\la}$, $R$ and $R^u$.   Explicitly, this equivalence of sites sends an \'etale map $R^{\la}\to B$ to the base change along $R^{\la}\to R$, and its inverse sends a Berkovich-\'etale map $R\to R'$ (given as a composite of $\dagger$-rational localizations and finite \'etale maps) to the affinoid Gelfand stack  $Y$ given by the pullback 
\[
\begin{tikzcd}
Y \ar[r] \ar[d]& \GSpec (R')^{\dR} \ar[d] \\ 
\GSpec (R^{\la}) \ar[r] & \GSpec (R)^{\dR}.
\end{tikzcd}
\]
 \end{remark}

The following lemma shows that Berkovich-\'etale maps  naturally decomplete for a classically decompleted Gelfand ring. 

\begin{lemma}\label{LemmaLocallyAnalyticDecompletionSites}
Let $G$ be a compact $p$-adic Lie group and $R$ a $G$-equivariant Gelfand ring endowed with an action of $G^{\ob{cont}}$, and let $R\to R'$ be a Berkovich-\'etale  $G$-equivariant morphism. If $R$ is classically decompleted, then $R'$ admits a classical decompletion relative to $R$. More precisely, let $B$ be the Berkovich-\'etale $R^{\la}$-algebra associated to $R\to R'$ via \Cref{RemarkClassicalDecompletionGelfandRingEtaleSites}. Then there is a natural locally analytic action of $G$ on $B$  and a $G$-equivariant isomorphism $B\otimes_{R^{\la}} R= R'$. 
\end{lemma}
\begin{proof}

Consider the natural map of Gelfand stacks 
\[
\GSpec(R)/G^{\ob{cont}}\to \GSpec (R^{\la})/G^{\la}.
\]
As $R$ is classically decompleted, it induces an isomorphism on de Rham stacks as both de Rham stacks identify with $\GSpec(R)^{\dR}/G^{\sm}$.
Let $X$ denote the pullback 
\[
\GSpec (R^{\la})/G^{\la} \times_{\GSpec(R)^{\dR}/G^{\sm}} \GSpec(R')^{\dR}/G^{\sm}
\] induced by the $G$-equivariant \'etale map $R\to R'$. Then, $X$ fits in a commutative diagram with cartesian squares
\[
\begin{tikzcd}
X  \ar[r] \ar[d]& \GSpec(R')^{\dR}/G^{\la} \ar[r] \ar[d] & \GSpec(R')^{\dR}/G^{\sm} \ar[d] \\
\GSpec(R^{\la})/G^{\la} \ar[r] & \GSpec(R)^{\dR}/G^{\la}\ar[r] & \GSpec(R)^{\dR}/G^{\sm}.
\end{tikzcd}
\]
In particular, the pullback $X\times_{BG^{\la}} \GSpec(\Q_p)$ is equal to $\GSpec(R^{\la})\times_{\GSpec (R)^{\dR}} \GSpec(R')^{\dR}$ which, by \Cref{RemarkClassicalDecompletionGelfandRingEtaleSites}, is precisely $\GSpec(B)$.  The previous discussion produces a natural $G^{\la}$-action on the \'etale $R^{\la}$-algebra $B$ such that the natural equivalence $B\otimes_{R^{\la}} R\cong R'$ is $G$-equivariant. This finishes the proof of the lemma. 
\end{proof}

\begin{corollary}\label{CorollaryPart2MainDecompletionTheorem}
  Let $G$ be a compact $p$-adic Lie group and $X=\GSpec (R)$ a Berkovich space endowed with an action of $G^{\ob{cont}}$. Let $h\colon X/G^{\ob{cont}}\to |X|_{\Betti}/G^{\la}$ and write $\s{O}_{X}^{\la}:= h_* \s{O}_{X/G^{\ob{cont}}}$. Suppose that $X$ is  classically decompleted.
  Then the natural map 
\[
\underline{\GSpec}_{|X|_{\Betti}/G^{\la}}(\s{O}^{\la}_{X}) \xrightarrow{\sim}  \GSpec (R^{\la})/G^{\la}
\]
is an isomorphism.
\end{corollary}
\begin{proof}
By construction, $\s{O}^{\la}_{X}$ is the sheaf on the topological space $|X|$  sending a $\dagger$-rational localization $Z\subset X$ to $\s{O}(Z)^{\la}$. Thanks to \Cref{LemmaLocallyAnalyticDecompletionSites}, if $R\to R_Z$ is the $\dagger$-rational localization associated to $Z$, then $R^{\la}\to (R_Z)^{\la}$ is the $\dagger$-rational localization obtained via the equivalence of analytic sites between $R^{\la}$ and $R$. This proves that the sheaf on $|X|$ associated to $R^{\la}$ via the  localization along its Berkovich space is precisely $\s{O}^{\la}_{X}$. This finishes the proof of the corollary.  
\end{proof}

Berkovich spaces admitting classical decompletions by locally analytic vectors satisfy the following stability properties:

\begin{lemma}\label{Lemma-Gelfand-connective cover}
Let $G$ be a $p$-adic Lie group acting on a Berkovich space $X$.  The following hold:

\begin{enumerate}

\item Let $X=\bigcup_i Z_i'$ be a strict closed cover of $X$, and let $G_i\subset G$ be compact open subgroups acting on $Z_i$ respectively.
  Suppose that the $Z_i$ are locally classically decompleted (with respect to the $G_i$-action), then $X$ is locally classically decompleted.

\item Suppose that $X$  locally  classically decompleted.
  Let $X'\to X$ be a Berkovich (finite) \'etale map such that there is an open subgroup $G_0\subset G$ and an extension of $G_0$ to an action on $X'$ (this extension exists for affinoid $X'$ by \Cref{Lemma-Etale-localization}).
  Then $X'$ is locally  classically decompleted, moreover, the map $X^{\prime,\la}\to  X^{\la}$ is Berkovich (finite) \'etale and we have $X'=X^{\prime,\la}\times_{X^{\la}} X$.

\item Conversely, let $f\colon X'\to X$ be a Berkovich-\'etale map as in (2) with the extension of an action of $G_0\subset G$. Suppose that $f$ is a cover and that $X'$ is locally  classically decompleted, then $X$ is locally  classically decompleted. 

\item Suppose that $X$ is locally classically decompleted and let $Y\to X^{\la}$ be a map of $G^{\la}$-equivariant Berkovich spaces.
  Then the fiber product $X\times_{X^{\la}} Y$ is endowed with a natural action of $G^{\rm cont}$, is locally classically decompleted, and we have $(X\times_{X^{\la}} Y)^{\la}= Y$.

\item Let $H$ be a $p$-adic Lie group acting on a Berkovich space $Y$ which is locally  classically decompleted for $H$.
  If $X$ is locally  classically decompleted for $G$, then $X\times Y$ is locally  classically decompleted  for $G\times H$ with decompletion given by $(X\times Y)^{(G\times H)-\la}=X^{G-\la}\times Y^{H-\la}$.

\item Let $G=H\ltimes N$ be a semidirect product of $p$-adic Lie groups. Let $X$ be a Berkovich space endowed with an action of $G^{\cont}$. Suppose that $X$ admits locally a classical  decompletion $X^{N-\la}$ with respect to $N$. Then $X^{N-\la}$ carries a natural action of $H^{\cont}\ltimes N^{\la}$. Moreover, suppose that $X^{N-\la}$ admits locally a classical decompletion  $(X^{N-\la})^{H-\la}$ with respect to $H$. Then $X$ admits locally a classical decompletion $X^{G-\la}=(X^{N-\la})^{H-\la}$ with respect to $G$.

\end{enumerate}
\end{lemma}
\begin{proof}

\begin{enumerate}

\item This is obvious from the local nature of the definition of being  locally classically decompleted.

\item By passing to a strict $\dagger$-analytic cover,  we can assume that $G$ is compact, $X=\GSpec (R)$ and that $R$ is classically decompleted. In this case the statement follows from \Cref{LemmaLocallyAnalyticDecompletionSites}.

\item We can assume without loss of generality that $G=G_0$. 
Covering both $X'$ and $X$ by strict closed covers by affinoids, and using part (1), we can assume without loss of generality that both $X'$ and $X$ are affinoid subspaces. 
Moreover, by taking a further analytic localization of $X$, and an \'etale localization of $X'$ and part (2), we can assume without loss of generality that either $X'\to X$ arises from a strict closed cover or a finite \'etale cover.
 Since the case of strict closed covers is precisely part (1), it remains to prove part (3) of the lemma for finite \'etale covers. 
 Furthermore, by taking a further finite \'etale cover of $X'$, we can  assume without loss of generality that $X'\to X$ is an $H$-torsor for a finite group $H$, and that the action of $G$ on $X'$ commutes with that of $H$. 
 Let $R=\s{O}(X)$ and $R'=\s{O}(X')$,  by assumption $R^{\prime,\la}$ is connective and the natural map $R^{\prime,\la}\to R'$ is an isomorphism on uniform completions. 
 Since $X'\to X$ is an $H$-torsor,  we get that 
\[
R^{\la}= \Gamma(H, R^{\prime,\la})
\]
is the $H$-cohomology of a Gelfand ring over $\bb{Q}_p$, and since $H$ is finite, $R^{\la}$ is itself a Gelfand ring.  We claim that $R^{\la}$ and $R$ have the same uniform completion.  For that, we note that  if $A\to B$ is an $H$-torsor of Gelfand rings, then $A^u$ is equal to the $H$-invariants of $B^u$. Indeed, the normalized trace
\[
\frac{1}{\#(H)}\ob{tr}\colon B\to A
\] 
is a retract that sends  $B^{\leq 1}$ to $A^{\leq c}$ where $c= | \#(H)^{-1} |_p$. Thus, passing to uniform completions, the normalized trace on $B^u$ produces a retract $\frac{1}{\#(H)}\ob{tr}\colon B^u\to A^u$ which yields $A^{u}=B^{u,H}$.
  We deduce that 
\[
R^{\la,u}= R^{\prime,\la,u,H} = R^{\prime,u,H}= R^{u},
\]
proving that $R$ admits a classical decompletion. 

\item This follows formally from the projection formula of the functor of locally analytic vectors.

\item This follows from the fact that taking locally analytic vectors for $G\times H$ can be computed component-wise. More precisely,  from the following commutative diagram 
\[
\begin{tikzcd}
B(G^{\rm cont}\times H^{\rm cont}) \ar[r] \ar[d] & B(G^{\la}\times H^{\rm cont}) \ar[d] \\ 
B(G^{\rm cont}\times H^{\la}) \ar[r] & B(G^{\la}\times H^{\la})
\end{tikzcd}
\]
one has that at the level of rings that  $X^{\la}\times Y^{\la}$  is the decompletion of $X\times Y$ by locally analytic vectors for the action of $G\times H$. Since, locally in the analytic topology, the uniform completion of  $X^{\la}$ and $Y^{\la}$ is $X$ and $Y$ respectively, the uniform completion of $X^{\la}\times Y^{\la}$ is $X\times Y$.

\item  By part (1), we can assume that $G$ is compact and that $X$ is an affinoid Berkovich space. The statement is then a consequence of \Cref{LemmaDecompletionsSemidirectProducrs}. 
\qedhere 
\end{enumerate}
\end{proof}

\subsection{Generalities on decompletions by locally analytic vectors}
\label{sec:gener-lemma-decompl}

In this subsection we discuss a technical  decompletion criterion for locally analytic representations which generalizes \cite[Proposition 3.4.1]{jacinto2023solid}. This criterion, \Cref{PropLocAnCriterioV2}, will be a key technical input in the discussion of \Cref{TheoDecompletion} (3) in \Cref{sec:decompl-perf-compl} below. 

\Cref{PropLocAnCriterioV2} is an abstract version of the following observation (\cite[Lemma 2.1.5]{camargo2022geometricsen}): if $V$ is a $\Q_p$-Banach space with a continuous $G$-action with $G$ a compact $p$-adic Lie group, $V^0\subseteq V$ a $G$-stable $\Z_p$-lattice and the action of $G$ on $V^0/p$ factors over a \emph{finite} quotient of $G$, then $V$ is locally analytic.
If $G=\Z_p=\gamma^{\Z_p}$, this is elementary, cf.\ \cite[Example 2.1.9.]{MR4390302}: assume that $(\gamma-1)(V^0)\subseteq p^2V^0$, then $a\cdot v=\sum_{n=0}^\infty \binom{a}{n}(\gamma-1)^n v$ for $v\in V$ showing that the action is locally analytic.

 In the rest of the section we  follow the conventions of \Cref{sss:BinomReps}, that is, we will assume that $G$ is an uniform pro-$p$-group, and we fix a chart $\psi\colon G\cong \bb{Z}_p^d$ defined by elements $g_{1}, \ldots, g_d\in G$. The group $G$ is endowed with a (rational) $p$-valuation $w\colon G\backslash \{1\}\to (1/(p-1),\infty)$ and we set $w_i=w(g_i)$.
For $I\subset \{1,\ldots, d\}$ we  denote $w_I:=\sum_{i\in I} w_i$, and write $\omega:=\sum_{I,J\subset \{1,\ldots, d\}}|w_I-w_J|$.  
We let $\Lambda$ denote an uniform Banach $\bb{Q}_p$-algebra with a pseudo-uniformizer $\varpi$, and let $|-|_{\varpi}$ be the spectral norm defined by $\varpi$, so that in particular $|\varpi|_{\varpi}=1/2$. We suppose that $|p|_{\varpi}\leq 1/2$, and also that the element $\varpi^{\epsilon}\in \Lambda$ exists for $\epsilon>0$ a rational number.\footnote{This condition is added for pure convenience in the exposition, it would suffice to assume that $\varpi^{\epsilon}$ is defined for a dense subset of $(0,1)$.}   We highlight that since both $p$ and $\varpi$ are pseudo-uniformizers of $\Lambda$, the $\varpi$- and $p$-adic topologies on $\Lambda^{\circ}$ agree and $\Lambda^{\circ}[\frac{1}{p}]=\Lambda^{\circ}[\frac{1}{\varpi}]$.

\subsubsection{Decompletion under the condition $\epsilon$-\textbf{TORS}} We consider the $\varpi$-adically complete Hopf algebras $C_{r,\varpi}(G, \Lambda^{\circ})$ of \Cref{DefinitionAnalyticAlgebrasBinomial}, and let
 \[
   \f{G}_{\ob{Bin}_r,\varpi}:=\Spf ( C_{r,\varpi}(G, \Lambda^{\circ}) )
 \]
 denote their formal spectrum  with respect to the $\varpi$-adic filtration.  
  
   Before stating the main proposition, we make the following ad-hoc definition of integral locally analytic vectors.

\begin{definition}\label{DefinitionrIntegralAnalyticVectors}

Keep the previous notation.  Let $M$ be a $\Lambda^{\circ}$-linear solid representation of $G$. The \textit{integral $\ob{Bin}_{r,\varpi}$-analytic vectors of $M$} are defined as 
\begin{equation}\label{equationIntegralLocAnVectors}
M^{\ob{Bin}_{r,\varpi}}:=\Gamma(G, M\otimes_{\Lambda^{\circ}} C_{r,\varpi}(G, \Lambda^{\circ})), 
\end{equation}
where $G$ acts diagonally on $M$ and with respect to the left regular action of $C_{r,\varpi}(G, \Lambda^{\circ})$. 
\end{definition}

Consider $\f{G}_{\ob{Bin}_{r,\varpi}}$ as an analytic stack over $\AnSpec (\Lambda^{\circ})$ where $\Lambda^{\circ}$ is endowed with the induced solid analytic ring structure from $\bb{Z}_p$.
It is an affinoid group over  $\AnSpec (\Lambda^{\circ})$ with the induced analytic ring structure, and by \cite[Proposition 3.3.13 (2)]{RCBBCourse2026} the category $\ob{D}(B\f{G}_{\ob{Bin}_{r,\varpi}})$ is naturally equivalent to the category of right $\f{G}_{\ob{Bin}_{r,\varpi}}$-comodules on $\Lambda^{\circ}$-modules.
We have a morphism of Hopf algebras $C_{r,\varpi}(G, \Lambda^{\circ})\to  C(G, \Lambda^{\circ})$ which produces a natural map of group objects in analytic stacks $G^{\ob{cont}}_{\Lambda^{\circ}}\to \f{G}_{\ob{Bin}_{r,\varpi}}$, where  $G^{\ob{cont}}_{\Lambda^{\circ}}=\AnSpec (C(G, \Lambda^{\circ}))$.
The $*$-pullback produces natural symmetric monoidal functors
\[
\ob{D}(B\f{G}_{\ob{Bin}_{r,\varpi}})\to \ob{D}(BG^{\ob{cont}}_{\Lambda^{\circ}})\to \cat{Rep}^{\sol}_{\Lambda^{\circ}}(G).
\]

\begin{lemma}\label{LemmaRightAdjointIntegralrAnVectors}
The symmetric monoidal functor  
\[
F\colon \ob{D}(B \f{G}_{\ob{Bin}_{r,\varpi}} ) \to \cat{Rep}_{\Lambda^{\circ}}^{\sol}(G)
\]
admits  \eqref{equationIntegralLocAnVectors} as a right adjoint. Furthermore, the functor $(-)^{\ob{Bin}_{\varpi,r}}\colon \cat{Rep}_{\Lambda^{\circ}}^{\sol}(G)\to \ob{D}(B \f{G}_{\ob{Bin}_{r,\varpi}} )$ is $\ob{D}(B \f{G}_{\ob{Bin}_{r,\varpi}} )$-linear. 
\end{lemma}

\begin{proof} 
We use the notation of the proof of \Cref{LemmaContinuousRepsvsSolid}. Let $\ob{D}_{\proet}(-,\Lambda^{\circ})$ denote the $0$-\'etale  $6$-functor formalism of solid $\Lambda^{\circ}$-linear sheaves on the  pro\'etale  site of condensed anima, that is, the  $6$-functor formalism whose values in  a light condensed anima  $T$ is $\ob{D}_{\proet}(T,\Lambda)= \cat{CondAni}_{/T}\otimes_{\cat{CondAni}} \ob{D}(\Lambda^{\circ})$.  We let $S=[*]_{\ob{D}_{\proet}(-,\Lambda^{\circ})}$ be the Gestalt associated to the point $*$ with respect to $\ob{D}_{\proet}(-,\Lambda^{\circ})$, and for a condensed anima $T$ write $T^{\sol}=[T]_{\ob{D}_{\proet}(-,\Lambda^{\circ})}$. Following the proof of  \Cref{LemmaContinuousRepsvsSolid}, we have a natural morphism of $S$-pointed Gestalten
\[
g\colon B G^{\sol}\to B(G^{\ob{cont}}_S).
\]
  We note that, as $G$ is uniform pro-$p$, the trivial representation admits its Lazard resolution by finite free $\bb{Z}_{p,\sol}[G]$-modules, cf.\ \cite[Theorem 5.7]{jacinto2021solid}.  Thus, the exact same argument of the lemma shows that $g$ is a prim map of Gestalten and that $g_{1,*}1=1\in  \s{O}(B(G^{\ob{cont}}_S))_1= \ob{D}(BG^{\cont}_{\Lambda^{\circ}})$. We deduce that the symmetric monoidal  functor $g_{1}^*\colon \ob{D}(BG^{\cont}_{\Lambda^{\circ}})\to \cat{Rep}^{\sol}_{\Lambda^{\circ}}(G)$ is fully faithful and has a linear right adjoint. Thus, to show the linearity of the right adjoint of $F$, it suffices to show that the map of analytic stacks $h\colon BG^{\cont}_{\Lambda^{\circ}}\to B\f{G}_{\ob{Bin}_{r,\varpi}}$ is $0$-proper. This follows from the properness of   $BG^{\cont}_{\Lambda^{\circ}}\to \ob{AnSpec}(\Lambda^{\circ})$,  the $0$-properness of the diagonal of  $B\f{G}_{\ob{Bin}_{r,\varpi}}$ as  $\f{G}_{\ob{Bin}_{r,\varpi}}$ is affinoid with induced structure over $\ob{AnSpec}(\Lambda^{\circ})$, and \Cref{LemmaDiagonalPrimGestalten}.

Finally, it  is left to see that the right adjoint of $F$ is given by the formula \eqref{equationIntegralLocAnVectors}. For this, consider the prim map of Gestalten 
\[
f\colon BG^{\sol}_{S}\to B(\f{G}_{\ob{Bin}_{r,\varpi},S})
\]
where $\f{G}_{\ob{Bin}_{r,\varpi},S}$ is the base change to $S$ of the incarnation of $\f{G}_{\ob{Bin}_{r,\varpi}}$ as a Gestalt. Then the formula \eqref{equationIntegralLocAnVectors} for the right adjoint follows formally by proper base change along the map $S\to B(\f{G}_{\ob{Bin}_{r,\varpi},S})$. 
  \end{proof}

\begin{remark}\label{RemarkAdjointsAlgebras}
  Let $A$ be a solid $G$-equivariant commutative $\Lambda^{\circ}$-algebra. Then, by \Cref{LemmaRightAdjointIntegralrAnVectors}, $A^{\ob{Bin}_{r,\varpi}}$ has a natural structure as a commutative algebra in $\ob{D}(B\f{G}_{\ob{Bin}_{r},\varpi})$, possibly non-connective.
  Moreover, if $M$ is a $G$-equivariant $A$-module then $M^{\ob{Bin}_{r,\varpi}}$ carries a natural structure as an $A^{\ob{Bin}_{r,\varpi}}$-module.

\end{remark}

\begin{remark}
  \label{sec:gener-decompl-locall-t-structure}
  We note that, by flatness of $C_{r,\varpi}(G,\Lambda^\circ)\cong \widehat{\bigoplus}_{\N} \Lambda^\circ$ as a $\Lambda^\circ$-module,  $\ob{D}(B\f{G}_{\ob{Bin}_{r,\varpi}})$ has a natural $t$-structure such that the pullback to $\AnSpec(\Lambda^\circ)$ is $t$-exact, cf.\ \cite[Lemma 4.3.5.(3)]{jacinto2023solid}.
\end{remark}

\begin{remark}
  \label{sec:gener-decompl-locall-independent-version-on-r}
  We also set $\f{G}_{\ob{Bin}_{\varpi}}:=\varprojlim_{r>1} \f{G}_{\ob{Bin}_{r,\varpi}}$, with classifying stack $B\f{G}_{\ob{Bin}_{\varpi}}$, and
  $M^{\ob{Bin}_{\varpi}}:=\varinjlim_{r\to 1} M^{\ob{Bin}_{r,\varpi}}\in \ob{D}(B\f{G}_{\ob{Bin}_{\varpi}})$ for any $\Lambda^\circ$-linear solid representation of $G$. 
\end{remark}

We now state the main result of this section.
Given $a\in \Lambda$ and $M$ a $\Lambda^{\circ}$-module we let $M_a:=M/^{\bb{L}}a$ denote the quotient of $M$ by $a$.

\begin{proposition}\label{PropLocAnCriterioV2}

Let $G$ be a uniform pro-$p$-group   of dimension $d$ as before.    Let $A$ be a $G$-equivariant  $\varpi$-complete solid animated $\Lambda^{\circ}$-algebra. 
Let $\epsilon \in \bb{Q}_{>0}$  and suppose that the following condition holds: 
\begin{itemize}
\item[($\epsilon$-\textbf{TORS})]  There is a sequence of numbers $2>r_{0}>r_1>\cdots > 1$ with $r_n\to 1$ as $n\to \infty$ such that  for all $n\geq 0$ we have  $\pi_i(A^{\ob{Bin}_{r_n,\varpi}})$ is $\varpi^{\epsilon}$-torsion for $i\leq -1$.\footnote{The proof of \Cref{PropLocAnCriterioV2} only needs the existence of some $2>r>1$ satisfying this vanishing, but it will be convenient later to assume the existence of such sequence $r_0>r_1>\cdots$.}
\end{itemize}
Let $N\in \cat{Rep}^{\sol}_{A}(G)_{\geq 0}$ be a  connective $G$-equivariant $\varpi$-complete $A$-module.
Suppose that there is a  rational number $0\leq \delta< \epsilon$,    a (non-necessarily connective)  $ \Gamma(G, (\tau_{\geq 0}A^{\ob{Bin}_{r_0,\omega}})_{\varpi^{(d+1)\epsilon}})$-module $M$ endowed with the  trivial $G$-action, and $G$-equivariant maps of $ A_{\varpi^{(d+1)\epsilon}}$-modules
\[
N_{\varpi^{(d+1)\epsilon}} \to M\otimes_{ \Gamma(G, (\tau_{\geq 0}A^{\ob{Bin}_{r_0,\omega}})_{\varpi^{(d+1)\epsilon}}) }  A_{\varpi^{(d+1)\epsilon}} \to N_{\varpi^{(d+1)\epsilon}} 
\]
such that the composite is given by multiplication by $\varpi^{\delta}$.
Then, for any $r_n$ as in the condition $\epsilon$-\textbf{TORS},   $N[\frac{1}{\varpi}]$ admits a  decompletion by $\ob{Bin}_{r_n,\varpi}$-analytic vectors  as a $A[\frac{1}{\varpi}]$-module, that is, the natural map 
\begin{equation}\label{eqDecompletionAnalyticVectors}
N[\frac{1}{\varpi}]^{\ob{Bin}_{r_n,\varpi}}\otimes_{A[\frac{1}{\varpi}]^{\ob{Bin}_{r_n,\varpi}}} A[\frac{1}{\varpi}] \to N[\frac{1}{\varpi}]
\end{equation}
is an equivalence.  In particular,   $N[\frac{1}{\varpi}]$ admits a decompletion by locally analytic vectors.
\end{proposition}

\begin{remark}\label{Remark-Same-Psuedo-unif}
  Note that by the description of $\Lambda^{\circ}$, it is $(\varpi,p)$-complete and we have  $\Lambda^{\circ}[\frac{1}{\varpi}]=\Lambda^{\circ}[\frac{1}{p}]$, i.e.  both $p$ and $\varpi$ are pseudo-uniformizers.
  Thus, it is well defined to talk about the locally analytic vectors of $N[\frac{1}{\varpi}]=N[\frac{1}{p}]$. 
\end{remark}

\begin{proof}[Proof of \Cref{PropLocAnCriterioV2}] It suffices to prove the proposition for $r:=r_0$ as  $\mathrm{Bin}_{r,\varpi}$-analytic implies $\mathrm{Bin}_{r',\varpi}$-analytic for $r>r'$ and $\Lambda$-linear representations, cf. \Cref{LemmaGoodShierkPropertiesGBinom}.

Set $B:=\tau_{\geq 0}(A^{\ob{Bin}_{r,\varpi}})$, which is a connective $\varpi$-complete solid $G$-equivariant $\Lambda^\circ$-algebra.
  
We want to show that  the natural map  \eqref{eqDecompletionAnalyticVectors} is an equivalence. We first make  the following claim: 

\begin{claim*}
The multiplication map $\varpi^{d\epsilon}\colon A^{\ob{Bin}_{r,\varpi}}\to A^{\ob{Bin}_{r,\varpi}}$ factors through $B=\tau_{\geq 0}(A^{\ob{Bin}_{r,\varpi}})$ as a map of $B$-modules. 
\end{claim*}

Indeed, by definition we have  $A^{\ob{Bin}_{r,\varpi}}=\Gamma(G, A\otimes_{\Lambda^{\circ}} C_{r,\varpi}(G,\Lambda^{\circ}))$.
By Lazard's resolution, solid $G$-cohomology has cohomological dimension $\leq d=\ob{dim}(G)$, and this implies that  $A^{\ob{Bin}_{r,\varpi}}$ is a $(-d)$-connective ring.
The claim follows as each homotopy group $\pi_{i}(A^{\ob{Bin}_{r,\varpi}})$ for $i\in [-d,-1]$ is killed by $\varpi^{\epsilon}$, which implies that the cofibre of $B\to A^{\ob{Bin}_{r,\varpi}}$ is killed by $\varpi^{d\epsilon}$ as $G$-equivariant  $B$-module.

To  prove the proposition, it will suffice to prove the following fact:

\begin{claim*}
 The cofiber $Q$ of the map 
\begin{equation}\label{eqKeyMapOldProp521}
 N^{\ob{Bin}_{r,\varpi}}\otimes_{B} A \to N
\end{equation}
has homotopy groups of $\varpi^{\delta+d\epsilon}$-torsion. 
\end{claim*}

Indeed, if the claim holds, then $Q[\frac{1}{\varpi}]=0$.
By the condition ($\epsilon$-\textbf{TORS}), we have that $B[\frac{1}{\varpi}] = A^{\ob{Bin}_{r,\varpi}}[\frac{1}{\varpi}]$.
Thus, \eqref{eqDecompletionAnalyticVectors} becomes an isomorphism by passing to generic fibres on \eqref{eqKeyMapOldProp521}.

To show the claim, first note that $Q$ is a $\varpi$-complete module,  namely, $N$ is $\varpi$-complete, and the left term of \eqref{eqKeyMapOldProp521} is a tensor product $\varpi$-complete eventually connective modules over a connective ring and hence $\varpi$-complete thanks to \cite[Proposition 2.12.10]{mann2022p}.  By \Cref{LemmaTorsion} (applied with $f=\varpi^{\epsilon-\delta}$, $g=\varpi^{\delta+d\epsilon}$, and using $\delta<\epsilon$ to get $f$-completeness of $Q$), it suffices to show that $Q_{\varpi^{(d+1)\epsilon}}$ has $\varpi^{\delta+d\epsilon}$-torsion cohomology groups.
By assumption, $N_{\varpi^{(d+1)\epsilon}}$  is a $\varpi^{\delta}$-retract of $K:=M\otimes_{ \Gamma(G, B_{\varpi^{(d+1)\epsilon}}) }  A_{\varpi^{(d+1)\epsilon}} $.
Therefore it suffices to show that the cofibre of the map 
\[
\alpha\colon K^{\ob{Bin}_{r,\varpi}}\otimes_{B_{\varpi^{(d+1)\epsilon}}} A_{\varpi^{(d+1)\epsilon}} \to K
\]
has $\varpi^{d\epsilon}$-torsion homotopy groups. 

The projection formula for solid $G$-cohomology and the fact that the $G$-action on $M$ is trivial yield natural equivalences of $B_{\varpi^{(d+1)\epsilon}}$-modules
\[
\begin{aligned}
K^{\ob{Bin}_{r,\varpi}} & = \Gamma(G, K \otimes_{\Lambda^{\circ}} C_{r,\varpi}(G,\Lambda^{\circ})) \\ 
& =\Gamma\bigg(G, \big(M\otimes_{ \Gamma(G, B_{\varpi^{(d+1)\epsilon}})}  A_{\varpi^{(d+1)\epsilon}} \big) \otimes_{\Lambda^{\circ}} C_{r,\varpi}(G,\Lambda^{\circ})\bigg) \\ 
& =M\otimes_{ \Gamma(G,B_{\varpi^{(d+1)\epsilon}}) }   \Gamma(G, A_{\varpi^{(d+1)\epsilon}}\otimes_{\Lambda^{\circ}} C_{r,\varpi}(G,\Lambda^{\circ})) \\ 
& =M\otimes_{ \Gamma(G, B_{\varpi^{(d+1)\epsilon}}) }  A^{\ob{Bin}_{r,\varpi}}_{\varpi^{(d+1)\epsilon}}.
\end{aligned}
\]
Applying $-\otimes_{B_{\varpi^{(d+1)\epsilon}}}A_{\varpi^{(d+1)\epsilon}}$ yields the isomorphism of $B_{\varpi^{(d+1)\epsilon}}$-modules
\[
  K^{\ob{Bin}_{r,\omega}}\otimes_{B_{\varpi^{(d+1)\epsilon}}}A_{\varpi^{(d+1)\epsilon}} \cong K\otimes _{B_{\varpi^{(d+1)\epsilon}}}A_{\varpi^{(d+1)\epsilon}}^{\ob{Bin}_{r,\varpi}},
\]
so the cofibre of $\alpha$ identifies with the cofibre of the map
\[
  K\otimes_{B_{\varpi^{(d+1)\epsilon}}} A_{\varpi^{(d+1)\epsilon}}^{\ob{Bin}_{r,\varpi}}\to K
\]
induced by the $A_{\varpi^{(d+1)\epsilon}}$-linear, hence $A^{\ob{Bin}_{r,\varpi}}_{\varpi^{(d+1)\epsilon}}$-linear, structure of $K$.
This map admits the natural map
\[
 \beta\colon K=K\otimes_{B_{\varpi^{(d+1)\epsilon}}} B_{\varpi^{(d+1)\epsilon}}\to K\otimes_{B_{\varpi^{(d+1)\epsilon}}} A_{\varpi^{(d+1)\epsilon}}^{\ob{Bin}_{r,\varpi}}
\]
as a $B_{\varpi^{(d+1)\epsilon}}$-linear retract, so it suffices to show that the cofibre of $\beta$ is $\varpi^{d\epsilon}$-torsion. 

As was noted before, the cofibre of $B_{\varpi^{(d+1)\epsilon}}\to A_{\varpi^{(d+1)\epsilon}}^{\ob{Bin}_{r,\varpi}}$ is killed by $\varpi^{d\epsilon}$ as $B_{\varpi^{(d+1)\epsilon}}$-module, which then implies the same for the cofibre of $\beta$.
\end{proof}

We used the following lemma.

\begin{lemma}\label{LemmaTorsion}
  Let $R$ be an animated ring and $f,g\in R$.
  Let $M$ be a connective  $f$-adically   complete $R$-module. Then the  following are equivalent:
\begin{enumerate}

\item  The homotopy groups $\pi_i(M/fg)$ are $g$-torsion for $0\leq i\leq n$.

\item  The homotopy groups $\pi_i(M)$ are $g$-torsion for $0\leq i\leq n$.

\end{enumerate}  
  In particular,  the homotopy groups of $M/fg$ are $g$-torsion if and only if  the homotopy groups of  $M$ are $g$-torsion.
\end{lemma}

\begin{proof}
  Without loss of generality we can assume that $R=\bb{Z}[X,Y]$ and that $f=X$ and $Y=g$. We use the decalage functor from \cite{Bhatt2018}.   We have that $ L\eta_{g} (M/fg)= (L\eta_g M)/f$ by  \cite[Lemma 7.9]{Bhatt2018} (and the connectivity of $M$ is used to find a resolution by a connective complex of $g$-torsion free objects). 
  Since the operator $L\eta_g$ preserves $f$-adically complete modules (\cite[Lemma 6.19]{Bhatt2018}), $L\eta_g M$ is $f$-adically complete, and by Nakayama's lemma we know that $L\eta_g M$ is $(n+1)$-connective if and only if $L\eta_{g} (M/fg)$ is $(n+1)$-connective.  This translates in the equivalence of (1) and (2), proving what we wanted. 
\end{proof}

For \Cref{PropLocAnCriterioV2} to be useful, we must supply a good amount of  examples verifying the condition $\epsilon$-\textbf{TORS} that occur in practice.
We start with the case of the trivial action $A$.

\begin{example}[Case of trivial action]\label{ExampleTechnicalPropositionHoldsForTrivial}
  Let $2>r>1$.
  The ring $A=\Lambda^{\circ}$ satisfies the condition $\epsilon$-\textbf{TORS} for all $\epsilon>0$.
  Indeed, the  $\ob{Bin}_{r,\varpi}$-analytic vectors of $\Lambda^{\circ}$ are given by $\Gamma(G, C_{r,\varpi}(G,\Lambda^{\circ}))$ by definition. The fiber of the map $\Lambda^{\circ}\to \Gamma(G, C_{r,\varpi}(G,\Lambda^{\circ}))$ is represented by the Lazard complex $\ob{LZ}^{\bullet}(C_{r,\varpi}(G,\Lambda^{\circ}))$.   Thanks to \Cref{LemmaExtensionLazardResolutionDist}, multiplication by $\varpi^{\epsilon}$ on  $\ob{LZ}^{\bullet}(C_{r,\varpi}(G,\Lambda^{\circ}))$ is homotopic to zero for all $\epsilon$ such that  $2^{-\epsilon}< r^{-\omega}$. In particular, the higher cohomology groups of  $\Gamma(G, C_{r,\varpi}(G,\Lambda^{\circ}))$ are killed by $\varpi^{\epsilon}$.  Letting $r\to 1$, we can let $\epsilon\to 0$.
\end{example}

We deduce the following corollary from \Cref{PropLocAnCriterioV2} and \Cref{ExampleTechnicalPropositionHoldsForTrivial}, which is essentially the locally analytic criterion of \cite[Proposition 3.4.1]{jacinto2023solid}.

\begin{corollary}\label{LemmaCriterionv1}
Let $A$ be a connective $p$-complete solid $\Z_p$-algebra endowed with a trivial action of a compact  $p$-adic Lie group $G$. Let $P\in \cat{Rep}_{A}^{\sol}[G]$ be a connective and $p$-complete  $A$-linear solid $G$-representation such that the following holds: there is $s\geq 1$ and an open subgroup $G_0\subset G$ such that the action  of  $G_0$ on $P/p^{s}$ arises from the trivial representation. 
 Then $P[\frac{1}{p}]$ admits a decompletion by locally analytic vectors over $A[\frac{1}{p}]$, in other words, the natural map
\begin{equation}\label{eq0jpmpqwnoqwpkq3de}
(P[\frac{1}{p}])^{\la}\otimes_{(A[\frac{1}{p}])^{\la}} A[\frac{1}{p}]\to P[\frac{1}{p}]
\end{equation}
is an isomorphism. 
\end{corollary}
\begin{proof}
  By \Cref{Lemma-Decompletions} (1) and \Cref{ExampleTechnicalPropositionHoldsForTrivial}, we can assume without loss of generality that $A=\bb{Z}_p$.  By \Cref{Lemma-Decompletions} (3) it suffices to prove the claim after base change along $\bb{Z}_p\to \Lambda^{\circ}$, and  we can assume  that  $A=\Lambda^{\circ}$. We can shrink $G$  and assume it is uniform pro-$p$ and that $G=G_0$. In that case, $\tau_{\geq 0}(A^{\ob{Bin}_{r,\varpi}})= A$ as $A$ is static with trivial $G$-action, and we can apply \Cref{PropLocAnCriterioV2} with $N=P$ and $\delta=0$, and $\epsilon$ such that $\varpi^{(d+1)\epsilon}| p^s$ (which can be chosen arbitrarily small after enlarging $r$ thanks to \Cref{ExampleTechnicalPropositionHoldsForTrivial}). 
\end{proof}

 \subsubsection{Stability properties of the condition $\epsilon$-\textbf{TORS}} The following lemma extends \Cref{ExampleTechnicalPropositionHoldsForTrivial} to algebras which have trivial $G$-action on special fibre.

 \begin{lemma}
   \label{sec:stab-prop-cond-locally-analytic-satisfy-eps-tors}
   Let $A$ be a $\varpi$-complete $G$-equivariant $\Lambda^{\circ}$-algebra. Suppose that there is some $\delta>0$ such that the action of $G$ on $A_{\varpi^{\delta}}$ arises from a trivial representation.
   Then $A$  satisfies the condition $\epsilon$-\textbf{TORS} for all $\epsilon>0$. 
\end{lemma}
\begin{proof}
  Fix $\epsilon>0$.
  We claim that for all $2>r>1$ close enough to $1$ the ring $A^{\ob{Bin}_{r,\varpi}}$ has $\varpi^{\epsilon}$-torsion higher cohomology groups.
  We may assume $\epsilon<\delta$, as the condition $\epsilon$-\textbf{TORS} only gets stronger.
  Then by \Cref{LemmaTorsion}, it suffices to show this property after reduction modulo $\varpi^{\delta}$, in which case we have that 
\[
\begin{aligned}
A^{\ob{Bin}_{r,\varpi}}_{\varpi^{\delta}} & = \Gamma(G,   A_{\varpi^{\delta}}\otimes_{\Lambda^{\circ}}C_{r,\varpi}(G,\Lambda^{\circ})) \\ 
& = A_{\varpi^{\delta}}\otimes_{\Lambda^{\circ}} \Gamma(G, C_{r,\varpi}(G,\Lambda^{\circ})),
\end{aligned} 
\]
where the second equivalence holds by the projection formula (\Cref{sec:local-locally-analyt-nuclear-projection-formula-for-uniform-groups}) since $A_{\varpi^{\delta}}$ arises from a trivial $G$-representation.
The claim then follows from \Cref{ExampleTechnicalPropositionHoldsForTrivial} as the higher cohomology groups of $ \Gamma(G, C_{r,\varpi}(G,\Lambda^{\circ}))$ are concentrated in degrees $[-d,-1]$, and are $\varpi^{\epsilon}$-torsion for $r\to 1$. 
\end{proof}

The following lemma shows that the property $\epsilon$-\textbf{TORS} is transitive in a suitable sense. 

\begin{lemma}\label{LemmaTransitivePropertyTorsionCondition}
  Let $A$ be a $\varpi$-complete $G$-equivariant $\Lambda^{\circ}$-algebra satisfying condition $\epsilon$-\textbf{TORS}.
  Let $A\to B$ be a $G$-equivariant morphism of connective $\varpi$-complete algebras, and suppose that there is $\delta>  d\epsilon >0$, a connective $\tau_{\geq 0}\Gamma(G, A_{\varpi^{\delta}})$-module $M$ endowed with the trivial $G$-action, and an isomorphism of $G$-equivariant $A_{\varpi^{\delta}}$-modules
\begin{equation}\label{equationNaturalMapInvariants0j03j3}
M \otimes_{\tau_{\geq 0}\Gamma(G,A_{\varpi^{\delta}})} A_{\varpi^{\delta}}\cong  B_{\varpi^{\delta}}.
\end{equation}
 Then $B$ satisfies the condition $d\epsilon$-\textbf{TORS}. 
\end{lemma}
\begin{proof}
  Let $\delta>d\epsilon>0$ and let $r>1$ be such that   the higher cohomology groups of $A^{\ob{Bin}_{r,\varpi}}$ are $\varpi^{\epsilon}$-torsion.  We claim that the higher cohomology groups of  $B^{\ob{Bin}_{r,\varpi}}$ are $\varpi^{d\epsilon}$-torsion.
  By \Cref{LemmaTorsion} it suffices to prove this claim modulo $\varpi^{\delta}$.
  In that case, we have that 
\[
\begin{aligned}
B^{\ob{Bin}_{r,\varpi}}_{\varpi^{\delta}} & = \Gamma(G, B_{\varpi^{\delta}}\otimes_{\Lambda^{\circ}} C_{r,\varpi}(G,\Lambda^{\circ}))\\ 
& = M\otimes_{\tau_{\geq 0}\Gamma(G, A_{\varpi^{\delta}})} \Gamma(G, A_{\varpi^{\delta}}\otimes_{\Lambda^{\circ}} C_{r,\varpi}(G,\Lambda^{\circ})) \\ 
& =    M\otimes_{\tau_{\geq 0}\Gamma(G, A_{\varpi^{\delta}})} A^{\ob{Bin}_{r,\varpi}}_{\varpi^{\delta}}. 
\end{aligned}
\] 
Then, applying \Cref{LemmaTorsion}, we know that   the higher cohomology groups of $A_{\varpi^{\delta}}^{\ob{Bin}_{r,\varpi}}$ are $\varpi^{\epsilon}$-torsion.
This implies that the cofibre of $\tau_{\geq 0}(A_{\varpi^{\delta}}^{\ob{Bin}_{r,\varpi}})\to A_{\varpi^{\delta}}^{\ob{Bin}_{r,\varpi}}$ is killed by $\varpi^{d\epsilon}$ as a $\tau_{\geq 0} \Gamma(G,A_{\varpi^{\delta}})$-module (as this last algebra is connective).   Since $M\otimes_{\tau_{\geq 0} \Gamma(G, A_{\varpi^{\delta}})} \tau_{\geq 0}(A_{\varpi^{\delta}}^{\ob{Bin}_{r,\varpi}})$ remains connective, one deduces that the   higher cohomology groups of $B_{\varpi^{\delta}}^{\ob{Bin}_{r,\varpi}}$ are killed by $\varpi^{d\epsilon}$. This proves the lemma.
\end{proof}

From \Cref{LemmaTransitivePropertyTorsionCondition} we can deduce the following stability properties of the $\epsilon$-\textbf{TORS} condition. 

\begin{corollary}(Stability under relatively locally analytic base change)\label{CoroStabilityBaseChange}
  Consider morphisms of $G$-equivariant $\Lambda^{\circ}$-algebras $A\to B$ and $A\to C$.
  Suppose that there is $\delta> d\epsilon>0$ and a connective  $\tau_{\geq 0}\Gamma(G,A_{\varpi^{\delta}})$-module $M$ with trivial $G$-action  together with a $G$-equivariant isomorphism of $A_{\varpi^{\delta}}$-modules as in  \eqref{equationNaturalMapInvariants0j03j3}.
  Suppose that $C$ satisfies the condition $\epsilon$-\textbf{TORS}.
  Then $B\otimes_A C$ also satisfies  $d\epsilon$-\textbf{TORS}. 
\end{corollary}
\begin{proof}
This follows from \Cref{LemmaTransitivePropertyTorsionCondition} applied to the morphism $C\to B\otimes_A C$ as we can write 
\[
(B\otimes_A C)_{\varpi^{\delta}}= (M\otimes_{\tau_{\geq 0}\Gamma(G, A_{\varpi^{\delta}})} \tau_{\geq 0}\Gamma(G, C_{\varpi^{\delta}}))\otimes_{\tau_{\geq 0}\Gamma(G, C_{\varpi^{\delta}})} C_{\varpi^{\delta}}.\qedhere
\]
\end{proof}

In \Cref{ExampleTechnicalPropositionHoldsForTrivial} and \Cref{sec:stab-prop-cond-locally-analytic-satisfy-eps-tors} we saw examples of algebras that satisfy the condition $\epsilon$-\textbf{TORS} for all $\epsilon>0$.
In the applications of Sen theory this will not be quite the situation, instead we will have the following uniform control of $\epsilon$ as the subgroup $G_0\subset G$ shrinks:

\begin{definition}\label{DefinitionUniformlyEpsilonBla}
  Let $G$ be a uniform pro-$p$-group of dimension $d$ as before.
  Let $A$ be a $G$-equivariant $\varpi$-adically complete $\Lambda^{\circ}$-algebra and let $\epsilon>0$. We say that $A$ satisfies the condition \textit{$\epsilon$-\textbf{TORS} uniformly on $G$} if for a basis of compact open subgroups  $G_0\subset G$ the restriction of $A$ to a $G_0$-equivariant ring satisfies the property $\epsilon$-\textbf{TORS}. 
\end{definition}

\begin{remark}\label{RemarkStabilityUniformPropertyEpsilonTors}
  Let $A$, $B$ and $C$ be as in \Cref{CoroStabilityBaseChange}. Suppose in addition that $C$ satisfies the condition $\epsilon$-\textbf{TORS} uniformly on $G$.
  Then $B\otimes_A C$ satisfies the condition $d\epsilon$-\textbf{TORS} uniformly on $G$.
  In fact, if $G_0\subset G$ and $B_{\varpi^{\delta}}\cong M\otimes_{\tau_{\geq 0}\Gamma(G, A_{\varpi^{\delta}})} A_{\varpi^{\delta}}$, then, letting $N:= M\otimes_{\tau_{\geq 0}\Gamma(G, A_{\varpi^{\delta}})} \tau_{\geq 0}\Gamma(G_0, A_{\varpi^{\delta}}) $, we have that 
\[
B_{\varpi^{\delta}}\cong N\otimes_{\tau_{\geq 0}\Gamma(G_0, A_{\varpi^{\delta}})} A_{\varpi^{\delta}} 
\]
and \Cref{CoroStabilityBaseChange} applies.
\end{remark}

The next two lemmas  combined with \Cref{CoroStabilityBaseChange} imply that the  condition ``$\epsilon$-\textbf{TORS} uniformly on $G$''  is stable under integral models of Tate-rational localizations.

\begin{lemma}\label{LemmaEtaleLocalizationDecomplete}
  Let $A$ be a  $p$-complete solid animated ring endowed with a solid action of a uniform pro-$p$ group $G$.
  Let $s> k\geq 1$ be integers, and let $f\in A$ be an element such that its reduction modulo $p^s$ arises from a $G$-invariant element.
  More precisely, such that the element $f\in A/p^s$ arises from a $G$-equivariant map $(\Z/p^s)[T]\to A/p^s$ where $G$ acts trivially on $(\Z/p^s)[T]$. The following holds:

\begin{enumerate}

\item The algebra $A\langle \frac{f}{p^k} \rangle:=A\otimes_{\Z_p} \Z_p\langle T \rangle/(p^kT-f)$ has a natural  $G$-action making $A\to A\langle \frac{f}{p^k} \rangle$ a $G$-equivariant map.
  The same holds for the algebra $A\langle \frac{p^k}{f} \rangle:=A\otimes_{\Z_p} \Z_p\langle T\rangle/ (fT-p^k)$. 

\item Let $B$ denote $A\langle \frac{f}{p^k} \rangle$ or $A\langle \frac{p^k}{f} \rangle$ and suppose that $A$ is a $\Lambda^{\circ}$-algebra.
  Then $B$ satisfies the hypothesis of \Cref{LemmaTransitivePropertyTorsionCondition}.
  In particular, $B$ satisfies the condition $d\epsilon$-\textbf{TORS}  uniformly on $G$ if $A$ does so (cf. \Cref{RemarkStabilityUniformPropertyEpsilonTors}). 

\end{enumerate}

\end{lemma}
\begin{proof}
The solid algebra $A\langle \frac{f}{p^k} \rangle$ is the solidification of the sheafification of the presheaf sending a light profinite set $S$ to the classical $p$-complete algebra $A(S)\widehat{\otimes}_{\Z_p} \Z_p\langle T \rangle/ (p^kT-f)$, similarly for $A\langle \frac{p^k}{f} \rangle$. Therefore, in order to construct the $G$-action on $A\langle \frac{p^k}{f}\rangle$, it suffices to construction  a, functorial in $S$, action of the group $G(S)$ on  $A(S)\widehat{\otimes}_{\Z_p} \Z_p\langle T \rangle/ (p^kT-f)$ for all $S$ (similarly for $A\langle \frac{p^k}{f}\rangle$). 

Hence, it suffices to consider the following abstract situation: let $H$ be an abstract group and let $C$ be a $p$-complete $H$-equivariant animated ring endowed with an element $f\in C$ and a $H$-equivariant map $(\Z/p^s)[T]\to C/p^s$ where $H$ acts trivially on $(\Z/p^s)[T]$. We want to construct a natural action (both in $H$ and $C$) of $H$ on the rings $C\langle \frac{f}{p^k} \rangle$ and $C\langle \frac{p^k}{f} \rangle$. To prove this, let $X\to BH$ be the $p$-complete formal algebraic stack given by the relative spectrum of $C$ over $BH$, where $BH=(\Spf \Z_p)/H$. The datum of $f$ and the $H$-equivariant map from $\Z/p^s[T]$ to $C/p^s$ is equivalent to a map of $p$-complete algebraic stacks over  $BH$
\[
(f\mod p^s)\colon X\to BH \times \mathbb{A}^{1}_{\Z_p}/ p^s\mathbb{A}^{1}_{\Z_p}. 
\]
and where the affine line $\A^1_{\Z_p}$ is taken to be $p$-complete, i.e., to be the formal spectrum of $\Z_p\langle T\rangle=\Z_p[T]^{\wedge}_p$. 

\begin{itemize}

\item[Case 1.] We do the construction for $C\langle \frac{f}{p^k} \rangle$. For that, consider the pullback of algebra stacks over $BH$
\[
\begin{tikzcd}
X'  \ar[r] \ar[d] & BH\times p^k \A^1_{\Z_p}/p^s\A^1_{\Z_p} \ar[d] \\
X \ar[r] & BH\times \A^1_{\Z_p}/p^s\A^1_{\Z_p}. 
\end{tikzcd}
\]
Then $X'\to BH$ is represented by an affine scheme and its pullback to $\Spf (\Z_p)$ is the algebra $C\langle \frac{f}{p^k} \rangle$.
This produces the desired $H$-action, it is clear from the construction that it is functorial on $H$ and $C$.

\item[Case 2.] Next, we do the construction for $C\langle \frac{p^k}{f} \rangle$. For that, let $\bb{D}(|p^k|, 1)$ denote the rigid analytic annulus of radius $[|p^k|,1]$ given by the Banach algebra $\Q_p\langle T, \frac{p^k}{T} \rangle$.
  This is an open adic subspace of the analytic affine line over $\Q_p$, and it is stable under the action of the disc $\bb{D}(|p^s|)$ by addition (because $s>k$).
  Passing to power-bounded functions, the formal scheme $\widetilde{\mathbb{D}}(|p^k|,1):=\Spf( \Z_p\langle T,\frac{p^k}{T} \rangle)$ carries a natural action of the scaled affine line $p^s\mathbb{A}^1_{\Z_p}=\Spf ( \Z_p\langle \frac{T}{p^s} \rangle )$.
  Thus, the natural projection $\widetilde{\mathbb{D}}(|p^k|,1)\to \mathbb{A}^1_{\Z_p}$ is $p^s\mathbb{A}^1_{\Z_p}$-equivariant for the addition.
  The construction of the $H$-action follows in this case analogously by using the pullback diagram 
\[
\begin{tikzcd}
X' \ar[d]  \ar[r]& BH\times  \widetilde{\mathbb{D}}(|p^k|,1)/p^s \mathbb{A}^1_{\Z_p} \ar[d] \\ 
X \ar[r] &  BH\times  \mathbb{A}^1_{\Z_p}/p^s\mathbb{A}^1_{\Z_p}
\end{tikzcd}
\]
Again, the functoriality in $H$ and $C$ is clear. 
\end{itemize}

Finally, it is clear by construction that  $A\langle \frac{f}{p^k}\rangle$ and $A\langle \frac{p^k}{f} \rangle$ satisfy the hypothesis of \Cref{LemmaTransitivePropertyTorsionCondition}. That is, we have 
\[
  A\langle \frac{f}{p^k}\rangle/p^s= A/p^s\otimes_{ (\bb{Z}/p^s)[T]} (\bb{Z}/p^s)[T,S]/(p^kS-T),
\]
which is isomorphic to $A/p^s\otimes_{\tau_{\geq 0}\Gamma(G,A_{p^s})}(\tau_{\geq 0}\Gamma(G,A_{p^s})\otimes_{ (\bb{Z}/p^s)[T]} (\bb{Z}/p^s)[T,S]/(p^kS-T))$,
and similarly
\[
A\langle \frac{p^k}{f}\rangle/p^s= A/p^s\otimes_{ (\bb{Z}/p^s)[T]} (\bb{Z}/p^s)[T,S]/(TS-p^k),
\]
where in both cases $T\mapsto f$ in $A$.
\end{proof}

\begin{lemma}\label{Lemma-Colimit-Solid}
  Let $R$ be a Gelfand $\Q_p$-algebra endowed with a solid action of a light profinite group $\Pi$ with finite $p$-cohomological dimension.
  Let $r_1,\ldots, r_s\in R^{\circ}$ be power bounded elements, see \cite[Definition 2.1.1]{drff_stacks}. The following hold: 
\begin{enumerate}

\item  We can write  $R^{\circ}=\varinjlim_{i} A_i$ as a filtered colimit of $\Pi$-equivariant $p$-complete solid animated rings. 

\item Let $R^{\circ}=\varinjlim_i A_i$ be a presentation of $R^{\circ}$ as a filtered colimit of $\Pi$-equivariant $p$-complete connective rings.   Let $k\geq 1$ and let $\Pi_0\subset \Pi$ be an open subgroup such that  $\Pi_0$ acts trivially on the reduction  modulo $p^k$ of  $r_1,\ldots, r_s$. Then there is $i$ large enough such that the elements $r_1,\ldots, r_s$ lift to elements $a_1,\ldots, a_s \in A_i$  such that their mod $p^k$ fiber arise from a $\Pi_0$-equivariant map 
\[
(\Z/p^k)[T_1,\ldots, T_s]\to A_i/p^k,
\]
where the left term has the trivial $\Pi_0$-action. 
\end{enumerate}
\end{lemma}
\begin{proof}
  We first  prove part (1).
  The rings $\Z_{p,\sol}[ \N[S] ]:= \Sym_{\Z_{\sol}} \Z_{\sol}[S]$ with $S$ a light profinite set form a family of compact projective generators of the category of solid animated $\Z_p$-algebras.
  Since $\Z_{p,\sol}[G\times S]\cong \Z_{p,\sol}[G]\otimes_{\Z_{p,\sol}} \Z_{p,\sol}[S]$ is the free $G$-equivariant solid $\Z_p$-module generated by $S$, we deduce that $\Z_{p,\sol}[\N[G\times S]]$ is a family of compact projective  $G$-equivariant solid $\Z_p$-algebras.
  Thus, we can write $R^{\circ}=\varinjlim_{I} P_i$ as a sifted colimit of $G$-equivariant animated rings with $P_i$ of the form $\Z_{p,\sol}[\N[G\times S]]$.  By the definition of $R^{\circ}$, any map $\Z_{p,\sol}[\N[G\times S]]\to R^{\circ}$ will factor uniquely through the $p$-completion $\Z_{p,\sol}\langle \N[G\times S] \rangle \to R^{\circ}$.
  Therefore, as the map $\Z_{p,\sol}[\N[G\times S]]\to \Z_{p,\sol}\langle\N[G\times S]\rangle $  is idempotent by \cite[Lemma 2.4.7]{RCAndeRham},  the diagram $\{P_j\to R^{\circ}\}_{J}$ will factor through the $p$-completion $\{P_{j,p}^{\wedge}\to R^{\circ}\}_{J}$ and the natural map $\varinjlim_{J}P_{j,p}^{\wedge}\to R^{\circ}$ is an equivalence. Indeed, $\varinjlim_{J}P_{j,p}^{\wedge}$ is an idempotent $R^{\circ}$-algebra being a colimit of idempotent morphisms, so that the $R^{\circ}$-linear map  $\varinjlim_{J}P_{j,p}^{\wedge}\to R^{\circ}$ must be an equivalence.

  Writing $J$ as a filtered colimit of finite subdiagrams, we have the desired expression $R^{\circ}=\varinjlim_{i} A_i$ as a filtered colimit of $G$-equivariant $p$-complete animated solid $\Z_p$-algebras.  

Next we prove part (2). Since $R$ is Gelfand, $R^{\circ}/p^k$ is a discrete static ring.  In particular, there is an open subgroup $\Pi_0$ leaving stable the reduction of the elements $r_1,\ldots, r_s\in R^{\circ}/p^k$. Thus, we have a morphism of $\Pi_0$-equivariant animated rings 
\begin{equation}\label{eao9qhjn3kllawkdfa}
\Z/p^k[T_1,\ldots, T_s]\to R^{\circ}/p^k. 
\end{equation}

By assumption,  we have   $R^{\circ}/p^{k}=\varinjlim_{I} A_i/p^k$. Now, since $\Pi_0$ has finite $p$-cohomological dimension, $\Z/p^k[T]$ is a compact $\Pi_0$-equivariant animated ring since it correpresents the functor sending a $\Pi_0$-equivariant ring $B$ to the connective cover of its invariants $\tau_{\geq 0} \Gamma(\Pi_0, B)$ which commutes with filtered  colimits. Thus, there is some index $i\in I$ such that the map \eqref{eao9qhjn3kllawkdfa} factors through a $\Pi_0$-equivariant map $\Z/p^k[T_1,\ldots, T_s]\to A_i/p^k$ proving what we wanted. 
\end{proof}

\subsection{Decompletions of perfect complexes, and part (3) of \Cref{TheoDecompletion}}
\label{sec:decompl-perf-compl}

In this section we prove part  (3) of \Cref{TheoDecompletion}, i.e., that the functor $\mathrm{Mod}_{\s{O}^\la_X}\ob{D}(|X|_\Betti/G^\la)\to \ob{D}(X/G^{\rm cont})$ induces an equivalence of perfect complexes under a suitable technical condition on $X$. The main results are \Cref{TheoremDecompletionPerfectComplexesII} and \Cref{sec:decompl-perf-compl-1-decompletion-etale-coverings}.

Since the statement is local for the analytic topology of $X$ (by descent), and only depends on open compact subgroups $G_0\subset G$,   we can assume without loss of generality that  $X=\GSpec (R)$ is affinoid and that $G$ is a uniform pro-$p$-group as in the conventions of \Cref{sec:gener-lemma-decompl}, with an action of $G^\cont$ on $\GSpec(R)$.
In particular, the $\Z_{p,\sol}[G]$-module $R$ arises from a $C(G,\Z_p)$-comodule.
We note that $R^\circ\subseteq R$ is a sub $C(G,\Z_p)$-comodule.\footnote{This is a condition, which can be checked on the uniform completion, but $R^u$ is a classical Banach algebra, where the statement is true.}

The following lemma shows that it suffices to prove the descent of perfect complexes locally in the target.

\begin{lemma}\label{LemmaDescentCoefficientsPerfect}
Consider a cartesian square of Gelfand stacks 
\[
\begin{tikzcd}
Y \ar[r,"f^\prime"]  \ar[d,"g^\prime"] & Y' \ar[d,"g"] \\ 
X\ar[r,"f"] & X' 
\end{tikzcd}
\]
Suppose that the following conditions hold:
\begin{enumerate}

\item  $f$ is prim and $1\xrightarrow{\sim} f_* 1$ is an equivalence. 

\item  The functor $f^{\prime,*}\colon \Perf(Y')\to \Perf(Y)$ is an equivalence. 

\item The map $g$ is an epimorphism of Gelfand stacks.

\end{enumerate}

Then the pullback $f^{*}\colon \Perf(X')\to \Perf(X)$ is an equivalence. 
\end{lemma}

\begin{proof}
This follows formally by proper base change, namely, condition (1) implies that $1\xrightarrow{\sim} f_{*}^{\prime} 1$ is an equivalence, and hence both $f^*$ and $f^{\prime,*}$ are fully faithful. It is left to check that any perfect complex over $X$ arises as pullback of a perfect complex from $X'$. Let $P\in \Perf(X)$, since $f$ is a $!$-cover of Gelfand stacks (it is prim and descendable), it suffices to show that the natural map $f^*f_*P\to P$ is an isomorphism. Since by condition (3), $g$ is an epimorphism of Gelfand stacks, then so is $g^\prime$, and we can prove this fact after base change along $g^\prime$. By proper base change it suffices to show that the natural map 
\[
f^{\prime,\ast}f'_*g^{\prime,\ast} P \to g^{\prime,\ast} P
\]
is an equivalence. This follows from condition (2). 
\end{proof}

\begin{lemma}\label{CorollaryLocAnDecompletion}
Let $G$ be a $p$-adic Lie group and $X$ a Berkovich space endowed with a $G^{\ob{cont}}$-action.  Suppose that we are given with a $G^{\ob{cont}}$-equivariant  morphism of Berkovich spaces $ X\to X'$  which is, locally in the analytic topology on $X'$, represented by affinoid Gelfand stacks with same uniform completion.  Suppose that the action of $G$ on $X'$ is locally analytic. Let $G_0\subset G$ be an open  subgroup,  and let $Y'\to X'$ be a $G^{\la}_0$-equivariant morphism of Berkovich spaces with pullback $Y:= Y'\times_{X'} X$. Then the following holds:

\begin{enumerate}

\item  The following commutative square is Cartesian
\[
\begin{tikzcd}
Y/G^{\cont}_0 \ar[d,"g'"] \ar[r,"f'"]& Y'/G^{\la}_0 \ar[d,"g"] \\ 
X/G^{\cont} \ar[r,"f"] & X'/ G^{\la}.
\end{tikzcd}
\]
Furthermore, if $X$ admits locally a classical decompletion given by $X'$, then $Y$ admits locally a classical decompletion given by $Y'$.  Conversely, if $Y'\to X'$ is an epimorphism and $Y$ admits locally a classical decompletion given by $Y'$, then $X$ admits locally a classical decompletion given by $X'$. 

\item  Suppose that $X$ admits $X'$ as the locally analytic decompletion, that    $Y'\to X'$ is  an epimorphism of Gelfand stacks, and that the pullback $f^{\prime,*}\colon \Perf(Y'/G^{\la}_0)\xrightarrow{\sim} \Perf(Y/G^{\cont}_0)$ is an equivalence. Then the pullback $f^*\colon \Perf( X'/ G^{\la})\xrightarrow{\sim} \Perf(X/G^{\cont})$ is an equivalence. 

\end{enumerate}
 
\end{lemma}
\begin{proof}
\begin{enumerate}
\item  We have a commutative diagram with left and right cartesian squares
\[
\begin{tikzcd}
Y/G^{\cont}_0 \ar[r]\ar[d]    & Y'/G^{\cont}_0 \ar[r] \ar[d] & Y'/G^{\la}_0 \ar[d] \\
 X/G^{\cont} \ar[r] & X'/G^{\cont} \ar[r] & X'/G^{\la}
\end{tikzcd}
\] 
This implies that the outer square is cartesian by the pasting lemma. Suppose that $X'$ is the locally analytic decompletion of $X$, the fact that $Y'$ is the locally analytic decompletion  of $Y$ follows from \Cref{Lemma-Gelfand-connective cover} (4). Conversely, suppose that $Y'\to X'$ is an epimorphism and that $Y$ admits a locally analytic decompletion given by $Y'$. We first claim that the map $X\to X'$ is $0$-proper, namely, this follows from the fact that they have the same underlying topological space and that the map $X\to |X|_{\Betti}$ is $0$-proper by \Cref{sec:deriv-berk-spac-berkovich-spaces-as-relative-spectrum}. Hence, we see that $f_{*} \s{O}_X = \s{O}_X^{\la}$ is the sheaf of locally analytic vectors of $\s{O}_X$. Thus, to show that $X$ admits a locally  analytic decompletion given by $X'$, it suffices to show that the natural map $\s{O}_{X'}\to \s{O}_X^{\la}$ is an isomorphism (as the equivalence on uniform completions, locally on the analytic topology, holds by assumption).  Since the map $Y'\to X'$ is an epimorphism, we can prove this fact after base change along the map  $g$  of part (1), but then the claim boils down the the natural equivalence $\s{O}_{Y'}\xrightarrow{\sim} \s{O}_{Y}^{\la}$ which holds since $Y'$ is the locally analytic decompletion of $Y$.

\item  Part (2) follows from \Cref{LemmaDescentCoefficientsPerfect}, noticing  that the morphism $f$ is $0$-proper with $1\cong f_\ast 1$ by \Cref{TheoDecompletion} (1).  \qedhere
  
  \end{enumerate}
\end{proof}

Let $\Lambda$ be a uniform Banach $\bb{Q}_p$ algebra with pseudo-uniformizer $\varpi$. Let $|-|_{\varpi}$ be the spectral norm defined by $\varpi$ so that $|\varpi|_{\varpi}=1/2$. We suppose that $|p|_{\varpi}\leq 1/2$ and that for $\epsilon>0$ a rational number $\varpi^{\epsilon}$ exits in $\Lambda$  (see the conventions of \Cref{sec:gener-lemma-decompl}).
From now on we will suppose that $R$ is a $\Lambda$-algebra and that the action of $G$ is $\Lambda$-linear.

Before proving the theorem, we need the following technical fact.
Almost mathematics is always referring to almost mathematics with respect to the ideal $\bigcup_{\delta>0}\varpi^\delta\Lambda^\circ=\Lambda^{\circ\circ}\subseteq \Lambda^\circ$.

\begin{lemma}\label{LemmaApproximationLocAnvsCont}
  \begin{enumerate}
  \item The natural map of $\Lambda^{\circ}$-algebras  $\varinjlim_{r} C_{r,\varpi}(G, \Lambda^{\circ}) \to C(G,\Lambda^{\circ})$ is an almost isomorphism  after $\varpi$-completion.
  \item If $M$ is a $C(G,\Lambda^\circ)$-comodule, then the natural map of $C(G,\Lambda^\circ)$-comodules
\[
M^{\ob{Bin}_{\varpi}}=\varinjlim_{r\to 1} M^{\ob{Bin}_{r,\varpi}}\to M
\]
is an almost isomorphism after $\varpi$-completion.
\item If $M$ is a connective $C(G,\Lambda^\circ)$-comodule, there are $\epsilon>0$ and a sequence $r_0>r_1>r_2>\cdots 1$  converging to $1$ such that $\pi_{i}(M^{\ob{Bin}_{r_n,\varpi}})$ is $\varpi^\epsilon$-torsion for any $n\geq 0$ and $i\leq -1$, then the natural map
  \[
    \varinjlim_{r\to 1} \tau_{\geq 0}M^{\ob{Bin}_{r,\varpi}}\to M
  \]
  is an almost isomorphism after $\varpi$-completion.
  \end{enumerate}

\end{lemma}
\begin{proof}
  If $r=|\varpi^{\epsilon}|$ with $\epsilon>0$, then the inclusion $C_{r,\varpi}(G,\Lambda^0)\to C(G,\Lambda^\circ)=C(G,\Z_p)\otimes_{\Z_p}\Lambda^{\circ}$ identifies with the injective map
  \[
    \widehat{\bigoplus_{\alpha\in \N^d}} \Lambda^\circ |\varpi^\epsilon|^{\tau \alpha} \binom{\underline{x}}{\alpha}\to \widehat{\bigoplus_{\alpha\in \N^d}} \Lambda^\circ \binom{\underline{x}}{\alpha} 
  \]
  on Mahler expansions.
  Taking the colimit over $r$ and modding out $\varpi$ this map becomes the direct sum over $\alpha\in \N^d$ of the almost isomorphisms
  \[
    \varinjlim\limits_{\epsilon\to 0}\Lambda^\circ_\varpi |\varpi^{\epsilon}|^{\tau \alpha}\binom{\underline{x}}{\alpha}\to \Lambda^\circ_\varpi \binom{\underline{x}}{\alpha}
  \]
  with $\Lambda^\circ_\varpi=\Lambda^\circ/^{\mathbb{L}} \varpi$.
  
  This implies the first claim.
  The second follows from (1) by taking $G$-invariants on $M\otimes_{\Lambda^\circ} \varinjlim_{r} C_{r,\varpi}(G, \Lambda^{\circ}) \to M\otimes_{\Lambda^\circ}C(G,\Lambda^{\circ})$, using that $\Gamma(G,M\otimes_{\Lambda^\circ} C(G,\Lambda^\circ))\cong M $.
  This implies that the map
  \[
    \varinjlim_{r\to 1} M^{\ob{Bin}_{r,\varpi}}\to M
  \]
  is an almost isomorphism after $\varpi$-completion.

  For the third one, set $K:=M^{\ob{Bin}_\varpi}=\varinjlim_{r\to 1} M^{\ob{Bin}_{r,\varpi}}$, so that $\tau_{\geq 0}K=\varinjlim_{r\to 1} \tau_{\geq 0}M^{\ob{Bin},{r,\varpi}}$.
  Consider the fiber sequence
  \[
    (\tau_{\geq 0} K)^\wedge_\varpi\overset{\alpha}{\to} K^\wedge_\varpi\to (\tau_{\leq -1}K)^{\wedge}_\varpi,
  \]
  where $(-)^\wedge_\varpi$ denotes $\varpi$-adic completion.
  We know that $K^\wedge_\varpi$ is almost isomorphic to $M$ by (2), and thus almost connective. We also know that $(\tau_{\geq 0} K)^\wedge_\varpi$ is connective, which implies that  $(\tau_{\leq -1}K)^{\wedge}_\varpi$ is almost connective. 
  In particular, it suffices to show that $\tau_{\leq -1} K$ is already $\varpi$-adically complete for a cofinal sequence of $r$'s as then it is $(-1)$-coconnective and thus almost zero. 
  But the assumption implies that each homotopy group of $\tau_{\leq -1}K$ is $\varpi^\epsilon$-torsion for $r_n\to 1$, and thus $\varpi$-adically complete.
\end{proof}

In the following we call animated algebras in $\ob{D}(\AnSpec(\Lambda^\circ)/G^\cont_{\Lambda^\circ})$ simply $G^\cont_{\Lambda^\circ}$-equivariant animated algebras. 

\begin{proposition}
  \label{sec:decompl-perf-compl-1-decompletion-for-pi-complete}
  Let $A$ be a $G^\cont_{\Lambda^\circ}$-equivariant animated algebra.
  Assume
  \begin{enumerate}
  \item $A$ is $\varpi$-adically complete and $A/\varpi$ is discrete, and 
  \item $A$ satisfies the condition $\epsilon$-\textbf{TORS} uniformly on $G$.
  \end{enumerate}
  Let $P$ be a $G$-equivariant perfect module over $A[1/\varpi]$.
  Then $P$ admits a decompletion, i.e., the natural map
  \[
    P^\la\otimes_{A[1/\varpi]^\la}A[1/\varpi]\to P
  \]
  is an equivalence.
\end{proposition}

Condition (1) implies that $A[1/\varpi]=A[1/p]$ is Gelfand.

\begin{proof}
  We note that the $G$-action on $A$ is assumed to be induced by a $C(G,\Lambda^\circ)$-coaction, however, we will view the object $P$ as a $\Lambda_{A[1/\varpi]}(G)$-module, i.e., as a solid $G$-equivariant $A[1/\varpi]$-module.
  The proof will rest on various claims.
  \begin{claim*}[A] \label{sec:decompl-perf-compl-1-claim-a}
    There exists a compact $\Lambda_{A[1/\varpi]}(G)$-module $C$ and a  trace class endomorphism $f\colon C\to C$ in $\cat{Rep}^\sol_{A[1/\varpi]}(G)$ such that $P$ is a retract of $\mathrm{cofib}(C\overset{f-\mathrm{Id}}{\to}C)$.
  \end{claim*}
  Indeed, we note that each solid $G$-equivariant $A[1/\varpi]$-module that is compact as an $A[1/\varpi]$-module is compact as a $\Lambda_{A[1/\varpi]}(G)$-module (this uses that the unit in $\cat{Rep}^\sol_{A[1/\varpi]}G$ is compact by the Lazard resolution, and that compact modules on $\cat{Rep}^\sol_{A[1/\varpi]}G$ are stable under tensor products over $A[\frac{1}{\varpi}]$).
  Moreover, the identity of $P$ is trace class in the category of solid $G$-equivariant $A$-modules (by dualizability of $P$).
  In particular, we can set $C=P$ and $f=\mathrm{Id}_P$ as clearly $P$ is a retract of $\mathrm{cofib}(P\overset{0}{\to}P)\cong P\oplus P[1]$.
\\

  Now, we approximate $C$ and $f$ integrally.
  \begin{claim*}[B] \label{sec:decompl-perf-compl-2-claim-b}
    There exists a compact $\Lambda_{A}(G)$-module $C_0$, a trace class endomorphism $g\colon C_0\to C_0$ in $\cat{Rep}^\sol_{A}(G)$, and $k\geq 0$ such that $P$ is a retract of $\mathrm{cofib}(C_0\overset{g-\omega^k}{\to} C_0)[1/\varpi]$.
  \end{claim*}

  Take $C,f$ as in the previous claim.
  Let $\widetilde{f}\colon A[1/\varpi]\to C^\vee\otimes_{A[1/\varpi]} C$ be a trace-class representative for $f$.
  Writing $C=\varinjlim_{i} M_i$ as a filtered colimit of compact $\Lambda_{A}(G)$-modules, we get $C=C[\frac{1}{\varpi}]= \varinjlim_{i} M_i[\frac{1}{\varpi}]$. 
  As $C$ is compact, it is a retract of the $\varpi$-generic fiber of $M_i$ for some $i\in I$; set $C_0:=M_i$. 
   Moreover, we can assume that $f$ is a retract of a trace class endomorphism $f_0\colon C_0[1/\varpi]\to C_0[1/\varpi]$, with trace class representative $\widetilde{f}_0$: namely, if $c\colon C_0[1/\varpi]\to C$ and $\iota\colon C\to C_0[1/\varpi]$ satisfy $c\circ \iota=\mathrm{Id}_C$, then one can set $f_0=\iota\circ f\circ c$ and $\widetilde{f}_0=(c^\vee\otimes \iota)(\widetilde{f})\in (C_0[1/\varpi])^\vee\otimes_{A[1/\varpi]}C_0[1/\varpi]$.

  Now, each compact $\Lambda_{A}(G)$-module is compact as an $A$-module (this reduces to the compact generators $A\otimes_{\Z_p}\Z_p[G\times S]$ of $\cat{Rep}_{A}^\sol(G)$ for light profinite sets $S$, where the compactness over $A$ is clear).
  This implies that
  \[
    (C_0[1/\varpi])^\vee=\Hom_{A[1/\varpi]}(C_0[1/\varpi],A[1/\varpi])\cong \Hom_{A}(C_0,A)[1/\varpi]=C_0^\vee[1/\varpi],
  \]
  and hence there exists some $k\geq 0$ and a $G$-equivariant morphism
  \[
    \widetilde{g}\colon A\to C_0^\vee\otimes_{A} C_0
  \]
  such that $\widetilde{f}_0=\frac{\widetilde{g}}{\varpi^k}$.
  Let $g\colon C_0\to C_0$ be the endomorphism induced by $\widetilde{g}$.
  Then $\mathrm{cofib}(C_0\overset{g-\varpi^k}{\to}C_0)[1/\varpi]\cong \mathrm{cofib}(C_0[1/\varpi]\overset{f_0-\mathrm{Id}}{\to}C_0[1/\varpi])$ admits $\mathrm{cofib}(C\overset{f-\mathrm{Id}}{\to} C)$ as a retract (cf.\ \Cref{sec:decompl-perf-compl-3-categorical-lemma} with $\lambda=a=1$), and hence $P$ as a retract as well. This proves claim (B).
\\

    Now, let $\epsilon>k$ be  such that $A$ satisfies $\epsilon$-\textbf{TORS} uniformly in $G$ (such an $\epsilon$ exists as $A$ satisfies $\epsilon$-\textbf{TORS} for all sufficiently large $\epsilon$), and define $s:=(d+1)\epsilon$.
    We set
    \[
      D_r:=\tau_{\geq 0}A^{\ob{Bin}_{r,\varpi}}\in \cat{Rep}^\sol_{\Lambda^\circ}(G),\   D:=\varinjlim\limits_{r\to 1} D_r
    \]
    and
    \[
      D_{r,H}:=\Gamma(H,D_r)\in \cat{Rep}^\sol_{\Lambda^\circ}(G)
    \]
    for $r>1$ and $H\subseteq G$ a compact open normal subgroup.
    We continue to use the notation $(-)_{\varpi^s}$ for the derived quotient $(-)/\varpi^s$.

    \begin{claim*}[C]\label{sec:decompl-perf-compl-3-claim-c}
      We have $D_{r,\varpi^s}\cong \varinjlim_{H}D_{r,H,\varpi^s}$ for each $r>1$, $D_{\varpi^s}\cong \varinjlim_{r,H}D_{r,H,\varpi^s}$ and the natural map
      \[
        \alpha\colon D_{\varpi^s}\to A_{\varpi^s}
      \]
      is an almost isomorphism.
    \end{claim*}
    We note that $(\tau_{\geq 0} A^{\ob{Bin}_{r,\varpi}})/\varpi^s$ and $A_{\varpi^s}$ are discrete $G$-equivariant $\Lambda^\circ$-modules, and as such they are the colimit of their derived $H$-invariants for $H$ running through open, normal subgroups in $G$ (e.g.,  from \cite[Remark 3.4.18 and Lemma 3.4.19]{mann2022p}).
    Using our $\epsilon$-\textbf{TORS} on $A$ and \Cref{LemmaApproximationLocAnvsCont} (3), we know that $A^{\ob{Bin}_\varpi}_{\varpi^s}=^a\varinjlim_{r\to 1} (\tau_{\geq 0} A^{\ob{Bin}_{r,\varpi}})_{\varpi^s}$ is almost isomorphic to $A_{\varpi^s}$.
   This altogether implies the claim. 
\\

   Set $\delta>0$ such that $k+2\delta<\epsilon$.
   We will show that $Q(C_0):=\mathrm{cofib}(C_0\overset{\varpi^\delta g-\varpi^{k+\delta}}{\to}C_0)$ satisfies the hypothesis of \Cref{PropLocAnCriterioV2}.

    \begin{claim*}[D]\label{sec:decompl-perf-compl-4-claim-d}
      There exists a compact $G$-equivariant $D_{\varpi^s}$-module $M$ such that $C_{0,\varpi^s}$ is a $G$-equivariant retract of $ A_{\varpi^s}\otimes_{D_{\varpi^s}}M$, and such that $\varpi^\delta g$ is a $\varpi^\delta$-retract of $\mathrm{Id}_{A_{\varpi^s}}\otimes h$ for a $G$-equivariant trace class endomorphism $h\colon M\to M$ in $\cat{Rep}_{D_{\varpi^s}}^{\sol}(G)$.
      In particular,
      $\mathrm{cofib}(C_{0,\varpi^s}\overset{\varpi^\delta g-\omega^{k+\delta}}{\to} C_{0,\varpi^s})$
      is a $\varpi^\delta$-retract of $\mathrm{cofib}(M\otimes_{D_{\varpi^s}}A_{\varpi^s}\overset{h\otimes \mathrm{Id}_{A_{\varpi^s}}-\omega^{k+\delta}}{\to} M\otimes_{D_{\varpi^s}} A_{\varpi^s})$.
    \end{claim*}
    
We note that as $g$ is a trace class endomorphism of $C_{0,\varpi^s}$ as $G$-equivariant  $A_{\varpi^s}$-module, it is a trace class endomorphism in the almost category  $\cat{Rep}^{\sol,a}_{A_{\varpi^s}}(G)=\cat{Rep}^{\sol,a}_{D_{\varpi^s}}(G)$, where the previous equivalence of categories is clear from the almost equivalence of rings of Claim $C$. Therefore, $\Hom_{D_{\varpi^s}}(C_{0,\varpi^s}, D_{\varpi^s})\otimes_{D_{\varpi^s}} C_{0,\varpi^s}$ is almost isomorphic as $G$-equivariant $D_{\varpi}$-representation to $\Hom_{A_{\varpi^s}}(C_{0,\varpi^s}, A_{\varpi^s})\otimes_{A_{\varpi^s}} C_{0,\varpi^s}$.  Therefore, $\varpi^{\delta/3}g$ lifts to a $G$-equivariant trace class endomorphism of $C_{0,\varpi^s}$ as $D_{\varpi^s}$-module.

 The object $C_{0,\varpi^s}$ is almost compact (as in \cite[Definition 2.2.28]{mann2022p}) in the almost category $\cat{Rep}_{A}^{\sol,a}(G)$ of $G$-equivariant semilinear solid $A$-representations, with respect to the almost structure of $(\Lambda^{\circ\circ},\Lambda^{\circ})$. By the previous claim, and by  \cite[Proposition 2.2.29 (ii)]{mann2022p} (with $\n{C}_0$ the class of compact $\Lambda_{D_{\varpi^s}}(G)$-modules), there is a   compact $G$-equivariant $D_{\varpi^s}$-module $M_0$ and $G$-equivariant morphisms $\iota\colon C_{0,\varpi^s}\to M_0\otimes_{D_{\varpi^s}}A_{\varpi^s},\ c\colon M_0\otimes_{D_{\varpi^s}}A_{\varpi^s}\to C_{0,\varpi^s}$ such that $c\circ \iota= \varpi^{\delta/3}\cdot \mathrm{Id}_{C_{0,\varpi^s}}$.

  Moreover,  by the previous claim, there exists a $G$-equivariant $D_{\varpi^s}$-semilinear morphism $\beta\colon A_{\varpi^s}\to D_{\varpi^s}$ such that $\alpha\circ \beta=\varpi^{\delta/3}\cdot \mathrm{Id}_{A_{\varpi^s}}$ and $\beta\circ \alpha= \varpi^{\delta/3}\cdot \mathrm{Id}_{D_{\varpi^s}}$.
    Set
    \[
      \iota^\prime:=(\mathrm{Id}\otimes \beta)\circ \iota\colon C_{0,\varpi^s}\to M_0
    \]
    and
    \[
      c':=c\circ (\mathrm{Id}\otimes \alpha)\colon M_0\to C_{0,\varpi^s}.
    \]
    Then
    \[
      c'\circ \iota'=\varpi^{2\delta/3}\cdot \mathrm{Id}_{C_{0,\varpi^s}} 
    \]
    as $\alpha\circ \beta=\varpi^{\delta/3}\cdot \mathrm{Id}_{D_{\varpi^s}}$ and $c\circ \iota = \varpi^{\delta/3}\cdot \ob{Id}_{C_{0,\varpi^s}}$.

    We use now \Cref{sec:decompl-perf-compl-3-categorical-lemma} in the category $\n{C}=\cat{Rep}_{D_{\varpi^s}}^{\sol}(G)$, with modules $M=M_0$ and $C=C_{0,\varpi^s}$,   morphisms   $\iota=\iota',c= c'$ and $f=\varpi^{\delta/3}g$, and constants $ \lambda=\omega^{2\delta/3}$ and $a=\varpi^{k+\delta}$.
    We conclude that $\varpi^\delta g$ is a $\omega^\delta$-retract of the trace class endomorphism $h:=\iota'\circ (\varpi^{\delta/3} g)\circ c'$ of $M_0$ and that $\mathrm{cofib}(C_{0,\varpi^s}\overset{\varpi^\delta g-\omega^{k+\delta}}{\to} C_{0,\varpi^s})$
      is a $\varpi^\delta$-retract of $\mathrm{cofib}(M_0\otimes_{D_{\varpi^s}}A_{\varpi^s}\overset{h\otimes\mathrm{Id}_{A_{\varpi^s}}-\omega^{k+\delta}}{\to} M_0\otimes_{D_{\varpi^s}}A_{\varpi^s})$. This proves claim (D). 
\\

    Fix now $h\colon M\to M$ as in the previous claim.
    By \Cref{sec:decompl-perf-compl-2-trace-class-morphisms-in-discrete-case} (2) applied to $B=D_{\varpi^s}=\varinjlim_{r,H}D_{r,H,\varpi^s}$, there exists $r>1$, an open normal subgroup $H\subseteq G$ and an $H$-equivariant factorization
    \[
      M\overset{\gamma}{\to} N:=N_0\otimes_{D_{r,H,\varpi^s}} D_{\varpi^s}\overset{\rho}{\to} M
    \]
    of $h$, where $N_0$ is a $D_{r,H,\varpi^s}$-module with trivial $H$-action and $\gamma$ is of trace class in $\cat{Rep}^{\sol}_{D_{\varpi^s}}(H)$.
    By \Cref{sec:decompl-perf-compl-1-approximation-of-trace-class-endomorphisms} (1)  there exists some $1<r'<r$ and $H'\subseteq H$ open and normal in $G$ such that the trace class morphism
    \[
      \gamma\circ \rho\colon N\to N
    \]
    is given by $k\otimes_{D_{r',H',\varpi^s}} D_{\varpi^s}$ for a trace class endomorphism
    \[
      k\colon N_1:=N_0\otimes_{D_{r,H,\varpi^s}}D_{r',H',\varpi^s}\to N_1.
    \]
    We can now replace $G$ by $H$, and assume that $H=H'$, $r=r'$ and thus $N_0=N_1$.
       
    \begin{claim*}[E]\label{sec:decompl-perf-compl-5-claim-e}
      We have that $Q(M):=\mathrm{cofib}(M\overset{h-\omega^{k+\delta}}{\to} M)$ is a $\varpi^{k+\delta}$-retract of $Q(N):=\mathrm{cofib}(N\overset{\gamma\circ \rho-\varpi^{k+\delta}}{\to} N)$.
      By base change, the same holds after tensoring over $D_{\varpi^s}$ with $A_{\varpi^s}$.
    \end{claim*}
    This follows from the commutative diagram
  \[\begin{tikzcd}[column sep=large]
	M & M \\
	N & N \\
	M & M
	\arrow["h-\varpi^{k+\delta}", from=1-1, to=1-2]
	\arrow["\gamma", from=1-1, to=2-1]
	\arrow["\gamma", from=1-2, to=2-2]
	\arrow["{\gamma\circ \rho-\varpi^{k+\delta}}", from=2-1, to=2-2]
	\arrow["\rho", from=2-1, to=3-1]
	\arrow["\rho", from=2-2, to=3-2]
	\arrow["h-\varpi^{k+\delta}", from=3-1, to=3-2]
      \end{tikzcd}\]
    by passing to cofibers and using that $\rho\circ \gamma=h=\varpi^{k+\delta}$ on $\mathrm{cofib}(M\overset{h-\varpi^{k+\delta}}{\to} M)$.

    \begin{claim*}[F]\label{sec:decompl-perf-compl-6-claim-f}
      \Cref{PropLocAnCriterioV2} applies to $Q(C_0)=\mathrm{cofib}(C_0\overset{\varpi^\delta g-\varpi^{k+\delta}}{\to} C_0)$.
      In particular, $P$ admits a decompletion and the proof of \Cref{sec:decompl-perf-compl-1-decompletion-for-pi-complete} is finished.
    \end{claim*}

    We first note that $P$ is a retract of $Q(C_0)[1/\varpi]$ using Claim (B), thus it suffices to show the first part of the claim.
    
    Now, by Claim (D) $Q(C_{0,\varpi^s})$ is a $\varpi^\delta$-retract of $Q(M)\otimes_{D_{\varpi^s}} A_{\varpi^s}$.
    By Claim (E) $Q(M)$ is a $\varpi^{k+\delta}$-retract of $Q(N)=\mathrm{cofib}(N\overset{\gamma\circ \rho-\omega^{k+\delta}}{\to}N)$.
    Now, $\gamma\circ\rho=k\otimes_{D_{r,G,\varpi^s}} D_{\varpi^s}$
    for a $D_{r,G,\varpi^s}$-linear endomorphism $k\colon N_0\to N_0$ by our previous reductions.
    This implies that $Q(N)\cong \mathrm{cofib}(N_0\overset{k-\omega^{k+\delta}}{\to}N_0)\otimes_{D_{r,G,\varpi^s}}D_{\varpi^s}$, where $\mathrm{cofib}(N_0\overset{k-\omega^{k+\delta}}{\to}N_0)$ has the trivial $G$-action.
    Altogether, $Q(C_{0,\varpi^s})$ is a $\omega^{k+2\delta}$-retract of $Q(N)$, and as $k+2\delta<\epsilon$ we see that the hypotheses of \Cref{PropLocAnCriterioV2} are fulfilled.
\end{proof}

\begin{remark}
  \label{sec:decompl-perf-compl-3-categorical-lemma}
  Let $\mathcal{C}$ be a $\Lambda^\circ$-linear category.
  Let $C,M\in \mathcal{C}$ be a objects and let $c\colon M\to C,\ \iota\colon C\to M$ and $f\colon C\to C$ be morphisms in $\mathcal{C}$.
  Assume we are given an isomorphism $c\circ \iota=\lambda\cdot \mathrm{Id}_C$ for some $\lambda\in \Lambda^\circ$.
  Then we have a natural commutative diagram
\[\begin{tikzcd}
	C & C \\
	M & M \\
	C & C
	\arrow["{\lambda\cdot f}", from=1-1, to=1-2]
	\arrow["\iota"', from=1-1, to=2-1]
	\arrow["\iota", from=1-2, to=2-2]
	\arrow["{f\circ c}", from=2-1, to=1-2]
	\arrow["{\iota \circ f\circ c}", from=2-1, to=2-2]
	\arrow["c"', from=2-1, to=3-1]
	\arrow["c", from=2-2, to=3-2]
	\arrow["{\iota\circ f}", from=3-1, to=2-2]
	\arrow["{\lambda\cdot f}", from=3-1, to=3-2]
      \end{tikzcd}\]
    in $\mathcal{C}$.
        If $a\in \Lambda^\circ$, we can conclude that the diagram
\[\begin{tikzcd}
	C & {} & C \\
	M && M \\
	C && C
	\arrow["{\lambda\cdot f-a\cdot \mathrm{Id}_C}", from=1-1, to=1-3]
	\arrow["\iota", from=1-1, to=2-1]
	\arrow["\iota", from=1-3, to=2-3]
	\arrow["{\iota\circ f\circ c-a\cdot \mathrm{Id}_M}", from=2-1, to=2-3]
	\arrow["c", from=2-1, to=3-1]
	\arrow["c", from=2-3, to=3-3]
	\arrow["{\lambda\cdot f-a\cdot \mathrm{Id}_C}", from=3-1, to=3-3]
      \end{tikzcd}\]
    naturally commutes, and thus that $\mathrm{cofib}(C\overset{\lambda\cdot f-a\cdot \mathrm{Id}_C}{\to}C)$ is a $\lambda$-retract of $\mathrm{cofib}(M\overset{\iota\circ f\circ c-a\cdot \mathrm{Id}_M}{\to}M)$.
\end{remark}

We used the following lemmata.

\begin{lemma}
  \label{sec:decompl-perf-compl-1-approximation-of-trace-class-endomorphisms}
  Let $B_i\in \ob{CAlg}(\cat{Rep}^\sol_{\Lambda^\circ}(G)), i\in I,$ be a filtered diagram with colimit $B$.
  \begin{enumerate}
  \item Let $i_0\in I$ and let $M_{i_0}\in \cat{Rep}^\sol_{B_{i_0}}(G)$ be compact.
    If $h\colon M':=B\otimes_{B_{i_0}}M_{i_0}\to M'$ is a trace class morphism in $\cat{Rep}^\sol_{B}(G)$, then there exists $j\geq i_0$ and a trace class morphism $h_j\colon B_j\otimes_{B_{i_0}}M_{i_0}\to B_j\otimes_{B_{i_0}}M_{i_0}$  in $\cat{Rep}^\sol_{B_{j}}(G)$ such that $h=B\otimes_{B_j}h_j$.
  \item Assume that $M\in \cat{Rep}^\sol_B(G)$ is compact and that $f\colon M\to M$ is a trace class endomorphism in $ \cat{Rep}^\sol_{B}(G)$.
   Then there exists some $j\in I$, a compact object $M_j\in \cat{Rep}^\sol_{B_j}(G)$ and a trace class endomorphism $f_j\colon M_j\to M_j$ in $ \cat{Rep}^\sol_{B_j}(G)$ such that $f$ is a retract of $B\otimes f_j\colon B\otimes_{B_j}M_j\to B\otimes_{B_j}M_j$.

  \end{enumerate}
  
\end{lemma}
\begin{proof}
  We note that for each $i\in I$ the category $\cat{Rep}^\sol_{B_i}(G)$ is compactly generated by $B_i\otimes_{\Lambda^\circ} K$ for compact objects $K\in \cat{Rep}^\sol_{\Lambda^\circ}(G)$, so that in (2) we deduce that there exists some $i_0\in I$ and a  compact object $M_{i_0}\in \cat{Rep}_{B_{i_0}}(G)$ such that $M$ is a retract of $B\otimes_{B_{i_0}}M_{i_0}$.
  Following \Cref{sec:decompl-perf-compl-3-categorical-lemma} we can assure that $f$ is a retract of a trace class endomorphism $h$ of $M':=B\otimes_{B_{i_0}}M_{i_0}$ in $\cat{Rep}^{\sol}_{B}(G)$.
  Thus, it suffices to show (1).
  Now, we have $G$-equivariant isomorphisms
  \[
    \Hom_{B}(M',B)\otimes_{B} M'\cong \Hom_{B_{i_0}}(M_{i_0},B)\otimes_{B_{i_0}} M_{i_0}\cong \varinjlim_{i} \Hom_{B_i}(M_i,B_i)\otimes_{B_i}M_i
  \]
  with $M_i:=B_i\otimes_{B_{i_0}}M_{i_0}$, so that each trace class representative for $h$ is defined over some $M_j$ with $j\geq i_0$ in the category $\cat{Rep}^{\sol}_{B_j}(G)$.
 In particular, we deduce the existence of $h_j$ as desired.
\end{proof}

\begin{lemma}
  \label{sec:decompl-perf-compl-2-trace-class-morphisms-in-discrete-case}
  Let $B\in \ob{CAlg}(\mathrm{Rep}^\sol_{\Lambda^\circ}(G))$ be discrete.
  Let $M\in \mathrm{Rep}^\sol_B(G)$ be compact and let $h\colon M\to M$ be a trace class morphism in $\mathrm{Rep}^\sol_B(G)$.
  \begin{enumerate}
  \item The morphism $h$ factors as a $G$-equivariant $B$-semilinear map into a composition
  \[
    M\to M^\sm\to M,
  \]
  where $M^\sm=\varinjlim_{H\subseteq G}\Gamma(H,M)\cong \varinjlim_{H\subseteq G} \Gamma(H,M)\otimes_{\Gamma(H,B)}B$, and the map $M\to M^{\sm}$ is trace class in $\cat{Rep}^{\sol}_{B}(G)$. 
  \item There exists some open normal subgroup $G_0\subseteq G$, a compact $\Gamma(G_0,B)$-module $N$ with trivial $G_0$-action and a factorization of $h$ as a $G_0$-equivariant $B$-linear map into
  \[
    M\to N\otimes_{\Gamma(G_0,B)}B\to M
  \]
  with $M\to N\otimes_{\Gamma(G_0,B)} B$ a trace class morphism in $\cat{Rep}_{B}^{\sol}(G_0)$.
  \end{enumerate} 
\end{lemma}
\begin{proof}
  As $B$ is discrete and $M$ compact, the $G$-equivariant $B$-module $M^\vee=\Hom_B(M,B)$ is discrete and hence $M^\vee\cong \varinjlim_{H\subseteq G} \Gamma(H,M^\vee)$ in $\mathrm{Rep}^\sol_{\Lambda^\circ}(G)$, where $H\subseteq G$ runs through all open normal subgroups.
  We can deduce  that
  \[
    M^\vee\otimes_B M\cong \varinjlim_{H\subseteq G} \Gamma(H,M^\vee)\otimes_{\Gamma(H,B)}M
  \]
  (using as well that $B$ is a discrete $G$-module).
  By compactness of $\Lambda^\circ\in \mathrm{Rep}^\sol_{\Lambda^\circ}(G)$ (implied by compactness of $\Z_p\in \mathrm{Rep}^\sol_{\Z_p}(G)$, using the Lazard resolution), we can deduce that $h$ admits a trace class representative given by a $G$-equivariant morphism
  \[
    \widetilde{h}\colon \Lambda^\circ\to \Gamma(G_0,M^\vee)\otimes_{\Gamma(G_0,B)} M
  \]
  for some open normal subgroup $G_0\subseteq G$.
  In particular, $\widetilde{h}$ factors over a morphism
  \[
    \Lambda^\circ\to \Gamma(G_0,\Gamma(G_0,M^\vee)\otimes_{\Gamma(G_0,B)} M)\cong \Gamma(G_0,M^\vee)\otimes_{\Gamma(G_0,B)} \Gamma(G_0,M)
  \]
  by the projection formula for $\Gamma(G_0,-)$. Thus, the map $h$ arises from a $G$-equivariant morphism  of $B$-modules
  \[
  \begin{aligned}
  B & \to (\Gamma(G_0,M^\vee)\otimes_{\Gamma(G_0,B)} \Gamma(G_0,M))\otimes_{\Gamma(G_0,B)} B  \\ & \to  M^{\vee}\otimes_{B}(B\otimes_{\Gamma(G_0,B)} \Gamma(G_0,M)) \\ & \to M^{\vee}\otimes_B M^{\sm} \\ & \to M^{\vee}\otimes_{B} M. 
  \end{aligned}
  \]
  where  in the third arrow we used that $B=\varinjlim_{H\subset G} \Gamma(H,B)$ as it is discrete. 
  This shows part (1). 
  
  For part (2),  by the previous discussion the endomorphism $h$ factors as a composite in   $\cat{Rep}^{\sol}_{B}(G)$
  \[
  M\xrightarrow{\gamma} B\otimes_{\Gamma(G_0, B)} \Gamma(G_0,M) \to M 
  \]
  where the left map $\gamma$ is trace class. Writing $\Gamma(G_0,M)=\varinjlim_j N_j$ as a filtered colimit of  compact  $\Gamma(G_0,B)$-modules,  the $G_0$-equivariant map  of $B$-modules (with $N_j$ endowed with the trivial $G_0$-action)
  \[
  B\to M^{\vee} \otimes_B(B\otimes_{\Gamma(G_0,B)}  \Gamma(G_0,M) )=\varinjlim_{j} M^{\vee} \otimes_B (B\otimes_{\Gamma(G_0,B)} N_j)
  \]
   inducing the map $\gamma$ factors through a map $B\to M^{\vee} \otimes_B (B\otimes_{\Gamma(G_0,B)} N_j)$ for some $j$. Denoting $N=N_j$ we deduce part (2) as desired. 
   \end{proof}

We now state our main theorem about decompletions of perfect complexes.
We need a technical assumption on affinoids
in order to reduce to \Cref{sec:decompl-perf-compl-1-decompletion-for-pi-complete}.
Let $\Lambda$ be a uniform Banach algebra  with pseudo-uniformizer $\varpi$  as before, that is, we have  $|p|_{\varpi}\leq |\varpi|_{\varpi}=1/2$ and  $\varpi^{\epsilon}\in \Lambda$  exists for $\epsilon$ rational.

\begin{definition}\label{AssumptionKeyTheoremDecompletionPerfectComplexes}
Let $G$ be a $p$-adic Lie group of dimension $d$, $\Lambda$ a uniform Banach $\bb{Q}_p$-algebra as above.
\begin{enumerate}
\item Let $X=\AnSpec (A)$ be an affinoid Berkovich space over $\Lambda$ endowed with a $\Lambda$-linear action by $G^\cont$.   We say that $A$ is \textit{strongly decompleted} if:
 
 \begin{enumerate}
 
 \item[(a)]  there is an open uniform pro-$p$-subgroup  $G_0\subset G$,
 
 \item[(b)]  a direct system $\{A_i\}_{i\in I}$ of $G_0$-equivariant  $\varpi$-adically   complete animated $\Lambda^{\circ}$-algebras with $A_i/\varpi$ discrete, and  such that $A=\varinjlim_i A_i[\frac{1}{\varpi}]$, and
\item[(c)] this datum satisfies the following property:  for each $i\in I$ there is some $\epsilon_i>0$ such that $A_i$ satisfies $\epsilon_i$-\textbf{TORS} uniformly on $G_0$.

\end{enumerate} 
\item Let $X$  be a locally classically decompleted Berkovich space endowed with an action of $G^{\cont}$. Let $X^{\la}$ denote its decompletion by locally analytic vectors.  We say that $X$ is \textit{locally strongly decompleted} if
 \begin{enumerate}
 \item[(a)]  there is $G_0\subset G$ an open compact subgroup,
 
 \item[(b)]  a $G^{\la}_0$-equivariant Berkovich space $Y'$ over $\GSpec (\Lambda)$, and

  \item[(c)]  a $G^{\la}_0$-equivariant epimorphism $Y'\to X^{\la}$ with pullback $Y:=Y'\times_{X^{\la}} X$, such that $Y=\bigsqcup_j \GSpec (R_j)$ is a disjoint union of affinoid Gelfand stacks which are strongly decompleted as in (1).
  \end{enumerate}
\end{enumerate}
 \end{definition}

\begin{remark}\label{Rem-Strong-Decompleted-Implies-Classically-Decompleted}
Let $A$ be a strongly decompleted $\Lambda$-algebra  as in \Cref{AssumptionKeyTheoremDecompletionPerfectComplexes}. Then $A$ is classically decompleted in the sense of \Cref{Definition-Classically-decompleted}. Indeed, the vanishing of of higher locally analytic vectors follow by the condition $\epsilon$-\textbf{TORS}, and the equivalence on uniform completions of the map $A^{\la}\to A$ follows from \Cref{LemmaApproximationLocAnvsCont} (3). Indeed, by \textit{loc. cit.} the natural map   $ \varinjlim_r \tau_{\geq 0}A_i^{\ob{Bin}_{r,\varpi}} \to A_i$ is an almost equivalence after $\varpi$-completion. Passing to colimits on $i$, inverting $p$, and taking uniform completions we get the desired equivalence $A^{\la,u}\cong A^u$. Therefore, if $X$ is a  $G^{\cont}$-equivariant Berkovich space which  locally in the analytic topology is strongly decompleted, then  it is locally classically decompleted. 
\end{remark}

\begin{theorem}\label{TheoremDecompletionPerfectComplexesII}
 Let $X$  be a locally classically decompleted Berkovich space endowed with an action of $G^{\cont}$ for $G$ a $p$-adic Lie group, with decompletion $X^\la$. If $X$ is locally strongly decompleted, the pullback of perfect complexes
 \[
 \Perf(X^{\la}/G^{\la})\to \Perf(X/G^{\cont})
 \]
 is an equivalence of categories. In particular, \Cref{TheoDecompletion} (3) holds true. 
 \end{theorem}

\begin{proof}
By \Cref{CorollaryLocAnDecompletion} it suffices to prove that $G^{\cont}$-equivariant perfect complexes on  a strongly decompleted affinoid Gelfand stack $X=\GSpec(R)$ admit locally analytic  decompletions. 
Write $R=\varinjlim_{i}A_i[1/\varpi]$ as in our assumption on $X$.
  By \Cref{sec:decompl-perf-compl-1-decompletion-for-pi-complete} it suffices to show that each perfect complex $P\in \ob{D}(\GSpec(R)/G^\cont)$ is a retract of the pullback of a perfect complex $P_i\in \ob{D}(\GSpec(A_i[1/\varpi])/G^\cont)$ for some $i$.
  As the rings $A_i$ are nuclear over $\Lambda^{\circ}$ (having discrete mod $\varpi$-fibre), the ring $R$ is nuclear over $\Lambda$ and hence over $\bb{Q}_p$. We deduce that the perfect complex $P$ is nuclear over $\bb{Q}_p$, so we may equivalently argue for $P\in \mathrm{Rep}^\sol_{R}(G)$, cf.\ \Cref{Remark-Difference-Actions}.
  As in Claim A of the proof of \Cref{sec:decompl-perf-compl-1-decompletion-for-pi-complete}, we can write $P$ as a retract of the cofiber of $f-\mathrm{Id}_C$, where $C\in \mathrm{Rep}^\sol_{R}(G)$ is compact and $f\colon C\to C$ is trace class (one can just take $C=P$ and $f=\mathrm{id}_P$).
  By \Cref{sec:decompl-perf-compl-1-approximation-of-trace-class-endomorphisms}, there exists some $i\in I$ such that $f$ is the base change of a trace class morphism $f_i\colon C_i\to C_i$ for some compact object $\mathrm{Rep}^\sol_{A_i[1/\varpi]}(G)$.
  By \Cref{sec:decompl-perf-compl-3-categorical-lemma}, we can deduce that $P$ is a retract of a base change of $Q_i:=\mathrm{cofib}(C_i\overset{f_i-\mathrm{Id}_{C_i}}{\to} C_i)$, and it suffices to see that $Q_i$ is perfect.
  But $A_i[1/\varpi]$ is Gelfand, so it suffices to see that $Q_i$ is dualizable.
  But $Q_i$ is compact and its identity is induced by $f_i$, and thus of trace class.
  This finishes the proof (note that fully faithfulness of $\Perf(X^\la/G^\la)\to \Perf(X/G^{\cont})$ was already proven in \Cref{TheoDecompletion} (1)).  
\end{proof}

We aim to extend the decompletion of perfect complexes to quotients $Y/G\to X/G$, where $Y\to X$ is Berkovich-\'etale.

\begin{corollary}
  \label{sec:decompl-perf-compl-1-decompletion-etale-coverings}
  Let $X$ and $G$ be as in \Cref{TheoremDecompletionPerfectComplexesII}.
  Assume that $Y\to X$ is a Berkovich-\'etale $G^\cont$-equivariant morphism of Berkovich spaces. Then $Y$ is locally classically decompleted and  the pullback  functor $\Perf(Y^\la/G^\la)\to \Perf(Y/G^\cont)$ is an equivalence.
\end{corollary}
\begin{proof}
  We may assume that $G$ is a uniform pro-$p$ group.
  We note that $Y$ admits locally a classical decompletion by \Cref{RemarkClassicalDecompletionGelfandRingEtaleSites}.
  We may also localize on $X$ and assume that $X=\GSpec(R)$ (after possibly shrinking $G$ again).
  Similarly, we may assume that $Y=\GSpec(R')$ is affinoid.
  By \Cref{CorollaryLocAnDecompletion} and our assumption on $X$, we may assume that $R$ is a $\Lambda$-algebra, which is a filtered colimit $R=\varinjlim_{i\in I} A_i[1/\varpi]$ as in our assumptions on $X$, i.e., each $A_i$ satisfies the assumptions in \cref{sec:decompl-perf-compl-1-decompletion-for-pi-complete} (for varying $\epsilon$).

  First assume that $Y\to X$ is a $G^\cont$-stable $\dagger$-rational localization, defined as the locus $\{|f_1|\leq |g|,\ldots, |f_n|\leq |g|\}$ for some elements $f_1,\ldots, f_n,g\in \pi_0(R)$ generating the unit ideal.
  There exists some $i\in I$ such that $f_1,\ldots, f_n,g$ lift to elements in $\pi_0(A_i[1/\varpi])$, we call these lifts $f_1,\ldots, f_n,g$ again.
  Multiplying by a power of $\varpi$ we may even assume that $f_1,\ldots, f_n,g$ lie in $A_i$.
  There exists some $k\geq 0$ such that changing the $f_1,\ldots, f_n,g$ by some element in $\varpi^kA_i$ does not change the locus $\{|f_1|\leq |g|,\cdots, |f_n|\leq |g|\}$.
  As $A_i/\varpi$ is discrete, we may assume that the residue classes of $f_1,\ldots, f_n,g\in A_i$ in $A_{i}/\varpi^{k+1}$ are fixed by $G$ (by shrinking $G$ if necessary).  Moreover, since the locus $\{|f_i|\leq |g|\}$ is compact and contains $\{|g|\neq 0\}$, it contains the locus $\{|\varpi^b|\leq |g|\}$ for some  $b\gg 0$. Thus, we can suppose without loss of generality that $f_{1}= \varpi^b$.  In particular, the $\dagger$-rational localization is a composite of $\dagger$-rational localizations of the form $\{|\varpi^k|\leq |g|\}$ and $\{|f_i/g| \leq 1\}$ for $i=2,\ldots, n$.

  Now, for $\delta>0$, let $Y_{i,\delta}=\GSpec(A_i\langle \frac{f_1}{\varpi^{-\delta}g}, \cdots, \frac{f_n}{\varpi^{-\delta}g}\rangle[1/\varpi])$ be the Tate-rational localization of $\GSpec(A_i[1/\varpi])$ defined by $\{|f_1|\leq |\varpi^{-\delta}g|,\cdots, |f_n|\leq |\varpi^{-\delta}g|\}$, in particular $B_{i,\delta}:=A_i\langle \frac{f_1}{\varpi^{-\delta}g},\cdots, \frac{f_n}{\varpi^{-\delta}g}\rangle$ is $\varpi$-complete and discrete mod $\varpi$. 
  We note that $Y_{i,\delta}$ is stable under $G$ if $\delta$ is sufficiently small, and that $Y$ is the inverse limit over $j\geq i$ and $\delta>0$ of the base changes $Y_{i,\delta}\times_{\GSpec(A_i[1/\varpi])} \GSpec(A_{j}[1/\varpi])=\GSpec(B_{j,\delta}[1/\varpi])$, with $B_{j,\delta}=A_j\langle \frac{f_1}{\varpi^{-\delta}g},\cdots, \frac{f_n}{\varpi^{-\delta}g}\rangle$.  
  By (iterating) \Cref{LemmaEtaleLocalizationDecomplete} (2), we know that $B_{j,\delta}$ carries a natural action of $G$ and  satisfies the assumption $\epsilon$-\textbf{TORS} uniformly on $G$ for some $\epsilon'>0$ as $A_j$ satisfies it for some $\epsilon>0$.
  As clearly $B_{j,\delta}$ is $\varpi$-complete and discrete mod $\varpi$, we see that perfect complexes on $Y/G^\cont$ decomplete by \Cref{TheoremDecompletionPerfectComplexesII}.

  Next, assume that $Y=\GSpec(R')\to X=\GSpec(R)$ is Berkovich-\'etale.
  After $\dagger$-rational localizations on $X$ and $Y$ (using what we have already shown) we may assume that $R\to R'$ is finite \'etale and that $\Perf(X/G)\cong \Perf(X^\la/G^\la)$.
  By \Cref{LemmaLocallyAnalyticDecompletionSites} we know that $R'\cong R^{\prime,\la}\otimes_{R^\la}R$ for a finite \'etale $R^\la$-algebra $R^{\prime,\la}$.
  If $M$ is a $G$-equivariant perfect complex over $R'$, its restriction to an $R$-module is a $G$-equivariant perfect complex over $R$.
  By assumption, $M$ admits a decompletion as an $R$-module.
  By \Cref{Lemma-Decompletions} we can deduce that $M$ admits a decompletion as an $R'$-module.
  This finishes the proof.  
\end{proof}

A direct corollary of \Cref{TheoremDecompletionPerfectComplexesII} is a suitable version in families. In the following we fix a Banach field extension $K$ of $\bb{Q}_p$, and assume that the ring $\Lambda$  is a $K$-algebra.

\begin{corollary}\label{CorollaryDescentCoefficients}
Let $X$ be a Berkovich space  over a Banach field  $K$  endowed with an action of $G^{\cont}$ for $G$ a $p$-adic Lie group. Let $\{B_{j,0}\}_{j\in J}$ be a cofiltered famility of $p$-adically complete animated $\n{O}_K$-algebras with $B=\varinjlim_j B_{j,0}[\frac{1}{p}]$.  Suppose that $X$ is locally strongly decompleted, then so is $X_{B}:=X\times_{\GSpec (K)} \GSpec (B)$. In particular, the natural map 
\[
\Perf(X_{B}^{\la}/G^{\la})\xrightarrow{\sim} \Perf(X_B/G^{\cont})
\]
is an equivalence. 
\end{corollary}

\begin{proof}
By applying \Cref{CorollaryLocAnDecompletion}, we can replace $X$ by the cover $Y'\to X$ and  assume without of generality that $G$ is uniform pro-$p$ and that $X=Y'=\GSpec (R)$ as in  \Cref{AssumptionKeyTheoremDecompletionPerfectComplexes}.  Then, we can write $R =\varinjlim_{i} A_i[\frac{1}{\varpi}]$  as  the generic fibre of a filtered diagram of  $\varpi$-complete $\Lambda^{\circ}$-algebras. Then, we see that  $R\otimes_{K}B = \varinjlim_{i,j} (A_i\otimes_{\n{O}_K} B_{j,0})[\frac{1}{\varpi}]$, and it suffices to show that the rings $B_i\otimes_{\n{O}_K} C_{j,0}$ satisfy the condition $\epsilon_{i,j}$-\textbf{TORS} for some $\epsilon_{i,j}>0$ uniformly on $G$.  This follows from  \Cref{CoroStabilityBaseChange} as the action of $G$ is trivial both in $\n{O}_K$ and $C_{j,0}$. 
\end{proof}

\begin{remark}\label{RemarkFamilyVersionNuclearGelfandAlgebras}
We expect \Cref{CorollaryDescentCoefficients} to hold more generally for  families of nuclear Gelfand $\bb{Q}_p$-algebras. We will not pursue this question in this paper. 
\end{remark}

\subsection{Decompletion and classical Sen theory} \label{subsec-Sen-theory-revisited}

We have proven a  decompletion of perfect complexes by locally analytic vectors in \Cref{TheoDecompletion}, see \Cref{TheoremDecompletionPerfectComplexesII}. 
 In this section, we prove that the standard frameworks of Sen theory, the cyclotomic tower and the toric tower, are locally strongly decompleted.

\subsubsection{Sen theory: abstract framework} \label{sss:SenTheoryAbstractFramework}

We begin our discussion of Sen theory with the appropriate abstract framework  that applies both in the arithmetic and geometric set up.
The axioms to be presented down below are stronger assumptions than those presented in \cite{camargo2022geometricsen}.

Let $\Lambda$ be a uniform Banach $\bb{Q}_p$-algebra with a pseudo-uniformizer $\varpi$, and let $\Lambda^\circ\subseteq \Lambda$ be the $\varpi$-complete subring of powerbounded elements.
We assume that $\varpi^{\epsilon}$ exists in $\Lambda$ for all $\epsilon\in \bb{Q}$ (this is again just for convenience: it would be enough to assume that $\varpi^\epsilon$ exists for a dense subset of $\Q$), and we endow $\Lambda$ with the $\varpi$-norm $|-|_{\varpi}$ such that for all $x\in \n{M}(\Lambda)$ we have $|\varpi(x)|_{\varpi}= \frac{1}{2}$. From now on, all elements $r\in \bb{R}_{>0}$ will be of the form  $r=|\varpi^{\epsilon}|$ for some  $\epsilon\in \bb{Q}$.
Let $\Gamma\cong \bb{Z}_p$ be a $1$-parameter $p$-adic Lie group  with generator $\gamma$, and for $k\geq 0$  we denote $\Gamma_k= \gamma^{p^k\bb{Z}_p} $ and set $\gamma_k:= \gamma^{p^k}$.
Let $\{A_{n}\}_{n\in \bb{N}}$ be a direct system of static classically $\varpi$-complete and $\varpi$-torsion free  animated $\Lambda^\circ$-algebras with $A_{\infty}=(\varinjlim_{n} A_n)^{\wedge}_{\varpi}$.
Assume that $A_n/\varpi$ is discrete for any $n\geq 0$. We consider the following assumption:

\begin{assumption}\label{AssumptionSenTheory}
Given   $k\geq 0$ there is $n_0\geq 0$  such that for all  $n\geq n_0$ there are  constants  $0<\rho_{n,k}<\epsilon_{n,k}$ with $\epsilon_{n,k}\to 0$ as $n\to \infty$\footnote{If there exists some $\delta>0$ such that $\epsilon_{n,k}\geq \delta>0$ for all $n$, then $A_\infty$ satisfies $\epsilon$-\textbf{TORS} by \Cref{sec:stab-prop-cond-locally-analytic-satisfy-eps-tors}, so this case is already settled.} such that:
\begin{enumerate}

\item  The action of $\Gamma_{k}$ on $A_{n}/\varpi^{\epsilon_{n,k}}$ is trivial. 

\item  The cohomology  groups of $\Gamma( \Gamma_k, A_{\infty}/A_n)$ are  concentrated in degree $1$ and of $\varpi^{\rho_{n,k}}$-torsion. 

\end{enumerate}
\end{assumption}

Assumption (1) implies that $A_n$ is a $C_{r}(\Gamma_k,\Lambda^\circ)$-comodule for any $r<r_{n,k}:=|\varpi^{-\epsilon_{n,k}}|$.
Indeed, we have the following general lemma.

\begin{lemma}\label{LemmaClaimStaticComodulesCrvarpi} Let $2>r>1$ and let $M$ be a static $\Lambda^{\circ}$-module which is classically $\varpi$-complete  and $\varpi$-torsion free endowed with a continuous action of $\Gamma\cong \bb{Z}_p$. Suppose that $M/\varpi$ is discrete. 
Let $\delta>0$ be such that $r<|\varpi^{\delta}|^{-1}$ and suppose that the action of $\Gamma$ on $M/\varpi^{\delta}$ is trivial. Then the orbit map $\n{O}\colon M\to C(\Gamma,M)= C(G,\Lambda^{\circ})\otimes_{\Lambda^{\circ}} M$ factors through the subspace $C_{r,\varpi}(\Gamma,M):=C_{r,\varpi}(\Gamma, \Lambda^{\circ})\otimes_{\Lambda^{\circ}} M \subset C(\Gamma,M)$. In particular, $M$ has a natural structure as a $C_{r,\varpi}(\Gamma, \Lambda^{\circ})$-comodule.
\end{lemma}  
\begin{proof}
  To show that the orbit map $\n{O}\colon M\to C(\Gamma,M)$ factors trough $C_{r,\varpi}(\Gamma,M)$, we can equivalently prove that the multiplication $\bb{Z}_{p,\sol}[\Gamma]\otimes_{\bb{Z}_{p,\sol}} M\to M$ extends (necessarily uniquely by density) to a map $\n{D}_{r,\varpi}(\Gamma,\Lambda^{\circ})\otimes_{\Lambda^{\circ}} M\to M$ where  $\n{D}_{r,\varpi}(\Gamma,\Lambda^{\circ})$ is the $\Lambda^{\circ}$-dual of $C_{r,\varpi}(\Gamma,M)$.
  Let $\gamma\in \Gamma$ be a generator. Let $\delta>\delta'>0$ be such that $|\varpi^{\delta}|<|\varpi^{\delta'}|<r^{-1}$  and denote $s:=|\varpi^{\delta'}|^{-1}$.   By  the power series description of \eqref{eqBinomDistAlgebrasFormulas33} we have a natural morphism of algebras $\n{D}_{r,\varpi}(\Gamma,\Lambda^{\circ})\to \n{D}_{s,\varpi}(\Gamma,\circ) = \Lambda^{\circ}[[\frac{\gamma-1}{\varpi^{\delta'}}]]$. Thus, it suffices to see that the action of $\bb{Z}_{p,\sol}[\Gamma]\cong \bb{Z}_p[[\gamma-1]]$ on $M$ extends to   $ \Lambda^{\circ}[[\frac{\gamma-1}{\varpi^{\delta'}}]]$.  By assumption, the action of $\gamma-1$ on $M$ is divisible by $\varpi^{\delta}$, this implies that the action of $\frac{\gamma-1}{\varpi^{\delta'}}$ on $M$ is topologically nilpotent and therefore  it naturally extends to an action of the algebra $\Lambda^{\circ}[[\frac{\gamma-1}{\varpi^{\delta'}}]]$, proving the claim. 
\end{proof}

\begin{proposition}\label{PropositionTechnicalLemmaTraces}
Let $\Gamma$ and $A_{\infty}$ be as in \Cref{AssumptionSenTheory}. Then the ring  $A_{\infty}$ satisfies   the property $\epsilon$-\textbf{TORS} uniformly on $\Gamma$ for all $\epsilon>0$. Furthermore, $(A_{\infty}[\frac{1}{\varpi}])^{\la}=\varinjlim_n A_n[\frac{1}{\varpi}]$.
\end{proposition}
\begin{proof}
  Let $k\geq 0$ and $\epsilon>0$.
  It suffices to show the following: for $n\geq n_0$ and for each $|\varpi|^{\epsilon_{n,k}}<r^{-1}<|\varpi|^{\rho_{n,k}}$ we have

\begin{enumerate}

\item the cohomology of  $(A_{\infty}/A_n)^{\ob{Bin}_{r,\varpi}(\Gamma_k)}$ is of $\varpi^{\rho_{n,k}}$-torsion, and

\item the $H^1$-cohomology of $A_n^{\ob{Bin}_{r,\varpi}(\Gamma_k)}$ is $\varpi^{\epsilon_{n,k}}$-torsion. 

\end{enumerate}

Indeed, from the short exact sequence $0\to A_n\to A_{\infty} \to A_{\infty}/A_n\to 0$ we have a long exact sequence of homotopy groups

\[
\cdots  \to \pi_{-1}(A_n^{\ob{Bin}_{r,\varpi}(\Gamma_k)})\to \pi_{-1}(A_{\infty}^{\ob{Bin}_{r,\varpi}(\Gamma_k)}) \to \pi_{-1}((A_{\infty}/A_n)^{\ob{Bin}_{r,\varpi}(\Gamma_k)})\to 0.
\]
This implies that $\pi_{-1}(A_{\infty}^{\ob{Bin}_{r,\varpi}(\Gamma_k)})$ is $\varpi^{\rho_{n,k}+\epsilon_{n,k}}$-torsion.  Since $\rho_{n,k}+\epsilon_{n,k}\to 0$ as $n\to \infty$ (for $r\in (\rho_{n,k}, \epsilon_{n,k})$), we deduce that $A_{\infty}$ satisfies $\epsilon$-\textbf{TORS} uniformly on $\Gamma$.

Let us first prove the claim (1).
For this, let $\delta\in (\rho_{n,k}, \epsilon_{n,k})$ such that $ r^{-1} < s^{-1}:=|\varpi^{\delta}| <|\varpi^{\rho_{n,k}}|$.
We first show that the action of $\gamma_k$ on $C_{r,\varpi}(\Gamma_k,\Lambda^{\circ})$ is trivial modulo $\varpi^{\delta}$.
We have an inclusion of Hopf algebras $C_{r,\varpi}(\Gamma_k,\Lambda^{\circ})\subset C_{s,\varpi}(\Gamma_k,\Lambda^{\circ})$  which dualizes to a morphism 
\[
\n{D}_{s,\varpi}(\Gamma_k,\Lambda^{\circ})\to \n{D}_{r,\varpi}(\Gamma_k,\Lambda^{\circ}).
\]
By the power-series description of \eqref{eqBinomDistAlgebrasFormulas33} we find that $\n{D}_{s,\varpi}(\Gamma_k,\Lambda^{\circ})= \Lambda^{\circ}[[\frac{\gamma_k-1}{\varpi^{\delta}}]]$. In particular, $\gamma_k-1$ is divisible by $\varpi^{\delta}$ on  $\n{D}_{s,\varpi}(\Gamma_k,\Lambda^{\circ})$ and then so is in $\n{D}_{r,\varpi}(\Gamma_k,\Lambda^{\circ})$, taking duals one deduces that the action of $\gamma_k$ is trivial on $C_{r,\varpi}(\Gamma_k,\Lambda^{\circ})$ modulo $\varpi^{\delta}$.  By \Cref{LemmaTorsion}, to show that $\pi_{-1}((A_{\infty}/A_n)^{\ob{Bin}_{r,\varpi}(\Gamma_k)})$ is $\varpi^{\rho_{n,k}}$-torsion, it suffices to prove it modulo $\varpi^{\delta}$.  But then we have 
\[
\begin{aligned}
(A_{\infty}/A_n)^{\ob{Bin}_{r,\varpi}(\Gamma_k)}_{\varpi^{\delta}} &  = \Gamma(\Gamma_k, A_{\infty}/A_n\otimes_{\Lambda^{\circ}} C_{r,\varpi}(\Gamma_k,\Lambda^{\circ})_{\varpi^{\delta}})  \\ 
& =  \Gamma(\Gamma_k, A_{\infty}/A_n)\otimes_{\Lambda^{\circ}}  C_{r,\varpi}(\Gamma_k,\Lambda^{\circ})_{\varpi^{\delta}} \\ 
\end{aligned}
\]  
where in the second equivalence we used the projection formula for $\Gamma_k$-cohomology.
By \Cref{AssumptionSenTheory} (2),  $H^*(\Gamma_k, A_{\infty}/A_n)$ is concentrated in degree $1$ and of  $\varpi^{\rho_{n,k}}$-torsion, in particular multiplication by $\varpi^{\rho_{n,k}}$ on $\Gamma(\Gamma_k, A_{\infty}/A_n)$ is homotopic to zero. This implies that  multiplication by $\varpi^{\rho_{n,k}}$ on $((A_{\infty}/A_n)^{\ob{Bin}_{r,\varpi}(\Gamma_k)}_{\varpi^{\delta}})$ is homotopic to zero, and   in particular of $\varpi^{\rho_{n,k}}$-torsion as wanted.

Next, we prove (2).
By \Cref{LemmaClaimStaticComodulesCrvarpi}, we know that $A_n$ is a $C_{r,\varpi}(\Gamma_k,\Lambda^{\circ})$-comodule, and thus in particular an object in $B\f{G}_{\ob{Bin}_{r,\varpi}}$.
Consider the natural map of analytic stacks $f\colon BG^{\cont}_{\Lambda^{\circ}}\to B\f{G}_{\ob{Bin}_{r,\varpi}}$; this map is $0$-proper as $\f{G}_{\ob{Bin}_{r,\varpi}}$ is affinoid with the induced structure from $\Lambda^{\circ}$, and since  $BG^{\cont}_{\Lambda^{\circ}}$ is $0$-proper (using that the map $\AnSpec (\Lambda^{\circ})\to  BG^{\cont}_{\Lambda^{\circ}}$ is $0$-proper and descendable thanks to the Lazard resolution).
By \Cref{LemmaRightAdjointIntegralrAnVectors} and the projection formula we deduce that 
\[
A_{n}^{\ob{Bin}_{r,\varpi}(\Gamma_k)}=  f_*f^* A_n = f_* \Lambda^{\circ} \otimes_{\Lambda^{\circ}} A_n=  \Lambda^{\ob{Bin}_{r,\varpi}(\Gamma_k)} \otimes_{\Lambda^{\circ}} A_n. 
\]
By \Cref{ExampleTechnicalPropositionHoldsForTrivial}, the $H^1$-cohomology group of  ${\Lambda^{\circ}}^{\ob{Bin}_{r,\varpi}(\Gamma_k)}$ is $\varpi^{\epsilon_{n,k}}$-torsion (because $r<|\varpi^{-\epsilon_{n,k}}|$).
This implies the same torsion bound for $H^1(A_{n}^{\ob{Bin}_{r,\varpi}(\Gamma_k)})$ proving what we wanted. 

The final claim about the locally analytic vectors follows from the fact that $(A_{\infty}/A_n)^{\ob{Bin}_{r,\varpi}(\Gamma_k)}$ is torsion by taking colimits as $r\to 1$ and $n\to \infty$. 
\end{proof}

\subsubsection{Sen theory: arithmetic case} \label{sss:SenTheoryArithmetic} In this section we consider Sen theory for the Fargues--Fontaine curve over the cyclotomic extension of $\Q_p$. Classically, this is considered in  \cite[Section III]{CherbonnierColmez}. Let us first state the main result of this subsection. 

Let $\bb{Q}_p^{\cyc}$ be the cyclotomic extension of $\bb{Q}_p$ and let $\bb{Q}_{p}^{\cyc,\flat}$ its tilt. Let $\epsilon=(\zeta_{p^n})_{n\geq 0}\in \bb{Q}_p^{\cyc,\flat}$ denote a  sequence of $p$-power roots of unity, and denote $\pi=\epsilon-1$. We have an isomorphism of perfectoid $\bb{F}_p$-algebras 
\begin{equation}\label{eqIsoPerfectoidFpAlgebrasCyc}
\bb{Q}_p^{\cyc, \flat} = \bb{F}_p((\pi^{1/p^{\infty}})).
\end{equation}
 Let $\Gamma:=\ob{Gal}(\bb{Q}_p(\zeta^{p^{\infty}})/\bb{Q}_p)$ denote the cyclotomic Galois extension and let $\chi^{\ob{cyc}}\colon \Gamma\xrightarrow{\sim} \bb{Z}_p^{\times}$ be the cyclotomic character. The isomorphism \eqref{eqIsoPerfectoidFpAlgebrasCyc} is $\Gamma$-equivariant for the natural action on $\bb{Q}_p^{\cyc}$, and the (unique) $\bb{Z}_p^{\times}$-action on $\bb{F}_p((\pi^{1/p^{\infty}}))$ determined by 
 \[
 a \cdot \pi = \epsilon^{a}-1= (\pi+1)^a-1. 
 \]
 for $a\in \bb{Z}_p^{\times}$.

Let $1+\bb{D}_{n,\bb{Q}_p}^{<1}$ be the open unit disc centered at $1$ with coordinate $q^{1/p^n}$. It is endowed with a $\bb{Z}_p^{\times,\la}$-action given by $a\cdot q^{1/p^n}= q^{a/p^n}$, and Frobenius given by $\varphi(q^{1/p^n})=q^{1/p^{n-1}}$.    Let $1+\widetilde{\bb{D}}^{<1}_{\bb{Q}_p}= \varprojlim_{n} (1+\bb{D}_{n,\bb{Q}_p}^{<1})$ denote the uncompleted perfectoid disc around $1$, and let  $1+\widetilde{\bb{D}}^{<1,u}_{\bb{Q}_p}$ be its uniform completion.

We want to show the following theorem.

\begin{theorem}[Sen theory: arithmetic case]\label{PropositionSenTheoryCyc}
Consider $X= (1+\widetilde{\bb{D}}^{<1,u}_{\bb{Q}_p} )\backslash \{1\}$ and $G=\bb{Q}_p^{\times}$.  Then $X$ is locally strongly decompleted. Furthermore, $X^{\la}= ( 1+\widetilde{\bb{D}}^{<1}_{\bb{Q}_p}  )\backslash \{1\}$. 

\end{theorem}
The proof of \Cref{PropositionSenTheoryCyc} will require a number of preliminary lemmas. By definition of being locally strongly decompleted, and \Cref{CorollaryLocAnDecompletion} (1),  to prove the theorem it suffices   to find a suitable $\bb{Z}_p^{\times,\la}$-equivariant  cover $Y'$    of $( 1+\widetilde{\bb{D}}^{<1}_{\bb{Q}_p}  )\backslash \{1\}$ with pullback $Y= (( 1+\widetilde{\bb{D}}^{<1,u}_{\bb{Q}_p}  )\backslash \{1\})\times_{( 1+\widetilde{\bb{D}}^{<1}_{\bb{Q}_p}  )\backslash \{1\}}Y'$ and $Y^{\prime}=Y^{\la}$,  where $Y$ consists of a disjoint union of  strongly decompleted affinoids as in  \Cref{AssumptionKeyTheoremDecompletionPerfectComplexes}.  Write $Y = \bigsqcup_i \GSpec (R_i)$. The key idea is to construct $Y$  such that $R_i = A_{\infty,i}[\frac{1}{\varpi}]$ where $A_{\infty,i}$ is  a  static, $\varpi$-torsion free and $\varpi$-complete  $\Lambda^{\circ}$-algebra of the form $A_{\infty,i}= (\varinjlim_n A_{n,i})^{\wedge}_{\varpi}$  verifying the axioms of \Cref{AssumptionSenTheory}. By \Cref{PropositionTechnicalLemmaTraces}, the rings $R_i$ will be strongly decompleted as desired and therefore classically decompleted (cf. \Cref{Rem-Strong-Decompleted-Implies-Classically-Decompleted}).

To construct such affinoids and to establish the key torsion properties of \Cref{AssumptionSenTheory}, we need to have a rather precise control of the action of $\Gamma$ on affinoid subspaces of $(1+\widetilde{\bb{D}}_{\bb{Q}_p}^{<1,u})\backslash \{1\}$.
The key idea behind is the classical construction of traces for the cyclotomic extension that we explain in the following lemma. 

\begin{lemma}\label{LemmaClassicalTracesCyclotomicSenTheory}
Denote $q=[\epsilon]\in \bb{A}_{\inf}(\bb{Q}_p^{\cyc})$. Then the natural map 
\[
\widehat{\bigoplus}_{r\in [0,1)\cap \bb{N}[\frac{1}{p}]}  \bb{Z}_p[[q-1]] \cdot q^r \to  \bb{A}_{\inf}(\bb{Q}_p^{\cyc})
\]
of $(p,q-1)$-adically complete modules is an equivalence. 
\end{lemma}
\begin{proof}
  By Nakayama's lemma, it suffices to prove the claim after taking reductions modulo $p$.
  Note that the reduction  of $q-1$ modulo $p$ is $\pi\in \bb{Q}_p^{\cyc,\flat}\cong \bb{F}_p((\pi^{1/p^{\infty}}))$ with $\pi=\epsilon-1$. Thus,   we want to show that the natural map  of $\pi$-adically complete modules
\[
\widehat{\bigoplus}_{r\in  [0,1)\cap \bb{N}[\frac{1}{p}]} \bb{F}_p[[\pi]] \epsilon^{r} \to \bb{F}_p[[\epsilon^{1/p^{\infty}}-1]]
\]
is an equivalence. By passing to $\pi$-completions, it suffices to show that the natural map 
\[
\widehat{\bigoplus}_{r\in  [0,1)\cap \bb{N}[\frac{1}{p}]} \bb{F}_p[\epsilon-1] \epsilon^r\to  \bb{F}_p[\epsilon^{1/p^{\infty}}-1]
\]
is an equivalence, but $\bb{F}_p[\epsilon-1]=\bb{F}_p[\epsilon]$ and  $\bb{F}_p[\epsilon^{1/p^{\infty}}-1]=\bb{F}_p[\epsilon^{1/p^{\infty}}]$ as these are just polynomial algebras. It is then clear that $\epsilon^{r}$ for $r\in  [0,1)\cap \bb{N}[\frac{1}{p}]$ is a basis of $\bb{F}_p[\epsilon^{1/p^{\infty}}]$ over $\bb{F}_p[\epsilon]$, finishing the proof of the lemma. 
\end{proof}

\begin{remark}\label{RemarkSenTracesCyclotomicTrace}
By taking Frobenius inverses, we see that \Cref{LemmaClassicalTracesCyclotomicSenTheory} holds for $q^{1/p^n}$ instead of $q$, namely, we have a natural isomorphism of $(p,q^{1/p^n}-1)$-complete modules
\[
\widehat{\bigoplus}_{r\in [0,1)\cap \bb{N}[\frac{1}{p}]} \bb{Z}_p[[q^{1/p^n}-1]] \cdot q^{r/p^n}\xrightarrow{\sim} \bb{A}_{\inf}(\bb{Q}_p^{\cyc}). 
\]
\end{remark}

Let $a\in \bb{Z}_p^{\times}$, by construction one has that $a\cdot q^{r}= q^{ar}$ for any $r\in \bb{N}[\frac{1}{p}]$.
Let $I_{m}\subset (0,1)\cap \bb{N}[\frac{1}{p}]$ the subset of rational numbers  with $p$-adic valuation $v_p(r)\leq -m$, and let $I_m^*:=I_m\backslash I_{m+1}$ be the subset of rational numbers with $v_p(r)=-m$.
Then, for all $n\geq 0$ and $m\geq 1$ the action of $\bb{Z}_p^{\times}$ on $\bb{A}_{\inf}(\bb{Q}_p^{\cyc})$ leaves the subspaces 
\begin{equation}\label{eqSubSpacesStableAction}
\widehat{\bigoplus}_{r\in I_m} \bb{Z}_p[[q^{1/p^n}-1]] \cdot q^{r/p^n} \mbox{ and } \widehat{\bigoplus}_{r\in I_m^*}  \bb{Z}_p[[q^{1/p^n}-1]] \cdot q^{r/p^n}
\end{equation}
stable.
We also have for $m\geq n$ that  
\[
\bb{Z}_p[[q^{1/p^m}-1]]= \bigoplus\limits_{ i=0}^{p^{m-n}-1} \bb{Z}_p[[q^{1/p^n}-1]] \cdot  q^{i/p^{m}}. 
\]

Hence, for $n\geq 0$, the natural projection map $R_n\colon \bb{A}_{\inf}(\bb{Q}_p^{\cyc})\to \bb{Z}_p[[q^{1/p^n}-1]]$  of $\bb{Z}_p[[q^{1/p^n}-1]]$-modules is $\bb{Z}_p^{\times}$-equivariant and decomposes $\bb{A}_{\inf}(\bb{Q}_p)$ as a $\bb{Z}_p^{\times}$-equivariant direct sum of $\bb{Z}_p[[q^{1/p^n}-1]]$-modules
\begin{equation}\label{eqDirectSumSenTracesAinfn}
\bb{A}_{\inf}(\bb{Q}_p^{\cyc}) = Q_n\oplus \bb{Z}_p[[q^{1/p^n}-1]]
\end{equation}
where, for any $0\leq k\leq n$, we have as a description of $(p,q^{1/p^k}-1)$-complete  $\bb{Z}_p[[q^{1/p^k}-1]]$-modules
\begin{equation}\label{eqDescriptionDirectSummandSenTraces}
Q_n= \widehat{\bigoplus}_{r\in I_{n-k+1}} \bb{Z}_p[[q^{1/p^k}-1]]\cdot q^r. 
\end{equation}
The maps $R_n$ are called the \textit{normalized Sen--Tate traces} of $\bb{A}_{\inf}(\bb{Q}_p^{\cyc})$. We will need the following decomposition of $\bb{A}_{\inf}(\bb{Q}_p^{\cyc})$.

\begin{lemma}\label{LemmaDecompositionAinfLargeDirectSum}
For $n\geq 0$ we let $Q_{n}^*:= R_{n+1}(Q_n)$. Then we have 
\[
Q_{n}^*=\bigoplus_{i=1}^{p-1} \bb{Z}_p[[q^{1/p^n}-1]] \cdot q^{i/p^{n+1}}
\]
Moreover, set $Q_{-1}^*:=\bb{Z}_p[[q-1]]$, then   the natural maps of $(p,q-1)$-complete modules
\[
 \widehat{\bigoplus}_{n\geq -1} Q_{n}^* \to  \bb{A}_{\inf}(\bb{Q}_p^{\cyc}) \mbox{ and } \widehat{\bigoplus}_{n\geq m} Q_{n}^* \to Q_{m}
\]
are  isomorphisms for $m\geq 0$. 
\end{lemma}
\begin{proof}
  By definition, $R_n$ projects a power series $f=\sum_{r\in [0,1)\cap \bb{N}[\frac{1}{p}] } a_r q^{r}$ with $a_r\in \bb{Z}_p[[q-1]]$ to the power series $R_n(f)= \sum_{v_p(r)\geq -n} a_r q^r$.
  In particular, $Q_{n}\subset \bb{A}_{\inf}(\bb{Q}_p^{\cyc})$ consists of those  power series  of the form $\sum_{v_{p}(r)\leq -n-1} a_r q^r$.
  Therefore, we see that $R_{n+1}(Q_n)$ is the subspace of power series of the form $\sum_{v_p(r)= -n-1} a_r q^r$, which can be equally written as the subspace of power series $\sum_{v_p(r)=-1} b_r q^{r/p^n}$ with $b_r\in \bb{Z}_p[[q^{1/p^n}-1]]$. This proves the first claim about the description of $Q_n^*$.
 
  The second claim follows from the fact that any power series $f=\sum_{r} a_r q^r$ in $\bb{A}_{\inf}(\bb{Q}_p^{\cyc})$ can be written uniquely as the sum 
 \[
 f = a_{0} + \sum_{n \geq 0} \big( \sum_{v_p(r)=-n-1} a_rq^r \big)
 \]
 with $a_0\in \bb{Z}_p[[q-1]]= Q_{-1}^*$, and  $\sum_{v_p(r)=-n-1} a_rq^r \in Q_n^{*}$.  The same discussion works for $Q_m$. 
\end{proof}

The action of $\Gamma$ has the following behaviour on the elements $q^r$ in  $\bb{A}_{\inf}(\bb{Q}_p^{\cyc})$. For $k\geq 1$ we denote $\Gamma_{k}:= 1+p^k \bb{Z}_p$. 

\begin{lemma}\label{LemmaKeyComputationsSenTheory}
The following hold:

\begin{enumerate}

\item  For $m\geq n$ we have $(q^{1/p^m}-1)^{p^{m-n}}\equiv q^{1/p^n}-1 \mod p\bb{Z}_{p}[[q^{1/p^m}-1]]$. 

\item  For $a\in \bb{Z}_p^{\times}$ and $n\geq 0$, we have that $q^{a/p^n}-1 = u_{a,n}  (q^{1/p^n}-1)$ with $u_{a,n}\in \bb{Z}_p[[q^{1/p^n}-1]]^{\times}$. Moreover, $u_{a,n}\equiv a \mod q^{1/p^n}-1 $. 

\item Let $\gamma_k=1+ap^k$, then for $r\in \bb{N}[\frac{1}{p}]$ we have that 
\[
(\gamma_k-1)(q^{r})= q^{r}(q^{ap^kr}-1). 
\]
In particular, if $v_p(r)=-n$ and $(a,p)=1$, we have that 
\[
(\gamma_k-1)(q^{r}) = q^{r} (q^{1/p^{n-k}}-1) v_{k,r}
\]
with $v_{k,r}\in \bb{Z}_p[[q^{1/p^{n-k}}-1]]^{\times}$. 
\end{enumerate}
\end{lemma}
\begin{proof}
\begin{enumerate}

\item  This is obvious using the Frobenius homomorphism modulo $p$. 

\item We have that 
\[
q^{a/p^n}-1 = (1+(q^{1/p^n}-1))^a-1=  \sum_{s\geq 1} \binom{a}{s} (q^{1/p^n}-1)^s =  (q^{1/p^n}-1) \bigg(  \sum_{s\geq 0}  \binom{a}{s+1} (q^{1/p^n}-1)^s \bigg)
\]
in $\bb{Z}_p[[q^{1/p^n}-1]]$. The term  $u_{a,n}= \sum_{s\geq 0}  \binom{a}{s+1} (q^{1/p^n}-1)^s$ is congruent to $a$ modulo $(q^{1/p^n}-1)$, which makes it a unit in $\bb{Z}_p[[q^{1/p^n}-1]]$. 

\item By definition we have that 
\[
(\gamma_{k}-1)(q^r) = q^{r(1+ap^k)} -q^r = q^r(q^{rap^k}-1),
\]
which yields the desired formula. The assertion when $v_p(r)=-n$ follows from (2) as then we can write $rap^{k} = bp^{-{n-k}}$ for some $b\in \bb{Z}_p^{\times}$. \qedhere
\end{enumerate}
\end{proof}

Up to this point all the previous lemmas and power series computations are standard in the theory of Tate--Sen traces for the cyclotomic extension, e.g.\ as explained in \cite[Section III]{CherbonnierColmez}.
However, we are interested in localizing such power series descriptions on the Fargues--Fontaine curve of $\bb{Q}_p^{\cyc}$.

Consider the radius map $\kappa\colon \FFdisc{\bb{Q}_p^{\cyc}}^{\times} = (1+\widetilde{\bb{D}}_{\bb{Q}_p}^{<1,u})\backslash \{1\}\to (0,\infty)$ with respect to the pseudo-uniformizer $[\pi]$ with $\pi=\epsilon-1$.

Given $I=[a,b]\subset (0,\infty)$ an interval with rational ends we let $\FFdisc{\bb{Q}_p^{\cyc},I}\subset \FFdisc{\bb{Q}_p^{\cyc}}^{\times}$ denote the affinoid subspace whose functions are given by the Tate-rational localization $\bb{B}_{I}(\bb{Q}_p^{\cyc})$ with respect to the locus $\{|\pi|_x^{1/a}\leq |p|_x \leq |\pi|_x^{1/b}\}$ in $\Spa (\bb{A}_{\inf}(\bb{Q}_p))$.
Since $|\FFdisc{\bb{Q}_p^{\cyc}}^\times| = \varprojlim_{n} | 1+\bb{D}_{n,\bb{Q}_p}^{<1}\setminus \{1\}| $ (for the Berkovich but also the adic spectrum), there is some $N\geq 0$ and  a Tate-rational localization $U_{N,I}\subset 1+\bb{D}_{N,\bb{Q}_p}^{<1}$  such that $\FFdisc{\bb{Q}_p^{\cyc},I} = \FFdisc{\bb{Q}_p^{\cyc}}^{\times} \times_{  1+\bb{D}_{N,\bb{Q}_p}^{<1}} U_{N,I} $.  We note that, thanks to \Cref{LemmaKeyComputationsSenTheory} (2), the ideal generated by $\epsilon^a-1$  for $a\in \bb{Z}_p^{\times}$ in $\bb{F}_p[[\pi^{1/p^{\infty}}]]$ is independent of $a$, so that   the affinoid $\FFdisc{\bb{Q}_p^{\cyc},I}$ is stable under the action of $\bb{Z}_p^{\times}$. Thus, after enlarging $N$ if necessary, one can  choose $U_{N,I}$ also $\bb{Z}_p^{\times}$-stable.  For $n\geq N$ we denote $U_{n,I}:= (1+\bb{D}_{n,\bb{Q}_p}^{<1}) \times_{1+\bb{D}_{N,\bb{Q}_p}^{<1}} U_{N,I}$. 

To simplify our discussion, we will use the following notation: 

\begin{itemize}

\item   We let $X_{\infty}^u:= 1+\widetilde{\bb{D}}^{<1,u}_{\bb{Q}_p}$,  $X_n := 1+\bb{D}_{n,\bb{Q}_p}^{<1}$ and $X_{\infty}= \varprojlim_{n} X_n$. 

\item We write $X_{\infty}^{u,\times}:= X_{\infty}^{u} \backslash \{1\}$, $X_{\infty}^{\times}:= X_{\infty} \backslash \{1\}$ and $X_n^{\times}:= X_{n} \backslash \{1\}$.  

\item  We let $q_n:=q^{1/p^n}$ denote the natural coordinate  of $X_n$ which satisfies $|q_n-1|<1$.  

\item  For $I=[a,b]\subset (0,\infty)$ we let  $X_{\infty,I}\subset X_{\infty}^*$ denote the Tate-rational localization with respect to $\{|\pi|_x^{1/a}\leq  |p|_x\leq |\pi|_x^{1/b}\}$. 

\end{itemize}

\begin{remark}\label{RemarkUnionDescription}
It is not true that $X_{n+1}\to X_{n}$ sends $X_{n+1}^{\times}$ to $X_{n}^{\times}$, namely, the map is given by  $x\mapsto x^p$, and so the $p$-th roots of unity of the disc $X_{n+1}$ are all sent to $1\in X_n$. Nevertheless, we have that 
\[
X_{\infty}^{\times} =  \bigcup_{n} f_{n}^{-1}(X_{n}^{\times})
\]
where $f_n\colon X_{\infty}\to X_n$. 
\end{remark}

We want to have a reasonably explicit description of the affinoid $U_{N,I}\subset X_{N}$ in terms of the coordinates $\pi^{1/p^N}$, at least for $N\gg 0$.  This is done in the following lemma:

\begin{lemma}\label{LemmaExplicitDescription}
  Let $[a,b]\in [0,\infty)$ be an interval with rational ends.
  The following hold:

\begin{enumerate}

\item For $n>\log_p b$ the quotient $\frac{q_{n}-1}{[\pi]^{1/p^n}}$  is in  $\bb{B}_{I}(\bb{Q}_p^{\cyc})^{\circ,\times}$, i.e., it is a norm $1$ unit.   In particular we have that $\{|[\pi]|^{1/a} \leq  |p|\leq |[\pi]|^{1/b}\}= \{|q_n-1|^{p^n/a}\leq  |p|\leq |q_n-1|^{p^n/b}\}$ as Tate-rational subspaces of $\Spa (\bb{A}_{\inf}(\bb{Q}_p^{\cyc}))$.

\item  For $m\geq n> \log_p b$  we have that 
\[
(q_{m}-1)^{p^{m-n}}= (q_n-1) u
\]
with $u\in \bb{B}_{I}(\bb{Q}_p^{\cyc})^{\circ,\times}$. In particular, $(q_{m}-1)^{p^{m-n}}$ and $q_n-1$ define the same spectral norm normalized such that $|q_m-1|^{p^{m-n}}_x=|q_n-1|_x=2^{p^{-n}}$ for all $x\in \n{M}(\bb{B}_{I}(\bb{Q}_p^{\cyc}))$.

\item Let $[a,b]\subset (0,\infty)$,  for $n>\log_p b$ let $X_{n,I}$ denote the Tate-rational localization $\{|q_n-1|^{p^n/a}\leq  |p|\leq |q_n-1|^{p^n/b}\}$ of $X_n$. Then, we have that 
\[
X_{m,I}= X_m\times_{X_n} X_{n,I} \mbox{ and } X_{\infty,I}^u = X_{\infty}^u \times_{X_n} X_{n,I} 
\]
for $m\geq n$. We shall write $X_{\infty,I}:= \varprojlim_{m\geq n} X_{m,I}$.  

\end{enumerate}
\end{lemma}
\begin{proof}
\begin{enumerate}

\item By definition, we know that $q_n-1\equiv \epsilon^{1/p^n}-1 = \pi^{1/p^n}$ modulo $p\bb{A}_{\inf}(\bb{Q}_p^{\cyc})$.
  Thus, $q_{n}-1= [\pi]^{1/p^n} + p c$ with $c\in \bb{A}_{\inf}(\bb{Q}_p^{\cyc})$.
  This produces an equality of open subspaces in the adic space
\[
U:= \{|p|_x< |q_{n}-1|_x\}= \{|p|_x<|[\pi]|^{1/p^n}_x \}
\] 
as $[\pi]^{1/p^n}$ and $q_{n}-1$ differ by a unit of norm $1$ in both the $>|p|$ locus of $|q_n-1|$ and $|[\pi]^{1/p^n}|$.
Since $p^n>b$, we know that $|[\pi]|^{1/b}= |[\pi]^{1/p^n}|^{ p^n/b}< |[\pi]^{1/p^n}|$, this implies that  $\{ |\pi|^{1/a} \leq |p|\leq |\pi|^{1/b}\}\subset U$. Since in the locus $U$ we have that $q_{n}-1= [\pi]^{1/p^n} u$ with $u$ a unit of norm $1$, we deduce that 
\[
\{ |\pi|^{1/a} \leq |p|\leq |\pi|^{1/b}\} = \{|q_{n}-1|^{p^n/a} \leq |p| \leq |q_{n}-1|^{p^n/b} \}
\] 
proving part (1).

\item Thanks to \Cref{LemmaKeyComputationsSenTheory} we have for $m\geq n>\log_p b$ that $(q_m-1 )^{p^{m-n}}\equiv q_n-1\mod p \bb{Z}_p[[q_m-1]]$. Thus, we have equivalent loci
\[
\{|p|< |q_{m}-1|^{p^{m-n}}\} = \{|p|<|q_n-1|\}
\]
in $X_k$ for all $\infty \geq k\geq m$, and in those loci $(q_{m}-1)^{p^{m-n}}$ and $q_{n}-1$ differ by a unit of norm $1$.
In this locus we have that 
\begin{equation}\label{eqLocalizationRationalLocalizationsFiniteLevel}
\{ |q_m-1|^{p^m/a} \leq |p| \leq  |q_m-1|^{p^m/b} \} = \{ |q_n-1|^{p^n/a} \leq  |p|\leq |q_n-1|^{p^n/b}\}
\end{equation}
which implies part (2).

\item This follows formally from part (1) and (2), and the identification of Tate-rational localizations \eqref{eqLocalizationRationalLocalizationsFiniteLevel}. 
\qedhere
\end{enumerate}
\end{proof}

Let $I=[a,b]\subset (0,\infty)$ be a closed interval with rational ends as before and let $n>\log_p (b)$.
We will write $B_{I,n}:= \s{O}(X_{n,I})$ and $B_{I,\infty}:= \s{O}(X_{\infty}^u)=\bb{B}_{I}(\bb{Q}_p^{\cyc})$.
From \Cref{LemmaDecompositionAinfLargeDirectSum} we can deduce the following easy explicit description of the ring $B_{I,\infty}$ in terms of the $B_{I,n}$'s:

\begin{lemma}\label{LemmaDescriptionOnAffinoids}
Let $n>\log_p(b)$, then there is a natural $\bb{Z}_p^{\times}$-equivariant and $B_{I,n}$-linear decomposition
\[
B_{I,\infty} = \widehat{\bigoplus}_{r\in [0,1)\cap \bb{N}[\frac{1}{p}]} B_{I,n} \cdot q^{r}_n.
\]
\end{lemma}
\begin{proof}
This follows from \Cref{RemarkSenTracesCyclotomicTrace} and \Cref{LemmaExplicitDescription} by taking the base change along $\bb{Z}_p[[q_n-1]]\to B_{I,n}$. 
\end{proof}

It would be very convenient if the power series description of \Cref{LemmaDescriptionOnAffinoids} also hold for the power bounded rings.
Unfortunately this might not be the case due to the passage to Tate-rational localizations of $(p,q-1)$-complete rings.
In order to use \Cref{PropositionTechnicalLemmaTraces} we need to measure the difference between the rings $B_{I,\infty}^{\circ}$ and $\bb{A}_{\inf}(\bb{Q}_p^{\cyc})\otimes_{\bb{Z}_p[[q_n-1]]} B_{I,n}^{\circ}$.
Furthermore, even if we can appropriately control this difference, we also need to base change our rings to a suitable uniform Banach algebra $\Lambda$ with enough pseudo-uniformizers.
In order to have a reasonable control of both of these problems, we will need the following lemmas:
 
\begin{lemma}\label{LemmaAlgebraPresentation}
  Consider the given $\bb{Z}_p^{\times}$-action on $\bb{F}_p((\pi^{1/p^{\infty}}))$ with $\pi=\varepsilon-1$.
  Let $\varpi$ be a variable and consider the perfectoid field $\bb{F}_p((\varpi^{1/p^{\infty}}))$.
  Finally, let $R:= \bb{F}_p((\varpi^{1/p^{\infty}}))\big\langle \big(\frac{\pi}{\varpi} \big)^{\pm 1/p^{\infty}} \big\rangle$. The following holds:

\begin{enumerate}

\item  We have  $R=\bb{F}_p((\pi^{1/p^{\infty}})) \big\langle \big(\frac{\pi}{\varpi} \big)^{\pm 1/p^{\infty}} \big\rangle$.

\item We have $R^{\circ}= \bb{F}_p[[\varpi^{1/p^{\infty}}]] \big\langle \big(\frac{\pi}{\varpi} \big)^{\pm 1/p^{\infty}} \big\rangle$. 

\item The action of $\bb{Z}_p^{\times}$  on $\bb{F}_p((\pi^{1/p^{\infty}}))$ extends uniquely to a $\bb{F}_p((\varpi^{1/p^{\infty}}))$-linear action on $R$.

\end{enumerate}
\end{lemma}
\begin{proof}

Part (1) is clear as $R$ is the free perfectoid ring over $\bb{F}_p$ with two pseudo-uniformizers $\pi$ and $\varpi$ such that $\pi=u\varpi$ with $u$ a unit of norm $1$. Part (2) is also clear as $R$ is just a perfectoid torus relatively to either $\bb{F}_p((\varpi^{1/p^{\infty}}))$ or  $\bb{F}_p((\pi^{1/p^{\infty}}))$. Finally, part (3) is a consequence of the fact that for $a\in \bb{Z}_p^{\times}$ the element $\frac{\epsilon^a-1}{\epsilon-1}\in \bb{F}_p[[\pi]]^{\times}$ is a unit of norm $1$ which converges to $1$ as $a\to 1$. In particular, $a\cdot \frac{\pi}{\varpi}= \frac{a\cdot \pi}{\pi} (\frac{\pi}{\varpi})$  has norm $1$, and converges to $ \frac{\pi}{\varpi}$ as $a\to 1$. 
\end{proof}

From now on, we shall write $E:=\bb{F}_p((\varpi^{1/p^{\infty}}))$ and let $\Lambda_{I}:= \bb{B}_{I}(E)$.
By an abuse of notation we also write $\varpi$ for the Teichm\"uller lift of $\varpi$, and similarly for $\pi$.
The ring $\Lambda_I$  is a uniform Banach $\bb{Q}_p$-algebra  with pseudo-uniformizer $\varpi$ admitting all $p$-th power roots. 

\begin{lemma}\label{LemmaTechnicalBaseChange}
  Let $I:=[a,b]\subset (0,\infty)$ be a closed interval with rational ends and let  $n> \log_p(b)$.
  Let $R$ be as in \Cref{LemmaAlgebraPresentation}. The following hold:

\begin{enumerate}

\item The quotient $\frac{q_{n}-1}{\varpi^{1/p^n}}$  is a unit in $\bb{B}_{I}(R)^{\circ}$ congruent to $\frac{\pi^{1/p^n}}{\varpi^{1/p^n}}$ modulo $p/\varpi^{1/p^n}$.

\item For $r\in [0,1)\cap \bb{N}[\frac{1}{p}]$ let $r= a_r/p^{k_r}$ be the unique expression with $ 0 \leq a_r< p^{k_r}$ and $(a,p^{k_r})=1$ (so that if $a=0$ then $k_r=0$).  We have a natural isomorphism  of $\varpi$-complete algebras
\begin{equation}\label{eqBasisRingAdaptedPseudoUniformizers}
\bb{B}_{I}^{\circ}(R) = \widehat{\bigoplus}_{r\in [0,1)\cap \bb{N}[\frac{1}{p}]} \Lambda_{I}^{\circ}  \bigg\langle \big(\frac{q_n-1}{\varpi^{1/p^{n}}} \big)^{\pm 1} \bigg\rangle  \cdot \big(\frac{q_{n+k_r}-1}{\varpi^{1/p^{n+k_r}}}\big)^{a_r}.
\end{equation}

\end{enumerate}

\end{lemma}
\begin{proof}
  By construction, we know that $\frac{[\pi]}{[\varpi]}$ is a unit  in $\bb{A}_{\inf}(R)$ and therefore a unit of norm $1$ in $\bb{B}_{I}(R)$.
  On the other hand, by \Cref{LemmaExplicitDescription} (1), the quotient $\frac{q_n-1}{[\pi]^{1/p^n}}$ is a norm $1$ unit in $B_{I,\infty}$.
  This implies that $\frac{q_n-1}{[\varpi]^{1/p^n}}$ is a norm $1$ unit in $\bb{B}_I(R)$ as desired.
  We have $q_n-1\equiv \pi^{1/p^n}\mod p\bb{A}_{\inf}(\bb{Q}_p^{\cyc})$, i.e. $q_n-1= \pi^{1/p^n}+px$, with $x\in \bb{A}_\inf(\Q_p^\cyc)\subset \bb{B}_I(R)^\circ$. By the previous, both sides of this equality are divisible by $\varpi^{1/p^n}$ in $\bb{B}_I(R)^\circ$ (using the hypothesis $n>\log_p(b)$); dividing by $\varpi^{1/p^n}$ (note that $\mathbb{B}_I(R)^\circ$ is $\varpi^{1/p^n}$-torsion free), we deduce the desired congruence modulo $p/\varpi^{1/p^n}$. 
  This proves   part (1).

For part (2), by \Cref{LemmaAlgebraPresentation} (2) we have an isomorphism of $(p,[\varpi])$-complete algebras
\[
\bb{A}_{\inf}(R)= \bb{A}_{\inf}(E) \big\langle  \big(\frac{[\pi]}{[\varpi]}\big)^{\pm 1/p^{\infty}}\big\rangle.
\] 
Localizing in the radius $[a,b]$ (which is independent of $\pi$ or $\varpi$ as both have the same norm), we have that 
\[
\bb{B}_I(R)=\Lambda_{I}  \big\langle  \big(\frac{[\pi]}{[\varpi]}\big)^{\pm 1/p^{\infty}}\big\rangle
\]
and thus that 
\[
\bb{B}_I(R)^{\circ}= \Lambda^{\circ}_I \langle \big(\frac{[\pi]}{[\varpi]}\big)^{\pm 1/p^{\infty}}\big\rangle. 
\]
In particular,  $\bb{B}_I(R)^{\circ}$ has a Banach basis over  $\Lambda^{\circ}_I \langle (\frac{[\pi]}{[\varpi]})^{\pm 1/p^n} \rangle$ given by $(\frac{[\pi]}{[\varpi]})^{ r/p^n}$ with  $r\in [0,1)\cap \bb{N}[\frac{1}{p}]$.

The statement then follows from Nakayama's lemma, since by part  (1) we have that 
\[
(\frac{q_{n+k_r}-1}{\varpi^{1/p^{n+k_r}}} )^{a_r} \equiv (\frac{\pi^{1/p^{n+k_r}}}{\varpi^{1/p^{n+k_r}}})^{a_r} = \frac{\pi^{rp^{-n}}}{\varpi^{rp^{-n}}} \mod \frac{p}{\varpi^{1/p^n}} \bb{B}_{I}(R)^{\circ} 
\]
and since both rings in \eqref{eqBasisRingAdaptedPseudoUniformizers} are $p/\varpi^{1/p^n}$-adically complete. 
\end{proof}

We shall denote $Y_{\infty}^{\times,u}:= \FFdisc{R}^{\times}$ and for $I=[a,b]\subset (0,\infty)$ write $Y_{I,\infty}^u\subset Y_{\infty}^{\times,u}$ for the appropriate Tate-rational localization for the radius map defined by $\varpi$ (or equivalently by $\pi$).
For $n>\log_p b$,  we define $Y_{n}^{\times}\to X_{n}^{\times}$ to be the analytic adic space over $\FFdisc{E}^{\times}$ defined by the Tate-rational localization  $Y_{n}^{\times} \subset  X_n\times_{\bb{Q}_p} \FFdisc{E}^{\times} $ making $|q_n-1|=|\varpi^{1/p^n}|$.
Finally, we define the affinoid subspace $Y_{I,n}:= Y_n^{\times}\langle  |\varpi|^{1/b}\leq |p|\leq |\varpi|^{1/b} \rangle$.

By construction we have that $Y_{\infty}^{\times,u}= \bigcup_{I\subset (0,\infty)} Y_{I,\infty}^u$ and that $Y_{I,\infty}^u$ is the uniform completion of the limit $Y_{I,\infty}:=\varprojlim_{n>\log_p(b)} Y_{I,n}$. 

\begin{lemma}\label{LemmaShierkCoverDescentAction}
\begin{enumerate}
\item The  analytic $\bb{Z}_p^{\times}$-action on $X_{I,n}$ for $n>\log_p (b)$ extends uniquely to a $\Lambda_I$-linear analytic action on $Y_{I,n}$.
  Moreover, these actions are compatible for the tower $\{Y_{I,n}\}_{n>\log_p(b)}$ and produce a locally analytic action on the affinoid space $Y_{I,\infty}$. 

\item  The collection of maps $\bigsqcup_{I\subset (0,\infty)} Y_{I,\infty}\to X_{\infty}^{\times}$ form a $\bb{Z}_p^{\times, \la}$-equivariant $!$-cover of $X_{\infty}^{\times}$ with base change 
\begin{equation}\label{eqBaseChangeOverIIsUniformCompletion}
Y_{I,\infty}^{\times,u}\cong X^{\times,u}_{\infty}\times_{X_{\infty}^{\times}} Y_{I,\infty}.
\end{equation}

\end{enumerate}
\end{lemma}
\begin{proof}

\begin{enumerate}

\item We can equally describe $Y_{I,n}$ as the adic space $Z$ over $\GSpec \Lambda_I$ with pseudo-uniformizer $|q_n-1|$ satisfying $ |q_n-1|_x= |\varpi^{1/p^n}|_x $ for all $x\in Z$.
  Thus, $Y_{I,n}$ is affinoid with algebra of  functions $\Lambda_{I}\langle \frac{q_n-1}{\varpi^{1/p^n}} \rangle$.
  By \Cref{LemmaKeyComputationsSenTheory} (1), for $a\in \bb{Z}_p^{\times}$ 
\[
q^{a}_n-1 = u_{a,n}(q_n-1)
\]
with $u_{a,n}\equiv a \mod (q_n-1)\bb{Z}_p[[q_n-1]]$. In particular, the element $\frac{q_n^a-1}{\varpi^{1/p^n}}$ is defined in  $\s{O}(Y_I)^{\circ}$. 
On the other hand, by \Cref{LemmaKeyComputationsSenTheory} (3) we know that 
\[
\frac{q^a_n-1}{\varpi^{1/p^n}} - \frac{q_n-1}{\varpi^{1/p^n}} = \frac{q_n (q_{n}^{a-1}-1) }{\varpi^{1/p^n}}.
\]
If $a\to 1$ this element converges to zero in $ \s{O}(Y_I)^{\circ}$. This produces a continuous  action of  $\bb{Z}_p^{\times}$ on $\s{O}(Y_I)^{\circ}$ such that for all $k>0$ the action of $\bb{Z}_p^{\times}$ on  $\s{O}(Y_I)^{\circ}/\varpi^k$  factors through a finite quotient. By \Cref{LemmaCriterionv1} the action of $\bb{Z}_p^{\times}$ on   $\s{O}(Y_I)$ is analytic as desired.

\item The morphism $\bigsqcup_{I\subset (0,\infty)} Y_{I,\infty}\to  X_{\infty}^{\times}$ is clearly $\bb{Z}_p^{\times,\la}$-equivariant by construction. To see that it is a cover, note that we have by construction 
\[
Y_{I,\infty} = \varprojlim_{n>\log_p(b)} Y_{I,n}.
\]
Thus, it suffices to see that $Y_{I,n}\to X_{I,n}$ is descendable of  index of descendability independent of $n$. For this, we note that we have a factorization
\begin{equation}\label{eqCompositionMapsYnXnDescendable}
Y_{I,n}\to X_{I,n}\times_{\GSpec(\Q_p)} \bb{Q}_p\langle T^{\pm 1}_n\rangle \to X_{I,n}
\end{equation}
where $T_n= \frac{q_{n}-1}{\varpi^{1/p^n}}$ is a unit of norm $1$ in $Y_{I,n}$. The second map in \eqref{eqCompositionMapsYnXnDescendable}  is the projection onto $X_{I,n}$ and it admits a section and therefore it is descendable of index $\leq 1$.
The first map is the base change along the profinite \'etale map $\FFdisc{E}^{\times}\to \bb{D}^{<1,\times}_{\bb{Q}_p}$  with coordinate $\varpi^{1/p^n}$.
This map is affine and the map $\s{O}_{\bb{D}^{<1,\times}_{\bb{Q}_p}}\to \s{O}_{\FFdisc{E}^{\times}}$  of sheaves on $\bb{D}^{<1,\times}_{\bb{Q}_p}$ admits a section as  $\s{O}_{\bb{D}^{<1,\times}_{\bb{Q}_p}}$-modules, which makes it descendable of index $\leq 1$.
In total, we see that $Y_{I,n}\to X_{I,n}$ is desendable of index $\leq 2$, and by \cite[Proposition 2.7.2]{mann2022p} the map $Y_{I,\infty}\to X_{I,\infty}$ is descendable of index $\leq 4$.

Finally, the claim about the base change   \eqref{eqBaseChangeOverIIsUniformCompletion} follows from the fact that $Y_{I,\infty} = X_{\infty}^{\times}\times_{X_n^{\times}} Y_{I,n}$ so that 
\[
X_{\infty}^{\times,u}\times_{X^{\times}_{\infty}} Y_{I,u}= X_{\infty}^{\times,u} \times_{X_n^{\times}} Y_{I,n} \subset X_{\infty}^{\times,u}\times_{\GSpec(\Q_p)} \FFdisc{E,I}^{\times} 
\]
is the Tate-rational localization of the Fargues--Fontaine curve  of $\Spd E\times_{\Spd \bb{F}_p} \Spd \bb{Q}_p^{\cyc}$ to the locus where $\{|\varpi|^{1/a} \leq  |p|\leq |\varpi|^{1/b}\}$ and $\{|q_n-1|=|\pi|^{1/p^n}=|\varpi|^{1/p^n}\}$. This is precisely $Y_{I,\infty}^{\times,u}$.\qedhere
\end{enumerate}

\end{proof}

Thus, to prove \Cref{PropositionSenTheoryCyc} it suffices to show that for $I\subset (0,\infty)$ the  affinoid algebra $\s{O}(Y_{I,\infty}^u)^{\circ}$ satisfies the hypothesis of \Cref{PropositionTechnicalLemmaTraces}, and that its locally analytic vectors are given precisely by $\s{O}(Y_{I,\infty})=\varinjlim_{n>\log_p(b)} \s{O}(Y_{I,n})$.
In order to keep the same notation of our previous set up and that of \Cref{PropositionTechnicalLemmaTraces}, we will fix $I$ and denote   $\Lambda^{\circ}=\Lambda^{\circ}_I$,   $A_{\infty}:= \s{O}(Y_{I,\infty}^u)^{\circ}$, $A_n:=\s{O}(Y_{I,n})^{\circ}$.

 Recall that for $k\geq 1$ we denote $\Gamma_k:=1+p^k \bb{Z}_p$ and we set $\Gamma:=\Gamma_1$. 

\begin{lemma}\label{LemmaFinalComputationSenTheory}
Write $\bb{A}_{\inf}(\bb{Q}_p^{\cyc})= Q_n\bigoplus \bb{Z}_p[[q_n-1]]$ as a decomposition by Sen traces as in \eqref{eqDirectSumSenTracesAinfn}.  The following hold:

\begin{enumerate}

\item  Let $n>\log_p(b)$, then the map 
\begin{equation}\label{eqCokernelDifferentCoordinatesCyclotomicExt}
\bb{A}_{\inf}(\bb{Q}_p^{\cyc})\otimes_{\bb{Z}_{p}[[q_n-1]]} A_n\to A_{\infty} 
\end{equation}
is an inclusion with cokernel of $\varpi^{1/p^{n+1}}$-torsion.

\item Let $k\geq 1$ and let $\gamma_k\in \Gamma_{k}$ be a generator.
  Then the action of $1-\gamma_k$ on $A_{n}$ is divisible by $\varpi^{\frac{p^k-1}{p^n}}$.

\item  Let $k\geq 2$, then $H^*(\Gamma_k, Q_n\otimes_{\bb{Z}_p[[q_n-1]]} A_n)$ is concentrated in degree $1$ and  of $ \varpi^{1/p^{n+1-k}}$-torsion. 

\item Let $k\geq 2$, then  $H^{*}(\Gamma_k, A_{\infty}/A_n)$ is concentrated in degree $1$ and  of $\varpi^{1/p^{n+1-k}+1/p^{n+1}}$-torsion.

\end{enumerate}

In particular, the ring $A_{\infty}=(\varinjlim_{n} A_n)^{\wedge}_{\varpi}$ satisfies  \Cref{AssumptionSenTheory}. 
\end{lemma}
\begin{proof}

\begin{enumerate}

\item   By construction, we have that $A_n=\Lambda_{I}^{\circ}\big\langle \big(\frac{q_n-1}{\varpi^{1/p^n}} \big)^{\pm 1} \big\rangle$, thus, \Cref{LemmaTechnicalBaseChange} implies that $A_{\infty}$ has a Banach basis over $A_n$ with elements given by $e_{r}:=\big(\frac{q_{n+k_r}-1}{\varpi^{1/p^{n+k_r}}}\big)^{a_r}$ with $r=\frac{a_r}{p^{k_r}}\in [0,1)\cap \bb{N}[\frac{1}{p}]$ and $(a_r,p^{k_r})=1$.   In particular, for $v_p(r)<0$ we see that 
\[
\varpi^{r/p^{n}}e_r = (q_{n+k_r}-1)^{a_r}\in  \bb{A}_{\inf}(\bb{Q}_p^{\cyc})\otimes_{\bb{Z}_{p}[[q_n-1]]} A_n.
\]
This implies that the cokernel of \eqref{eqCokernelDifferentCoordinatesCyclotomicExt} is $\varpi^{1/p^{n+1}}$-torsion, proving part (1).

\item  Let $\gamma_k=1+p^k a\in \Gamma_k$ be an element with $(a,p)=1$.
  By \Cref{LemmaKeyComputationsSenTheory} (3)  we have that 
\[
(\gamma_k-1)\frac{q_n-1}{\varpi^{1/p^n}} = \frac{q_n(q_{n-k}^{a}-1) }{\varpi^{1/p^n}}=(q_n \frac{q^a_{n-k}-1}{q_{n-k}-1} ) \varpi^{\frac{p^k-1}{p^n}} \frac{q_{n-k}-1}{\varpi^{1/p^{n-k}}} = v \varpi^{\frac{p^k-1}{p^n}} \frac{q_{n-k}-1}{\varpi^{1/p^{n-k}}}. 
\]
By \Cref{LemmaKeyComputationsSenTheory} (2) the element $v$ is a unit in $A_n$, and by construction $ \frac{q_{n-k}-1}{\varpi^{1/p^{n-k}}}\in A_{n-k}\subset A_n$. This shows that the action of $\gamma_k$ on $A_{n}/\varpi^{\frac{p^k-1}{p^n}}$ is trivial, proving what we wanted.

\item  By \Cref{LemmaDecompositionAinfLargeDirectSum} we can write 
\[
Q_n = \widehat{\bigoplus}_{m\geq n} Q_{m}^* 
\] 
as  $(p,q_n-1)$-complete $\bb{Z}_p[[q_n-1]]$-modules, with 
\[
Q_m^{*}= \bigoplus_{i=1}^{p-1}  \bb{Z}_p[[q_m-1]] \cdot q_{m+1}^{i}.
\]
Thus, we have that 
\[
Q_n\otimes_{\bb{Z}_p[[q_n-1]] }A_n = \widehat{\bigoplus}_{m\geq n} (Q_m^*\otimes_{\bb{Z}_p[[q_n-1]] }A_n)
\]
and 
\[
Q_m^*\otimes_{\bb{Z}_p[[q_n-1]] }A_n =  \bigoplus_{i=1}^{p-1} (A_n[q_m]) \cdot q_{m+1}^{i}. 
\]

It suffices to show that for all $m\geq n$ the cohomology $H^1(\Gamma_k, Q_m^*\otimes_{\bb{Z}_p[[q_n-1]] }A_n)$ is $\varpi^{1/{p^{n+1-k}}}$-torsion.
We are going to prove a stronger statement:

\begin{itemize}

\item[i.] For $n\leq m<n+k$ the cohomology $H^*(\Gamma_k, Q_m^*\otimes_{\bb{Z}_p[[q_n-1]] }A_n)$ is $\varpi^{1/p^{m+1-k}}$-torsion, and hence $\varpi^{1/p^{n+1-k}}$-torsion. 

\item[ii.] For $m\geq n+k$ the cohomology  $H^*(\Gamma_k, Q_m^*\otimes_{\bb{Z}_p[[q_n-1]] }A_n)$ is $\varpi^{1/p^{m+1-k}+ 1/p^{n+1}}$-torsion, hence $\varpi^{2/p^{n+1}}$-torsion and $\varpi^{1/p^{n+1-k}}$-torsion (using $k\geq 1$).

\end{itemize}

As the modules $Q_m^*\otimes_{\bb{Z}_p[[q_n-1]] }A_n$ are $\varpi$-torsion free, the claims (i) and (ii) imply in particular that $\Gamma(\Gamma_k, Q_n\otimes_{\bb{Z}_p[[q_n-1]] } A_n)$ is concentrated in degree $1$.

Let $\gamma_k:=1+p^k$ be a generator of $\Gamma_k$.
As $\frac{1}{p^{m-k}}>\frac{1}{p^{m+1-k}}$,  by \Cref{LemmaTorsion}, it suffices to prove (i) after taking $\varpi^{\frac{1}{p^{m-k}}}$-quotients. For $m \leq n+k-1$ we have that 
\[
(\gamma_k-1)(q_m)= q_m(q_{m-k}-1) = \varpi^{1/p^{m-k}} q_m \frac{q_{m-k}-1}{ \varpi^{1/p^{m-k}}}
\]
with $\frac{q_{m-k}-1}{ \varpi^{1/p^{m-k}}}\in A_n^{\times}$.
This implies that $\gamma_{k}$ acts trivially on $A_n[q_m]/\varpi^{1/p^{m-k}}$.  On the other hand,  for $i=1,\ldots, p-1$,  we have that 
\[
(\gamma_k-1)(q_{m+1}^i) =  q_{m+1}^i( q_{m+1-k}^i -1) = q_{m+1}^i \varpi^{1/p^{m+1-k}} \frac{q_{m+1-k}^i -1}{\varpi^{1/p^{m+1-k}}}. 
\]
Since $m\leq n+k-1$ the element  $\frac{q_{m+1-k}^i -1}{\varpi^{1/p^{m+1-k}}}\in A_n^{\times}$ is a unit.
Thus, the previous calculation implies that the operator  $\gamma_k-1$ on $Q_m^*\otimes_{\bb{Z}_p[[q_n-1]]} A_n$ is divisible by $\varpi^{1/p^{m+1-k}}$, and Nakayama's lemma shows that $\varpi^{-1/p^{m+1-k}}(\gamma_k-1)$ is an isomorphism of   $Q_m^*\otimes_{\bb{Z}_p[[q_n-1]]} A_n=\bigoplus_{i=1}^{p-1} (A_n[q_m]) \cdot q_{m+1}^{i}$.
This yields the desired torsion statement in claim (i).

Let us now discuss the claim (ii).
Let $m\geq n+k$.
By part (1), the cokernel of $A_n[q_{m-k+1}]\to A_{m-k+1}$ is $\varpi^{1/p^{n+1}}$-torsion.
Thus, $\Gamma(\Gamma_k, Q_m^*\otimes_{\bb{Z}_p[[q_{m-k+1}-1]]} A_{n}[q_{m-k+1}])$ is a $\varpi^{1/p^{n+1}}$-retract of $\Gamma(\Gamma_k, Q_m^*\otimes_{\bb{Z}_p[[q_{m-k+1}-1]]} A_{m-k+1})$. By claim (i),  the cohomology groups of  $\Gamma(\Gamma_k, Q_m^*\otimes_{\bb{Z}_p[[q_{m-k+1}-1]]} A_{m-k+1})$  are $\varpi^{1/p^{m+1-k}}$-torsion, which implies that the cohomology groups of  $\Gamma(\Gamma_k, Q_m^*\otimes_{\bb{Z}_p[[q_{m-k}-1]]} A_{n}[q_{m-k}])$ are $\varpi^{1/p^{m+1-k}+ 1/p^{n+1}}$-torsion as desired. 

\item By part (1), $A_{\infty}/A_n$  is a $\varpi^{1/p^{n+1}}$-retract of   the base change $Q_n\otimes_{\bb{Z}_p[[q_n-1]]} A_n$.  Thus, $H^*(\Gamma_k, A_{\infty}/A_n)$ is a $\varpi^{1/p^{n+1}}$-retract of  $H^*(\Gamma_k,Q_n\otimes_{\bb{Z}_p[[q_n-1]]} A_n)$ which is concentrated in degree $1$ and of  $\varpi^{1/p^{n+1-k}}$-torsion by part (3). We deduce that  $H^*(\Gamma_k, A_{\infty}/A_n)$ is  concentrated in degree $1$ and of $\varpi^{1/p^{n+1-k}+1/p^{n+1}}$-torsion as desired.

\end{enumerate}

Finally, to see that the ring $A_{\infty}$ satisfies \Cref{AssumptionSenTheory},   it suffices to observe that 
\[
  1/p^{n+1-k}+1/p^{n+1}<  \frac{p^k-1}{p^n}=\frac{1}{p^{n-k}}-\frac{1}{p^n}
\]
for $k\geq 2$ fixed and $n\geq 1$, and that they converge to $0$ as $n\to \infty$. 
\end{proof}

\begin{proof}[Proof of \Cref{PropositionSenTheoryCyc}]
Thanks to \Cref{LemmaFinalComputationSenTheory}, the ring $A_{\infty}$ satisfies \Cref{AssumptionSenTheory} so that,  by \Cref{PropositionTechnicalLemmaTraces}, it satisfies $\epsilon$-\textbf{TORS}  uniformly on $\Gamma$ for some $\epsilon>0$. Thus, the affinoids $Y_{I,\infty}^u$ are strongly decompleted as in  \Cref{AssumptionKeyTheoremDecompletionPerfectComplexes} with decompletion given by $Y_{I,\infty}=\varprojlim_{n} Y_{I,n}$.  Then, thanks to \Cref{LemmaShierkCoverDescentAction}, we deduce that  the map 
\[
((1+\widetilde{\bb{D}}^{<1,u}_{\bb{Q}_p})\backslash \{1\})/ \bb{Q}_p^{\times,\cont}= X_{\infty}^{\times,u}/\bb{Q}_p^{\cont}  \to  X_{\infty}^{\times}/\bb{Q}_p^{\la}= ( (1+\widetilde{\bb{D}}^{<1}_{\bb{Q}_p})\backslash \{1\})/ \bb{Q}_p^{\times,\la}
\]
satisfies the hypothesis of \Cref{TheoremDecompletionPerfectComplexesII},   finishing the proof of \Cref{PropositionSenTheoryCyc}.
\end{proof}

\subsubsection{Sen theory: geometric case}\label{sss:SenTheotyGeometric}

Let $$S=(1+\widetilde{\bb{D}}^{<1,u}_{\bb{Q}_p})\backslash \{1\} = \FFdisc{\bb{Q}^{\cyc}_p}^{\times}$$ be the punctured Fargues--Fontaine disc of $\bb{Q}_p^{\cyc}$. We denote $q=[\epsilon]$ with $\epsilon=(\zeta_{p^n})_{n\geq 0}$, we write $\pi=\epsilon-1$ and we set $\varpi:=[\pi]$ as its Teichm\"uller lift. 

 For $n\geq 0$  let $\bb{D}_{S,n}= S\times_{\GSpec(\Q_p)} \GSpec (\bb{Q}_p\langle T^{1/p^n} \rangle)$ be the affinoid unit disc over $S$  in the variable $T^{1/p^n}$. 
 We let $\bb{D}_{S,\infty}=\varprojlim_n \bb{D}_{S,n}$ be its limit in Gelfand stacks, and let $\bb{D}_{S,\infty}^{u}$ be its uniform completion.
  Equivalently,  $\bb{D}_{S,\infty}^{u}$ is the punctured Fargues--Fontaine disc of $\bb{F}_p((\pi^{1/p^{\infty}}))\langle T^{1/p^{\flat,\infty}} \rangle$, where $T^{\flat}=(T^{1/p^n})_{n\geq 0}$. 
   We let $\Gamma\cong \bb{Z}_p$ be a $1$-parameter group with generator $\gamma$.  We endow $\bb{D}_{S,n}$ with the unique action of $\Gamma$ such that $\gamma(T^{a/p^n})= q^{a/p^n} T^{a/p^n}$ for $a/p^n\in p^{-n} \bb{N}$. 
    This endows $\bb{D}_{S,n}$ with a locally analytic action of $\Gamma$ which is analytic on affinoids, and linear over $S$. 
    By passing to limits and taking uniform completions, one obtains a $\Gamma^{\cont}$-action on the perfectoid disc  $\bb{D}_{S,\infty}^{u}$.
     Finally, we let $I=[a,b]\subset (0,\infty)$ be a closed interval with rational ends (with respect to the radius map defined by $\pi$) and consider the affinoid  subspaces $\bb{D}_{S_I,n}$ and $\bb{D}_{S_I,\infty}^{u}$. Thanks to \Cref{LemmaExplicitDescription}, for $n>\log_p(b)$ we have that $|q^{1/p^n}-1|_x=|\varpi^{1/p^n}|$ for all $x\in S_{I}$, in particular they define the same norm. 
 
 To ease notation we shall write $B_I:=\s{O}(S_I)$, $A_{n}:=B_{I}^{\circ}\langle T^{1/p^n}\rangle$ and $A_{\infty}= B_{I}^{\circ} \langle T^{1/p^{\infty}} \rangle$. We also  write $\Gamma_k=\gamma^{p^k\bb{Z}_p}$ for $k\geq 0$, and denote $\gamma_k:=\gamma^{p^k}$.   We have the following proposition:

\begin{proposition}\label{PropositionSenTheoryGeometricCase}
  With the previous notations, the $\Gamma$-action on $A_\infty$ satisfies \Cref{AssumptionSenTheory}. More precisely, let $k \geq 0$ and $n>\log_p(b)$, then  the following holds:
  
  \begin{enumerate}
  
\item The action of $\gamma_k-1$ on $A_{n}$ is trivial modulo $\varpi^{1/p^{n-k}}$. 

\item  The cohomology $\Gamma(\Gamma_{k}, A_{\infty}/A_n)$ is concentrated in degree $1$ and of $\varpi^{1/p^{n+1-k}}$-torsion.  
  
  \end{enumerate}
\end{proposition}

\begin{remark}\label{RemarkTorusCase}
Using \Cref{LemmaEtaleLocalizationDecomplete} one can deduce the $\epsilon$-\textbf{TORS} property uniformly on $\Gamma$ for the case of a torus $A_{\infty}=B_I^{\circ}\langle T^{\pm 1/p^{\infty}} \rangle$. One can also just directly check that the exact argument of \Cref{PropositionSenTheoryGeometricCase} applies in this case. 
\end{remark}

\begin{proof}[Proof of \Cref{PropositionSenTheoryGeometricCase}]
  Let $k\geq 0$.
  We have
  \begin{equation}
    \label{eq:2}
    (\gamma_k-1)T^r=(q^{p^kr}-1)T^r  
  \end{equation}
  for any $r\in \N[1/p]$.
 By  \Cref{LemmaKeyComputationsSenTheory} (2), for $v_p(r)=-n$ one has that $q^{r}-1= u_{r}(q^{1/p^n}-1)$ for $u$ a unit of norm $1$. In particular, 
\[
(\gamma_k-1) T^{r} = (q^{rp^k}-1) T^r = (q^{1/p^{n-k}}-1) u_{rp^k} T^{r}
\] 
with $u_{rp^k}\in B_I^{\circ,\times}$. 
In particular, the action of  $\gamma_k$ on $T^{1/p^n}$ is trivial modulo $\varpi^{1/p^{n-k}}$ for $n>\log_p(b)$. This yields part (1).

For part (2), we can write 
\[
A_{\infty} = \bigoplus_{r\in [0,1)\cap \bb{N}[\frac{1}{p}]} A_n \cdot T^{r},
\]
thus 
\[
A_{\infty}/A_n =  \bigoplus_{r\in (0,1)\cap \bb{N}[\frac{1}{p}]} A_n \cdot T^{r}. 
\]
One deduces that 
\[
\Gamma(\Gamma_k,A_{\infty}/A_n) = \bigoplus_{r\in (0,1)\cap \bb{N}[\frac{1}{p}]} A_n/ \varpi^{p^{v_p(r)+k}} \cdot T^{r},
\]
and thus that it is $\varpi^{1/p^{n+1-k}}$-torsion, since $v_p(r)\leq -n-1$. This finishes the proof of the proposition. 
\end{proof}

\begin{corollary}[Sen theory: geometric case] \label{cor:decompletion-geometric-case}
Let $d\geq 1$ and let $D$ be a perfectoid space over $\GSpec(\Q_p^{\cyc})$. 
 Consider the action of $G=\Z_p(1)^d$ on the $d$-dimensional perfectoid multiplicative group $\widetilde{\mathbb{G}}^{d,\diamond}_{m,D}$ with coordinates $T_1,\dots,T_d$, and its associated action at the level of perfect $\delta$-spaces
 $$
 X:=\bb{D}^{\times}_{\FF,\widetilde{\mathbb{G}}^{d,\diamond}_{m,D}}\cong \widetilde{\mathbb{G}}_{m,\Q_p}^{d,  u}\times_{\GSpec (\Q_p)} \bb{D}^{\times}_{\FF,D},
 $$
  where $\widetilde{\mathbb{G}}_{m,\Q_p}^{d, u}$ is the $d$-dimensional pre-perfectoid multiplicative group over $\Q_p$ with coordinates $[T_1^{\flat}], \dots, [T_d^\flat]$ where $T_i^{\flat}=(T_i^{1/p^{n}})_{n\geq 1}$, for $1\leq i\leq d$. Then $X$ is locally strongly decompleted.  Furthermore, 
  \[
X^{\la}= \widetilde{\bb{G}}_{m,\Q_p}^d \times_{\GSpec (\Q_p)} \mathbb{D}^{\times}_{\FF,D}.
\]
\end{corollary}
\begin{proof}
We deal, for simplicity of notation, with the case $d=1$, the general case reduces to this by iterating the argument. We have
\[
X \cong \widetilde{\mathbb{G}}_{m,\Q_p}^{ u}\times_{\GSpec (\Q_p)} \bb{D}^{\times}_{\FF,D},
\]
with trivial $\Z_p(1)$-action on the second factor. To prove that $X$ with its $G$-action is locally classically decompleted, it suffices, cf. \Cref{Lemma-Gelfand-connective cover} (1), to show that $\widetilde{\mathbb{G}}_{m,\Q_p}^{ u}\times_{\GSpec (\Q_p)} U$ is classically decompleted, for any affinoid subspace $U=\GSpec(\Lambda) \subset \bb{D}^{\times}_{\FF,D}$. In this case, it follows from \Cref{PropositionSenTheoryGeometricCase,PropositionTechnicalLemmaTraces}, which shows that
\[ (\widetilde{\mathbb{G}}_{m,\Q_p}^{ u}\times_{\GSpec (\Q_p)} U)^{G-\la} = \widetilde{\mathbb{G}}_{m,\Q_p} \times_{\GSpec (\Q_p)} U.
\]
 Moreover, this proposition also shows $X$ is locally strongly decompleted (by choosing $Y'$ to be the pullback to $X^\la$ of an affinoid cover of $\bb{D}_{\FF,D}^\times$).
\end{proof}

\subsubsection{Summary}\label{sss:SenTheoryMixedCase}

We summarize our results in one statement, for use in the next chapter. 
\begin{theorem}\label{thm:summary-main-results-sen}
\begin{enumerate}
\item    Consider the perfectoid space $\GSpec(\Q_p^\cyc)$ with its natural action of $G_{\rm arith}=\Z_p^{\times}$, and the induced action on the associated perfect $\delta$-space
$$X_{\rm arith}=\mathbb{D}^{\times}_{\FF, \Q_p^{\rm cyc}}= 1+\widetilde{\mathbb{D}}^{<1,u} \backslash \{1\}.
$$
   Then $X_{\rm arith}$ is locally classically decompleted with decompletion given by the Banach-Colmez space 
\[
X_{\rm arith}^{\la}=\BC(\s{O}(1))\backslash \{0\}=1+\widetilde{\mathbb{D}}^{<1}\setminus \{1\}.
\]  
\item Let $d\geq 1$ and let $D$ be a perfectoid space over $\GSpec(\Q_p^{\cyc})$. 
 Consider the action of $G_{\rm geo}=\Z_p(1)^d$ on the $d$-dimensional perfectoid multiplicative group $\widetilde{\mathbb{G}}^{d,\diamond}_{m,D}$, and its induced action at the level of perfect $\delta$-spaces
 $$
 X_{\rm geo}=\bb{D}^{\times}_{\FF,\widetilde{\mathbb{G}}^{d,\diamond}_{m,D}}\cong \widetilde{\mathbb{G}}_{m,\Q_p}^{d,  u}\times_{\GSpec (\Q_p)} \bb{D}^{\times}_{\FF,D},
 $$
  where $\widetilde{\mathbb{G}}_{m,\Q_p}^{d, u}$ is the $d$-dimensional pre-perfectoid multiplicative group over $\Q_p$.  Then $X_{\rm geo}$ is locally classically decompleted, with
  \[
X_{\rm geo}^{\la}= \widetilde{\bb{G}}_{m,\Q_p}^d \times_{\GSpec (\Q_p)} \mathbb{D}^{\times}_{\FF,D}.
\]
\item For $?\in \{\mathrm{arith},\mathrm{geo}\}$,   let $Y\to X_?$ be a Berkovich-\'etale $G_?^\cont$-equivariant morphism of Berkovich spaces. Then $Y$ is locally classically decompleted and the pullback  functor 
$$\Perf(Y^\la/G_?^\la)\to \Perf(Y/G_?^\cont)
$$
 is an equivalence.
 
 \item Let $d\geq 0$. Consider the perfectoid space $\widetilde{\mathbb{G}}_{m,\Q_p^{\cyc}}^{d,u}$, with the natural continuous action of $G_{\rm mixed}=G_{\rm arith}\ltimes G_{\rm geo}$ on it, and the induced action on its perfect $\delta$-space
 $$
 X_{\rm mixed}= \bb{D}^{\times}_{\FF,\widetilde{\mathbb{G}}^{d,u}_{m,\Q_{p}^{\cyc}}}= \widetilde{\bb{G}}_{m,\Q_p}^{d,u} \times_{\GSpec (\Q_p)} X_{\rm arith}.
 $$
  Then $X_{\rm mixed}$ is locally classically decompleted with 
\[
X_{\rm mixed}^{\la}=\widetilde{\bb{G}}_{m,\Q_p}^d \times_{\GSpec (\Q_p)} (\BC(\s{O}(1))\backslash \{0\}),
\]
where the multiplicative group has coordinate $[T_1^{\flat}],\dots,[T_d^\flat]$ as in (2). Moreover, for any $Y\to X_{\rm mixed}$ Berkovich-\'etale $G_{\rm mixed}^\cont$-equivariant morphism of Berkovich spaces, $Y$ is locally classically decompleted for the action of $G_{\rm mixed}^\cont$ and the pullback  functor $\Perf(Y^\la/G_{\rm mixed}^\la)\to \Perf(Y/G_{\rm mixed}^\cont)$ is an equivalence.

\end{enumerate}
\end{theorem}
\begin{proof}
\begin{enumerate}
\item This follows directly from \Cref{PropositionSenTheoryCyc}.
\item This follows directly from \Cref{cor:decompletion-geometric-case}.
\item By \Cref{PropositionSenTheoryCyc} or \Cref{cor:decompletion-geometric-case}, the hypothesis of \Cref{TheoremDecompletionPerfectComplexesII} hold. Applying this theorem and \Cref{sec:decompl-perf-compl-1-decompletion-etale-coverings}, we get the result.
\item We first observe that $X_{\rm mixed}$ is locally classically decompleted for the action of $G_{\rm geo}^\cont$, with decompletion
$$
X_{\rm mixed}^{G_{\rm geo}-\la} = \widetilde{\bb{G}}_{m,\Q_p}^d \times_{\GSpec (\Q_p)}  X_{\rm arith}.
$$
Indeed, this follows from the decomposition $X_{\rm mixed}= \widetilde{\bb{G}}_{m,\Q_p}^{d,u} \times_{\GSpec (\Q_p)} X_{\rm arith}$ of $X_{\rm mixed}$ as a product and \Cref{Lemma-Gelfand-connective cover} (5) (applied with $Y=X_{\rm arith}$, $G=G_{\rm geo}$ and $H$ the trivial group). Moreover, from this description we see (applying \Cref{Lemma-Gelfand-connective cover} (5) again, now with $G=G_{\rm arith}$, $H$ trivial) that $X_{\rm mixed}^{G_{\rm geo}-\la}$ is locally classically decompleted for the $G_{\rm arith}^\cont$-action, with decompletion $\widetilde{\bb{G}}_{m,\Q_p}^d \times_{\GSpec (\Q_p)} (\BC(\s{O}(1))\backslash \{0\})$. By \Cref{Lemma-Gelfand-connective cover} (6), we deduce that $X_{\rm mixed}$ is locally classically decompleted for the action of $G_{\rm mixed}$, with decompletion 
$$
X_{\rm mixed}^{\la}=\widetilde{\bb{G}}_{m,\Q_p}^d \times_{\GSpec (\Q_p)} (\BC(\s{O}(1))\backslash \{0\}).
$$
It remains to prove the equivalence of perfect complexes. Let $Y\to X_{\rm mixed}$ be a Berkovich-\'etale $G_?^\cont$-equivariant morphism of Berkovich spaces. By \Cref{Lemma-Gelfand-connective cover} (2), $Y$ is locally classically decompleted for the action of $G_{\rm mixed}^\cont$. Let $Y^\la$ be its decompletion. The morphism $Y/G_{\rm mixed}^\cont \to Y^\la/G_{\rm mixed}^\la$ factors as a composition
$$
Y/G_{\rm mixed}^\cont \to Y^{G_{\rm geo}-\la}/(G_{\rm arith}^\cont \ltimes G_{\rm geo}^\la) \to Y^\la/G_{\rm mixed}^\la.
$$
Pullback along the first map induces an equivalence on perfect complexes: indeed, this can be checked after pullback along $\GSpec(\Q_p) \to BG_{\rm arith}^\cont$ by \Cref{LemmaDescentCoefficientsPerfect} and then follows from part (2) (with $D=\GSpec(\Q_p^\cyc)$).  To finish, we prove that pullback along the second map also induces an equivalence on perfect complexes. Again by \Cref{LemmaDescentCoefficientsPerfect}, it suffices to prove that pullback along 
$$
Y^{G_{\rm geo}-\la}/G_{\rm arith}^\cont \to Y^\la/G_{\rm arith}^\la
$$
is an equivalence on perfect complexes. 
First assume that $Y=X_{\rm mixed}$. In this case, as observed above, 
$$
X_{\rm mixed}^{G_{\rm geo}-\la} = \widetilde{\bb{G}}_{m,\Q_p}^d \times_{\GSpec (\Q_p)}  X_{\rm arith},
$$
with trivial $G_{\rm arith}$-action on the first factor. Part (1) shows that $X_{\rm arith}$ is locally strongly decompleted for the action of $G_{\rm arith}^\cont$. By \Cref{CorollaryDescentCoefficients} (applied by taking $C$ to be rings of functions of affinoid subspaces of $\widetilde{\bb{G}}_{m,\Q_p}^d$), the same holds for the fiber product, hence $X_{\rm mixed}^{G_{\rm geo}-\la}$ is locally strongly decompleted for the action of $G_{\rm arith}^\cont$, giving the equivalence of perfect complexes by \Cref{TheoremDecompletionPerfectComplexesII}. In general, note that by \Cref{Lemma-Gelfand-connective cover} (2), $Y^{G_{\rm geo}-\la}$ is Berkovich-\'etale $G_{\rm arith}^\cont$-equivariant over $X_{\rm mixed}^{G_{\rm geo}-\la}$, and thus \Cref{sec:decompl-perf-compl-1-decompletion-etale-coverings} applies.
\qedhere
\end{enumerate}
\end{proof}

\subsection{Decompletions of torsors}
\label{sec:decompl-tors-1}
Let $X$ be a Berkovich space with an action by $G^{\rm cont}$, where $G$ is a $p$-adic Lie group. Assume that $X$ is locally classically decompleted for the action of $G^\cont$. We have the natural morphism
\[
  h\colon X/G^{\rm cont}\to X^\la/G^\la
\]
of Gelfand stacks.
For any Gelfand stack $Z$ over $X^\la/G^\la$ we get a natural map
\[
  \beta_{X,Z}\colon \mathrm{Map}_{X^\la/G^\la}(X^\la/G^\la,Z)\to \mathrm{Map}_{X^\la/G^\la}(X/G^{\rm cont},Z),\ f\mapsto f\circ h.
\]
We are interested in conditions on $Z$ that guarantee that  $\beta_{X,Z}$ is an equivalence. 

In order to have a good compatibility for \'etale localizations, we actually study a stronger property.
Let $\mathcal{C}$ be the Berkovich-\'etale site of $X/G^\cont$ given by Berkovich-\'etale maps $X'/G^\cont\to X/G^\cont$ with $X'\to X$ Berkovich-\'etale and $G^\cont$-equivariant (by \Cref{Lemma-Gelfand-connective cover} this is equivalent to the Berkovich-\'etale site $\n{C}^{\la}$  of $X^\la/G^\la$). Given a morphism $Z\to X^{\la}/G^{\la}$ we let $\n{F}_{Z}$ denote the \'etale sheaf on $\n{C}$ defined by 
\[
\n{F}_{Z}:= \ob{Map}_{X^{\la}/G^{\la}}(-,Z).
\]
Similarly, we let $\n{G}_Z$ denote the \'etale sheaf on $\n{C}^{\la}$ defined in the same way. Under the natural equivalence of sites $\n{C}\cong \n{C}^{\la}$, we have a natural morphism of \'etale sheaves
\[
\beta_{X,Z}^{\et}\colon \n{G}_{Z}\to \n{F}_{Z}
\] 
inducing $\beta_{X,Z}$ on global sections and we are interested in the question when $\beta_{X,Z}^{\et}$ is an equivalence of sheaves.

We have the following simple permanence properties:

\begin{lemma}\label{LemmaPermanencePropertiesBetaXZ}
Let $X$ be a Berkovich space, with a $G^\cont$-action, for $G$ a $p$-adic Lie group. Assume that $X$ is locally classically decompleted for this action.
The following hold:

\begin{enumerate}
\item The class of Gelfand stacks $Z$  over $X^{\la}/G^{\la}$ for which $\beta^{\et}_{X,Z}$ is an equivalence is stable under limits. 

\item Suppose that $Z\to Z^{\dR/(X^{\la}/G^{\la})}$ is an equivalence, then $\beta_{X,Z}^{\et}$ is an equivalence. 

\item Let $Z\to X^{\la}/G^{\la}$ be a Gelfand stack and let $Z'\to Z$ be a Berkovich-\'etale map.
  If $\beta_{X,Z}^{\et}$ is an equivalence, then $\beta_{X,Z'}^{\et}$ is an equivalence.
  Conversely, if $Z'\to Z$ is a Berkovich-\'etale cover and $\beta_{X,Z'}^{\et}$ is an equivalence, then $\beta_{X,Z}^{\et}$ is an equivalence.

\end{enumerate}
\end{lemma}
\begin{proof}
\begin{enumerate}

\item  This is obvious from the definition as the formation of $\n{F}_Z$ and $\n{G}_Z$ commutes with limits in $Z$.  

\item Since $X$ is classically decompleted, we have that $X^{\diamond}= X^{\la,\diamond}$ as the uniform completions of $X$ and $X^{\la}$ agree locally on a strict cover.
  Since $G^{\ob{cont},\diamond}=G^{\la,\diamond}= \underline{G}$, we deduce that
\[
(X/G^{\cont})^{\dR}= (X^{\la}/G^{\la})^{\dR}
\]
and therefore that the relative de Rham stack of $X/G^{\cont}\to X^{\la}/G^{\la}$ is an isomorphism.
Thus, if $X'/G^{\cont}\to X/G^{\cont}$ is a Berkovich-\'etale map with decompletion $X^{\prime,\la}/G^{\la}$, we have that 
\[
\begin{aligned}
  \ob{Map}_{X^{\la}/G^{\la}}(X'/G^{\cont}, Z^{\dR/(X^{\la}/G^{\la})})  
  = & \ob{Map}_{(X^\la/G^{\la})^{\dR}}(X'/G^{\cont}, Z^{\dR}) \\
  =& \ob{Map}_{(X^{\la}/G^{\la})^{\dR}}((X'/G^{\cont})^{\dR}, Z^{\dR}) \\
  = &\ob{Map}_{(X^{\la}/G^{\la})^{\dR}}((X^{\prime,\la}/G^{\la})^{\dR}, Z^{\dR}) \\
  =  &\ob{Map}_{(X^{\la}/G^{\la})^{\dR}}(X^{\prime,\la}/G^{\la}, Z^{\dR}) \\
= & \ob{Map}_{X^{\la}/G^{\la}}(X^{\prime,\la}/G^{\la}, Z^{\dR/(X^{\la}/G^{\la})}), 
\end{aligned}
\]
which proves part (2).

\item  Suppose that $Z'\to Z$ is a Berkovich  \'etale map over $X^{\la}/G^{\la}$.
  Then we have an equivalence $Z'= Z^{\prime,\dR}\times_{Z^{\dR}} Z$ which yields
\[
\n{F}_{Z'}= \n{F}_{Z^{\prime,\dR/(X^{\la}/G^{\la})}}\times_{\n{F}_{Z^{\dR/(X^{\la}/G^{\la})}}} \n{F}_{Z}
\] 
and similarly for $\n{G}_{Z'}$. By part (2), we have natural equivalences 
\[
\n{G}_{Z^{\prime,\dR/(X^{\la}/G^{\la})}} =\n{F}_{Z^{\prime,\dR/(X^{\la}/G^{\la})}} \mbox{ and } \n{G}_{Z^{\dR/(X^{\la}/G^{\la})}} =\n{F}_{Z^{\dR/(X^{\la}/G^{\la})}}.
\]
Therefore, we have a Cartesian diagram
\[
\begin{tikzcd}
\n{G}_{Z'} \ar[r,"{\beta^{\et}_{X,Z'}}"] \ar[d] & \n{F}_{Z'} \ar[d] \\
\n{G}_Z \ar[r,"{\beta^{\et}_{X,Z}}"] & \n{F}_Z.
\end{tikzcd}
\]
This implies that if $\beta^{\et}_{X,Z}$ is an equivalence, then so is $\beta_{X,Z'}^{\et}$.
Conversely, if $Z'\to Z$ is a Berkovich-\'etale cover, then the map $\n{F}_{Z'}\to \n{F}_Z$ is a surjection  of \'etale sheaves on $\n{C}$. Thus, by descent, if $\beta_{X,Z'}^{\et}$ is an equivalence, then so is $\beta^{\et}_{X,Z}$ as wanted. \qedhere

\end{enumerate}
\end{proof}

Under perfectoidness assumptions, one can deduce slightly stronger statements.

\begin{lemma}
  \label{sec:decompl-tors-3-decompletion-of-torsors-perfectoid-case}
  Let $X$ be a Berkovich space, with a $G^\cont$-action, for $G$ a $p$-adic Lie group. Assume that $X$ is locally classically decompleted for this action and that the natural map $\widehat{X^{\diamond}}\to X$ is an isomorphism, e.g., $X$ is perfectoid.
  Then for any Gelfand stack $Z$ with $Z':=X^\la/G^\la\times Z$ such that $\beta_{X,Z'}$ is an isomorphism, the natural map
  \[
    \mathrm{Map}_{\GelfStk}(X^\la/G^\la,Z)\to \mathrm{Map}_{\ArcStk_{\Q_p}}(X^\diamond/\underline{G},Z^\diamond)
  \]
  induced by $(-)^\diamond$ is an equivalence. 
\end{lemma}
\begin{proof}
  This follows from the adjunction between $\widehat{(-)}$ and $(-)^\diamond$.
  Indeed, we get
  \[
    \mathrm{Maps}_{\GelfStk}(X^\la/G^\la,Z)\cong \mathrm{Maps}_{\GelfStk}(X/G^{\rm cont},Z)\cong \mathrm{Maps}_{\ArcStk_{\Q_p}}(X^\diamond/\underline{G},Z^\diamond),
  \]
  where the first isomorphism holds by our assumption on $\beta_{X,Z'}$, and the second by our assumption on $X$ and the fact that $\widehat{(X^{\diamond}/G^{\cont})}= \widehat{X^{\diamond}}/G^{\cont}$.
\end{proof}

\begin{remark}
Let $Y$ be a Gelfand stack and $\Pi$ a locally profinite group, we note that the adjunction between de Rham stacks and perfectoidization yields natural equivalences
\[
\ob{Map}_{\cat{GelfStk}}(Y, B\Pi^{\sm}) =  \ob{Map}_{\cat{arcStk}}(Y^{\diamond}, B\underline{\Pi}) = \ob{Map}_{\cat{GelfStk}}(\widehat{Y^{\diamond}}, B\Pi^{\cont})
\]
 where in the last equality we used that $ \widehat{B\Pi}=  B\Pi^{\cont}\times_{\GSpec(\Q_p)} \widehat{\bb{Q}_p^{\diamond}}$.  In particular, any $\Pi^{\cont}$-torsor of $Y$ has a natural \textit{decompletion} to a $\Pi^{\sm}$-torsor, such decompletion makes sense for arbitrary $\Pi$ and $Y$, and hence it is not given by taking locally analytic vectors in any sense. 
\end{remark} 

\begin{remark}
  \label{sec:decompl-tors-4-berkovich-space-over-decompletion}
  Assume that $W\to X^\la$ is a $G^\la$-equivariant morphism with $W$ a Berkovich space.
  Set $Y:=X\times_{X^\la} W$ with its $G^\cont$-action.
  Then $Y$ admits a classical decompletion for its $G^\cont$-action with $Y^\la=W$.
  This follows from \Cref{Lemma-Gelfand-connective cover} (4).
  Alternatively, we have to see that $h_{Y,\ast}(\s{O})\in \ob{D}(|Y|_\Betti/G^\la)$ is connective for $h_Y\colon Y/G^\cont\to |Y|_\Betti/G^\la$.
  Using $|W|=|Y|$ (as can be checked on arc-stacks as $X^\diamond=X^{\la,\diamond}$), we have the following commutative diagram, whose left square is cartesian (but in general not the right one!):
 \[\begin{tikzcd}
	{Y/G^\cont} \ar[r,"g'"] \ar[rr,"h_Y", bend left] \ar[d] & {W/G^\la}  \ar[r ,"f'"] \ar[d] & {|Y|_\Betti/G^\la} \ar[d] \\
	{X/G^\cont} \ar[r,"g"] \ar[rr,"h_X"', bend right]& {X^\la/G^\la}  \ar[r,"f"]& {|X|_\Betti/G^\la}
      \end{tikzcd}\]
      
    By definition of $X^\la$, we have $g_\ast(\s{O})=\s{O}$ and by primness of $g$ we can deduce $g'_\ast(\s{O})=\s{O}$.
    The morphism $f'$ is the relative spectrum of a connective Gelfand ring (because $W$ is a Berkovich space), and hence $f_\ast^\prime$ preserves connective objects.
    We can conclude that $h_{Y,\ast}(\s{O})=f'_\ast(\s{O})$ is connective, and in fact given by $\s{O}_W\in \ob{D}(|W|_\Betti/G^\la)$.
    This shows that $Y$ admits a classical decompletion by $W$.

    As a concrete example, $W$ could be the geometric $G^\la$-equivariant vector bundle over $X^\la$ for some vector bundle $\mathcal{E}$ on $X^\la/G^\la$, and then $Y$ is given by the geometric vector bundle associated with the pullback of $\mathcal{E}$ to $X/G^\cont$. Another example could be an analytic $G^{\la}$-equivariant $H$-torsor over $X^{\la}$ where $H$ is an affine algebraic group over $\bb{Q}_p$.  
\end{remark}

It is natural to ask about the behaviour of locally analytic decompletions with respect to two commuting actions of $p$-adic Lie groups.  To address this question we need some preliminary lemmas.

\begin{lemma}\label{LemmaAffinoid}
Let $S$ be a Gelfand stack and let $\n{H}\to S$ be a group in Gelfand stacks over $S$, which $!$-locally on $S$ is induced by a morphism of Gelfand rings.  
Let $Y\to S$ be a Gelfand stack over $S$ and $f\colon \n{T}\to Y$ an $\n{H}$-torsor.
Then $f$ is $!$-locally on the target a morphism of Gelfand rings.
 \end{lemma}
\begin{proof}
This follows from the fact that torsors are $!$-locally trivial. 
\end{proof}

\begin{remark}
  \label{sec:decompl-tors-2-question-on-representability-of-g-la-torsors}
  In \Cref{LemmaAffinoid} assume that $\mathcal{H}=G^\la$ for a compact $p$-adic Lie group and that $S$ is a Berkovich space.
  Then each $G^\la$-torsor $Y\to S$ is $!$-locally a morphism of Gelfand rings, but it is unfortunately not clear that $Y$ is a Berkovich space.
  For example, if $G=\mathrm{GL}_n(\Z_p)$, then $Z:=Y\times^{\mathrm{GL}_n(\Z_p)^\la} \mathrm{GL}_n$ is a Berkovich space and $Z$ identifies with the overconvergent neighborhood of the closed subspace $|Y|\subseteq Z$, but it is is not clear that this overconvergent neighborhood is a Berkovich space.
  
\end{remark}

\begin{remark}\label{RemarkaffinenessTorsor}
  By \Cref{LemmaAffinoid}, if $H$ is a compact $p$-adic Lie group and $Y$ a Gelfand stack, we see that $H^{\cont}$ and $H^{\la}$-torsors over $Y$ are $!$-locally on the target induced by morphisms of Gelfand rings.
  Moreover, if $\n{T}\to Y$ is a $H^{\cont}$-torsor we can form its \textit{stacky locally analytic decompletion} $\n{T}^{\ob{stla}}$  induced by the pushout along the morphism of groups $H^{\ob{cont}}\to H^{\la}$.
  Geometrically,  this decompletion is given by 
\[
\n{T}^{\ob{stla}}=(\n{T}\times H^{\la})/H^{\cont},
\]
where $H^{\cont}$ acts antidiagonally via its actions on $\n{T}$ and $H^{\la}$.
Since $\n{T}^{\ob{stla}}$ is $!$-locally over $Y$ induced by a morphism of Gelfand rings by \Cref{LemmaAffinoid} (being an $H^\la$-torsor over $Y$),  we deduce by $\ast$-pushforward to $Y$ that (as objects in $\ob{D}(Y)$) we have isomorphisms 
\[
\s{O}(\n{T}^{\ob{stla}}) = \Gamma(H^{\cont}, \s{O}(\n{T})\otimes C^{\la}(H,\bb{Q}_p))= \s{O}(\n{T})^{H-\la}.
\]
That is, the (derived) ring of functions of $\n{T}^{\ob{stla}}$ (over $Y$) is nothing but the  locally analytic functions of $\n{T}$ with respect to the action of $H$. In this way, we can think of $\n{T}^{\ob{stla}}$ as the decompletion of $\n{T}$ by $H$-locally analytic vectors, even though they are not stacks associated to Berkovich spaces as in \Cref{Definition-Classically-decompleted}. Of course, the price to pay to allow this level of generality is that $\n{T}$ and $\n{T}^{\ob{stla}}$ are only Gelfand stacks, so that a priori we do not have a precise control on their geometry via their underlying Berkovich spaces. 
\end{remark}

\begin{lemma}\label{LemmaHlaTorsorDecompletion}
  Let $X$ be a Berkovich space endowed with an action of a $p$-adic Lie group $G^{\cont}$. Suppose that $X$ is locally classically decompleted for the action of $G$, with decompletion $Y^\la$.
  Let $H$ be a compact $p$-adic Lie group and let  $\n{T}\to X$ be a $G^{\cont}$-equivariant $H^{\la}$-torsor.
  Assume that $\beta_{X,Z}$ is an equivalence for $Z=BH^\la\times X^\la$.
  Let $\n{T}^{\ob{stla}}\to X^{\la}$ be the $G^{\la}$-equivariant $H^{\la}$-torsor over $X^{\la}$ associated to $\n{T}$ via $\beta_{X,Z}$.
  Then, as sheaves on the topological space $|X|$, we have a natural equivalence 
\[
\s{O}(\n{T})^{G-\la}= \s{O}(\n{T}^{\ob{stla}}).
\]
\end{lemma}
\begin{proof}
We have that $\n{T}=\n{T}^{\ob{stla}}\times_{X^{\la}} X$ as $H^{\la}$-torsors. Thus, we have a diagram with pullback squares 
\[
\begin{tikzcd}
 \n{T} /G^{\cont}\ar[r] \ar[d] & \n{T}^{\ob{stla}}/G^{\la } \ar[d]  \ar[r] & \GSpec (\bb{Q}_p )\ar[d] \\
X/G^{\cont} \ar[r] & X^{\la}/G^{\la} \ar[r]  & BH^{\la}. 
\end{tikzcd}
\]
The lemma follows formally from the fact that $ \n{T} /G^{\cont} \to  X/G^{\cont}$ and $\n{T}^{\ob{stla}}/G^{\la}\to X^{\la}/G^{\la}$ are $!$-locally induced by morphisms of Gelfand rings by \Cref{LemmaAffinoid}, and  that $f\colon X/G^{\cont} \to X^{\la}/G^{\la}$ is $0$-proper with $f_* 1= 1$. 
\end{proof}

\begin{proposition}\label{PropCommutingActionOnBerkovichSpace}
  Let $G$, $H$ be two $p$-adic Lie groups. 
  Let $Z$ be a Gelfand stack endowed with an action of $(G\times H)^\cont$.
  Denote $Y:=Z/G^{\cont}$ and assume the following:
  \begin{enumerate}
  \item $Y$ is a Berkovich space which is locally classically decompleted for the $H^\cont$-action, with decompletion $Y^{H-\la}$,
  \item $\beta_{Y,BG^\la}$ is an equivalence, in particular, each $H^\cont$-equivariant $G^\la$-torsor $\mathcal{P}$ over $Y$ decompletes to a $H^\la$-equivariant $G^\la$-torsor over $Y^{H-\la}$.
  \end{enumerate}
  Then the following holds:
  \[
    Y^{H-\la}/H^\la\cong \underline{\GSpec}_{|Z|_\Betti}(\mathscr{O}_{Z}^{(G\times H)-\la})/(G\times H)^\la.
  \]
  In particular, if $X:=Z/H^\cont$ satisfies assumptions (1) and (2) as well (with the roles of $G$ and $H$ swapped), then $Y^{H-\la}/H^\la\cong X^{G-\la}/G^\la$.

\end{proposition}
\begin{proof}
We can assume without loss of generality that both $G$ and $H$ are compact $p$-adic Lie groups. 

We note that $Z$ is a $H^\cont$-equivariant $G^\cont$-torsor over $Y$, and that therefore $\mathcal{P}:=Z\times^{G^\cont}G^\la$ is a $H^\cont$-equivariant $G^\la$-torsor over $Y$.
By our assumption, there exists a $H^\la$-equivariant $G^\la$-torsor $\mathcal{T}$ over $Y^{H-\la}$ with an $H^\cont$-equivariant isomorphism $\mathcal{P}\cong \mathcal{T}\times_{Y^{H-\la}} Y$.

By \Cref{LemmaAffinoid}, we know that
\[
  \mathcal{T}\cong \underline{\GSpec}_{Y^{H-\la}}(\mathcal{A})
\]
for a $H^\la$-equivariant commutative algebra in $\mathcal{A}\in\ob{D}(Y^{H-\la})$, which $!$-locally on $Y^{H-\la}$ is represented by a Gelfand ring.
Let $f\colon Y\to Y^{H-\la}$ be the natural $H^\cont$-equivariant morphism.
Then
\[
  \mathcal{T}\times_{Y^{H-\la}}Y\cong \underline{\GSpec}_{Y}(f^\ast \mathcal{A}),
\]
and we deduce from fully faithfulness in \Cref{TheoDecompletion} (1) that $\mathcal{A}\cong (f^\ast \mathcal{A})^{H-\la}$.
By \Cref{RemarkaffinenessTorsor} we see that $f^\ast \mathcal{A}$ is $G^\la\times H^\cont$-equivariant isomorphic to $\mathcal{B}^{G-\la}$, where $\mathcal{B}=g_\ast(\mathscr{O}_Z)\in \ob{D}(Y)$ for $g\colon Z\to Y$ the natural map.
Using \Cref{LemmaAffinoid}, we have $Z\cong \underline{\GSpec}_{Y}(\mathcal{B})$.
But as $Z$ is a Berkovich space we have $Z\cong \underline{\GSpec}_{|Z|_\Betti}(a_\ast \mathscr{O}_Z)$ for $a\colon Z\to |Z|_\Betti$ the natural map (cf.\ \Cref{sec:deriv-berk-spac-berkovich-spaces-as-relative-spectrum}).
This implies the claim, using \Cref{LemmaDecompletionsSemidirectProducrs} to see that $(-)^{(G\times H)-\la}\cong ((-)^{G-\la})^{H-\la}$.
For the ``in particular statement'' note that under the given assumptions on $X$ and $Y$, we get the isomorphisms
\[
  Y^{H-\la}/H^\la\cong \underline{\GSpec}_{|Z|_\Betti}(\mathscr{O}_{Z}^{(G\times H)-\la})\cong X^{G-\la}/G^\la
\]
as the term in the middle is completely symmetric.
\end{proof}

Finally, the next proposition proves several useful instances  where $\beta^{\et}_{X,Z}$ is an equivalence. Given an algebraic stack $\n{T}$ over $\Q_p$, we let $\n{T}^\an$ be its associated analytification as a Gelfand stack over $\Q_p$, see \Cref{sec:analyt-tann-stacks}. 

\begin{proposition}
  \label{sec:decompl-tors-2-beta-x-z-an-equivalence}
 Let $X$ be a Berkovich space, with a $G^\cont$-action, for $G$ a $p$-adic Lie group. Assume that $X$ is  locally classically decompleted for the action of $G^\cont$, and for any $X'\to X$ Berkovich-\'etale and $G^{\cont}$-equivariant morphism, perfect complexes on $X'/G^\cont$ admit a locally analytic decompletion.  The map $\beta_{X,Z}^{\et}$ is an equivalence in the following cases: 
   \begin{enumerate}
    \item $Z=\mathcal{T}^\an\times_{\GSpec(\Q_p)} X^\la/G^\la$, where $\mathcal{T}$ is a perfect strongly Tannakian stack over $\Spec(\Q_p)$, as in \Cref{sec:analyt-tann-stacks} (for example, the classifying stack of an affine algebraic group $H$ over $\Q_p$),
  \item $Z=BU\times X^\la/G^\la$, where $U$ is an open subgroup of $H(\Q_p)^\la$ for $H$ an affine algebraic group over $\Q_p$, and $H(\Q_p)^\la=H^\an\times_{H^{\an,\dR}}H(\Q_p)^\dR$.
  \item $Z=B\Gamma^\la \times X^\la/G^\la$, where $\Gamma$ is a $p$-adic Lie group.
    \end{enumerate}
\end{proposition}

For example, the assumption of the proposition holds if $X$ is locally strongly decompleted (using \Cref{sec:decompl-perf-compl-1-decompletion-etale-coverings}).

\begin{proof}

\begin{enumerate}

\item Recall from \Cref{PropositionFunctorOfPointsAnalytification} that for any Gelfand stack $Y$, morphisms $Y\to \n{T}^\an$ are the same  as symmetric monoidal functors from perfect complexes on $\n{T}$ to perfect complexes on $Y$, sending connective pseudo-coherent complexes in the ind-extension to connective objects. Since by assumption, pullback induces a symmetric monoidal equivalence between perfect complexes on $X^{\prime \la}/G^\la$ and $X'/G^\cont$ for any $X'\to X$ Berkovich-\'etale, we get the result (the assumption on preservation of connectivity can be checked using \cite[Lemma 3.3.3]{drff_stacks}).

\item   The second statement follows from the first and \Cref{LemmaPermanencePropertiesBetaXZ}: because $U\subseteq H(\Q_p)^\la$ is open, the map $BU\to H(\Q_p)^\la$ is a Berkovich-\'etale cover. Hence, by  \Cref{LemmaPermanencePropertiesBetaXZ} (3), it suffices to prove the claim where $U= H(\Q_p)^{\la}$. In that case, we have a pullback diagram 
\[
B H(\Q_p)^{\la} = BH(\Q_p)^{\sm} \times_{ BH^{\an,\dR}} BH^{\an},
\]
where we notice that $BH(\Q_p)^{\sm} = (BH(\Q_p)^{\la})^{\dR}$. Then,  by part (1) the map $\beta_{X,Z}^{\et}$ is an isomorphism for $Z=BH^{\an}$, and by \Cref{LemmaPermanencePropertiesBetaXZ} (2) it is an equivalence for the stacks  $Z=BH(\Q_p)^{\sm},BH^{\an,\dR}$. By \Cref{LemmaPermanencePropertiesBetaXZ} (2) it is an equivalence for $B H(\Q_p)^{\la}$ as desired.

\item  For (3), we first need to make some reductions.
  Given $\Gamma_0\subset \Gamma$ an open  subgroup,  the map $B\Gamma_0\to B\Gamma$ is a Berkovich  \'etale  cover.
  Thus,   thanks to \Cref{LemmaPermanencePropertiesBetaXZ} (3), we are allowed to replace $\Gamma$ by $\Gamma_0$ and suppose that $\Gamma$ is a compact $p$-adic Lie group, as small as necessary.
  Let $\f{g}=\ob{Lie}(\Gamma)$ denote the Lie algebra of $\Gamma$.
  Taking $\n{L}\subset \f{g}$ a free $\bb{Z}_p$-submodule stable under the Lie bracket with $[\n{L},\n{L}]\subset p \n{L}$, we can use the Baker-Campbell-Hausdorﬀ formula to define an affinoid rigid analytic group $\mathbb{G}(\n{L})$  as in \cite[Part II. Ch V.4]{MR1176100}.
  We denote $\exp(\n{L}):=\mathbb{G}(\n{L})(\bb{Q}_p)$, and call it the exponential group of $\n{L}$. Furthermore, after taking $\n{L}$ small enough, we have an open immersion of $p$-adic Lie groups    $\exp(\n{L})\subset \Gamma$ by \cite[Part II. Ch V.7 Corollary 2]{MR1176100}. Hence, we can assume without loss of generality that $\Gamma=\exp(\n{L})$ for a suitable $\bb{Z}_p$-submodule $\n{L}\subset \f{g}$ with $[\n{L},\n{L}]\subset p \n{L}$.

Now, if $\f{g}$ is the Lie algebra of an affine algebraic group $H$, we have that $\Gamma \subset H(\bb{Q}_p)$ is an open subgroup for $\n{L}$ small enough,\footnote{This follows from \cite[Chapter III.\$7.1, Th\'eor\`eme 2, Th\'eor\`eme 3]{bourbaki1975lie}.}  and the lemma follows from part (2). Moreover, we note that the previous holds if $\f{g}$ is a semisimple Lie algebra over $\bb{Q}_p$  thanks to \cite[Theorem  7.9]{borel2012linear}, or if $\f{g}$ is abelian by taking $H\cong \bb{G}_a^n$.

 We prove the general case by induction on the dimension of $\f{g}$, the basic cases being when $\f{g}$ is abelian or semisimple.  Let $\f{g}$ be a Lie algebra over $\bb{Q}_p$ which is not semisimple.   We have an extension 
\[
0\to \f{a} \to \f{g}\to \f{f} \to 0
\]
 where $\f{a}$ is a non-zero abelian Lie algebra and $\ob{dim}_{\bb{Q}_p}(\f{f})<\ob{dim}_{\bb{Q}_p}(\f{g})$. We can take a $\bb{Z}_p$-lattice  $\n{L}\subset \f{g}$ stable under the bracket as before, fitting in a short exact sequence of finite free $\bb{Z}_p$-modules
 \[
 0\to \f{a}\cap \n{L}\to \n{L}\to \n{L}/(\f{a}\cap \n{L})\to 0.
 \]
 We let 
 \begin{equation}\label{eqExtensionLieGroups}
 1\to A\to \Gamma\to F\to 1
 \end{equation}
 denote the associated morphism of compact $p$-adic Lie groups obtained by taking the exponential group of the lattices.
 The extension $\Gamma$ of $F$ produces a natural continuous linear action of $F$ on $A$, and $\Gamma$ arises from a continuous group cohomology class $\eta_{\Gamma}\in H^2(F, A)$, cf.\ \cite[Chapter IV.3 Theorem 3.12]{MR1324339}.
 Explicitly, let $\n{X}$ denote the topos generated by the site of $p$-adic Lie manifolds with topological covers.
 We can see $A$ as an  $F$-equivariant commutative group, and then the class $\eta_{\Gamma}$ corresponds to a map $\eta_{\Gamma}\colon BF \to (B^2A)/F$ over $BF$, where $F$ acts on $B^2A$ via the twofold-delooping of its additive action on $A$.
 The morphism $\eta_\Gamma$ is naturally pointed along the trivial $F$-torsor $*\to BF$ (as $H^2(*,A)=0$), and therefore it is equivalent to a morphism of $F$-augmented associative groups  in $\n{X}$
 \[
 F\xrightarrow{s} F\ltimes BA.
 \] 
The extension $\Gamma$ is then recovered as the pullback 
\[
\begin{tikzcd}
\Gamma \ar[r] \ar[d] & F  \ar[d,"\iota_0"] \\ 
F \ar[r, "s"] & F\ltimes BA,
\end{tikzcd}
\] 
where $s$ is the group homomorphism depicted by $\eta_{\Gamma}$, and $\iota_0$ the trivial map which is zero on the component of $BA$. Passing to classifying stacks we deduce a pullback diagram 
\[
B\Gamma= BF\times_{B(F\rtimes BA)} BF,
\]
where the right arrow is the trivial quotient and the left arrow the extension $\eta_{\Gamma}$.
Passing to the realization on locally analytic $p$-adic Lie groups in Gelfand stacks, in order to show that $\beta_{X,Z}^{\et}$ is an equivalence for $Z=B\Gamma^{\la} \times X^{\la}/G^{\la}$, by \Cref{LemmaPermanencePropertiesBetaXZ} (1)  it suffices to show that $\beta_{X,Z}^{\et}$  is an equivalence for $Z=BF^{\la}\times X^{\la}/G^{\la} $ and $Z=B(F^{\la}\ltimes BA^{\la}) \times X^{\la}/G^{\la}$. The case of $Z=BF^{\la}\times X^{\la}/G^{\la}$ follows by induction on the dimension of the Lie groups.
It is left to show the case $Z= B(F^{\la}\ltimes BA^{\la}) \times X^{\la}/G^{\la}$. For $A$  small  enough,  one has an $F$-equivariant open immersion of groups $\log\colon A\to \f{a}$  produced by the logarithm map.
Let $\bb{V}(\f{a})$ be the geometric vector bundle over $\bb{Q}_p$ associated to $\f{a}$.
The action of $F$ on $\f{a}$ produces a natural action of $F^{\la}$ on $\bb{V}(\f{a})$ and we have an $F^{\la}$-equivariant pullback diagram
 \[
 A^{\la} = A^{\dR}\times_{\bb{V}(\f{a})^{\dR}} \bb{V}(\f{a}). 
 \] 
 Taking suspensions we get a pullback diagram of Gelfand stacks
 \[
 B(F^{\la}\ltimes BA^{\la}) =  B(F^{\sm}\ltimes BA^{\dR} ) \times_{B(F^{\sm}\ltimes B\bb{V}(\f{a})^{\dR} )} B(F^{\la}\ltimes B\bb{V}(\f{a}) ).
 \]
 Therefore, by applying \Cref{LemmaPermanencePropertiesBetaXZ} (1) and (2), it suffices to show that $\beta_{X,Z}^{\et}$ is an equivalence when $Z=B(F^{\la}\ltimes B\bb{V}(\f{a}))\times X^{\la}/G^{\la}$.   Let us denote $Y:=BF^{\la}\times  X^{\la}/G^{\la}$, and consider the natural map of sheaves on the \'etale site $\n{C}$ of $X^{\la}/G^{\la}$
 \[
 \begin{tikzcd}
  \n{G}_{Z} \ar[d]\ar[r] & \n{F}_Z  \ar[d] \\ 
  \n{G}_{Y} \ar[r,"\sim"] & \n{F}_{Y}. 
 \end{tikzcd}
 \]
 By induction on the dimension, the map $\n{G}_{Y}\xrightarrow{\sim} \n{F}_Y$ is an isomorphism.
 Thus, in order to show that  $\beta_{X,Z}^{\et}$ is an equivalence, it suffices to argue for the fibers over $\n{G}_{Y}\cong \n{F}_Y$. In other words, after replacing $X^{\la}/G^{\la}$ by an \'etale map $X^{\prime,\la}/G^{\la}\to X^{\la}/G^{\la}$, we can fix a map $\nu\colon X^{\la}/G^{\la}\to BF^{\la}$. This morphism and the $F^{\la}$-representation $\f{a}$  give rise to a vector bundle $\s{V}$ on $X^{\la}/G^{\la}$.
 In this situation, the fibers $\n{G}_{Z/\nu}\to \n{F}_{Z/\nu}$ are the  \'etale sheaves given by the  $(-2)$-connective cover  of the cohomology  of $\s{V}$  in the Berkovich-\'etale site of $X/G^{\cont}$ or $X^{\la}/G^{\la}$.
 This map of cohomologies  is an isomorphism thanks to the equivalence of perfect complexes under locally analytic decompletion of \Cref{TheoDecompletion} (3). This finishes the proof of the lemma. \qedhere
 \end{enumerate}  
\end{proof}

\newpage
\section{Pro-\'etale comparison and applications to $p$-adic Hodge theory}
\label{sec:comp-proet-cohom}

In this section we prove \Cref{thm:equivalence-perfect-complexesintro} in its various incarnations (absolute and relative).
We also discuss some consequences.
The general theme is that the analytic prismatization of partially proper smooth rigid analytic varieties over $\Q_p$, or over a perfectoid space, allows the decompletion of various period sheaves on the pro-\'etale site, or vector bundles under them, and that this decompletion often reveals some extra structure.
For example, and most notably, the existence of the Sen operator of a semilinear $\C_p$-representation or of the canonical Higgs field of a v-vector bundle are direct consequences of the gerbe property of the analytic Hodge--Tate stack (\Cref{TheoremPropertiesHodgeTateStack}, \Cref{Theo:Geometric-Case-Hodge-TateStack-Decompletion}).

Before starting, let us fix notations. Recall from \eqref{composition-for-twistor} that for each Gelfand stack $X$ there exists a natural morphism
\[
  \alpha_X\colon \FF_{X^\diamond}\to X^{\Div^1}.
\]
Restricting to the Hodge-Tate locus it induces a morphism
\[
  \alpha_X^{\HT}\colon \widehat{X^\diamond}\to X^{\HT}.
\]
Moreover, we have the canonical morphism $\pi_X\colon X^\HT\to X$.

If $D$ is a perfectoid space over $\Q_p$, and a morphism $X\to D$ is fixed, then we also get the relative variants:
\[
   \alpha_{X/D}\colon \FF_{X^\diamond}\to X^{\Div^1/\FF_D}, ~~   \alpha^{\HT}_{X/D}\colon \widehat{X^\diamond}\to X^{\HT/D}=X^\HT
 \]
(where we used \Cref{prop:ht-stack-of-perfectoids} to identify the relative and absolute Hodge--Tate stacks of $X$) and $\pi_{X/D}\colon X^{\HT}\to X.$
 
\subsection{The arithmetic case: $p$-adic fields}
\label{sec:arithm-case:-q_ppr}

In this section, let $X=\GSpec(K)$ for a finite extension $K$ of $\Q_p$.\footnote{One can generalize the results of this subsection to finite extensions $K$ of $W(k)[1/p]$, where $k$ is a perfect field of finite transcendence degree, by replacing $X^{\Div^1}$ with the base change $X^{\Div^1}\times_{\GSpec(\Q_p)/\varphi^\Z} \GSpec(W(k)[1/p])/\varphi^\Z$.
For the identification of $\mathbb{L}_{X^\HT/\Q_p}$ in \Cref{sec:arithmetic-case:-p-adic-fields} one can use \Cref{sec:cotang-compl-morph-cotangent-complex-via-p-completion} to see that still $\mathbb{L}_{K/\Q_p}=0$ (because $\mathbb{L}_{k/\F_p}=0$ by perfectness of $k$).} The following theorem is the arithmetic incarnation of the pro-\'etale and prismatic cohomological comparison.

\begin{theorem}
  \label{sec:arithmetic-case:-p-adic-fields}
  The following hold true:
  \begin{enumerate}
  
\item The natural map $\alpha_X\colon \FF_{X^{\diamond}}\to X^{\Div^1}$ is $0$-proper. Moreover, the natural map $1\xrightarrow{\sim} \alpha_{X,*}1$   is an equivalence. In particular, for any map of Gelfand stacks $Y\to X^{\Div^1}$ the $*$-pullback map 
\[
\ob{D}(Y)\to \ob{D}(\FF_{X^{\diamond}}\times_{X^{\Div^1}} Y)
\]
is fully faithful. 
 
  \item Let $U\to X^{\Div^1}$ be a Berkovich-\'etale map, then the natural pullback functor  
    \[
      \alpha^\ast_X\colon \Perf(U)\xrightarrow{\sim} \Perf(\FF_{X^\diamond}\times_{X^{\Div^1}} U)
    \]
    is an equivalence of categories. 
    
  \end{enumerate}
\end{theorem}
\begin{proof}
  Since $K$ is a finite extension of $\Q_p$ (hence \'etale), we have a cartesian square 
  \[
\begin{tikzcd}
 \FF_{X^{\diamond}}  \ar[r] \ar[d] & \FF_{\Q_p^{\diamond}}  \ar[d] \\ 
X^{\Div^1} \ar[r]&  \Div^1. 
\end{tikzcd}
\] 
Therefore, by base change, to prove (1) it suffices to consider $X=\GSpec (\bb{Q}_p)$.

We have $\FF_{\Q_p^\diamond}=\FF_{\Q_p^\cyc}/\Z_p^{\times, \cont}$. By \Cref{thm:summary-main-results-sen} (1), $\FF_{\Q_p^\cyc}$ is locally classically decompleted for the action of $\Z_p^{\times}$ and by the explicit description of $\Div^1$, provided in \Cref{Section:Analytic-Prism}, we have an isomorphism
 $$
 \Div^1=\FF_{\Q_p^\cyc}^\la/G^\la
 $$
 such that the map $\FF_{\Q_p^\diamond}\to \Div^1$ identifies with the natural map $\FF_{\Q_p^\cyc}/\Z_p^{\times,\cont} \to \FF_{\Q_p^\cyc}^\la/\Z_p^{\times,\la}$.

Hence, part (1) follows from \Cref{TheoDecompletion} (1) (and \Cref{rk:reformulation-theo-decompletion}) and part (2) from \Cref{thm:summary-main-results-sen} (3). 
\end{proof}

\begin{remark}\label{RemarkTateRationalLocalizationsDecompletions}
The equivalence of perfect complexes of \Cref{sec:arithmetic-case:-p-adic-fields} is stated for Berkovich-\'etale  localizations in order to have a clean statement and proofs that only involve Gelfand stacks which are quotients of Berkovich spaces  groups on Berkovich spaces. Thanks to the general decompletion results of perfect modules of \Cref{TheoremDecompletionPerfectComplexesII}, and to arithmetic Sen theory (cf. \Cref{PropositionSenTheoryCyc}),  one can deduce similar equivalences of perfect complexes for Tate-rational localizations of $\FF_{X^{\diamond}}$, for instance, for the stacks associated to the period rings $\bb{B}_{I}$ where $I\subset (0,\infty)$ is a closed interval with rational boundary points. 
\end{remark}

We continue the discussion of the arithmetic case with the pro-\'etale and Hodge--Tate cohomological comparison. Before that, observe that the $\bb{G}_m^{\dagger}\cong \bb{G}_a^{\dagger}$-gerbe structure of the map $X^{\HT}\to X=\GSpec (K)$ induces an action of $(\ob{D}(\bb{G}_m^{\dagger}),*)$ on $\ob{D}(X^{\HT})$ via $!$-functors, where $*$ is the convolution   symmetric monoidal structure of $\ob{D}(\bb{G}_m^{\dagger})$. By Cartier duality, e.g. as in \cite[Theorem 3.2.37]{RCCartierGerbe},  we have a symmetric monoidal action of $\ob{D}(\bb{G}_a^{\an})$ on $\ob{D}(X^{\HT})$, where  $\ob{D}(\bb{G}_a^{\an})$ has its natural symmetric  monoidal structure.   Let $\theta$ denote the coordinate of $\bb{G}_a^{\an}$, taking the restriction of $\ob{D}(X^{\HT})$ to a module category over $\ob{D}(K[\theta])$, we see that any object $M\in \ob{D}(X^{\HT})$ is  endowed with a natural endomorphism $\theta_M\colon M\to M$ that we call the \textit{arithmetic Sen operator of $M$.}

\begin{theorem}\label{TheoremPropertiesHodgeTateStack}
The following hold true:
\begin{enumerate}

\item The map $\alpha_{X}^{\HT}\colon \widehat{X^{\diamond}}\to X^{\HT}$ is $0$-proper, and the natural map $1\xrightarrow{\sim} \alpha_{X,*}^{\HT} 1$ is an equivalence.  Furthermore, the pullback map induces an equivalence on perfect modules
\[
\Perf(X^{\HT})\xrightarrow{\sim} \Perf(\widehat{X^{\diamond}}).
\]

\item The cotangent complex $\mathbb{L}_{X^\HT/\Q_p}\cong \mathbb{L}_{X^\HT/X}$ is canonically isomorphic to $\s{O}_{X^\HT}[-1]$.

\item Let $M\in \ob{D}(\widehat{X^{\diamond}})$ be a (possibly infinite dimensional and solid) $\bb{C}_p$-semilinear Galois representation in the essential  image of $\alpha^{\HT,*}_X$. Then $M$ carries a natural $\bb{C}_p$-linear and Galois-equivariant  endomorphism  $\theta_M$ obtained as the pullback of the arithmetic Sen operator of $\alpha_{X,*}^{\HT}(M)$.  Furthermore, the arithmetic Sen operator satisfy the following properties:  

\begin{enumerate}

\item Let $K^{\cyc}$ be the completed cyclotomic extension of $K$ and  let $\Gamma=\Gal(K(\zeta_{p^{\infty}})/K)$ denote the cyclotomic Galois group. Then, under the natural equivalence of stacks
$$
\GSpec (\bb{C}_p)/\Gal_K^{\cont}= \widehat{X^{\diamond}}=\GSpec (K^{\cyc})/\Gamma^{\cont},
$$
 the Sen operator $\theta_{M}$ on $M$ is obtained as the base change along $K(\zeta_{p^{\infty}})\to \bb{C}_p$ of the derivation of the action of $\Gamma$ on the locally analytic representation\footnote{If $M$ is a nuclear $\bb{C}_p$-vector space then $\Gal_{K}^{\cont}$-cohomology agrees with solid $\Gal_{K}$-cohomology by \Cref{Remark-Difference-Actions}.}
\[
\alpha_{X,*}^{\HT}(M)|_{K(\zeta_{p^{\infty}})}=\Gamma(\Gal_{K^{\cyc}}^{\cont},M)^{\Gamma-\la}. 
\]

\item Suppose that $M\in \ob{D}(\GSpec (\bb{C}_p)/\Gal_{K}^{\cont})^{\heartsuit}$ is static and  that it belongs to the essential image of $\alpha^{\HT,*}_{X}$. Then $\ob{Sen}(M)$ is static and we have an isomorphism of cohomologies
\begin{equation}\label{eqIsomorphismCohomologySenOperators}
H^*(\Gal_{K}^{\cont}, M)\cong H^*(\theta_M,M)^{\Gal_K},
\end{equation}
where  $H^*(\theta_M,M)$ is the cohomology of the complex $\theta_M\colon M\to M$, and $(-)^{\Gal_K}$ refers to the classical (i.e. underived) Galois invariants in solid $\bb{Q}_p$-vector spaces.

\end{enumerate}

\end{enumerate}
\end{theorem}
\begin{proof}
We have a cartesian diagram 
\[
\begin{tikzcd}
\widehat{X}^{\diamond} \ar[r] \ar[d] & X^{\HT}  \ar[d] \\ 
 \FF_{\widehat{X}^{\diamond}} \ar[r] &  X^{\Div^1},
\end{tikzcd}
\]
thus \Cref{sec:arithmetic-case:-p-adic-fields} (1) implies the $0$-properness and the equivalence $1\xrightarrow{\sim} \alpha_{X,\ast}^{\HT} 1$.  To see that $\alpha^{\HT}_{X}$ induces an equivalence on perfect complexes, one can just apply \Cref{sec:arithmetic-case:-p-adic-fields} (2) as the $*$-pushforward functor along $\widehat{X^{\diamond}}\to \FF_{X^{\diamond}}$ is conservative and preserves perfect complexes, being a Cartier divisor. This shows part (1). 

Part (2) is a consequence of the fact that the map $X^{\HT}\to B\bb{G}_m$ is $\dagger$-formally \'etale, and that the cotangent complex of $B\bb{G}_m$ is $\s{O}[-1]$. 

For part (3), the description of the Sen operator as in (a) follows from the natural morphism of Gelfand stacks
\[
\alpha_{X}^{\HT}\colon \widehat{X^{\diamond}}= \GSpec( \bb{C}_p)/\Gal_{K}^{\cont}\cong  \GSpec (K^{\cyc})/\Gamma^{\cont} \to \GSpec (K(\zeta_{p^{\infty}}))/\Gamma^{\la}\cong X^{\HT}.  
\]
Indeed, from this presentation it is clear that, at the level of modules, the  $*$-pushforward along the  first equivalence  is  $\Gamma(\Gal(\C_p/K^{\cyc})^{\cont}, -)$, and that  the $*$-pushforward along the second map is the functor of $\Gamma$-locally analytic vectors. The action of the Sen operator appears, by applying Cartier duality, from the natural $B\bb{G}_a^{\dagger}$-torsor structure on $X^{\HT}$ obtained from the pullback square
\[
\begin{tikzcd}
X^{\HT} \ar[r] \ar[d] & \GSpec (K) \ar[d,"s"] \\ 
\GSpec (K)/\bb{G}_m \ar[r] & \GSpec (K)/\bb{G}_m^{\dR},
\end{tikzcd}
\]
where $s$  classifying the Tate-twist over $\Spd K$. Over $K_{\infty}:= K(\zeta_{p^{\infty}})$ the map $s$ trivializes as the uniform completion $K^{\cyc}\cong K_{\infty}^u$ is perfectoid. This produces an equivalence of Gelfand stacks
\[
\GSpec (K_{\infty})/\Gamma^{\dagger} \cong  \GSpec (K_{\infty})\times_{\GSpec (K)} X^{\HT}
\]
for which  $\GSpec (K_{\infty})/\Gamma^{\dagger}\to X^{\HT}$ is the $\Gamma^{\sm}$-torsor isomorphic to the map $\GSpec (K_{\infty})/\Gamma^{\dagger}\to \GSpec (K)/\Gamma^{\la}$. The previous discussion shows that the $B\bb{G}_a^{\dagger}$-torsor structure on $X^{\HT}$ induces the trivial action of $B\Gamma^{\dagger}$ on $\GSpec (K_{\infty})$. We deduce that for  $N\in \ob{D}(X^{\HT})$,  the restriction to $\GSpec (K_{\infty})$ of the  aritmetic Sen operator of the module $N$ is nothing but the derivation of the locally analytic action of $\Gamma$.  This finishes the proof of point (a). 

Next, we show point (b) of part (3). For the claim about the $t$-exactness of the decompletion functor, it suffices to observe that the base change along $K(\zeta_{p^{\infty}})\to \bb{C}_p$ is $t$-exact for the solid tensor product, namely, it suffices to see that $K\to \bb{C}_p$ is $t$-exact, but this follows from \cite[Lemma 3.21]{jacinto2021solid}.  Finally, for the natural isomorphism \eqref{eqIsomorphismCohomologySenOperators} one can argue as in \cite[Proposition 2.5.5]{camargo2022geometricsen}. We add the argument for completeness. Let $N\in \ob{D}(X^{\HT})^{\heartsuit}$ be a static module for the presentation  $X^{\HT}=\GSpec (K_{\infty})/ \Gamma^{\la}$,  and let $M=N\otimes_{K_{\infty}} \bb{C}_p$. By projection formula we have that 
\[
\Gamma(\Gal_{K}^{\cont},M) = \Gamma(X^{\HT}, N)= \Gamma(\Gamma^{\sm},\Gamma(\theta_N, N)),
\]
where in the second equivalence we used the group cohomology comparisons of  \cite[Theorem 6.3.4]{jacinto2023solid}. Therefore, since $\Gamma^{\sm}$-cohomology is $t$-exact as we are working in characteristic zero,  to deduce \eqref{eqIsomorphismCohomologySenOperators} it suffices to show that 
\[
H^*(\theta_N,N) = ((H^*(\theta_M,M) )^{\Gal_{K^{\cyc}}})^{\Gamma-\sm}.
\]
By the equivalence $\GSpec (\bb{C}_p)/\Gal_{K^{\cyc}}\cong \GSpec (K^{\infty})$, and the fact that $M=N\otimes_{K_{\infty}} K^{\cyc}$, we have that 
\[
\Gamma(\Gal_{K^{\cyc}}^{\cont}, M) = N\otimes_{K_{\infty}} K^{\cyc}
\]
and similarly for the representations $H^*(\theta_M,M)$. In particular, we have that 
\[
H^*(\theta_N , N)\otimes_{K_{\infty}} K^{\cyc} \cong (H^*(\theta_M, M))^{\Gal_{K^{\cyc}}}.
\]
Taking $\Gamma$-smooth vectors, and noticing that the action of $\Gamma$ on  $H^*(\theta_N , N)$ is locally constant as the Lie algebra acts trivially, we obtain $H^*(\theta_N , N)\cong ((H^*(\theta_M, M))^{\Gal_{K^{\cyc}}} )^{\Gamma-\sm}$, as desired. 
\end{proof}

\begin{example}
  \label{sec:arithmetic-case:-p-new-proof-of-tates-theorem}
  If $W=\s{O}_{X^\HT}\{n\}$, then one calculates that $\theta_{W}=n$, e.g., by pulling back along $X^\HT\to B\Gm$.
Therefore \Cref{TheoremPropertiesHodgeTateStack} (1)  implies
  \[
    H^i(\mathrm{Gal}_K,\C_p(n))\cong H_{\proet}^i(X,\widehat{\s{O}}(n))\cong H^i(X^\HT,\s{O}\{n\})\cong H^i(K\xto{n}K)
  \]
  for $i\geq 0$, which is Tate's fundamental calculation of $H^\ast(\mathrm{Gal}(\overline{K}/K),\C_p)$ in $p$-adic Hodge theory.
  While logically this argument is not circular, Tate's argument appears  in the proof of  \Cref{PropositionSenTheoryCyc}.
  \end{example}

\begin{remark}
  \label{sec:arithmetic-case:-p-unramified-case}
  If $K/\Q_p$ is unramified, then $X^{\Div^1}\cong \GSpec(K)\times \Div^1$.
  Indeed, this follows from $X^\HK\cong X\times_{\GSpec(\Q_p)} \Q_p^\HK$ and $X^{\Div^1}\cong X^{\cone}\times_{X^{\cone,\dR}}X^\HK$.
\end{remark}

\begin{example}
  \label{sec:arithmetic-case:-p-pullback-of-b-dr-and-construction-of-overconvergent-period-ring}
  Localizations on $\Div^1$ allow the construction of various period sheaves. 
  For example, let $\Div^{1,\dagger}:=\Div^1\times_{\A^{1,\dR}/\Gm^{\dR}}B\Gm^\dR$ be the overconvergent neighborhood of the Hodge--Tate point, and $\iota\colon \Div^{1,\dagger}\to \Div^1$ the inclusion.
  Given any perfectoid space $Z$  over $\Q_p$, let $\rho\colon \FF_{Z^\flat}\to \Div^1$ be the natural morphism.
  Then $\rho^\ast(\iota_\ast(\s{O}_{\Div^{1,\dagger}}))$ identifies with the overconvergent $B^+_\dR$-period sheaf that is constructed and studied in \cite{wiersig2023solution}.
\end{example}

Next, we explain how the classical functors $D_{\HT}, D_\dR, D_{\mathrm{pdR}},\ldots $ fit into our perspective.

\begin{example}
  \label{sec:arithmetic-case:-p-d-ht-d-dr}
  By \Cref{sec:arithmetic-case:-p-adic-fields} we have a fully faithful embedding $V\mapsto \mathcal{E}_V$ of the category of continuous finite dimensional (over $\Q_p$) $\mathrm{Gal}_K:=\mathrm{Gal}(\overline{K}/K)$-representations to vector bundles on $\Div^1$ (the essential image consists of those vector bundles whose pullback to $\FF_{\Q_p^\cycl}$ has slope $0$).
  In this remark, we explain how one can use geometry to reinterpret various constructions of their $p$-adic Hodge theory.
  \begin{enumerate}
  \item ($D_\HT$)
   Vector bundles on $B\Gm$ identify with graded finite dimensional $\Q_p$-vector spaces and the $H^0$ of the  pushforward of $\mathcal{E}_{V|K^\HT}$ along $K^{\HT}\to B\Gm$ identifies with $D_{\HT}(V)$ (up to changing the grading by a sign, depending on the exact conventions).
  \item ($D_\dR$) Vector bundles on $\A^1/\Gm$ identify with filtered finite dimensional $\Q_p$-vector spaces and the $H^0$ of the pushforward of $\mathcal{E}_{V|K^{\Div^{1,\dagger}}}$ along $K^{\Div^{1,\dagger}}\to \A^1/\Gm$ identifies with $D_\dR(V)$ (cf.\ \cite{gao2026galois} for related considerations).
  \item ($D_{\mathrm{pdR}}$) Vector bundles on $B\Ga$ identify with finite dimensional $\Q_p$-vector spaces with a nilpotent endomorphism.
    Let $\gamma\colon K^{\Div^{1,\dagger}}\to B\Ga$ be the morphism classified by the class $\mathrm{log}(\chi_{\cycl})\in H^1(K^{\Div^{1,\dagger}},\s{O})$.
    Then pushforward along the map $\gamma':=(\mathrm{can},\rho)\colon K^{\Div^1,\dagger}\to \A^1/\Gm\times B\Ga$ of $\mathcal{E}_{V|K^{\Div^{1,\dagger}}}$ identifies with $D_{\mathrm{pdR}}(V)$.
  \item ($D_{\mathrm{pst}}$) Vector bundles on $K^{\HK}$ identify with $(\varphi,N,\mathrm{Gal}_K^{\sm})$-modules over $\Q_p^{\mathrm{un}}$ via pullback from $Z:=B_{\GSpec(\Q_p)/\varphi^\Z}(\mathbb{V}(-1))/\mathrm{Gal}_K^\sm$ (cf.\ \cite[\S 7]{drff_stacks}), and the pushforward of $\mathcal{E}_{V}[1/t]$ along the map $K^{\Div^1}\to Z$ identifies with $D_{\mathrm{pst}}(V)$ (after forgetting the filtration on $D_{\mathrm{pst}}(V)$). Conversely, pulling back a $(\varphi,N,\mathrm{Gal}_K^{\sm})$-module over $\Q_p^{\mathrm{un}}$ along $K^{\Div^1}\to K^\HK$ produces a $(\varphi,\Gamma)$-module (in general not \'etale) over the Robba ring, often called ``Berger's $p$-adic differential equation''. 
    \item ($D_{\mathrm{Sen}}$) Let $K_\infty=K(\mu_{p^\infty})$ be the (uncompleted) cyclotomic extension of $K$, and $X_\infty=\GSpec(K_{\infty})$.
    Then the morphism $X_\infty^{\HT}\to X_{\infty}$ admits a natural section $X_{\infty}\to X^{\HT}$ (induced by the rational cyclotomic prism).
    In particular, we get a morphism $\beta\colon X_\infty\times B\Ga^\dagger\to X^\HT\to \Div^1$.
    Vector bundles on $X_\infty\times B\Ga^\dagger$ identify with pairs $(M,\theta_M)$ with $M$ a finite dimensional $K_\infty$-vector space and $\theta_M\colon M\to M$ a $K_\infty$-linear endomorphism.
    Now, $D_{\mathrm{Sen}}(V)$ identifies with the pullback of $\mathcal{E}_{V}$ along $\beta$.
  \item ($D_{\mathrm{dif}}^+$) Similarly to (5), one gets a morphism 
  \[
  \beta^\prime\colon  \Spf K_{\infty}[[t]]/\Ga^{\dagger}=\bigcup_{n} \GSpec(K_\infty[[t]]/t^n)/\Ga^\dagger\to \Div^1,
  \] where $t=\log([\varepsilon])$ is Fontaine's element and $\Ga^\dagger$ acts on $\GSpec(K_\infty[[t]]/t^n)$ via the derivation $t \frac{\partial}{\partial t}$.
    Vector bundles on $\bigcup_{n} \GSpec(K_\infty[[t]]/t^n)/\Ga^\dagger$ identify with finite free $K_\infty[[t]]$-modules $(M,\nabla_M)$ with a $t$-connection $\nabla_M$, i.e., a connection lying over the derivation $t \frac{\partial}{\partial t}$.
    Pullback along $\beta^\prime$ identifies with Fontaine's lift $D_{\mathrm{dif}}^+$ of $D_{\mathrm{Sen}}$ for $B^+_\dR$-representations. 
  \end{enumerate}
\end{example}

\subsection{The geometric case: smooth varieties over perfectoid spaces}
\label{sec:geometric-case}

Let $D$ be a perfectoid space over $\Q_p$,\footnote{We assume this only for simplicity of exposition, one could also work relative to a rational prism $(S,\mathcal{X},D)$, e.g., the rational cyclotomic prism, where $D=\GSpec(\Q_p(\mu_{p^\infty}))$ is not perfectoid.} and let $X\to D$ be a rigid smooth morphism of Berkovich spaces.

We have the following analogue of \Cref{sec:arithmetic-case:-p-adic-fields} in this relative situation.

\begin{theorem}
  \label{sec:geom-case:-smooth-relative-case}
  The following hold true:
  \begin{enumerate}
  \item The map $\alpha_{X/D}\colon \FF_{X^\diamond}\to X^{\Div^1/\FF_D}$ is $0$-proper,  and the natural map  $ 1\xrightarrow{\sim} \alpha_{X/D,\ast}1$ is an equivalence.  In particular, for any map of Gelfand stacks  $Y\to X^{\Div^1/\FF_D}$ the $*$-pullback 
  \[
  \ob{D}(Y)\to \ob{D}(Y\times_{X^{\Div^1/\FF_D}} \FF_{X^{\diamond}})
  \]
  is fully faithful.

  \item Let  $U\subseteq X^{\Div^1/\FF_D}$ be a Berkovich-\'etale map, then the natural pullback functor
  \[
  \Perf(U)\xrightarrow{\sim} \Perf(\FF_{X^{\diamond}}\times_{X^{\Div^1/\FF_D}} U)
  \]
  is an equivalence of categories. 
    \end{enumerate}
\end{theorem}

\begin{proof}
First we will prove that the map $\alpha_{X/D}$ is $0$-proper, and that $1\xrightarrow{\sim} \alpha_{X/D,*} 1$ is an equivalence which yields part (1).

All the previous statements are local on target, hence we can assume that $D=\GSpec(A)$ is affinoid.
  By descent we may furthermore assume that there exists a morphism $D\to \GSpec(\Q_p^\cycl)$, that we fix from now on.
  Thanks to the descent of the analytic prismatization for $\dagger$-formally \'etale covers, see \Cref{TheoremDescentPrismatization} (1), we can furthermore localize on $X$ in the rigid \'etale topology, in particular, for the analytic topology (we also note that $X\mapsto \FF_{X^\diamond}$ preserves open immersions and finite \'etale maps\footnote{Though not $\dagger$-rational localizations, because these get completed on $X^\diamond$.}).
  Thus, we can assume without loss of generality that there is a rigid \'etale map $X\to \mathbb{G}_{m,A}^{d}$ factoring as open localizations and finite \'etale maps.
 
  In this case, we have a cartesian diagram
\[
\begin{tikzcd}
 \FF_{X^{\diamond}}  \ar[r] \ar[d] & \FF_{\bb{G}_{m,A}^{d,\diamond}}  \ar[d] \\ 
X^{\Div^1/\FF_{D}} \ar[r]&  \bb{G}_{m,A}^{d,\Div^1/\FF_{D}}. 
\end{tikzcd}
\] 
so that it suffices to prove part (1) for $X=\bb{G}_{m,A}^d$. We can even, by K\"unneth, assume $d=1$. Moreover, by passing to the $\bb{Z}$-covers before taking by Frobenius, it suffices to show the claim for the map 
\begin{equation}\label{EqKeyMapGmProetComparison}
\FFdisc{\bb{G}_{m,A}}^{\times}\to \bb{G}_{m,A}^{\widetilde{\Div^1}/\FFdisc{A}^{\times}}. 
\end{equation}
By \Cref{sec:analyt-prism-gm-description-over-cyclotomic-prism}, the map \eqref{EqKeyMapGmProetComparison} is equivalent to 
\[
\widetilde{\bb{G}}_{m,\FFdisc{A}^{\times}}^{u}/ \Gamma^{\cont}\to \widetilde{\bb{G}}_{m,\FFdisc{A}^{\times}}/\Gamma^{\la}
\]

where:
\begin{enumerate}
\item  $\Gamma\cong \bb{Z}_p(1)$ is the Galois group of $\varprojlim_{x\mapsto x^p} \bb{G}_{m,A}\to \bb{G}_{m,A}$.  We fix  $\bb{Z}_p(1)\cong \bb{Z}_p$ a trivialization determined by a compatible sequence of roots of unit $(\zeta_{p^n})_{n\geq 0}$, and denote $\gamma\in \Gamma$ its associated the generator.

\item $\widetilde{\bb{G}}_m=\varprojlim_{x\mapsto x^p} \bb{G}_m$ is the decompleted preperfectoid multiplicative group with coordinate $[T^{\flat}]$.

\item The action of $\Gamma$ on  $\widetilde{\bb{G}}_m$ is via  $\gamma\cdot [T^{\flat}]^r = [\epsilon]^r [T^{\flat}]^r$ where $\epsilon=(\zeta_{p^n})_{n\geq 0}$, and $r\in \bb{Z}[\frac{1}{p}]$. 
\end{enumerate} 
Then the $0$-properness   and the equivalence $1\xrightarrow{\sim} \alpha_{X/D,*}1$ follow from \Cref{TheoDecompletion} (1), \Cref{rk:reformulation-theo-decompletion} and the decompletions  of geometric Sen theory \Cref{cor:decompletion-geometric-case}.   This finishes the proof of (1). 

For part (2) one reduces as above to the case where $X$ is rigid \'etale over $\bb{G}_m^d$. In this case, the result follows from \Cref{thm:summary-main-results-sen} (3). 
\end{proof}

For the Hodge--Tate stack, we have the following analogue of \Cref{TheoremPropertiesHodgeTateStack}. We recall that $X^{\HT}\to  X$ has the structure of a $B\bb{T}_{X/D}^{\dagger}\{1\}$-torsor, cf. \Cref{PropRelativeDescriptionPrismatization}. In particular, passing to sheaves, we get an action of $\ob{D}(B\bb{T}_{X/D}\{1\}^{\dagger})$ on $\ob{D}(X^{\HT})$ by convolution along $!$-functors, and by passing to Cartier duals a $\ob{D}(\bb{T}^*_{X/D}\{-1\})$-linear structure on $\ob{D}(X^{\HT})$.   Let $M\in \ob{D}(X^{\HT})$. This  produces a natural structure of $\Sym_{\s{O}_X}(\bb{T}_{X/D}\{1\})$-module on $M$, or equivalently, a structure of $\Omega^1_{X/D}\{-1\}$-Higgs field $\theta_M$ on $M$ over $\ob{D}(X^{\HT})$. This Higgs field structure  on $M$ is called the \textit{geometric Sen operator of $M$}. 

\begin{theorem}\label{Theo:Geometric-Case-Hodge-TateStack-Decompletion}
The following hold true:

\begin{enumerate}

\item   The map $\alpha_{X/D}^{\HT}\colon \widehat{X^{\diamond}}\to X^{\HT}$ is $0$-proper  and the natural map $1\xrightarrow{\sim} \alpha_{X/D,*}^{\HT} 1$ is an equivalence. Furthermore, for all Berkovich-\'etale map $U\to X$ the pullback map induces an equivalence on perfect modules
\[
\Perf(X^{\HT}\times_{X} U ) \xrightarrow{\sim}\Perf(\widehat{X^{\diamond}}\times_X U). 
\]

\item The cotangent complex $\mathbb{L}_{X^{\HT}/X}$ is canonically isomorphic to $\Omega^1_{X/D}\{-1\}[-1]$, where
\[
\Omega^1_{X/D}\{-1\}=\pi_X^\ast \Omega^1_{X/D}\otimes_{\s{O}_{X^{\HT}}}\s{O}_{X^{\HT}}\{-1\}.
\]

\item Let $M\in \ob{D}(\widehat{X^{\diamond}})$ be a solid $\widehat{\s{O}}_X$-sheaf over $X_{\proet}$ in the essential image of $\alpha_{X/D}^{\HT,\ast}$. Then $M$ carries a natural $\widehat{\s{O}}_X$-linear structure of $\Omega^1_{X/D}\{-1\}$-Higgs module structure obtained as the pullback of the geometric Sen operator of $\alpha_{X/D,*}^\HT M$.

\item Let $\pi_X\colon X^{\HT}\to X$. Then, the kernel $1\in \ob{D}(B\bb{T}_{X/D}^{\dagger}\{-1\})$ acts on $\ob{D}(X^{\HT})$ via the functor $\pi_X^*\pi_{X,*}$. Furthermore, the relative de Rham complex of the map $X\to B\bb{T}_{X/D}^{\dagger}\{-1\}$ produces a natural filtration on $\pi_X^*\pi_{X,*}$ depicted as the Higgs complex for $M\in \ob{D}(X^{\HT})$
\[
\pi_X^*\pi_{X,\ast}M \cong [M \xrightarrow{\theta_M}  M\otimes_{\s{O}_X} \Omega^{1}_{X/D}\{-1\}\to \cdots \to M \otimes_{\s{O}_X} \Omega^{\dim X}_{X/D}\{-1\} \to 0 ].  
\]

\item Let $M$ be a relative locally analytic ON Banach $\widehat{\s{O}}_X$-sheaf  as in \cite[Definition 1.0.1]{camargo2022geometricsen}, then $M$ belongs to the essential image of the functor $\alpha^{\HT,*}_{X/D}$, and the geometric Sen operator of $M$ as in (2) agrees with the geometric Sen operator of \cite[Theorem 1.0.3]{camargo2022geometricsen}.
\end{enumerate}

\end{theorem}
\begin{proof}
Part (1) follows the exact same argument as in \Cref{TheoremPropertiesHodgeTateStack}, namely, it follows from \Cref{TheoremPropertiesHodgeTateStack} (1) and the pullback diagram $\widehat{X^{\diamond}}=\FF_{X^{\diamond}}\times_{X^{\Div^1/\FF_D}}X^{\HT}$. 

For part (2), recall the Cartesian square 
\begin{equation}\label{eqCartesianDiagramHTStack}
\begin{tikzcd}
X^{\HT}  \ar[d] \ar[r,"\pi_X"] & X \ar[d,"s"] \\ 
X/\bb{T}_{X/D}\{1\} \ar[r,"f"] &  X/\bb{T}_{X/D}^{\dR}\{1\}
\end{tikzcd}
\end{equation}
From this we deduce that $X^{\HT}\to X/\bb{T}_{X/D}\{1\}$ is $\dagger$-formally \'etale. It is  clear that the cotangent complex of $X/\bb{T}_{X/D}\{1\}$ is precisely $\Omega^1_{X/D}\{-1\}[-1]$, taking a pullback to $X^{\HT}$ we obtain the desired claim.

Part (3) is a formal consequence of the fully faithfulness of the map $\alpha^{\HT*}_{X/D}$ and the $\bb{T}_{X}^{\dagger}\{1\}$-gerbe structure of $X^{\HT}\to X$. 

For part (4), recall that in the cartesian diagram \eqref{eqCartesianDiagramHTStack}, $s$ arises from the natural isomorphism of pro-\'etale cohomology $ R^1\nu_* \widehat{\s{O}}_{X} \cong \Omega^1_{X/D}\{-1\}$ where $\nu \colon X_{\proet}\to X_{\an}$ is the projection of sites. The $B\bb{T}_{X/D}^{\dagger}\{1\}$-torsor structure of $\pi_X$ arises from the similar torsor structure of the map  $f$. The kernel $1\in \ob{D}(B\bb{T}_{X}^{\dagger}\{1\})$ arises from $*$-pullback of the unit  along the map $B\bb{T}_{X}^{\dagger}\{1\}\to X$ (where $X$ is seen as the trivial group over $X$), thus, from the cartesian diagram 
\[
\begin{tikzcd}
B\bb{T}_{X}^{\dagger}\{1\}\times_X X^{\HT }\ar[r,"m"] \ar[d,"\mathrm{pr}_{X^{\HT}}"]  &  X^{\HT} \ar[d,"\pi_X"] \\ 
X^{\HT}  \ar[r,"\pi_X"] & X
\end{tikzcd}
\]
with $m$ the projection map, the kernel $1_{B\bb{T}_{X}^{\dagger}\{1\}}$ induces the functor $ m_!\circ \mathrm{pr}_{X^{\HT}}^*\cong \pi_X^*\pi_{X,!} =\pi_X^*\pi_{X,\ast}$ where we have used the natural trivialization of the codualizing sheaf of $\pi_X$ as this map is $0$-proper. In terms of the natural geometric Sen action of $\Sym_{\s{O}_X} (\bb{T}_{X/D}\{1\})$ on objects in $\ob{D}(X^{\HT})$, one has by construction and Cartier duality that 
\[
\pi_X^*\pi_{X,*} M = \Hom_{\Sym_{\s{O}_X} (\bb{T}_{X/D}\{1\})}(\s{O}_{X^{\HT}}, M)\in \ob{D}(X^{\HT}).
\]
The final claim about the filtration of the functor $\pi_X^*\pi_{X,\ast}$ follows from the fact that the relative de Rham cohomology of the map $X\to X/\bb{T}_X^{\dagger}\{1\}$ is the unit as  $\bb{T}_X^{\dagger}\{1\}$ is an overconvergent group, this implies that    the relative de Rham cohomology produces a filtration of $1_{B\bb{T}_X^{\dagger}\{1\}}$ with graded pieces given by $\bigwedge^{\bullet} (\bb{T}_X^{\dagger}\{1\})$. By Cartier duality, this filtration corresponds to the $\s{O}_X$-linear dual of the standard Koszul complex of $\s{O}_X$ as $\Sym_{\s{O}_X} (\bb{T}_{X/D}\{1\})$-module, which produces by construction the Higgs complex  of $M$.

Finally for part (5), we can prove that $M$ belongs to the essential image of $\alpha^{\HT,*}_{X/D}$ locally on $X$. On the other hand, by  construction,  the geometric Sen operators of \cite{camargo2022geometricsen} are obtained by recollement along the equivalence $a\colon R^1\nu_* \widehat{\s{O}}_{X}\cong \Omega^1_{X/D}\{-1\}$ of the local geometric Sen operators obtained by differentiating toroidal charts. Since the $B\bb{T}_{X/D}^{\dagger}\{1\}$-torsor structure of $X^{\HT}$ arises from the map $s$ which is induced by   the equivalence  $a$, it suffices to identify the geometric Sen operators local in coordinates. Thus, for both statements we can assume without loss of generality that $X$ admits an \'etale chart $\psi\colon X\to \bb{G}_{m,D}^d$. Let $\widetilde{X}\to X $ be the $\Gamma^{\sm}=\bb{Z}_p(1)^{d,\sm}$-torsor obtained by taking pullbacks along the uncompleted perfect multiplicative group $\widetilde{\bb{G}}^d_{m,D}\to \bb{G}^d_{m,D}$, the preservation of $\dagger$-formally \'etaleness of \Cref{sec:defin-first-prop-5-prismatization-for-dagger-formally-etale-maps} and the explicit presentation of the prismatization of $\bb{G}_m$ of \Cref{sec:analyt-prism-gm-description-over-cyclotomic-prism} yields 
\[
\widehat{X^{\diamond}}\cong \widetilde{X}^u/\Gamma^{\cont}\to \widetilde{X}/\Gamma^{\la}\cong X^{\HT} \to \widetilde{X}/\Gamma^{\sm}=X.
\]
In particular, the $\bb{T}_{X/D}^{\dagger}\{1\}$-gerbe structure of $\pi_X$ arises from the $\Gamma^{\dagger}$-gerbe structure by base change along $\bb{Z}_p\to \s{O}_X$  by considering the dual of the  isomorphism 
\[
\Lie (\Gamma)^{\vee}\otimes_{\bb{Z}_p} \s{O}_X \cong  R^1 \nu_* \widehat{\s{O}}_X \cong \Omega^1_{X/D}\{-1\},
\]
where the first equivalence follows from \cite[Proposition 3.3.1]{camargo2022geometricsen}, and the second from the map $a$.

 The first claim about decompletions of relatively locally analytic ON Banach $\widehat{\s{O}}_X$-sheaves $M$ then follows from \cite[Theorem 2.4.4 (2)]{camargo2022geometricsen} by taking a Tate-\'etale covering of $X$ where $M$ admits a Banach basis. Finally, by construction, the geometric Sen operators of part (3) of sheaves on  $\widetilde{X}^{\HT}\cong \widetilde{X}/\Gamma^{\dagger}$ arise from the natural derivations of along the action of $\Gamma^{\dagger}$, these are precisely the local geometric Sen operators constructed in \cite{camargo2022geometricsen}.
\end{proof}

\begin{remark}
  \label{sec:geom-case:-smooth-non-perfectoid-base}
  Instead of a perfectoid base $D$ we could work with \emph{any} rigid smooth morphism $f\colon X\to Y$ of Berkovich spaces, by replacing $\FF_D$ by $\FF_{Y^\diamond}$, and $X^{\Div^1/\FF_D}$ by the fiber product $\FF_{Y^\diamond}\times_{Y^{\Div^1}}X^{\Div^1}$.  The analogues of \Cref{sec:geom-case:-smooth-relative-case,Theo:Geometric-Case-Hodge-TateStack-Decompletion} for this map will follow from the case of perfectoid spaces by descent along an arc-cover of $Y^{\diamond}$.
  \end{remark}

\begin{remark}
\label{conclusion-ht-stack-valid-smooth-not-pp}
 \Cref{Theo:Geometric-Case-Hodge-TateStack-Decompletion}   also applies verbatim to a smooth rigid space over $D$, not necessarily partially proper (even though the associated Berkovich space may not be rigid smooth).  However,  it  does not hold for arbitrary smooth $\dagger$-rigid spaces (such as the overconvergent unit disc), as $\widehat{X^\diamond}$ is independent of the $\dagger$-structure.
\end{remark}

\begin{remark}[Analogy with the work of Ogus--Vologodsky]
The theory of vector bundles on the Hodge--Tate stack of a smooth rigid space $X$ over a perfectoid field can be compared to the theory of the de Rham stack of a smooth scheme $Y$ over a perfect field $k$ of characteristic $p$. The relative de Rham stack $Y^{\dR/k}$ (defined by transmutation from $\Ga/\Ga^\sharp$) is a $\mathbb{T}_{Y/k}^\sharp$-gerbe over $Y'$, the Frobenius twist of $Y$. 

Therefore, a vector bundle $V$ on $Y^{\dR/k}$ carries a natural $\Omega^1_{Y'/k}$-Higgs module structure which is nothing but the \textit{$p$-curvature of $V$}. The vector bundle $V$ descends to $Y'$ if and only if the $p$-curvature vanishes. In an  analogous fashion, a vector bundle $V$ on $X^{\HT}$ carries the natural structure of a  $\Omega^1_{X}\{-1\}$-Higgs module given by the geometric Sen operator. The vector bundle descends to $X$ if and only it the geometric Sen operator vanishes.   In both cases, the Higgs field and the $p$-curvature arise from the gerbe structure.

Finally, by a result of Terentiuk \cite{terentiuk2026ogus} (who uses this to reprove a theorem of Ogus--Vologodsky), the pushout of the $\mathbb{T}_{Y/k}^\sharp$-gerbe $Y^{\dR/k}\to Y'$ along $\mathbb{T}_{Y/k}^\sharp \to \mathbb{T}_{Y/k}$ is the gerbe of liftings of $Y$ to $W_2(k)$, which can be compared to the fact that the pushout of $X^\HT$ along $\mathbb{T}_{X/C}^\dagger(1) \to \mathbb{T}_{X/C}(1)$ is the gerbe of liftings to $B_\dR^+(C)/t^2$ (cf.\ \Cref{sec:p-adic-twistor-stack}).\footnote{To mention yet another example, the relation between the $p$-curvature of a vector bundle on $Y^{\dR/k}$ coming from the syntomification of $Y$ and the Kodaira-Spencer class (a generalization of Katz's formula) should be compared to the relation between the Higgs field of a vector bundle on $X^\HT$ coming from the relative syntomification of $X$ (i.e.\ a de Rham relative $(\varphi,\Gamma)$-module, by \cite{hauck}) and the Kodaira-Spencer class, cf. e.g.\ \cite{howe2025geometric}.}   

In fact, this analogy can be made more precise using the integral version of the analytic prismatization, which in particular leads to an integral version of the analytic Hodge--Tate stack (defined for $p$-adic analytic spaces over $\Z_p$), interpolating between the analytic Hodge--Tate stack of the generic fiber and a refined, Frobenius-untwisted, version of the de Rham stack of the special fiber. 
\end{remark}

\subsection{The mixed case: smooth varieties over $p$-adic fields}
\label{sec:mixed-case}
In this subsection, we combine our results from \Cref{sec:arithm-case:-q_ppr} and \Cref{sec:geometric-case}. Let $K$ be a finite extension of $\Q_p$.\footnote{Similar generalizations to the ones mentioned in \Cref{sec:arithm-case:-q_ppr} can be done here.}

\begin{theorem}
\label{sec:mixed-case:-smooth-mixed-case1}
Let $X$ be a smooth partially proper rigid space over $K$.  The map $\alpha_{X}\colon \FF_{X^\diamond}\to X^{\Div^1}$ is $0$-proper, and the natural map $ 1\xrightarrow{\sim}  \alpha_{X,\ast}(1)$ is an equivalence. Furthermore, for all Berkovich-\'etale morphism $Y\to X^{\Div^1}$ the natural map 
  \[
  \Perf(Y)\xrightarrow{\sim}  \Perf(\FF_{X^{\diamond}}\times_{X^{\Div^1}} Y)
  \]
  is an equivalence.
\end{theorem}
\begin{proof}
The proof is entirely analogous the proof of \Cref{sec:geom-case:-smooth-relative-case}, using \Cref{thm:summary-main-results-sen} (4). 
\end{proof}

 We turn to the Hodge--Tate stack. Recall that the absolute Hodge--Tate stack sits in a pullback square
\begin{equation}\label{eqSquareAbsoluteHTStackSection6}
\begin{tikzcd}
X^{\HT}  \ar[r] \ar[d] & X \ar[d,"s"] \\  
 X/\bb{T}_X\rtimes \bb{G}_m \ar[r]  &  X/(\bb{T}_{X}\rtimes \bb{G}_m)^{\dR/X}
\end{tikzcd}
\end{equation}
where $\bb{G}_m$ acts by multiplication by $-1$ on $\bb{T}_X$, and $s$ is the morphism of Gelfand stacks induced by the Faltings' extension, see \Cref{restriction-ht-overconvergent}. Let us denote  $\n{G}_{\mathrm{FE}}:= \bb{T}_X\rtimes \bb{G}_m $, it is a smooth partially proper rigid analytic group over $X$.  The pullback square \eqref{eqSquareAbsoluteHTStackSection6} produces a natural action of 
\[
\ob{End}_{B \n{G}_{\mathrm{FE}}^{\dR/X}} (B\n{G}_{\mathrm{FE}}) = \n{G}_{\mathrm{FE}} \backslash \n{G}_{\mathrm{FE}}^{\dR/X}  / \n{G}_{\mathrm{FE}} 
\]  on  $X^{\HT}$ by convolution in the category of kernels of $X$, and thus an action of the monoidal category $\ob{D}(\n{G}_{\mathrm{FE}} \backslash \n{G}_{\mathrm{FE}}^{\dR/X}  / \n{G}_{\mathrm{FE}} )$ on $\ob{D}(X^{\HT})$.  We have a natural equivalence 
\[
\n{G}_{\mathrm{FE}} \backslash \n{G}_{\mathrm{FE}}^{\dR}  / \n{G}_{\mathrm{FE}}  \cong B(\n{G}_{\mathrm{FE}} \subset \n{G}_{\mathrm{FE}} \times \n{G}_{\mathrm{FE}})^{\dagger}=B \n{G}_{\mathrm{FE}}^{\ob{Har}},
\] 
where  $\n{G}_{\mathrm{FE}}^{\ob{Har}}:=(\n{G}_{\mathrm{FE}} \subset \n{G}_{\mathrm{FE}} \times \n{G}_{\mathrm{FE}})^{\dagger}$ is the overconvergent neighbourhood of the diagonal map. Therefore,  $\ob{D}(\n{G}_{\mathrm{FE}} \backslash \n{G}_{\mathrm{FE}}^{\dR}  / \n{G}_{\mathrm{FE}} )$ is a suitable category of \textit{analytic Harish--Chandra bimodules} for the group $\n{G}_{\mathrm{FE}}$.

We wish to analyse the action of Harish--Chandra bimodules on $\ob{D}(X^{\HT})$, and from this to construct a natural structure of \textit{Higgs--Sen modules} on sheaves over the Hodge--Tate stack. As hinted from the cartesian square \eqref{eqSquareAbsoluteHTStackSection6}, this construction should only depend on the bottom map, so we will discuss this universal case. Let $X$ be a partially proper rigid space and $\n{G}/X$ a partially proper smooth rigid analytic group over $X$. Let $\n{G}^{\dagger}\subset \n{G}$ be the overconvergent neighbourhood at the unit section. The conjugate action of $\n{G}$ on itself leaves the subgroup $\n{G}^{\dagger}$ stable, and therefore we can form the semidirect product $\n{G}^{\dagger}\rtimes \n{G}$. We have a natural multiplication morphism of groups $m\colon \n{G}^{\dagger}\rtimes \n{G}\to \n{G}$, such that the map 
\[
(\ob{pr}_G,m)\colon  \n{G}^{\dagger}\rtimes \n{G}\to \n{G}\times_X \n{G}, \;\; (a,b) \mapsto (b,ab)
\]
induces an isomorphism of groups on Berkovich spaces
\begin{equation}\label{eqDifferentPresentationHarishChandraBiModules}
 \n{G}^{\dagger}\rtimes \n{G}\xrightarrow{\sim} (\n{G}\subset \n{G}\times_X \n{G})^{\dagger} \cong \n{G}\times_{\n{G}^{\dR/X}} \n{G}. 
\end{equation}
Under this equivalence, the projections $\ob{pr}_i\colon \n{G}\times_{\n{G}^{\dR/X}} \n{G}\to \n{G}$ identify as $\ob{pr}_1\cong \ob{pr}_G$ and $\ob{pr}_2 \cong m$.   Furthermore, the diagonal map $\n{G}\to \n{G}\times_{\n{G}^{\dR/X}} \n{G}$ is identified with the natural morphism of groups $\n{G}\to \n{G}^{\dagger}\rtimes \n{G}$ given by $(1,\ob{id}_G)$. Thus, the natural sequence of groups $\n{G}\xrightarrow{(1,\ob{id}_G)} \n{G}^{\dagger}\rtimes \n{G} \xrightarrow{m} \n{G}$ gives rise morphisms of classifying stacks over $X$
\begin{equation}\label{eqMapsOfClassifyingStacks}
B\n{G} \xrightarrow{i} B( \n{G}^{\dagger}\rtimes \n{G}) \xrightarrow{\phi} B \n{G}.
\end{equation}
Sheaves on  $ B( \n{G}^{\dagger}\rtimes \n{G}) $  consist of $(\f{g},\n{G})$-modules (with $\f{g}=\Lie (\n{G})$) for which the action of $\f{g}$ is locally analytic.  Indeed, we have a natural morphism of groups $\widehat{\n{G}}\rtimes \n{G} \to \n{G}^{\dagger}\rtimes \n{G}$, with $\widehat{\n{G}}$ the formal completion of $\n{G}$ at $1$, so that the morphism of classifying stacks
\begin{equation}\label{eqMorphismsGroupsdeRhamStacksSemidirectProduct}
B(\widehat{\n{G}}\rtimes \n{G})\to B(\n{G}^{\dagger}\rtimes \n{G})
\end{equation}
is fibered in the map $\n{G}^{\dagger}/\widehat{\n{G}}\to X$ which is the algebraic de Rham stack of $\n{G}^{\dagger}$, and in particular prim with contractible cohomology by the Poincar\'e lemma. This implies that the pullback along \eqref{eqMorphismsGroupsdeRhamStacksSemidirectProduct} is fully faithful as desired.

The pullback map $\phi^*\colon \ob{D}( B \n{G}^{\dagger})\to \ob{D}(B( \n{G}^{\dagger}\rtimes \n{G}) )$ sends a $\n{G}$-representation to the underlying analytic $(\f{g},\n{G})$-module obtained by taking derivations along the action of $\n{G}$, and  sheaves on the $BG$-pointed stack $B(\widehat{\n{G}}\rtimes \n{G})$  can be described  as $\ob{LMod}_{U(\f{g})}(\ob{D}(B\n{G}))$, where $\f{g}$ is endowed with the adjoint action.

The following theorem describes the relationship between the stack $\widehat{X^{\diamond}}$ and the Hodge--Tate stack in the mixed case.

\begin{theorem}
\label{sec:mixed-case:-smooth-mixed-case}
Let $X$ be a smooth partially proper rigid space over $K$.  The following hold true:
  \begin{enumerate}
  
  \item  Let $\pi_{X,1}\colon X^{\HT}\to X\times K^{\HT}$ be the natural projection map. Then $\pi_{X,1}$ is a $B\bb{T}_{X}^{\dagger}\{1\}$-torsor, where $\{1\}$ refers to a Tate twist over  $K^{\HT}$. In particular, applying Cartier duality, any module $M\in \ob{D}(X^{\HT})$ is naturally endowed with an $\s{O}_{X^{\HT}}$-linear action of $\Sym_{\s{O}_X}(\bb{T}_{X}\{1\})$, or equivalently, with the structure of a $\Omega^{1}_X\{-1\}$-Higgs module. Furthermore, this action is compatible with that of \Cref{Theo:Geometric-Case-Hodge-TateStack-Decompletion} after base change along perfectoids $D\to \bb{Q}_p^{\HT}$, and produces a natural descent of the  geometric Sen operator  to $X\times K^{\HT}$.

\item The cotangent complex $\mathbb{L}_{X^{\HT}/X\times K^\HT}$ is canonically isomorphic to $\Omega^1_{X/K}\{-1\}[-1]$, with $\Omega^1_{X/K}\{-1\}:=\pi_X^\ast \Omega^1_{X/K}\otimes_{\s{O}_{X^{\HT}}}\s{O}_{X^{\HT}}\{-1\}$ and $\pi_X\colon X^{\HT}\to X$ is the natural morphism. Furthermore, the cotangent complex $\mathbb{L}_{X^\HT/X}$ is canonically an extension
  \[
    0\to \s{O}_{X^{\HT}}[-1] \to \mathbb{L}_{X^\HT/X}\to \mathbb{L}_{X^\HT/(X\times K^\HT)}\cong \Omega^1_{X/K}\{-1\}[-1] \to 0,
  \]
  whose pullback along $\alpha_{X}^{\HT}$ to $\widehat{X^\diamond}$ identifies with the inverse $\mathrm{FE}_X^{-}$ of the Faltings' extension (up to twisting by $\{1\}$ and shifting). More precisely, under the equivalence of perfect complexes of part (1)  we have 
  \[
  \mathrm{FE}_X^{-}\{-1\} =  \mathbb{L}_{X^\HT/X}[1]. 
  \]

\item  Let $\f{g}_{\mathrm{FE}}$ denote the dual of the Lie algebra of $\n{G}_{\mathrm{FE}}$ over $X$ endowed with its adjoint action. Then pullback of the dual $\f{g}^*_{\mathrm{FE}}$ along the map $X^{\HT}\to B \n{G}_{\mathrm{FE}}$ is naturally isomorphic to $\bb{L}_{X^{\HT}/X}[1]$. We let $\f{g}_{X,\ob{FE}}$ denote the pullback of $\f{g}_{\mathrm{FE}}$ to $X^{\HT}$ considered as an $\s{O}_{X^{\HT}}$-linear Lie algebra.

\item   Consider the commutative diagram of Gelfand stacks with cartesian squares
\[
\begin{tikzcd}
X^{\HT}  \ar[r] \ar[d] & Y_1 \ar[r,"p"]  \ar[d] &  Y_0  \ar[r,"m'"] \ar[d] & X^{\HT} \ar[r,"\pi_X"] \ar[d] & X \ar[d] \\ 
B\n{G}_{\mathrm{FE}} \ar[r,"{(1,\ob{id})}"] & B(\widehat{\n{G}}_{\mathrm{FE}} \rtimes \n{G}_{\mathrm{FE}}) \ar[r] & B(\n{G}^{\dagger}_{\mathrm{FE}} \rtimes \n{G}_{\mathrm{FE}}) \ar[r,"m"] & B\n{G}_{\mathrm{FE}} \ar[r] & B\n{G}_{\mathrm{FE}}^{\dR/X}
\end{tikzcd}
\]
Then the following hold:

\begin{enumerate}

\item The map $p$ is prim and $1\xrightarrow{\sim}p_*1$, in particular the pullback $p^*$ is fully faithful. 

\item  There is a  natural equivalence of categories  
\[
\ob{D}(Y_1)\cong \ob{LMod}_{U_{\s{O}_{X^{\HT}}}(\f{g}_{X,\ob{FE}})}(\ob{D}(X^{\HT})).
\]
In particular, the functor $p^*\circ m^*\colon \ob{D}(X^{\HT})\to  \ob{LMod}_{U_{\s{O}_{X^{\HT}}}(\f{g}_{X,\ob{FE}})}(\ob{D}(X^{\HT}))$ naturally promotes a sheaf on $X^{\HT}$ to a $\f{g}_{X,\ob{FE}}$-Lie algebra representation. Furthermore, the restriction along the morphism of Lie algebra $\bb{T}_{X}\{1\}\to \f{g}_{X,\ob{FE}}$ recovers the structure of $\Omega^1_X\{-1\}$-Higgs module of part (2).

\item  There is a natural equivalence of endofunctors of $\ob{D}(X^{\HT})$
\[
\pi_X^*\pi_{X,*}\cong \Hom_{U_{\s{O}_{X^{\HT}}}(\f{g}_{X,\ob{FE}})}(\s{O}_{X^{\HT}},-).
\]
In particular, the Koszul resolution of  $\s{O}_{X^{\HT}}$ as a  $\f{g}_{X,\ob{FE}}$-representation produces a natural filtration of $\pi_X^*\pi_{X,*}$  depicted as a chain complex
\begin{equation}\label{eqFiltrationbetaStarMonadMixedCase}
\pi_X^*\pi_{X,*}\otimes M = [M\to M\otimes_{\s{O}_{X^{\HT}}} \mathrm{FE}_X^{-}\{-1\}\to \cdots \to  M\otimes_{\s{O}_{X^{\HT}}} (\bigwedge^{\dim_X +1}\mathrm{FE}_X^{-}) \{-\dim X-1\}]
\end{equation}

\item Let $M\in \widehat{X^{\diamond}}$ be a sheaf in the essential image of $\alpha^{\HT,*}_{X}$, then $M$ has a natural  structure of $\f{g}_{\mathrm{FE}_X}$-representation whose restriction to $\bb{T}_{X}\{1\}$ is its canonical Higgs field structure. 

\end{enumerate}

The underlying $\f{g}_{X,\ob{FE}}$-module structure of sheaves on $X^{\HT}$ or $\widehat{X^{\diamond}}$ is called the \textit{canonical Higgs--Sen module structure}.
    \end{enumerate}
\end{theorem}

\begin{proof}
  Part (1) follows from the pullback diagram \eqref{eqSquareAbsoluteHTStackSection6} and Cartier duality. Part (2) is a special case of \Cref{pullback-cotangent-complex-vb}, noticing that the map $s\colon X\to X/ (\bb{T}_{X}\rtimes \bb{G}_m)^{\dR/X}$ arises from the ($\{-1\}$-twist) of the Faltings' extension.  Part (3) follows from the compatibility of cotangent complexes with base change, and the fact that the cotangent complex of $B \n{G}_{\mathrm{FE}}\to B \n{G}_{\mathrm{FE}}^{\dR/X}$ is the same as the cotangent complex of $B\n{G}_{\mathrm{FE}}\to X$ which is $\f{g}_{\mathrm{FE}}^{\vee}[-1]$ (cf. \Cref{PropositionCotangentComplexesClassifyingStacks} (1)).  
It is left to show part (4). Point (a) follows from the analogous statement for the map $B(\widehat{\n{G}}_{\mathrm{FE}}\rtimes \n{G}_{\mathrm{FE}})\to B(\n{G}^{\dagger}_{\mathrm{FE}}\rtimes \n{G}_{\mathrm{FE}})$, which is fibered in the algebraic de Rham stack $\n{G}_{\mathrm{FE}}^{\dagger, \ob{alg}-\dR/X}= \n{G}_{\mathrm{FE}}^{\dagger}/\widehat{\n{G}}_{\mathrm{FE}}$ which is prim over $X$, and has contractible cohomology by the Poincar\'e lemma. For part (b), we note that, by base change, the map $X^{\HT}\to Y_1$ is suave as then so is the map $X\to B\widehat{\n{G}}_{\mathrm{FE}}$, hence by the monadicity theorem one formally deduces that the $*$-pull and $\sharp$-push adjunction along  $X^{\HT}\to Y_1$ produces an equivalence
\[
\ob{LMod}_{A}(\ob{D}(X^{\HT})) \cong \ob{D}(Y_1)
\]
for a suitable algebra $A\in \ob{Alg}(\ob{D}(X^{\HT}))$. Since this monad arises as the pullback of the monad associated to the suave map $B\n{G}_{\mathrm{FE}}\to B(\widehat{\n{G}}_{\mathrm{FE}}\rtimes \n{G}_{\mathrm{FE}})$, the equivalence  
\[
B(\widehat{\n{G}}_{\mathrm{FE}}\rtimes \n{G}_{\mathrm{FE}})\cong \ob{LMod}_{U(\f{g}_{\mathrm{FE}})}(\ob{D}(\n{G}_{\mathrm{FE}}))
\] yields the natural isomorphism $A\cong U_{\s{O}_{X^{\HT}}}(\f{g}_{X,FE})$. 
 Recall that the geometric Sen operators  of part (1) consists of a natural  action of the  abelian Lie algebra $\bb{T}_{X}\{1\}$ on  the objects  of $\ob{D}(X^{\HT})$, and that these  arise from the  $B\bb{T}^{\dagger}_X\{1\}$-torsor structure of the map $\pi_{X,1}$ via Cartier duality. Explicitly, the action of   $\ob{D}(B\bb{T}^{\dagger}_X\{1\})$ arise from convolution  and can be equivalently encoded from the pullback along the action map 
\[
a\colon  B\bb{T}^{\dagger}_X\{1\}\times_{X\times  K^{\HT}} X^{\HT} \to X^{\HT}.
\]
The map $a$ fits in a pullback square
\[
\begin{tikzcd}
B\bb{T}^{\dagger}_X\{1\}\times_{X\times K^{\HT}} X^{\HT} \ar[r,"a"]  \ar[d] & X^{\HT} \ar[d] \\ 
B ( \bb{T}^{\dagger}_X\{1\}\times_{X\times K^{\HT}}  \bb{T}_{X}\{1\}) \ar[r,"m"] & B \bb{T}_{X}\{1\}
\end{tikzcd}
\]
where the bottom map $m$ is the natural addition. This implies that the restriction of the $\f{g}_{X,\ob{FE}}$-module structure of the objects of $\ob{D}(X^{\HT})$ to a module over the Lie algebra $\bb{T}_{X}\{1\}\subset \f{g}_{X,\ob{FE}}$ is precisely the canonical Higgs structure,  proving point (b).

 For point (c), we need to look more explicitly at the action of (analytic) Harish--Chandra bimodules of $\n{G}_{\mathrm{FE}}$ on $\ob{D}(X^{\HT})$. By \eqref{eqDifferentPresentationHarishChandraBiModules} we have an equivalence
\[
B(\n{G}^{\dagger}_{\mathrm{FE}}\rtimes \n{G}_{\mathrm{FE}} )\cong B\n{G}_{\mathrm{FE}}\times_{B\n{G}_{\mathrm{FE}}^{\dR/X}} B \n{G}_{\mathrm{FE}}
\]
where we have made explicit the projections and diagonal map.  By proper base change, the kernel given by the unit $1\in B(\n{G}^{\dagger}_{\mathrm{FE}}\rtimes \n{G}_{\mathrm{FE}} )$ acts on $\ob{D}(X^{\HT})$ as the map $\pi_X^*\pi_{X,!}=\pi_X^*\pi_{X,\ast}$ where we identified $\pi_{X,!}=\pi_{X,*}$ thanks to the trivialization of the codualizing sheaf of the map $X^{\HT}\to X$ (as it is $0$-proper and the diagonal is $0$-affine, see \Cref{sec:deriv-berk-spac-f-lower-shriek-via-0-affine-stuff}). The relative de Rham complex of the map $a=(1,\ob{id})\colon B\n{G}_{\mathrm{FE}}\to B(\n{G}^{\dagger}_{\mathrm{FE}}\rtimes \n{G}_{\mathrm{FE}} )$ produces a resolution of the unit $1\in \ob{D}(B(\n{G}^{\dagger}_{\mathrm{FE}}\rtimes \n{G}_{\mathrm{FE}} ))$ of the form
\[
0\to 1\to a_* 1_{B\n{G}_{\mathrm{FE}}} \to a_* \f{g}_{\mathrm{FE}}^{\vee}\to \cdots \to a_* \bigwedge^{\dim X+1} \f{g}_{\mathrm{FE}}^{\vee}\to 0
\] 
which is nothing but the de Rham complex of  $\n{G}^{\dagger}_{\mathrm{FE}}$ seen as a representation of  $\n{G}^{\dagger}_{\mathrm{FE}}\rtimes \n{G}_{\mathrm{FE}}$. Acting by kernels on $\ob{D}(X^{\HT})$ one obtains a filtration  of $\pi_X^*\pi_{X,*}$ of the form \eqref{eqFiltrationbetaStarMonadMixedCase} which by construction is precisely the de Rham complex of a module $M\in \ob{D}(X^{\HT})$ with respect to the action of the Lie algebra $\f{g}_{X,\ob{FE}}$. This produces the desired identification $\pi_X^*\pi_{X,*}= \Hom_{U_{\s{O}_{X^{\HT}}}(\f{g}_{X,\ob{FE}})}(\s{O}_{X^{\HT}},-)$.

Finally, point (d) follows from point (b) and the fully faithful embedding of $\alpha^{\HT,*}_{X}\colon \ob{D}(X^{\HT}) \hookrightarrow \ob{D}(\widehat{X^{\diamond}})$ deduced from   (1) by the projection formula.  
  \end{proof}

\begin{lemma}\label{pullback-cotangent-complex-vb}
  Let $S$ be a Gelfand stack, $\mathcal{W}$ be a vector bundle over $S$, $\mathcal{V}=\mathcal{W}(1)$ on $S\times B\Gm$, and $Y$ be a Gelfand stack with a map $h\colon Y\to B_{S\times B\Gm} \mathcal{V}$.
  Let $\widetilde{Y}\to Y$ be the inverse of the $\mathcal{V}$-torsor over $Y$ coming via pullback along $h$ of the unit section $S\times B\Gm \to B_{S\times B\Gm} \mathcal{V}$.
  Then the sheaf of affine functions on $\widetilde{Y}$ is an extension of $\mathcal{V}^\vee$ by $\s{O}_{Y}$, which identifies with the pullback along $h$ of $\mathbb{L}_{B_{S\times B\Gm} \mathcal{V}/S}[1]$. 
\end{lemma}
As the proof will show, the extension
\[
  0\to \s{O}_Y\to h^\ast\mathbb{L}_{B_{S\times B\Gm}\mathcal{V}/S}[1]\to h^\ast\mathcal{V}^\vee\cong h^\ast\mathbb{L}_{B_{S\times B\Gm}(\mathcal{V})/S\times B\Gm}[1]\to 0\]
corresponds to the $\mathcal{V}$-torsor $\widetilde{Y}$ over $Y$.
\begin{proof}
  The sheaf of affine functions on $\widetilde{Y}$ will be constructed during the proof.
  We note that the statement can be pulled back from the case of (derived) schemes over $\Spec(\Q_p)$ (using the theory of \Cref{sec:analyt-tann-stacks}).
  For simplicity, we write $\Gm, \mathrm{GL}_n,\ldots$ now for the classical derived schemes over $\Spec(\Z)$.
  We can furthermore reduce to the universal case $S=B\mathrm{GL}_n$ and the universal $\mathcal{V}$-torsor $f\colon B(\GL_n\times \Gm)\to B_{B\mathrm{GL}_n\times \Gm}\mathcal{V}\cong B(V\rtimes (\mathrm{GL}_n\times \Gm))$, where $V$ is the standard representation of $\GL_n$ and $\Gm$ acts via multiplication.
  We also identify (classical) quasi-coherent sheaves on the occuring classifying stacks with comodules under the respective coalgebra, i.e., with algebraic representations.
  The morphism $f$ is To\"en-affine and the relative Spec of the algebra $\mathrm{Sym}^\bullet(V^\vee)$.
  Here, the $V\rtimes(\GL_n\times \Gm)$-action in the inverse of the universal torsor is as follows:
  $\GL_n\times \Gm$ act via the natural action coming from $V$, and the $V$-action comes from the map $\Sym^\bullet(V^\vee)\to \Sym^\bullet(V^\vee)\otimes \Sym^\bullet(V^\vee)\cong \Sym^\bullet(V^\vee\oplus V^\vee)$ induced by the antidiagonal map $V^\vee\to V^\vee\oplus V^\vee,\ v\mapsto (v,-v)$ (as we pass to the inverse torsor).
  The sheaf of affine functions of the inverse of the universal torsors is the subrepresentation $\Sym^{\leq 1}(V^\vee)=\Z\oplus V^\vee$, which is naturally the extension $0\to \Z\to \Sym^{\leq 1}(V^\vee)\to V^\vee\to 0$ (this sheaf of affine functions can then be pulled back along any map to $B(V\rtimes (\GL_n\times \Gm))$).
  We want to identify $\Sym^{\leq 1}(V^\vee)$ naturally with $M:=\mathbb{L}_{B(V\rtimes (\GL_n\times \Gm))/B\GL_n}[1]$, as this will imply the general case by pullback.
  Now, $M$ identifies naturally with the cokernel of the map
  \[
    \mathfrak{gl}_n^\vee\to V^\vee\oplus \mathfrak{gl}_n^\vee\oplus \mathfrak{h}^\vee
  \]
  with $\mathfrak{h}:=\mathrm{Lie}(\Gm)$.
  Via projection, this is isomorphic to $V^\vee\oplus \mathfrak{h}^\ast$, as a $\GL_n\times \Gm$-representation for the natural action.
  The $V$-action needs a bit more care.
  Identifying $V\oplus \mathfrak{h}$ with $\mathrm{Lie}(V\rtimes\Gm)$, one checks that the adjoint action of $v\in V$ on $(w,x)\in V\otimes \mathfrak{h}$ is given by $v\cdot (w,x)=(-xv+w,x)$, e.g., by a matrix calculation in the mirabolic.
  One checks directly that this implies that $\mathrm{Sym}^{\leq 1}(V^\vee)$ is isomorphic to $M$.
\end{proof}

We note that by decomposing a vector bundle on $S\times B\Gm$ into graded pieces and by pulling back the case of \Cref{pullback-cotangent-complex-vb} along the multiplication by $B\Gm$, one can generalize \Cref{pullback-cotangent-complex-vb} to any vector bundle $\mathcal{V}$ on $S\times B\Gm$.
    
    \begin{remark}\label{conclusion-ht-stack-valid-smooth-not-pp2}
    As in \Cref{conclusion-ht-stack-valid-smooth-not-pp}, \Cref{sec:mixed-case:-smooth-mixed-case} does not need the assumption that $X$ is partially proper.
    \end{remark}

    \begin{remark}
     \label{sec:mixed-case:-smooth-remark-atiyah-classes}
     It seems plausible to us that one can interpret the Sen operator, Higgs field or Higgs-Sen module structure on v-vector bundles arising from the results of this section and the last two as instances of a general theory of \textit{Atiyah classes} for morphisms of Gelfand stacks (we thank David Hansen for sharing this viewpoint).
     We do not develop a theory of Atiyah classes in this paper, and content ourselves with the following remarks:
     \begin{enumerate}
     \item Let $f\colon X\to S$ be a morphism of Gelfand stacks, whose diagonal $\Delta\colon X\to X\times_S X$ is prim.
       By primness of $\Delta$ the comonad $\Delta^\ast\Delta_\ast$ on $\ob{D}(X)$ is given by tensoring with the algebra $\mathcal{H}_f:=\Delta^\ast\Delta_\ast(1)$.
       Let $\mathrm{pr}_1\colon X\times_S X\to X$ be the first projection.
       For $M\in \ob{D}(X)$ we have $M=\Delta^\ast(\mathrm{pr}_1^\ast(M))$, and thus we can deduce a $\mathcal{H}_f$-coaction $M\to M\otimes_{\mathscr{O}_X} \mathcal{H}_f$ induced by the unit $\mathrm{pr}_1^\ast(M)\to \Delta_\ast(\Delta^\ast \mathrm{pr}_1^\ast(M))$.
       We note that $\mathcal{H}_f$ is naturally a bialgebra because $\Delta_\ast$ is lax-symmetric monoidal.
     \item Assume that $f\colon X=\GSpec(B)\to S=\GSpec(A)$.
       Then $\mathcal{H}_f=B\otimes_{B\otimes_A B}B$ is the (solid) Hochschild homology of $B$ relative to $A$.
       Now, $\mathcal{H}_f$ is equipped with a canonical decreasing filtration $\mathrm{Fil}^\bullet$, cf.\ proof of \Cref{sec:cotang-compl-morph-cotangent-complex-via-p-completion}, such that $\mathcal{H}_f/\mathrm{Fil}^2\cong \mathbb{L}_{B/A}[1]\oplus B$.
       In particular, the $\mathcal{H}_f$-coaction on some $M\in \ob{D}(X)$ induces a morphism $M\to M\otimes_{B}\mathbb{L}_{B/A}[1]$, as desired for Atiyah classes, cf.\ \cite{kuhn2024atiyah}.
     \item Assume that $f\colon X=BG\to S=\GSpec(\Q_p)$ for some affinoid Gelfand group stack $G=\GSpec(\mathscr{O}_G)$.
       Then $\mathcal{H}_f\cong \mathscr{O}_G$ with the conjugation action, and for $M\in \ob{D}(X)$, seen as a comodule under $\mathscr{O}_G$, its $\mathcal{H}_f$-coaction is the given coaction morphism $M\to M\otimes \mathcal{H}_f=M\otimes \mathscr{O}_G$, which is a morphism of $\mathscr{O}_G$-comodules for the diagonal coaction.
       If $G$ is the overconvergent neighborhood of the identity of a rigid-analytic group over $\Q_p$ (or we work in analytic stacks and $G$ is an affine algebraic group), then by evaluating at the first infinitesimal neighborhood of the identity we get a natural morphism $\mathscr{O}_G\to \mathfrak{g}^\vee\oplus \Q_p$ with $\mathfrak{g}:=\mathrm{Lie}(G)$ and thus a natural morphism
       \[
         M\to M\otimes \mathfrak{g}^\vee,
       \]
       which encodes the action of the Lie algebra on $M$.
       As $\mathbb{L}_{BG/\ast}\cong \mathfrak{g}^\vee[-1]$, this morphism should again be an instance of a general theory of Atiyah classes.
     \item Given a morphism $g\colon X'\to X$ of two Gelfand stacks $f\colon X\to S,\ f'\colon X' \to S$ over $S$, we get a natural morphism $g^\ast \mathcal{H}_f\to \mathcal{H}_{f'}$ of coalgebras (even bialgebras) and for $M\in \ob{D}(X)$ the natural $\mathcal{H}_{f'}$-coaction on $g^\ast M$ is given by the composition
       \[
         g^\ast M\to g^\ast(M\otimes_{\mathscr{O}_Y}\mathcal{H}_f)\cong g^\ast M\otimes_{\mathscr{O}_{X'}} g^\ast \mathcal{H}_f\to g^\ast M\otimes_{\mathscr{O}_{X'}} \mathcal{H}_{f'}, 
       \]
       where the first morphism is the pullback of the canonical $\mathcal{H}_f$-coaction on $M$.
     \item If $f\colon X\to S$ is a gerbe banded by some abelian group Gelfand stack $A$ over $S$ such that $BA\to S$ is prim, we can use the isomorphism
       \[
         X\times_S X\cong X\times BA
       \]
       to see that the construction of the canonical Sen operator/the canonical Higgs field in \Cref{TheoremPropertiesHodgeTateStack,Theo:Geometric-Case-Hodge-TateStack-Decompletion,sec:mixed-case:-smooth-mixed-case} are particular instances of the general construction in (1).
     \end{enumerate}
\end{remark}

\subsection{Decompleting torsors on the $p$-adic twistor stack}
\label{sec:decompl-tors}    
Let $G$ be a $p$-adic Lie group. In this section, we apply the general formalism of decompletions of torsors of \Cref{sec:decompl-tors-1} to prove that pro-\'etale $\underline{G}$-torsors on a partially proper smooth rigid analytic variety decomplete uniquely into $G^{\la}$-torsors of their associated Hodge--Tate stacks, and furthermore any such torsor produces a natural $G^{\la}$-torsor over the  $p$-adic twistor stack.

\begin{theorem}\label{CoroEquivalenceTorsor}
  Let  $X$ be a smooth partially proper rigid space over $\Q_p$.  Let $G$ be a  $p$-adic Lie group.
  The following hold.
\begin{enumerate}

\item  There are  natural equivalences of anima
\begin{equation}\label{eqiho3ampnpoaef}
\ob{Map}_{\GSpec(\Q_p)}(X^{\HT}, BG^{\la}) \xrightarrow{\sim}  \ob{Map}_{/\GSpec(\Q_p)}  (\widehat{X^{\diamond}}, BG^{\la})   \xrightarrow{\sim}  \ob{Map}_{/\Marc(\Q_p)}  (X^{\diamond}, B\underline{G}),
\end{equation}
where the first map arises from the pullback along $f\colon \widehat{X^{\diamond}}\to X^{\HT}$ and the second map is perfectoidization. 

\item  For any  Berkovich-\'etale map $U \subset \FF_{X^{\diamond}}^{\diamond}$, let $X^{\Div^1}|_U \to X^{\Div^1}$ and $\FF_{X^{\diamond}}|_U\to  \FF_{X^{\diamond}}$ be the associated  Berkovich-\'etale map of  Gelfand stacks.  Then there is a natural equivalence of anima 
\[
\ob{Map}_{\GSpec(\Q_p)}(X^{\Div^1}|_U, BG^{\la}) \xrightarrow{\sim}  \ob{Map}_{/\GSpec(\Q_p)}  (\FF_{X^{\diamond}}|_U, BG^{\la}) 
\]
given by pullback along $\FF_{X^{\diamond}}|_U\to X^{\Div^1}|_U$.

\item There is a natural morphism of anima 
\begin{equation}\label{eqkwnih3ianwlfdasf}
\ob{Map}_{/\Marc(\F_p)}  (X^{\diamond}, B\underline{G})\to  \ob{Map}_{/\GSpec(\Q_p)}(\FF_{X^{\diamond}}, BG^{\la}),
\end{equation}
whose composition with perfectoidization is the functor sending  a $\underline{G}$-torsor $\widetilde{X}^{\diamond}\to X^{\diamond}$ to the $\underline{G}$-torsor $\FF_{\widetilde{X}^{\diamond}}^\diamond \to \FF_{X^{\diamond}}^\diamond$. 
\end{enumerate}
A similar result holds for $X$ smooth partially proper over a perfectoid space $D$ over $\Q_p$, using the relative $p$-adic twistor stack $X^{\Div^1/\FF_D}$ instead.
\end{theorem}

\begin{proof}
The statements to be proved are local for the rigid \'etale topology on $X$, so we can assume that $X=\Gm^d$, for some $d\geq 0$. Part (1) follows from \Cref{sec:decompl-tors-2-beta-x-z-an-equivalence}, \Cref{sec:decompl-tors-3-decompletion-of-torsors-perfectoid-case} together with the explicit description of $X^\HT$ in terms of locally analytic vectors as in \Cref{thm:summary-main-results-sen}. Similarly, the explicit description of $X^\prism$ and \Cref{sec:decompl-tors-2-beta-x-z-an-equivalence} give part (2). 

We now prove part (3).
Let $f\colon X^\diamond\to B\underline{G}$ be a map.
Applying the functor $\FF_{-}$ from arc-stacks over $\F_p$ to Gelfand stacks, we get $g\colon\FF_{X^\diamond}\to \FF_{B\underline{G}}\cong B\underline{G}\times_{\GSpec(\Q_p)}\FF_{\Marc(\F_p)}$.
Composing with the projection to $B\underline{G}$, and the morphism $B\underline{G}\to BG^\la$, we get the desired natural map
\[
\ob{Map}_{/\Marc(\F_p)}  (X^{\diamond}, B\underline{G})\to  \ob{Map}_{/\GSpec(\Q_p)}(\FF_{X^{\diamond}}, BG^{\la}),
\]
and by construction it has the desired description.
\end{proof}

\begin{remark}
  \label{sec:decompl-tors-p-non-partially-proper}
  Following \Cref{conclusion-ht-stack-valid-smooth-not-pp2}, we note that for the Hodge--Tate stack, i.e.\ point (1) of \Cref{CoroEquivalenceTorsor}, the hypothesis that $X$ is partially proper is not needed.
\end{remark}

 \begin{remark}\label{TheoLocallyAnalyticComparisonGeneralbis}
 Let $X$ be a partially proper smooth rigid space over $\Q_p$, let $G$ be a  $p$-adic Lie group and let $\widetilde{X}^{\diamond}\to X^{\diamond}$ be a   pro-\'etale $\underline{G}$-torsor.  By \Cref{CoroEquivalenceTorsor}, this gives rise to a commutative diagram
\[\begin{tikzcd}
	{\FF_{X^\diamond}}  & {X^{\Div^1}}   \\
	{BG^\cont} & {BG^\la},
	\arrow["\alpha_X", from=1-1, to=1-2]
	\arrow["f", from=1-1, to=2-1]
	\arrow["g",from=1-2, to=2-2]
	\arrow["h", from=2-1, to=2-2]
      \end{tikzcd}\]
      so that the $G^\cont$-torsor on $\FF_{X^\diamond}$ determined by the left vertical map is $\FF_{\widetilde{X}}$. 
 By \Cref{sec:mixed-case:-smooth-mixed-case1}, $\alpha_{X,\ast} 1=1$. Recall also that $h$ is cohomologically proper and that $h_\ast=h_!$ computes derived $G$-locally analytic vectors (\cite[Proposition 6.1.5]{jacinto2023solid}).  Let  $\n{Y}\to X^{\Div^1}$ be the $G^{\la}$-torsor obtained from the right map.  Then $f_! 1 = \Gamma_c(\FF_{\widetilde{X}^{\diamond}},\s{O})$ computes the compactly supported cohomology of $\FF_{\widetilde{X}^{\diamond}}$, endowed with its natural $G^{\cont}$-action. On the other hand, $g_! 1=\Gamma_c(\n{Y},\s{O}_{\n{Y}})$ computes the compactly supported cohomology of $\n{Y}$ seen as a locally analytic representation of $G$. The commutativity of the diagram yields
 \[
\Gamma_c(\n{Y},\s{O}_{\n{Y}})= g_!1=g_! \alpha_{X,*}1 = h_!f_! 1 = (\Gamma_c(\FF_{\widetilde{X}^{\diamond}},\s{O}_{\n{T}}))^{G-\la}.
 \] 
In other words, $!$-pushforward along $g$ computes the locally analytic vectors of the $!$-pushforward along $f$. When applied to torsors arising from (toroidal compactifications of) Shimura varieties, we see that the cohomology with compact supports of the  natural $G^{\la}$-torsor over  $X^{\Div^1}$  will compute locally analytic vectors of completed cohomology as studied in \cite{MR4390302}. This  property was one of the major expectations for the theory of  analytic prismatization. 
 
Taking $G^{\la}$-cohomology,  we also deduce from this commutative square the isomorphism
 \[ \Gamma(X^{\Div^1},\s{O})= \Gamma(G^\la, \Gamma(\FF_{\widetilde{X}^{\diamond}},\s{O})^{G-\la}),
 \]
 expressing the cohomology of the $p$-adic twistor stack of $X$ in terms of the $G$-locally analytic vectors of the (derived) global sections of $\FF_{\widetilde{X}^{\diamond}}$. Of course, an analogous statement can be formulated for the Hodge--Tate stack, replacing $\FF_{\widetilde{X}^{\diamond}}$ by $\widehat{\widetilde{X}^{\diamond}}$, and for partially proper smooth rigid spaces over perfectoid spaces, by taking the relative $p$-adic twistor and Hodge--Tate stacks.
\end{remark}

In the rest of the section we explain in more detail the relationship of decompletions of torsors of \Cref{CoroEquivalenceTorsor}, and the work of Berger--Colmez \cite{BergerColmezII} and Pan  \cite{MR4390302}.  Following \Cref{sec:arithm-case:-q_ppr,sec:geometric-case,sec:mixed-case}, we will separate the discussion in the arithmetic, geometric, and mixed cases. We begin with the construction  of the \textit{Sen operator}  of a torsor. 

\begin{definition}\label{DefinitionSenOperatorTorsor}
Let $D$ be a perfectoid space and $X$ a partially proper smooth rigid space over $D$. Let $G$ be a $p$-adic Lie group and  $\beta^{\diamond}\colon \widetilde{X}^{\diamond}\to X^\diamond$ a pro-\'etale $\underline{G}$-torsor. Let $\beta\colon X^{\HT/D}\to BG^{\la}$ be the $G^{\la}$-torsor associated to $\widetilde{X}^{\diamond}$ by \Cref{CoroEquivalenceTorsor}. The \textit{geometric Sen operator of the torsor $\beta$} is the induced morphism of (shifted) cotangent complexes $ \beta^* \bb{L}_{BG^{\la}}[1] \to \bb{L}_{X^{\HT}/X}[1]$, written as 
\[
\theta_{\beta}^{\ob{geo}}\colon  \f{g}^{\vee}_{\beta}:=\beta^* \f{g}^{\vee}\to \Omega^{1}_{X/D}\{-1\},
\]
where we identified  $\bb{L}_{X^{\HT}/X}[1]\cong \Omega^{1}_{X/D}\{-1\}$ via \Cref{Theo:Geometric-Case-Hodge-TateStack-Decompletion} (2). 

 Similarly, if $X$ is a partially proper smooth rigid space over $K$, a finite extension of $\bb{Q}_p$, the \textit{mixed Sen operator of the torsor $\beta$} is the morphism
\[
\theta_{\beta}^{\ob{mix}}\colon \f{g}^{\vee}_{\beta}\to \ob{FE}^{-}_X\{-1\}
\]
obtained as the $1$-shift of maps of cotangent complexes of $\beta$ (where we have identified $\bb{L}_{X^{\HT}/X}[1]\cong \ob{FE}_X^{-}\{-1\}$ via \Cref{sec:mixed-case:-smooth-mixed-case} (2)).   

In the case $X=\GSpec(K)$, we have  $\ob{FE}^{-}_X\{-1\}=\s{O}_{\HT}$ and we write $\theta^{\ob{arith}}_X$ for $\theta_{\beta}^{\ob{mix}}$. We call the map 
\[
\theta_{\beta}^{\ob{arith}}\colon \f{g}^{\vee}_{\beta}\to \s{O}_{X^{\HT}}
\]
the \textit{arithmetic Sen operator of the torsor $\beta$}. 
\end{definition}

One of the main properties of the Sen operators of torsors is that their duals are morphisms of Lie algebras, and that they determine the action of the Sen operators of \Cref{Theo:Geometric-Case-Hodge-TateStack-Decompletion,sec:mixed-case:-smooth-mixed-case} on sheaves that arise from pullbacks along the map $\beta\colon X^{\HT}\to BG^{\la}$, cf. \cite[Theorem 1.0.4]{camargo2022geometricsen}. Proving this fact requires a reinterpretation of the Sen operators, which can be explained in a general framework as follows:

\begin{construction}\label{Construction:SenOperator-Torsor}
Let $S$ be a Gelfand stack and let $G\to S$ be a $\dagger$-formally smooth group over $S$. Let $s\colon S\to B_SG^{\dR/S}$ be a $G^{\dR/S}$-torsor over $S$ and consider the pullback 
\[
\begin{tikzcd}
 S^{\ob{gerbe}} \ar[r] \ar[d,"s'"] & S \ar[d,"s"] \\ 
 B_S G \ar[r] & B_S G^{\dR/S}. 
\end{tikzcd}
\]
Let $G^{\dagger}=\ob{fib}(G\to G^{\dR/S})$ denote the overconvergent neigbourhood at $1$ in $G$. We have a sequence of groups over $S$ given by $G\xrightarrow{(1,\ob{id}_G)} G^{\dagger}\rtimes G\xrightarrow{m} G$ (with $m$ the multiplication). Consider the following commutative  diagram with pullback squares
\begin{equation}\label{eqBigCartesianSquareGerbesTorsors111}
\begin{tikzcd}
S^{\ob{gerbe}}  \ar[r,"\iota"] \ar[d,"s'"] & \n{Y} \ar[r] \ar[d] & S^{\ob{gerbe}} \ar[r] \ar[d,"s'"] & S \ar[d,"s"] \\ 
 B_S G  \ar[r]& B_S(G^{\dagger}\rtimes G) \ar[r] & B_S G \ar[r] & B_S G^{\dR/S}. 
\end{tikzcd}
\end{equation}
The pullback along  $s'$ induces a group $\n{G}^{\dagger}$ over $S^{\ob{gerbe}}$  such that $\n{Y} = B_{S^{\ob{gerbe}}} \n{G}^{\dagger}$.  Let $H$ be now another group over $S$ and suppose we are given a map $\beta\colon S^{\ob{gerbe}}\to BH$. Taking the  precomposition with the morphism $\n{Y}\to S^{\ob{gerbe}}\to BH$,  one obtains a map of $S^{\ob{gerbe}}$-pointed stacks over $S^{\ob{gerbe}}$ 
\begin{equation}\label{eq:Morphism-Torsors-Gerbes}
\Phi_{\beta}\colon  \n{Y}\to BH\times_{S} S^{\ob{gerbe}},
\end{equation}
and therefore, by taking loops with respect to the  $S^{\ob{gerbe}}$-point, a natural morphism of group objects over $S^{\ob{gerbe}}$
\[
\Omega (\Phi_{\beta})\colon \n{G}^{\dagger}\to H_{\beta},
\]
where $H_{\beta}\to S^{\ob{gerbe}}$ is  the group of automorphisms of the $H$-torsor $\beta$. Since  $\n{G}^{\dagger}$ has trivial de Rham stack over $S^{\ob{gerbe}}$, the map $\Omega (\Phi_{\beta})$ factors uniquely as a morphism of groups
\[
\Omega (\Phi_{\beta})\colon \n{G}^{\dagger}\to H^{\dagger}_{\beta}. 
\]
\end{construction}

\begin{remark}\label{RemarkMapsCotangentComplexesAbstractFramework}
Keep the notation of \Cref{Construction:SenOperator-Torsor}, and suppose that all stacks involved admit a cotangent complex relative to $S$. Let $\beta\colon S^{\ob{gerbe}}\to B_S H$ be an $H$-torsor over $S^{\ob{gerbe}}$. The map $\Phi_{\beta}$ of \eqref{eq:Morphism-Torsors-Gerbes} can be obtained as the base change of $\beta$ along the map $S^{\ob{gerbe}}\to S$. Indeed since $G$ is $\dagger$-formally smooth over $S$ the map $G\to G^{\dR/S}$ is surjective, and the natural morphism of groups
\[
G^{\dagger}\rtimes G \to G\times_{G^{\dR/S}} G, \; \; (h,g)\mapsto (g, hg)
\]
is an isomorphism. This implies that 
\[
\n{Y} \cong  B(G\times_{G^{\dR/S}} G) \times_{B_S G^{\dR}} S = S^{\ob{gerbe}}\times_{S} S^{\ob{gerbe}}
\]
where the map $\n{Y}\to S^{\ob{gerbe}}$ of \eqref{eqBigCartesianSquareGerbesTorsors111} is identified  with the second projection.   From the previous, we can deduce that the morphism of cotangent complexes $\beta^* \bb{L}_{H/S}\to \bb{L}_{S^{\ob{gerbe}}/S}$ is equal to the pullback along $\iota\colon S^{\ob{gerbe}}\to \n{Y}$ of the map of cotangent complexes 
\begin{equation}\label{eqMapCotangentComplexesClassGroups223}
\Phi_{\beta}^*\bb{L}_{B_SH\times_S S^{\ob{gerbe}}/S^{\ob{gerbe}}}\to \bb{L}_{\n{Y}/S^{\ob{gerbe}}}
\end{equation}
 on $\n{Y}$.  Since $\Phi_{\beta}$ arises by passing to classifying stack in a group homomorphism $\n{G}^{\dagger}\to H_{\beta}$, \Cref{PropositionCotangentComplexesClassifyingStacks} (1) implies that the pullback along $\iota$ of the  map \eqref{eqMapCotangentComplexesClassGroups223} is a $(-1)$-shift of the dual of the morphism of Lie algebras  $\bb{T}_{\n{G}^{\dagger}/S^{\ob{gerbe}}}\to \bb{T}_{H_{\beta}/S^{\ob{gerbe}}}$. 
\end{remark}

\begin{remark}\label{RemarkSenOperatorTorsorConceptualConstruction}
Let $D$ be a perfectoid space and $X$ a partially proper smooth rigid space over $D$. Let $G$ be a $p$-adic Lie group and  $\widetilde{X}^{\diamond}\to X^{\diamond}$ a pro-\'etale $\underline{G}$-torsor and let $\beta\colon X^{\HT}\to BG^{\la}$ be its locally analytic decompletion.    \Cref{Construction:SenOperator-Torsor}, applied with $S=X$, $G=\bb{T}_{X/D}\{1\}$, $s$ the morphism from \eqref{eqCartesianDiagramHTStack} (so that $S^{\mathrm{gerbe}}=X^\HT$, cf. \eqref{eqCartesianDiagramHTStack}) and $H=G^\la$ produces from $\beta$ a natural  morphism of groups over $X^{\HT}$
\[
\Omega(\Phi_{\beta})\colon \bb{T}_{X/D}\{1\}^{\dagger} \to G^{\dagger}_{\beta}. 
\]
Passing to tangent spaces, this produces  a  morphism  of  Lie algebras 
\[
\ob{Sen}_{\beta}\colon \bb{T}_{X/D}\{1\} \to \f{g}_{\beta}, 
\] where   $\f{g}_{\beta}$ is the vector bundle on $X^{\HT}$  obtained as pullback along $\beta$ of the adjoint action of $G^{\la}$.  
\end{remark} 

\begin{lemma}\label{LemmaIdentificationSenOperators}
The dual of the morphism $\ob{Sen}_{\beta}\colon \bb{T}_{X/D}\{1\} \to \f{g}_{\beta}$ is the geometric Sen operator of \Cref{DefinitionSenOperatorTorsor}. Furthermore, $\theta^{\ob{geo}}_{\beta}$ naturally identifies with the geometric Sen operator of \cite[Theorem 1.0.4]{camargo2022geometricsen}. 
\end{lemma}
\begin{proof}
The first statement identifying the dual of $\ob{Sen}_{\beta}$ with $\theta^{\ob{geo}}_{\beta}$ is a special case of \Cref{RemarkMapsCotangentComplexesAbstractFramework}. It is left to show that the geometric Sen operator $\theta^{\ob{geo}}_{\beta}$ agrees with that of \cite[Theorem 1.0.4]{camargo2022geometricsen}. Let $\theta_{\beta}'$ denote the geometric Sen operator from \textit{loc. cit.}.  We observe that by \cite[Theorem 1.0.4]{camargo2022geometricsen} the map $\theta_{\beta}'\colon \f{g}_{\beta}^{\vee}\to \Omega^1_{X/D}\{-1\}$ is completely determined by the fact that if $M\in \ob{D}(BG^{\la})$, then the geometric Sen operator of $\beta^*M\in \ob{D}(X^{\HT})\subset \ob{D}(\widehat{X^{\diamond}})$ (with inclusion given by $*$-pullback along $\widehat{X^{\diamond}}\to X^{\HT}$) is given by the map $\theta_{\beta}'$ composed with the derivations of  $\f{g}_{\beta}$ on $M$.  This property is satisfied by $\theta^{\ob{geo}}_{\beta}$  being the dual of the morphism of Lie algebras $\ob{Sen}_{\beta}\colon \bb{T}_{X/D}\{1\}\to \f{g}_{\beta}$ obtained from \Cref{Construction:SenOperator-Torsor}. 
\end{proof}

\begin{remark}\label{RemarkConstructionGeometricSenOperatorsTorsorsMixed}
The analogues of \Cref{RemarkSenOperatorTorsorConceptualConstruction} and \Cref{LemmaIdentificationSenOperators}  apply in the mixed case for a smooth partially proper smooth rigid space $X$ over a finite extension $K$ over $\bb{Q}_p$. Indeed, consider a pro-\'etale $\underline{G}$-torsor $\widetilde{X}^{\diamond}\to X^{\diamond}$ with decompletion $\beta\colon X^{\HT}\to BG^{\la}$. By \Cref{Construction:SenOperator-Torsor}, one has a map of Lie algebras
\[
\ob{Sen}_{\beta}\colon (\ob{FE}_X\{-1\})^{\vee} \to \f{g}_{\beta}.
\]
By the same argument as in \Cref{LemmaIdentificationSenOperators}, the dual of $\ob{Sen}_{\beta}$ is the mixed Sen operator $\theta^{\ob{mix}}_{\beta}$ of \Cref{DefinitionSenOperatorTorsor}. In the case where $X=K$ and $\widetilde{X}^{\diamond}$ arises from a Galois extension of $X$,   the dual of the arithmetic  Sen operator $\theta^{\ob{arith}}_{\beta}$ is precisely the classical Sen operator of  \cite[\S 3.2]{SenContinuousCoho}  and \cite[Remarque 1.4]{BergerColmezII}. 
\end{remark}

A final property that is relevant to compare with  \cite{BergerColmezII}, \cite{MR4390302} and \cite{camargo2022geometricsen} is the description of the natural locally analytic decompletions of $\underline{G}$-torsors in terms of locally analytic vectors at infinite level.  

\begin{proposition}[Arithmetic case]\label{pro:Relation-SenTheory-BergerColmez}
Let $X=\GSpec (K)$ with $K$ a finite extension of $\bb{Q}_p$, and let $\widetilde{X}=\GSpec (K_{\infty})$ with $K_{\infty}$ a an algebraic Galois extension of $K$ with Galois group a $p$-adic Lie group $G$, and such that $K^u_{\infty}$ is perfectoid. Let $\beta\colon X^{\HT}\to BG^{\la}$ be the  $G^{\la}$-torsor  obtained from the  pro-\'etale extension $K_{\infty}^u/K$ by \Cref{CoroEquivalenceTorsor}. Then the map $\beta$ is $0$-affine and $\beta_* 1 = K^{u,\la}_{\infty}\subset  K^{u}_{\infty}$ is the subspace of locally analytic vectors of $K^u_{\infty}$. In particular, we have a natural equivalence of Gelfand stacks
\[
X^{\HT}\xrightarrow{\sim}\GSpec (K^{u,\la}_{\infty})/ G^{\la}.
\]  
\end{proposition}
\begin{proof}
Note that the composite map $f\colon \widehat{X^{\diamond}}\to X^{\HT}\to BG^{\la}$ is induced by the $G^{\cont}$-torsor $\widehat{\widetilde{X}^{\diamond}}\to \widehat{X^{\diamond}}$. In particular, $f_* 1  = \Gamma_{v}(\widetilde{X}^{\diamond},\widehat{\s{O}}_X )^{G-\la} = K_{\infty}^{u,G-\la}$, where as before $(-)^{G-\la}$ refers to the derived locally analytic vectors. Then, to prove the statement we need to show that  $K_{\infty}^{u,G-\la}$ is in degree $0$, and that the map $X^{\HT}\to BG^{\la}$ is $0$-affine. The   maps $\GSpec (K)\to BG^{\la}$ and $\GSpec (K_{\infty}^u)\to \GSpec (K)$ are  epimorphisms of Gelfand stacks with $t$-exact pullback; the first map is clear, for the second map we observe that $K\to K_{\infty}^u$ is descendable as it admits a  section as $K$-module, and flatness follows from \cite[Lemma 3.21]{jacinto2021solid}. Therefore,  the $0$-affineness and the vanishing  of higher locally analytic vectors  can be shown   after base change along $\GSpec (K_{\infty}^u)\to BG^{\la}$ (cf. \cite[Proposition 7.9]{GestaltenScholze}).

We first consider the base change of the map $\beta\colon X^{\HT}\to \GSpec (K)/G^{\la}$ along $\GSpec (K_{\infty}^u)\to \GSpec (K)$. Recall from \Cref{prop:ht-stack-of-perfectoids} that  for a perfectoid $D$ one has $D^{\HT}=D$, thus, we have a natural section $\GSpec(K_\infty^^u)=K^{u,\HT}_{\infty}\to K^{\HT}\to \GSpec (K)$ and this produces a natural splitting of the $\bb{G}_m^{\dagger}$-gerbe $K^{\HT}\times_{\GSpec (K)} \GSpec (K^{u}_{\infty}) \to  \GSpec (K^{u}_{\infty})$, showing that  the base change of $\beta$ to $K^{u}_{\infty}$ is equivalent to the morphism of classifying stacks over $K^u_{\infty}$
\[
\beta_{K^u_{\infty}}\colon \GSpec (K_{\infty}^u)/ \bb{G}_m^{\dagger}\to  \GSpec (K_{\infty}^u)/G^{\la}. 
\]
On the other hand,   the composite map 
\[
\GSpec (K_{\infty}^u) \to \GSpec (K_{\infty}^u)/ \bb{G}_m^{\dagger}\to  \GSpec (K_{\infty}^u)/G^{\la}
\]
gives rise to the trivial  $G^{\la}$-torsor over $K^{u}_{\infty}$ by construction, namely, this follows from the fact that the diamond of $\beta$ corresponds to the $\underline{G}$-torsor $\Spd K^u_{\infty} \to \Spd K$. In other words, $\beta_{K_\infty^u}$ is a pointed map. Therefore, by taking loops, the map  $\beta_{K^u_{\infty}}$  is equivalent to the morphism of groups over  $K^u_{\infty}$
\[
\Omega(\beta_{K^u_{\infty}})\colon \bb{G}_m^{\dagger}\times \GSpec (K^u_{\infty}) \to G^{\la}\times \GSpec (K^u_{\infty}),
\]
and the pullback of $\beta_{K^u_{\infty}}$ along the map $\GSpec (K_{\infty}^u) \to \GSpec (K_{\infty}^u) /G^{\la}$ is the morphism
\[
G^{\la}_{K^u_{\infty}}/\bb{G}_{m,K^{u}_{\infty}}^{\dagger}\to \GSpec (K^{u}_{\infty}). 
\]
In total, we deduce a pullback square
\[
\begin{tikzcd}
G^{\la}_{K^u_{\infty}}/\bb{G}_{m,K^{u}_{\infty}}^{\dagger} \ar[r] \ar[d]  &  \GSpec(K_{\infty}^u) \ar[d] \\ 
X^{\HT}  \ar[r] & BG^{\la}\times X.
\end{tikzcd}
\]
Now, thanks to \Cref{RemarkConstructionGeometricSenOperatorsTorsorsMixed} the morphism of groups $\bb{G}_{m,K^u_{\infty}}^{\dagger}\to G^{\dagger}_{K^u_{\infty}} \to G^{\la}_{K^u_{\infty}}$ arises from the exponential in overconvergent groups of the morphism of Lie algebras $\ob{Sen}\colon K^u_{\infty}\to \f{g}\otimes K^u_{\infty}$ produced by the arithmetic Sen operator as in \cite[\S 6]{BergerColmezII}. Thanks to a theorem of Sen \cite[\S 5]{zbMATH03407789}, and since $K^u_{\infty}$ is perfectoid,  the map $\ob{Sen}$ is non-zero and hence injective. Fixing a splitting of the map $\ob{Sen}$, one deduces that the quotient stack  $G^{\la}_{K^u_{\infty}}/\bb{G}_{m,K^u_{\infty}}^{\dagger}$ is a Berkovich space over $K^u_{\infty}$, whose sheaves of functions given by the dual nuclear Fr\'echet $K^u_{\infty}$-vector space 
\[
\s{O}(G^{\la}_{K^u_{\infty}}/\bb{G}_{m,K^u_{\infty}}^{\dagger})= C^{\la}(G, K^u_{\infty})^{\theta^{\ob{arith}}=0} \subset  C^{\la}(G, K^u_{\infty})
\]
of $K^u_{\infty}$-locally analytic functions of $G$ that are killed by the arithmetic Sen operator acting by derivations.  Indeed, the map $\beta_{K^u_{\infty}}$ factors through a morphism 
\[
 \GSpec(K^u_{\infty})/\bb{G}^{\dagger}_{m,K^u_{\infty}}\xrightarrow{\widetilde{\beta}_{K^u_{\infty}}} \GSpec(K^u_{\infty})/G^{\dagger}\cong G^{\sm}_{K^u_{\infty}}/G^{\la}\to \GSpec(K^u_{\infty})/G^{\la}.
\]
Since the map  $\GSpec(K^u_{\infty})/G^{\dagger}\to \GSpec(K^u_{\infty})/G^{\la}$ is $0$-affine and the $*$-push is $t$-exact as it is an induction,  it suffices to see that $\widetilde{\beta}_{K^u_{\infty}}$ is $0$-affine and that $\widetilde{\beta}_{K^u_{\infty},*}1$ is in degree $0$. This last claim can be proven after base change along $\GSpec(K^u_{\infty})\to \GSpec(K^u_{\infty})/G^{\dagger}$, and hence it suffices to show that $G^{\dagger}_{K^u_{\infty}}/\bb{G}_m^{\dagger}$ is represented by the ring $\s{O}(G^{\dagger}_{K^u_{\infty}})^{\theta^{\ob{arith}}=0}$. Now, if $\f{g}_{K^u_{\infty}}\cong \langle \theta^{\ob{arith}} \rangle \oplus V$ is a splitting of the arithmetic Sen operator, we have a $\theta$-equivariant isomorphism of commutative algebras
\[
\s{O}(G^{\dagger}_{K^u_{\infty}})\cong \s{O}(\bb{G}_{m,K^u_{\infty}}^{\dagger})\otimes_{K^u_{\infty}} \s{O}(V^{\dagger})
\]
which proves that the action of $\bb{G}_m^{\dagger}$ on $G^{\dagger}_{K^u_{\infty}}$ is free, and hence that $G^{\dagger}_{K^u_{\infty}}/\bb{G}_m^{\dagger}$ is the spectrum of $\s{O}(G^{\dagger}_{K^u_{\infty}})^{\theta^{\ob{arith}}=0}\cong \s{O}(V^{\dagger})$.

Finally, the $0$-affineness of the map $\beta\colon X^{\HT}\to BG^{\la}$ and the fact that $\beta_* 1$ is a static Gelfand ring (hence connective) imply that $\beta$ is To\"en-affine and therefore that the natural map 
\[
X^{\HT}\xrightarrow{\sim} \GSpec (K^{u,\la}_{\infty})/G^{\la}
\]
is an equivalence of Gelfand stacks, see \Cref{sec:deriv-berk-spac-different-notions-of-affineness} (4). 
\end{proof}

\begin{remark}
The proof of \Cref{pro:Relation-SenTheory-BergerColmez} establishes the formula
$$
K_{\infty}^{u,\la} \otimes_{K} \ K_{\infty}^{u} \cong C^{\la}(G, K^u_{\infty})^{\theta^{\ob{arith}}=0},
$$
which is exactly \cite[Theorem 6.1]{BergerColmezII}. 
\end{remark}

To state, the geometric and mixed analogues of the last proposition, let us first introduce some notation.  Let $D$ be a perfectoid space and $X$ a partially proper smooth rigid space over $D$. Let $\widetilde{X}^{\diamond}\to X^{\diamond}$ be a  profinite \'etale Galois extension   with Galois group a compact $p$-adic Lie group $G$ and let $\beta\colon X^{\HT}\to B_X G^{\la}$ be the decompletion obtained from \Cref{CoroEquivalenceTorsor}. Let $\widetilde{X}\to X$ be the $G^{\sm}$-torsor obtained from $\widetilde{X}^{\diamond}\to X^{\diamond}$. In particular, we have a natural map $X^{\HT}\to \widetilde{X}/G^{\la}$, whose composition with $\widetilde{X}/G^{\la}\to X/G^\la$ is $\beta$. We denote by $\widetilde{\beta}$ the composition
\[
\widetilde{\beta}\colon X^{\HT}\to \widetilde{X}/G^{\la} \to  |\widetilde{X}|_{\Betti}/G^{\la}. 
\]

\begin{proposition}[Geometric and mixed cases]\label{pro:Relation-SenTheory-Pan}
Keep the assumptions and notations from above. Suppose that the geometric Sen operator $\theta_{\widetilde{X}}\colon \f{g}^{\vee}_{\beta}\to \Omega^1_{X/D}\{-1\}$ is   surjective, then the map $\widetilde{\beta}\colon X^{\HT}\to |\widetilde{X}|_{\Betti}/G^{\la}$ is $0$-affine, and $\widetilde{\beta}_{*} 1 = \s{O}_{\widehat{\widetilde{X}^{\diamond}}}^{\la}\subset \widehat{\s{O}}_{\widehat{X}^{\diamond}}$ is the subsheaf on $|\widetilde{X}|$ of locally  analytic sections of the structural sheaf of $\widehat{\widetilde{X}^{\diamond}}$. Furthermore, there is a strict closed basis   $\{Y\subset |\widetilde{X}|\}$ such that  $\s{O}_{\widehat{\widetilde{X}^{\diamond}}}^{\la}(Y)$ is a  static Gelfand ring. Moreover, the natural morphism of Gelfand stacks
 \begin{equation}\label{eqEquivalenceTorsorInfiniteLecelHTStack}
 X^{\HT}\to \GSpec_{|\widetilde{X}|_{\Betti}}(\s{O}_{\widehat{\widetilde{X}^{\diamond}}}^{\la})/G^{\la}
 \end{equation}
 is an equivalence. 

An analogous statement holds  if $X/K$ is a partially proper smooth rigid space over a finite extension $K$ of $\bb{Q}_p$,  after replacing the geometric Sen operator with the mixed Sen operator of \Cref{RemarkConstructionGeometricSenOperatorsTorsorsMixed}. 
\end{proposition}
\begin{proof}
We will prove the proposition only in the geometric case, the mixed case being analogous.  

By construction, the map $\widetilde{\beta}\colon  X^{\HT}\to |\widetilde{X}|_{\Betti}/G^{\la}$ factors as 
\[
X^{\HT}\to \widetilde{X}/G^{\la} \to  |\widetilde{X}|_{\Betti}/G^{\la}. 
\]
The map  $\widetilde{X}\to  |\widetilde{X}|_{\Betti}$ is $0$-affine by \Cref{sec:deriv-berk-spac-different-notions-of-affineness} (4) and the proof of \Cref{sec:deriv-berk-spac-berkovich-spaces-as-relative-spectrum}. Thus, by \cite[Proposition 6.6 (iii)]{GestaltenScholze}, it suffices to show that    $\widetilde{\beta}'\colon X^{\HT}\to \widetilde{X}/G^{\la}$ is $0$-affine and that $\widetilde{\beta}'_* 1$ is static when restricted to  a suitable strict closed basis of affinoids. 

 The rest of the proof follows essentially the same steps as that of \Cref{pro:Relation-SenTheory-BergerColmez}, with some minor technical modifications. The $0$-affineness can be proven locally on the target by  \cite[Proposition 7.9]{GestaltenScholze}.    Therefore, by localizing in the \'etale topology of $X$, we can assume without loss of generality that $X$ admits a chart $\psi\colon X\to \bb{G}_{m,D}^d$ with coordinates $T_1,\ldots, T_d$, and that $G$ is a compact $p$-adic Lie group. 
 By  taking a cover by affinoid rigid spaces, we can even assume that $X=\GSpec (R)$ is affinoid with $R$ a smoothoid Banach algebra over $D$. Hence, it suffices to show that the map $\beta\colon X^{\HT}\to X/G^{\la}$ is $0$-affine  and $\beta_* 1$ is static.

 We  prove the $0$-affineness of  the map $\beta$ altogether with the fact that $\beta_* 1$ is static.  Let $\widetilde{\bb{G}}_{m,D}^{u} \to \bb{G}_{m,D}$ be the perfectoid multiplicative group over $D$, passing to arc-stacks it is a  $\bb{Z}_p(1)$-pro-\'etale cover of the multiplicative group. Taking $d$-th fold products and pullback along the map $\psi$, we have a $\Gamma=\bb{Z}_p(1)^d$-pro-\'etale torsor $X_{\infty}^{\diamond}\to X$ given by taking perfectoid coordinates $T_1^{\flat},\ldots, T_d^{\flat}$.  Let $\widetilde{X}_{\infty}^{\diamond}:= \widetilde{X}^{\diamond}\times_{X^{\diamond}} X^{\diamond}_{\infty}$, it is a $\underline{G\times \Gamma}$-pro-\'etale torsor over $X^{\diamond}$, an $\underline{G}$-torsor over $X_{\infty}^{\diamond}$ and in particular an  affinoid perfectoid space. In particular, $\widehat{\widetilde{X}_{\infty}^{\diamond}}\to X$ is a $!$-cover, and by \cite[Proposition 7.9]{GestaltenScholze}, the map $\beta$ is $0$-affine if and only if the base change
  \[
\beta|_{\widehat{\widetilde{X}^{\diamond}_{\infty}}} \colon X^{\HT}\times_X \widehat{\widetilde{X}^{\diamond}_{\infty}} \to X/G^{\la} \times_X  \widehat{\widetilde{X}^{\diamond}_{\infty}}. 
 \]
 is $0$-affine. Furthermore, the map $\s{O}(X)\to \s{O}(\widehat{X_{\infty}^{\diamond}})$ admits an $\s{O}(X)$-linear retract by the explicit description of the perfectoid torus, and the map  $\s{O}(\widehat{X_{\infty}^{\diamond}})\to \s{O}(\widehat{\widetilde{X}_{\infty}^{\diamond}})$ admits an  $\s{O}(\widehat{X_{\infty}^{\diamond}})$-linear retract by \Cref{LemmaSectionCountableProFiniteEtalePerfectoid}. Therefore, to show that $\beta_* 1$ is static, by \Cref{LemmaSectionAndStatic} and proper base change 
  it suffices to show that $\beta|_{\widehat{\widetilde{X}^{\diamond}_{\infty}},*} 1$ is static.

  To ease notation let us write $S:= \widehat{\widetilde{X}^{\diamond}_{\infty}} $. Since $S$ is perfectoid, we have $S^{\HT}=S$ by \Cref{prop:ht-stack-of-perfectoids}, in particular, the map $S\to X$ induces a splitting $S=S^{\HT}\to X^{\HT}$ of the Hodge--Tate stack of $X$. Therefore, we have an equivalence
 \[
 X^{\HT}\times_X S \cong S/\bb{T}_{X/D}^{\dagger}\{1\}.
 \] 
 On the other hand,  the $G^{\la}$-torsor $\beta$ is trivial when base changed to $S$. Therefore, the map $\beta|_S$ produces a morphism of classifying stacks
 \[
\beta|_S \colon S/\bb{T}_{X/D}^{\dagger}\{1\} \to S/G^{\la}
 \]
 whose loops is a map of groups over $S$
 \[
 \bb{T}_{X/D}^{\dagger}\{1\}\times_X S\to G^{\la}\times S. 
 \]
The previous yields a Cartesian square 
 \[
 \begin{tikzcd}
  (G^{\la}\times S)/\bb{T}_{X/D}^{\dagger}\{1\}\times_X S \ar[r] \ar[d]&   S \ar[d]\\ 
X^{\HT} \ar[r,"\beta"] & X/G^{\la}.
 \end{tikzcd}
 \]
Now, by assumption the geometric Sen operator $\theta^{\ob{geo}}_{\beta}\colon \f{g}_{\beta}\to \Omega^1_{X/D}\{-1\}$ is surjective.  Fixing a splitting of this map, one deduces that  $(G^{\la}\times S)/\bb{T}_{X/D}^{\dagger}\{1\}\times_X S$ is the Berkovich space associated to the algebra 
\[
C^{\la}(G, \s{O}(S))^{\theta^{\ob{geo}}=0}\subset  C^{\la}(G, \s{O}(S))
\]
of $\s{O}(S)$-valued locally analytic functions of $G$ that are killed by the geometric Sen operator. The $0$-affineness follows and staticness of $\beta|_{\widehat{\widetilde{X}^{\diamond}_{\infty}},*}1$.

 The last claim about the equivalence \eqref{eqEquivalenceTorsorInfiniteLecelHTStack} follows from \Cref{sec:deriv-berk-spac-different-notions-of-affineness} (4). 
\end{proof}

We used the following lemmas. 

\begin{lemma}\label{LemmaSectionCountableProFiniteEtalePerfectoid}
Let $f\colon A\to B$ be a profinite \'etale and surjective on Berkovich spectra morphism of separable Tate-perfectoid rings with pseudo-uniformizer $\varpi$ with $\varpi^{p}|p$.  Then $A\to B$ admits a retract as $A$-module.  
\end{lemma}
\begin{proof}
We can write $B=(\varinjlim_{n\in \bb{N}} A_n)^u$ as the uniform completion of a sequential  colimit along finite \'etale faithfully flat transition maps with $A=A_0$. This produces an almost isomorphism with respect to the almost structure $(A^{\circ},A^{\circ\circ})$ 
\[
B^u=^a (\varinjlim_n A_n^{\circ})^{\wedge}_{\varpi}.
\]
By \cite[Theorem 5.5.9]{kedlaya_liu_relative_p_adic_hodge_theory_foundations}, given $n\in \bb{N}$ the map $A_n^{\circ}\to A_{n+1}^{\circ}$ is almost finite \'etale and faithfully flat, and by \cite[Theorem 4.1.14]{gabber2003almost}, given a pseudo-uniformizer $\pi_n\in A^{\circ}$, there is an $A_{n}^{\circ}$-linear  $\pi_{n}$-section $r_n\colon A^{\circ}_{n+1}\to A^{\circ}_n$ of the inclusion $s_{n}\colon A_{n}^{\circ} \to A^{\circ}_{n+1}$ given by a trace, i.e. a map $r_n$ such that $r_n\circ s_n=\pi_n \ob{id}_{A_n^{\circ}}$.  For $n\geq 1$, let  $\varpi_n\in A^{\circ}$ be a pseudo-uniformizer such that $\varpi_n\equiv \varpi^{1/p^{n+1}}\mod \varpi A^{\circ}$, then set $\pi_0=\varpi$ and $\pi_{n} = \pi_{n-1}\varpi^{-1}_{n}$ for $n\geq 1$. By construction, we have that $\pi_n\in A^{\circ}$ is a pseudo-uniformizer with $\pi_{n}\equiv \varpi^{1-\sum_{k=2}^{n+1} p^{-k}} \mod \varpi A^{\circ}$ for $n\geq 1$. Denote $R_{n}:=\pi_n^{-1}r_n\colon A_{n+1}\to A_n$ and for $m\geq n$ set $R_{n}^m\colon A_{m}\to A_n$ the composite  $R_n\circ R_{n+1}\cdots R_{m-1}$. By construction the maps $R_{n}^m$ are compatible under composition. Furthermore, by the choice of the pseudo-uniformizers $\pi_n$ it satisfies  $R^m_n(A_m^{\circ}) \subset \varpi^{-1} A_n^{\circ}$. Therefore, the sequence of maps $R^{m}_0$ as $m\to \infty$ extends to a natural $A$-linear section $R^{\infty}_0\colon B\to A$ proving what we wanted.
\end{proof}

\begin{lemma}\label{LemmaSectionAndStatic}
Let $A\to B$ be a morphism of Gelfand rings that admits a retract $B\to A$ as $A$-module. Let $M\in \ob{D}(A)$ be such that $B\otimes_A M$ is static, then $M$ is static. 
\end{lemma}
\begin{proof}
This follows  from the fact that $M\to B\otimes_{A} M$ admits a retract as a morphism of $A$-modules. 
\end{proof}

\subsection{$p$-adic twistor cohomology}
\label{sec:p-adic-twistor}

One of the most powerful aspects of the stacky approach to (rational) $p$-adic Hodge theory is that we automatically get a full fledged cohomology theory, that we call \emph{$p$-adic twistor cohomology}, through the six-functor formalism
\[
  X\mapsto \ob{D}(X^{\Div^{1}})
\]
on Gelfand stacks.

In this subsection, we analyse some aspects of this cohomology theory.
We first discuss $!$-ability.

\begin{definition}
  \label{sec:p-adic-twistor-1-prism-shriekable}
  Let $f\colon Y\to X$ be a morphism of  Gelfand stacks.
  We say that $f$ is $\prism$-$!$-able if the morphism $f^{\Div^1}\colon Y^{\Div^1}\to X^{\Div^1}$ is $!$-able with respect to the $6$-functor formalism $\ob{D}(-)$ on Gelfand stacks. 
\end{definition}

$\prism$-$!$-ability can in practice equivalently be asked with $(-)^{\Div^1}$ replaced by $(-)^\prism$ or $(-)^{\widetilde{\Div^1}}$ because $!$-ability for morphisms on Hyodo-Kato stacks has been discussed in \cite{drff_stacks}.

In order to produce sufficient examples of $\prism$-$!$-able maps  we will need a suitable notion of finite presentation for morphisms of Berkovich spaces. Since our aim is just to produce enough examples, and not to develop a general theory of finite presentation in Berkovich geometry, we will content ourselves with a suitable definition that combines the key cases we would like to capture, namely, $\dagger$-rational localizations and analytifications of morphisms of finite presentation of classical schemes.

\begin{definition}\label{DefinitionLocalDaggerPresentation}
\begin{enumerate}

\item A morphism $A\to B$ of connective $\bb{Q}_p$-algebras is said to be \textit{finitely presented} if one of the following equivalent properties hold (cf. \cite[Theorem 7.4.3.18]{lurie_higher_algebra}):\footnote{In \cite{lurie_higher_algebra} what we call ``finitely presented'' is named ``locally finitely presented''.} 
\begin{enumerate}

\item $\pi_0(B)$ is a finitely presented $\pi_0(A)$-algebra in the classical sense, and $\bb{L}_{B/A}$ is a perfect complex. 

\item $B$ sits in the smallest category  of $A$-algebras spanned by $A[T]$ and stable under finite colimits and retracts.  

\end{enumerate}

\item A morphism $Y\to X$ of Berkovich spaces is said to be \textit{locally $\dagger$-finitely presented} if, locally after passing to strict closed covers in both $Y$ and $X$, $f$ is a composite of $\dagger$-rational localizations, and analytifications of finitely presented morphisms of affine schemes. 

\end{enumerate} 
\end{definition}

\begin{proposition}
  \label{sec:p-adic-twistor-2-properties-for-prism-shriekable-maps}
  Let $f\colon Y\to X$ be a morphism of  Gelfand stacks.
  \begin{enumerate}
  \item If $f$ is $\dagger$-formally \'etale and such that the induced map on perfectoidization $f^\diamond: Y^\diamond \to X^\diamond$ is locally on the arc-topology of $X^\diamond$ lqfd (cf.\ \cite[Definition 5.13]{drff_stacks}), $f$ is $\prism$-$!$-able. For example, pro-(Berkovich)\'etale maps are $\prism$-$!$-able.
  
  \item The affine line $\bb{A}^1_{\bb{Q}_p}\to \GSpec (\bb{Q}_p)$ is $\prism$-$!$-able with  $0$-affine diagonal. 
  
\item Let $f\colon Y=\GSpec (B) \to \GSpec (A)=X$ be a morphism of affinoid Gelfand stacks such that $A\to B$ is surjective and  of finite presentation as animated rings.  Then  $f^{\ob{Div}^1}$ is  $0$-affine. In particular,    there is a natural trivialization of the codualizing sheaf $\delta_{f^{\ob{Div}^1}}\cong 1$.

\item Assume that $X$ and $Y$ are Berkovich spaces, and that $f$ is locally of $\dagger$-finite presentation. Then $f$ is $\prism$-$!$-able.
    \end{enumerate}
\end{proposition} 

\begin{proof}
  \begin{enumerate}
  \item By \Cref{sec:defin-first-prop-5-prismatization-for-dagger-formally-etale-maps} and stability of $!$-ability under base change, the statement reduces to $!$-ability of the induced map on de Rham stacks of Fargues--Fontaine curves, which is proved in \cite[Theorem 6.3.1]{drff_stacks}.

\item   By covering $\bb{A}^1_{\bb{Q}_p}= (\bb{A}^1_{\bb{Q}_p}\backslash \{1\})\cup (\bb{A}^1_{\bb{Q}_p}\backslash \{0\})$, it suffices to show that $\bb{G}_{m,\bb{Q}_p}\to \GSpec (\bb{Q}_p)$ is $\prism$-$!$-able. This follows from the explicit presentation of $\bb{G}_m^{\Div^1}$ of \Cref{sec:analyt-prism-gm-description-over-cyclotomic-prism}, the fact that morphisms between Berkovich spaces are $!$-able, and that classifying stacks of locally analytic $p$-adic Lie groups are $!$-able. To see that $\bb{A}^1_{\bb{Q}_p}$ has $0$-affine diagonal after passing to the $p$-adic twistor stacks, since we have a pullback square
\[
\begin{tikzcd}
\A^1_{\bb{Q_p}} \ar[r,"\Delta"] \ar[d] & \A^1_{\bb{Q}_p} \times \A^1_{\bb{Q}_p} \ar[d,"{(x,-y)}"] \\ 
\GSpec (\bb{Q}_p) \ar[r,"0"] & \A^1_{\bb{Q}_p},  
\end{tikzcd}
\]
it suffices to show that the zero section $0\to \bb{A}^1$ induces a $0$-affine map after passing to $p$-adic twistor stacks, this is a particular case of  (3).

\item  Consider the map of Gelfand rings $A\to B$, since it is of finite presentation as classical animated rings, we have that $B(*)\otimes_{A(*)}A\xrightarrow{\sim} B$, and therefore that $B$ is discrete relative to $A$. Since $A\to B$  is surjective and of finite presentation, we have a presentation  of $A$-algebras
\[
B=\varinjlim_{n\geq -1} A_n
\]
such that:
\begin{itemize}

\item[i.]  $A_n$ is an $A$-algebra of finite presentation.
 
\item[ii.]  $A_n\to B$ is  $n+1$-connective for $n\geq -1$.  

\item[iii.]  There is a finite free $A$-module $M_n$ and a map of $A$-modules $M_n[n+1]\to A_n$ such that $A_{n+1} = A_n\otimes_{\Sym_{A} M_{n}[n+1]} A$.  

\end{itemize}

The map $\pi_0(A)\to \pi_0(B)$ is classically of finite presentation and hence we can write  $\pi_0(B)= \pi_0(A)/(a_1,\ldots, a_d)$. Taking lifts of $a_1,\ldots, a_n$ to $A$, we can set $A_{-1}:=A/^{\bb{L}} (a_0,\ldots, a_{d})$. The  rest of the sequence is constructed  as in the proof of \cite[Theorem 7.4.3.18]{lurie_higher_algebra}.

Now, since $B$ is of finite presentation, it is compact as $A$-algebra and hence there is some $n\geq 0$ and a retract $B\to A_n\to B$.  Write $X_m:=\GSpec (A_m)$, to prove part (3) it will suffice to show the following claim:

\begin{claim*} The following  hold true:
\begin{itemize}

\item[(a)] The morphisms $X_m^{\Div^1}\to X_{m-1}^{\Div^1}$ are represented, locally in the $!$-topology of $X_{m-1}$, by morphisms of affinoid Gelfand stacks. 

\item[(b)] Let $Z\to X$ be a $0$-affine morphism of Gelfand stacks, and let $Y\to Z \to Y$ be a retract over  $X$. Then $Y\to X$ is $0$-affine. 

\end{itemize}
\end{claim*}

Indeed, (a) will imply that the maps $X_m^{\Div^1}\to X_{m-1}^{\Div^1}$ are $!$-able, $0$-proper and  $0$-affine, and hence that the composite $X_m$. Since $Y$ is a retract of $X_{n}^{\Div^1}\to X^{\Div^1}$ is also $0$-proper and $0$-affine. Condition (b) implies that, since   $Y^{\Div^1}\to X_{n}^{\Div^1}$ admits a retract, the map $Y^{\Div^1}\to X^{\Div^1}$ is $0$-proper and $0$-affine, and in particular $!$-able as desired.

\begin{proof}[Proof of Claim]
For (a), denote $X_{-2}:=X$. The map $X_{m}\to X_{m-1}$ for $m\geq -1$ arises as a pullback 
\[
\begin{tikzcd}
X_{m-1} \ar[r] \ar[d] & \GSpec (\bb{Q}_p) \ar[d,"0"] \\
X \ar[r] & \bb{A}^{k_{m}}_{\bb{Q}_p}[-(m+1)],
\end{tikzcd}
\]
where   $\bb{A}^{k_{m}}_{\bb{Q}_p}[-(m+1)]$ is the $(m+1)$-loops of the affine space $\bb{A}^{k_{m}}_{\bb{Q}_p}$. We note that if $m\geq 0$ then  $\bb{A}^1_{\bb{Q}_p}[-(m+1)]=\GSpec(\bb{Q}_p[X_{m+1}])$ where $X_{m+1}$ is a free variable in homological degree $m+1$.  Therefore, by base change, it suffices to show that the zero section $\GSpec (\bb{Q}_p)\to \bb{A}_{\bb{Q}_p}[-(m+1)]$ is, locally in the $!$-topology of the target, represented by affinoid Gelfand stacks after passing to $p$-adic twistor stacks. By induction on $m$, it suffices to deal with the case $m= -1$, and by moving the base point from $0$ to $1$, it suffices to prove this  property for the map 
\[
g\colon \GSpec (\bb{Q}_p)  \to \bb{G}_{m,\bb{Q}_p}
\]
given by the unit. Consider the map 
\[
g^{\Div^1}\colon \Div^1\to \bb{G}_{m,\bb{Q}_p}^{\Div^1}.
\]
It suffices to show that $g$ is represented by affinoid Gelfand stacks after taking a base change along  the cover $S:=\FFdisc{\bb{Q}_p}^{\times}\to \Div^1$.  Therefore, it suffices to show that the map 
\[
g^{\widetilde{\Div^1}/S}\colon S\to \bb{G}_{m,\bb{Q}_p}^{\widetilde{\Div^1}/S}
\]
is, locally on the target, represented by affinoid Gelfand stacks.  This follows from  the fact that $g$ is the unit map, and from the explicit presentation of the prismatization of $\bb{G}_m$ of \Cref{sec:analyt-prism-gm-description-over-cyclotomic-prism}. Indeed, this is a consequence of the fact that the pullback of $g^{\widetilde{\Div^1}/S}$ along $\widetilde{\bb{G}}_{m}\times S\to \bb{G}_{m,\bb{Q}_p}^{\widetilde{\Div^1}/S}$ is isomorphic to the group $ \n{G}_{\widetilde{\Div^1}}\times_{\n{X}_{\ob{cyc}}} S\cong q^{\bb{Z}_p^{\la}}$ which is affinoid relative to $\widetilde{\bb{G}}_{m}\times S$. Thus, (a) is proved.

For (b), one simply applies \cite[Proposition 6.6]{GestaltenScholze} which implies that $0$-affine morphisms of Gestalten are stable under retracts. 
\end{proof}

Finally, the trivialization of the codualizing sheaf follows from \Cref{sec:deriv-berk-spac-f-lower-shriek-via-0-affine-stuff}. This proves (b).

  \item Being $\prism$-$!$-able can be checked $\dagger$-formally \'etale locally on $Y$ and $X$ (again by \Cref{sec:defin-first-prop-5-prismatization-for-dagger-formally-etale-maps}, the fact that arc-covers give $!$-covers on de Rham stacks of Fargues--Fontaine curves, and the fact that $!$-ability is $!$-local on source and target). Hence, since $\dagger$-rational localizations are $\dagger$-formally \'etale,  we can assume without loss of generality that $X=\GSpec (A)$ is affinoid, and that $Y$ is the  analytification of an affine scheme of  finite presentation over $X$. Any such analytification factors through 
  \[
  Y\to \bb{A}^{d}_X\to X
  \]
  where $Y\to \bb{A}^{d}_X$ is locally on affinoids of  $\bb{A}^{d}_X$  given by a surjective map $B\to C$ of Gelfand rings of finite presentation as animated rings.   It follows by parts (2) and (3) that $f$ is $\prism$-$!$-able as desired.   \qedhere
    \end{enumerate}
\end{proof}

Next, we establish several cohomological smoothness and properness results. We start with the case $X=\GSpec(\Q_p)$, which is a version of local Tate duality (with $\Q_p$-coefficients).

\begin{theorem}
  \label{sec:div1-analyt-prism-properness-smoothness-for-div-1}
  \begin{enumerate}
  \item The morphism $\Div^{1,\dR}\to \GSpec(\Q_p)$ is cohomologically proper and cohomologically smooth with dualizing complex $\s{O}\langle 1 \rangle$, where $\s{O}\langle 1\rangle$ is the pullback along $\Div^{1,\dR}=\Q_p^\HK \to \Marc(\F_p)^\HK=\GSpec(\Q_p)/\varphi^\Z$ of the shifted line bundle $\s{O}(1)[2]$ (where $(1)$ stands for the Frobenius twist, i.e.\ the line bundle attached to the rank $1$ isocrystal of slope $-1$).
  \item The morphism $\beta \colon\Div^{1}\to \Div^{1,\dR}$ is $0$-proper with $0$-affine diagonal. 
    \item The morphism $\Div^1\to \GSpec(\Q_p)$ is $0$-proper and cohomologically smooth with dualizing complex $\s{O}_{\Div^1}\{1\}[2]$ given by the  $(\varphi,\Gamma)$-module associated with the cyclotomic character shifted in cohomological degree $-2$. 
  \end{enumerate}
\end{theorem}
\begin{proof}
\begin{enumerate}
\item Cohomological smoothness is deduced from the presentation 
\begin{equation}\label{eqPresentationDiv1dR}
\Div^{1,\dR}=( (1+\widetilde{\bb{D}}^{<1,\dR}_{\bb{Q}_p}) \backslash 1)/\bb{Q}_p^{\times,\sm}
\end{equation} cohomological smoothness of classifying stacks of smooth locally profinite sets, and of the perfectoid open disc \cite[Proposition 5.7.1]{drff_stacks}.  Cohomological properness follows from \cite[Proposition 5.1.4]{drff_stacks} as $\Div^{1,\diamond}\to \Spd \bb{Q}_p$ is a proper morphism of $v$-stacks by \cite[Proposition II.1.21]{fargues2021geometrization}. 

From the presentation of \eqref{eqPresentationDiv1dR}, we deduce from \cite[Example 5.3.21]{heyer20246functorformalismssmoothrepresentations} (for $B\bb{Q}_p^{\times,\sm}$) and \cite[Proposition 5.7.1]{drff_stacks} (for the perfectoid open disc) that  the dualizing complex  of $\Div^{1,\dR}$ is invertible in cohomological degree $-2$.  By \cite[Theorem 7.1.1]{drff_stacks}, $\omega_{\Div^{1,\dR}}[-2]$ a rank $1$ $(\varphi,N,\mathrm{Gal}_{\Q_p})$-module over $\Q_p^{\rm un}$. The calculation of the cohomology of rank $1$ $(\varphi,N,\mathrm{Gal}_{\Q_p})$-modules over $\Q_p^{\rm un}$ shows that the only possibility for the dualizing complex is $\s{O}\langle 1 \rangle$. For a more conceptual computation of the dualizing sheaf of $\Div^{1,\dR}$ see \Cref{alternative-computation-dualizing-Div1dR}. 

\item By \Cref{prop:pullback-description-twistor-stack}  and the preservation of  $0$-properness and $0$-affineness by base change, it suffices to show the analogous properties for the map $f\colon \A^1/\Gm\to \A^{1,\dR}/\Gm^\dR$. Since $0$-proper maps  are stable under composition, it suffices to show that each of the maps 
\[
\A^1/\Gm\to \bb{A}^{1,\dR}/\bb{G}_m \to \A^{1,\dR}/\Gm^\dR
\]
are $0$-proper. The first map is $0$-proper as $\A^1\to \A^{1,\dR}$ is so, being represented in affinoid Gelfand stacks. The second map is $0$-proper since the fibre of $B\bb{G}_m\to B\bb{G}_m^{\dR}$ is $B\bb{G}_m^{\dagger}$, and $B\bb{G}_{m}^{\dagger}$ is prim with diagonal represented in affinoid Gelfand stacks.  It is left to see that the diagonal of $f$ is $0$-affine, we have a commutative diagram with pullback square
\[
\begin{tikzcd}
\A^1 /\Gm \ar[r] & \A^1/\Gm\times_{\bb{A}^{1,\dR}/\bb{G}_m}  \A^1/\Gm  \ar[r] \ar[d] &  \A^1/\Gm \times_{\A^{1,\dR}/\Gm^\dR} \A^1/\Gm \ar[d]  \\  
 &  \bb{A}^{1,\dR}/\bb{G}_m \ar[r] &   \bb{A}^{1,\dR}/\bb{G}_m\times_{ \A^{1,\dR}/\Gm^\dR}  \bb{A}^{1,\dR}/\bb{G}_m.
\end{tikzcd}
\]
The map  $\A^1 /\Gm \to \A^1/\Gm\times_{\bb{A}^{1,\dR}/\bb{G}_m}  \A^1/\Gm $ is $0$-affine as $\A^1\to \A^{1,\dR}$ is so. As $0$-affine maps are stable under composition and base change  (\cite[Proposition 6.6 (ii) and (iii)]{GestaltenScholze}), it is left to show that the diagonal of $B\bb{G}_m\to B\bb{G}_m^{\dR}$ is $0$-affine, this follows from the following equivalent presentation of the diagonal 
\[
B\bb{G}_m^{\dR} = \bb{G}_m\backslash \bb{G}_m/\bb{G}_m \to   \bb{G}_m\backslash \bb{G}_m^{\dR} / \bb{G}_m
\]
and the fact that $\bb{G}_m\to \bb{G}_m^{\dR}$ is $0$-affine. 

\item Since $0$-properness can be checked $0$-proper-locally on the source, the $0$-properness of $\Div^1$ follows from (1) and (2) and surjectivity of $\beta$. Cohomological smoothness is proven in \cite{mikami2025finiteness}, equivalently, one can use the cohomologically smoothness of $\FF_{\Q_p^\diamond}$ of \cite[Theorem 4.5.1]{ALBM-FF} and the invertibility of $\alpha_* 1$ in $\ob{D}(\Div^1)$ for the $0$-proper map $\alpha\colon \FF_{\Marc(\Q_p)}\to \Div^1$ (so that $\alpha_*1$ is suave, and in fact invertible). For the identification of the dualizing complex of $\Div^1$, by \cite[Corollary 4.5.2]{ALBM-FF} the dualizing complex $\omega$ of $\FF_{\bb{Q}_p^{\diamond}}$ is the perfect $(\varphi,\Gamma)$-module associated to the cyclotomic character  sitting in cohomological degree $-2$. The map  $\alpha\colon \FF_{\bb{Q}_p^{\diamond}}\to \Div^1$ is $0$-proper  (and therefore has trivial codualizing  sheaf) by \Cref{sec:arithmetic-case:-p-adic-fields} (1), hence \cite[Lemma 4.4.9 (ii)]{heyer20246functorformalismssmoothrepresentations}
 yields 
 \[
 \omega_{\Div^1} \cong \alpha_* \omega
 \]
proving what we wanted. For a more conceptual computation of the dualizing sheaf of $\Div^1$ see \Cref{alternative-computation-dualizing}.  \qedhere
 \end{enumerate} 
 \end{proof}

\begin{remark}\label{alternative-computation-dualizing-Div1dR}
The dualizing sheaf of $\Div^{1,\dR}$ can be computed purely in terms of representation theory. We have an open immersion $\Div^{1,\dR}\subset \BC(\s{O}(1))^{\dR}/\bb{Q}_p^{\times,\sm}$, thus it suffices to compute the dualizing sheaf of $\BC(\s{O}(1))^{\dR}/\bb{Q}_p^{\times,\sm}$.   Consider the composite 
\[
\BC(\s{O}(1))^{\dR}/\bb{Q}_p^{\times,\sm} \xrightarrow{a} \GSpec (\bb{Q}_p)/\bb{Q}_p^{\times,\sm} \xrightarrow{b} \GSpec \bb{Q}_p.
\]
We have $\omega_{\BC(\s{O}(1))^{\dR}/\bb{Q}_p^{\times,\sm}}  = a^*b^!1\otimes a^!1$, by \cite[Lemma 4.4.9(i)]{heyer20246functorformalismssmoothrepresentations} (using that $a$ and $b$ are suave). By \cite[Example 5.3.21]{heyer20246functorformalismssmoothrepresentations}, 
 the dualizing sheaf of classifying stacks of locally profinite groups are modulus characters. Since $\bb{Q}_p^{\times}$ is abelian, we have that $b^!1\cong 1$ and thus $\omega_{\BC(\s{O}(1))^{\dR}/\bb{Q}_p^{\times,\sm}}\cong a^! 1$. 

 Since  $\BC(\s{O}(1))^{\dR}/\bb{Q}_p^{\times,\sm}$ has a natural structure of a group over $\GSpec (\bb{Q}_p)/\bb{Q}_p^{\times,\sm}$, its dualizing sheaf is  $\BC(\s{O}(1))^{\dR}$-equivariant and equal to the pullback along  $\BC(\s{O}(1))^{\dR}/\bb{Q}_p^{\times,\sm}\to \GSpec (\bb{Q}_p)/\bb{Q}_p^{\times,\sm}$ of the dualizing sheaf of the map $\GSpec (\bb{Q}_p)/\bb{Q}_p^{\times,\sm}\to  \GSpec (\bb{Q}_p)/  \BC(\s{O}(1))^{\dR} \rtimes \bb{Q}_p^{\times,\sm}$. 

Then, consider the diagram
\[
\GSpec (\bb{Q}_p) / \bb{Q}_p^{\times,\sm}\xrightarrow{\alpha} \GSpec (\bb{Q}_p) / ( \bb{Q}_p^{\sm}\rtimes  \bb{Q}_p^{\times,\sm}  ) \xrightarrow{\beta} \GSpec (\bb{Q}_p)/   \BC(\s{O}(1))^{\dR}\rtimes  \bb{Q}_p^{\times,\sm}.
\]
The map $\beta$ is fibered in the morphism 
\[
\bb{A}^{1,\dR}_{\bb{Q}_p}/\bb{Q}_p^{\times,\sm}\to \GSpec (\bb{Q}_p)/\bb{Q}_p^{\times,\sm},
\]
this shows that it is represented in de Rham stacks of smooth rigid spaces and by Poincar\'e duality \cite[Theorem 1.3.2 (2)]{drff_stacks} one has $\beta^!1=1[2]$. 

By \cite[Example 5.3.21]{heyer20246functorformalismssmoothrepresentations},  again, we have that $\alpha^! 1\cong \alpha^\ast \delta^{-1}$ where $\delta$ is the dualizing sheaf of the mirabolic  $\bb{Q}_p\rtimes \bb{Q}_p^{\times} $.  We deduce that 
\[
(\alpha\circ \beta)^! 1 = \alpha^! \beta^! 1 =   (\alpha^! 1)[2] \cong  \alpha^\ast \delta^{-1}[2]. 
\] 
Finally, to finish the computation of $\omega_{\BC(\s{O}(1))^{\dR}/\bb{Q}_p^{\times,\sm}}$, one observes that  $\alpha^\ast \delta$ is the character sending $a\in \bb{Q}_p^{\times}$ to $|a|_p^{-1}= p^{v_p(a)}$, and thus that $\alpha^\ast \delta^{-1}$  is the character $\delta^{-1}(a)=p^{-v_p(a)}$, this is precisely the character associated to the isocrystal of slope $-1$.  
\end{remark}

\begin{remark}\label{alternative-computation-dualizing}
The dualizing sheaf of $\Div^1$ can be computed in terms of the cartesian square \Cref{prop:pullback-description-twistor-stack}.  In this computation, we will write $\omega_g$ for the dualizing complex of a suave morphism $g$ and we will use several times the fact that if $g$ is a morphism represented in smooth rigid spaces of pure dimension $d$, then one has an isomorphism $\omega_g\cong \mathbb{L}_g[d]$ (see \cite[\S 7.7]{solidnotes}).

From \Cref{prop:pullback-description-twistor-stack} we have a  fiber square
\[\begin{tikzcd}
	{\Div^1} & {\Div^{1,\dR}} \\
	{\A^1/\Gm} & {\A^{1,\dR}/\Gm^\dR}.
	\arrow["h'", from=1-1, to=1-2]
	\arrow["{f}"', from=1-1, to=2-1]
	\arrow["{f^\dR}", from=1-2, to=2-2]
	\arrow["h", from=2-1, to=2-2]
      \end{tikzcd}\]
Since $\Div^1 \to \GSpec(\Q_p)$ is suave and since the diagonal of $\A^1/\Gm$ is suave (using that a map between smooth schemes is suave, being l.c.i.), we deduce from \cite[Lemma 4.5.9(ii)]{heyer20246functorformalismssmoothrepresentations} that $f$ is suave. In particular, the map $f^{\dR}$ is also suave.\footnote{This suggests that $\Div^{1,\diamond}\to \bb{A}^{\diamond}/\bb{G}_m^{\diamond}$ is cohomologically smooth for the six-functor formalism of Berkovich motives, but a proof of this fact seems to be non-obvious.}  We deduce from \cite[Lemma 4.5.17(i)]{heyer20246functorformalismssmoothrepresentations} that the dualizing sheaf $\omega_{\Div^1}$ is isomorphic to
\begin{equation}\label{eqComputationDualizingSheafDiv1Square}
\omega_f \otimes f^\ast \omega_{\A^1/\Gm}
\end{equation}
We claim that the dualizing sheaf $\omega_{\A^1/\Gm}$ is $\s{O}(-1)$. Indeed, the map $\A^1/\Gm\to \GSpec(\Q_p)$ factors as the composition of the suave maps $\A^1/\Gm\to B\Gm$ and $B\Gm\to \GSpec(\Q_p)$. Using \cite[Lemma 4.5.17(i)]{heyer20246functorformalismssmoothrepresentations}  again for this composition, we see that we only have to compute the dualizing sheaves for the two maps. For $B\Gm\to \GSpec(\Q_p)$, using Verdier diagonal's trick, the dualizing sheaf is isomorphic to $\Delta^\ast \omega_p$, for $p: B\Gm \times B\Gm \to B\Gm$ the first projection and $\Delta$ the diagonal. By \cite[Lemma 4.5.17(i)]{heyer20246functorformalismssmoothrepresentations} again, $\Delta^\ast \omega_p \otimes \omega_{\Delta}= \omega_{p\circ \Delta}=\s{O}$. The morphism $\Delta$ is equivalent to 
\[
\bb{G}_m\backslash  \bb{G}_m / \bb{G}_m\to \bb{G}_m \backslash  (\GSpec \bb{Q}_p)/\bb{G}_m
\] 
and in particular represented in smooth rigid spaces. Hence, $\omega_{\Delta}= \mathbb{L}_{\Delta}[1] = e^* \bb{L}_{\bb{G}_m}[1]$ where $e\colon 1\to \bb{G}_m$ is the unit section, and $\omega_{\Delta}$ is given by the dual of the adjoint   representation of $\bb{G}_m$ sitting in cohomological degree $-1$, that is,  $\omega_{\Delta}=\s{O}[1]$. We deduce that $\omega_{B\bb{G}_m} = \s{O}[-1]$.

 The morphism $g:\A^1/\Gm\to B\Gm$ is represented in smooth rigid spaces of dimension $1$ and so $\omega_g=\mathbb{L}_g[1]$, the line bundle  $\mathbb{L}_g$ is nothing but the cotangent bundle of $\bb{A}^1$ endowed with the natural action of $\bb{G}_m$, namely, $\omega_{g}=\s{O}(-1)[1]$. Putting everything together, we deduce that $\omega_{\A^1/\Gm}=\s{O}(-1)$, proving our claim. 
 
 From \eqref{eqComputationDualizingSheafDiv1Square} one gets
$$
\omega_{\Div^1} =  f^! 1 \otimes f^\ast \s{O}(-1).
$$
By \cite[Lemma 4.5.13(i)]{heyer20246functorformalismssmoothrepresentations} (suave base change),
$$
f^!1= h^{\prime,\ast} f^{\dR,!} 1.
$$
The dualizing sheaf $\omega_{\A^{1,\dR}/\Gm^\dR}$ is the unit,\footnote{Indeed, this follows from Poincar\'e duality for de Rham stacks on smooth rigid spaces, and its obvious extension to Artin stacks.} hence $f^{\dR,!}1= \omega_{\Q_p^\HK} = \s{O}_{\Q_p^\HK} \langle 1 \rangle[2]$. Thus we get
$$
\omega_{\Div^1} = h^{\prime, \ast} \s{O}_{\Q_p^\HK}\langle 1\rangle[2]  \otimes f^\ast \s{O}(-1),
$$
which is (by definition) the (shifted by $2$) Breuil--Kisin twist on $\Div^1$.

We note that the reasoning applies to $\Q_p^\prism$ as well, showing that its dualizing complex is given by the (shifted by $2$) Breuil-Kisin twist.
\end{remark}

\begin{remark}
From \Cref{sec:div1-analyt-prism-properness-smoothness-for-div-1}, one deduces that $\Q_p^\prism$ is cohomologically smooth, with dualizing complex given by $\s{O}_{\Q_p^\prism}\{1\}[2]$. Indeed, since $(0,\infty)$ is cohomologically smooth, it suffices to show that the map $\Q_p^\prism \to (0,\infty)$ is cohomologically smooth and to describe the dualizing complex for it. Both claims are local on target, and thus by \Cref{cor:prismatization-push-out}, the question reduces to the analogous question for the Hyodo--Kato stack and for the $p$-adic twistor stack. Now the result follows from the fact the rational shtuka attached to $\s{O}_{\Q_p^\prism}\{1\}$ is $\s{O}_{\Div^1}\{1\} = \mathcal{I}_{\Div^1}(1) \to \s{O}_{\Div^1}\langle 1 \rangle$ (cf.\ \Cref{sec:breuil-kisin-twist-1-modification-associated-to-bk-twist} and \Cref{rk:line-bundles-on-qp-prism-adn-div1} (4)) together with the description of the dualizing sheaf of $\Div^1$ and $\Div^{1,\dR}$ from \Cref{sec:div1-analyt-prism-properness-smoothness-for-div-1}. Analogously, one can just use the computation of \Cref{alternative-computation-dualizing}.
\end{remark}

From the main results of \cite{ALBM-FF} we can deduce Poincar\'e duality for the analytic prismatization.
For any Gelfand stack $Z$ over $\Div^1$ and any $d\in \Z$, we denote in the following theorem the pullback of the line bundle $\mathscr{O}_{\Div^1}\{d\}$ to $Z$ by $\mathscr{O}_Z\{d\}$. 

\begin{theorem}
  \label{CorollaryPoincar\'e duality} 
  Let $f\colon Y\to X$ be a rigid smooth morphism of Berkovich spaces of relative dimension $d$.
  \begin{enumerate}
  \item The morphism $f^{\Div^1}\colon Y^{\Div^1}\to X^{\Div^1}$ is cohomologically smooth.
  \item Let $f^{\Div^1/\FF_{X^\diamond}}\colon Y^{\Div^1/\FF_{X^\diamond}}=Y^{\Div^1}\times_{X^{\Div^1}}\FF_{X^\diamond}\to \FF_{X^\diamond}$ be the base change of $f^{\Div^1}$ along the natural morphism $h_X\colon \FF_{X^\diamond}\to X^{\Div^1}$.
    Then there is a natural isomorphism $f^{\Div^1/\FF_{X^\diamond},!}\s{O}_{\FF_{X^\diamond}}\cong \s{O}_{Y^{\Div^1/\FF_{X^\diamond}}}\{d\}[2d]$.
  \item Assume that $X$ is a partially proper smooth rigid space over $\Q_p$. 
    Then there is a natural isomorphism $f^{\Div^1,!}\mathscr{O}_{Y^{\Div^1}}\cong \mathscr{O}_{Y^{\Div^1}}\{d\}[2d]$.
  \end{enumerate}
\end{theorem}
\begin{proof}
  By \Cref{sec:defin-first-prop-5-prismatization-for-dagger-formally-etale-maps} we can prove cohomological smoothness of $f^{\Div^1}$  by taking a  rigid \'etale cover on the source. Hence,   we can assume that $f$ is the composition of a rigid \'etale map and the projection $\mathbb{G}_{m,X}^d\to X$, and then in fact that $Y=\mathbb{G}_{m,X}^d$. 
  By stability of cohomologically smooth maps by base change and compositions,    we can assume $X=\GSpec(\Q_p)$ and $Y=\Gm$. 
  In this case, the morphism $h_X\colon \FF_{X^\diamond} \to X^{\Div^1}$ is an epimorphism (cf. \Cref{sec:geom-case:-smooth-relative-case}, \Cref{sec:mixed-case:-smooth-mixed-case1}).
  Hence it suffices to show the cohomological smoothness of the morphism
  \[
    g:=f^{\Div^1/\FF_{X^\diamond}}\colon Y^{\Div^1/\FF_{X^\diamond}}\to \FF_{X^\diamond}.
  \]
  By \cite[Theorem 4.4.4]{ALBM-FF} we know that the morphism $f^{\FF}\colon \FF_{Y^{\diamond}}\to \FF_{X^{\diamond}}$ induced by $f$ is cohomologically smooth.\footnote{Strictly speaking, we do not work here with the $6$-functor formalism from \cite{ALBM-FF}, but with its light variant from \cite[\S 4.4]{drff_stacks}; as explained there, one can however compare the $!$-functors.}
  On the other hand, by \Cref{sec:geom-case:-smooth-relative-case}, \Cref{sec:mixed-case:-smooth-mixed-case1} (and descent to $\FF_{X^\diamond}$ from the case that $X$ is perfectoid), the map $a \colon \FF_{Y^{\diamond}}\to Y^{\Div^1/\FF_{X^\diamond}}$ is prim and satisfies $a_*\s{O}_{\FF_{Y^\diamond}}=\s{O}_{Y^{\Div^1/\FF_X}}$.
  This implies that the unit of $Y^{\Div^1/\FF_X}$ is suave over $\FF_{X^{\diamond}}$. 
  It remains to compute the dualizing sheaf of $f^{\Div^1/\FF_{X^\diamond}}$ resp.\ $f^{\Div^1}$.
  This will imply then cohomological smoothness in (1) instead of mere suaveness, i.e., invertibility of the dualizing complex, as we can again reduce to the case $Y=\mathbb{G}_{m,\Q_p}$, where we can apply (3).
  
  The idea of the proof below is, using the assumption on $X$, to reduce to the analogous statement for relative Fargues--Fontaine curves, where the computation was already performed in \cite{ALBM-FF}.
  We start with the case for $f^{\Div^1/\FF_{X^\diamond}}$, i.e., we prove part (2).
  For this, we apply \cite[Lemma 4.4.9.(ii)]{heyer20246functorformalismssmoothrepresentations} to the diagram
\[\begin{tikzcd}
	{\FF_{Y^\diamond}} & \\
	{Y^{\Div^1/\FF_{X^\diamond}}} & {\FF_{X^\diamond}}
	\arrow["a"', from=1-1, to=2-1]
	\arrow["{f^{\FF}}", from=1-1, to=2-2]
	\arrow["g", from=2-1, to=2-2]
      \end{tikzcd}\]
    with $f^{\FF}, g$ suave, $a$ $0$-proper (and thus with trivial codualizing sheaf), 
    and the objects $P=Q=\mathscr{O}_{\FF_{Y^\diamond}}$.
    This implies that
    \[
      g^!\mathscr{O}_{\FF_{X^\diamond}}\cong a_\ast(f^{\FF,!}\mathscr{O}_{\FF_{X^\diamond}})\cong a_\ast(\mathscr{O}_{\FF_{Y^\diamond}}(d)[2d])\cong \mathscr{O}_{Y^{\Div^1/\FF_{X^\diamond}}}\{d\}[2d]
    \]
    using \cite[Theorem 4.4.6]{ALBM-FF}, $a^\ast(\mathscr{O}_{Y^{\Div^1/\FF_{X^\diamond}}}\{d\}[2d])\cong \mathscr{O}_{\FF_{Y^\diamond}}(d)[2d]$, the projection formula for $a_\ast$ and $a_\ast(\mathscr{O}_{\FF_{Y^\diamond}})=\mathscr{O}_{Y^{\Div^1/\FF_{X^\diamond}}}$.

    We now prove (3), and thus assume that $X$ is a partially proper smooth rigid variety over $\Q_p$.
    By \Cref{sec:mixed-case:-smooth-mixed-case1} the morphism $h_X\colon \FF_{X^\diamond}\to X^{\Div^1}$ is $0$-proper, and therefore  with trivial codualizing sheaf,  such that $h_{X,\ast}(\mathscr{O}_{\FF_{X^\diamond}})\cong \mathcal{O}_{X^{\Div^1}}$.
    
  By suave base-change (\cite[Lemma 4.5.13]{heyer20246functorformalismssmoothrepresentations}) applied to the cartesian square (and the established suaveness of $f^{\Div^1}$)
  \[
    \begin{tikzcd}
	{Y^{\Div^1/\FF_{X^\diamond}}} & {\FF_X^{\diamond}} \\
	{Y^{\Div^1}} & {X^{\Div^1}},
	\arrow["g", from=1-1, to=1-2]
	\arrow["h_X'", from=1-1, to=2-1]
	\arrow["h_X", from=1-2, to=2-2]
	\arrow["f^{\Div^1}", from=2-1, to=2-2]
      \end{tikzcd}\]
    one has using (2) that 
    \begin{equation}
      \label{eq:1}
      h_X^{\prime \ast} f^{\Div^1,!} \s{O}_{X^{\Div^1}} \cong g^! \s{O}_{\FF_{X^\diamond}}\cong \mathscr{O}_{Y^{\Div^1/\FF_{X^\diamond}}}\{d\}[2d]\cong h_X^{\prime,\ast}\mathscr{O}_{Y^{\Div^1}}\{d\}[2d].
    \end{equation}
As $h_{X}^\prime$ is a base change of $h_X$, we deduce that $h_X'$ is $0$-proper  and satisfies $h_{X,\ast}^\prime \s{O}_{Y^{\Div^1/\FF_{X^\diamond}}} = \s{O}_{Y^{\Div^1}}$.
Applying $h_{X,*}^\prime$ and the projection formula to \eqref{eq:1}, we deduce that
\[
  f^{\Div^1,!}\mathscr{O}_{X^{\Div^1}}\cong \mathscr{O}_{Y^{\Div^1}}\{d\}[2d]
\]
as desired.
 \end{proof}

 \begin{remark}
   \label{sec:p-adic-twistor-1-remark-general-case}
   The assumption in \Cref{CorollaryPoincar\'e duality} (3) that $X$ is a partially proper smooth rigid space  over $\bb{Q}_p$ is there only for simplicity and should not be necessary.
   To get the result in general, one should run a deformation to the normal bundle argument to identify the dualizing sheaf.  Similarly, an analogue technique involving the cartesian square of the cone as in \Cref{alternative-computation-dualizing} should  give rise to a different proof of \Cref{CorollaryPoincar\'e duality} (3) that is independent of the six-functor formalism of Fargues--Fontaine stacks.   We do not pursue this here. 
\end{remark}

\begin{proposition}\label{berkovich-smooth-qc-weakly-cohomologically-proper}
Let $f:Y\to X$ be a morphism of Berkovich spaces   locally of $\dagger$-finite presentation (in the sense of \Cref{DefinitionLocalDaggerPresentation})  whose underlying morphism of arc-stacks is proper. Then $f^{\prism}: Y^{\prism} \to X^{\prism}$ is  $0$-proper (hence $!$-able).
 In particular, it has trivial codualizing sheaf.
\end{proposition}
\begin{proof}
  Let us first observe that the statement to be proved is $!$-local on the target.
  Localizing over $(0,\infty)$, we thus see that it is enough to prove the result separately for the Hyodo--Kato stack and for the $p$-adic twistor stack, using \Cref{cor:prismatization-push-out}. For the Hyodo--Kato stack the statement follows from \cite[Proposition 5.1.4]{drff_stacks}), noticing that the locally $\dagger$-finite presentation assumption on $f$ implies that $f^{\diamond}\colon Y^{\diamond}\to X^{\diamond}$ is lqfd.   We give the argument for the $p$-adic twistor stack.

$0$-properness can be proved locally on the target (by \cite[Lemma 4.4.8]{heyer20246functorformalismssmoothrepresentations}), we can then suppose that $X$ is affinoid. Since $Y^{\diamond}\to X^{\diamond}$ is proper, there is a finite closed analytic cover by affinoids $Y_i\to Y$ for which the base change $Y_i\times_{Y} Y_j \to Y_{i}$ is a closed immersion of analytic stacks, in particular, the map $\bigsqcup_{i\in I} Y_i\to Y$ is  a cohomologically proper  and $0$-affine $!$-cover (as closed maps of analytic stacks are $0$-affine), and the same holds for $\bigsqcup_i Y_i^{\Div^1}\to Y^{\Div^1}$. By \cite[Remark 4.4.11(iii)]{heyer20246functorformalismssmoothrepresentations},  $0$-proper maps are local along cohomologically $0$-proper $!$-covers on the source.\footnote{The reference is only stated for cohomologically proper maps but the exact same argument holds for $0$-proper maps.}   
\end{proof}

  \begin{corollary}\label{cor:poincare-duality}
    Let $f\colon Y \to X$ be a smooth proper morphism of partially proper smooth rigid spaces over $\bb{Q}_p$, of relative dimension $d$.
    Then pushforward along $f^{\Div^1}: Y^{\Div^1}\to X^{\Div^1}$ preserves perfect complexes. If $M$ is a perfect complex on $Y^{\Div^1}$, we have
  $$
  f^{\Div^1}_\ast M^\vee\{d\}[2d] \cong (f_\ast^{\Div^1} M)^\vee,
  $$
  where on the left $(-)^\vee$ refers to linear duality on $Y^{\Div^1}$ and on the right on $X^{\Div^1}$.
\end{corollary}
\begin{proof}
Since rigid smooth maps are Berkovich-smooth, the first assertion is a direct consequence of  \Cref{CorollaryPoincar\'e duality} and \Cref{berkovich-smooth-qc-weakly-cohomologically-proper}, together with the general fact that in a $6$-functor formalism, pushforward along a cohomologically smooth and prim morphism preserves dualizable objects, i.e.\ (on Gelfand stacks) perfect complexes. The second assertion is also an immediate corollary of \Cref{CorollaryPoincar\'e duality} and \Cref{berkovich-smooth-qc-weakly-cohomologically-proper} by applying \cite[Lemma 4.4.9 (ii)]{heyer20246functorformalismssmoothrepresentations}.
\end{proof}

\begin{remark}
When $X=\GSpec(\Q_p)$, the above corollary can be combined with \Cref{sec:div1-analyt-prism-properness-smoothness-for-div-1} to get the following Poincar\'e duality statement. Let $Y$ be a  smooth proper rigid space  over $\Q_p$ and let $g: Y^{\Div^1}\to \GSpec(\Q_p)$. Then for any perfect complex $M$ on $Y^{\Div^1}$, 
$$
g_\ast M^\vee\{(d+1)\}[2(d+1)] \cong (g_\ast M )^\vee,
$$
 where on the left $(-)^\vee$ refers to linear duality on $Y^{\Div^1}$ and on the right on $\GSpec(\Q_p)$ (i.e.\ in solid $\Q_p$-vector spaces). If we only assume that $Y$ is smooth rigid (not necessarily quasi-compact), a similar result holds using $g_!$ on the right hand side. 

For example, this can be applied when $M$ is a finite rank pro-\'etale $\Q_p$-local system $\mathbb{L}$ on a smooth partially proper rigid space $Y$, in which case $g_\ast M^\vee$, resp.\ $g_! M$ identifies with the pro-\'etale cohomology of $\mathbb{L}^\vee$, resp.\ the compactly supported pro-\'etale cohomology of $\mathbb{L}$, and moreover the Breuil--Kisin twist on $Y^{\Div^1}$ is the descent of the Tate twist (cf.\ \Cref{rk:line-bundles-on-qp-prism-adn-div1} (4)). Such results have been proved already in the literature, cf.\ \cite{colmez2023arithmetic} (in dimension $1$, for Tate twists), \cite{zhenghui_li_2024} (for Tate twists), \cite{ALBM-FF}.
\end{remark}

\subsection{The first Chern class for rational prismatic cohomology}
\label{sec:prism-chern-class}

In this subsection, we introduce the first Chern class.
For its construction we follow the general discussion in \cite[Remark 5.1.10(3)]{drff_stacks} with base $\mathcal{S}=\Q_p^\prism$, ring stack $\mathcal{R}=\A^{1,\prism}$, and ``Tate twist'' $\mathcal{L}$ given by the first Breuil--Kisin twist $\s{O}_{\Q_p^\prism}\{1\}$ on $\Q_p^\prism$.
In this general setup, a (weak) theory of Chern classes can be constructed by providing a morphism (of animated abelian group objects)
\[
  c_1^{\mathcal{R}}\colon (B\Gm)^{\mathcal{R}}\to B^2 \mathbb{V}(\mathcal{L}). 
\]
Given a line bundle $\mathcal{M}$ on some $X$, classified by a morphism $f_{\mathcal{M}}\colon X\to B\Gm$, its first Chern class $c_1^{\mathcal{R}}(\mathcal{M})\in H^2(X^{\mathcal{R}},\mathcal{L})=\pi_0(\mathrm{Map}(X^{\mathcal{R}},B^2 \mathbb{V}(\mathcal{L})))$
is then defined as the composition of $f_{\mathcal{M}}^{\mathcal{R}}$ with $c_1^{\mathcal{R}}$.
Here, $X^{\mathcal{R}}$ is the transmutation of $X$ with respect to $\mathcal{R}$.

We note that we do not claim to construct the analytic prismatization as an example of a general theory of transmutation for a ring stack over some base (as this necessitates a discussion of ``Gelfand ring stacks'', that we avoid), but the same construction carried out with our examples of interest (analytic prismatization, Hyodo-Kato stacks, $p$-adic twistor stacks) works.

\begin{construction}
  \label{sec:chern-classes-p-construction-weak-theory-of-first-chern-classes}
  We will now construct the first Chern class as a morphism
  \[
   c_1^\prism\colon B\Gm^\prism\to B^2\mathbb{V}(\s{O}_{\Q_p^\prism}\{1\})
  \]
  of animated abelian group Gelfand stacks over $\Q_p^\prism$, where $\s{O}_{\Q_p^\prism}\{1\}$ is the Breuil-Kisin twist from \Cref{sec:breuil-kisin-twist}.
  Delooping, it suffices to construct
  \[
    c_1^\prism\colon \Gm^\prism\to B\mathbb{V}(\s{O}_{\Q_p^\prism}\{1\})
  \]
  as a morphism of animated abelian group Gelfand stacks over $\Q_p^\prism$.

  Recall from \Cref{sec:analyt-prism-ga} the explicit description
\[
\Gm^{\prism} =\widetilde{\mathbb{G}}_{m,\Q_p^\prism}/\mathcal{G}_{\Q_p^\prism}.
\]
In particular, we get the connecting map
\[
\Gm^{\prism} \to \mathcal{G}_{\Q_p^\prism}[1]
\]
of animated abelian group stacks over $\Q_p^\prism$, that we can compose with the map of abelian group stacks over $\Q_p^\prism$
\[
\log_{\prism}\colon \mathcal{G}_{\Q_p^\prism}[1]  \to \mathbb{V}(\s{O}_{\Q_p^\prism}\{1\})[1],
\]
the prismatic logarithm, to get the desired map $c_1^\prism$.
\end{construction}

\begin{remark}
  \label{sec:chern-classes-p-weak-theory-of-chern-classes}
  \Cref{sec:chern-classes-p-construction-weak-theory-of-first-chern-classes} implies that the $6$-functor formalism $X\mapsto \ob{D}(X^\prism)$ on Berkovich spaces (equipped with the Berkovich-\'etale topology) admits a weak theory of first Chern classes in the sense of \cite[Definition 5.2.4]{zavyalov2023poincaredualityabstract6functor}.
  Indeed, this amounts to constructing a natural transformation
 \begin{equation}\label{weakchern}
  \Gamma(-,\Gm[1])\to \Gamma((-)^{\prism},\s{O}\{1\}[2])
 \end{equation}
 of Berkovich-\'etale sheaves with values in $\ob{D}(\Z)$, where the left hand side is cohomology on the Berkovich-\'etale site of Berkovich spaces.
 The right hand side satisfies Berkovich-\'etale-descent, cf.\ \Cref{TheoremDescentPrismatization}, so it suffices to construct a natural transformation
\[
\mathrm{Maps}_{\cat{GelfStk}}(-,B\Gm)= \tau_{\leq 0}\Gamma(-,\Gm[1])\to \tau_{\leq 0}\Gamma((-)^{\prism},\s{O}\{1\}[2])=\mathrm{Maps}_{\cat{GelfStk}}(-,\mathbb{V}(\s{O}_{\Q_p^\prism}\{1\}[2]))
\]
of sheaves of animated abelian groups (namely, the \'etale-sheafification in $\ob{D}(\Z)$ of $B\Gm=\tau_{\leq 0}\Gamma(-,\Gm[1])$ as an animated abelian group valued presheaf yields $\Gamma(-,\Gm[1])$), and this map can be given by composition with the map $c_1^\prism\colon B\Gm\to B^2\mathbb{V}(\s{O}_{\Q_p^\prism}\{1\})$ from \Cref{sec:chern-classes-p-construction-weak-theory-of-first-chern-classes}.
\end{remark}

\begin{remark}
  \label{sec:first-chern-class-compatibility-general}
  The general discussion in \cite[Remark 5.1.10.(3)]{drff_stacks} is natural in the following sense: given base Gelfand stacks $\mathcal{S}',\mathcal{S}$, a morphism $\mathcal{S}'\to \mathcal{S}$, ring stacks $\mathcal{R}'\to \mathcal{S}', \mathcal{R}\to \mathcal{S}$ and a morphism $\mathcal{R'}\to \mathcal{S}'\times_{\mathcal{S}}\mathcal{R}$ of ring stacks, the pullback along $(-)^{\mathcal{R}'}\to (-)^{\mathcal{R}}$ on transmutations yields a pullback on cohomology
  \[
    \Phi\colon \Gamma((-)^{\mathcal{R}},\mathcal{L})\to \Gamma((-)^{\mathcal{R}'},\mathcal{L})
  \]
  for any $\mathcal{L}\in \ob{D}(\mathcal{S})$ (we denote the pullback of $\mathcal{L}$ to stacks over $\mathcal{S}$ again by $\mathcal{L}$).
  Assume now that $\mathcal{L}$ is a line bundle on $\mathcal{S}$, and that we are given a commutative diagram
  \[
    \begin{tikzcd}
	{\Gm^{\mathcal{R}'}} & {\Gm^{\mathcal{R}}\times_{\mathcal{S}}\mathcal{S}'} \\
	{\mathbb{V}(\mathcal{L})\times_{\mathcal{S}}\mathcal{S}'} & {\mathbb{V}(\mathcal{L})\times_{\mathcal{S}}\mathcal{S}'}
	\arrow[from=1-1, to=1-2]
	\arrow["{c_1^{\mathcal{R}'}}"', from=1-1, to=2-1]
	\arrow["{c_1^{\mathcal{R}}\times_{\mathcal{S}}\mathcal{S}'}", from=1-2, to=2-2]
	\arrow[from=2-1, to=2-2]
      \end{tikzcd}\]
    of animated abelian group Gelfand stacks over $\mathcal{S}'$.
    Then the pullback $\Phi$ preserves first Chern classes of line bundles, i.e., if $X$ is a Gelfand stack with a line bundle $\mathcal{M}$, then
    \[
      \Phi(c^{\mathcal{R}}_1(\mathcal{M}))=c_1^{\mathcal{R}'}(\mathcal{M}).
    \]
  \end{remark}
  
    \begin{remark}[Compatibility with pro-\'etale Chern classes]
      \label{sec:first-chern-class-compatibility-with-pro-etale-chern-classes}
      From the reasoning in \Cref{sec:first-chern-class-compatibility-general} we can deduce compatibility of $c_1^\prism$ with pro-\'etale Chern classes by the following facts:
      \begin{enumerate}
      \item For an arc-stack $Z$ we set $\FFdisc{Z,(0,p)}^\times:=\FFdisc{Z}^\times\times_{(0,\infty)_{\Betti}}(0,p)_{\Betti}$ (using the radius map from \Cref{rk:radius-map} via $\FFdisc{Z}^\times\to \FFdisc{\Q_p}^\times\to \Q_p^\prism$).
        Then the pullback of $\s{O}_{\Q_p^\prism}\{1\}$ to $\FFdisc{\Q_p,(0,p)}^\times$ identifies canonically with the pullback along $\FFdisc{\Q_p,(0,p)}\to \FF_{\Q_p}$ of the Tate twist, seen as a line bundle on $\FF_{\Q_p}$.
        Indeed, this follows from the fundamental exact sequence of $p$-adic Hodge theory and \Cref{ConstructionBKtwist}, cf.\ \Cref{rk:line-bundles-on-qp-prism-adn-div1}.
      \item Given a Gelfand stack $X$, we have a natural map $\alpha_X\colon \FFdisc{X}^\times\to X^\prism$ (which can be thought of as the morphism $X^{\mathcal{R}'}\to X^{\mathcal{R}}$ with $\mathcal{R}'$ the ring stack $\FFdisc{\A^{1,\diamond}}^\times$ over $\FFdisc{\Q_p}^\times$ and $\mathcal{R}$ the ring stack $\A^{1,\prism}\to \Q_p^\prism$).
        The pullback along $\alpha_X$ is the pro-\'etale realization of rational analytic prismatic cohomology.
      \item We have (by the construction in \Cref{ConstructionPrismaticLog}) a commutative diagram
\[\begin{tikzcd}
	{\Z_p(1)^\cont\times{\FFdisc{\Q_p}^\times}} & {\widetilde{\bb{G}}_m^{u}\times \FFdisc{\Q_p}^\times} & {\FFdisc{\Gm^\diamond}^\times} \\
	{\mathcal{G}_{\Q_p^\prism}} & {\widetilde{\bb{G}}_m\times \Q_p^\prism} & {\Gm^\prism,}
	\arrow[from=1-1, to=1-2]
	\arrow[from=1-1, to=2-1]
	\arrow[from=1-2, to=1-3]
	\arrow[from=1-2, to=2-2]
	\arrow[from=1-3, to=2-3]
	\arrow[from=2-1, to=2-2]
	\arrow[from=2-2, to=2-3]
      \end{tikzcd}\]
    where $\widetilde{\bb{G}}_m^u$ is the uniform completion (locally on affinoids) of $\widetilde{\bb{G}}_m:=\varprojlim_{x\mapsto x^p}\Gm$.
    We can deduce a commutative diagram
  \[\begin{tikzcd}
	{\FFdisc{\Gm^\diamond}^\times} & {B_{\FFdisc{\Q_p}^\times}(\Z_p(1)^\cont\times{\FFdisc{\Q_p}^\times})} \\
	{\Gm^\prism} & {B_{\Q_p^\prism}(\mathcal{G}_{\Q_p^\prism})}
	\arrow[from=1-1, to=1-2]
	\arrow[from=1-1, to=2-1]
	\arrow[from=1-2, to=2-2]
	\arrow[from=2-1, to=2-2]
\end{tikzcd}\]
    \item Using (1), (3) (with $\Gm$ replaced by $1+\mathbb{D}_{\Q_p}$) and the diagram in \Cref{ConstructionTwists} we can deduce that we have commutative diagram
\[\begin{tikzcd}
	{\Z_p(1)^\cont\times \FFdisc{\Q_p,(0,p)}^\times} & {\mathcal{G}_{\Q_p^\prism}} \\
	{\mathbb{V}(\s{O}_{\mathrm{Tate},\FFdisc{\Q_p,(0,p)}^\times}(1))} & {\mathbb{V}(\mathcal{O}_{\Q_p^\prism}\{1\}),}
	\arrow[from=1-1, to=1-2]
	\arrow["{\mathrm{can}}"', from=1-1, to=2-1]
	\arrow["{\mathrm{log}_\prism}", from=1-2, to=2-2]
	\arrow[from=2-1, to=2-2]
      \end{tikzcd}\]
    where $\s{O}_{\mathrm{Tate},\FFdisc{\Q_p,(0,p)}^\times}(1)$ is the Tate twist (or rather its associated line bundle on $\FFdisc{\Q_p,(0,p)}^\times$, i.e, $\s{O}_{\mathrm{Tate},\FFdisc{\Q_p,(0,p)}^\times}(1)\cong \s{O}_{\mathrm{Tate},\FFdisc{\Q_p,(0,p)}^\times}\otimes_{\Z_p}\Z_p(1)$), and $\mathrm{can}$ is the canonical map from $\Z_p(1)^\cont$ to the global sections of $\s{O}_{\mathrm{Tate},\FFdisc{\Q_p,(0,p)}^\times}(1)$.
  \item The first Chern class in pro-\'etale cohomology is constructed using the Kummer sequence.
    In other words, it agrees with the (weak) theory of first Chern classes (for the $6$-functor formalism $\ob{D}(\FFdisc{(-)}^\times)$ on arc-stacks) constructed as in \Cref{sec:chern-classes-p-construction-weak-theory-of-first-chern-classes} via the morphism
    \[
      c_1^{\mathbb{D}_{\mathrm{FF}}^\times}\colon \FFdisc{\Gm^\diamond}^\times\to \Z_p(1)^\cont\times \FFdisc{\Q_p}^\times[1]\overset{\mathrm{can}}{\to} \mathbb{V}(\s{O}_{\mathrm{Tate},\FFdisc{\Q_p}^\times}(1))[1].
    \]
    We note that $c_1^{\mathbb{D}_{\mathrm{FF}}^\times}$ descents naturally to a morphism
    \[
      c_1^\proet\colon \FF_{\Gm^\diamond}\to \mathbb{V}(\s{O}_{\mathrm{Tate},\FF_{\Q_p}}(1))[1].
    \]
  \item From the previous points, we can conclude that after localizing to $\Q_{p,(0,p)}^\prism=\widetilde{\Div^1}_{(0,p)}$ and pulling back to the restriction of ${\FFdisc{\Q_p,(0,p)}^\times}$ the image of the prismatic first Chern class $c_1^\prism(\mathcal{M})$ for a line bundle $\mathcal{M}$ on a Gelfand stack $X$ agrees with the pullback of the pro-\'etale first Chern class $c_1^{\proet}(\mathcal{M})$ in $H^{2}_{\proet}(X^{\diamond},\bb{Q}_p(1))=H^2(\FF_{X^\diamond},\s{O}_{\ob{Tate},\FF_{X^\diamond}}(1))$.
   \item We note that over $\Q_{p,(0,p)}^\prism$ the first Chern classes descent naturally to the $p$-adic twistor stack $\Div^1$ (using \Cref{quotient-perfect}), and the same reasoning as before show that they pullback to the first pro-\'etale Chern classes along $\FF_{X^\diamond}\to X^\prism$. 
   \end{enumerate}

\begin{remark}
  \label{sec:first-chern-class-relation-to-pan-petrov}
  In ongoing work, Pan and Petrov define pro-\'etale Chern classes of $p$-adic local systems that should agree with the ones defined in \Cref{sec:first-chern-class-compatibility-with-pro-etale-chern-classes}. More precisely,  these Chern classes should be constructed from the cohomology of the map $(\ob{Vect})^{\n{R}}\to (\GSpec \bb{Q}_p)^{\n{R}}$ for the different rings stacks $\n{R}$ discussed in \Cref{sec:first-chern-class-compatibility-with-pro-etale-chern-classes}, where $\ob{Vect}$ the Gelfand stack of vector bundles. 
 A discussion of what happens in the de Rham case would be interesting, but seems to fit more naturally in the context of syntomic Chern classes (\cite{hauck}), which is beyond the scope of our paper.
\end{remark}

   We also have compatibility with Chern classes in Hyodo-Kato cohomology.
 \end{remark}
 \begin{remark}[Compatibility with Chern classes in Hyodo-Kato cohomology]
   \label{sec:first-chern-class-compatibility-with-hk-chern-classes}
   \cite[Lemma 6.2.9]{drff_stacks} constructed Chern classes in Hyodo-Kato cohomology (as well for the $6$-functor formalism $\ob{D}(\FFdisc{-}^{\times,\dR})$ on arc-stacks over $\Marc(\Z_p)$), and $c_1^\prism$ restricts to these in the following sense:
   after localizing to $\Q_{p,(1,\infty)}^\prism={\FFdisc{\Q_p}^{\times,\dR}}_{(1,\infty)}$, the prismatic first Chern class agrees with the pullback of the first Chern class in Hyodo-Kato cohomology along the natural morphism $X^{\prism}_{(1,\infty)}\to X^{\HK}$.

   This statement is easier to check than \Cref{sec:first-chern-class-compatibility-with-pro-etale-chern-classes}, namely it suffices to note the following:
 over $(1,\infty)_{\Betti}$ the Breuil-Kisin twist is isomorphic to $\s{O}_{\Q_p^\prism}\langle1\rangle$, $\mathcal{G}_{\Q_p^\prism}$ gets isomorphic to $\Z_p(1)^\sm\times \Q_{p,(1,\infty)}^\prism$ and $c_1^\prism$ agrees exactly with the map constructed in \cite[Lemma 6.2.9]{drff_stacks}.
 \end{remark}

 We can deduce that prismatic first Chern classes satisfy the projective bundle formula.
 
\begin{proposition}\label{LemComputationP1Div1}
  The 6-functor formalism of analytic prismatization on Berkovich spaces has a strong theory of first Chern classes $c_1^{\prism}$, in the sense of \cite[Definition 5.2.8]{zavyalov2023poincaredualityabstract6functor}.
  In fact, for any integer $d\ge 1$, and for any Gelfand stack $X$ over $\Q_p$, denoting $f\colon\mathbb{P}_X^{d}\to X$ the natural projection, we have that the morphism
\begin{equation}\label{projbundleformula}
 \sum_{k=0}^d c_1^k\{d-k\}[2d-2k]: \bigoplus_{k=0}^d \s{O} \{d-k\}[2d-2k] \longrightarrow f_{*}^{\prism} \s{O}\{d\}[2d]
\end{equation}
is an isomorphism; here, $c_1$ denotes the first Chern class of the universal line bundle.
\end{proposition}
\begin{proof}
  By \Cref{sec:chern-classes-p-construction-weak-theory-of-first-chern-classes} and \Cref{sec:chern-classes-p-weak-theory-of-chern-classes} we have a weak theory of first Chern classes on Berkovich spaces.

  To prove that it is strong, let us first unravel what \eqref{projbundleformula} is.
  The first Chern class $c_1^\prism$ gives a map $H^1(\mathbb{P}_X^d,\Gm) \to H^0(\mathbb{P}_X^{d,\prism}, \s{O}\{1\}[2])$, and hence the universal line bundle on $\mathbb{P}_X^d$ defines a class in $H^0(\mathbb{P}_X^{d,\prism}, \s{O}\{1\}[2])$; equivalently, a map $\s{O}\to  \s{O}\{1\}[2]$ on $\mathbb{P}_X^{d,\prism}$. By adjunction this is the same as a map $\s{O} \to f_\ast^\prism  \s{O}\{1\}[2]$ on $X^\prism$, which is $c_1$ and the $c_1^k$ are its tensor powers. To prove that \eqref{projbundleformula} is an isomorphism, since $\mathbb{P}_X^d\to X$ is Berkovich-smooth and the underlying morphism of arc-stacks is proper, the induced morphism on prismatizations is $0$-proper by \Cref{berkovich-smooth-qc-weakly-cohomologically-proper}, so we can by proper base change assume that $X=\GSpec(\Q_p)$. Using \Cref{cor:prismatization-push-out}, \Cref{sec:first-chern-class-compatibility-with-pro-etale-chern-classes} and \Cref{sec:first-chern-class-compatibility-with-hk-chern-classes}, proving that  \eqref{projbundleformula} is an isomorphism reduces to the analogous statement for the Hyodo--Kato stack and the $p$-adic twistor stack.
  For the Hyodo--Kato stack, the claim is proved in \cite[Lemma 6.2.9]{drff_stacks}. For the $p$-adic twistor stack, the right hand side of  \eqref{projbundleformula} is the pro-\'etale $\Q_p$-cohomology of $\mathbb{P}^d$ seen as a Galois-equivariant vector bundle on the Fargues-Fontaine curve, by \Cref{sec:mixed-case:-smooth-mixed-case1}, where the result is classical.
  Finally, the case where $X$ is a general Gelfand stack follows by base change.
\end{proof}

\subsection{Relation to the Riemann--Hilbert functor} \label{sec:relation-to-rh}
 Having discussed consequences for cohomology, we turn to a (brief) discussion of coefficients. Indeed, since the $p$-adic twistor stack geometrizes (for smooth partially proper rigid spaces over $\Q_p$) $p$-adic pro-\'etale cohomology, quasi-coherent sheaves on the $p$-adic twistor stack are coefficients for this cohomology theory and it is natural to wonder what they look like. We already observed (and this motivated our choice of terminology) that vector bundles on the $p$-adic twistor stack look like ``$p$-adic variations of twistor structures'', in that away from the Hodge--Tate locus, the $p$-adic twistor stack is isomorphic to the Hyodo-Kato stack (hence there one sees analytic $D$-modules on the relative Fargues--Fontaine curve without the de Rham-Hodge-Tate locus) and it is a gerbe for the overconvergent Tate-twisted tangent bundle on the Hodge--Tate locus (hence there one sees \emph{twisted} Higgs bundles). A natural next step is to try to unravel further what happens on the overconvergent (or formal) neighborhood of the Hodge--Tate locus. We now explain some results in this direction, in relation to the Riemann--Hilbert functor in $p$-adic Hodge theory (in both contexts, arithmetic and geometric). 

Let $X$ be a Gelfand stack. We have a natural morphism $X^{\Div^1} \xrightarrow{f} X^{\cone} \times_{ \bb{A}^/\bb{G}_m} \Div^1$, and  base  changes along $\widehat{\A}^1/\Gm \to \A^1/\Gm$  and  $\A^{1,\dagger}/\Gm\to \A^1/\Gm$ given by 
 \[
   X^{\widehat{\HT}}  \to X^{\widehat{\cone}} \mbox{ and }  X^{\dagger-\HT} \to X^{\dagger-\cone}
 \]
  respectively. 
   We call these stacks the \textit{formal} and \textit{overconvergent} Hodge--Tate and cone stacks respectively, see  
   \Cref{sec:defin-first-prop-1-definition-cone-construction} and \Cref{OverconvergentHTLocus}. Finally, we will consider the further base changes on the target
   \[
   \widehat{f}\colon   X^{\widehat{\HT}}     \to X^{\widehat{\cone}}\times_{\widehat{\A}^1/\Gm} \Q_p^{\widehat{\HT}} \mbox{ and } f^{\dagger}\colon X^{\dagger-\HT}\times_{\bb{G}^{\dagger}_a/\Gm} \Q_p^{\dagger-\HT}. 
   \]

\begin{proposition}\label{fg-suave}
  Let $X$ be a smooth rigid space over $\Q_p$.
  In the notations above, $\widehat{f}$ is a $\bb{T}_{X}^{\dR}\{1\}$-torsor where $\{1\}$ refers to a Breuil--Kisin twist, in particular $\widehat{f}$ is suave and the left adjoint $f_\natural$ to the pullback $f^\ast$ satisfies $f_\natural 1=1$. The same holds for $f^{\dagger}$ if $X$ is a partially proper smooth rigid space over $\Q_p$. 
\end{proposition}
\begin{proof}
  These claims are a rather direct consequences of \Cref{PropOverconvergentHT}.  For example, in the formal case for smooth rigid spaces one gets from \Cref{restriction-ht-overconvergent}  a cartesian square 
\[\begin{tikzcd}
	{X^{\widehat{\HT}}} & {X^\dR} \\
	{X^{\widehat{\cone}}} \times_{\widehat{\A}^1/\Gm} \Q_p^{\widehat{\HT}} & {X^{\widehat{\cone}},\dR} \times_{B\Gm^\dR} \GSpec(\Q_p)
	\arrow[from=1-1, to=1-2]
	\arrow[from=1-1, to=2-1]
	\arrow[from=1-2, to=2-2]
	\arrow[from=2-1, to=2-2]
      \end{tikzcd}\]
One has that $X^{\widehat{\cone},\dR} = \bb{T}_{X}^{\dR}(-1)[1]$  is the $1$-shifted de Rham stack of the tangent bundle  over $X^{\dR}\times B\bb{G}_m^{\dR}$ and twisted with weight $-1$. Taking a base change along $\Q_p^{\widehat{\HT}}\to \GSpec \bb{Q}_p$ and $X\to X^{\dR}$, one obtains a cartesian square
\[
\begin{tikzcd}
X^{\widehat{\HT}}  \ar[r] \ar[d] &  X\times \Q_p^{\widehat{\HT}} \ar[d] \\
X^{\widehat{\cone}} \times_{\widehat{\A}^1/\Gm} \Q_p^{\widehat{\HT}}\ar[r] & (X\times\Q_p^{\widehat{\HT}} ) / \bb{T}_X^{\dR}\{1\}
\end{tikzcd}
\] 
where $\{1\}$ now refers to a Breuil--Kisin twist (which corresponds to weight $-1$ as it corresponds to the ideal defined by the Cartier divisor).  One deduces that $\widehat{f}$ is a $\bb{T}_{X}^{\dR}\{1\}$-torsor, and therefore that it is suave with contractible homology.

The exact same argument holds for $f^{\dagger}$ for partially proper smooth rigid spaces by applying \Cref{PropOverconvergentHT} instead.
\end{proof}

\begin{remark}
  The above proof shows that the result holds as well for any $\dagger$-cohesive and formally smooth Gelfand stack. In the following we will focus in the case of formal completions as these are the ones leading to a functor that can directly be related to the existing $p$-adic Riemann--Hilbert functor.  However, see \Cref{Remark-relation-Wiersig-RH} for a comment about   overconvergent  Hodge--Tate stacks. 
 \end{remark}

Let $X$ be a smooth rigid space over $\bb{Q}_p$.  As a direct consequence of \Cref{fg-suave}, the functor
\[ 
\widehat{f}^\ast\colon \ob{D}(X^{\widehat{\cone}} \times_{\widehat{\A}^1/\Gm} \Q_p^{\widehat{\HT}}) \to \ob{D}(X^{\widehat{\HT}})
\]
is fully faithful.  Let $\Spf(B_\dR^+(\Q_p^\cyc)) =  \underset{n} \colim \GSpec(B_\dR^{+}(\Q_p^\cyc)/t^n)$, where $t$ is Fontaine's $2i\pi$. Passing to locally analytic vectors on  $B_\dR^{+}(\Q_p^\cyc)/t^n$ for the action of $\bb{Z}_p^{\times}$, one obtains from \Cref{thm:summary-main-results-sen} (1) that  
\[
(B_\dR^{+}(\Q_p^\cyc)/t^n)^{\la} = K_{\infty}[t]/t^n
\]
with $K_{\infty}= \bb{Q}_p(\zeta_{p^{\infty}})$. Hence, we have the presentation
\[
\Q_p^{\widehat{\HT}} =\Spf(K_{\infty}[[t]])/\Z_p^{\times,\la}.
\]

\begin{proposition}\label{Proposition-RH-Geometric-Case1}
Let $X$ be a smooth rigid space over $\Q_p$ and write $B_{\dR}^+=  B_{\dR}^{+}(\bb{Q}_p^{\cyc})$. After base change along $\Spf(B_\dR^+) \to \Q_p^{\widehat{\HT}}$ the pullback induces a fully faithful functor 
 \[
 \widehat{f}^{\ast}\colon \ob{D}(X^{\widehat{\cone}/ (B_{\dR}^+,t)})\to \ob{D}(X^{\widehat{\HT}}\times_{\Q_p^{\widehat{\HT}}} \Spf(B_{\dR}^+))
 \] 
with essential image  consisting on those sheaves in $\ob{D}(X^{\widehat{\HT}}\times_{\bb{Q}_p^{\HT}} \Spf(B_{\dR}^+))$ whose pullback to $X^\HT \times_{\Q_p^\HT} \GSpec(\Q_p^{\rm cyc})=X_{\Q_p^\cyc}^\HT$ have $\dagger$-nilpotent geometric Sen operators. Furthermore,  under the equivalence of perfect complexes
\[
\Perf(\FF^{\wedge}_{X_{\Q_p^\cyc}}) \cong \Perf( X^{\widehat{\HT}}\times_{\Q_p^{\HT}} \Spf (B_{\dR}^+))
\] deduced from  \Cref{sec:geom-case:-smooth-relative-case}, we have a natural identification between vector bundles on $X^{\widehat{\cone}/ (B_{\dR}^+,t)}$ and $\bb{B}_{\dR}^+$-local systems on $X_{\bb{Q}_p^{\cyc}}$ whose underlying $\widehat{\s{O}}$-vector bundle has nilpotent geometric Sen operators. 
\end{proposition}
\begin{proof}
Fully faithfulness is a consequence of \Cref{fg-suave} by proper base change. To describe the essential image of this functor, we  use the $\bb{T}_{X}^{\dR}\{1\}$-torsor structure of $\widehat{f}$. By descent along the map $\GSpec (\bb{Q}_{p}^{\cyc})\to \Spf B_{\dR}^+$, it  suffices to show that the essential image of the pullback along 
\[
X_{\Q_p^\cyc}^{\HT}\to X_{\bb{Q}_p^{\cyc}}/ \bb{T}_{X}\{1\},
\]
(which is the special fiber $\widehat{f}_s$ of $\widehat{f}$)   consists of the sheaves on   $X_{\Q_p^\cyc}^{\HT}$ with $\dagger$-nilpotent Sen operator.   Indeed, consider the sequence of maps 
\[
X_{\Q_p^\cyc}^{\HT} \to  X_{\bb{Q}_p^{\cyc}}/ \bb{T}_{X}\{1\}\to X_{\bb{Q}_p^{\cyc}}.
\]
The left morphism is a $\bb{T}_{X}^{\dR}\{1\}$-torsor, the right morphism a $B\bb{T}_{X}\{1\}$-torsor and the composite a $B\bb{T}_{X}^{\dagger}\{1\}$-torsor.  Now, the essential image along $\widehat{f}^*$ consist precisely of those sheaves  for which the action of $B\bb{T}_{X}^{\dagger}\{1\}$ by convolution integrates to an action of  $B\bb{T}_{X}\{1\}$, namely, by descent this can be proven locally on $X$ and hence it suffices for the case of trivial gerbes where it is a simple application of Cartier duality for vector bundles as in \cite{RCCartierGerbe}.  This integration property is equivalent, by Cartier duality, to the fact that the geometric Sen operators are $\dagger$-nilpotent.

For the final claim on vector bundles, we observe that sheaves on  $X^{\widehat{\cone}} \times_{\widehat{\A}^1/\Gm} \Spf(B_\dR^{+})$ identify with $t$-connections on $X\times_{\GSpec(\Q_p)} \Spf(B_\dR^{+})$, which are nilpotent modulo $t$ (\Cref{sec:cone-constr-sheaf-relation-to-t-connections}). Taking pullbacks along $\widehat{f}$ and then to $\FF^{\wedge}_{X_{\bb{Q}_p^{\cyc}}}$ we see that the associated $\bb{B}_{\dR}^+$-local system has nilpotent Sen operators. Conversely, if  $\bb{V}$ is a $\bb{B}_{\dR}^+$-de Rham local system with nilpotent Sen operators, the associated vector bundle $\bb{V}^{\la}$ on $ X^{\widehat{\HT}}\times_{\Q_p^{\widehat{\HT}}} \Spf (B_{\dR}^+)$ has nilpotent Sen operators, and by the first statement it belongs to the essential image of $\widehat{f}^*$. This finishes the proof of the proposition. 
\end{proof}

\begin{remark}\label{RemarkGeometricDiscussionGeneralPerfectoid}
 The proof of \Cref{Proposition-RH-Geometric-Case1} holds verbatim for the base change of $X$ to any perfectoid base $D$.  On the other hand, one can prove that $\bb{G}_a^{\dR}$ is Cartier self dual in Gestalten. This makes the Cartier dual of $\bb{T}^{\dR}_{X}\{1\}$ equal to $\bb{T}^{*,\dR}_{X}\{-1\}$ endowing $\ob{D}(X^{\widehat{\HT}}\times_{\Q_p^{\widehat{\HT}}} \Spf (\B_{\dR}^+))$ with the structure of a  $\ob{D}(\bb{T}^{*,\dR}_{X}\{-1\})$-linear category.  This action defines a notion of \textit{singular support} for sheaves on the completed Hodge--Tate stack. One can rephrase \Cref{Proposition-RH-Geometric-Case1} by saying that the essential image of $\widehat{f}^*$ is formed by the sheaves with zero singular support.

\end{remark}

Next, we discuss the relation with the $p$-adic Riemann--Hilbert functor of \cite{scholze2013p}. Recall that this functor sends a vector bundle $\mathcal{E}$ with a filtration by subbundles and a flat connection satisfying Griffiths transversality over $X$ to
\begin{equation}
\label{functor-rh-scholze}
F(\mathcal{E}):= \mathrm{Fil}^0(\mathcal{E}\otimes_{\s{O}} \s{O}\mathbb{B}_{\dR,X})^{\nabla=0},
\end{equation}
a $\mathbb{B}_\dR^+$-local system on $X$, cf.\ \cite[Theorem 7.6]{scholze2013p}. The functor $F$ is a fully faithful functor of abelian categories, but not fully faithful at the derived level. In order to have a clean  comparison between the functors $F$ and $\widehat{f}^*$, it is convenient to fix this failure of fully faithfulness. In the following we shall denote $\Spf(\bb{B}_{\dR,X}^+):=\FF^{\wedge}_{X^{\diamond}}$. 

\begin{construction}\label{Construction-pdR-variants}
Let $X$ be a smooth rigid space over $\bb{Q}_p$,  and consider the morphisms of Gelfand stacks
\[
\Spf (\bb{B}_{\dR,X}^+) \xrightarrow{\alpha_{X,\dR}} X^{\widehat{\HT}}\xrightarrow{\widehat{f}} X^{\widehat{\cone}}\times_{\widehat{\A}^1/\Gm} \bb{Q}_p^{\widehat{\HT}} \xrightarrow{h}  X^{\widehat{\cone}}
\]
All the maps involved are $0$-proper, and we have that $\alpha_{X,\dR,*}=1$, $\widehat{f}_{*} 1=1$ and $h_*1 = \bb{Q}_p\oplus \bb{Q}_p[1]$ where the last computation arises from the Galois cohomology of the filtered ring $\bb{B}_{\dR}$, see \Cref{sec:arithmetic-case:-p-d-ht-d-dr} (3).  The $H^1$-class of the cohomology of  $h$ gives rise to a $0$-proper morphism of analytic stacks
\[
h'\colon X^{\widehat{\cone}}\times_{\widehat{\A}^1/\Gm} \bb{Q}_p^{\widehat{\HT}} \xrightarrow{h}  X^{\widehat{\cone}} \times B\bb{G}_a^{\ob{alg}}. 
\]
satisfying $h'_{*}1 =1$ by design.  By taking pullbacks along this $\AnSpec (\bb{Q}_p)\to B\bb{G}_{a}^{\ob{alg}}$, we obtain a diagram of analytic  stacks
\[
\Spf (\bb{B}_{\ob{pdR},X}^+) \xrightarrow{\alpha_{X,\ob{pdR}}}  X^{\widehat{\ob{pHT}}} \xrightarrow{\widehat{f}_{\ob{pHT}}} X^{\widehat{\cone}}.
\]
All the previous arrows remain $0$-proper and they satisfy $\alpha_{X,\ob{pdR},*}1=1$ and $\widehat{f}_{\ob{pHT},*} 1=1$, in particular, $*$-pullback is fully faithful at the level of  categories of modules.   Finally, we shall write 
\[
\widetilde{F}\colon \Spf (\bb{B}_{\ob{pdR},X}^+)\to X^{\widehat{\cone}}.
\]
\end{construction}

\begin{proposition}\label{prop:comparison-with-rh}
Let $X$ be a smooth rigid space over $\Q_p$ and consider the map  $\widetilde{F}\colon \Spf (\bb{B}_{\ob{pdR},X}^+)\to X^{\widehat{\cone}}$. Then the functor $\widetilde{F}^*$ restricted to vector bundles on  $X^{\widehat{\cone}}$ is naturally equivalent to the base change of the Riemann--Hilbert functor of \eqref{functor-rh-scholze} along the map $\bb{B}_{\dR,X}^+\to \bb{B}_{\ob{pdR},X}^+$. 
\end{proposition}

\begin{proof}
Note that the natural map 
\[
X^{\widehat{\cone}} \to X^{\widehat{\mathrm{cone-alg}}}
\]
makes the source the ``Gelfandification'' of the target, and moreover, pullback along this map induces an equivalence on perfect complexes (this reduces, by completeness, to the special cones $X\times B \bb{T}_{X/\Q_p}(1)$ and $X\times B \bb{T}_{X/\Q_p}^{\mathrm{alg}}(1)$, where it can for example be seen using Cartier duality). Therefore, to continue with the comparison of both functors, we will continue working with the algebraic cone of $X$.

Consider the natural map 
\[
\widetilde{F}' \colon \Spf(\bb{B}_{\ob{pdR},X}^+)\to X^{\widehat{\ob{cone-alg}}},
\]
 we want to show that the functors $\widetilde{F}^{\prime,*}$ and  \eqref{functor-rh-scholze} are naturally equivalent on vector bundles, after base change to $B_{\ob{pdR}}$. By   \Cref{RemarkAlgebraicConeGeneralAnalyticStacks} (1) the map $X^{\widehat{\ob{cone-alg}}}\to \widehat{\A}^1/\Gm$ is $0$-affine given by the relative spectrum of the de Rham complex seen as a cdga over $|X|_{\Betti}$.  By \Cref{Construction-pdR-variants} we know that $\widetilde{F}^{\prime}_* 1 =1\in \ob{D}(X^{\widehat{\ob{cone-alg}}})$. Considering both objects as sheaves over $|X|_{\Betti}$, this equivalence is nothing but the natural equivalence of cdga
 \begin{equation}\label{eqNaturalIsoCohomologyBdR}
 \dR_{X/\Q_p} \cong R\nu_* \bb{B}_{\ob{pdR},X}
 \end{equation}
 which arises from the filtered algebra  $\bb{B}_{\ob{pdR},X}$, where $\nu\colon X_{\proet}\to X_{\an}$ is the projection of sites.  Let us recall how the classical computation of $R\nu_* \bb{B}_{\ob{pdR},X}$ goes. First, one has an equivalence of filtered algebras in chain complexes
 \[
\ob{Fil}^{\bullet}( \bb{B}_{\ob{pdR},X})=\ob{Fil}^{\bullet} [\s{O}\bb{B}_{\ob{pdR},X} \xrightarrow{\nabla}\s{O}\bb{B}_{\ob{pdR},X} \otimes_{\s{O}_X} \Omega^1_X\to \cdots \to  \s{O}\bb{B}_{\ob{pdR},X} \otimes_{\s{O}_X} \Omega^{\ob{dim X}}_X], 
 \]
 where $\Omega^i_X$ has degree $i$. Taking derived pushforwards along $\nu$  one gets an explicit representative of the filtered cdga (i.e. an explicit representative in $\ob{D}(|X|_{\Betti}\times \widehat{\A}^1/\bb{G}_m)$)
 \[
 \ob{Fil}^{\bullet} R\nu_{*} \bb{B}_{\ob{pdR},X} = \ob{Fil}^{\bullet}(\s{O}_X \xrightarrow{d_X}  \Omega^{1}_X\{-1\}\xrightarrow{d_X} \Omega^{2}_X\{-2\} \to  \cdots \to \Omega^{\ob{dim}(X)}_X\{-\ob{dim}(X)\}).  
 \]
 
 The previous shows that the natural isomorphism \eqref{eqNaturalIsoCohomologyBdR} arises from an automorphism of the de Rham complex $(\Omega^{\bullet}_{X},d_X)$ which is determined by the isomorphism of $\Omega^1_X$  sitting in degree $1$. This isomorphism is pinned down by its compatibility with first Chern classes discussed in  \Cref{sec:first-chern-class-compatibility-with-pro-etale-chern-classes} and the $\ob{dlog}$ map by \cite[Lemma 3.24]{scholze2013perfectoid}, hence it must be the identity.  
 
 From the previous, we deduce that the functor $\widetilde{F}^{\prime,*}$ sends a filtered vector bundle with flat connection $(\n{E},\nabla)$ to the base change along the morphism of  filtered  cdga algebras 
 \[
\dR_{X/\Q_p} \xrightarrow{\sim} R\nu_{*} \bb{B}_{\ob{pdR},X}  \to (\dR_{X/\Q_p}\otimes_{\s{O}_X} \s{O} \bb{B}_{\ob{pdR},X}, \nabla) \xrightarrow{\sim} \bb{B}_{\ob{pdR},X} 
 \]
 where the last equivalence is given  by taking horizontal sections.  This is precisely the functor \eqref{functor-rh-scholze} after base change to $B_{\ob{pdR}}$, finishing the proof of the proposition.  
 \end{proof}

Another way to get rid of the base change to $\Q_p^{\widehat{\HT}}$ is to consider the geometric situation, where $X$ is a smooth rigid space over $D$ perfectoid. Assume that $X$ lifts to a smooth  rigid space $\widetilde{X}$ over $\Spf(B_\dR^{+}(D))$. Let
$$
X^{\widehat{\HT}/B_\dR^{+}(D)} := X^{\Div^1/\FF_D} \times_{\FF_D} \Spf(B_\dR^{+}(D)).
$$
The relative version of \Cref{restriction-ht-overconvergent} gives us a cartesian square 
\[\begin{tikzcd}
	{X^{\widehat{\HT}/B_\dR^{+}(D)}} & {X^\dR \times_{D^\dR}  \Spf(B_\dR^{+}(D))} \\
	{X^{\widehat{\cone}} \times_{D^{\widehat{\cone}}}  \Spf(B_\dR^{+}(D))} & {X_s^{\cone,\dR} \times_{D_s^{\cone,\dR}}  \Spf(B_\dR^{+}(D)). }
	\arrow[from=1-1, to=1-2]
	\arrow[from=1-1, to=2-1]
	\arrow[from=1-2, to=2-2]
	\arrow[from=2-1, to=2-2]
      \end{tikzcd}\]
  
Since we have also a cartesian square (using that the cone construction commutes with pullbacks)
\[\begin{tikzcd}
	{X^{\widehat{\cone}}} & {D^{\widehat{\cone}}} \\
	{\widetilde{X}^{\widehat{\cone} }} &  {\Spf(B_\dR^{+}(D))^{\widehat{\cone}},}
	\arrow[from=1-1, to=1-2]
	\arrow[from=1-1, to=2-1]
	\arrow[from=1-2, to=2-2]
	\arrow[from=2-1, to=2-2]
      \end{tikzcd}\]
we see that 
$$
X^{\widehat{\cone}} \times_{D^{\widehat{\cone}}}  \Spf(B_\dR^{+}(D)) = \widetilde{X}^{\widehat{\cone}} \times_{\Spf(B_\dR^{+}(D))^{\widehat{\cone}}}  \Spf(B_\dR^{+}(D)) = X^{\widehat{\cone}/(D,B_\dR^{+}(D))}
$$
and thus we can rewrite our cartesian square as
\[\begin{tikzcd}
	{X^{\widehat{\HT}/B_\dR^{+}(D)}} & {X^\dR \times_{D^\dR}  \Spf(B_\dR^{+}(D)) = \widetilde{X}^{\dR/B_\dR^{+}(D)}} \\
	{X^{\widehat{\cone}/(D,B_\dR^{+}(D))}} & {(\widetilde{X}_s^{\cone})^{\dR/B_\dR^{+}(D)}= B_{\widetilde{X}^\dR} T_{\widetilde{X}/B_\dR^{+}(D)}^\dR. }
	\arrow[from=1-1, to=1-2]
	\arrow[from=1-1, to=2-1]
	\arrow[from=1-2, to=2-2]
	\arrow[from=2-1, to=2-2]
      \end{tikzcd}\]

We deduce from these considerations that the left vertical map
$$
f\colon X^{\widehat{\HT}/B_\dR^{+}(D)} \to X^{\widehat{\cone}/(D,B_\dR^{+}(D))}
$$
is suave, with $f_\natural 1=1$.

\begin{proposition}\label{prop:rh-geometric-case}
Let $X$ be a smooth rigid space over $D$ perfectoid. Assume that $X$ lifts to a smooth rigid space $\widetilde{X}$ over $\Spf(B_\dR^{+}(D))$. In the notations above, the pullback functor $f^\ast$ is a fully faithful functor from the category of flat $t$-connections (relative to $B_\dR^{+}(D)$), nilpotent modulo $t$, on $\widetilde{X}$ to the category of $\mathbb{B}_\dR^+$-local systems on $X$. 
\end{proposition}
\begin{proof}
Fully faithfulness along $f^*$ follows from suaveness of $f$ and the fact that $f_{\sharp}1=1$. It is left to identify vector bundles on  $X^{\widehat{\ob{cone}}/(D,B^+_{\dR}(D))}$ with $t$-connections nilpotent modulo $t$,   this follows from \Cref{sec:cone-constr-sheaf-relation-to-t-connections} (more precisely, one can pass to Tate-rational localizations to assume that $X$ is a Berkovich space, e.g., $X$ affinoid).
\end{proof}

\begin{remark}
  \label{sec:relat-riem-hilb-in-terms-of-a-period-sheaf}
  It is also possible to unravel \Cref{prop:rh-geometric-case} in terms of period sheaves.
  But the (filtered) structural $\mathbb{B}_\dR$-de Rham period sheaf $\mathscr{O}\mathbb{B}_{\dR,X}$, which geometrically arises through \Cref{prop:fiber-product-in-terms-of-faltings-extension} as the fiber product
  \[
    \Spf(\mathbb{B}^+_{\dR,X})\times_{X^{\widehat{\cone-\ob{alg}}}}(X\times \widehat{\mathbb{A}}^1/\mathbb{G}_m^{\ob{alg}}),
  \]
  has to be replaced by a new ($t$-complete) period sheaf $\mathscr{O}\mathbb{B}_{\dR,X}^{\widetilde{X}}$, depending on $\widetilde{X}$ and carrying a $t$-connection, whose formal spectrum is given by 
  \[
      \Spf(\mathbb{B}^+_{\dR,X})\times_{X^{\widehat{\cone-\ob{alg}}/(D,B^+_\dR(D))}}\widetilde{X}.
\]
\end{remark}

\begin{remark}\label{Remark-relation-Wiersig-RH}
Let $X$ be a smooth partially proper rigid space over $\Q_p$. Let $\FF_{X^\diamond}^\dagger$ be the overconvergent neighborhood of $\widehat{X^\diamond}$ in $\FF_{X^\diamond}$. Noting that $(\FF_{X^\diamond}^\dagger)^\dR = X^\dR$, we can also form the fiber product $\FF_{X^\diamond}^\dagger \times_{X^\dR} X$. By arguments similar to the ones of \Cref{LemmaReconstructionOBdR}, one sees that 
$$
\FF_{X^\diamond}^\dagger=\GSpec_{\FF_{X^{\diamond}}}(\B_{\dR,X}^{\dagger,+}), ~~ \FF_{X^\diamond}^\dagger \times_{X^\dR} X=\GSpec_{\FF_{X^\diamond}^\dagger}(\s{O}\B_{\dR,X}^{\dagger,+}),
$$
for the period sheaves $\B_{\dR,X}^{\dagger,+}, \s{O}\B_{\dR,X}^{\dagger,+}$ introduced in \cite{wiersig2023solution}. Let  $\GSpec_{\FF_{X^\diamond}^\dagger}(\B_{\dR,X}^{\dagger})$ be the closed subanalytic stack of $\FF_{X^\diamond}^\dagger$ obtained as the pullback of the $\Gm$-equivariant idempotent algebra $\Q_p\langle T \rangle_{\leq 0}[1/T]$ over $\Q_p\langle T \rangle_{\leq 0}$, along the natural map $\FF_{X^\diamond}^\dagger \to \A^{1,\dagger}/\Gm$. Let $\GSpec_{\FF_{X^\diamond}^\dagger}(\s{O}\B_{\dR,X}^{\dagger})$ be the pullback of $\GSpec_{\FF_{X^\diamond}^\dagger}(\B_{\dR,X}^{\dagger})$ to $\GSpec_{\FF_{X^\diamond}^\dagger}(\s{O}\B_{\dR,X}^{\dagger,+})$. Moreover, we have natural $\mathbb{G}_a^{\rm alg}$-torsors over them, induced by the composition
$$
\FF_{X^\diamond}^\dagger \to X^{\dagger-\HT} \to \Div^{1,\dagger} \to \A^{1,\dagger}/\Gm \times B\mathbb{G}_a^{\rm alg},
$$
where the last map is as discussed in \Cref{sec:arithmetic-case:-p-d-ht-d-dr} (3) (coming from Galois cohomology of $B_\dR$ in degree $1$). These are given by replacing in the above the various de Rham period sheaves by their \textit{presque} de Rham variants. 

The composite map
$$
w: \GSpec_{\FF_{X^\diamond}^\dagger}(\B_{\mathrm{pdR},X}^{\dagger}) \to \GSpec_{\FF_{X^\diamond}^\dagger}(\B_{\mathrm{pdR},X}^{\dagger,+})\to X^\dR
$$
is prim, as a composite of prim maps (the first map is a closed immersion, for the second, the statement being local on $X$, it can be reduced by smoothness and descent to the case of the perfectoid torus). Hence, to check that $w^\ast$ is fully faithful, it will be enough to show that the natural map $1\to w_\ast 1$ is an isomorphism. Since $X\to X^\dR$ is an epimorphism, this can be be checked after pullback along this map. By proper base change and the above formula, this amounts to showing that the pushforward of the unit along
$$
\GSpec_{\widehat{X^\diamond}}(\s{O}\B_{\mathrm{pdR},X}^{\dagger}) \to X 
$$
is $1$, i.e. that $R\nu_\ast \s{O}\B_{\mathrm{pdR},X}^{\dagger}=\s{O}_X$.  This fact is proven in \cite[Corollary VII.4.4]{WiersigProetalePeriodSheaves}, and the fully faithfulness of the map $w^*$ is an extension to analytic $D$-modules of the Riemann--Hilbert functor of \cite{WiersigFullyFaithfulRH}.    
\end{remark}

\begin{remark}
  \label{sec:relat-riem-hilb-hauck}
In \cite[\S 7]{hauck}, Hauck provides another, very interesting, geometric perspective on the Riemann--Hilbert functor from \cite{scholze2013p}, using the geometry of the (rational) analytic syntomification.
\end{remark}

\begin{remark}[Comparison with \cite{RHRobbaLiuZhu}] \label{relationrobaarh}
  We finish this section by comparing our work with the $p$-adic Riemann--Hilbert correspondence discussed in \cite[Theorem 1.1]{RHRobbaLiuZhu}. Let $C$ be a perfectoid field and set $A= \mathbb{A}_{\inf}(C)$.
  Let $\f{X}$ be a semistable $p$-adic formal  scheme over $\Spf(\n{O}_C)$, with rigid generic fiber $\mathcal{X}_\eta$.
  Assume given $\widetilde{\f{X}}\to \Spf(A)$ a semistable lift to $A$.  For an interval $I\subset (0,\infty)$ with rational ends, let  $\FFdisc{\f{X}_{\eta}^{\diamond},I}\subset \FFdisc{\f{X}_{\eta}}^{\times}$ be the Tate--rational subspace defined by the interval $I$, and let $\widetilde{\f{X}}_{I}:= \widetilde{\f{X}} \times_{\GSpec(A)} \FFdisc{C^{\flat},I}$.  Assume that $I\subset (\frac{1}{p-1}, p)$. The Riemann--Hilbert correspondence of  \cite[Theorem 1.1]{RHRobbaLiuZhu} produces an equivalence between vector bundles on the Gelfand stacks $\FFdisc{\f{X}^{\diamond}_{\eta},I}$ and $(\widetilde{\f{X}}_{I})^{\dR,+/\FFdisc{C^{\flat},I}}$ over $\FFdisc{C^{\flat},I}$ (with respect to the Hodge--Tate divisor of $\FFdisc{C^{\flat},I}$), under some smallness assumption that we do not recall here.  In the following we discuss a conceptual way to obtain a very similar result as \textit{loc.\ cit.}, at least for smooth formal schemes.\footnote{The extension to semistable schemes requires an extension to log-schemes that we will not pursue in this paper.}  We only sketch the main ideas; the details will appear elsewhere. 

 Following the proof of  \Cref{Proposition-RH-Geometric-Case1}, the key idea will be  to realize such a Riemann--Hilbert correspondence as induced by a pullback along a morphism of Gelfand stacks.  This morphism of stacks should arise from transmutation along a morphism of ring stacks, hence we will focus our discussion on the case $ \f{X}=\mathbb{G}_{a,\n{O}_{C}}$. There is also a mixed variant  of this construction of the Riemann--Hilbert correspondence, and in fact we will only discuss this one, i.e.\ take $\f{X}=\mathbb{G}_{a,\bb{Z}_p}$ instead. The geometric situation of \cite{RHRobbaLiuZhu} described above would follow from a relative variant of what is outlined below, in the exact same way as our above discussion of Riemann--Hilbert had an absolute and a relative incarnation.

 On one side of the correspondence, we want to see a full subcategory  of vector bundles over an affinoid on the Fargues--Fontaine disc of $\f{X}^{\diamond}_{\eta}$.
 Thanks to \Cref{sec:geom-case:-smooth-relative-case}, we can decomplete such vector bundles and consider instead $(\f{X}_{\eta})^{\widetilde{\Div^1}}= (\bb{G}_{a,\bb{Z}_p,\eta})^{\widetilde{\Div}^1}$; this step of the discussion should be compared with the  decompletion result \cite[Theorem 1.2]{RHRobbaLiuZhu}.
 Since we have seen throughout the paper that it is more convenient to replace affinoids by overconvergent spaces when possible, we shall work instead with the ring Gelfand stack $\bb{G}_{a}^{\leq 1, \widetilde{\Div^1}}\to \widetilde{\Div}^1$. For an  interval $I\subset (0,\infty)$, we let  $(\bb{G}_{a}^{\leq 1, \widetilde{\Div}^1})_I$ denote the pullback to $I\subset (0,\infty)$ along the radius map $\kappa\colon \widetilde{\Div^1}\to (0,\infty)_{\Betti}$.\footnote{In particular, for closed intervals we are considering $\dagger$-rational localizations and not Tate-rational localizations.}  The ring Gelfand stack  $(\bb{G}_{a}^{\leq 1, \widetilde{\Div}^1})_I$ (for a well-chosen interval $I$) will play the role of the left hand side of the Riemann--Hilbert correspondence in  \cite[Theorem 1.1]{RHRobbaLiuZhu}. 
 
Moving to the other side, we define the ring stack corresponding to (a variant of) small $t$-connections, namely, the ring stack that will replace the right hand side of \cite[Theorem 1.1]{RHRobbaLiuZhu}. For that, we make use of the arithmetic ring stack $\bb{G}_a^{\mathrm{arith}}=\bb{G}_a^{\leq 1}/\bb{G}_a^{<1}$ of \cite{qin2026arithmetic}. We recall the functor of points description of   $\bb{G}_a^{\mathrm{arith}}$. Given a nilperfectoid ring $R$,  we write $\overline{R}:= R^{\leq 1} /R^{<1}$, this ring can be equivalently described as 
\[
\overline{R}= R^{u,\circ}/R^{u,<1}= R^{u,\flat,\circ}/R^{u,\flat, <1}.
\] 
In particular, we have $\bb{G}_a^{\mathrm{arith}}(R)= \overline{R}$ and it is a perfect $\bb{F}_p$-algebra.  The morphism of ring stacks $\bb{G}_a^{\leq 1}\to \bb{G}_a^{\mathrm{arith}}(R)$ gives rise to a natural deformation  $\bb{G}_a^{\mathrm{arith},+}\to \bb{D}^{<1}/\bb{T}^{\dagger}$ (with $\bb{T}^{\dagger}=\GSpec(\bb{Q}_p\langle T^{\pm 1}\rangle_{\leq 1})$ the overconvergent torus), defined as follows:
\begin{equation}\label{eqRingStack}
\bb{G}_a^{\mathrm{arith},+}(\GSpec(R)\to \bb{D}^{<1}/\bb{T}^{\dagger}) =  R^{\leq 1}/ (I\otimes_{R^{\leq 1}} R^{<1}),
\end{equation}
where the map $I\to R^{\leq 1}$ is the (topologically nilpotent) Cartier divisor produced by the natural map to  $\bb{D}^{<1}/\bb{T}^{\dagger}$. There is a natural morphism of ring stacks 
$$
f\colon \bb{G}_a^{\leq 1, \widetilde{\dR}, +}:=\bb{G}_{a}^{\leq 1, \dR,+}\times_{\bb{A}^1/\bb{G}_m} \bb{D}^{<1}/\bb{T}^{\dagger} \to \bb{G}_a^{\ob{arith},+}
$$ 
which is suave and satisfies $f_{\sharp} 1 =1$. As a consequence, $f^{*}$ is fully faithful (and similarly for the transmutation along  smooth schemes of finite type over $\bb{Z}_p$).  This allows to describe vector bundles on $\bb{G}_a^{\ob{arith},+}$ as filtered connections satisfying a suitable convergence condition (or equivalently, a smallness condition).\footnote{If we were considering the geometric case, over a perfectoid base, we would consider a relative variant of the $+$-arithmetic de Rham stack, on which, under a liftability assumption, vector bundles are small $t$-connections, recovering the set up of \cite{RHRobbaLiuZhu}. Compare with the discussion of \Cref{sss:ConeFlatConnection}.}  

Finally, one can relate the two sides. Consider the map 
\[
\mu\colon \widetilde{\Div^1}\to \bb{A}^1/\bb{G}_m
\] 
defined by the universal degree $1$ divisor over perfect $\delta$-spaces. As shown in \cite[Lemma 3.20]{hauck}, the map $\mu$ has a natural lift
\[
\widetilde{\mu}\colon  \widetilde{\Div^1}\to \bb{D}^{<1}/\bb{T}^{\dagger}.
\]
 One can show using the functor of points description of both sides that there is a natural morphism of ring Gelfand stacks
\[
g\colon (\bb{G}_{a}^{\leq 1})^{\widetilde{\Div^1}} \to \bb{G}_a^{\mathrm{arith},+}\times_{\bb{D}^{<1}/\bb{T}^{\dagger},\widetilde{\mu}} \widetilde{\Div^1}
\]
over $\widetilde{\Div^1}$. Moreover, this morphism is suave and if we set $I:= (1/(p-1),p)$,  the restriction to 
\[
g_{I}\colon  (\bb{G}_{a}^{\leq 1})^{\widetilde{\Div^1}}_{I}\to \bb{G}_a^{\mathrm{arith},+}\times_{\bb{D}^{<1}/\bb{T}^{\dagger},\widetilde{\mu}} \widetilde{\Div^1_I}
\] 
satisfies $g_{I,\sharp} 1\xrightarrow{\sim} 1$. In particular, $g_{I}^*$ is fully faithful. This provides the desired stacky reformulation (of a variant) of the smooth case of \cite{RHRobbaLiuZhu}.\footnote{We observe that this variant will not directly produce the case of $p$-adic formal schemes as the ring stack $\bb{G}^{\leq 1}_a$ is not $p$-adically complete. There are two possible ways around  this problem: first, one could work instead with some theory of weak formal schemes as in \cite{MR330167} or \cite{ErtlYamadaRigidHK}, where the building blocks are integral overconvergent algebras $\bb{Z}_p\langle  T_1,\ldots, T_n \rangle_{\leq 1}$. A second possible  solution is to consider a \textit{solid variant} of all the previous constructions replacing $\bb{G}_{a}^{\leq 1}$ by $\bb{G}_{a,\sol}$. This corresponds to replacing all instances of $R^{\leq 1}$ by $R^{+}$  in the previous constructions. This second approach will allow the transmutation of $p$-adic formal schemes, but requires an extension of the theory of Gelfand stacks  containing the ring $\bb{Q}_{p}\langle T  \rangle_{\sol}$. }
\end{remark}

\newpage
\appendix

\section{Deformation theory for Gelfand stacks}\label{app:DefTheoryGelfand}\label{appendix-def-theory}

This appendix develops (a bit of) deformation theory in Gelfand analytic geometry. 

\subsection{Cotangent complexes of Gelfand rings}
\label{sec:cotang-compl-gelf-gelfand-rings}

In this subsection, we provide some results on cotangent complexes for morphisms of Gelfand rings (following the discussions in \cite[Section 3.4]{RCAndeRham} and \cite[Section 7.4]{solidnotes}).\footnote{The theory can be generalized to analytic rings, cf.\ \cite[Section 3.4]{RCAndeRham}, but we don't need this generality.}

Given a morphism $A\to B$ of solid animated $\Z$-algebras, we recall that its cotangent complex $\mathbb{L}_{B/A}\in \ob{D}_{\geq 0}(B)$ is the connective complex corepresenting the functor
\[
  M\in \ob{D}_{\geq 0}(B)\mapsto \Der_A(B,M),
\]
where $\Der_A(B,M):=\Hom_{A}(B,B\oplus M)\times_{\Hom_A(B,B)}\{\Id_B\}$ is the $B$-module of $A$-linear derivations, cf.\ \cite[Definition 7.4.1]{solidnotes}.

\begin{lemma}
  \label{sec:cotang-compl-gelf-1-properties-of-cotangent-complexes}
  \begin{enumerate}
  \item Let $A\to A', A\to B$ be morphisms of solid animated $\Z$-algebras with pushout $B'$.
    Then the natural morphism $\mathbb{L}_{B/A}\otimes_{B}B'\to \mathbb{L}_{B'/A'}$ is an equivalence.
  \item If $A\to B$ is the colimit of a diagram of morphisms $\{A_i\to B_i\}$ of solid animated $\Z$-algebras, then the natural map $\varinjlim_i B\otimes_{B_i}\mathbb{L}_{B_i/A_i}\to \mathbb{L}_{B/A}$ is an equivalence.
  \item Let $A\to B\to C$ be morphisms of solid animated $\Z$-algebras.
    Then the natural morphism $C\otimes_{B}\mathbb{L}_{B/A}\to \mathbb{L}_{C/A}\to \mathbb{L}_{C/B}$ form a natural fiber sequence of $C$-modules.
  \item If $A\to B$ is idempotent, then $\mathbb{L}_{B/A}=0$.
  \end{enumerate}
\end{lemma}
\begin{proof}
  This is \cite[Lemma 7.4.3]{solidnotes}.
\end{proof}

The cotangent complex has strong descent properties.

\begin{lemma}\label{LemmaDescentLAB}
  Let $R$ be a solid animated $\Z$-algebra.
  Let $A\to A_0$ be a descendable morphism of solid animated $R$-algebras.
  Let $A\to A_{\bullet}$ be its \v{C}ech co-nerve.
  Then for all $m\in \N$ the natural map 
\[
\wedge^{m}_{A} \mathbb{L}_{A/R}\to \ob{Tot}(\wedge^m_{A_{\bullet}} \mathbb{L}_{A_{\bullet}/R})
\]
is an equivalence. 
\end{lemma}

We observe that  the same proof applies when $A\to A_0$ is a descendable morphism of analytic rings, with $A_0$ having the induced analytic ring structure.
\begin{proof}
  For simplicity of notation we suppress $R$ and just write $\mathbb{L}_{(-)}:=\mathbb{L}_{(-)/R}$.
  By \cite[Proposition 3.20]{mathew2016galois}, the map $A\to A_0$ of induced analytic rings is descendable if and only if the tower $\big((\ob{Tot}_{\leq n}(A_{\bullet})\big))_{n\in \N}$ is pro-constant and isomorphic to $A$.
  Equivalently, if and only if there is some $n\geq 0$ such that, if $K_n=\ob{fib}(A\to \ob{Tot}_{\leq n}(A_{\bullet}))$, the map $K_n\to A$ is homotopic to zero as $A$-module.
  Equivalently, the map $A\to \ob{Tot}_{\leq n}(A_{\bullet})$ admits a retract $\ob{Tot}_{\leq n}(A_{\bullet})\to A$ as an $A$-module. 

  We first prove the case of $m=1$, i.e., of $\mathbb{L}_{(-)}$.
  We have a fiber sequence of connective  $A_{\bullet}$-modules
\begin{equation}\label{eqFiberSequenceCotangentComplex}
A_{\bullet}\otimes_{A} \mathbb{L}_A\to \mathbb{L}_{A_{\bullet}}\to \mathbb{L}_{A_{\bullet}/A}.
\end{equation}
By descent, we know that 
\begin{equation*}
A_{\bullet}\otimes_{A} \mathbb{L}_A=\ob{Tot}(A_{\bullet}\otimes_{A}\mathbb{L}_{A}),
\end{equation*}
thus, it suffices to show that $\ob{Tot}(\mathbb{L}_{A_{\bullet}/A})=0$.  To prove this, it suffices to consider the base change 
\[
\mathbb{L}_{A_{\bullet}/A}\otimes_A \ob{Tot}_{|[m]|\leq n}(A_{m})= \ob{Tot}_{|[m]|\leq n}(\mathbb{L}_{A_{\bullet}/A} \otimes_A A_m)=\ob{Tot}_{|[m]|\leq n}(\mathbb{L}_{A_{\bullet+m+1}/A_m})
\]
for $[m]\in \Delta_{\leq n}$ and $n\in \N$. Choosing $n$ such that $A$ is a retract of $\ob{Tot}_{\leq n}(A_{\bullet})$, we see that $\mathbb{L}_{A_{\bullet}/A}$ is a retract of $\ob{Tot}_{|[m]|\leq n}(\mathbb{L}_{A_{\bullet+m+1}/A_m})$ as a cosimplicial $A$-module, and taking limits that $\ob{Tot}(\mathbb{L}_{A_{\bullet}/A})$ is a retract of $\ob{Tot}_{|[m]|\leq [n]} (\ob{Tot}(\mathbb{L}_{A_{\bullet+m+1}/A_m}))$, but $\ob{Tot}(\mathbb{L}_{A_{\bullet+m+1}/A_m})=0$ for $m\geq 0$ since the augmented simplicial diagram $A_m\to A_{\bullet+m+1}$ is split. It follows that $\ob{Tot}(\mathbb{L}_{A_{\bullet}/A})=0$ proving what we wanted.

To handle the general case of $\wedge^m_A \mathbb{L}_A$, we use the length  $m$ filtration induced by the fiber sequence \eqref{eqFiberSequenceCotangentComplex} (cf.\ \cite[Example 2.2]{bhatt2019topological}) and induction on $m$ to reduce to the descent for $m=1$. \end{proof}

We now show that cotangent complexes can be calculated as desired.

\begin{lemma}
  \label{sec:cotang-compl-morph-cotangent-complex-via-p-completion}
  
  Let $A$ be a solid animated $\Z$-algebra.
  \begin{enumerate}
  \item If $A$ is discrete, then $\mathbb{L}_{A/\Z}$ is discrete and $\mathbb{L}_{A/\Z}(\ast)=\mathbb{L}_{A(\ast)/\Z}$.
  \item If $A$ is $p$-complete, then $\mathbb{L}_{A/\Z_p}$ is $p$-complete.
  In particular, $\mathbb{L}_{A[1/p]/\Q_p}\cong (\varprojlim_{n} \mathbb{L}_{A/p^n/\Z/p^n})[1/p]$.
  \end{enumerate}
\end{lemma}
\Cref{sec:cotang-compl-morph-cotangent-complex-via-p-completion} implies that the solid cotangent complex agrees with previous constructions of cotangent complexes, e.g., for $p$-complete nuclear animated $\Z_p$-algebras or Tate-Huber rings (\cite{guo2021period}).
\begin{proof}
  For the first statement see \cite[Example 7.4.4.1]{solidnotes}.
  For the second, we follow an argument that arose from conversations with Akhil Mathew.  Assume that $A$ is $p$-complete.
  We first make the following observations:
  \begin{enumerate}
  \item The cotangent complex of $A$ is the solidification of the sheafification of the presheaf of $A$-modules $S\mapsto \mathbb{L}_{A(S)/\Z(S)}$ on light profinite sets $S$.
  \item We can consider the (derived) Hochschild homology $\mathrm{HH}(A/\Z)$ of $A$, which is the solid tensor product $A\otimes_{A\otimes_{\Z}A}A$.
    In particular, $\mathrm{HH}(A/\Z)$ is $p$-complete as $A$ is $p$-complete and connective (and implicitly we use that $\Z\to \Z_p$ is idempotent, which implies $\mathrm{HH}(A/\Z)=\mathrm{HH}(A/\Z_p)$).
  \item Equivalently, $\mathrm{HH}(A/\Z)$ is the solidification of the sheafification of the presheaf $S\mapsto \mathrm{HH}(A(S)/\Z(S))$ on light profinite sets.
    In particular, there exists as a consequence of the HKR-theorem an $\N$-indexed decreasing filtration $\Fil^\bullet_{\mathrm{HKR}} \mathrm{HH}(A/\Z)$ of $\mathrm{HH}(A/\Z)$ by solid $A$-modules whose $i$-th graded piece $\mathrm{gr}^i_{\mathrm{HKR}} \mathrm{HH}(A/\Z)$ is $\wedge^i_{A}\mathbb{L}_{A/\Z}[i]$ (\cite[Theorem 6.2.6]{raksit2020hochschild}).
  \item For each integer $\ell\geq 2$ prime to $p$ the spectrum $\mathrm{HH}(A/\Z[1/\ell])$ carries a functorial operator $\psi^\ell\colon \mathrm{HH}(A/\Z[1/\ell])$, the $\ell$-th Adams operation, \cite[Construction 6.4.3]{raksit2020hochschild}.
    These operators respect the HKR-filtration (\cite[Proposition 6.4.4]{raksit2020hochschild}), and act by $\ell^i$ on the $i$-th graded object (\cite[Proposition 6.4.12]{raksit2020hochschild}).
    Indeed, the operator $\psi^\ell$ exists already on the presheaf $S\mapsto \mathrm{HH}(A(S)/\Z[1/\ell](S))$ of spectra on light profinite sets, with the desired preservation of the HKR-filtration and description on the associated gradeds.
    We can moreover conclude that $\psi^\ell$ is an automorphism of $\mathrm{HH}(A/\Z[1/\ell])$.
  \end{enumerate}

  Fix some integer $\ell\geq 2$ prime to $p$.
  We want to deduce the $p$-completeness of $C:=\mathbb{L}_{A/\Z}\cong \mathbb{L}_{A/\Z[1/\ell]}$ from the $p$-completeness of $M:=\mathrm{HH}(A/\Z)\cong \mathrm{HH}(A/\Z[1/\ell])$ by the use of Adams operations.
  Let $F^\bullet$ be the HKR-filtration on $M$ and fix $i\geq 0$.     We need to see that $\pi_i(C)$ is $p$-complete. The case $i=0$ is clear since $\pi_0(\bb{L}_{A/\bb{Z}})$ is isomorphic to the quotient of derived $p$-complete static modules $I/I^{2}$ where $I=\ob{ker}(\pi_0(A\otimes_{\Z_{p}} A)\to \pi_0(A))$. 
  Using induction and $\pi_0(M)=A$ we may assume that $\pi_{i'}(C)$ is $p$-complete for $0\leq i'<i$.
  We note that by $p$-completeness of $A$, the object $F^1$ is $p$-complete as well. 
  Let $j\gg i$ such that $\pi_{i'+2}(F^1)\cong \pi_{i'+2}(F^1/F^j)$ for $0\leq i' \leq i$ (this is possible by the connectiveness of the cotangent complex and its wedge powers).
  We consider the complex $N:=\mathrm{fib}(F^1/F^j\xto{\psi^\ell-\ell\mathrm{Id}} F^1/F^j)$, i.e., the eigenspace of $\psi^\ell$ for the multiplication by $\ell$-action on $F^1$.
  The filtration $F^\bullet$ induces a filtration $G^j=0\to G^{j-1}\to \ldots G^1=N$ on $N$ (by stability of the HKR-filtration under the Adams operations), whose graded piece $\mathrm{gr}^k$ is $\mathrm{fib}(F^k/F^{k+1}\xto{\ell^k-\ell} F^k/F^{k+1})$ for $1\leq k\leq j-1$.
  As $\ell^k-\ell\neq 0$ for $k\geq 2$, we see that the graded pieces are uniformly $p^n$ torsion for $2\leq k\leq j-1$ (for some $n$ only depending on $j$).
  In fact, we can conclude (by the finiteness of the filtration $G^\bullet$) that there exists some $m$ such that $\pi_k(G^2)$ is $p^m$-torsion for any $k\geq 0$.
  However, $G^1/G^2\cong C[1]\oplus C[2]$ as $\psi^\ell-\ell \mathrm{Id}=0$ on $G^1/G^2$, and then the $p$-completeness of $\pi_{i'+2}(N)$ for $0\leq i'\leq i$ implies the $p$-completeness of $\pi_i(C)$. 
\end{proof}

\begin{example}
  \label{sec:cotang-compl-morph-b-dr-as-universal-deformation}
  \Cref{sec:cotang-compl-morph-cotangent-complex-via-p-completion} and Fontaine's classical computation of $\Omega^1_{\s{O}_{\C_p}/\Z_p}$ imply that $\mathbb{L}_{\C_p/\Q_p}\cong \C_p(1)[1]$, and hence that there exists\footnote{See \url{https://mathoverflow.net/questions/407484/witt-vectors-the-cotangent-complex-and-a-solid-construction-of-b-dr} for a discussion.} a universal solid pro-nilpotent thickening $B$ of $\C_p$ as a $\Q_p$-algebra with ideal $I:=\mathrm{fib}(B\to \C_p)$ and $I/I^2\cong \C_p(1)$.
  This is of course exactly Fontaine's ring $B_\dR^+$!
  \end{example}

The existence of a well-behaved theory of cotangent complexes for Gelfand stacks, instead of rings, is not automatic. We leave its development for the following section that also discusses a reasonable deformation theory context for Gelfand analytic geometry.
We shall follow the deformation theory formalism of \cite[Chapter 7.3]{lurie_higher_algebra} and \cite[Chapter V.17]{Luriea}.

\subsection{Deformation theoretic properties}\label{s:DefTheoryProperties}
There are three key properties in any deformation theory context: nilcompleteness, infinitesimal cohesivity, and the existence of the cotangent complex. 

Nilcompleteness compares the value of a stack $X$ at a ring $A$ with the limit of its values at the Postnikov truncations $\tau_{\leq n} A$.
Similarly, infinitesimal cohesivity compares the value of $X$ on a square-zero extension $\widetilde{A}$ arising from a derivation $d\colon A\to M[1]$ with the fiber product $X(A)\times_{X(A\oplus M[1])} X(A)$, where the left map  is induced by the zero derivation $A\to M[1]$ and the right map by the derivation $d$. The existence of the cotangent complex roughly means that the deformation theory of a Gelfand stack is controlled by linear data given by the cotangent complex. In the more familiar setting of derived algebraic geometry, these hypotheses together allow one to control via the cotangent complex a derived stack from its classical truncation, using Postnikov tower induction, e.g.\ in Lurie's version of Artin's representability theorem.

 One technical  problem in our setup of Gelfand stacks compared with classical algebraic geometry is that we use the $!$-topology, which is the strongest that guarantees descent of quasi-coherent sheaves after arbitrary base change. Unfortunately, this topology  may lose deformation theoretic information.

\begin{example}\label{ExampleDeformationTheoryhatGa}
The following example illustrates why the $!$-topology is not well suited to do deformation theory in a classical sense. The map $\GSpec (\bb{Q}_p)\to \widehat{\mathbb{A}}^1$ is an epimorphism for the $!$-topology inducing a natural equivalence
\begin{equation}\label{eqLoopshatGa}
B\bb{G}_a[-1]=\widehat{\mathbb{A}}^1,
\end{equation}
where $\bb{G}_a[-1]=\GSpec( \bb{Q}_p[\epsilon_1])$ where $\epsilon_1$ is a free variable in homological degree $1$ (indeed, $\bb{G}_a[-1]=\GSpec(\Q_p)\times_{\widehat{\mathbb{A}}^1}\GSpec(\Q_p)$).
However, in the standard context of deformation theory as in \cite[Chapter IV]{Luriea},  it is not true that $\GSpec (\bb{Q}_p)\to \widehat{\mathbb{A}}^1$ is a surjection, namely, in this context a surjection would mean a surjection of presheaves on artinian $\mathbb{\Q}_p$-algebras.

We must highlight that the identification \eqref{eqLoopshatGa} is not a bug but a feature of the theory! Indeed, when applied to the quotient stack $\widehat{\bb{A}}^1/\bb{G}_m$ we end up with the equivalence
\[
\widehat{\bb{A}}^1/\bb{G}_m=B(\bb{G}_a[-1]\rtimes \bb{G}_m),
\]
and the right term can be directly related with the category of differential graded modules (up to a shearing, \cite[Example 5.1.4]{RCBBCourse2026}). 
Moreover, this equivalence explains geometrically Beilinson's $t$-structure in filtered modules (as it identifies with the natural  $t$-structure on $B(\bb{G}_a[-1]\rtimes \bb{G}_m)$ defined by $\ast$-descent and after shearing).

\end{example}

We introduce the following deformation-theoretic notions.
Here, $\cat{GelfStk}^{\ob{std}}$ denotes the category of affinoid strictly totally disconnected Gelfand stacks.
 
 \begin{definition}  \label{sec:cotang-compl-morph-definition-having-a-cotangent-complex}
   Let $f\colon Y\to X$ be a map of presheaves (of anima) on affinoid Gelfand rings.
   We make the following definitions.
 \begin{enumerate}
 
\item We define the functor $Y^{\ob{nil}}\colon \cat{GelfStk}^{\ob{std}}\to \cat{Ani}$ to be the one sending a strictly totally disconnected\footnote{This assumption is convenient, e.g., for the $\dagger$-quasi-pro-\'etale topology, cf.\ \Cref{sec:deform-theor-prop-flatness-for-dagger-quasi-pro-etale-maps}.} Gelfand ring $A$ to 
\[
Y^{\ob{nil}}(A)=\varprojlim_n Y(\tau_{\leq n} A).
\]
We define $Y^{\ob{nil}/X}:=Y^{\ob{nil}}\times_{X^{\ob{nil}}} X$.  We say that $f$ is \textit{nilcomplete} if the natural map  $Y\to Y^{\ob{nil}/X}$ of functors $\cat{GelfRing}^{\ob{std}}\to \cat{Ani}$  is an equivalence.  We say that $Y$ is \textit{nilcomplete} if $Y\to \GSpec(\Q_p)$ is nilcomplete.

\item Let $\cat{Sqz}$ denote the category of square-zero extensions, i.e., the data of square-zero extensions $\widetilde{A}_d\to A$ with associated derivation $d\colon A\to M[1]$ with $A$ a Gelfand ring and $M\in \ob{D}_{\geq 0}(A)$ the ideal of the square-zero extension.
  We define the functors $\iota, (-)^{\ob{coh}}\colon \cat{Sqz}\to \cat{Ani}$
\[
\iota(Y)(d\colon A\to M[1]):= Y(\widetilde{A}_d) \mbox{ and } Y^{\ob{coh}}(d\colon A\to M[1]):= Y(A)\times_{Y(A\oplus M[1])} Y(A) ,
\]
where the left map $Y(A)\to Y(A\oplus M[1])$ corresponds to the zero derivation $0\colon A\to M[1]$ and the right map to the derivation $d$. We have a natural transformation of functors $\iota(Y)\to Y^{\ob{coh}}$, and we denote $Y^{\ob{coh}/X}=Y^{\ob{coh}}\times_{X^{\ob{coh}}} \iota(X)$. We say that $f$ is \textit{infinitesimally cohesive} if the natural map $\iota(Y)\to Y^{\ob{coh}/X}$ is an equivalence.  We say that $Y$ is \textit{infinitesimally cohesive} if $Y\to \GSpec(\Q_p)$ is infinitesimally cohesive.

\item  We say that $f$ admits a cotangent complex if there is a locally eventually connective object\footnote{This condition is there to guarantee that if the cotangent complex exists then it is unique up to a contractible anima.}, cf. \Cref{Subsection:Gelfand-Stacks}, $\bb{L}_{Y/X}\in \ob{D}_+(Y)$\footnote{Here, $\ob{D}_+(Y)$ is the value at $Y$ of the right Kan extension of $\ob{D}_+$ on affinoid Gelfand stacks to all presheaves on Gelfand rings.} such that for all maps $\eta\colon \GSpec (R)\to Y$ with $R$ a Gelfand ring, and all connective modules $M\in \ob{D}_{\geq 0}(M)$, the fiber 
\begin{equation}\label{eqTangentSpaceDefinitionPullback}
\begin{tikzcd}
F_f(\eta,R,M) = \ob{fib} (Y(R\oplus M) \to Y(R)\times_{X(R)} X(R\oplus M))
\end{tikzcd}
\end{equation}
 with respect to the point induced by $\eta\in Y(R)$ and the composite with the zero section $\GSpec (R\oplus M)\to \GSpec (R)\to X$, is corepresented by $\eta^*\bb{L}_{Y/X}$. That is, if there is a natural on $M$  equivalence of anima
 \[
 F_f(\eta,R,M)=\tau_{\geq 0} \Hom_{R}(\eta^*\bb{L}_{Y/X}, M). 
 \]
We say that $Y$ admits a cotangent complex if $Y\to \GSpec(\Q_p)$ admits a cotangent complex, in that case we write $\bb{L}_Y=\bb{L}_{Y/\GSpec(\Q_p)}$.

\item Suppose that $f$ admits a cotangent complex.
  We say that $f$ is \textit{formally \'etale} if $\bb{L}_f=0$, we say that $f$ is \textit{formally smooth} if $\bb{L}_f$ is the dual of a connective perfect complex, that is, a perfect complex of non-positive perfect amplitude.

\item Finally, we say that $f$ \textit{admits a deformation theory} if it is nilcomplete, infinitesimally cohesive, and admits a cotangent complex. 
 \end{enumerate}
 \end{definition}
 
 \begin{example}
   \label{sec:deform-theor-prop-example-admitting-a-deformation-theory}
   Using the results of \Cref{sec:cotang-compl-gelf-gelfand-rings} it is easy to see that each morphism $\GSpec(B)\to \GSpec(A)$ admits a deformation theory
 \end{example}

A natural source of ($\dagger$-)formally \'etale morphisms satisfying all deformation theoretic properties introduced above are de Rham stacks. 

\begin{lemma}\label{LemmadeRhamDefProperty}
  Let $X$ be a Gelfand stack.
  Then the de Rham stack $X^{\ob{dR}}$ is nilcomplete, infinitesimally cohesive, and admits a cotangent complex with $\bb{L}_{X^{\dR}/\GSpec(\Q_p)}=0$.
\end{lemma}
\begin{proof}
  By definition, $X^\dR$ is the functor $A\mapsto X^\diamond(\Marc(A^u))$ on Gelfand rings, where $A^u$ is the uniform completion.
  Since $A^u$ only depends on the reduction of the connected components $\pi_0(A)$ of $A$, one deduces that $X^{\dR}$ satisfies all properties of \Cref{sec:cotang-compl-morph-definition-having-a-cotangent-complex}.
\end{proof}

The properties of \Cref{sec:cotang-compl-morph-definition-having-a-cotangent-complex} satisfy the following standard stability properties.

\begin{lemma}\label{LemmaStabilityCotangentComplexes}
  Let $Z\xrightarrow{f} Y\xrightarrow{g} X$ be morphisms of presheaves on affinoid Gelfand stacks.
  The following hold.
\begin{enumerate}

\item  Suppose that $g$ is nilcomplete (resp.\ infinitesimally cohesive).
  Then $f$ is nilcomplete (resp.\ infinitesimally cohesive) if and only if $g\circ f$ is so.

\item Suppose that $g$ admits a cotangent complex. Then $f$ admits a cotangent complex if and only if $g\circ f$ does so, in which case we have a fiber sequence of sheaves on $Z$
\[
f^*\bb{L}_{Y/X}\to \bb{L}_{Z/X}\to \bb{L}_{Z/Y}.
\]

\item Suppose that $g$ is nilcomplete (resp.\ infinitesimally cohesive) and that $h\colon X\to  X'$ is a map of presheaves on Gelfand rings with base change $g'\colon Y'\to X'$ and $h'\colon Y\to Y$.
  Then $g'$ is nilcomplete (resp.\ infinitesimally cohesive).

\item  Keep the same notation of (3). Suppose that $g$ admits a cotangent complex, then $g'$ admits a cotangent complex and there is a natural equivalence 
\[
h^{\prime,*}\bb{L}_{g}=\bb{L}_{g'}.
\]

\end{enumerate}
\end{lemma}
\begin{proof}

\begin{enumerate}

\item  We give the proof for the nilcomplete  property, that of infinitesimal cohesivity being exactly the same.
  We have a commutative diagram (of presheaves on strictly totally disconnected affinoid Gelfand stacks) with cartesian squares
\[
\begin{tikzcd} 
Z \ar[d]& & \\
Z^{\ob{nil}/Y} \ar[r] \ar[d]& Y \ar[d] & \\
Z^{\ob{nil}/X}\ar[r] \ar[d]& Y^{\ob{nil}/X} \ar[d] \ar[r] & X\ar[d] \\
Z^{\ob{nil}} \ar[r] & Y^{\ob{nil}}\ar[r] & X^{\ob{nil}} 
\end{tikzcd}
\]
By hypothesis the map $Y\to Y^{\ob{nil}/X}$ is an equivalence, which makes the map $Z^{\ob{nil}/Y}\to Z^{\ob{nil}/X}$ an equivalence. From this we deduce that $Z$ is nilcomplete over $X$ if and only if it is so over $Y$.

\item  This is proven in the exact same way as (1), see \cite[Proposition 17.2.5.2]{Luriea}.

\item We give the proof for the nilcomplete  property, that of infinitesimal cohesivity being exactly the same. By assumption, we have that $Y\to Y^{\ob{nil}/X}$ is an equivalence. The claim then follows by noticing that $Y'=Y\times_X X'$ and 
\[
\begin{aligned}
Y^{\ob{nil}/X}\times_{X} X' &  = ( Y^{\ob{nil}}\times_{X^{\ob{nil}}} X^{\prime,\ob{nil}} )\times_{X^{\prime,\ob{nil}}} X^{'} \\ 
 & = Y^{\prime,\ob{nil}}  \times_{X^{\prime,\ob{nil}}} X' \\ 
 & = Y^{\prime,\ob{nil}/X'}. 
\end{aligned}
\]

\item  This is proven in the exact same way as (3), see \cite[Remark 17.2.4.6]{Luriea}. \qedhere

\end{enumerate}
\end{proof}

 In order to find a good balance between having a reasonable context for deformation theory, and staying close to Gelfand stacks, we will use the notions \Cref{sec:cotang-compl-morph-definition-having-a-cotangent-complex} for morphisms $f\colon Y\to X$ of sheaves for some topology $\tau$ on Gelfand rings (we recall that implicitly we work with qfd Gelfand rings throughout).\footnote{On a more technical level it is probably false that $Y^{\ob{nil}}$ and $Y^{\ob{coh}}$ are Gelfand stacks, i.e., invert $!$-equivalences if $Y$ is a Gelfand stack.}
 There are many possible choices for $\tau$.
 We make the running assumption that each covering $\{\GSpec(A_i)\to \GSpec(A)\}_{i\in I}$ in $\tau$ is a $!$-cover.\footnote{We also assume the convenient propertes that $I$ is finite and that finite disjoint unions of spectra of Gelfand rings are covers in $\tau$ to reduce to single-element coverings. This is only added for simplicity of exposition, and can easily be avoided.}
 In addition, we will use the following possible properties of $\tau$:
 \begin{enumerate}
 \item[(formsm)] Given a square-zero extension $\widetilde{A}\to A$ each covering $\{\GSpec(A_i)\to \GSpec(A)\}_{i\in I}$ lifts to a covering of $\GSpec(\widetilde{A})$.

 \item[(flat)] Each covering $\{\GSpec(A_i)\to \GSpec(A)\}_{i\in I}$ with $A$ strictly totally disconnected\footnote{This is convenient for the $\dagger$-quasi-pro-\'etale topology, but could be replaced by some other class of Gelfand rings as well. Consequently, one should restrict to this class as well in the definition of $Y^{\ob{nil}}$ in \Cref{sec:deform-theor-prop-example-admitting-a-deformation-theory}.} consists of morphisms $\GSpec(A_i)\to \GSpec(A)$ with $A_i$ strictly totally disconnected and $A\to A_i$ flat for the solid tensor product.
 \end{enumerate}

 Property (formsm) will allow us to study local properties of infinitesimal cohesity and the existence of cotangent complexes, while the property (flat) is needed only when discussing nilcompleteness.

Our framework for deformation theory in analytic geometry is the following category. 

\begin{definition}\label{DefinitionQuasiProetStacks}
  We define the category $\cat{GelfStk}_{\tau}$ of $\tau$-Gelfand stacks to be the category of sheaves for the topology $\tau$ on Gelfand rings. 
\end{definition}

\begin{remark}\label{RemarkDescent}
By the assumption that $\tau$-covers are $!$-covers, there is a left exact and colimit preserving functor of $\infty$-topoi
\[
\cat{GelfStk}_{\tau}\to \cat{GelfStk}
\]
which is the sheafification under $!$-equivalences. In particular, by passing to right adjoints we have a fully faithful functor $\cat{GelfStk}\subset \cat{GelfStk}_{\tau}$, which allows us to apply all considerations below to Gelfand stacks.
 We let $\ob{D}$ denote the induced functor of quasi-coherent sheaves on $\cat{GelfStk}_{\tau}$.
\end{remark}

The following lemma clarifies the sheaf properties of the constructions of \Cref{sec:cotang-compl-morph-definition-having-a-cotangent-complex}.

\begin{lemma}\label{LemmaQuasiproetTop}
Let $Y$ be a $\tau$-Gelfand stack. 

\begin{enumerate}

\item Assume that $\tau$ satisfies (flat).
  Then the functor $Y^{\ob{nil}}\colon \cat{GelfRing}^{\ob{std}}\to \cat{Ani}$ is a sheaf for the topology induced by $\tau$ on $\cat{GelfRing}^{\ob{std}}$.
\item Endow the category $\ob{Sqz}$  of square-zero extensions $\widetilde{A}\to A$ with the $\tau$-topology on $\widetilde{A}$.
  The functors $\iota(Y)$ and $Y^{\ob{coh}}$ from $\ob{Sqz}\to \cat{Ani}$ are $\tau$-Gelfand stacks.
\end{enumerate}
\end{lemma}
\begin{proof}
\begin{enumerate}

\item Let $A$ be a strictly totally disconnected  Gelfand ring and $A\to B$ a $\tau$-cover with \v{C}ech nerve $B^{\bullet}$.
  By assumption it is in particular a flat cover and we have 
\[
\tau_{\leq n} B^{\bullet}= B^{\bullet}\otimes_{A} \tau_{\leq n} A. 
\]
Then, $ \tau_{\leq n} A\to \tau_{\leq n} B$ is a $\tau$-cover with \v{C}ech nerve  $\tau_{\leq n} B^{\bullet}$ and the natural map
\[
Y(\tau_{\leq n} A)\to \ob{Tot}(Y(\tau_{\leq n} B^{\bullet}))
\]
is an equivalence. Taking limits as $n\to \infty$ we deduce that $Y^{\ob{nil}}$ is a sheaf for $\tau$. 
\item Let $\widetilde{A}\to A$ be a square-zero extension associated to a derivation $d\colon A\to M[1]$ with $M$ a connective $A$-module.
  Let $\widetilde{A}\to \widetilde{B}$ be a $\tau$-cover with base change $\widetilde{B}\to B:=A\otimes_{\widetilde{A}} \widetilde{B}$.
  Then $\widetilde{B}\to B$ is a square-zero extension written as a pullback 
\[
\widetilde{B} = B\times_{B\oplus B\otimes_A M[1]} B,
\]
where the left map is the trivial derivation and the right map the base change along $\widetilde{A}\to \widetilde{B}$ of the derivation $(\ob{id}_A,d)\colon A\to A\oplus M[1]$. Thus, the \v{C}ech nerve $\widetilde{B}^{\bullet}$ of $\widetilde{A}\to \widetilde{B}$  satisfies 
\[
\widetilde{B}^{\bullet} = B^{\bullet}\times_{B^{\bullet}\oplus B^{\bullet}\otimes_A M[1]} B^{\bullet}. 
\]
It follows that both $\iota(Y)$ and $Y^{\ob{coh}}$ are sheaves for the $\tau$-topology on $\widetilde{A}$. \qedhere
\end{enumerate}
\end{proof}

Given $\tau$ the existence of the cotangent complex can conveniently be checked on a restricted class of rings, in fact locally on a basis for $\tau$.
The following lemma discusses this.
We note that it does \emph{not} claim that one can check the existence of the cotangent complex on a $\tau$-cover on $Y$ (this needs more difficult arguments and property (formsm), cf.\ \Cref{LemmaDescentSourceKey}).
  
\begin{lemma}  \label{sec:cotang-compl-morph-weakening-of-definition}
  Let $f\colon Y\to X$ be a map of $\tau$-Gelfand stacks and $\s{C}\subset \cat{GelfRing}_{/Y}$ be a basis for $\tau$ over $Y$ in the sense that any Gelfand ring $A$ over $Y$ admits a $\tau$-cover $A\to B$ with $B\in \s{C}$, and such that the co-\v{C}ech nerve $B^{\bullet}$ of $A\to B$ lies in $\s{C}$.
  Suppose that $f$ admits a cotangent complex when restricted to $\s{C}$, that is, that the hypothesis of \Cref{sec:cotang-compl-morph-definition-having-a-cotangent-complex} holds when restricted to Gelfand rings in $\s{C}$.
  Then $f$ admits a cotangent complex as in \textit{loc.\ cit.}. 
\end{lemma}
\begin{proof}
By assumption, there is an object $\bb{L}_f\in \ob{D}_+(Y)$ satisfying the hypothesis of \Cref{sec:cotang-compl-morph-definition-having-a-cotangent-complex} when restricted to the basis $\s{C}$.  We want to see that $\bb{L}_f$ gives rise to a cotangent complex for $f$.

More precisely,  let $A$ be a Gelfand ring and $\eta \colon \GSpec(A)\to Y$ a map. Consider the functor 
\[
F_f(A,-,\eta)\colon \ob{D}_{\geq 0}(A)\to \ob{Ani}
\]
given by the fiber \eqref{eqTangentSpaceDefinitionPullback}.  We want to see that $F_f(A,-,\eta)$ is corepresented by $\eta^*\bb{L}_{f}$.  For that, let $A\to B$ be a cover with $B\in \s{C}$ and let $B^{\bullet}$ be its co-\v{C}ech nerve. 
 By \Cref{LemmaDescentlocallyAlmostConnective} we have an equivalence of categories 
\[
\ob{D}_+(A)\cong \ob{Tot}(\ob{D}_+(B^{\bullet}))
\]
because we assume that every $\tau$-cover is a $!$-cover.
In particular, we have an equivalence of functors 
\[
F_f(A,-,\eta)=\ob{Tot}(F_f(B^{\bullet},(-)\otimes_A B^{\bullet},\eta_{B^{\bullet}})),
\]
where 
\[
\eta_{B^{\bullet}}\colon \GSpec (B^{\bullet})\to \GSpec(A)\to Y
\]
 is the composite map.  Indeed, this follows from the definition of $F_f$ as a pullback so it commutes with totalizations,  and the  fact that if $A\to B$ is a $!$-cover with co-\v{C}ech nerve $B^{\bullet}$ then so is 
 $A\oplus M\to (A\oplus M)\otimes_A B= B\oplus M\otimes_A B
 $ with co-\v{C}ech nerve $B^{\bullet}\oplus M\otimes_A B^{\bullet}$.

By hypothesis, we have natural equivalences of cosimplicial diagrams of functors $\ob{D}_+(A)\to \cat{Ani}$
\[
\begin{aligned}
F_f(B^{\bullet},(-)\otimes_A B^{\bullet},\eta_{B^{\bullet}})  & = \tau_{\geq 0}\ob{Hom}_{B^{\bullet}}(\eta_{B^{\bullet}}^* \bb{L}_{f}, (-)\otimes B^{\bullet} ) \\ 
&  = \tau_{\geq 0}\ob{Hom}_{A}(\eta^* \bb{L}_{f}, (-)\otimes B^{\bullet} ).
\end{aligned}
\]
Passing to totalizations we deduce the natural equivalence
\[
F_f(A,-,\eta)=\tau_{\geq 0}\ob{Hom}_{A}(\eta^* \bb{L}_{f}, -)
\]
proving what we wanted. 
\end{proof}

The following lemma tells us that one can also check infinitesimal cohesivity  when testing on a restricted class of rings, in fact a basis for the topology.
 
\begin{lemma}\label{LemmadeRhamStacks}
Let $f\colon Y\to X$ be a map of $\tau$-Gelfand stacks.
Let $\s{C}\subset \ob{Sqz}$ be a basis for the $\tau$-topology on square-zero extensions in the sense that given $\widetilde{A}\to A$ in $\ob{Sqz}$ there exists a $\tau$-cover $\widetilde{A}\to \widetilde{B}$ with \v{C}ech nerve $\widetilde{B}^{\bullet}\to A\otimes_{\widetilde{A}} \widetilde{B}^{\bullet}$ in $\s{C}$.
Suppose that $f$ is infinitesimally cohesive when restricted to functors in $\s{C}$, that is, \Cref{sec:cotang-compl-morph-definition-having-a-cotangent-complex} (2)  holds as functors  $\s{C}\to \cat{Ani}$.
Then $f$ satisfies \Cref{sec:cotang-compl-morph-definition-having-a-cotangent-complex} (2) when tested on all Gelfand rings.
\end{lemma}
\begin{proof}
  The map $f\colon Y\to X$ is infinitesimally cohesive if and only if $Y\to Y^{\ob{coh}/X}$ is an equivalence of $\tau$-sheaves on $\cat{Sqz}^{\ob{std}}$ (using \Cref{LemmaQuasiproetTop}).
  For testing whether it is an equivalence, as they are both $\tau$-sheaves,  it suffices to check when evaluated on a basis for $\tau$.
  The lemma follows. 
\end{proof}

Next, we want to address  the descent properties for the deformation theoretic conditions of \Cref{sec:cotang-compl-morph-definition-having-a-cotangent-complex}.  We start with the easier case of descent on the target.

\begin{lemma}\label{DescentDefTheoryTarget}
  Let $f\colon Y\to X$ be a morphism of $\tau$-Gelfand stacks and let $X'\to X$ be an epimorphism of $\tau$-Gelfand stacks with base change $f'\colon Y':=X'\times_X Y\to X'$.
\begin{enumerate}
\item Assume that $\tau$ satisfies (flat) and that $f'$ is nilcomplete.
  Then $f$ is nilcomplete.
\item Assume that $f'$ is infinitesimally cohesive.
  Then $f$ is infinitesimally cohesive.
\item Assume that $f'$ admits a cotangent complex.
  Then $f$ admits a cotangent complex.
\end{enumerate}
\end{lemma}
\begin{proof}
First, we argue for nilcompleteness.
Suppose that $f\colon Y'\to X'$ is nilcomplete.
We have a commutative diagram of $\tau$-sheaves on $\cat{GelfStk}^{\ob{std}}$ with cartesian squares
\[
\begin{tikzcd}
Y' \ar[rr] \ar[d]& & Y \ar[d] \\
Y^{\prime, \ob{nil}/X'} \ar[r] \ar[d] & Y^{\prime,\ob{nil}/X} \ar[r]\ar[d] & Y^{\ob{nil}/X}\ar[d] \\
X'\ar[r] & X^{\prime,\ob{nil}/X}\ar[r] & X
\end{tikzcd}
\]
By assumption, the map $X'\to X$ is a $\tau$-cover.
Thus, by descent, the map $Y\to Y^{\ob{nil}/X}$ is an isomorphism if and only if $Y'\to Y^{'\ob{nil}/X'}$ is so, and this is precisely saying that $f'$ is nilcomplete if and only if $f$ is so.

Suppose now that $f'$ is infinitesimally cohesive, in this case we have a commutative diagram of cartesian squares of $\tau$-sheaves on $\cat{Sqz}$
\[
\begin{tikzcd}
\iota(Y') \ar[rr] \ar[d]& & \iota(Y) \ar[d] \\
Y^{\prime, \ob{coh}/X'} \ar[r] \ar[d] & Y^{\prime,\ob{coh}/X} \ar[r]\ar[d] & Y^{\ob{coh}/X}\ar[d] \\
\iota(X')\ar[r] & X^{\prime,\ob{coh}/X}\ar[r] & \iota(X)
\end{tikzcd}
\]
The fact that $X'\to X$ is a $\tau$-cover implies that $\iota(X')\to \iota(X)$ is as well.
Indeed, given a square-zero extension $\widetilde{A}\to A$ and $x\in X(\widetilde{A})$ we find a $\tau$-cover $\widetilde{A}\to B$ with a lift of $x$ to $X'$, but then we also lift $x\in \iota(X)(\widetilde{A}\to A)=X(\widetilde{A})$ to $X'$ over the square-zero extension $(B\to A\otimes_{\widetilde{A}}B)$.
The same argument as before shows that $f$ is infinitesimally cohesive. 
  
Next, we argue for the cotangent complex. Suppose that $f'\colon Y'\to X'$ admits a cotangent complex.
Let $\s{C}_{Y}\subset \cat{GelfRing}_{/Y}$ be the full subcategory of affinoid Gelfand rings over $Y$ that admit a morphism to $X'$ over $X$.
Then $\s{C}_{Y}$ is a basis for $\tau$ and by \Cref{sec:cotang-compl-morph-weakening-of-definition} it suffices to show that $f$ admits a cotangent complex when restricted to $\s{C}_{Y}$.
Thus, it suffices to show that given $\GSpec(R)\to X$ a map admitting a lift $\GSpec(R) \to X'$, the base change $Y\times_{X} \GSpec(R)\to \GSpec(R)$ admits a cotangent complex, namely, the same argument as in \cite[Proposition 17.2.4.7]{Luriea} will show that $f$ admits a cotangent complex when restricted to the basis $\s{C}_X$.
But now, we have that 
\[
Y\times_X \GSpec(R) = Y'\times_{X'} \GSpec(R)
\]
and $Y'\to X'$ admits a contangent complex by assumption.
Thus $Y\times_{X} \GSpec(R)\to \GSpec(R)$ admits a cotangent complex by \Cref{LemmaStabilityCotangentComplexes} (4) and the lemma follows.
\end{proof}

\begin{example}
  \label{sec:deform-theor-prop-general-berkovich-spaces}
  Assume that covers by $\dagger$-rational localizations are coverings in $\tau$.
  From \Cref{sec:cotang-compl-morph-weakening-of-definition} and \Cref{LemmadeRhamStacks} we can deduce that a general morphism $f\colon Y\to X$ of Berkovich spaces is infinitesimally cohesive and admits a cotangent complex.
  Indeed, by \Cref{sec:deform-theor-prop-general-berkovich-spaces} we can first reduce to the case that $X$ is affinoid, and then in \Cref{sec:cotang-compl-morph-weakening-of-definition} we can take $\mathcal{C}_Y$ as those affinoid Gelfand stacks $\GSpec(R)$ over $Y$, which factor over some affinoid $\dagger$-rational subspace in $Y$.  
\end{example}

\begin{example}
  \label{sec:deform-theor-prop-more-examples-for-infinitesimal-cohesivity}
  Assume that covers by $\dagger$-rational localizations are coverings in $\tau$.
  From \Cref{DescentDefTheoryTarget} and \Cref{sec:deform-theor-prop-general-berkovich-spaces} we can deduce that weakly representable morphisms $f\colon Y\to X$ of Gelfand stacks, i.e., there exists a cover $X'\to X$ with $X'$ a Berkovich space such that $Y\times_X X'$ is a Berkovich space, are infinitesimally cohesive and admit a cotangent complex.
\end{example}

Next, we want to address the descent properties on the source.
This is a more subtle question and requires additional assumptions.

\begin{lemma}\label{LemmaDescentSourceKey}
Let $Z\xrightarrow{f} Y \xrightarrow{g} X$ be a diagram of $\tau$-Gelfand stacks. The following hold.
\begin{enumerate}

\item Assume that $\tau$ satisfies (flat).
  Suppose that $f$ and $g\circ f$ are nilcomplete and that the natural map $Z^{\ob{nil}/X}\to Y^{\ob{nil}/X}$ is an epimorphism of $\tau$-sheaves on $\cat{GelfStk}^{\ob{std}}$.
  Then $g$ is nilcomplete.  

\item  Suppose that $f$ and $g\circ f$ are infinitesimally cohesive and that the map of functors $Z^{\ob{coh}/X}\to Y^{\ob{coh}/X}$ is an epimorphism for the $\tau$-topology on square-zero extensions.
  Then $g$ is infinitesimally cohesive. 

\item  Suppose that $f$ and $g\circ f$ admit a cotangent complex, $g$ is infinitesimally cohesive and $f$ is an epimorphism.
  Then $g$ admits a cotangent complex. 

\end{enumerate}
\end{lemma}
\begin{proof}
\begin{enumerate}

\item We have a commutative diagram of $\tau$-sheaves on $\cat{GelfStk}^{\ob{std}}$ (by \Cref{LemmaQuasiproetTop}) with cartesian squares
\[
\begin{tikzcd}
Z \ar[d] &  & \\
Z^{\ob{nil}/Y} \ar[r] \ar[d] & Y \ar[d]  & \\ 
Z^{\ob{nil}/X} \ar[r] \ar[d] & Y^{\ob{nil}/X} \ar[d] \ar[r] & X \ar[d] \\
Y^{\ob{nil}} \ar[r] & Y^{\ob{nil}} \ar[r] & X^{\ob{nil}}.  
\end{tikzcd}
\]
By assumption, the map $Z^{\ob{nil}/X}\to Y^{\ob{nil}/X}$ is an epimorphism. Thus, to show that $Y\to Y^{\ob{nil}/X}$ is an equivalence (i.e.\ $g$ is nilcomplete) it suffices to show that $Z^{\ob{nil}/Y}\to Z^{\ob{nil}/X}$ is an equivalence. By assumption $Z$ is both nilcomplete over $X$ and $Y$, this implies that all the arrows
\[
Z\to Z^{\ob{nil}/Y}\to Z^{\ob{nil}/X}
\]
are equivalences, proving what we wanted. 

\item This is proven in the exact same way as part (1).

\item Let $\s{C}_{Y}\subset \cat{Sch}^{\ob{aff}}_{/Y}$  be the basis consisting of those maps $\GSpec(R)\to Y$ that admit a lift to $Z$. By \Cref{sec:cotang-compl-morph-weakening-of-definition}  it suffices to show that $Y\to X$ admits a cotangent complex when restricted to $\s{C}_{Y}$.
  Let $\eta\colon  \GSpec(R)\to Y$ be a fixed map with composition $\nu\colon \GSpec(R)\to X$. Consider the functor 
\[
F_g(\eta,R,-)\colon \ob{D}_{\geq 0}(R)\to \cat{Ani}
\]
sending $M$ to the fiber $\ob{fib}(Y(R\oplus M)\to Y(R)\times_{X(R)}X(R\oplus M))$ of  \Cref{sec:cotang-compl-morph-definition-having-a-cotangent-complex} (3). Since $g$ is cohesive, and $R\oplus M= R\times_{R\oplus M[1]} R$ where the maps are   the zero derivations $R\to M[1]$, we have a pullback diagram of fibers over $\nu$
\[
\begin{tikzcd}
Y(R\oplus M)_{/\nu} \ar[r] \ar[d] & Y(R)_{/\nu}  \ar[d] \\ 
Y(R)_{/\nu} \ar[r] & Y(R\oplus M[1])_{/\nu}.
\end{tikzcd}
\]
 The previous implies that $F_{g}(\eta,R,-)$ sends pushouts in $\ob{D}_{\geq 0}(R)$ to pullbacks in anima, and therefore that it defines a spectrum object. In particular, we have that $F_g(\eta,R,M)\to \Omega F_{g}(\eta,R, M[1]) $ is an equivalence. Let $\tau\colon \GSpec(R)\to Z$ be a lift of $\eta$ (which exists as $\eta$ is in $\s{C}_{/Y}$). We have a fiber sequence of functors $\ob{D}_{\geq 0}(R)\to \cat{Ani}$
 \begin{equation}\label{eq003ko3smasasf}
 F_{f}(\tau, R, -)\to F_{g\circ f}(\tau,R,-)\to F_{g}(\eta,R,-). 
 \end{equation}
Since $f$ and $g\circ f$ admit a cotangent complex, there are natural equivalences for $M\in \ob{D}_{\geq 0}(R)$
\begin{equation}\label{eq0ko3klamslassssss}
 F_{f}(\tau, R, M)= \tau_{\geq 0} \ob{Hom}_{R}(\tau^* \bb{L}_{f},M) \mbox{ and } F_{g\circ f}(\tau, R, M)= \tau_{\geq 0} \ob{Hom}_{R}(\tau^* \bb{L}_{g\circ f},M)
\end{equation}
In particular, they form spectrum objects from $\ob{D}_{\geq 0}(R)\to \cat{Ani}$, i.e., upgrade naturally to functors to spectra (in fact, these are given by $\ob{Hom}_{R}(\tau^* \bb{L}_{f},-)$ resp.\ $\ob{Hom}_{R}(\tau^* \bb{L}_{g\circ f},-)$).  

 It follows that \eqref{eq003ko3smasasf} is a fiber-cofiber sequence of spectra-valued functors, and we deduce that 
\[
F_{g}(\eta,R,-)=\tau_{\geq 0} \Hom_{R}(\tau^*\ob{fib}(\bb{L}_{g\circ f}\to \bb{L}_f), M)
\] 
is corepresented by an object in $\ob{D}_{+}(R)$ compatible with base change along maps $\Spec R'\to \Spec R$ in $\s{C}_{Y}$.
We deduce from \Cref{sec:cotang-compl-morph-weakening-of-definition} that $g$ admits a cotangent complex as wanted. \qedhere
\end{enumerate}
\end{proof}

The additional assumptions in \Cref{LemmaDescentSourceKey} can be guaranteed under formally smoothness hypothesis on the cover of the source, and by assuming the property (formsm) on $\tau$.

\begin{proposition}\label{PropFinalDescentDeformationTheory}
Let $Z\xrightarrow{f} Y \xrightarrow{g} X$ be a diagram of $\tau$-Gelfand stacks such that $f$ is an epimorphism, infinitesimally cohesive and formally smooth. The following hold: 
\begin{enumerate}
\item Assume that $\tau$ satisfies (flat) and that each affinoid Gelfand stacks admits a $\tau$-cover by a strictly totally disconnected affinoid Gelfand stack.
  The map $Z^{\ob{nil}}\to Y^{\ob{nil}}$ is an epimorphism.
  In particular, if $f\circ g$ and $f$ are nilcomplete, then $g$ is nilcomplete. 
\item Assume that $\tau$ satisfies (formsm).
  The map $Y^{\ob{coh}}\to X^{\ob{coh}}$ is an epimorphism for the $\tau$-topology of functors $\cat{Sqz}\to \cat{Ani}$.
  In particular, if $g\circ f$ is infinitesimally cohesive then so is $g$. 
\item Assume that $\tau$ satisfies (formsm).
  Suppose that $g\circ f$ is infinitesimally cohesive and admits a cotangent complex, then $g$ admits a cotangent complex. 
\end{enumerate}
\end{proposition}
\begin{proof}
\begin{enumerate}
\item Let $\GSpec(R)\to Y^{\ob{nil}}$ be a map with $R$ strictly totally disconnected given by compatible maps $\eta_n\in Y(\tau_{\leq n}R)$. Since $Z\to Y$ is an epimorphism in $\tau$, we can find an $\tau$-cover $R\to R'$ and a lift $\GSpec(\pi_0(R'))\to Z$ of $\GSpec \pi_0(R)\to Y$.
  By exchanging $\eta_n$ by the composite $\GSpec \tau_{\leq n}R'\to \GSpec \tau_{\leq n} R\to Y$, we can assume without loss of generality that $R=R'$ and that we are given with a lift $\nu_0\in Z(\pi_0(R))$ of $\eta_0\in Y(\pi_0(R))$

  We claim that there are compatible elements $\nu_{n}\in Z(\tau_{\leq n} R)$ lifting $\eta_{n}$.
  Indeed, suppose we have such a lift for $n\geq 0$, we want to construct the lift for $n+1$. The map $\tau_{\leq n+1} R\to \tau_{\leq n} R$ is a square-zero extension associated to a derivation $\alpha\colon \tau_{\leq n} R \to \tau_{\leq n} R \oplus \pi_{n+1}(R) [n+2]$. The infinitesimal  cohesivity  of $f$ yields a cartesian square
\[
\begin{tikzcd}
Z(\tau_{\leq n+1}(R))_{/\eta_{n+1}} \ar[r] \ar[d] & Z(\tau_{\leq n} R)_{/\eta_{n}} \ar[d] \\ 
 Z(\tau_{\leq n} R)_{/\eta_{n}} \ar[r] &  Z(\tau_{\leq n} R\oplus  \pi_{n+1}(R) [n+2])_{/\eta_{n}},
\end{tikzcd}
\]
where to morphisms from $Z(\tau_{\leq n})_{\eta_n}$ are induced by the $0$-derivation resp.\ $\alpha$.
Thus, the obstruction to find a lift of $\nu_{n}\in Z(\tau_{\leq n} R)_{/\eta_{n}}$ to $\GSpec( \tau_{\leq n+1} R)$ over $\eta_{n+1}\colon \GSpec ( \tau_{\leq n+1} R)\to Y$ is given by an element in 
\[
\pi_0(\ob{Hom}(\nu^*_n \bb{L}_{Z/Y}, \pi_{n+1}(R) [n+2])=0
\]
using the formal smoothness of $f$ (which means that $\bb{L}_{Z/Y}$ is the dual of a connective perfect complex).
Hence the obstruction vanishes, and we can always find such lift $\nu_{n+1} \in Z(\tau_{\leq n+1}R)_{/\eta_{n+1}}$. This proves the first claim. The second claim about the nilcompleteness of $g$ follows from \Cref{LemmaDescentSourceKey} (1).

\item

 We follow the argument of \cite[Proposition 7.10 (1)]{MR4176782}.  Let $\widetilde{R}\to R$ be a square-zero extension of $R$ by $M$ with derivation $d\colon R\to M[1]$ and let $x\in Y(R)\times_{Y(R\oplus M[1])} Y(R)$  an object with first projection $x_1$. We want to show that, locally in the $\tau$-topology on $\cat{Sqz}$, the point $x$ admits a preimage along the map 
\[
Z(R)\times_{Z(R\oplus M[1])} Z(R) \xrightarrow{f}  Y(R)\times_{Y(R\oplus M[1])} Y(R). 
\]
By localizing on $R$ for the $\tau$-topology, we can assume that there exists a covering $h\colon R\to R'$ such that the point $h^\ast(x_1)\in Y(R')$ admits a lift $y_1\in Z(R')$, which we denote $\eta\colon \GSpec(R')\to Z$.
As $\tau$ satisfies (formsm) the covering $R\to R'$ lifts  to a $\tau$-cover $\widetilde{R}\to \widetilde{R'}$.
Then we may replace the square-zero extension $\widetilde{R}\to R$ by $\widetilde{R'}\to R'$ and assume $R'=R$.
Consider now the following commutative square 
\[
\begin{tikzcd}
  \ob{fib}(g) \ar[r] \ar[d]& \ob{fib}(p) \ar[r,"g"]  \ar[d]& \ob{fib}(q) \ar[d] \\
 \ob{fib}(f) \ar[r] & Z(R)\times_{Z(R\oplus M[1])} Z(R)  \ar[r,"f"]   \ar[d,"p"] & Y(R)\times_{Y(R\oplus M[1])} Y(R) \ar[d,"q"] \\ 
 & Z(R) \ar[r] & Y(R), 
\end{tikzcd}
\]
where the fiber of $f$ and $g$ is at the point $x$, the fiber of $p$ at $y_1$ and the fiber of $q$ at $x_1$.  It suffices to show that $\ob{fib}(g)$ is non empty.

Then, in the terminology of \cite[Definition 7.5]{MR4176782}, the map $g$ is the morphism  of anima of paths  from  $0 \rightsquigarrow d$
\[
\Omega_{0,d} \ob{Der}_{Z}(R, M[1]) \xrightarrow{g} \Omega_{0,d} \ob{Der}_{Y}(R, M[1]),
\]
where for a functor $F\colon \cat{GelfRing}\to \cat{Ani}$ and a map $\GSpec(R)\to F$ we set $\ob{Der}_F(R,M):=\ob{fib}(F(R\oplus M)\to F(R))$ (and define similarly the relative versions). Hence, we find that 
\[
\ob{fib}(g)=\Omega_{0,d} \ob{Der}_{Z/Y}(R, M[1])=\Omega_{0,d} \tau_{\geq 0}\Hom_{R}(\eta^*\bb{L}_{Z/Y}, M[1]). 
\]
But, as $Z\to Y$ is  formally smooth,  $\tau_{\geq 0}\Hom_{R}(\eta^*\bb{L}_{Z/Y}, M[1])$ is a connective anima, so that $0$ and $d$ are in the same connected component, and $\ob{fib}(g)$ is non-empty as wanted.

\item  By part (2), the map $Z^{\ob{coh}}\to Y^{\ob{coh}}$ is an epimorphism. Since both $f$ and $g\circ f$ are infinitesimally cohesive, \Cref{LemmaDescentSourceKey} (3) implies that $g$ is infinitesimally cohesive. Then, part (4) of \textit{loc.\ cit.} implies that $g$ admits a cotangent complex. \qedhere
\end{enumerate}
\end{proof}

Next, we discuss several choices for $\tau$.

\begin{example}
  \label{sec:deform-theor-prop-shriek-formally-etale}
  \begin{enumerate}
  \item One can consider the $!$-formally \'etale topology, where a $\tau$-covering $\{\GSpec(A_i)\to \GSpec(A)\}_{i\in I}$ is given by a $!$-cover by formally \'etale morphisms $\GSpec(A_i)\to \GSpec(A)$.
    This topology satisfies (formsm) as formally \'etale morphisms lift (uniquely) along square-zero extensions.
  \item One can consider just the $!$-topology on Gelfand rings.
    This topology satisfies neither (flat) nor (formsm).
  \item One can consider the $\dagger$-formally \'etale topology, where a $\tau$-covering $\{\GSpec(A_i)\to \GSpec(A)\}_{i\in I}$ is given by $\dagger$-formally \'etale maps (\Cref{sec:prop-morph-gelf-1-equivalent-characterizations-dagger-formally-etale-morphisms}) such that $\{\Marc(A_i)\to \Marc(A)\}_{i\in I}$ is an arc-cover.
    Note that $\dagger$-formally \'etale covers  are $!$-covers by \Cref{sec:prop-morph-gelf-1-weakly-formally-smooth-is-epimorphism}.
    This topology satisfies (formsm), as $\dagger$-formally \'etale maps are formally \'etale, and has the advantage that (geometric realizations of) $\dagger$-formally \'etale hypercovers are $!$-equivalences (by the argument of \Cref{sec:prop-morph-gelf-1-weakly-formally-smooth-is-epimorphism} as arc-hypercovers realize to $!$-equivalences).
  \end{enumerate}
\end{example}

An example for $\tau$ satisfying (formsm) and (flat) is given by the $\dagger$-quasi-pro-\'etale topology.

\begin{definition}\label{DefinitionGelfandRingsDaggerSmooth}
  A map $A\to B$ of Gelfand rings is said to be \textit{$\dagger$-quasi-pro-\'etale} if the map $\Marc(B)\to \Marc(A)$ is quasipro-\'etale and $\GSpec(B)\to \GSpec(A)$ is $\dagger$-formally \'etale.
  A \textit{$\dagger$-quasipro-\'etale cover} of $A$ is a finite collection of $\dagger$-quasipro-\'etale maps $\{A\to B_i\}_i$ such that $\{\n{M}(B_i)\to \n{M}(A)\}$ is jointly surjective. 
  \end{definition}

\begin{remark}\label{RemarkSTDCase}
\begin{enumerate}
\item Let $A\to B$ be a map of strictly totally disconnected Gelfand rings such that $\Marc(B)\to \Marc(A)$ is quasipro-\'etale. Then, we have a pullback diagram of arc-stacks $\Marc(B)=\n{M}(B)_{\Betti}\times_{\n{M}(A)_{\Betti}} \Marc(A)$ (cf.\ \cite[Lemma 4.5]{scholze2025berkovich}). Thus, in this situation, \Cref{DefinitionGelfandRingsDaggerSmooth} is equivalent to 
\[
\GSpec(B) = \n{M}(B)_{\Betti}\times_{\n{M}(A)_{\Betti}} \GSpec(A). 
\]
\item Let $f\colon A\to B$ be a $\dagger$-rational localization, then we have a commutative diagram of Gelfand stacks
\[
\begin{tikzcd}
\GSpec(B) \ar[r] \ar[d]& \GSpec(A) \ar[d]\\ 
\Marc(B)^{\dR} \ar[r] \ar[d]& \Marc(A)^{\dR} \ar[d] \\ 
\n{M}(B)_{\Betti} \ar[r] & \n{M}(A)_{\Betti}.
\end{tikzcd}
\]
In particular, $f$ is $\dagger$-quasipro-\'etale.

\item By the previous point, a map $A\to B$ is $\dagger$-quasipro-\'etale if and only if it is so in a $\dagger$-rational cover of $\GSpec(A)$. Moreover, \cite[Lemma 4.7.9]{drff_stacks} implies that standard Berkovich-\'etale maps are $\dagger$-quasipro-\'etale. We deduce that any Berkovich-\'etale map is $\dagger$-quasipro-\'etale. 

\end{enumerate}

\end{remark}

$\dagger$-formally \'etale and $\dagger$-quasipro-\'etale maps satisfy the following  properties.
In particular, Gelfand rings admit $\dagger$-quasi-pro-\'etale covers by strictly totally disconnected Gelfand rings (recall that we implicitly only consider qfd Gelfand rings).

\begin{lemma}\label{LemmaDaggerQuasiProet}
The following hold.

\begin{enumerate}

\item  Let $A\xrightarrow{f} B\xrightarrow{g} C$ be morphisms of  Gelfand rings and suppose that $f$ is $\dagger$-formally \'etale. Then $g$ is $\dagger$-formally \'etale if and only if $g\circ f$ is so. 

\item Let $A\to B$ be a $\dagger$-formally \'etale morphism of  Gelfand rings, and let $A\to A'$ be any morphism with base change $B\to B'$. Then $A'\to B'$ is $\dagger$-formally \'etale. Conversely, if $A'\to B'$ is $\dagger$-quasiproetale and $A\to A'$ is a $!$-cover, then $A\to B$ is $\dagger$-formally \'etale. 

\item Let $\{A_i\to B_i\}$ be a filtered diagram of $\dagger$-formally \'etale maps with colimit $f\colon A\to B$, then $f$ is $\dagger$-formally \'etale.

\item Let $A$ be a Gelfand ring. Then there is a $\dagger$-quasipro-\'etale (in particular, $\dagger$-formally \'etale) cover $A\to B$ with $B$ a strictly totally disconnected Gelfand ring. In particular, the \v{C}ech nerve $A\to B^{\bullet}$ consists of strictly totally disconnected Gelfand rings. 

\item Let $A\to B$ be a morphism of Gelfand rings and let $B\to C$ be a  $\dagger$-formally \'etale cover such that $A\to C$ is $\dagger$-formally \'etale, then $A\to B$ is $\dagger$-formally \'etale.

\item  Let $A\to B^{\bullet}$ be a $\dagger$-formally \'etale hypercover of a Gelfand ring $A$, then the hypercover of Gelfand stacks
\[
\GSpec (B^{\bullet})\to \GSpec(A)
\]
is a $!$-cover.

\end{enumerate}
The same statements are true if $\dagger$-formally \'etale is replaced by $\dagger$-quasipro-\'etale throughout.
\end{lemma}
\begin{proof}

\begin{enumerate}

\item  We have a commutative diagram with cartesian squares
\[
\begin{tikzcd}
\GSpec(C) \ar[d] &  & \\ 
\Marc(C)^{\dR}\times_{\Marc(B)^{\dR}}\GSpec(B) \ar[r]\ar[d] & \GSpec(B)  \ar[d,"\alpha"] & \\ 
 \Marc(C)^{\dR}\times_{\Marc(A)^{\dR}} \GSpec(A) \ar[r] \ar[d] &  \Marc(B)^{\dR}\times_{\Marc(A)^{\dR}} \GSpec (A)  \ar[r] \ar[d]& \GSpec (A) \ar[d]\\
\Marc(C)^{\dR} \ar[r,"\Marc(g)"]& \Marc(B)^{\dR}\ar[r,"\Marc(f)"]  & \Marc(A)^{\dR}.
\end{tikzcd}
\]
As $\GSpec(B)\to \GSpec(A)$ is $\dagger$-formally \'etale, the morphism $\alpha$ is an isomorphism.
 We deduce that the  map 
\[
\Marc(C)^{\dR}\times_{\Marc(B)^{\dR}}\GSpec(B)\to  \Marc(C)^{\dR}\times_{\Marc(A)^{\dR}} \GSpec(A)
\]
is an isomorphism, and therefore that $g$ is $\dagger$-quasipro-\'etale if and only if $f\circ g$ is so.
Let us now finish the claim for the additional claim about $\dagger$-quasipro-\'etaleness. By \cite[Definition 4.3]{scholze2025berkovich} one deduces that $\Marc(g)$ is quasipro-\'etale if and only if $\Marc(f\circ g) $ is so, namely, being quasi-pro-\'etale is equivalent to the fact that the completion of the fibers of the maps at geometric points $\Marc(C)$ are given by the spectrum of $C(S,C)$ with $S$ a profinite set. 

\item Base change along $\GSpec(A')\to \GSpec(A)$  yields  $\GSpec (B')=\GSpec(B)\times_{\GSpec(A) }\GSpec (A')$ and 
\[
\begin{gathered}
 (\Marc(B)^{\dR}\times_{\Marc(A)^{\dR}} \GSpec(A))\times_{\GSpec(A)} \GSpec (A')  \\  =  \Marc(B)^{\dR}\times_{\Marc(A)^{\dR}} \GSpec (A') \\  
  =  (\Marc(B)^{\dR}\times_{\Marc(A)^{\dR}} \Marc(A')^{\dR})\times_{\Marc(A')^{\dR}} \GSpec (A') \\ 
  = \Marc(B')^{\dR}\times_{\Marc(A')^{\dR}} \GSpec (A').
 \end{gathered}
\]
Now, if $A\to B$ is $\dagger$-formally \'etale (resp.\ $\dagger$-quasipro-\'etale) the previous computations shows that $A'\to B'$ is also $\dagger$-formally \'etale (resp.\ $\dagger$-quasipro-\'etale).
Conversely, if $A'\to B'$ is $\dagger$-formally \'etale (resp.\ $\dagger$-quasipro-\'etale) and $\GSpec (A')\to \GSpec(A)$ is a $!$-cover, by descent we deduce that $f\colon A\to B$ is $\dagger$-formally \'etale (resp.\ $\dagger$-quasipro-\'etale).

\item  Let $\{A_i\to B_i\}$ is a filtered  system of $\dagger$-formally \'etale (resp.\ $\dagger$-quasipro-\'etale) maps  with colimit $A\to B$. Then, we have as Gelfand stacks
\[
\Marc(A)^{\dR}=\varprojlim_i \Marc(A_i)^{\dR} \mbox{ and } \GSpec(A)= \varprojlim_i \GSpec(A_i)
\]
and similarly for $B$.  Indeed, the limit property holds for arc stacks and Gelfand stacks,  and the de Rham stack functor commutes with limit being a right adjoint. Since the map $\GSpec (B_i)\to \Marc(B_i)^{\dR}\times_{\Marc(A_i)^{\dR}} \GSpec (A_i)$ is an isomorphism, taking limits we deduce that $A\to B$ is $\dagger$-formally \'etale (resp.\ $\dagger$-quasipro-\'etale, using that cofiltered inverse limits of quasi-pro-\'etale maps of arc-stacks remain quasi-pro-\'etale).

\item  To construct this cover we use the descent results of \cite[Section 4.6]{drff_stacks}, particularly of Proposition 4.6.4 of \textit{loc.\ cit.}. First, consider the map $\GSpec(A)\to \n{M}(A)_{\Betti}$ to the Betti stack (which exists as $A$ is qfd, see \cite[Proposition 3.2.10]{drff_stacks}) and let $S\to  \n{M}(A)_{\Betti}$ be a surjection from a light profinite set. The map $S_{\Betti}\to  \n{M}(A)_{\Betti}$ is descendable by \cite[Lemma 3.2.4]{drff_stacks}. Then the pullback $S_{\Betti}\times_{\n{M}(A)_{\Betti}} \GSpec(A)$ is represented by a Gelfand ring $B$ which is $\dagger$-quasi-pro-\'etale over $A$ and such that $\n{M}(B)=S$. Hence, we can assume without loss of generality that $\n{M}(A)$ is profinite, i.e.\ that $A$ is totally disconnected. Then, by \cite[Proposition 4.6.4 (4)]{drff_stacks}, we can find a $!$-cover $A\to B=\varinjlim_i B_i$ with each $A\to B_i$ a finite \'etale map. By \Cref{RemarkSTDCase} (3)  finite \'etale maps are $\dagger$-quasipro-\'etale, and the previous argument yields a $\dagger$-quasipro-\'etale cover by part (3) with $B$ strictly totally disconnected as wanted.

\item Consider the diagram of (1). Since $B\to C$ is a $\dagger$-formally \'etale (resp.\ $\dagger$-quasipro-\'etale) cover the map $\Marc(B)\to \Marc(C)$ is a quasipro-\'etale cover of arc-stacks, and if $\Marc(C)\to \Marc(A)$ is quasipro-\'etale  then so is $\Marc(B)\to \Marc(A)$.
  Passing to de Rham stacks, to see that $\GSpec(B)\to \Marc(B)^{\dR}\times_{\Marc(A)^{\dR}} \GSpec(A)$ is an isomorphism (i.e.\ that $A\to B$ is $\dagger$-formally \'etale), it suffices to prove it after base change along $\Marc(C)^{\dR}\to \Marc(B)^{\dR}$, in which case we claim  that
\[
\Marc(C)^{\dR}\times_{\Marc(B)^{\dR}} \GSpec(B)\to \Marc(C)^{\dR}\times_{\Marc(A)^{\dR}} \GSpec(A)
\]
is an isomorphism. But $B\to C$ and $A\to C$ are $\dagger$-formally \'etale, so that the previous stacks are naturally isomorphic to $\GSpec C$, proving what we wanted.

\item  This follows from the fact that we have a cartesian diagram
\[
\GSpec (B^{\bullet})=\Marc(B^{\bullet})^{\dR}\times_{\Marc(A)^{\dR}} \GSpec(A),
\]
and that $|\Marc(B^{\bullet})|\to \Marc(A)^{\dR}$ is an  equivalence of Gelfand stacks by \cite[Theorem 5.6.1]{drff_stacks}. \qedhere
\end{enumerate}

\end{proof}

The major advantage of the $\dagger$-quasi-pro-\'etale topology is the following flatness (implying property (flat) for the $\dagger$-quasi-pro-\'etale topology).

\begin{lemma}
  \label{sec:deform-theor-prop-flatness-for-dagger-quasi-pro-etale-maps}
  Let $A$ be a strictly totally disconnected Gelfand ring, then any $\dagger$-quasipro-\'etale map $A\to B$ can be written as $B=\varinjlim_i B_i$ where $B_i$ is a finite \'etale $A$-algebra.
  In particular, it is flat for the solid tensor product.   
\end{lemma}
\begin{proof} By  \Cref{RemarkSTDCase} (1)  we have a pullback diagram $\GSpec(B)=\GSpec(A) \times_{\n{M}(A)_{\Betti}} \n{M}(B)_{\Betti}$. The claim follows from the fact that $\n{M}(B)_{\Betti}\to \n{M}(A)_{\Betti}$ is a map of light profinite sets and hence can be written as a countable limit of finite disjoint unions of clopen subspaces of $\n{M}(A)_{\Betti}$.   
\end{proof}

\begin{remark}
  \label{sec:deform-theor-prop-quasi-proetale-version}
  Assume that $\tau$ is the $\dagger$-quasi-pro-\'etale topology.
  The functor $\ob{D}$ on $\tau$-Gelfand stacks from \Cref{RemarkDescent} has a well defined $t$-structure (by \Cref{LemmaDaggerQuasiProet} (5)) such that for $X$  a  quasi-pro-\'etale hypersheaf and $M\in \ob{D}(X)$, the sheaf $M$ is connective if and only if the pullback to all strictly totally disconnected affinoid Gelfand stacks $\GSpec(R)\to X$ is connective.
  Indeed, in this setup the pullback functors for morphisms of strictly totally disconnected affinoid Gelfand stacks is flat.
\end{remark}

\subsection{Tangent and cotangent spaces of Gelfand stacks}\label{ss:GelfandStacksDeformationTheory}

In this section we shall discuss some examples of morphisms of Gelfand stacks that admit a (co)tangent complex. Thanks to the descent of the deformation property conditions of \Cref{DescentDefTheoryTarget} and \Cref{PropFinalDescentDeformationTheory}, we shall focus ourselves to $\tau$-sheaves (for a topology as in \Cref{s:DefTheoryProperties}), even though all  the relevant  examples are already honest Gelfand stacks.
First, we begin with the definition of the tangent space.

\begin{construction}[Tangent space]\label{ConstructionTangentSpaces}
  Let $f\colon Y\to X$ be a morphism of $\tau$-Gelfand stacks which is infinitesimally cohesive and admits a cotangent complex.
  Let $\cat{Mod}_{\geq 0, /Y}\to \cat{GelfRing}_{/Y}$ be the cocartesian fibration over the category of Gelfand rings over $Y$ associated to the functor $\ob{D}_{\geq 0}\colon \cat{GelfRing}_{/Y}  \to \cat{Cat}_{\infty}$. The \textit{tangent space of $f$} is the functor
\[
\underline{\bb{T}}_{f}\colon \cat{Mod}_{\geq 0, /Y}\to \cat{Ani}
\]
sending a triple $(R,\eta,M)\in \cat{Mod}_{\geq 0, /Y}$ with $M\in \ob{D}_{\geq (R)}$ and $\eta\colon \GSpec(R)\to Y$ to the fiber
\[
\underline{\bb{T}}_f(\eta,R,M):=\ob{fib}(Y(R\oplus M)\to Y(R)\times_{X(R)} X(R\oplus M) ),
\]
where the fiber is taken at the point $x\in Y(R)\times_{X(R)} X(R\oplus M) $ defined by the map $\eta$ and the map $\GSpec (R\oplus M)\to \GSpec (R)\to X$ induced by the trivial derivation $R\to R\oplus M$.\footnote{To be completely formal, one first needs to construct a functor $\cat{Mod}_{\geq 0}\to \cat{Sqz}$ from the cocartesian fibration associated to $\ob{D}_{\geq 0}(-)\colon \cat{GelfRing}\to \cat{Cat}_{\infty}$ to the category of square-zero extensions of Gelfand rings, this is done in the exact same way as for discrete animated rings, eg.\ using the deformation theory  formalism of \cite[Chapters 7.3 and 7.4]{lurie_higher_algebra}.}  

Since $f$  is infinitesimally cohesive, $\bb{T}_f$ is a spectrum  functor $\cat{Mod}_{\geq 0,/Y}\to \cat{Ani}$, that is, for $(\eta,R,M)$ as before the natural map 
\[
\underline{\bb{T}}_f(\eta,R,M)\to \Omega\underline{\bb{T}}_f(\eta,R,M[1])
\]
is an isomorphism. Furthermore, the fact that $f$ admits a cotangent complex means precisely that $\underline{\bb{T}}_f$ is corepresented by the locally eventually connective $\bb{L}_{f}$, that is, there are natural equivalences of anima
\[
\underline{\bb{T}}_f(\eta,R,M)=\tau_{\geq 0} \Hom_R(\eta^*\bb{L}_{f}, M). 
\]

We want to consider the tangent complex as a $\tau$-Gelfand stack, as well as all its possible shifts.
For that, given $i\geq 0$, we define the \textit{$i$-th shifted geometric tangent space} to be the $\tau$-Gelfand stack
\[
\bb{T}_f[i]\colon \cat{GelfRing}_{/Y}\to  \cat{Ani}
\]
sending $\eta\colon \GSpec(R)\to Y$ to 
\[
\bb{T}_f[i](\GSpec(R)\to Y):= \underline{\bb{T}}_f(\eta,R,R[i]). 
\]
We note that, since $f$ admits a cotangent complex, we have that 
\[
\bb{T}_f[i](\GSpec(R)\to Y)=  \tau_{\geq -i} \Hom_{R}(\eta^*\bb{L}_{f}, R). 
\] 
In particular, since quasi-coherent sheaves satisfy descent along $!$-equivalences, if $f\colon Y\to X$ is a morphism of Gelfand stacks the shifted tangent spaces $\bb{T}_f[i]$ are also Gelfand stacks over $Y$. To keep the definitions consistent, for $i\leq 0$ we define the $i$-th  shift of the geometric cotangent complex $\bb{T}_f[i]$ to be the $(-i)$-th iterated loop $\Omega^{-i} \bb{T}_f$. 
\end{construction}

\begin{remark}\label{RemarkdeRhamStackTangentComplex}
Let $f\colon Y\to X$ be a map of Gelfand stacks admitting a cotangent complex and $i\geq 0$.   The de Rham stack of $\bb{T}_f$ is given by the object $(\bb{T}_f[i])^{\dR}\in \cat{GelfStk}_{/Y}$ sending a nilperfectoid ring $\eta\colon \GSpec(A)\to Y$ to 
\[
(\bb{T}_f[i])^{\dR}(A)=\tau_{\geq -i}\ob{Hom}_{A}(\eta^*\bb{L}_{f}, A^{u}) = \tau_{\geq -i}\ob{Hom}_{A^u}(\overline{\eta}^* \bb{L}_f, A^u)
\]
where $\overline{\eta}\colon \GSpec(A^u)\to Y$. We define $(\bb{T}_{f}[i])^{\dagger}:=\ob{fib}(\bb{T}_{f}[i]\to (\bb{T}_f[i])^{\dR})$. Thus, for a nilperfectoid affinoid Gelfand stack $\eta\colon \GSpec(A)\to Y$ one has 
\[
(\bb{T}_f[i])^{\dagger}(A)=\tau_{\geq -i}\ob{Hom}_{A}(\eta^*\bb{L}_{f}, \ob{Nil}^{\dagger}(A))
\]
and  there is a fiber sequence of  Gelfand stacks over $Y$
\[
(\bb{T}_f[i])^{\dagger}\to \bb{T}_f[i]\to (\bb{T}_{f}[i])^{\dR}.
\]
This sequence naturally upgrades to a fiber sequence of spectrum objects.
\end{remark}

The following list gives some basic  examples of Gelfand stacks admitting a deformation theory as in \Cref{sec:cotang-compl-morph-definition-having-a-cotangent-complex}. In the following we shall use the descent statements of \Cref{DescentDefTheoryTarget} and \Cref{PropFinalDescentDeformationTheory}.

In the following, we assume that $\tau$ is the $!$-formally \'etale topology from \Cref{sec:deform-theor-prop-shriek-formally-etale}.

\begin{example}
  \label{sec:cotang-compl-morph-examples-of-morphisms-admitting-a-cotangent-complex}
  \begin{enumerate}
  
\item By \Cref{LemmadeRhamDefProperty} any map between de Rham stacks admits a deformation theory and is formally \'etale.  
  
\item Any morphism $\GSpec(B)\to \GSpec(A)$ of affinoid Gelfand stacks admits a deformation theory, cf.\ \Cref{sec:deform-theor-prop-example-admitting-a-deformation-theory}; the nilcomplete and infinitesimal cohesivity is immediate from the definitions.
  The existence of cotangent complex is explained  in \cite[Remark 7.4.2]{solidnotes}.   
  
    \item Each morphism $Y\to X$ of Berkovich spaces admits a deformation theory. Indeed, applying analytic descent on the source  we can assume that $X$ is affinoid. Then, applying descent on the target $\dagger$-rational localizations, we can assume that $Y$ is itself affinoid which then admits a cotangent complex from  (2).

    \item More generally,  let $f\colon Y\to X$ which is \emph{weakly representable in Berkovich spaces} in the following sense\footnote{We caution the reader that there does not exist a good notion of an ``affine'' morphism in non-archimedean geometry, cf.\  \cite[Example 9.1.2]{scholze2020berkeley}.}: there exists a covering $X'\to X$ such that $Y\times_X X'$ is a  morphism of Berkovich spaces.   Then $f$ is infinitesimally cohesive and admits a cotangent complex, if in addition $X'\to X$ is an epimorphism in the $\dagger$-quasipro-\'etale topology, then $f$ is nilcomplete and thus admits a deformation theory.

\end{enumerate}

\end{example}

The vanishing of cotangent complexes for morphisms between perfectoid rings implies the following.

\begin{lemma}
  \label{sec:tang-cotang-spac-cotangent-complexes-of-diamonds}
  Let $Y\to X\in \ArcStk_{\Q_p}$ be a morphism of arc-stacks, which is weakly representable in arc-diamonds, i.e., there exists an arc-cover $X'\to X$ with $X'$ a perfectoid space such that the base change $Y\times_X X'$ admits a quasi-pro-\'etale surjection from a perfectoid space.
  Then $\widehat{Y}\to \widehat{X}$ is infinitesimally cohesive and admits a cotangent complex, in fact $\mathbb{L}_{\widehat{Y}/\widehat{X}}=0$. 
\end{lemma}
\begin{proof}
  By \Cref{DescentDefTheoryTarget} we may assume that $X$ is a perfectoid space such that $Y$ admits a quasi-pro-\'etale surjection $f\colon Z\to Y$ by a perfectoid space $Z$.
  We note that $f$ is weakly representable in perfectoid spaces because quasi-pro-\'etale spaces over a strictly totally disconnected spaces are represented by perfectoid spaces.
  In particular, $\widehat{f}$ is infinitesimally cohesive and admits a cotangent complex.
  From the vanishing of the cotangent complex for morphisms between perfectoid spaces (a standard calculation via $\mathbb{A}_{\inf}$ using \Cref{sec:cotang-compl-morph-cotangent-complex-via-p-completion}) we moreover see that $\widehat{f}$ is formally \'etale.
  By \Cref{PropFinalDescentDeformationTheory} we can deduce that $\widehat{Y}\to \widehat{X}$ admits a cotangent complex.
  The fiber sequence
  \[
    f^\ast \mathbb{L}_{\widehat{Y}/\widehat{X}}\to \mathbb{L}_{\widehat{Z}/\widehat{X}}\to \mathbb{L}_{\widehat{Z}/\widehat{Y}}
  \]
  and the vanishing of $\mathbb{L}_{\widehat{Z}/\widehat{Y}}, \mathbb{L}_{\widehat{Z}/\widehat{X}}$ with the fact that $\widehat{f}$ is a $!$-cover imply that $\mathbb{L}_{\widehat{Y}/\widehat{X}}=0$. 
\end{proof}
The following proposition  gives some mild conditions for classifying stacks to admit a deformation theory.  

\begin{proposition}\label{PropositionCotangentComplexesClassifyingStacks}
Let $S$ be a Gelfand stack and $G$ a group object in Gelfand stacks over $S$ such that the structural map $G\to S$ admits a deformation theory and is formally smooth.  

\begin{enumerate}

\item   Suppose that $BG\to S$ admits a deformation theory (this holds for example if $S\to BG$ is an epimorphism in the $\dagger$-quasi-pro-\'etale topology).
  Let $e\colon BG\to G/^{\ob{conj}}G$ denote the unit section where the right quotient is with respect to the conjugation action and  $f\colon  G/^{\ob{conj}}G\to BG$ the projection map. Then there is a natural equivalence of cotangent complexes
\[
\bb{L}_{BG/S}\cong  e^* \bb{L}_{f} [-1]. 
\]
In other words, $\bb{L}_{BG/S}$ is the $-1$-shift of the adjoint action of $G$ on the cotangent complex of $G$ at $1$.

\item Let $S$ be Berkovich space an and suppose that, locally on a strict closed cover by affinoid subspaces $\GSpec(R)=Z\subset S$, the group $G|_Z$ admits a $\dagger$-quasipro-\'etale group homomorphism $G|_Z\to H^{an}_Z$ to the analytification of an affine  smooth algebraic group $H_Z$ over $R(*)$. Then   $BG\to S$ admits a deformation theory.

\item Let $H\to S$ be a group such that $H\to H^{\dR/S}$ is an epimorphism and that $S\to BH^{\dR/S}$ is an epimorphism in the $\tau$-quasipro-\'etale topos.  Let $H^{\dagger}\subset H$ be the overconvergent neighbourhood of the unit section $e\colon S\to H$.   Then, if  $BH^{\dagger}\to S$ admits a deformation theory  so does $BH\to S$.  Furthermore,  let $k\geq 1$,  if $H$ is abelian, $B^nH^{\dagger}\to S$ admits a deformation theory,  and $S\to B^nH^{\dR}$ is a surjection in the  $\dagger$-quasipro-\'etale topology for all $1\leq n\leq k$, then  $B^nH\to S$ admits a deformation theory for all $1\leq n\leq k$.  

\item Let $S$ be a Gelfand stack, $V$ a vector bundle over $S$ and $H\to V$ a $\dagger$-quasipro-\'etale homomorphism of commutative groups over $S$. Then for all $n\geq 0$ the map $f\colon B^nH\to S$  admits a deformation theory. Furthermore, we have a natural isomorphism  of vector bundles. 
\[
\bb{L}_{f}= f^*V^{\vee}[-n]. 
\]

\item Let $S$ be a Gelfand stack and $H$ a commutative group stack over $S$  as in (4).  Let $X\to S$ be a Gelfand stack over $X$ and $\eta: X\to B^2 H$ a map depicting a gerbe $\widetilde{X}\to X$ banded by $H$. Then $\widetilde{X}\to X$ admits a deformation theory and there is a natural isomorphism of  sheaves in $\ob{D}_+(\widetilde{X})$
\[
\bb{L}_{\widetilde{X}/X}= V^{\vee}\otimes_{\s{O}_S} \s{O}_{\widetilde{X}} [-1]
\]
where $V^{\vee}=e^*\bb{L}_{H/S}$ is the fiber at $1$ of the cotangent complex of $H$. 
\end{enumerate}

\end{proposition}

\begin{proof}
In the following we will implicitly use the descent  of the deformation theory properties of \Cref{DescentDefTheoryTarget} and \Cref{PropFinalDescentDeformationTheory}.. 

\begin{enumerate}

\item   By assumption $BG\to S$ admits a deformation theory, we want to  compute $\bb{L}_{BG/S}$. Let $g\colon S\to BG$, and $l\colon S\to G$, it suffices to construct a natural equivalence 
\[
s^*\bb{L}_{BG/S}=l^*\bb{L}_{G/S}[-1].
\] 
Indeed, the statement will hold by applying the previous  construction to the classifying stack  of the group $G/^{\ob{conj}}G\to BG$ over $BG$ which is 
\[
BG/ (G/^{\ob{conj}}G) = BG\times BG,
\]
namely, the group $G/^{\ob{conj}}G$ is the \textit{inertia stack} which can be described as the loops of the diagonal map $BG\to BG\times BG$.
Note that the cotangent complex $\bb{L}_{BG\times BG/BG}$ of the second projection is the pullback of $\bb{L}_{BG/S}$ along the first projection, and pulling back along $s=\Delta\colon BG\to BG\times BG$ yields $\bb{L}_{BG/S}$. 

Let us now compute $s^*\bb{L}_{BG/S}$.  The diagram $S\to BG\to S$ produces a natural equivalence 
\[
e^*\bb{L}_{BG/S}[1]=\bb{L}_{S/BG}.
\] 
Now, since the map $S\to BG$ can be written as $G/G\to S/G$, we have that $\bb{L}_{S/BG}=l^* \bb{L}_{G/S}$ is  the fiber at $1$ of the cotangent complex of $G\to S$.  This yields the natural equivalence
\[
e^*\bb{L}_{BG/S}=l^* \bb{L}_{G/S} [-1]
\]
proving what we wanted.

\item All claims are local on the $\dagger$-quasipro-\'etale topology of $S$. Hence, we can assume without loss of generality that $S=\GSpec(R)$ is affine, and that  $G$ admits a $\dagger$-quasipro-\'etale morphism of groups to $H^{an}$ with $H$ a smooth algebraic group over $R(*)$. Since $H^{an}\to H^{\dR,an}$ is an epimorphism of Gelfand stacks as $H$ is smooth, we have a pullback diagram of Gelfand stacks
\[
\begin{tikzcd}
BG \ar[r] \ar[d]& BH^{\an} \ar[d] \\ 
BG^{\dR} \ar[r] & BH^{\an,\dR}.
\end{tikzcd}
\]
The lower bottom map admits a deformation theory being a morphism of de Rham stacks, this implies that $BG\to BH^{\an}$ admits a deformation theory. Hence, it suffices to show that $BH^{\an}$ admits a deformation theory. By the Tannakian  description of $H$-torsors (see \Cref{PropositionFunctorOfPointsAnalytification}), an $H$-torsor on a Gelfand affinoid stack is trivial locally in the \'etale topology as this is the case for schemes. This implies that $S\to BH^{\an}$ is an epimorphism in the Berkovich-\'etale topology, and in particular in the $\dagger$-quasipro-\'etale topology, and by (1) that $BH^{\an}\to S$ admits a deformation theory.
\item Since $H\to H^{\dR/S}$ is an epimorphism, we have a cartesian square of Gelfand stacks
  \[
  \begin{tikzcd}
 BH^\dagger  \ar[r] \ar[d] & BH \ar[d] \\ 
  S \ar[r] & BH^{\dR/S}
  \end{tikzcd}
  \] 
The bottom horizontal arrow   is surjection for the  $\dagger$-quasipro-\'etale topos, since it admits a deformation theory and is formally \'etale being a map of de Rham stacks, the same holds for the base change $BH^{\dagger}\to BH$. Thus, by descent of having a deformation theory on the source, to guarantee that $BH$ admits a deformation theory it suffices that $BH^{\dagger}$ does so, which is the case by assumption. 
  The case of the higher classifying stacks follow from the same argument as before and induction in $n$.

\item  We claim that for a vector bundle $V$ over $S$, the maps $S\to B^nV^{\dR/S}$ and $S\to B^nV$ are surjections in the $\dagger$-quasipro-\'etale topos.   Indeed, let $\eta\colon \GSpec(R)\to S$ a map from a Gelfand ring over $S$. For $B^nV$, maps $S\to B^nV$ are the same as cohomology classes $H^n(\GSpec(R), \eta^* V)$ which vanish for $n\geq 1$.  For $B^nV^{\dR/S}$, there is a $\dagger$-quasipro-\'etale cover $R\to R'$ with $R'$ strictly totally disconnected. Hence, after enlarging $R$, we have assume that $\Marc(R)$ is perfectoid. Thus, the maps $\GSpec(R)\to B^nV^{\dR/S}$ over $S$ are the same as the $v$-cohomology of vector bundles   $H^n(\Marc(R), \eta^*V\otimes_R \widehat{\s{O}})$ which vanishes for $n>0$ as $R^u$ is perfectoid. By an inductive argument using descent on the source, we deduce that $B^nV\to S$ admits a deformation theory

Now, let $H\to V$ be a $\dagger$-quasipro-\'etale group. We have a cartesian diagram
\[
\begin{tikzcd}
B^nH \ar[r] \ar[d] & B^nV \ar[d]\\
B^n H^{\dR/S} \ar[r] & B^nV^{\dR/S}. 
\end{tikzcd}
\]
The bottom line admits a deformation theory being a map of relative de Rham stacks. Since $B^nV\to S$ admits a deformation theory by the previous discussion, the map $B^nH\to S$ will also admit a deformation theory by stability under composition.

It is left to compute the cotangent complex of $B^nH$.  By (1) this is a shift of the adjoint action $e^*\bb{L}_{H/S}[-n]$ where $e\colon S\to H$ is the unit section. Since the group is abelian, this adjoint action is trivial (that is, for an abelian group $G$ one has a natural equivalence $G/^{\ob{conj}} G = G \times BG$, so that the action in the cotangent complex has a natural trivialization), and we get that $\bb{L}_{B^nH/S} = f^* e^*\bb{L}_{H/S}[-n]$ where $f\colon B^nH\to S$. It is left to compute $\bb{L}_{H/S}$, since $H\to V$ is $\dagger$-quasipro-\'etale, it is formally \'etale   and $\bb{L}_{H/S}=\bb{L}_{V/S}|_H$. Hence, we can assume that $H=V$, in that case it is classical (eg. by reducing to the universal situation over $S=B\GL_n$) that $e^*\bb{L}_{V/S} = V^{\vee}$ as a vector bundle over $S$. Part (4) follows.  

\item We have a pullback diagram
  \[
  \begin{tikzcd}
\widetilde{X}  \ar[r] \ar[d] & S \ar[d,"e"] \\ 
  X \ar[r, "\eta"] & B^2 H
  \end{tikzcd}
  \] 
  where the right vertical map is the unit section of $B^2 H$. From (4) above, we deduce that $\widetilde{X}\to X$ admits a deformation theory, in particular a cotangent complex.  Moreover, the cotangent complex of  $\widetilde{X} \to X$  is given by the pullback along $\widetilde{X} \to S$ of $\mathbb{L}_{S /B^2 H}$. We have a fiber sequence
  $$
  e^\ast \bb{L}_{B^2H/S} \to \bb{L}_{S/S} \to \bb{L}_{S/B^2 H}
  $$
so $\bb{L}_{S/B^2 H}=  e^\ast \bb{L}_{B^2H/S}[1]$. By (4) above, this is $V^{\vee}[-1]$, whence the result.
\qedhere
\end{enumerate}
\end{proof}

The following lemma shows that classifying stacks of smooth rigid groups admits a deformation theory. 

\begin{lemma}\label{LemmaRigidGroups}
Let $K$ be a Banach field over $\bb{Q}_p$ and $G$ a partially proper rigid analytic group over $K$. The following hold:
\begin{enumerate}

\item The map $\GSpec (K)\to BG^{\dR/K}$ is an epimorphism for the $\dagger$-quasipro-\'etale topology.

\item Let $G^{\dagger}\subset G$ be the overconvergent neighbourhood at $1$. Then  $BG^{\dagger}$ admits a deformation theory over $\Spec K$.

\item  The map $BG^{\dagger}\to BG$ is an epimorphism for the $\dagger$-quasipro-\'etale topology, in particular $BG$ admits a deformation theory.

\item Let $H\subset G^{\diamond}$ be a locally closed subgroup of the diamond of $G$, and let   $(H\subset G)^{\dagger}\subset G$ be the overconvergent neighbourhood at $H$. Then $B(H\subset G)^{\dagger}$ admits a deformation theory.

\end{enumerate}

In particular, the stacks $BG$ and $BG^{\dagger}$ admit a deformation theory over $\GSpec (K)$. 

\end{lemma}
\begin{proof}
The maps $\GSpec (K)\to BG$ and $\GSpec (K)\to BG^{\dagger}$ admit deformation theory  and are formally smooth being represented in locally closed subspaces of smooth Berkovich spaces. Since the maps are epimorphisms for the $\dagger$-quasipro-\'etale topology, the classifying stacks admit a deformation theory thanks to \Cref{PropFinalDescentDeformationTheory}.   Thus, it is left to show (1)-(4).

\begin{enumerate}

\item Let $\eta\colon \GSpec(R)\to BG^{\dR/K}$ be a torsor over a Gelfand ring $R$, we want to show that locally in the $\dagger$-quasipro-\'etale topology the torsor is trivial. We can assume without loss of generality that $R$ is strictly totally disconnected. Then the map $\eta$ is the same as an $G$-torsor on $\Marc(R)$ for the arc-topology (equivalently for the $v$ topology as $G$ is partially proper), and \cite[Lemma 4.10]{heuer2026g} guarantees that such torsors always have a section \'etale locally on $\Marc(R)$, i.e., a trivialization \'etale locally. 

\item  Let $\f{g}=\ob{Lie} (G)$ be the Lie algebra of $G$. We have a short exact sequence of Lie algebras $\f{g}^{\ob{rad}}\to \f{g}\to \f{g}^{\ob{ss}}$ where $\f{g}^{\ob{rad}}$ is the radical of the Lie algebra and $\f{g}^{\ob{ss}}$ its semisimplification. The Lie algebras have attached overconvergent groups (with group laws produced by the Baker-Campbell-Hausdorff formula).  In particular, if $G^{\dagger,\ob{rad}}\to G^{\dagger}\to G^{\dagger,\ob{ss}}$ is the short exact sequence of the  corresponding overconvergent groups, we have a pullback diagram 
\[
\begin{tikzcd}
BG^{\dagger, \ob{rad}} \ar[r] \ar[d] & BG^{\dagger} \ar[d]\\
\Spec K \ar[r] & BG^{\dagger,\ob{ss}}.
\end{tikzcd}
\]
Hence, it suffices to prove the claim when either $\f{g}$ is semisimple or solvable. In the solvable case, a further inductive argument reduces the problem to $\f{g}=K$ and therefore $G^{\dagger}=\bb{G}_a^{\dagger}$.  This case  is handled in  \Cref{PropositionCotangentComplexesClassifyingStacks} (2) and (4). We can then assume that $\f{g}$ is semisimple, and by extending $K$ to some finite extension that $\f{g}$ is the Lie algebra of a semisimple reductive group $\bb{G}$ over $K$. Then $G^{\dagger}\to \bb{G}^{an}$ is $\dagger$-quasipro-\'etale and this case is handled in \Cref{PropositionCotangentComplexesClassifyingStacks} (2).

\item By (1), the map $\Spec K\to BG^{\dR/K}$ is an epimorphism in the $\dagger$-quasipro-\'etale topology, and by base change the same holds for $BG^{\dagger}\to BG$. Hence, since $BG^{\dagger}$ admits a deformation theory and the previous map is formally \'etale, we deduce that $BG$ admits a deformation theory.

\item This follows from the fact that $B(H\subset G)^{\dagger}= BH^{\dR/K}\times_{BG^{\dR/K}} BG$, that maps between de Rham stacks admit a deformation theory, and that $BG\to \Spec K$ admits a deformation theory by part (3) above. \qedhere 
\end{enumerate}

\end{proof}

\begin{example}\label{ExampleCotangentComplexGla}
Let $G$ be a $p$-adic Lie group and let  $G^{\la}$ be the locally analytic incarnation of $G$.  Then $BG^{\la}\to \GSpec(\bb{Q}_p)$ admits a deformation theory and its cotangent complex is  given by the adjoint action on $\f{g}^{\vee}[-1]$. Indeed,  the formula for the cotangent complex is \Cref{PropositionCotangentComplexesClassifyingStacks} (1)  provided it exists.  To show that $BG^{\la}$ admits a deformation theory, we can argue locally \'etale on $BG^{sm}=(BG^{\la})^{\dR}$ and assume that $G$ is compact equal to the $\bb{Q}_p$-points of a partially proper rigid analytic group $\n{G}$. Then, we have a pullback square
\[
BG^{\la}=BG^{\sm}\times_{B\n{G}^{\dR}} B\n{G}.
\]
This implies that $BG^{\la}\to B\n{G}$ admits a deformation theory, and by \Cref{LemmaRigidGroups} (3)  $B\n{G}\to \GSpec(\Q_p)$ does as well. We deduce that $BG^{\la}\to \GSpec(\Q_p)$ admits a deformation theory as wanted. 
\end{example}

\bibliographystyle{alpha}
\bibliography{biblio.bib}

\newcommand{\etalchar}[1]{$^{#1}$}
\begin{thebibliography}{ABLB{\etalchar{+}}25}

\bibitem[ABLB{\etalchar{+}}25]{drff_stacks}
Johannes Ansch\"utz, Guido Bosco, Arthur-C\'esar Le~Bras, Juan~Esteban
  Rodr\^iguez~Camargo, and Peter Scholze.
\newblock Analytic de {R}ham stacks of {F}argues--{F}ontaine curves.
\newblock 2025.
\newblock arXiv:2510.15196.

\bibitem[AGT16]{abbes2016p}
Ahmed Abbes, Michel Gros, and Takeshi Tsuji.
\newblock {\em The $p$-adic {S}impson {C}orrespondence (AM-193)}, volume 193.
\newblock Princeton University Press, 2016.

\bibitem[AHLB22]{anschutz2022v}
Johannes Ansch{\"u}tz, Ben Heuer, and Arthur-C{\'e}sar Le~Bras.
\newblock v-vector bundles on $p$-adic fields and {S}en theory via the
  {H}odge-{T}ate stack.
\newblock {\em arXiv preprint arXiv:2211.08470}, 2022.

\bibitem[ALB20]{anschutz2020p}
Johannes Ansch{\"u}tz and Arthur-C{\'e}sar Le~Bras.
\newblock The $p$-completed cyclotomic trace in degree $2$.
\newblock {\em Annals of K-Theory}, 5(3):539--580, 2020.

\bibitem[ALBM26]{ALBM-FF}
Johannes Ansch{\"u}tz, Arthur-C{\'e}sar Le~Bras, and Lucas Mann.
\newblock A 6-functor formalism for solid quasi-coherent sheaves on the
  {Fargues}-{Fontaine} curve.
\newblock Preprint, {arXiv}:2412.20968 [math.{AG}] (2026), 2026.

\bibitem[AM24]{anschutz2024descent}
Johannes Ansch{\"u}tz and Lucas Mann.
\newblock Descent for solid quasi-coherent sheaves on perfectoid spaces.
\newblock {\em arXiv preprint arXiv:2403.01951}, 2024.

\bibitem[BC16]{BergerColmezII}
Laurent Berger and Pierre Colmez.
\newblock Th\'eorie de {S}en et vecteurs localement analytiques.
\newblock {\em Ann. Sci. \'Ec. Norm. Sup\'er. (4)}, 49(4):947--970, 2016.

\bibitem[Bha23]{bhatt2023algebraic}
Bhargav Bhatt.
\newblock Algebraic geometry in mixed characteristic.
\newblock In {\em International Congress of Mathematicians}, pages 712--748.
  European Mathematical Society-EMS-Publishing House GmbH, 2023.

\bibitem[Bha25]{bhatt-aspects}
Bhargav Bhatt.
\newblock Aspects of $p$-adic {H}odge theory.
\newblock 2025.

\bibitem[Bha26]{bhatt_takagi}
Bhargav Bhatt.
\newblock Non-abelian $p$-adic {H}odge theory and stacks.
\newblock {\em Takagi lectures, in preparation}, 2026.

\bibitem[BL22a]{bhatt2022absolute}
Bhargav Bhatt and Jacob Lurie.
\newblock Absolute prismatic cohomology.
\newblock {\em arXiv preprint arXiv:2201.06120}, 2022.

\bibitem[BL22b]{bhatt2022prismatic}
Bhargav Bhatt and Jacob Lurie.
\newblock Prismatic {F}-gauges.
\newblock {\em Lecture notes available at \url{https://www. math. ias. edu/\~{}
  bhatt/teaching/mat549f22/lectures. pdf}}, 2022.

\bibitem[BL22c]{bhatt2022prismatization}
Bhargav Bhatt and Jacob Lurie.
\newblock The prismatization of $ p $-adic formal schemes.
\newblock {\em arXiv preprint arXiv:2201.06124}, 2022.

\bibitem[BMS18]{Bhatt2018}
Bhargav Bhatt, Matthew Morrow, and Peter Scholze.
\newblock Integral $p$-adic {H}odge theory.
\newblock {\em Publications math{\'e}matiques de l'IH{\'E}S}, 128(1):219--397,
  2018.

\bibitem[BMS19]{bhatt2019topological}
Bhargav Bhatt, Matthew Morrow, and Peter Scholze.
\newblock Topological {H}ochschild homology and integral $p$-adic {H}odge
  theory.
\newblock {\em Publications math{\'e}matiques de l'IH{\'E}S}, 129(1):199--310,
  2019.

\bibitem[Bor12]{borel2012linear}
Armand Borel.
\newblock {\em Linear algebraic groups}.
\newblock Springer Science \& Business Media, 2012.

\bibitem[Bou75]{bourbaki1975lie}
Nicolas Bourbaki.
\newblock {\em Lie groups and Lie algebras: Chapters 1-3}, volume~1.
\newblock Springer Science \& Business Media, 1975.

\bibitem[Bro94]{MR1324339}
Kenneth~S. Brown.
\newblock {\em Cohomology of groups}, volume~87 of {\em Graduate Texts in
  Mathematics}.
\newblock Springer-Verlag, New York, 1994.
\newblock Corrected reprint of the 1982 original.

\bibitem[BS22]{bhatt_scholze_prisms_and_prismatic_cohomology}
Bhargav Bhatt and Peter Scholze.
\newblock Prisms and prismatic cohomology.
\newblock {\em Annals of Mathematics}, 196(3):1135--1275, 2022.

\bibitem[BZFN10]{BenZviFrancisNadlerIntegralTransforms}
David Ben-Zvi, John Francis, and David Nadler.
\newblock Integral transforms and {D}rinfeld centers in derived algebraic
  geometry, 2010.

\bibitem[BZSV24]{benzvi2024relativelanglandsduality}
David Ben-Zvi, Yiannis Sakellaridis, and Akshay Venkatesh.
\newblock Relative {Langlands} {Duality}.
\newblock Preprint, {arXiv}:2409.04677 [math.{RT}] (2024), 2024.

\bibitem[CC98]{CherbonnierColmez}
F.~Cherbonnier and P.~Colmez.
\newblock Repr\'esentations {$p$}-adiques surconvergentes.
\newblock {\em Invent. Math.}, 133(3):581--611, 1998.

\bibitem[CGN23]{colmez2023arithmetic}
Pierre Colmez, Sally Gilles, and Wies{\l}awa Nizio{\l}.
\newblock Arithmetic duality for $ p $-adic pro-{\'e}tale cohomology of
  analytic curves.
\newblock 2023.
\newblock Preprint. Available at \url{https://arxiv.org/abs/2308.07712}.

\bibitem[CGN25]{colmez2024duality}
Pierre Colmez, Sally Gilles, and Wies{\l}awa Nizio{\l}.
\newblock Duality for {$p$}-adic geometric pro-{\'e}tale cohomology.
\newblock Preprint, {arXiv}:2411.12163 [math.{AG}] (2025), 2025.

\bibitem[CLW{\etalchar{+}}26]{RHRobbaLiuZhu}
Zekun Chen, Ruochuan Liu, Yupeng Wang, Jiahong Yu, and Xinwen Zhu.
\newblock A {$p$}-adic {Riemann}-{Hilbert} correspondence over {Robba} rings.
\newblock Preprint, {arXiv}:2610.21065 [math.{AG}] (2026), 2026.

\bibitem[Col94]{colmez1994resultat}
Pierre Colmez.
\newblock Sur un r{\'e}sultat de {S}hankar {S}en.
\newblock {\em Comptes rendus de l'Acad{\'e}mie des sciences. S{\'e}rie 1,
  Math{\'e}matique}, 318(11):983--985, 1994.

\bibitem[CS23]{clausenscholze_analytic_stacks}
Dustin Clausen and Peter Scholze.
\newblock Analytic stacks.
\newblock {\em URl:
  https://people.mpim-bonn.mpg.de/scholze/AnalyticStacks.html}, 2023.

\bibitem[{\v{C}}S24]{cesnavicius2019purity}
K{\k{e}}stutis {\v{C}}esnavi{\v{c}}ius and Peter Scholze.
\newblock Purity for flat cohomology.
\newblock {\em Ann. Math. (2)}, 199(1):51--180, 2024.

\bibitem[Dri24a]{drinfeld20211}
Vladimir Drinfeld.
\newblock A 1-dimensional formal group over the prismatization of
  {$\mathrm{Spf}\,\mathbb{Z}_p$}.
\newblock {\em Pure Appl. Math. Q.}, 20(1):233--305, 2024.

\bibitem[Dri24b]{drinfeld2020prismatization}
Vladimir Drinfeld.
\newblock Prismatization.
\newblock {\em Sel. Math., New Ser.}, 30(3):150, 2024.
\newblock Id/No 49.

\bibitem[Dri26]{drinfeld2025ring}
Vladimir Drinfeld.
\newblock Ring stacks conjecturally related to the stacks {$BT_n^{G,\mu}$}.
\newblock Preprint, {arXiv}:2510.04958 [math.{AG}] (2026), 2026.

\bibitem[EY19]{ErtlYamadaRigidHK}
Veronika Ertl and Kazuki Yamada.
\newblock Rigid analytic reconstruction of {Hyodo}--{Kato} theory.
\newblock Preprint, {arXiv}:1907.10964 [math.{NT}] (2019), 2019.

\bibitem[FC90]{faltings_chai_degenerations_of_abelian_varieties}
Gerd Faltings and Ching-Li Chai.
\newblock {\em Degeneration of abelian varieties}, volume~22 of {\em Ergebnisse
  der Mathematik und ihrer Grenzgebiete (3) [Results in Mathematics and Related
  Areas (3)]}.
\newblock Springer-Verlag, Berlin, 1990.
\newblock With an appendix by David Mumford.

\bibitem[FS26]{fargues2021geometrization}
Laurent Fargues and Peter Scholze.
\newblock {\em Geometrization of the local {Langlands} correspondence (to
  appear)}, volume 466 of {\em Ast{\'e}risque}.
\newblock Paris: Soci{\'e}t{\'e} Math{\'e}matique de France (SMF), 2026.

\bibitem[GL21]{guo2021period}
Haoyang Guo and Shizhang Li.
\newblock Period sheaves via derived de {R}ham cohomology.
\newblock {\em Compositio Mathematica}, 157(11):2377--2406, 2021.

\bibitem[GR03]{gabber2003almost}
Ofer Gabber and Lorenzo Ramero.
\newblock {\em Almost ring theory}, volume 1800 of {\em Lect. Notes Math.}
\newblock Berlin: Springer, 2003.

\bibitem[Gro68]{grothendieck1968crystals}
Alexander Grothendieck.
\newblock Crystals and the de {R}ham cohomology of schemes.
\newblock {\em Dix expos{\'e}s sur la cohomologie des sch{\'e}mas}, 3:306--358,
  1968.

\bibitem[GW26]{gao2026galois}
Hui Gao and Yupeng Wang.
\newblock Galois representations over convergent de {R}ham period ring.
\newblock {\em arXiv preprint arXiv:2604.21605}, 2026.

\bibitem[Hau26]{hauck}
Maximilian Hauck.
\newblock Rational analytic syntomic cohomology.
\newblock Preprint, {arXiv}:2604.15193 [math.{AG}] (2026), 2026.

\bibitem[Heu25]{heuer2025p}
Ben Heuer.
\newblock A $p$-adic {S}impson correspondence for smooth proper rigid
  varieties.
\newblock {\em Inventiones mathematicae}, 240(1):261--312, 2025.

\bibitem[Heu26a]{heuer2026g}
Ben Heuer.
\newblock {$ G $-torsors on perfectoid spaces}.
\newblock {\em {\'E}pijournal de G{\'e}om{\'e}trie Alg{\'e}brique}, 10, 2026.

\bibitem[Heu26b]{MR5022120}
Ben Heuer.
\newblock Moduli spaces in {$p$}-adic non-abelian {H}odge theory.
\newblock {\em J. Algebraic Geom.}, 35(2):225--318, 2026.

\bibitem[HKR25]{hekking2025deformationnormalbundleblowups}
Jeroen Hekking, Adeel~A. Khan, and David Rydh.
\newblock Deformation to the normal bundle and blow-ups via derived {Weil}
  restrictions.
\newblock Preprint, {arXiv}:2511.19412 [math.{AG}] (2025), 2025.

\bibitem[HM24]{heyer20246functorformalismssmoothrepresentations}
Claudius Heyer and Lucas Mann.
\newblock 6-{Functor} {Formalisms} and {Smooth} {Representations}.
\newblock Preprint, {arXiv}:2410.13038 [math.{CT}] (2024), 2024.

\bibitem[How25]{howe2025geometric}
Sean Howe.
\newblock The geometric {S}en morphism is the unique lift of the
  {K}odaira--{S}pencer morphism.
\newblock {\em arXiv preprint arXiv:2512.09669}, 2025.

\bibitem[How26]{howe2025inscriptiontwistorspadicperiods}
Sean Howe.
\newblock Inscription, twistors, and {$p$}-adic periods, 2026.

\bibitem[Jia26]{jiang2026differential}
Yuanyang Jiang.
\newblock Differential operators on locally analytic {S}himura varieties.
\newblock {\em arXiv preprint arXiv:2604.09116}, 2026.

\bibitem[KL13]{kedlaya_liu_relative_p_adic_hodge_theory_foundations}
K.~S. {Kedlaya} and R.~{Liu}.
\newblock {Relative $p$-adic Hodge theory: Foundations}.
\newblock {\em ArXiv e-prints}, January 2013.

\bibitem[Koh11]{MR2774315}
Jan Kohlhaase.
\newblock The cohomology of locally analytic representations.
\newblock {\em J. Reine Angew. Math.}, 651:187--240, 2011.

\bibitem[Kuh24]{kuhn2024atiyah}
Nikolas Kuhn.
\newblock The {A}tiyah class on algebraic stacks.
\newblock In {\em Forum of Mathematics, Sigma}, volume~12, page e100. Cambridge
  University Press, 2024.

\bibitem[Laz65]{zbMATH03225421}
Michel Lazard.
\newblock {{\(p\)}}-adic analytic groups.
\newblock {\em Publ. Math., Inst. Hautes {\'E}tud. Sci.}, 26:389--603, 1965.

\bibitem[Li25]{zhenghui_li_2024}
Zhenghui Li.
\newblock Duality for {Arithmetic} {$p$}-adic {Pro}-{\'e}tale {Cohomology} of
  {Analytic} {Spaces}.
\newblock Preprint, {arXiv}:2412.11786 [math.{AG}] (2025), 2025.

\bibitem[Lur11]{lurie_dag_IX}
Jacob Lurie.
\newblock Derived {A}lgebraic {G}eometry {IX}: {C}losed {I}mmersions.
\newblock Available at \url{https://www.math.ias.edu/~lurie/}, 2011.

\bibitem[Lur17]{lurie_higher_algebra}
Jacob Lurie.
\newblock Higher algebra.
\newblock Available at \url{http://www.math.harvard.edu/~lurie/}, 2017.

\bibitem[Lur18]{Luriea}
Jacob Lurie.
\newblock Spectral algebraic geometry.
\newblock Available at \url{http://www.math.harvard.edu/~lurie/}, 2018.

\bibitem[Man22]{mann2022p}
Lucas Mann.
\newblock A $ p $-adic 6-functor formalism in rigid-analytic geometry.
\newblock {\em arXiv preprint arXiv:2206.02022}, 2022.

\bibitem[Man26]{mann2022zl}
Lucas Mann.
\newblock A 6-functor formalism for {{\(\mathbb{Z}_\ell\)}}- and
  {{\(\mathbb{Q}_\ell\)}}-sheaves on diamonds.
\newblock {\em Res. Math. Sci.}, 13(1):72, 2026.
\newblock Id/No 19.

\bibitem[Mat16]{mathew2016galois}
Akhil Mathew.
\newblock The {G}alois group of a stable homotopy theory.
\newblock {\em Advances in Mathematics}, 291:403--541, 2016.

\bibitem[Mer72]{MR330167}
David Meredith.
\newblock Weak formal schemes.
\newblock {\em Nagoya Math. J.}, 45:1--38, 1972.

\bibitem[Mik26]{mikami2025finiteness}
Yutaro Mikami.
\newblock Finiteness and duality of cohomology of {$(\phi,\Gamma)$}-modules and
  the 6-functor formalism of locally analytic representations.
\newblock Preprint, {arXiv}:2504.01780 [math.{NT}] (2026), 2026.

\bibitem[Pan22]{MR4390302}
Lue Pan.
\newblock On locally analytic vectors of the completed cohomology of modular
  curves.
\newblock {\em Forum Math. Pi}, 10:Paper No. e7, 82, 2022.

\bibitem[Pan26]{pan2022locally}
Lue Pan.
\newblock On locally analytic vectors of the completed cohomology of modular
  curves. {II}.
\newblock {\em Ann. Math. (2)}, 203(1):121--281, 2026.

\bibitem[Por24]{PoratLocAnVetorsFF}
Gal Porat.
\newblock Locally analytic vector bundles on the {F}argues-{F}ontaine curve.
\newblock {\em Algebra Number Theory}, 18(5):899--946, 2024.

\bibitem[Por25]{PoratSolidLocAn}
Gal Porat.
\newblock Solid locally analytic representations in mixed characteristic, 2025.

\bibitem[PY20]{MR4176782}
Mauro Porta and Tony~Yue Yu.
\newblock Representability theorem in derived analytic geometry.
\newblock {\em J. Eur. Math. Soc. (JEMS)}, 22(12):3867--3951, 2020.

\bibitem[Qin26]{qin2026arithmetic}
Junhui Qin.
\newblock The arithmetic de {R}ham stack.
\newblock {\em arXiv preprint arXiv:2608.20203}, 2026.

\bibitem[Rak26]{raksit2020hochschild}
Arpon Raksit.
\newblock Hochschild homology and the derived de {Rham} complex revisited.
\newblock {\em Sel. Math., New Ser.}, 32(3):96, 2026.
\newblock Id/No 61.

\bibitem[Ram12]{ramsey2009bsen}
Nick Ramsey.
\newblock {{\(B _{\text{Sen}}\)}} via distributions on weight space.
\newblock {\em Ramanujan J.}, 28(1):139--153, 2012.

\bibitem[RC24]{RCAndeRham}
Juan~Esteban Rodr\'iguez~Camargo.
\newblock The analytic de {Rham} stack in rigid geometry.
\newblock Preprint, {arXiv}:2401.07738 [math.{AG}] (2024), 2024.

\bibitem[RC26a]{RCCartierGerbe}
Juan~Esteban Rodr\'iguez~Camargo.
\newblock Cartier duality for gerbes of vector bundles.
\newblock Preprint, {arXiv}:2512.24967 [math.{AG}] (2026), 2026.

\bibitem[RC26b]{camargo2022geometricsen}
Juan~Esteban Rodr\'iguez~Camargo.
\newblock Geometric {S}en theory over rigid analytic spaces.
\newblock {\em J. Eur. Math. Soc. (JEMS)}, 28(8):3385--3451, 2026.

\bibitem[RC26c]{RCBBCourse2026}
Juan~Esteban Rodr\'iguez~Camargo.
\newblock Locally analytic representations and geometric representation theory.
\newblock {\em URl:
  https://sites.google.com/view/jerodriguezcamargo/courses-and-seminars?authuser=0},
  2026.

\bibitem[RC26d]{solidnotes}
Juan~Esteban Rodr\'iguez~Camargo.
\newblock Notes on {Solid} {Geometry}.
\newblock Preprint, {arXiv}:2603.03012 [math.{AG}] (2026), 2026.

\bibitem[RJRC22]{jacinto2021solid}
Joaqu\'in Rodrigues~Jacinto and Juan~Esteban Rodr\'iguez~Camargo.
\newblock Solid locally analytic representations of {$p$}-adic {L}ie groups.
\newblock {\em Represent. Theory}, 26:962--1024, 2022.

\bibitem[RJRC26]{jacinto2023solid}
Joaqu{\'\i}n Rodrigues~Jacinto and Juan~Esteban Rodr{\'\i}guez~Camargo.
\newblock Solid locally analytic representations.
\newblock Preprint, {arXiv}:2305.03162 [math.{RT}] (2026), 2026.

\bibitem[Sch13a]{scholze2013p}
Peter Scholze.
\newblock $ p $-adic hodge theory for rigid-analytic varieties.
\newblock In {\em Forum of Mathematics, Pi}, volume~1. Cambridge University
  Press, 2013.

\bibitem[Sch13b]{scholze2013perfectoid}
Peter Scholze.
\newblock Perfectoid spaces: a survey.
\newblock In {\em Current developments in mathematics 2012. Papers based on
  selected lectures given at the conference, Harvard University, Cambridge, MA,
  USA, November 2012}, pages 193--227. Somerville, MA: International Press,
  2013.

\bibitem[Sch18]{ScholzeLubinTate}
Peter Scholze.
\newblock On the {{\(p\)}}-adic cohomology of the {Lubin}-{Tate} tower.
\newblock {\em Ann. Sci. {\'E}c. Norm. Sup{\'e}r. (4)}, 51(4):811--863, 2018.

\bibitem[Sch24]{real_local_langlands}
Peter Scholze.
\newblock Geometrization of the real local {L}anglands correspondence.
\newblock Avaliable at
  \url{https://people.mpim-bonn.mpg.de/scholze/RealLocalLanglands.pdf}, 2024.

\bibitem[Sch25]{scholze2025berkovich}
Peter Scholze.
\newblock Berkovich motives.
\newblock {\em Journal of the American Mathematical Society}, 2025.

\bibitem[Sch26a]{scholze_etale_cohomology_of_diamonds}
Peter Scholze.
\newblock Etale cohomology of diamonds.
\newblock Preprint, {arXiv}:1709.07343 [math.{AG}] (2026), 2026.

\bibitem[Sch26b]{GestaltenScholze}
Peter Scholze.
\newblock Geometry and higher category theory.
\newblock \url{https://people.mpim-bonn.mpg.de/scholze/Gestalten.pdf}, 2026.

\bibitem[Sch26c]{scholze2022lectures}
Peter Scholze.
\newblock Six-{Functor} {Formalisms}.
\newblock Preprint, {arXiv}:2510.26269 [math.{AG}] (2026), 2026.

\bibitem[Sen73]{zbMATH03407789}
Shankar Sen.
\newblock Lie algebras of {Galois} groups arising from {Hodge}-{Tate} modules.
\newblock {\em Ann. Math. (2)}, 97:160--170, 1973.

\bibitem[Sen81]{SenContinuousCoho}
Shankar Sen.
\newblock Continuous cohomology and {$p$}-adic {G}alois representations.
\newblock {\em Invent. Math.}, 62(1):89--116, 1980/81.

\bibitem[Ser92]{MR1176100}
Jean-Pierre Serre.
\newblock {\em Lie algebras and {L}ie groups}, volume 1500 of {\em Lecture
  Notes in Mathematics}.
\newblock Springer-Verlag, Berlin, second edition, 1992.
\newblock 1964 lectures given at Harvard University.

\bibitem[Sim97]{simpson1996hodge}
Carlos Simpson.
\newblock The {Hodge} filtration on nonabelian cohomology.
\newblock In {\em Algebraic geometry. Proceedings of the Summer Research
  Institute, Santa Cruz, CA, USA, July 9--29, 1995}, pages 217--281.
  Providence, RI: American Mathematical Society, 1997.

\bibitem[ST03]{zbMATH01965478}
Peter Schneider and Jeremy Teitelbaum.
\newblock Algebras of {{\(p\)}}-adic distributions and admissible
  representations.
\newblock {\em Invent. Math.}, 153(1):145--196, 2003.

\bibitem[{Sta}17]{StacksProjectAuthors2017}
The {Stacks Project Authors}.
\newblock \textit{Stacks Project}.
\newblock \url{http://stacks.math.columbia.edu}, 2017.

\bibitem[Ste23]{stefanich2023tannakaduality1affineness}
Germ{\'a}n Stefanich.
\newblock Tannaka duality and 1-affineness.
\newblock Preprint, {arXiv}:2311.04515 [math.{AG}] (2023), 2023.

\bibitem[SW20]{scholze2020berkeley}
Peter Scholze and Jared Weinstein.
\newblock Berkeley lectures on {$p$}-adic geometry.
\newblock In {\em Berkeley Lectures on p-adic Geometry}. Princeton University
  Press, 2020.

\bibitem[Ter26]{terentiuk2026ogus}
Gleb Terentiuk.
\newblock Ogus--{V}ologodsky equivalence via stacks.
\newblock {\em arXiv preprint arXiv:2604.04317}, 2026.

\bibitem[Wan23]{wang2023p}
Yupeng Wang.
\newblock A p-adic {S}impson correspondence for rigid analytic varieties.
\newblock {\em Algebra \& Number Theory}, 17(8):1453--1499, 2023.

\bibitem[Wan26]{wang2026relativegagatheoremapplication}
Qixiang Wang.
\newblock The relative {GAGA} {Theorem} and an application to the analytic
  mapping stacks.
\newblock Preprint, {arXiv}:2601.12299 [math.{AG}] (2026), 2026.

\bibitem[Wie23]{wiersig2023solution}
Finn Wiersig.
\newblock On a solution functor for ${D}$-modules via $p$-adic {H}odge theory.
\newblock {\em arXiv preprint arXiv:2309.13769}, 2023.

\bibitem[Wie26a]{WiersigFullyFaithfulRH}
Finn Wiersig.
\newblock A fully faithful p-adic {Riemann}-{Hilbert} functor for coadmissible
  {D}-cap-modules.
\newblock Preprint, {arXiv}:2506.12601 [math.{AG}] (2026), 2026.

\bibitem[Wie26b]{WiersigProetalePeriodSheaves}
Finn Wiersig.
\newblock Galois and {Pro}-{\'e}tale {Cohomology} of {Overconvergent} de {Rham}
  {Period} {Rings}.
\newblock Preprint, {arXiv}:2309.13769 [math.{NT}] (2026), 2026.

\bibitem[Zav26]{zavyalov2023poincaredualityabstract6functor}
Bogdan Zavyalov.
\newblock Poincar{\'e} {Duality} in abstract 6-functor formalisms.
\newblock Preprint, {arXiv}:2301.03821 [math.{AG}] (2026), 2026.

\end{thebibliography}

\end{document}